# Extremes and Recurrence in Dynamical Systems

### Valerio Lucarini

Institute of Meteorology, CEN, University of Hamburg, Hamburg, Germany
Department of Mathematics and Statistics, University of Reading, Reading, UK

### Davide Faranda

LSCE, CEA CNRS UMR 8212, Gif-sur-Yvette France

### Ana Cristina Moreira Freitas

Centro de Matemática da Universidade do Porto, Porto, Portugal
Faculdade de Economia da Universidade do Porto, Porto, Portugal

### Jorge Milhazes Freitas

Centro de Matemática da Universidade do Porto, Porto, Portugal

### Mark Holland

School of Engineering, Computer Science & Mathematics, University of Exeter, Exeter, UK

### Tobias Kuna

Department of Mathematics and Statistics, University of Reading, Reading, UK

### Matthew Nicol

Department of Mathematics, University of Houston, Houston, USA

### Mike Todd

Mathematical Institute, University of St. Andrews, St. Andrews, UK

### Sandro Vaienti

Aix Marseille Université, CNRS, CPT, UMR 7332, Marseille, France
Université de Toulon, CNRS, CPT, UMR 7332, Toulon, France

*Ma l'impresa eccezionale,*
*dammi retta,*
*è essere normale*

*(L. Dalla, 1977)*

**4**

# Contents

**vi**

!

**2**

# 1
# Introduction

## 1.1
## A Transdisciplinary Research Area

The study of extreme events has long been a very relevant field of investigation at the intersection of different fields, most notably mathematics, geosciences, engineering and finance [1, 2, 3, 4, 5, 6, 7]. While extreme events in a given physical system obey the same laws as typical events, extreme events are rather special both from a mathematical point of view and in terms of their impacts. Often, procedures like mode reduction techniques, which are able to reliably reproduce the typical behaviour of a system, do not perform well in representing accurately extreme events and underestimate their variety. It is extremely challenging to predict extremes in the sense of defining precursors for specific events and, on longer time scales, to assess how modulations in the external factors (e.g. climate change in the case of geophysical extremes) impact on their properties.

Clearly, understanding the properties of the *tail of the probability distribution* of a stochastic variable attracts a lot of interest in many sectors of science and technology because extremes sometimes relate to situations of high stress or serious hazard, so that in many fields it is crucial to be able to predict their return times in order to cushion and gauge risks, such in the case of the construction industry, of the energy sector, of agriculture, of territorial planning, of logistics, and of financial markets, just to name a few examples. Intuitively, we associate the idea of an *extreme event* to either something which is very large, or something which is very rare, or, in more practical terms, to something with a rather abnormal impact with respect to an indicator (*e.g.* economic or environmental welfare) that we deem important. While overlaps definitely exist between such definitions, they are not equivalent.

An element of subjectivity is unavoidable when treating finite data - observational or synthetic - and when having a specific problem in mind: we might be interested in studying yearly or decadal temperature maxima in a given location, or the return period of river discharge larger than a prescribed value. Practical needs have indeed been crucial in stimulating the investigation of extremes, and most notably in the fields of hydrology [5] and of finance [2], which provided the first examples of empirical yet extremely powerful approaches.

Let's briefly consider, to take a relevant and instructive example, the case of geophysical extremes, which do not only cost many human lives each year, but also cause significant economic damages [4, 8, 9, 10]; see also the discussion and historical perspective given in [11]. For instance, freak ocean waves are extremely hard to predict and can have devastating impacts on vessels and coastal areas [12, 13, 14]. Windstorms are well-known to dominate the list of the costliest natural disasters, with many occurrences of individual events causing insured losses topping 1 Billion $ [15, 16]. Temperature extremes, like heat waves and cold spells, have severe impacts on society and ecosystems [17, 18, 19]. Notable temperature-related extreme events are the 2010 Russian heat wave, which caused 500 wild fires around Moscow, reduced grain harvest by 30% and was the hottest summer in at least 500 years [20], and the 2003 heat wave in Europe, which constituted the second hottest summer in this period [21]. The 2003 heat wave had significant societal consequences; *e.g.* it caused additional deaths exceeding 70000 [17]. On the other hand, recent European winters were very cold, with widespread cold spell hitting Europe during January 2008, December 2009 and January 2010. The increasing number of weather and climate extremes over the last few decades [22, 23, 24] has led to intense debates, not only amongst scientists but also policy makers and the general public, whether this increase is triggered by global warming.

Additionally, is some cases, we might be interested in exploring the spatial correlation of extreme events. See extensive discussion in [25, 26]. Intense precipitation events occurring at the same time within a river basin, which acts as spatial integrator of precipitations, can cause extremely dangerous floods. Large scale long-lasting droughts can require huge infrastructural investments to guarantee the welfare of entire populations and the production of agricultural goods. Extended wind storms can halt the production of wind energy in vast territories, changing dramatically the input of energy into the electric grid, with the ensuing potential risk of brown- or blackouts, or can seriously impact the air, land, and sea transportation networks. In general, weather and climate models need to resort to parametrizations for representing the effect of small-scale processes on the large scale dynamics. Such parametrizations are usually constructed and tuned in order to capture as accurately as possible the first moments (mean, variability) of the large scale climatic features. But it is indeed much less clear how spatially extended extremes could be affected. Going back to a more conceptual problem, one can consider the case where we have two or more versions of the same numerical model of a fluid, differing for the adopted spatial resolution. How can we compare the extremes of a local physical observable provided by the various versions of the model? Is there a coarse-graining procedure suited for upscaling to a common resolution the outputs of the models, such that we find a coherent representation of the extremes? At this regard, see in [27] a related analysis of extremes of precipitation in climate models.

When we talk about the impacts of geophysical extremes, a complex portfolio of aspects needs to be considered, so the study of extremes leads naturally to comprehensive transdisciplinary areas of research. The impacts of geohazards depend strongly not only on the magnitude of the extreme event, but also on the vulnerability of the affected communities. Some areas, *e.g.* coasts, are especially at risk of high-impact geophysical hazards, such as extreme floods caused by tsunami, storm surges and freak waves. Delta regions of rivers face additional risks due to flooding resulting from intense and extensive precipitation events happening upstream the river basin, maybe at thousands of Kms of distance. Sometimes, storm surges and excess precipitation act in synergy and create exceptional coastal flooding. Mountain areas are in turn, extremely sensitive to flash floods, landslides, and extreme solid and liquid precipitation events.

When considering the impacts of extreme events on the societal fabric, a primary role is played by the level of resilience and preparedness of the affected local communities. Such levels can vary enormously, depending on many factors including availability of technology, social structure, level of education, quality of public services, presence of social and political tensions, including conflicts, gender relations, and many others [28, 29, 30]. Geophysical extremes can wipe out or damage substantially the livelihood of entire communities, leading in some cases to societal breakdown and mass migration, as, *e.g.*, in the case of intense and persistent droughts. Prolonged and extreme climate fluctuations are nowadays deemed responsible for causing or accelerating the decline of civilisations - *e.g.* the rapid collapse of the Mayan empire in the XI century, apparently fostered by an extreme multidecadal drought event [31]. Cold spells can also have severe consequences. An important yet not so well known example is given by the dramatic impacts of the recurrent ultra cold winter *Dzud* events in the Mongolian plains, which can lead to the death of livestock due to starvation, and are deemed responsible for causing in the past the recurrent waves of migration of nomadic Mongolian populations and their clash with China, Central Asia, and Europe [32, 33]. The meteorological conditions and drivers of Dzud events are basically uninvestigated.

Nowadays public and private decision makers need support under great uncertainty from science and technology to optimally address how to deal with forecasts of extreme events in order to address questions such as: how to evacuate a coastal region forecasted to be flooded as a result of a storm surge; and how to plan for successive severe winter conditions affecting Europe's transportation networks? How to minimize the risk of drought-induced collapse in the availability of staple food in Africa? How to adapt to climate change? Along these lines, today, a crucial part of advising local and national governments is not only the prediction of natural hazards, but also the communication of risk for a variety of public and private stakeholders, as, e.g., in the health, energy, food, transport, and logistics sectors [23].

Other sectors of public and private interest where extremes play an important role are finance and (re-)insurance. Understanding and predicting events like the crash of the New York Stocks Exchange of October 1987 and the Asian crisis have become extremely important for investors and institutions. The ability to assess the very high quantiles of a probability distribution, and so delve into low-probability events is of great interest, because it translates into the ability to deal efficiently with extreme financial risks, as in the case of currency crises, stock market crashes, and large bond defaults, and, in the case of insurance, of low probability/high risk events [2].

The standard way to implement risk-management strategies has been, until recently, based on the value-at-risk (VaR) approach [34]. The VaR approach typically aims at estimating the worst anticipated loss over a given period with given probability assuming regular market conditions. The basic idea is to extract information from the typical events to predict the properties of the tails of the probability distribution. The VaR method has been recently criticised because of various limitations. In many applications, simple normal statistical models are used, resulting into an underestimation of the risk associated to the high quantiles in the common (in the financial sector) case where fat-tailed distributions are present. More sophisticated statistical models partially address this problem, but, since they are based on fitting distributions like the Student-$t$ or mixtures of normals, the properties of the non-typical events are poorly constrained. Non-parametric methods, instead, cannot be used beyond the range of observed values, so that it is virtually impossible to have any predictive power for assessing the probability of occurrence of out-of-sample events [35].

Intuitively, it seems impossible to be able to predict the probability of occurrence of events larger than what has already been observed, and, in general, of events which are extremely untypical. The key idea is to focus on the tail of the distribution, by constructing a statistical model aimed at describing only the extreme data, so that the fitting procedure is tailored to what one is interested in [36, 37]. In other terms, this requires separating clearly typical events from non-typical - extreme - events, disregarding entirely the former, and attributing to the latter *special* features. One needs to note that, by the very nature of the procedure, introducing spurious events in the group of selected extremes (the so-called *soft extremes*) may lead to substantial biases in the procedure of statistical fitting.

The goal of this introduction is to motivate the reader to delve into the mathematics of extremes by presenting some of the most interesting challenges in various areas of scientific research where extremes play a major role. In this sense, we stick to the history of the field, where mathematical findings and relevant applications have gone hand in hand since several decades. In the following sections we introduce the main themes of this book, we try to clarify the main ideas we will develop, and we underline the most problematic aspects related to the development of a rigorous theory of extremes for dynamical systems, and to its possible use for studying specific mathematical and more applied problems.

## 1.2
## Some Mathematical Ideas

One can safely say that in the case of extremes, as in many other sectors of knowledge, the stimuli leading to the mathematical theory have come from the need to systematize and understand current technological applications at a deeper level. A more complete and rigorous mathematical framework is also instrumental in defining more powerful and more robust algorithms and approaches to studying time series and to improve our ability to describe, study, and predict extremes, and, eventually, cushion their impacts. This book aims at providing a mathematical point of view on extremes able to relate their features to the dynamics of the underlying system.

Obviously, in order to develop a mathematical theory of extremes, it is necessary to carefully define the *rules of the game*, *i.e.* lay out criteria separating clearly extremes from regular events. Moreover, it is crucial to construct a conceptual framework which can be easily adapted to many different applications while, at the same time, is able to deliver as many results as possible which are universal in some suitable limits.

The quest for some level of universality, apart from being of obvious mathematical interest, is very relevant when using exact mathematical results for interpreting data coming from observations or from numerical simulations. In fact, the presence of universal mathematical properties gives more robustness to the procedures of statistical inference. It is clear that the tantalizing goal of constructive credible estimates for very high quantiles and for the probability of occurrence of yet unobserved events requires one to provide arguments to define the properties for the tails of distribution under very general conditions, and, possibly, of an explicit universal functional form describing them.

The classical construction of the mathematical theory of extremes and the definition of extreme value laws (EVLs) result from the generalization of the intuitive points of view - extreme as large and extreme as rare - introduced before. Following [38], one considers a random variable (r.v.) $X$ and investigates the conditions under which one can construct the properties of the r.v. $M_N$ given by the maximum of a sequence of $N$ independent and identically distributed (i.i.d.) r.v. $X_j, j = 1, \ldots, N$, in the limit $N \to \infty$. This is an extremely powerful and fruitful approach to the problem, which we will later discuss in detail. One can find that, under rather general conditions and using a suitable procedure of rescaling, it is possible to construct such a limit law for $M_N$. In practice, one finds a general three-parameter statistical model (Generalized Extreme Value or GEV distribution) for fitting the empirical probability distribution of the *block maxima* (BM), which are constructed from a time series of length $M = NK$ by chopping into $K$ (long) blocks of length $N$, and taking the maximum of each block.

The GEV distribution provides a generalization of the Frechét, Gumbel, and Weibull distributions, which have long been used for studying extremes events in many fields of science and technology. Nowadays, *GEV*-fitting is one of the most common methods for dealing with extremes in many applications. Giving a meaning to the statistics of *e.g* the annual maxima of surface temperature in a given location requires taking implicitly or explicitly the block maxima point of view [39]. The sign of the *shape parameter* determines the qualitative properties of the extremes, If the shape parameter is positive (Frechét distribution), the extremes are unlimited, if the sign is positive (Weibull distribution) the extremes are upper limited, with the Gumbel distribution (vanishing shape parameter) being the limiting case lying in-between. The *location parameter*, instead, describes how large, typically, extremes are, while the *scale parameter* indicates the variability of the extremes.

A crucial aspect is that under the same mathematical conditions allowing for the definition of the limit laws for the variables $M_N$, it is possible to look at the problem from a complementary point of view. One finds a one-to-one correspondence between the statistical properties of the block maxima in the limit of very large $N$, and those of the *above-threshold* events, described by the conditional probability of $X$ given that $X$ itself exceeds a certain threshold value $T$, for very large $T$, approaching infinity if the distribution of $X$ is not bounded from above, or approaching the maximum of the support of the probability distribution of $X$ in the opposite case. One can prove that *these* maxima are distributed according to the two-parameter Generalized Pareto Distribution (GPD) [40, 41, 1]. This point of view leads to looking at extremes using the *Peak-Over-Threshold* (POT) method. When one looks at the risk of occurrence of negative anomalies of input of wind energy into the electric grid larger than a given alert level $T$, the POT point of view is taken [42]. When performing POT statistical inference, one derives the values of the two parameters of the corresponding GPD, the *shape* and *scale* parameters, with a similar meaning as for the GEV case. See [43] for an elegant and comprehensive discussion of advanced use of the POT methods for the relevant case of rainfall data.

It is remarkable that while for a given series of i.i.d. r.v. $X_j$ the probability distributions of POT events and of the BM are different, also in the limit, because the procedure of selection of extremes is fundamentally different, the two points of view are in some sense equivalent. In other terms we have universal properties, which do not depend on the procedure of selection of the extremes. More specifically, if one learns the properties of extremes defined as (large) block maxima, the properties of extremes as events above a (high) threshold can be automatically deduced, and the other way around. In fact, the shape parameters of the GEV and GPD are the same, while one can write explicit relationships between linking the GEV's scale and location parameters on one side, with the GPD scale parameter and threshold $T$ (which acts as a hidden parameter) on the other side, so that a one-to-one correspondence between the two distributions can be found [44]. There has long been a very lively debate on whether the BM or POT method is better suited for treating finite time series coming from social, engineering, or natural systems. Most importantly, the existence of the corresponding well-defined and universal parametric probability distributions implies that if we are able to obtain a robust fit for the extremes of an observable given a set of observations, we are able to predict the probability of occurrence (or the return time) of events larger than anything observed, with an accuracy that depends on the quality of the fit. This clarifies why the existence of universality fosters predictive power.

## 1.3
## Some Difficulties and Challenges in Studying Extremes

### 1.3.1
### Finiteness of Data

It is important to note that, despite the powerful mathematics behind the EVLs, not all problems related to extreme events can be tackled following the approaches mentioned above. Difficulties emerge, *e.g.*, when considering finite data samples and attempting to derive via statistical inference the statistical models of the underlying

EVLs.

While the relationship between *very large* and *very rare* events is far from obvious, it seems hard to conceive - or provide any useful definition - of an extreme as a fairly common event, so that, by construction, the statistical inference of extremes always suffers from a relative (or serious) lack of data for fitting purposes [3]. The problem is more easily seen taking the BM point of view and the classical case of a time series $X_j$, $1 \leq j \leq NK \gg 1$. Assuming that each entry of the time series is the realization of i.i.d. stochastic variables, we need to divide the $NK$ entries into $K$ blocks each of size $N$, and perform our statistical inference over the $M_k$, $1 \leq k \leq K$. Since we are targeting extremes, we clearly have to consider the case $N \gg 1$. Moreover, since we need to perform statistical inference using the GEV model, we definitely need $K \gg 1$ in order to have sufficient robustness. In particular, fitting an $n-$parameter probability distribution requires $O(10^n)$ independent data [45]; hence, considering the GEV method, we need $K = O(10^3)$ candidates as true extremes. Assuming that yearly maxima are reasonable candidates for extremes of temperature at a given location, what said implies that we need time series covering $O(10^3)$ years to perform a reasonable GEV-based analysis of extremes. It is immediately apparent that available observational data - which cover at most three centuries for some meteorological stations (and neglecting all problems of homogenization)- are not appropriate, while one immediately sees the potential of using ultra-long numerical simulations with climate models.

If data are abundant, one can think at many possible options for dividing a time series of length $NK$ into $K$ blocks each of length $N$. One can proceed as follows: larger and larger values of $N$ are considered, until for $N > N_{crit}$ the estimates of the GEV parameters (and in particular of the shape parameter) are stable (and the goodness of fit is high). This allows us to say that we have reached - for all practical purposes - the asymptotic limit of the theory. Further increasing the value of $N$ makes our fitting procedure worse, because we reduce the value of $K$, thus increasing the uncertainty in the parameters' estimate [46]. The basic problem with the BM method is that many excellent candidates for extremes are removed from the analysis because only the maximum of each block is retained in the subsequent statistical inference procedure. Interestingly, in many applications, the POT approach is preferred to the BM approach for fitting time series because it provides more efficient use of data and has better properties of convergence when finite datasets are considered [3]. A comprehensive treatment of optimal threshold selection procedures for the POT method is presented in [43].

If data are relatively scarce, one is bound to relax the criteria for defining extremes (thus considering soft extremes (*e.g.* taking $N$ not too large, or choosing a moderate threshold $T$). As discussed in [46, 47], softening the criteria for choosing extremes leads to biases in the estimates of the EVL distribution parameters, the basic reason being that we quickly pollute the statistics with many bogus entries. Therefore, in some cases one needs to take a more heuristic point of view and construct empirical measures of extremes, defined, *e.g.*, as events above a given percentile - say $95^{th}$ - of the underlying probability distribution. This is, in fact, the standard point of view taken in most climate -related investigations [23]. Unfortunately, as soon as

these - pragmatic and useful - points of view are taken, we lose universality and, consequently, predictive power. This demonstrates that it is important to understand not only what the limiting EVLs of given stochastic processes are, but also evaluate how fast the convergence of finite-data estimates is [46]. The reader is encouraged to look into [48] for an in-depth analysis of extremes in simple atmospheric models, with a detailed discussion of how fast the statistics of block maxima converges to the GEV distribution family depending on the choice of observable.

## 1.3.2
### Correlation and Clustering

Apart from the elements of subjectivity and data requirements in defining how large or rare an event has to be in order to be called a *true extreme*, many additional aspects need to be dealt with when looking at many situations of practical interest. In particular, the property of *recurrence* of extremes, defining whether they come as statistically independent events, or whether they form clusters, *i.e.* groups of extreme events concentrated in time, is extremely important. The two scenarios relate to two different underlying dynamical processes, where the occurrence of extremes is related or not to the presence of persistent patterns of dynamics, and to memory effects, and, in terms of risk, imply entirely different strategies of adaptation, mitigation, and response.

The classical EVT is basically unaltered if the stochastic variables $X_j$'s are, instead of independent, weakly correlated, meaning that the correlation between the variables $X_j$ and $X_k$ decreases *rapidly enough* with $|j-k|$. In the presence of short range (*i.e.* small $|j-k|$) strong correlations between $X_j$ and $X_k$, the GEV and GPD based approaches are not equivalent, the basic reason being that the POT method is bound to select clusters of extremes, which are instead automatically neglected in the BM procedure [3]. As a result, one can prove that when estimating the shape parameter with the POT and BM method using the same time series, one expects to obtain different estimates, with the POT method giving biased results. At practical level, this may result is errors in the estimate of the return times of extreme events. The *extremal index* (EI), the inverse of the average cluster size [49], is the most important indicator of how important is clustering is a given time series, and various statistical methods have been devised to optimally estimate its value [50]. In order to use the POT approach, we need to post process the data by performing *declustering*. Commonly used declustering procedures are based on the idea of discarding all the elements of a cluster except the largest one. After this treatment of data, the POT method typically delivers the same estimates of the shape parameter as the BM method [3].

Taking a concrete example where these issues play a key role, one may want to accommodate situations where the occurrence of an extreme is not exclusively related to the occurrence of a large event, but to the persistence of the considered observable above a certain threshold for an extended period of time, so that clustering of individual events is crucial. This is exactly how heat stress indices are defined by the World Health Organization (WHO) in relation to the persistence of high temperature, be-cause the human body is well suited to resisting to short spells of high temperatures, but has instead great problems in dealing with extended period of physical stress due to heat. See also the definition of heat wave in [51]. Similarly, food security is strongly affected by prolonged heat and drought stress in some key regions in the world and contingency plans - based on risk reduction and insurance-based methods - are continuously updated to take into account the time scale of the persistence of extreme conditions [9, 52].

Some explicit results have been presented in the physical literature regarding extremes of time series resulting from stochastic variables $X_j$ featuring long-term correlations. In particular, one obtains that the recurrence times of occurrence of above-threshold events are not Poisson distributed (which, roughly speaking, implies that occurrence of one extreme does not influence the occurrence of another extreme), but rather follows a stretched exponential statistics, with ensuing implications on the possibility of predicting extremes [53, 54]. A detailed discussion of the performance and limitations of the POT and BM methods in the context of time series featuring substantial long-term correlations is given in [7].

### 1.3.3
### Time Modulations and Noise

Often, the parameters of a system, or of a model, or the boundary conditions of an actual system, change with time: what is the best way to analyze extremes in a context like this? The usual setting of EVT is based upon assuming stationarity in the stochastic variables. When performing statistical inference from data, is it more efficient to remove trends in the time series of the considered observables and then study the extremes of the obtained detrended time series? Or should we analyze the original time series, and use time as a covariate? How do we remove periodic components in the time-series of a process (*e.g.* energy consumption) influenced by *e.g.* the daily and seasonal cycle?

Some of these aspects are dealt with in [3, 55]: it is proposed to treat time as covariate and construct in this way a time-dependent description of extremes. See also [56], where this method is compared with what is obtained by dividing the time series to be studied in smaller blocks, and then performing standard EVT parameter inference in each block separately assuming stationarity, as proposed in [57]. Recently, [58] have proposed new statistical methods for constructing efficient time-dependent EVT models from non-stationary time series, while [59] have underlined the limitations of this approach when forcing terms have different time scales.

This issue is of extremely relevance and urgency, *e.g.* with regard to the very active field of investigation of the meteo-climatic extremes in the context of the changes in the properties of the climate system due to anthropogenic influence. In most cases, for the reasons outlined above, the investigation of changes in extremes is performed by looking at heuristic indicators, such as changes in the probability of occurrence of empirically defined above-thresholds events [60] or in value of very low and high quantiles [61]. Though, it is becoming more common in the geophysical literature to make explicit - *e.g.* [62] - or implicit -*e.g.* [63] - use of EVT. See also [5, 4, 64] for comprehensive reviews.

Another aspect to be mentioned is the role of noise or finite precision in observations. When taking into account real measuring devices, we also face the problem of errors - uncertainties and biases - in recording data. Therefore, observational noise needs to be addressed in order to connect mathematical results to inferences from data [65]. On a related note, [66, 67, 68] concluded that there is a substantial impact of finite precision (typically, 1 $mm$) on the rain gauge readings on the fitting procedures for reconstructing the statistical properties of rainfall data.

## 1.4
## Extremes, Observables, and Dynamics

The open issues and practical problems mentioned above clarify that it is necessary to develop a comprehensive mathematical view of extremes in order to achieve flexibility and predictive power in many real-life problems.

Traditionally, the theory of extremes events has been developed in the context of probability theory and stochastic processes, as suggested by the definitions provided above. Nonetheless, many of the examples we have hinted at suggest that one needs to move in the direction of extending and adapting the classical results of extreme value theory for investigating the extremes of observables of suitably defined dynamical systems. The reader is encouraged to look into [69] for a comprehensive presentation of the theory of dynamical systems, and into [70, 71] for a point of view specifically tailored for linking dynamical systems and (non-equilibrium) statistical mechanics.

Roughly speaking, the conceptual bridge relies on considering a given (continuous or discrete time) dynamical system as the generator of infinitely many stochastic variables, each defined by the time series of a given *observable*, *i.e.* a sufficiently regular function mapping the trajectory of the point representing the evolution of the system in the phase space into the real numbers, and then studying the extremes of such observables. Such point of view, first proposed in [72], has the crucial advantage of giving the possibility of relating the properties of the extremes to the underlying dynamics of the system generating them, thus providing a natural link between probability theory and dynamical systems theory and connecting, in practical terms, forward numerical simulations of - possibly complex - systems with the statistical inference of the properties of extremes from time series.

Moreover, it provides the perfect setting for studying the properties of extremes generated by numerical models, which are finite-precision (and thus noisy) implementations of systems ranging from simple low dimensional maps up to discretized (in time and space) representations of the evolution of continuum systems, such as fluids. It is clear that by considering dynamical systems with different properties of decay of correlations, one mirrors precisely the conditions of weak *vs.* strong correlations of stochastic variables described above. This substantiates the idea of clusters of extreme events, and can define cases where the one-to-one equivalence between BM and POT approaches is broken, so additional mathematical subtlety is required [73, 74, 49].

A key ingredient of a theory of extremes incorporating dynamics and recurrences is the choice of the observable. This provides a new dimension of the problem, and requires the scientist to define what a good, meaningful, useful, well-behaved observable is. Moreover, given a numerical model, one can practically explore many of the aspects related to temporal or spatial coarse graining just by constructing in rather simple ways the observable of interest. This issue naturally provides a more statistical mechanical, physical setting to the problem of extremes, paving the way for interesting applications where extremes can be used as indicators of the structural properties of a system, defining new, powerful methods to study its tipping points [75].

We shall see that, by looking at extremes, one can learn about qualitative properties of the dynamical system generating them, *e.g.* learning whether it features regular or chaotic motions [76], and, under suitable circumstances, understand basic information on the geometry of the attractor and on the Lyapunov exponents determining the growth or decay or errors in given directions. Therefore, extremes act as a probe of a dynamical system, and, when suitable observables are considered, they define a natural microscope to look at the details of the invariant measure, to the point of providing alternative ways to compute Hausdorff dimension of the attractor.

An especially important role is played by observables whose extremes correspond to close returns of the orbit near a reference point [73, 49, 77, 78]. Interestingly, perturbing systems with noise allows the establishment of EVLs for such observables also in the case of deterministic quasi-periodic motion and removes clusters of extreme events when strong correlations are otherwise present [79].

Recurrence-based techniques have also been shown to be directly usable for studying the properties of extremes in climatic time series [80]. Nonetheless, in many practical applications, one is interested in studying a different sort of observables, the so-called *physical observables* [81, 44], which *a priori* have nothing to do with the recurrence of an orbit near a given reference point, but rather describe macroscopic or anyway physically relevant properties of the system. As a simple example, one may consider the extremes of the energy [82] or of temperature [48] in a model of geophysical fluid . The extremes of physical observables permits the study of rather sophisticated aspects of the geometry of the attractor of the underlying system, providing a formidable tool for analyzing at the properties of the unstable and stable components of the tangent space.

One can formulate the problem of studying, at least heuristically, extremes for non stationary systems by taking the point of view of some recent results of non-equilibrium statistical mechanics and dynamical systems theory. This can involve the construction of a response theory for Axiom A dynamical systems to predict the change in the expectation value of suitably defined observables resulting from weak perturbations which are also time dependent, *e.g.* such as small changes in a parameter [83, 84].

In order to apply these results to the problem of assessing the time-dependent properties of extremes - see Sec. 1.3.3 - one needs to construct observables which can represent the extreme events of interest, and then apply response theory to compute their change as a result of the perturbation. Finally, the computed response can be reformulated in terms of time dependent correction to the value of the EVL param-eters [44]. An interesting aspect is that, since extremes are in this case described by the universal parametric family of EVLs, one could draw the rather counter-intuitive conclusions that, in some sense, the response of extremes to perturbations could be a better-posed problem than studying the response of the statistics of the *bulk* of the events [82, 56]. In practical terms, this gives a framework for rigorous questions mentioned before in this introduction, such as determining how extremes change when time-dependent perturbations are added to the system [56], as in the case of changes in climatic extremes in an overall changing climate [85, 86, 87]. Apart from the relevance in terms of applications, the mathematical interest in this regard is considerable.

A different yet related dynamical point of view on extremes of non-stationary systems is based upon considering the mathematical construction of the so-called *pullback* attractor [88, 89, 90], sometimes referred to as *snapshot* attractor [91], which is basically the time-dependent generalisation of the usual attractor appearing in many dissipative chaotic system [70, 71], and enjoys a lot of the useful properties of the latter, except time invariance. The time-dependent measure supported on the pullback attractor at time $t$ is constructed by evolving from an infinitely distant past a Lebesgue measure with compact support. This procedure, in practical numerical applications, corresponds to starting to run in the distant past a large number of simulations with different initial conditions, and observing the resulting ensemble at time $t$ [89]. The time-dependent properties of extremes can then be assessed studying the properties of such an ensemble [92, 93, 59, 94].

This setting suggests the possibility of achieving predictability of extremes in a statistical sense, *i.e.*, developing tools for understanding how their properties change as a result of changing parameters of the system that generates them, which is clearly a major scientific challenge in *e.g.* climate science [4, 23], where *big data* are being advocated [10]. Our ability to predict efficiently the occurrence of individual extreme events is still modest, see some examples in [4]. A crucial aspect is that it is not easy to anticipate (we are not talking about *ex-post* analysis) the dynamical causes leading to an extreme event. As clarified in [95, 96], (finite-time) Lyapunov exponents and related dynamical objects have an important role in assessing the potential predictability of a chaotic system *for typical conditions* [97], in terms of allowing for effectively extremes with a certain lead time. Some authors have proposed ingenious methods for detecting precursors [98, 99], but still a comprehensive understanding of this problem is missing.

## 1.5
## This Book

The scope for the lines of investigation described above is immense, and what we are proposing in this book is a limited perspective resulting from the collective effort of a group of authors coming together and joining forces from a rather diverse spectrum of scientific expertise, ranging from probability theory, to dynamical systems; from statistical mechanics, to geophysical fluid dynamics; from theoretical physics, to time series analysis. Without the hope nor the goal of being comprehensive or conclusive, this book aims at providing a new perspective on extreme events, a transdisciplinary field of research of interest for mathematicians, natural scientists, statisticians, engineers, and social scientists.

One can safely say that the main difference between this book and many other excellent monographs in the literature, ranging from rather abstract mathematical formulations of the theory of extremes [1, 100], to sophisticated presentation of algorithms for defining, detecting, and performing statistical inference of extremes [3, 50], to specific applications [2], is the focus on *dynamics*. In other terms, we do not take extremes as results of a *black box* - a stochastic process whose origin we might not necessarily be interested in *per se* - but rather explore the link between the (typically chaotic) system under investigation and the generating process leading to the extreme and the extreme event itself. The freedom of looking at different extremes is guaranteed by the possibility of choosing different observables, which may be tailored for looking at local, global, or recurrence properties of the system and for focusing on specific regions of its attractor. This perspective leads to considering extremes as a revealing source of information on the microscopic or macroscopic properties of the system, and so naturally suited for improving our understanding of its statistical mechanics. As opposed to many contributions in the scientific literature, our emphasis will not be on presenting ideas for optimizing statistical inference procedures, even if we will present examples and ideas in this direction. Though, we hope to contribute to providing useful guidance for statistical inference procedures by clarifying what they *should* find for given systems.

The main motivation for the approach we propose here comes - mostly but not exclusively - from a more general interest by the authors in exploring the fruitful and emerging *nexus* between mathematics and geophysical fluid dynamics, which has recently received global accolade with the Mathematics for Planet Earth international initiative (http://www.mpe2013.org), and in particular of the programme *Mathematics for the Fluid Earth* (http://www.newton.ac.uk/event/MFE) [101] held at the Newton Institute for Mathematical Science in Cambridge (UK); see also the recent review [102]. Moreover, the authors hope to contribute to stimulating the development of new effective and robust methods for studying extremes in a meteo-climatic context, thus contributing to the global effort for adapting to climate change and climate -related risk. This book tries to provide stimulations, hints, and new results having in mind a readership of (applied) mathematicians, statisticians, theoretical physicists, experts in probability, stochastic processes, statistical mechanics, time series analysis, and in (geophysical) fluid dynamics.

The structure of the book can be described as follows:

- In Chapter 2 we present an overview of general laws and concepts used for describing rare events and introduce some terminology.
- In Chapter 3 the basics of classical extreme value theory are introduced, concentrating on results that are useful for developing a dynamical systems perspective. In this part of the book, there is no reference to dynamics, whereas everything is formulated exclusively in terms of stochastic processes.

- Chapter 4 and Chapter 5 constitute the core of the book. In the former we introduce dynamical systems as generators of random processrandom processes and present a description of a variety of methods and approaches to derive EVLs for the so-called *distance observables*. In the latter, we construct the correspondence between EVLs and hitting/return time statistics in uniformly hyperbolic, non-uniformly hyperbolic, and quasi-periodic systems.

- In Chapter 6 we focus on specific dynamical systems of special interest, and study in detail the role of the decay of correlations for establishing extreme value laws and relate it to the chaotic nature of the dynamics, and investigate the rate of convergence of the statistics of extremes to the asymptotic EVLs. We also introduce and discuss the properties of the extremes of the so-called *physical observables*.

- In Chapter 7 we tackle the important problem of the impact of noise on the statistics of extremes in dynamical systems, treating the case of random perturbations to the dynamics and of observational noise.

- In Chapter 8 we take the point of view of statistical mechanics and, using Axiom A systems as mathematical framework, we discuss extremes in the context of high-dimensional dynamical systems, introducing the so-called *physical observables*, relating the properties of the extremes to the dynamical properties of the underlying system, and proposing a framework for a response theory of extremes.

- In Chapter 9 we move our focus to studying the procedures for statistical inference and present instructive applications of the theory in numerical simulations, investigating the role of finite size effects in the inference. We also present examples of how EVT can be used to derive relevant information on the geometrical and dynamical properties of the underlying system, and investigate how the presence of noise impacts the statistics of extremes.

- Chapter 10 focuses on physically-oriented applications of EVT, showing how it can be used for detecting tipping points in multistable systems and, additionally, for providing a rigorous characterisation of the properties of temperature fields in present and changing climate conditions.

- Chapter 11 contains the concluding remarks of the book and provides indications for future research activities in the field.

- Appendix A includes a few MATLAB© numerical codes used for producing some of the numerical results contained in the book, which are distributed for the benefit of the readers.

### Acknowledgments

The research activities that have led to the preparation of this book have received financial and in-kind support from national and international scientific institutions and funding agencies. Specific acknowledgements are reported below separately for each of the authors.

- VL: DFG (Germany) Cluster of Excellence CliSAP (Hamburg, Germany); Walker Institute for Climate System Research (Reading, UK); London Mathematical

Society (UK); EU FP7-ERC Starting Investigator Grant Thermodynamics of the Climate System - NAMASTE (No. 257106); EU EIT-Climate-KIC PhD Grant *Extreme Events in Southern Pakistan*.

- DF: EU FP7-ERC Advanced Investigator Grant A2C2 (No. 338965).
- ACMF: FCT (Portugal) grant SFRH/BPD/66174/2009, funded by the program POPH/FSE; FCT (Portugal) project PTDC/MAT/120346/2010, funded by national and European structural funds through the programs FEDER and COMPETE; CMUP (UID/MAT/00144/2013), funded by FCT (Portugal) with national (MEC) and European structural funds through the program FEDER, under the partnership agreement PT2020.
- JMF: FCT (Portugal) grant SFRH/BPD/66040/2009, funded by the program POPH/FSE; FCT (Portugal) project PTDC/MAT/120346/2010, funded by national and European structural funds through the programs FEDER and COMPETE; CMUP (UID/MAT/00144/2013), funded by FCT (Portugal) with national (MEC) and European structural funds through the program FEDER, under the partnership agreement PT2020.
- MH: EPSRC (PREDEX - Predictability of Extreme Weather Events): EP/I019146/1 (UK); London Mathematical Society Scheme 4 grants: 41325 and 41367 (UK).
- TK: Walker Institute for Climate System Research (Reading, UK); DFG (Germany) Cluster of Excellence CliSAP (Hamburg, Germany); London Mathematical Society (UK).
- MN: CNRS (France) through a poste d'accueil position at the Centre of Theoretical Physics in Luminy; NSF (USA) grant DMS 1101315; FCT (Portugal) project PTDC/MAT/120346/2010, which is funded by national and European structural funds through the programs FEDER and COMPETE.
- MT: CMUP (UID/MAT/00144/2013), funded by FCT (Portugal) with national (MEC) and European structural funds through the programs FEDER, under the partnership agreement PT2020.
- SV: ANR (France) Project Perturbations (Grant 10-BLAN-0106); the University of Houston (Houston, USA); CNRS (France) PICS Project *Propriétés statistiques des systémes dynamiques detérministes et aléatoires*, with the University of Houston (Grant PICS05968); CNRS (France) PEPS Project *Mathematical Methods of Climate Theory*, CONICYT (Chile) Project Atracciòn de Capital Humano Avanzado del Extranjero MEC 80130047, held at the CIMFAV, University of Valparaiso (Chile); Project MOD TER COM of the French Region PACA; FCT (Portugal) project PTDC/MAT/120346/2010, funded by national and European structural funds through the programs FEDER and COMPETE.

A special thank goes to the I. Newton Institute for Mathematical Sciences (Cambridge, UK) for having supported the scientific programme *Mathematics for the Fluid Earth*, which has provided invaluable opportunities for scientific exchange, cross-pollination, and discovery. This book can be considered a direct product of that stimulating experience.

**18**

This book would have not been possible without the many interactions the authors have had in these years from many friends, colleagues, and interesting people. We wish to acknowledge the direct and indirect support by S. Albeverio, V. Baladi, F. Bassani, R. Blender, T. Bódai, G. Canguilhem, E. Cartman, J.-R. Chazottes, M. Chekroun, P. Collet, R. Deidda, I. Didenkulova, K. Fraedrich, C. Franzke, G. Gallavotti. M. Ghil, G. Gioli, I. Gomes, P. Guiraud, S. Hallerberg, A. Hoffman, P. Imkeller, J. Lebowitz, J. Littell, C. Liverani, P. Manneville, A. Panday, K. Peiponen, B. Pelloni, R. Rajbhandari, D. Ruelle, O. Sacks, B. Saussol, A. Shulgin, A. Speranza, T. Tél, G. Turchetti, S. Tzu, S. Vannitsem, R. Vitolo, and J. Wouters.

This book is dedicated to friendship and to our friends.

# 2

# A Framework for Rare Events in Stochastic Processes and Dynamical Systems

## 2.1
## Introducing Rare Events

The first step of the study of rare events is to establish a mathematical framework where the meaning of rare event becomes precise. The wording already suggests us that we should have a notion of how to measure the frequency with which an event occurs. Hence a probabilistic framework is natural. Therefore, as a setup, we should start with a probability space $\mathcal{X}$ equipped with a $\sigma$-algebra $\mathcal{B}$ of events and a probability measure $\mathbb{P}$ designed to quantify the probability of the occurrence of any event $A \in \mathcal{B}$. We assume the reader has prior knowledge of probability theory but refer to [103, 104] as two examples of textbooks where every possibly missing definition and concept can be easily found and recalled. Any other standard textbook on probability theory should be equally suitable.

We can characterise rare events as those events $A \in \mathcal{B}$ such that $\mathbb{P}(A)$ is *small*. Of course what is meant by small needs to be quantified at some point. As mentioned in Chapter 1, one of the main reasons why we are interested in studying rare events is because they are usually associated with unwanted incidents of high potential environmental, social, financial, infrastructural or technological impacts. Therefore, while sometimes it is impossible to avoid such high risk scenarios, we would like at least to gauge the risk by assessing their probability of occurrence. Hence, another major ingredient in our consideration should be the time dimension of the process leading to the occurrence of the extremes. The classical statistical approach to this matter is to consider that $\mathcal{X}$ is just the space of realisations of a time series or stochastic process of random variables (r.v.) of interest, say $X_0, X_1, \ldots$, which stands as a theoretical model for treating time series of quantities of relevance for the phenomena under investigation

In most cases, along the book $X_0, X_1, \ldots$ is a stationary sequence of r.v. Stationary means that the joint distribution function (d.f.) of any finite collection of r.v. in the process is equal to the joint d.f. of that collection of r.v. displaced by the same time period, *i.e.*, for any $k, t \in \mathbb{N}$ and $i_1 < i_2 < \ldots < i_k \in \mathbb{N}$ the joint d.f. of $X_{i_1}, X_{i_2}, \ldots, X_{i_k}$ is equal to the joint d.f. of $X_{i_1+t}, X_{i_2+t}, \ldots, X_{i_k+t}$. We also denote by $F$ the common (cumulative) d.f. of the r.v. in the process, *i.e.*,

$$F(x) = \mathbb{P}(X_0 \leq x),$$

for every $x \in \mathbb{R}$.

If a r.v. is absolutely continuous, so that the Lebesgue-Stieltjes measure associated to $F$, $\mu_F$, is absolutely continuous with respect to the Lebesgue measure $\mathrm{Leb}$, we can define the Radon-Nikodym derivative function, $f = \frac{d\mu_F}{d\mathrm{Leb}}$, which is usually called a probability density function. If $f$ is Riemann integrable and continuous at $x$, then $f(x) = \frac{dF(x)}{dx}$. See [103, 104] for further details.

In this case, when provided with a finite sample $X_0, \ldots, X_{n-1}$, rare events are almost always tied with abnormal values of the observations in the sample, which means that we are interested in the extremal observations (either very large or very small) typically above (or below) some high (low) thresholds. For historical reasons, we will refer only to exceedances of high thresholds, although the study could be reframed in terms of analysis of extremely small observations, as well. We note that it is sufficient to change the sign of the variable of interest to go from one problem to the other. Hence, in this context, a rare event or extreme event (since it corresponds to an extreme value of the observations) corresponds to the occurrence of an exceedance of the threshold $u$,

$$U(u) := \{X_0 > u\}, \tag{2.1.1}$$

where $u$ is close to the right endpoint of the d.f. of $F$, *i.e.*,

$$u_F = \sup\{x : F(x) < 1\}. \tag{2.1.2}$$

which can be finite or infinite. Since we are interested in large values of the r.v., the behaviour of tail of its d.f. is of crucial importance. Hence, we introduce the following notation: let

$$\bar{F} = 1 - F.$$

the complementary distribution function. The speed at which $\bar{F}(u)$ approaches $0$ as $u \to u_F$ establishes the type of tail we have. In particular, informally, we say that $F$ has *heavy tails* if $u_F = \infty$ and $\bar{F}(u)$ vanishes polynomially fast in $u$, and we have *light tails* if $u_F < \infty$ or $\bar{F}(u)$ vanishes exponentially fast in $u$. We will see later in Chapter 3 how Extreme Value Theory (EVT) allows one to construct a very general framework for studying the tails of the distributions.

When considering time series produced by instruments or numerical models, one should keep in mind that our knowledge of the properties of the tail can depend critically on our observation window in time.

## 2.2
## Extremal Order Statistics

The study of the tails of the d.f. is connected to the behaviour of the extremal order statistics of the sample. Let $X_{1,n} \leq X_{2,n} \leq \ldots \leq X_{n,n}$ denote the order statis-tics of the sample $X_0, X_1, \ldots, X_{n-1}$, so that $X_{1,n}$ is the minimum and $X_{n,n}$ the maximum of such sample. We introduce the following notation:

$$M_n := \max\{X_0, \ldots, X_{n-1}\}, \tag{2.2.1}$$

so that $X_{n,n} = M_n$. Again, for historical reasons we will consider the maximum but note that $X_{1,n} = -\max\{-X_0, \ldots, -X_{n-1}\}$.

In classical statistical analysis, one is mostly concerned with the average behaviour of the sample, which is tied with its mean, $1/n \sum_{i=0}^{n-1} X_i$, or its central order statistics, such as the median $X_{\lfloor n/2 \rfloor, n}$, and the respective Gaussian type asymptotic nature. Here, since we are interested in the abnormal behaviour, described by the tail of $F$, we are concerned with the extreme order statistics such as $M_n$. In fact, observe that the knowledge of $M_n$ allows to conclude if an exceedance has or not occurred among the first $n$ observations since, if $\{M_n \leq u\}$ occurs, than there are no exceedances of $u$ up to time $n - 1$.

In the case of independent events, *i.e.*, , when the stochastic process $X_0, X_1, \ldots$ is a sequence of independent and identically distributed (i.i.d.) r.v., the first statement regarding $M_n$ is that $M_n$ converges almost surely (a.s.) to $u_F$. Then, the next natural question is whether we can find a distributional limit for $M_n$, when it is conveniently normalised. Hence, given a stochastic process $X_0, X_1, \ldots$ we define a new sequence of random variables $M_1, M_2, \ldots$ given by (2.2.1) and propose the following:

**Definition 2.2.1.** We say that we have an *Extreme Value Law* (EVL) for $M_n$ if there is a non-degenerate d.f. $H : \mathbb{R} \to [0, 1]$ with $H(0) = 0$ and, for every $\tau > 0$, there exists a sequence of levels $u_n = u_n(\tau)$, $n = 1, 2, \ldots$, such that

$$n\mathbb{P}(X_0 > u_n) \to \tau, \;\; \text{as } n \to \infty, \tag{2.2.2}$$

and for which the following holds:

$$\mathbb{P}(M_n \leq u_n) \to \bar{H}(\tau), \;\; \text{as } n \to \infty. \tag{2.2.3}$$

where the convergence is meant to hold at the continuity points of $H(\tau)$.

The motivation for using a normalising sequence $(u_n)_{n \in \mathbb{N}}$ satisfying (2.2.2) comes from the case when $X_0, X_1, \ldots$ are independent and identically distributed (i.i.d.). In this setting, it is clear that $\mathbb{P}(M_n \leq u) = (F(u))^n$, where $F$ is the d.f. of $X_0$. Hence, condition (2.2.2) implies that

$$\mathbb{P}(M_n \leq u_n) = (1 - \mathbb{P}(X_0 > u_n))^n \sim \left(1 - \frac{\tau}{n}\right)^n \to e^{-\tau}, \tag{2.2.4}$$

as $n \to \infty$. Moreover, the reciprocal is also true (see [1, Theorem 1.5.1] for more details). Note that in this case $H(\tau) = 1 - e^{-\tau}$ is the standard exponential d.f.

The study of EVLs is the starting point for a deeper analysis of the properties of rare events. But the main novelty we present here is the study of rare events for dynamical systems. We mentioned before that it is natural to use EVT to study the properties of time series reflecting the evolution in time of a random variable. But, in fact, one can bring together randomness and time evolution by considering purely deterministic systems featuring a so-called *chaotic* dynamics.

## 2.3
## Extremes and Dynamics

In many situations, the evolution of natural, engineering, and social phenomena can be represented by mathematical models able to incorporate how the properties of the system change in time. The research field focussing the study of how systems evolve with time is usually referred to as *dynamical systems* [70, 69]. In the last century, the discovery of the ubiquity of chaotic systems and the investigation of their erratic yet extremely structured behaviour has shed new light on the interpretation of stochastic processes.

A dynamical system is build up out of a phase space, where each point describes a given state of the system, a model for the passage of time, which can be continuous or discrete, and an evolution law that rules the unfolding of the system. In the continuous case, the evolution laws are usually given in terms of differential equations, while in the discrete case a map describe the transitions from one state to another. One of the most famous examples of a dynamical system is the celebrated Lorenz '63 model [105]. It consists of a system of three differential equations designed to provide a minimal representation of convection via a minimal representation of its dynamics. This system provides the prototype of chaotic behaviour and has been extensively studied in a variety of disciplines including mathematics, physics, geophysics, and time series analysis.

The erratic behaviour of this system makes its outputs hardly distinguishable from purely randomly generated numbers, at least on the long run. This led to the coining of the expression *butterfly effect*, or, in more precise terms, the sensitivity to initial conditions and the presence of impassable and intrinsic limits to deterministic prediction for chaotic system. See [97] for a detailed exposition of the whole conceptual and technological chain stretching from the Lorenz '63 model to the modern technological infrastructure and methods for weather forecast. Hence, chaotic dynamical systems provide another area to study rare events, the difference being that random behaviour is generated by a deterministic model, whose characteristics are shown later in the book to have a crucial impact on the properties of the extremes.

We briefly introduce some terminology for dynamical systems, taking the example of discrete time evolution. In this setting, $\mathcal{X}$ denotes the phase space which we endow some topological, or differentiable, or measure-theoretical structure, depending on the problem under consideration. We are particularly interested in the latter, and denote by $\mathcal{B}$ the associated $\sigma$-algebra of subsets of $\mathcal{X}$, which gives the measurable structure. The system itself is represented by a measurable map (with a structure compatible with that of $\mathcal{X}$) that we denote by $T$ (or sometimes $f$) such that $T : \mathcal{X} \to \mathcal{X}$ describes the time evolution of the system, *i.e.* determines the rule defining the transition from an initial state $x \in \mathcal{X}$ to the state $T(x)$, after one unit of time. Moreover, we require the existence of a probability measure $\mathbb{P}$ defined on $\mathcal{B}$ that is coherent with $T$, in the sense that, $\mathbb{P}(T^{-1}(A)) = \mathbb{P}(A), \forall A \in \mathcal{B}$, which we express by saying that $\mathbb{P}$ is $T$-invariant.

For an initial state $x \in \mathcal{X}$ we define its orbit to be the sequence of states through which the system will go if it is started at that particular state $x$, *i.e.*, the sequence

$x, T(x), T^2(x), \ldots$, where

$$T^n(x) = \underbrace{T \circ T \circ \ldots \circ T}_{\text{n times}}(x)$$

is the $n^{th}$ iterate of $T$. The main goal of dynamical systems theory is to study the long-term behaviour of the orbits of the system.

A rare event is given by some subset $A \subset \mathcal{X}$ of the phase space belonging to $\mathcal{B}$ and such that $\mathbb{P}(A)$ is small. The occurrence of a rare event is obtained whenever the orbit of some point hits the target set $A$. Poincaré's Recurrence Theorem guarantees that almost every point of $A$ will return to $A$ infinitely often (i.o.) (see [106]). Moreover, if a system is ergodic with respect to $\mathbb{P}$ and $\mathbb{P}(A) > 0$, then almost every (a.e.) orbit will hit $A$ at some point in time. A system is said to be ergodic with respect to $\mathbb{P}$ if for every $B \in \mathcal{B}$ such that $T^{-1}(B) = B$ then either $\mathbb{P}(B) = 0$ or $\mathbb{P}(B) = 1$. In other words, a system is ergodic if it is not decomposable into two parts (detected by the measure $\mathbb{P}$) that do not interact with each other, see [106].

We define the *first hitting time* $r_A : \mathcal{X} \to \mathbb{N} \cup \{\infty\}$ by

$$r_A(x) := \inf\{n \in \mathbb{N} : T^n(x) \in A\}. \tag{2.3.1}$$

If $x$ never enters $A$ under iteration by $f$ then $r_A(x) = \infty$. When the point $x \in A$ then we say that $r_A(x)$ is the first return of $x$ to $A$. Kac's theorem gives that the mean return time to $A$ is equal to the reciprocal of $\mathbb{P}(A)$. Hence, at least in average, $r_A(x)$ should go to $\infty$ as $\mathbb{P}(A)$ goes to zero. A natural question then is if we can find a limit distribution for $r_A$, when conveniently normalised. By Kac's theorem the natural candidate to normalise $r_A$ would be $\mathbb{P}(A)$. The study of such distributional limits are known as *Hitting Times Statistics* (HTS) and, when we start in $A$, as *Return Times Statistics*, see Chap. 5 for precise definitions.

But the same formalism used in the classical study of extreme values can be brought into consideration simply by using the system to generate a time series by evaluating a certain observable along its orbits. Let $\varphi : \mathcal{X} \to \mathbb{R} \cup \{\pm\infty\}$ be a measurable observable function and define the stochastic process $X_0, X_1, \ldots$ given by $X_0 = \varphi$ and

$$X_n = \varphi \circ T^n, \quad \text{for each } n \in \mathbb{N}.$$

A realisation of such process corresponds to picking an initial state $x \in \mathcal{X}$ at random, according to $\mathbb{P}$, and then evaluating the function $\varphi$ along all points of the orbit $x$, $T(x), T^2(x), \ldots$.

Having an exceedance of a certain threshold $u$ means that the orbit hits a region of the phase space described by $U(u) = \{X_0 > u\}$. Moreover, note that

$$T^{-1}(\{M_n \leq u\}) = \{r_{U(u)} > n\},$$

which means that there should be a connection between the existence of an EVL and the existence of HTS. This link will be fully described in Section 5.3.

The study can be deepened by considering for example point processes that keep track of the number of exceedances or hits to certain regions of the phase space and allow to paint a global picture of the extremal type behaviour of the systems. A main goal of this book is then performing such statistical study of chaotic systems and understanding how the structure (geometric, topological) and the properties of the system emerge and influence its extremal behaviour. Moreover, we will see that the study of the extremes for dynamical systems can be used to study and uncover geometric properties and other features of the system itself.

# 3
# Classical Extreme Value Theory

In this chapter we review some of the more emblematic results of classical EVT. It is not intended as a comprehensive exposition but the main motivation is to give a brief introduction to the techniques and notation, with the purpose of making the book as self contained as possible and to motivate the further developments in the following chapters.

The choice of subjects and issues addressed in this chapter does not indicate judgement of their importance from our side. Instead, we have the intent of providing the reader some background and perspective for a better understanding of the problems and advances described in the rest of the book.

Several excellent books on classical EVT are available and they should be seen both as complement and as reference for the material in this chapter. We mention some of them: [107, 1, 108, 3, 109, 110].

In this chapter we will start by recalling some by now classical results in EVT, first in the context of i.i.d. sequences of r.v. Then we recall some conditions introduced by Leadbetter to obtain EVL in the dependent case for stationary stochastic processes. Afterwards we introduce the concept of clustering and further developments concerning the dependence conditions needed to recover the existence of EVL in the presence of clustering.These mixing conditions are reviewed in detail in order to compare them with the new ones proposed in Chapter 4. One of the best ways to understand clusters is by studying point processes of exceedances or rare events. For that reason, we make a small detour to define such point processes and establish both notation and definitions (such as weak convergence) in order to make clearer the statements of some of the results appearing in the following chapters. Understanding clustering, the Extremal Index and how it affects statistical inference is further discussed in a brief description of the some well-known declustering procedures. Finally, a brief discussion of statistical inference based on EVT is presented.

## 3.1
## The i.i.d. Setting and the Classical Results

### 3.1.1
### Block Maxima and the Generalized Extreme Value Distribution

In the classical theory, the sequences of real numbers $u_n = u_n(\tau)$, $n = 1, 2, \ldots$, are usually taken to be one parameter linear families like

$$u_n = y/a_n + b_n, \tag{3.1.1}$$

where $y \in \mathbb{R}$ and $a_n > 0$, for all $n \in \mathbb{N}$. In fact, in the classical theory, one considers the convergence of probabilities of the form

$$\mathbb{P}(a_n(M_n - b_n) \leq y). \tag{3.1.2}$$

The main classical result of EVT is the so called *Extremal Types Theorem* usually credited to Gnedenko [38] although some previous version was already stated in the work of Fisher and Tippett [111].

**Theorem 3.1.1** ([111, 38]). *If $X_0, X_1, \ldots$ is a sequence of i.i.d. r.v. and there exist linear normalising sequences $(a_n)_{n \in \mathbb{N}}$ and $(b_n)_{n \in \mathbb{N}}$, with $a_n > 0$ for all $n$, such that*

$$\mathbb{P}(a_n(M_n - b_n) \leq y) \to G(y), \tag{3.1.3}$$

*where the convergence occurs at continuity points of $G$, and $G$ is non-degenerate [1], then $G(y) = \mathrm{e}^{-\tau(y)}$, where $\tau(y)$ is of one of the following three types (for some $\beta, \gamma > 0$):*

*1) $\tau_1(y) = \mathrm{e}^{-y}$ for $y \in \mathbb{R}$;*
*2) $\tau_2(y) = y^{-\beta}$ for $y > 0$;*
*3) $\tau_3(y) = (-y)^{\gamma}$ for $y \leq 0$.*

These types are usually called Gumbel (type 1), Fréchet (type 2) and Weibull (type 3) (cumulative) distributions.

**Definition 3.1.2.** The three types may be combined in a unified model called the Generalised Extreme Value (GEV) distribution:

$$G(y) = \mathrm{GEV}_\xi(y) = \begin{cases} \mathrm{e}^{-(1+\xi y)^{-1/\xi}}, 1 + \xi y > 0, & \text{if } \xi \neq 0 \\ \mathrm{e}^{-\mathrm{e}^{-y}}, & y \in \mathbb{R}, \quad \text{if } \xi = 0 \end{cases}. \tag{3.1.4}$$

When $\xi = 0$, the distribution corresponds to the Gumbel type; when the index is positive, it corresponds to a Fréchet; when the index is negative, it corresponds to a Weibull.

---

[1] A d.f. is non-degenerate if there is no $y_0 \in \mathbb{R}$ such that $G(y_0) = 1$ and $G(y) = 0$, for all $y < y_0$

Informally, one can sat that *exponential tail* corresponds to $\xi = 0$ or Gumbel type, *heavy tail* corresponds to $\xi > 0$ or Fréchet type, and *upper bounded* corresponds to $\xi < 0$ or Weilbull type.

We emphasise that as observed in [38], for i.i.d. sequences of r.v., the limiting distribution type of the partial maxima is completely determined by the tail of the d.f. $F$. We recall the definition of the right endpoint $u_F$ given in Eq. 2.1.2. As can also be found in [1, Theorem 1.6.2], in order to define the respective domain of attraction for EVT, we have the following sufficient and necessary conditions on the tail of the d.f. $F$:

**Type 1** $G(y) = \mathrm{e}^{-\tau_1(y)}$ (Gumbel) if and only if (iff) iff $u_F = \infty$ and there exists some strictly positive function $h : \mathbb{R} \to \mathbb{R}$ such that for all $y \in \mathbb{R}$

$$\lim_{s \to u_F} \frac{\bar{F}(s + yh(s))}{\bar{F}(s)} = \mathrm{e}^{-y}; \tag{3.1.5}$$

**Type 2** $G(y) = \mathrm{e}^{-\tau_2(y)}$ (Fréchet) iff $u_F = \infty$ and there exists $\beta > 0$ such that for all $y > 0$

$$\lim_{s \to u_f} \frac{\bar{F}(sy)}{\bar{F}(s)} = y^{-\beta}; \tag{3.1.6}$$

**Type 3** $G(y) = \mathrm{e}^{-\tau_3(y)}$ (Weibull) iff $u_F < \infty$ and there exists $\gamma > 0$ such that for all $y > 0$

$$\lim_{s \to 0} \frac{\bar{F}(u_F - sy)}{\bar{F}(u_F - s)} = y^{\gamma}. \tag{3.1.7}$$

*Remark* 3.1.3. It may be shown that $\int_0^\infty 1 - F(u)\,du < \infty$ when a Type 1 limit holds, and one appropriate choice of $h$ is given by

$$h(s) = \frac{\int_s^{u_F} 1 - F(u)\,du}{1 - F(s)},$$

for $s < u_F$.

*Remark* 3.1.4.

As we can see in [1, Corollary 1.6.3], the normalising constants $a_n$ and $b_n$ may be taken as follows:

**Type 1:**

$$a_n = [h(\gamma_n)]^{-1}, b_n = \gamma_n;$$

**Type 2:**

$$a_n = \gamma_n^{-1}, b_n = 0;$$

**Type 3:**

$$a_n = (u_F - \gamma_n)^{-1}, b_n = u_F,$$

with $\gamma_n = F^{-1}(1 - 1/n) = \inf\{y : F(y) \geq 1 - 1/n\}$.

*Remark* 3.1.5. An important aspect related to the study of the statistics of extremes is to provide inference methodologies appropriate to the tail of the d.f. $F$, in order to fit the data to the correct model and to estimate parameters connected with rare events, like high quantiles or the mean waiting time between the occurrence of extreme events. As mentioned in Chapter 1, the detailed exploration of the methods for performing accurate and efficient statistical inference of extremes is outside the scope of this book. While we discuss below and in subsequent chapters some aspects of this line of investigation, the reader is encouraged to look into, *e.g.*, [2, 3, 109, 50] for further reading and references on these topics.

In order to gather a more practical view on the problem, we can reformulate the main result of Theorem 3.1.1 as the fact that when the sample size $n \to \infty$,

$$\mathbb{P}(M_n \leq y) = F^n(y) \approx \text{GEV}_\xi \left( \frac{y - \mu}{\sigma} \right),$$

for some $\mu \in \mathbb{R}$ and $\sigma > 0$., where $\xi \in \mathbb{R}$ is the shape parameter, $\mu \in \mathbb{R}$ is the location parameter and $\sigma > 0$ the scale parameter. For $n$ large, the normalising constant $a_n \to \sigma^{-1}$ and $b_n \to \mu$. This motivates the parametric method suggested by Gumbel for studying the tail of a times series resulting from a stochastic process. It is the *block maxima* (BM) approach, whereby one basically divides the series of recorded data into $k \gg 1$ bins or blocks of length $n \gg 1$. Then the maximum of each block is retrieved and the empirical distribution of the $k$ maxima is then fitted by the *best* matching GEV type of d.f. This implies, in particular, estimating the shape $\xi$, location $\mu$, and scale $\sigma$ parameters. The obtained best estimate of $\mu$ and of $\sigma$ give the $n^{th}$ element of normalising sequence $b_n$ and the inverse of the $n^{th}$ element of normalising sequence $a_n$ given in Eq. 3.1.2, respectively. Note that, when $n$ is large enough, the estimate of $\xi$, as opposed to the case of $\sigma$ and of $\mu$, does not depend on $n$ [46]; see also Sect. 9.1.1.

The fitting procedure is often done using maximum likelihood estimation (MLE)[3]. Other approaches to the problem of fitting have been proposed in the literature, most notably the weighted moments method [112, 113] and the L-moments method [114]. [2] In Chapter 9 we describe in some detail some strengths and weaknesses of MLE and L-moments methods.

Of course, time series almost invariably feature non-trivial correlations, so that the whole mathematical construction provided above might seem irrelevant for time series analysis, as the stochastic variables $X_1, X_2, \ldots$ are assumed to be independent. In fact, we see shortly below in Section 3.2 how one can (partially) circumvent such an issue. The problem of serial dependence of data is also dealt with later in Section 4.1 when addressing observables of dynamical systems. Instead, in Chapter 9 we present various examples clarifying how to practically implement time series analysis.

---

[2] The literature in this field is very large and rapidly evolving. We have made here reference only to a very limited set of relevant contributions.

## 3.1.2
## Examples

In this subsection, we present one example for each of the three different types of limit laws for $M_n$ and compute, in each case, normalising sequences $a_n > 0$ and $b_n$ for which (3.1.3) holds.

**Example 3.1.1. - Exponential distribution**

For the exponential distribution with parameter 1, we have $F(x) = 1 - e^{-x}$, $x > 0$. For $\tau > 0$, we may consider $u_n$ such that $n(1 - F(u_n)) = \tau$, i.e.

$$u_n = \log n - \log \tau.$$

Thus,

$$\mathbb{P}(M_n - \log n \leq -\log \tau) \to e^{-\tau}.$$

Putting $\tau = e^{-y}$, we obtain

$$\mathbb{P}(M_n - \log n \leq y) \to e^{-e^{-y}}.$$

So, in this case, $M_n$ has an EVL of Type 1, with

$$a_n = 1 \quad \text{and} \quad b_n = \log n.$$

**Example 3.1.2. - Pareto distribution**

For the Pareto distribution, we have that $F(x) = 1 - kx^{-\alpha}$, for $\alpha > 0$, $k > 0$, $x > k^{1/\alpha}$. For $\tau > 0$, we may choose $u_n$ such that $n(1 - F(u_n)) = \tau$, i.e.

$$u_n = \left(\frac{kn}{\tau}\right)^{1/\alpha}.$$

Thus,

$$\mathbb{P}((kn)^{-1/\alpha} M_n \leq \tau^{-1/\alpha}) \to e^{-\tau}.$$

By writing $\tau = y^{-\alpha}$ for $y \geq 0$, we obtain

$$\mathbb{P}((kn)^{-1/\alpha} M_n \leq y) \to e^{-y^{-\alpha}}.$$

So, $M_n$ has an EVL of Type 2 with

$$a_n = (kn)^{-1/\alpha} \quad \text{and} \quad b_n = 0.$$

**Example 3.1.3. - Uniform distribution**

For the uniform distribution on $(0, 1)$, we have that $F(x) = x$, $0 \leq x \leq 1$. Let $\tau > 0$ and $u_n$ be such that $n(1 - F(u_n)) = \tau$, i.e.

$$u_n = 1 - \frac{\tau}{n}.$$

Consequently, we have that

$$\mathbb{P}(n(M_n - 1) \leq -\tau) \to e^{-\tau}.$$

By writing $\tau = -y$ for $y < 0$, we obtain

$$\mathbb{P}(n(M_n - 1) \leq y) \to e^y.$$

So, $M_n$ has an EVL of Type 3 with

$$a_n = n \quad \text{and} \quad b_n = 1.$$

### 3.1.3
### Peaks Over Threshold and the Generalised Pareto Distribution

One of the problems of the Gumbel's BM method is the fact that by looking only at the maxima in each bin, some high values (corresponding to *actual* extreme events) may be discarded.

The BM method is clearly very data hungry, so that one, in many practical situation, might want to be able to take into account any large value in a block not only the maximal one. Another approach to the analysis of extremes, called the Peaks Over Threshold (POT) method, relies on establishing a high threshold $T$ and retrieving all the data that exceed $T$, which is used as a tuning variable. This approach is based on the following theoretical result:

**Theorem 3.1.6** ([41, 40]). *Assume that $X_0, X_1, \ldots$ is a sequence of i.i.d. r.v.. Then, there exist linear normalising sequences $(a_n)_{n \in \mathbb{N}}$ and $(b_n)_{n \in \mathbb{N}}$, with $a_n > 0$ for all $n$, such that*

$$\mathbb{P}(a_n(M_n - b_n) \leq y) \to GEV_\xi(y)$$

*if and only if*

$$\lim_{T \to u_F} \sup_{0 \leq y < u_F - T} |F_T(y) - GPD_{\xi,\sigma}(y)| = 0,$$

*where, for $0 \leq y < u_F - T$,*

$$F_T(y) = \mathbb{P}(X_0 - T \leq y | X_0 > T) = \frac{F(y + T) - F(T)}{1 - F(T)},$$

*and, for $\sigma = \sigma(T) > 0$, $GPD_{\xi,\sigma}(y) = GPD_\xi(y/\sigma)$, where*

$$GPD_\xi(y) = \begin{cases} 1 - (1 + \xi y)^{-1/\xi}, & y \geq 0, 1 + \xi y > 0, & \text{if } \xi \neq 0 \\ 1 - e^{-y}, & y \geq 0, & \text{if } \xi = 0 \end{cases} \quad (3.1.8)$$

*is the so-called univariate Generalised Pareto Distribution (GPD). The parameters $\xi$ and $\sigma$ are referred to as shape and scale parameters, respectively.*

It is crucial to note that the shape parameter $\xi$ is the same for the corresponding GEV and GPD distributions. Furthermore, a simple functional relation connects the two distributions:

$$1 + \log(GEV_\xi(y)) = GPD_\xi(y).$$

Given a long time series, the POT approach to extremes boils down to fitting the selected data having value larger than the threshold $T$ data to the distributions $GPD_\xi(y/\sigma)$, and estimating the values of $\xi$ and $\sigma$.

As seen above, when analyzing extremes using the BM method, we need to find a good compromise between choosing very long bins and having many maxima to use for the statistical inference procedure. In the case of the the POT method, we have to find a good compromise between setting a very high threshold $T$, so that only true extremes are selected, and making sure that such selection is such that a sufficient number of data is retained for the purpose of data fitting. It is crucial to make sure that our estimates for $\xi$ and $\sigma$ are robust when $T$ is varied within the asymptotic regime of large values of $T$. In such a regime, the estimates of $\xi$ do not depend on $T$. Instead, as indicated above, $\sigma$ does depend on $T$ also in the asymptotic regime $T \to u_F$, while the quantity $\sigma - \xi T$ (*modified scale parameter*) does not. A successful fit is obtained when we can define a (high) threshold $\bar{T}$ such that $\forall T > \bar{T}$, the estimates of $\xi$ and $\sigma - \xi T$ are compatible within statistical uncertaininties [3, 115].

The estimation of the GPD tail index parameter $\xi$ is of great importance in several applications [3, 109, 50]. Plenty of estimators for the tail index have been proposed. We mention the popular Hill's estimator [116] and its enhanced version, where special weights are designed for the log differences in order to reduce the bias [117]. Many other estimators have been proposed in the literature; see, *e.g.*, [118, 119, 120]. In Chapter 9 we present examples of statistical inference of POT performed using standard GPD-based methods.

## 3.2
## Stationary Sequences and Dependence Conditions

After the success of the classical Extremal Types Theorem of Fisher-Tippet and Gnedenko in the i.i.d. setting, there has been a great deal of interest in studying the existence of EVL for dependent stationary stochastic processes. Building up on the work of Loynes and Watson, Leadbetter proposed in [121] two conditions on the dependence structure of the stochastic processes, which he called $D(u_n)$ and $D'(u_n)$, that guaranteed the existence of the same EVLs of the i.i.d. applied to the partial maxima of sequences of r.v. satisfying those conditions.

Let $X_0, X_1, X_2, \ldots$ be a stationary sequence of r.v..

Condition $D(u_n)$ is a sort of uniform mixing condition adapted to this setting of extreme values where the main events of interest are exceedances of the threshold $u_n$. Let $F_{i_1,\ldots,i_n}$ denote the joint d.f. of $X_{i_1}, \ldots, X_{i_n}$, and set $F_{i_1,\ldots,i_n}(u) = F_{i_1,\ldots,i_n}(u, \ldots, u)$.

**Condition** ($D(u_n)$)**.** We say that $D(u_n)$ holds for the sequence $X_0, X_1, \ldots$ if for any integers $i_1 < \ldots < i_p$ and $j_1 < \ldots < j_k$ for which $j_1 - i_p > t$, and any large $n \in \mathbb{N}$,

$$\left| F_{i_1, \ldots, i_p, j_1, \ldots, j_k}(u_n) - F_{i_1, \ldots, i_p}(u_n) F_{j_1, \ldots, j_k}(u_n) \right| \leq \alpha(n, t),$$

uniformly for every $p, k \in \mathbb{N}$, where $\alpha(n, t_n) \xrightarrow[n \to \infty]{} 0$, for some sequence $t_n = o(n)$.

Condition $D'(u_n)$ precludes the existence of clusters of exceedances of $u_n$. Let $(k_n)_{n \in \mathbb{N}}$ be a sequence of integers such that

$$k_n \to \infty, \quad \lim_{n \to \infty} k_n \alpha(n, t_n) = 0, \quad \text{and} \quad k_n t_n = o(n). \tag{3.2.1}$$

**Condition** $(D'(u_n))$**.** We say that $D'(u_n)$ holds for the sequence $X_0, X_1, X_2, \ldots$ if there exists a sequence $\{k_n\}_{n \in \mathbb{N}}$ satisfying (3.2.1) and such that

$$\lim_{n \to \infty} n \sum_{j=1}^{\lfloor n/k_n \rfloor} \mathbb{P}(X_0 > u_n, X_j > u_n) = 0.$$

**Theorem 3.2.1** ([122, Theorem 1.2])**.** *Let $X_0, X_1, \ldots$ be a stationary stochastic process and $(u_n)_{n \in \mathbb{N}}$ a sequence satisfying (2.2.2), for some $\tau > 0$. If $D(u_n)$ and $D'(u_n)$ hold, then $\bar{H}(\tau) = \mathrm{e}^{-\tau}$.*

*Remark* 3.2.2. Note that, as a consequence of the previous theorem, under the assumptions $D(u_n)$ and $D'(u_n)$ the statement of Theorem 3.1.6 remains valid. Therefore, in a variety of applications where one needs to study extremes of time series featuring a sufficiently fast decay of correlations, GEV- and GPD-based statistical inference methods are used almost interchangeably, under the overall consensus that following the POT approach is *more* efficient when the time series are not exceptionally long [3, 5, 43, 85, 4, 44]. Differences between the two methods emerge when extremes come in clusters; this is discussed below in Sections 3.2.2-3.4.

## 3.2.1
## The Blocking Argument

The core of the argument of the previous result is a blocking argument similar to the one introduced by Markov and used to prove a Central Limit Theorem for stationary stochastic processes. The idea is to break the $n$ observations into $k_n$ blocks (where $k_n$ is as in (3.2.1)) of size $\lfloor n/k_n \rfloor$ and then separate them by adding time gaps of size $t_n$ between each block. These time gaps introduce some sort of independence between the blocks on account of condition $D(u_n)$, which motivates coupling with i.i.d. r.v. with the same d.f. as that of $X_0$. Finally, condition $D'(u_n)$ is used to decorrelate the information contained in each block.

To be more precise, instead of $n$ r.v. of the process $X_0, X_1, \ldots$ we introduce the time gaps by considering $k_n(\lfloor n/k_n \rfloor + t_n)$, instead. We also introduce the following notation. For a subset $A \subset \mathbb{N}_0$ we let $M_A := \max_{k \in A} X_k$. The blocks are defined, for each $m = 1, \ldots, k_n$, by $B_m := \{(m-1)(\lfloor n/k_n \rfloor + t_n), \ldots, m\lfloor n/k_n \rfloor - 1 + (m-1)t_n)\}$. There are $k_n$ disjoint blocks of length $\lfloor n/k_n \rfloor$, which are separated by gaps defined, for each $m = 1, \ldots, k_n$, by $G_m := \{m\lfloor n/k_n \rfloor - 1 + (m-1)t_n), m(\lfloor n/k_n \rfloor + t_n) - 1\}$.

Then we proceed with the following sequence of approximations.

- $\left| \mathbb{P}(M_n \le u_n) - \mathbb{P}(M_{k_n(\lfloor n/k_n \rfloor + t_n)} \le u_n) \right| \le k_n t_n \mathbb{P}(X_0 > u_n)$. In this first step we notice that adding the r.v. that will be used to create the gaps does not change the asymptotic distributional limit. By choice of the sequences $(u_n)_{n \in \mathbb{N}}$, $(t_n)_{n \in \mathbb{N}}$ and $(k_n)_{n \in \mathbb{N}}$, satisfying (2.2.2) and (3.2.1), it follows easily that $\lim_{n \to \infty} k_n t_n \mathbb{P}(X_0 > u_n) = 0$

- $\left| \mathbb{P}(M_{k_n(\lfloor n/k_n \rfloor + t_n)} \le u_n) - \mathbb{P}(M_{\bigcup_{m=1}^{k_n} B_m} \le u_n) \right| \le k_n t_n \mathbb{P}(X_0 > u_n)$. In this step we see that disregarding the r.v. in the gaps leads to an error term equal to $k_n t_n \mathbb{P}(X_0 > u_n)$, which corresponds to the possibility of having at least one exceedance among the r.v. of the $k_n$ gaps of size $t_n$. The corresponding error term converges to zero due to the same argument as in the previous point.

- $\left| \mathbb{P}(M_{\bigcup_{m=1}^{k_n} B_m} \le u_n) - (\mathbb{P}(M_{B_1} \le u_n))^{k_n} \right| \le k_n \alpha(n, t_n)$. In this step, we use condition $D(u_n)$ to obtain an approximate independence between the blocks of r.v. separated by a time gap from each other. In fact, $k_n \alpha(n, t_n)$ goes to 0 by $D(u_n)$ and choice of $k_n$.

- $\left| \mathbb{P}^{k_n}(M_{\lfloor \frac{n}{k_n} \rfloor} \le u_n) - (1 - \lfloor \frac{n}{k_n} \rfloor \mathbb{P}(X_0 > u_n))^{k_n} \right|$. This final approximation is shown to be controlled by the anti-clustering condition $D'(u_n)$, which precludes the appearance of more than one exceedance in each block. The main error term in this approximation is bounded by a constant times $n \sum_{j=1}^{\lfloor n/k_n \rfloor} \mathbb{P}(X_0 > u_n, X_j > u_n)$, which vanishes, on account of $D'(u_n)$.

In order to apply the theory to dynamical systems, condition $D(u_n)$ has to be revised because it can hardly be verified in that context, except in some trivial situations. This means that the blocking argument needs to be refined in order to accommodate a weakening of the original $D(u_n)$ condition. This refinement is mostly concentrated on the third approximation and is accomplished by a more intensive use of condition $D'(u_n)$ to compensate the weakening of $D(u_n)$. This generalisation of the blocking argument is carried out in Section 4.1, where new conditions that, on one hand, allows us to weaken $D(u_n)$ and, on the other hand, allows us to treat both the presence and absence of clusters, are introduced and shown to be useful to prove the existence of EVL.

In Chapter 6, the blocking argument is revisited in a more dynamical setting and is used to devise, in particular, convergence rates for specific examples of dynamical systems, as in Section 4.1.

### 3.2.2
### The Appearance of clusters of Exceedances

In the previous result, condition $D'(u_n)$ plays a double role: together with $D(u_n)$, it guarantees the existence of a limit distribution for $M_n$ and also makes sure that the limit law was such that $\bar{H}(\tau) = e^{-\tau}$, which is exactly what appears in the i.i.d. setting. When $D'(u_n)$ does not hold but $D(u_n)$ does, although we cannot assert the existence of a limit distribution, in case this limit exists, we can still say something about the type of law one should expect for the limit. This is the content of the following theorem.

**Theorem 3.2.3** ([122, Theorem 2.2]). *Let* $X_0, X_1, \ldots$ *be a stationary stochastic process and* $(u_n)_{n \in \mathbb{N}}$ *a sequence satisfying* (2.2.2). *Suppose* $D(u_n)$ *holds for each choice of* $\tau$. *If the limit* (2.2.3) *exists then there exists* $0 \leq \theta \leq 1$ *such that* $\bar{H}(\tau) = e^{-\theta\tau}$ *for all* $\tau > 0$.

In certain circumstances, observed data clearly exhibit the existence of *clusters* of exceedances [123]. By clusters, we mean consecutive occurrences of an exceedance of a given threshold, which thus causes $D'(u_n)$ to fail. The presence of clusters is clearly related to the memory properties of the underlying process. This motivates the study of the effect of clustering on EVLs [122]. In fact, one can observe that clustering of exceedances essentially produces the same type of EVL but with a parameter $0 \leq \theta \leq 1$, the Extremal Index (EI), so that $\bar{H}(\tau) = e^{-\theta\tau}$: here $\theta$ quantifies the intensity of clustering [1].

**Definition 3.2.4.** We say that $X_0, X_1, \ldots$ has an *Extremal Index* (EI) $0 \leq \theta \leq 1$ if we have an EVL for $M_n$ with $\bar{H}(\tau) = e^{-\theta\tau}$ for all $\tau > 0$.

The notion of the EI was latent in the work of Loynes [124] but was established formally by Leadbetter in [122]. The parameter $\theta$ quantifies the strength of the dependence of $X_0, X_1, \ldots$, so that $\theta = 1$ indicates that the process has practically no memory while very low values of $\theta > 0$, conversely, reveals extremely long memory. In particular, when $\theta > 0$, one can interpret the inverse of the EI $\theta^{-1}$ as the mean number of exceedances of a high level in a cluster of large observations, *i.e.*, is the *mean size of the clusters*.

In order to show the existence of EVLs with a certain EI $\theta \leq 1$, new conditions (replacing $D'(u_n)$) were devised. We refer to condition $D''(u_n)$ of [125] and particularly the more general condition $D^{(k)}(u_n)$ of [126], which also includes the case of absence of clustering.

**Condition** ($D^{(k)}(u_n)$)**.** We say that condition $D^{(k)}(u_n)$ holds for the sequence $X_0$, $X_1, X_2, \ldots$ if there exist sequences $\{k_n\}_{n \in \mathbb{N}}$ and $\{t_n\}_{n \in \mathbb{N}}$ satisfying (3.2.1) and such that

$$\lim_{n \to \infty} n \mathbb{P}(X_0 > u_n \geq M_{1,k-1}, M_{k,\lfloor n/k_n \rfloor - 1} > u_n) = 0, \qquad (3.2.2)$$

where $M_{i,j} := -\infty$ for $i > j$, $M_{i,j} := \max\{X_i, \ldots, X_j\}$ for $i \leq j$.

Clearly (3.2.2) is implied by the condition

$$\lim_{n \to \infty} n \sum_{j=k+1}^{\lfloor n/k_n \rfloor} \mathbb{P}(X_0 > u_n \geq M_{1,k-1}, X_{j-1} > u_n) = 0 \qquad (3.2.3)$$

(see [126, Equation (1.2)]).

This last condition is equal to $D'(u_n)$, when $k = 1$ and to the aforementioned $D''(u_n)$ when $k = 2$.

Together with condition $D(u_n)$, the condition $D^{(k)}(u_n)$ gave an EVL, where $\bar{H}(\tau) = \mathrm{e}^{-\theta\tau}$, with an EI $\theta$ given by O'Brien's formula, whenever the following limit exists:

$$\theta = \lim_{n\to\infty} \theta_n = \lim_{n\to\infty} \frac{\mathbb{P}\left(X_0 > u_n, X_1 \leq u_n, \ldots, X_{k-1} \leq u_n\right)}{\mathbb{P}(X_0 > u_n)}. \qquad (3.2.4)$$

**Theorem 3.2.5** ([126])**.** *Let $X_0, X_1, \ldots$ be a stationary stochastic process and $(u_n)_{n\in\mathbb{N}}$ a sequence satisfying (2.2.2). Suppose $D(u_n)$ holds and $\liminf_{n\to\infty} \mathbb{P}(M_n \leq u_n) > 0$. If for each positive $k$, $D^{(k)}(u_n)$ holds, then*

$$\lim_{n\to\infty}(\mathbb{P}(M_n \leq u_n) - \mathrm{e}^{-\theta_n\tau}) = 0.$$

*Moreover, if the limit (3.2.4) exists then,*

$$\lim_{n\to\infty} \mathbb{P}(M_n \leq u_n) = \mathrm{e}^{-\theta\tau}.$$

## 3.3
## Convergence of Point Processes of Rare Events

A more sophisticated way of studying rare events consists in studying Rare Events Point Processes (REPP). These point processes keep record of the exceedances of the high thresholds $u_n$ by counting the number of such exceedances on a rescaled time interval. For every $A \subset \mathbb{R}$ we define

$$\mathcal{N}_{u_n}(A) := \sum_{i \in A \cap \mathbb{N}_0} \mathbf{1}_{X_i > u_n}.$$

Observe that $\mathcal{N}_{u_n}([0, n))$ counts the number of exceedances amongst the first $n$ observations of the process $X_0, X_1, \ldots, X_n$. Hence, clearly, in the i.i.d. case $\mathcal{N}_{u_n}([0, n))$ is binomial with parameters $(n, \mathbb{P}(X_0 > u_n))$, since it counts the number of successes (exceedances) among $n$ Bernoulli trials being the probability of each success equal to $\mathbb{P}(X_0 > u_n)$. The choice of the normalising sequences $(u_n)_{n\in\mathbb{N}}$ satisfying (2.2.2) means that the average number of successes is nearly constant (converging to $\tau \geq 0$, as $n \to \infty$), which implies that $\mathcal{N}_{u_n}([0, n))$ is asymptotically Poisson. Hence, it is natural to ask whether the exceedance instants form a Poisson Process, in the limit. Since the exceedances occur at integer time values we need to rescale time to obtain such sort of convergence. The REPP introduced describe the record of the times when an exceedance event was observed. The distance between consecutive events in the point processes is the waiting time for the next exceedance.

### 3.3.1
## Definitions and Notation

In order to provide a proper framework of the problem we introduce next the necessary formalism to state the results regarding the convergence of point processes. We recommend the books of Kallenberg [127] and Resnick ( [108, Section 3]) for further reading.

Consider the interval $[0, \infty)$ and its Borel $\sigma$-algebra $\mathcal{B}_{[0,\infty)}$. Let $x_1, x_2, \ldots \in [0, \infty)$ and define

$$\nu = \sum_{i=1}^{\infty} \delta_{x_i},$$

where $\delta_{x_i}$ is the Dirac measure at $x_i$, *i.e.*, for every $A \in \mathcal{B}_{[0,\infty)}$, we have that $\delta_{x_i}(A) = 1$ if $x_i \in A$ or $\delta_{x_i}(A) = 0$, otherwise. The measure $\nu$ is said to be a counting measure on $[0, \infty)$. Let $\mathcal{M}_p([0, \infty))$ be the space of counting measures on $([0, \infty), \mathcal{B}_{[0,\infty)})$. We equip this space with the vague topology, i.e., $\nu_n \to \nu$ in $\mathcal{M}_p([0, \infty))$ whenever $\nu_n(\psi) \to \nu(\psi)$ for any continuous function $\psi : [0, \infty) \to \mathbb{R}$ with compact support. A *point process* $N$ on $[0, \infty)$ is a random element of $\mathcal{M}_p([0, \infty))$, *i.e.*, let $(X, \mathcal{B}_X, \mu)$ be a probability space, then any measurable map $N : X \to \mathcal{M}_p([0, \infty))$ is a point process on $[0, \infty)$.

To give a concrete example of a point process, which in particular will appear as the limit of the REPP, we consider:

**Definition 3.3.1.** Let $T_1, T_2, \ldots$ be an i.i.d. sequence of r.v. with common exponential distribution of mean $1/\theta$. Let $D_1, D_2, \ldots$ be another i.i.d. sequence of r.v., independent of the previous one, and with d.f. $\pi$. Given these sequences, for $J \in \mathcal{B}_{[0,\infty)}$, set

$$N(J) = \int \mathbf{1}_J \, d\left(\sum_{i=1}^{\infty} D_i \delta_{T_1 + \ldots + T_i}\right) = \sum_{i=1}^{\infty} D_i \delta_{T_1 + \ldots + T_i}(J) = \sum_{i=1: \, T_1 + \ldots + T_i \in J}^{\infty} D_i,$$

where $\delta_t$ denotes the Dirac measure at $t > 0$. Let for example $X$ denote the space of all possible realisations of $T_1, T_2, \ldots$ and $D_1, D_2, \ldots$, equipped with a product $\sigma$-algebra and measure, then $N : X \to \mathcal{M}_p([0, \infty))$ is a point process which we call a compound Poisson process of intensity $\theta$ and multiplicity d.f. $\pi$.

*Remark* 3.3.2. Throughout the book, the multiplicity will always be integer valued which means that $\pi$ is completely defined by the values $\pi_k = \mathbb{P}(D_1 = k)$, for every $k \in \mathbb{N}_0$. Note that, if $\pi_1 = 1$ and $\theta = 1$, then $N$ is the standard Poisson process and, for every $t > 0$, the random variable $N([0, t))$ has a Poisson distribution of mean $t$.

As mentioned above, in order to define the REPP we need to rescale time. This rescaling is done so that the average number of exceedances is kept stabilised (converging to 1 on the unit interval $[0, 1]$), which means that, in the independent case, we are counting the number of successes of a Bernoulli trial with nearly constant expected number of successes, which in turn leads to a Poisson behaviour. Hence, we rescale time using the factor $v_n := 1/\mathbb{P}(X_0 > u_n)$, following what Kac's theorem would suggest, see Section 2.3 and detailed discussion in Chapter 4. However, before we give the definition, we need some more formalism. Let $\mathcal{S}$ denote the semi-ring of subsets of $\mathbb{R}_0^+$ whose elements are intervals of the type $[a, b)$, for $a, b \in \mathbb{R}_0^+$. Let $\mathcal{R}$ denote the ring generated by $\mathcal{S}$. Recall that for every $J \in \mathcal{R}$ there are $k \in \mathbb{N}$ and $k$ intervals $I_1, \ldots, I_k \in \mathcal{S}$ such that $J = \cup_{i=1}^k I_j$. In order to fix notation, let $a_j, b_j \in \mathbb{R}_0^+$ be such that $I_j = [a_j, b_j) \in \mathcal{S}$. For $I = [a, b) \in \mathcal{S}$ and $\alpha \in \mathbb{R}$, we denote $\alpha I := [\alpha a, \alpha b)$ and $I + \alpha := [a + \alpha, b + \alpha)$. Similarly, for $J \in \mathcal{R}$ define $\alpha J := \alpha I_1 \cup \cdots \cup \alpha I_k$ and $J + \alpha := (I_1 + \alpha) \cup \cdots \cup (I_k + \alpha)$.

**Definition 3.3.3.** We define the *rare event point process* (REPP) by counting the number of exceedances during the (rescaled) time period $v_n J \in \mathcal{R}$, where $J \in \mathcal{R}$. To be more precise, for every $J \in \mathcal{R}$, set

$$N_n(J) := \mathcal{N}_{u_n}(v_n J) = \sum_{j \in v_n J \cap \mathbb{N}_0} \mathbf{1}_{X_j > u_n}. \tag{3.3.1}$$

We will see that both in the classical theory and in the one developed in this book, the REPP converges either to a standard Poisson process or to a compound Poisson process. For completeness, we define here what we mean by convergence of point processes (see [127] for more details).

**Definition 3.3.4.** Let $(N_n)_{n \in \mathbb{N}} : X \to \mathcal{M}_p([0, \infty))$ be a sequence of point processes defined on a probability space $(X, \mathcal{B}_X, \mu)$ and let $N : Y \to \mathcal{M}_p([0, \infty))$ be another point process defined on a possibly distinct probability space $(Y, \mathcal{B}_Y, \nu)$. We say that $N_n$ converges in distribution to $N$ if $\mu \circ N_n^{-1}$ converges weakly to $\nu \circ N^{-1}$, *i.e.*, for every continuous function $\varphi$ defined on $\mathcal{M}_p([0, \infty))$, we have $\lim_{n \to \infty} \int \varphi d\mu \circ N_n^{-1} = \int \varphi d\nu \circ N^{-1}$. We write $N_n \overset{\mu}{\Longrightarrow} N$.

*Remark* 3.3.5. It can be shown that $(N_n)_{n \in \mathbb{N}}$ converges in distribution to $N$ if the sequence of vector r.v. $(N_n(J_1), \ldots, N_n(J_k))$ converges in distribution to $(N(J_1), \ldots, N(J_k))$, for every $k \in \mathbb{N}$ and all $J_1, \ldots, J_k \in \mathcal{S}$ such that $N(\partial J_i) = 0$ a.s., for $i = 1, \ldots, k$.

Note that

$$\{\mathcal{N}_{u_n}([0, n)) = 0\} = \{M_n \leq u_n\}, \tag{3.3.2}$$

hence the limit distribution of $M_n$ can be easily recovered from the convergence of the REPP.

### 3.3.2
### Absence of Clusters

In the independent case, the exceedances appear scattered on the time line and the waiting times between them are (asymptotically) exponentially distributed. Consequently, it is fairly easy to prove that in this case the REPP $N_n$ defined in (3.3.1) converges in distribution to a standard Poisson process $N$ of intensity 1.

Let's consider general stationary stochastic processes, where conditions $D(u_n)$ and $D'(u_n)$ hold. This implies that there is no clustering of exceedances. In this case, the waiting times between exceedances assume exactly the same exponential pattern as in the independent case. Therefore, not surprisingly, we obtain that the REPP $N_n$ defined in (3.3.1) converges in distribution to a standard Poisson process $N$ of intensity 1. This is the content of [1, Theorem 5.2.1]. The proof relies on a criterion, proposed by Kallenberg [127, Theorem 4.7], which reduces the study of the convergence of the REPP $N_n$ to the distributional properties of the maxima on disjoint multiple time intervals.

### 3.3.3
### Presence of Clusters

on many occasions, the exceedances appear concentrated in the time line in groups or clusters, and such a phenomenon is related to the presence of memory.

As an example, floods persist for days (so that a river gauge reading would give extremely high values for several days in a row); furthermore, having a flood creates a higher propensity for flood in the immediate future because the soil is saturated with water, so that even relatively weak rainfall events can cause inundations.

The time rescaling used to study the convergence of REPP is then responsible for collapsing all the exceedances in a cluster into a single point of mass in the time line, with the weight of the point mass encoding the size of the cluster. Hence, unsurprisingly, the limit of the REPP is a compound Poisson process, which can be thought of as having two components: one is the underlying asymptotic Poisson process governing the positions of the clusters of exceedances; and the other is the multiplicity distribution associated to each such Poisson event, which corresponds to the distribution of the cluster sizes.

Under the presence of clustering the convergence of the REPP is more complicated and it is usually obtained using Laplace transforms and a strengthening of the mixing conditions above is required. We will recall below [128, Theorem 4.2], which assumes the following strengthening of $D(u_n)$:

**Condition** $(\Delta(u_n))$. For $0 \le i \le j \le n$, let $\mathcal{F}_i^j(u_n)$ be the $\sigma$-field generated by the events $\{X_\ell \le u_n\}$, for $i \le \ell \le j$, where $(u_n)_{n \in \mathbb{N}}$ is a sequence of levels satisfying (2.2.2). Let

$$\alpha_{n,t} = \max\{|\mathbb{P}(A \cap B) - \mathbb{P}(A)\mathbb{P}(B)| : A \in \mathcal{F}_1^p(u_n), \ B \in \mathcal{F}_{p+t}^k(u_n), \ 0 \le p \le n-1-t\}.$$

We say that $\Delta(u_n)$ holds for the sequence $X_0, X_1, \ldots$ if $\alpha(n, t_n) \xrightarrow[n \to \infty]{} 0$, for some sequence $t_n = o(n)$.

Let $(k_n)_{n \in \mathbb{N}}$ be a sequence satisfying (3.2.1). We define the cluster size distribution, for each $j \in \mathbb{N}$

$$\pi_n(j) = \mathbb{P}\left(\mathcal{N}([0, n/k_n)) = j | \mathcal{N}([0, n/k_n)) > 0\right) = \frac{\mathbb{P}(\mathcal{N}([0, n/k_n)) = j)}{\mathbb{P}(M_{\lfloor n/k_n \rfloor} \le u_n)}.$$

The next theorem gives the convergence of the REPP to a compound Poisson process under the assumptions $\Delta(u_n)$, the existence of a distributional limit for $M_n$ and the existence of limit of the cluster size distribution $\pi_n$ defined above. The statement is an adjustment of the content of [128, Theorem 4.2], having in mind that the REPP $N_n$ defined in (3.3.1) has a slightly different time rescaling when compared to the one used in [128].

**Theorem 3.3.6.** *Assume that $X_0, X_1, \ldots$ is a stationary stochastic process satisfying condition $\Delta(u_n)$ and $\lim_{n\to\infty} \mathbb{P}(M_n \le u_n) = e^{-\eta}$ for some $\eta > 0$ and a sequence $(u_n)_{n\in n}$ as in (2.2.2). Suppose that there exists a probability distribution $\pi$ defined on $\mathbb{N}$ such that $\pi(j) = \lim_{n\to\infty} \pi_n(j)$ for every $j \in \mathbb{N}$. Then the REPP $N_n$ defined in (3.3.1) converges in distribution to a compound Poisson process $N$, with intensity $\theta = \eta/\tau$ and multiplicity distribution $\pi$.*

We note that by Theorem 3.2.3, the existence of $\eta > 0$ such that $\mathbb{P}(M_n \le u_n) = e^{-\eta}$ implies the existence of an EI $\theta = \eta/\tau$, where $\tau$ comes from the normalising sequence $(u_n)_{n\in n}$ satisfying (2.2.2). Moreover - excluding pathological cases of little interest here - it can be shown that the EI is equal to the inverse of the average cluster size or, in other words,

$$\theta^{-1} = \lim_{n\to\infty} \sum_{j=1}^{\infty} j\pi_n(j). \tag{3.3.3}$$

When clustering is present, the REPP $N_n$ converges to a compound Poisson process with intensity $\theta < 1$, which means that the cluster positions, which account for the Poisson events of the compound Poisson process, appear less frequently than in the absence of clustering, where the exceedances (which corresponded to the Poisson events) appeared at a frequency equal to 1. However, each cluster corresponding to a Poisson event has a multiplicity equal to the cluster size, whose average equals $\theta^{-1}$ so that, in the end, counting with the multiplicities, we still have an average frequency of exceedances equal to 1.

*Remark* 3.3.7. In the applications we will make to dynamical systems, the periodicity will impose some sort of Markovian property on the stochastic processes, which will translate in the fact that clustering $\pi$ is actually a geometric distribution of parameter $\theta \in (0, 1]$, *i.e.*, $\pi_k = \theta(1-\theta)^{k-1}$, for every $k \in \mathbb{N}_0$. This means that, as in [129] the random variable $N([0, t))$ follows a Pólya-Aeppli distribution, *i.e.*:

$$\mathbb{P}(N([0, t)) = k) = e^{-\theta t} \sum_{j=1}^{k} \theta^j (1-\theta)^{k-j} \frac{(\theta t)^j}{j!} \binom{k-1}{j-1},$$

for all $k \in \mathbb{N}$ and $\mathbb{P}(N([0, t)) = 0) = e^{-\theta t}$.

## 3.4
## Elements of Declustering

It is often recognised that the BM method presents some disadvantages when compared to the POT method to estimate the tail of the distribution, since, in the latter, there seems to a better use of the information available on the tail. On the other hand, when clustering is present, BM appears to be more robust because clustering does not change the type of limiting distribution for the maxima. Instead, the estimates of the GPD parameters, an in particular, of the shape parameter $\xi$, obtained using the POT approach can be seriously biased [3]. However, if the clusters are identified and suitable post-processing of the data is performed - the so-called *declustering* - the POT method can again become an efficient statistical inference tool.

Declustering methods are associated with corresponding estimators for the EI. The two main classical procedures to identify clusters are the so-called *blocks declustering* and *runs declustering* [50]:

- the blocks declustering consists of splitting data as described in Section 3.2.1 into $k_n$ blocks of size $\lfloor n/k_n \rfloor$ and considering that the exceedances within each block belong to the same cluster;
- the runs declustering method consists of choosing a run length $q_n$ and stipulating that any two exceedances separated from each other by a time gap smaller than $q_n$ belong to the same cluster.

After the declustering procedure, one retains only one peak (typically the largest one, but other choices are in principle possible) per cluster and performs the POT statistical inference procedure on such reduced dataset. As a result, the estimates of the GPD parameters are more uncertain than before the post-processing, but they are asymptotically unbiased. Additionally, these two declustering procedures can be used to estimate the EI as in [130]. The idea is to use (3.3.3), which says that the EI is the reciprocal of the average number of exceedances in a cluster. Hence, a natural estimation consists in taking the ratio of the number of clusters detected through the declustering procedure to the total number of exceedances. We recall first the so-called blocks estimator:

$$\hat{\theta}_n^B := \frac{\sum_{i=1}^{k_n} \mathbf{1}_{\mathcal{N}([(i-1)\lfloor n/k_n \rfloor, i\lfloor n/k_n \rfloor]) > 0}}{\sum_{i=1}^n \mathbf{1}_{X_i > u_n}}.$$

The runs estimator is given by

$$\hat{\theta}_n^R := \frac{\sum_{i=1}^n \mathbf{1}_{\{X_i > u_n, \mathcal{N}([i+1, i+q_n)) = 0\}}}{\sum_{i=1}^n \mathbf{1}_{X_i > u_n}}.$$

A problem of these two estimators is the sensitivity to the choices of $k_n$ and $q_n$, respectively. More recently, more robust declustering methods and more sophisticated estimation of the EI have been pursued, *e.g.*, in [131, 132, 133]. Some practical methodologies can be found in *e..g.* [3, 50].

# 4

# Emergence of Extreme Value Laws for Dynamical Systems

One of the main purposes of this book is the establishment of a theory of extreme values for dynamical systems. Chaotic dynamical systems exhibit erratic behaviour difficult to anticipate, which can be better understood from a probabilistic point of view. In many practical situations, such as in the now classical case of weather and climate dynamics, time series of observables can be modelled as resulting from a chaotic dynamical system, which describes its time evolution. For this reason, the development of a theory of extreme values for such systems opens up a huge pool of potential applications as well as embodying a new ground for the further development of the already existing theory.

The starting point of our analysis are stationary stochastic processes arising from a chaotic dynamical systems. One way to go in this direction is to evaluate an observable function along the orbit of the system. We focus on two qualitatively distinct kinds of functions, the *distance* and the *physical* observables. Then the extremal behaviour is analysed by studying the distributional limit of the partial maximum of such stochastic processes or by investigating the limit of point processes counting the number of exceedances of high thresholds during some time interval.

In the context of dynamical systems, the study of EVLs is a quite recent topic. It first appeared in the pioneering paper of Collet, [72], which has been an inspiration for plenty of the research on this issue. Then the subject has been further addressed and developed in many subsequent contributions including [134, 73, 135, 74, 136, 137, 138, 139, 49, 140, 81, 46, 76, 77, 78, 141, 142, 79, 44].

The classical theory for general stationary stochastic processes described in the previous chapter is based on the mixing condition $D(u_n)$ (or very similar ones), which can hardly be verified for processes arising from dynamical systems (other than trivial examples), because the most natural way to analyse the mixing properties is usually through the study of the rates of decay of correlations. Based on ideas introduced by Collet, a strengthening of the classical theory has been achieved first by introducing some weaker assumptions to replace the Leadbetter's $D(u_n)$ in order to obtain EVLs. We emphasise that these results apply to general stationary stochastic processes, although these efforts were motivated by application to the study of statistical properties of dynamical systems.

The structure of this chapter is as follows. In Section 4.1, we start by enhancing the EVT for general stationary stochastic processes under weaker mixing conditions that we will denote with a Cyrillic D: Д. These conditions are weaker than the original $D(u_n)$ condition of Leadbetter and and that will be essential in order to apply the condition to dynamical systems. Moreover, the asymptotic results are obtained with an indication of the respective convergence rates. In Section 4.2, we finally apply the theory developed in Section 4.1 to stationary processes arising from dynamical systems. In particular, we show that if a system has a sufficiently fast decay of correlations then there exists a dichotomy regarding the extremal behaviour, which establishes a connection between presence of clusters and periodicity. This is illustrated with a toy model of a uniformly expanding system: the special case of the Bernoulli shift maps known as doubling map. In Section 4.3, we enhance the study by considering the convergence of Rare Events Point Processes, which is much stronger than obtaining EVL. We also obtain a dichotomy in this case and apply this to the example mentioned before. We defer the more detailed discussion of the rate of convergence to the asymptotic EVL to Chap. 6 and specifically to Sect. 6.8.

## 4.1
## Extremes for General Stationary Processes – an Upgrade Motivated by Dynamics

The theory developed by Leadbetter and later by several other authors for stationary stochastic processes is not adequate for the study of the extremal behaviour of stochastic processes arising from dynamical systems. The reason is that the condition $D(u_n)$ requires a uniform mixing property on the dependence structure of the processes, which can hardly be met within a dynamical systems setting. Therefore, a substantial revision of the existing theory for stationary processes is needed in order to be able to prove the existence of EVL considering milder conditions on the dependence structure of the processes. This has been achieved first in the absence of clustering in [72, 134, 73, 139]. The presence of clustering brings some extra difficulty in pursuing the weakening of condition $D(u_n)$ of Leadbetter, eventually solved in [49]. In this series of papers, several mixing conditions have been devised in order to weaken $D(u_n)$ and still be able to obtain EVL for stochastic processes arising from dynamical systems. We present below the latest refinement of these conditions, which has been obtained in [143] and which, in particular, allows for addressing both the absence and presence of clusters at once.

### 4.1.1
### Notation

Let $X_0, X_1, \ldots$ be a stationary stochastic process, where each r.v. $X_i : \mathcal{Y} \to \mathbb{R}$ is defined on the measure space $(\mathcal{Y}, \mathcal{B}, \mathbb{P})$. We assume that $\mathcal{Y}$ is a sequence space with a natural product structure so that each possible realisation of the stochastic process corresponds to a unique element of $\mathcal{Y}$ and there exists a measurable map $T : \mathcal{Y} \to \mathcal{Y}$, the time evolution map, which can be seen as the passage of one unit of time, so that

$$X_{i-1} \circ T = X_i, \quad \text{for all } i \in \mathbb{N}.$$

*Note.* There is an obvious relation between $T$ and the *shift* map, but we do not want to stress this relation here, because the description considered here is definitely not the usual shift dynamics. The normal shift map acts on sequences from a finite or countable alphabet, while here $T$ acts on spaces like $\mathbb{R}^{\mathbb{N}}$, in the sense that the sequences can be thought as being obtained from an alphabet like $\mathbb{R}$. In fact, our applications include non uniformly hyperbolic systems and ball target sets. $T$ can also be thought as a discretisation of a continuous time system. Note that our framework is so general that any discrete time dynamical systems is automatically of this form, with the alphabet being the state space of the dynamical system.

Stationarity means that $\mathbb{P}$ is $T$-invariant. Note that $X_i = X_0 \circ T^i$, for all $i \in \mathbb{N}_0$, where $T^i$ denotes the $i$-fold composition of $T$, with the convention that $T^0$ denotes the identity map on $\mathcal{Y}$.

Following Sect. 2.1, we denote by $F$ the d.f. of $X_0$, *i.e.*, $F(x) = \mathbb{P}(X_0 \leq x)$. Given any d.f. $F$, let $\bar{F} = 1 - F$ and let $u_F$ denote the right endpoint of the d.f. $F$, *i.e.*, $u_F = \sup\{x : F(x) < 1\}$. We say we have an *exceedance* of the threshold $u < u_F$ at time $j \in \mathbb{N}_0$ whenever $\{X_j > u\}$ occurs. We define a new sequence of random variables $M_1, M_2, \ldots$ given by (2.2.1)

In what follows for every $A \in \mathcal{B}$, we denote the complement of $A$ as $A^c := \mathcal{X} \setminus A$. For some $u \in \mathbb{R}$, $q \in \mathbb{N}$, we define the events:

$$U(u) := \{X_0 > u\},$$

$$A^{(q)}(u) := U(u) \cap \bigcap_{i=1}^{q} T^{-i}(U(u)^c) = \{X_0 > u, X_1 \leq u, \ldots, X_q \leq u\}.$$

$$\text{(4.1.1)}$$

where $A^{(q)}(u)$ corresponds to the case where we have an extreme event at time zero that is not followed by another one up to time $t = q$. This is a condition clearly pointing to the absence of clusters. We also set for convenience $A^{(0)}(u) := U(u)$, $U_n := U(u_n)$ and $A_n^{(q)} := A^{(q)}(u_n)$, for all $n \in \mathbb{N}$ and $q \in \mathbb{N}_0$. Let

$$\theta_n := \frac{\mathbb{P}\left(A_n^{(q)}\right)}{\mathbb{P}(U_n)}. \quad \text{(4.1.2)}$$

Note that $0 \leq \theta_n \leq 1$. Let now $B \in \mathcal{B}$ be an event. For some $s \geq 0$ and $\ell \geq 0$, we define:

$$\mathcal{W}_{s,\ell}(B) = \bigcap_{i=\lfloor s \rfloor}^{\lfloor s \rfloor + \max\{\lfloor \ell \rfloor - 1, \, 0\}} T^{-i}(B^c). \quad \text{(4.1.3)}$$

The notation $T^{-i}$ is used for the preimage by $T^i$. We will write $\mathcal{W}_{s,\ell}^c(B) := (\mathcal{W}_{s,\ell}(B))^c$. Whenever is clear or unimportant which event $B \in \mathcal{B}$ applies, we will drop the $B$ and write just $\mathcal{W}_{s,\ell}$ or $\mathcal{W}_{s,\ell}^c$. Observe that

$$\mathcal{W}_{0,n}(U(u)) = \{M_n \leq u\} \qquad \text{and} \qquad T^{-1}(\mathcal{W}_{0,n}(B)) = \{r_B > n\}. \quad \text{(4.1.4)}$$

where $r_B(x)$ is the first hitting time defined in Eq. 2.3.1. Also observe that $\mathcal{W}_{0,n}$ has the following interpretation in terms of the sequence $X_0, X_1, X_2, \ldots$, namely $T^i(x) \notin A_n^{(q)}$ means that $X_i(x) \leq u_n$ or $X_{i+j}(x) > u_n$ for some $j = 1, \ldots, q$.

### 4.1.2
### The New Conditions

When the stochastic processes arise from dynamical systems as described in Section 4.2 below, condition $D(u_n)$ cannot be verified using the usual available information about mixing rates of the system except in some very special situations, and even then only for certain subsequences of $n$. This means that the theory developed by Leadbetter and others and discussed in Chapter 3 is not of practical utility in the dynamical systems context. For that reason, motivated by the work of Collet ([72]), Freitas and Freitas have proposed in [73] a new condition called $D_2(u_n)$ for general stationary stochastic processes, which imposes a much weaker uniformity requirement than $D(u_n)$, and, together with $D'(u_n)$, admits a proof of the existence of EVL in the absence of clustering (with $\theta = 1$).

The big advantage of $D_2(u_n)$ over $D(u_n)$ is that in $D(u_n)$ one has to establish this rate uniformly with respect to the size of both blocks whereas in $D_2(u_n)$ the rate has only to be uniform with respect to the size of one block. The great advantage is that using such a weaker condition the EVL can be established easily for dynamical systems and stochastic processes with sufficiently fast decay of correlations and hence EVLs can be easily shown to apply for a much larger class of systems. We remark that the original $D(u_n)$ has actually never been verified for any dynamical system, when considering the observables introduced later in Sect. 4.2.1. In the argument of [73], this weakening has been achieved by a fuller application of condition $D'(u_n)$.

In [49], the authors proved a connection between periodicity and presence of clusters. Motivated by the behaviour at periodic points, which led to the appearance of clusters of exceedances, the authors proposed new conditions in order to prove the existence of EVLs with EI less than 1. The main idea is that, under a condition $SP_{p,\theta}(u_n)$, which imposed some sort of periodic behaviour of period $p$ on the structure of general stationary stochastic processes, the sequences $\mathbb{P}(M_n \leq u_n)$ and $\mathbb{P}(\mathcal{W}_{0,n}(A^{(p)}(u_n)))$ share the same limit (see [49, Proposition 1]). Then the strategy has been to prove the existence of a limit for $\mathbb{P}(\mathcal{W}_{0,n}(A^{(p)}(u_n)))$, which has been achieved under conditions $D_p(u_n)$ and $D'_p(u_n)$. These latter conditions can be seen as being obtained from $D_2(u_n)$ in [134] and the original $D'(u_n)$, respectively, by replacing the role of exceedances $\{X_j > u_n\}$ by that of *escapes*, which correspond to the event $\{X_j > u_n, X_{j+p} \leq u_n\}$.

In [142] discontinuity points create two periodic types of behaviour (with possibly different periods) on the structure of the stochastic processes, so some further adjustments to conditions $D^p(u_n)$ and $D'_p(u_n)$ were needed.

We remark that in all the cases above the main advantage of the conditions $D_2(u_n)$, $D^p(u_n)$ is that they are much weaker than the original uniformity requirement imposed by $D(u_n)$ and, unlike $D(u_n)$, they all follow from sufficiently fast decay of correlations of the system.

While developing the techniques in [143] to sharpen the error terms, it has been necessary to improve the estimates in [49, Proposition 1] (this is done in Proposition 4.1.12 below). One consequence of this is that the authors were then able to essentially remove condition $SP_{p,\theta}(u_n)$. Whence we developed two conditions that refine all the previous conditions to obtain EVLs , in the way that they simultaneously include the cases of absence and presence of clustering and, on the other hand, to combine all the scenarios considered before with no periodic behaviour, with simple periodic behaviour or multiple types of periodic behaviour.

As seen in the historical discussion above, the notation of the condition $D$ is hampered with sub- and superscripts. In order to simplify the notation, here we follow [143] and employ instead a cyrillic D, i.e., Д.

**Condition ($Д_q(u_n)$).** We say that $Д_q(u_n)$ holds for the sequence $X_0, X_1, \ldots$ if for every $\ell, t, n \in \mathbb{N}$,

$$\left| \mathbb{P}\left( A_n^{(q)} \cap \mathcal{W}_{t,\ell}\left( A_n^{(q)} \right) \right) - \mathbb{P}\left( A_n^{(q)} \right) \mathbb{P}\left( \mathcal{W}_{0,\ell}\left( A_n^{(q)} \right) \right) \right| \leq \gamma(q, n, t), \quad (4.1.5)$$

where $\gamma(q, n, t)$ is decreasing in $t$ for each $n$ and, there exists a sequence $(t_n)_{n \in \mathbb{N}}$ such that $t_n = o(n)$ and $n\gamma(q, n, t_n) \to 0$ when $n \to \infty$.

*Remark* 4.1.1. Note that the new condition $Д_q(u_n)$ imposes a condition in the first block only on $q$ random variables. On the contrary, the original condition $D(u_n)$ does require a bound independent of the number of variables considered in the first block. This is the crucial difference that makes it possible to prove $Д_q(u_n)$ easily from decay of correlations of the underlying stochastic processes, in contrast to $D(u_n)$, which is not possible to be verified even in very simple situations. The weakening of the uniformity imposed by $D(u_n)$ came at price on the rate function: while for $D(u_n)$, we need $\lim_{n\to\infty} k_n \alpha(n, t_n) = 0$, for $Д_q(u_n)$, we need $\lim_{n\to\infty} n\gamma(q, n, t_n) = 0$. However, this is a very small price to pay because, when the stochastic processes arise from dynamical systems, the verification of $Д_q(u_n)$ means that we need decay of correlations having at least at a summable rate. But this is precisely the regime where one can prove Central Limit Theorems. Hence, even though we have a slight strengthening in the required mixing rate, when we compare $Д_q(u_n)$ to $D(u_n)$, the weakening on the uniformity is so much more important that we believe it is fair to say $Д_q(u_n)$ is considerably weaker than $D(u_n)$. This is substantiated by the fact that $Д_q(u_n)$ can be verified in a huge range of examples arising from dynamical systems, where condition $D(u_n)$ simply cannot be established.

*Remark* 4.1.2. Let's try to provide a more heuristic explanation of the condition $Д_q(u_n)$. What we are imposing is a sort of asymptotic independence between the occurrence of an event $A^{(q)}(u)$ (*i.e.* an exceedance of threshold at time zero *not* followed by any other one within the next $q$ time steps) and no occurrences of such an event for a long period of length $l$, starting from time $t$, which introduces a time gap. In other words, if one observes after a cluster for $q$ time steps no other exceedance then the time of the next exceedance is essentially an independent event. Instead, in the $D(u_n)$ condition, the role of $A^{(q)}(u)$ is played by simple occurrences $U(u_n)$, with the disadvantage that both blocks separated by the time gap $t$ may be arbitrarily large.

For some fixed $q \in \mathbb{N}_0$, consider the sequence $(t_n)_{n \in \mathbb{N}}$, given by condition $\unicode{x2A3F}(u_n)$ and let $(k_n)_{n \in \mathbb{N}}$ be another sequence of integers such that

$$k_n \to \infty \quad \text{and} \quad k_n t_n = o(n). \tag{4.1.6}$$

**Condition** $(\unicode{x2A3F}'_q(u_n))$**.** We say that $\unicode{x2A3F}'_q(u_n)$ holds for the sequence $X_0, X_1, X_2, \ldots$ if there exists a sequence $(k_n)_{n \in \mathbb{N}}$ satisfying (4.1.6) and such that

$$\lim_{n \to \infty} n \sum_{j=q+1}^{\lfloor n/k_n \rfloor - 1} \mathbb{P}\left( A_n^{(q)} \cap T^{-j}\left( A_n^{(q)} \right) \right) = 0. \tag{4.1.7}$$

*Remark* 4.1.3. Note that condition $\unicode{x2A3F}'_q(u_n)$ is very similar to condition $D^{(q+1)}(u_n)$ f[126, Equation (1.2)] (see (3.2.3)). Since $T^{-j}\left( A_n^{(q)} \right) \subset \{X_j > u_n\}$, it is slightly weaker than $D^{(q+1)}(u_n)$ [126], but in the applications considered that does not make any difference. Note that if $q = 0$ then we get back condition $D'(u_n)$ from Leadbetter.

### 4.1.3

### The Existence of EVL for General Stationary Stochastic Processes under Weaker Hypotheses

Since condition $\unicode{x2A3F}_q(u_n)$ is much weaker than the original $D(u_n)$ of Leadbetter, in the sense explained in Remark 4.1.1, then Theorem 4.1.4 can be seen, in particular, as a generalisation of [126, Corollary 1.3].

**Theorem 4.1.4.** *Let* $X_0, X_1, \ldots$ *be a stationary stochastic process and* $(u_n)_{n \in \mathbb{N}}$ *a sequence satisfying (2.2.2), for some* $\tau > 0$. *Assume that conditions* $\unicode{x2A3F}_q(u_n)$ *and* $\unicode{x2A3F}'_q(u_n)$ *hold for some* $q \in \mathbb{N}_0$, *and* $(t_n)_{n \in \mathbb{N}}$ *and* $(k_n)_{n \in \mathbb{N}}$ *are the sequences in those conditions. Then, there exists* $C > 0$ *such that for all* $n$ *large enough we have*

$$\left| \mathbb{P}(M_n \le u_n) - e^{-\theta_n \tau} \right| \le C \left[ k_n t_n \frac{\tau}{n} + n\gamma(q, n, t_n) \right.$$
$$+ n \sum_{j=q+1}^{\lfloor n/k_n \rfloor - 1} \mathbb{P}\left( A_n^{(q)} \cap T^{-j}\left( A_n^{(q)} \right) \right)$$
$$\left. + e^{-\theta_n \tau}\left( |\tau - n\mathbb{P}\left( U_n \right)| + \frac{\tau^2}{k_n} \right) + q\mathbb{P}\left( U_n \setminus A_n^{(q)} \right) \right],$$

*where* $\theta_n$ *is given by equation* (4.1.2).

In case the limit $\theta = \lim_{n \to \infty} \theta_n$ exists, where $\theta_n$ is as in (4.1.2), then we can use the previous result to obtain:

**Corollary 4.1.5.** *Let $X_0, X_1, \ldots$ be a stationary stochastic process and $(u_n)_{n \in \mathbb{N}}$ a sequence satisfying (2.2.2), for some $\tau > 0$. Assume that conditions $\underline{\underline{\mathrm{D}}}_q(u_n)$ and $\underline{\underline{\mathrm{D}}}'_q(u_n)$ hold for some $q \in \mathbb{N}_0$, and $(t_n)_{n \in \mathbb{N}}$ and $(k_n)_{n \in \mathbb{N}}$ are the sequences in those conditions. Moreover assume that the limit in (3.2.4) exists. Then, there exists $C > 0$ such that for all $n \in \mathbb{N}$ we have*

$$
\begin{aligned}
\left| \mathbb{P}(M_n \leq u_n) - \mathrm{e}^{-\theta \tau} \right| \leq C \Bigg[ & k_n t_n \frac{\tau}{n} + n \gamma(q, n, t_n) \\
& + n \sum_{j=1}^{\lfloor n/k_n \rfloor} \mathbb{P}\left( A_n^{(q)} \cap T^{-j}\left( A_n^{(q)} \right) \right) \\
& + \mathrm{e}^{-\theta \tau} \left( |\tau - n\mathbb{P}(U_n)| + \frac{\tau^2}{k_n} + |\theta_n - \theta| \, \tau \right) + q \mathbb{P}\left( U_n \setminus A_n^{(q)} \right) \Bigg],
\end{aligned}
$$

*where $\theta_n$ is given by equation (4.1.2) and $\theta = \lim_{n \to \infty} \theta_n$.*

This criterium can be effectively used for general stochastic processes as well as to a large class of dynamical systems such as those studied in [74, 136, 49, 144].

*Remark 4.1.6.* Note that the estimates of Theorem 4.1.4 and Corollary 4.1.5 hold under $\underline{\underline{\mathrm{D}}}_q(u_n)$ alone. However, if $\underline{\underline{\mathrm{D}}}'_q(u_n)$ does not hold, the upper bound is useless since the third term on the right hand side would not converge to 0 as $n \to \infty$.

Observe that a direct consequence of the previous result is the following:

**Corollary 4.1.7.** *Let $X_0, X_1, \ldots$ be a stationary stochastic process and $(u_n)_{n \in \mathbb{N}}$ a sequence satisfying (2.2.2), for some $\tau > 0$. Assume that conditions $\underline{\underline{\mathrm{D}}}_q(u_n)$ and $\underline{\underline{\mathrm{D}}}'_q(u_n)$ hold for some $q \in \mathbb{N}_0$, and $(t_n)_{n \in \mathbb{N}}$ and $(k_n)_{n \in \mathbb{N}}$ are the sequences in those conditions. Moreover assume that the limit in (3.2.4) exists. Then,*

$$
\lim_{n \to \infty} \mathbb{P}(M_n \leq u_n) = \mathrm{e}^{-\theta \tau}.
$$

*Remark 4.1.8.* The convergence result is based on the *blocking argument* as several quantities are constructed according to it. Namely, the number of blocks taken ($k_n$) and the size of the gaps between the blocks ($t_n$) have to satisfy (4.1.6). The first error term, depending on the choices for adequate $k_n$ and $t_n$, typically, decays like $n^{-\delta}$, for some $0 < \delta < 1$. The second term depends on the long range mixing rates of the process ($\underline{\underline{\mathrm{D}}}(u_n)$). The third term takes into account the short range recurrence properties ($\underline{\underline{\mathrm{D}}}'(u_n)$). The fourth has three components, the first depends on the asymptotics of relation (2.2.2), the third on the asymptotics of (3.2.4) and the second on the number of blocks, which must be traded off with the first term. Note that the term $\mathrm{e}^{-\theta \tau} \frac{\tau^2}{k_n}$ also appears in the i.i.d. case since it results from expansion (4.1.8) below. The fifth term results from replacing $U_n$ by $A_n^{(q)}$ (see Proposition 4.1.12) and should decay like $1/n$. The constant $C$ may depend on the rates just mentioned but not on $\tau$.

*Remark* 4.1.9. As observed above, in the special case $q = 0$, we have that condition $Ⅎ'_0(u_n)$ coincides with the original condition $D'(u_n)$ from Leadbetter, which means there are no clusters of exceedances and it is straightforward to check that by definition of $A_0^{(0)} = U_n$ that the limit in (3.2.4) trivially exists and is equal to 1.

*Remark* 4.1.10. In this approach, it is rather important to observe the prominent role played by condition $Ⅎ'_q(u_n)$. In particular, note that if condition $Ⅎ'_q(u_n)$ holds for some particular $q = q_0 \in \mathbb{N}_0$, then condition $Ⅎ'_q(u_n)$ holds for all $q \geq q_0$, which also implies that if the limit in (3.2.4) exists for such $q = q_0$ it will also exist for all $q \geq q_0$ and takes always the same value. This suggests that in trying to find the existence of EVL, one needs to find the smallest value of $q = q_0$ such that the condition $Ⅎ'_q(u_n)$ holds.

*Remark* 4.1.11. Also note that the verification of conditions $Ⅎ_q(u_n)$ and $Ⅎ'_q(u_n)$ for certain stochastic process and, in particular, for the ones arising from dynamical systems as described in the following sections, works as a sort of hunting license to use the block maxima method or the POT approach to make statistical inference on the tail of the considered distributions. This is clear in the case $q = 0$ but also when $q > 0$, since the fact that every cluster must end when there are at least $q$ observations with no exceedances, just after the latest exceedance. This helps in identifying clusters in order to perform a data declustering so that we can use the POT and consequent GPD analysis.

### 4.1.4
### Proofs of Theorem 4.1.4 and Corollary 4.1.5

In the following we construct the proofs of Theorem 4.1.4 and Corollary 4.1.5.

The following result gives a simple estimate but a rather important one. It is crucial in removing condition $SP_{p,\theta}(u_n)$ from [49], to present in a unified way the results under the presence and absence of clustering in Theorem 4.1.4 and guarantee that one can replace exceedances by the occurrence of $A_n^{(q)}$.

**Proposition 4.1.12.** *Given an event $B \in \mathcal{B}$, let $q, n \in \mathbb{N}$ be such that $q < n$ and define $A = B \setminus \bigcup_{j=1}^q T^{-j}(B)$. Then*

$$|\mathbb{P}(\mathcal{W}_{0,n}(B)) - \mathbb{P}(\mathcal{W}_{0,n}(A))| \leq \sum_{j=1}^q \mathbb{P}\left(\mathcal{W}_{0,n}(A) \cap T^{-n+j}(B \setminus A)\right).$$

*Proof.* Since $A \subset B$, then clearly $\mathcal{W}_{0,n}(B) \subset \mathcal{W}_{0,n}(A)$. Hence, we have to estimate the probability of $\mathcal{W}_{0,n}(A) \setminus \mathcal{W}_{0,n}(B)$ which corresponds to the set of points that at some time before $n$ enter $B$ but never enter its subset $A$.

Let $x \in \mathcal{W}_{0,n}(A) \setminus \mathcal{W}_{0,n}(B)$. Then $T^j(x) \in B$ for some $j = 1, \ldots, n-1$ but $T^j(x) \notin A$ for all such $j$. We will see that there exists $j \in \{1, \ldots, q\}$ such that $T^{n-j}(x) \in B$. In fact, suppose that no such $j$ exists. Then let $\ell = \max\{i \in \{1, \ldots, n-1\} : T^i(x) \in B\}$ be the last moment the orbit of $x$ enters $B$ during the time period in question. Then, clearly, $\ell < n - q$. Hence, if $T^i(x) \notin B$, for all $i = \ell + 1, \ldots, n - 1$ then we must have that $T^\ell(x) \in A$ by definition of $A$. But this contradicts the fact that $x \in \mathcal{W}_{0,n}(A)$. Consequently, we have that there exists $j \in \{1, \ldots, q\}$ such that $T^{n-j}(x) \in B$ and since $x \in \mathcal{W}_{0,n}(A)$ then we can actually write $T^{n-j}(x) \in B \setminus A$.

This means that $\mathcal{W}_{0,n}(A) \setminus \mathcal{W}_{0,n}(B) \subset \bigcup_{j=1}^{q} T^{-n+j}(B \setminus A) \cap \mathcal{W}_{0,n}(A)$ and then

$$\left| \mathbb{P}(\mathcal{W}_{0,n}(B)) - \mathbb{P}(\mathcal{W}_{0,n}(A)) \right| = \mathbb{P}(\mathcal{W}_{0,n}(A) \setminus \mathcal{W}_{0,n}(B))$$

$$\leq \mathbb{P} \left( \bigcup_{j=1}^{q} T^{-n+j}(B \setminus A) \cap \mathcal{W}_{0,n}(A) \right) \leq \sum_{j=1}^{q} \mathbb{P} \left( \mathcal{W}_{0,n}(A) \cap T^{-n+j}(B \setminus A) \right),$$

as required. □

In what follows we will need the error term of the limit expression $\lim_{n \to \infty} \left( 1 + \frac{x}{n} \right)^n = \mathrm{e}^x$, namely,

$$\left( 1 + \frac{x}{n} \right)^n = \mathrm{e}^x \left( 1 - \frac{x^2}{2n} + \frac{x^3(8 + 3x)}{24n^2} + O\left( \frac{1}{n^3} \right) \right), \tag{4.1.8}$$

which holds uniformly for $x$ on bounded sets. Also, by Taylor's expansion, for every $\delta \in \mathbb{R}$ and $x \in \mathbb{R}$ we have

$$\left| \mathrm{e}^{x+\delta} - \mathrm{e}^x \right| \leq \mathrm{e}^x \left( |\delta| + \mathrm{e}^{|\delta|} \delta^2 / 2 \right). \tag{4.1.9}$$

The strategy is to use a blocking argument, that goes back to Markov, which consists of splitting the data into blocks with gaps of increasing length. There are three main steps. The first step is to estimate the error produced by neglecting the data corresponding to the gaps. The second is to use essentially the mixing condition $Ð_q(u_n)$ (with some help from $Ð'_q(u_n)$) to show that the probability of the event corresponding to the global maximum being less than some threshold $u_n$ can be approximated by the product of the probabilities of the maxima within each block being less than $u_n$. The idea is that the gaps make the maxima in each block become practically independent from each other. The last step is to use condition $Ð'_q(u_n)$ to estimate the probability of the maximum within a block being smaller than $u_n$.

Next, we state a couple of lemmas and a proposition that give the main estimates regarding the use of a blocking argument. Let us first bound the effect of ignoring the exceedance events in the gaps.

**Lemma 4.1.13.** *For any fixed $A \in \mathcal{B}$ and $s, t', m \in \mathbb{N}$, we have:*

$$\left| \mathbb{P}(\mathcal{W}_{0,s+t'+m}(A)) - \mathbb{P}(\mathcal{W}_{0,s}(A) \cap \mathcal{W}_{s+t',m}(A)) \right| \leq t' \mathbb{P}(A).$$

*Proof.* Using stationarity we have

$$\mathbb{P}(\mathcal{W}_{0,s} \cap \mathcal{W}_{s+t',m}) - \mathbb{P}(\mathcal{W}_{0,s+t'+m}) = \mathbb{P}(\mathcal{W}_{0,s} \cap \mathcal{W}_{s,t'}^c \cap \mathcal{W}_{s+t',m})$$

$$\leq \mathbb{P}(\mathcal{W}_{0,t'}^c) = \mathbb{P}(\cup_{j=0}^{t'-1} T^{-j}(A))$$

$$\leq \sum_{j=0}^{t'-1} \mathbb{P}(T^{-j}(A)) = t' \mathbb{P}(A).$$

$\square$

**Lemma 4.1.14.** *For any fixed* $A \in \mathcal{B}$ *and integers* $s, t, m$, *we have:*

$$|\mathbb{P}(\mathcal{W}_{0,s}(A) \cap \mathcal{W}_{s+t,m}(A)) - (1 - s\mathbb{P}(A))\mathbb{P}(\mathcal{W}_{0,m}(A))| \leq$$

$$\left| s\mathbb{P}(A)\mathbb{P}(\mathcal{W}_{0,m}(A)) - \sum_{j=0}^{s-1} \mathbb{P}(A \cap \mathcal{W}_{s+t-j,m}(A)) \right| + 2s \sum_{j=1}^{s-1} \mathbb{P}(A \cap T^{-j}(A)).$$

*Proof.* Observe that the first term in the bound is measuring the mixing across the gap t and the second term is measuring the probability that two events $A$ appear in the first block. Adding and substracting and using the triangle inequality we obtain that

$$\left| \mathbb{P}(\mathcal{W}_{0,s} \cap \mathcal{W}_{s+t,m}) - \mathbb{P}(\mathcal{W}_{0,m})(1 - s\mathbb{P}(A)) \right| \leq$$

$$\left| s\mathbb{P}(A)\mathbb{P}(\mathcal{W}_{0,m}(A)) - \sum_{j=0}^{s-1} \mathbb{P}(A \cap \mathcal{W}_{s+t-j,m}(A)) \right| +$$

$$+ \left| \mathbb{P}(\mathcal{W}_{0,s} \cap \mathcal{W}_{s+t,m}) - \mathbb{P}(\mathcal{W}_{0,m}) + \sum_{j=0}^{s-1} \mathbb{P}(A \cap \mathcal{W}_{s+t-j,m}) \right|. \quad (4.1.10)$$

Regarding the second term on the right, by stationarity, we have

$$\mathbb{P}(\mathcal{W}_{0,s} \cap \mathcal{W}_{s+t,m}) = \mathbb{P}(\mathcal{W}_{s+t,m}) - \mathbb{P}(\mathcal{W}_{0,s}^c \cap \mathcal{W}_{s+t,m})$$
$$= \mathbb{P}(\mathcal{W}_{0,m}) - \mathbb{P}(\mathcal{W}_{0,s}^c \cap \mathcal{W}_{s+t,m}).$$

Now, since $\mathcal{W}_{0,s}^c \cap \mathcal{W}_{s+t,m} = \cup_{i=0}^{s-1} T^{-i}(A) \cap \mathcal{W}_{s+t,m}$, we have by Bonferroni's inequality that

$$0 \leq \sum_{j=0}^{s-1} \mathbb{P}(A \cap \mathcal{W}_{s+t-j,m}) - \mathbb{P}(\mathcal{W}_{0,s}^c \cap \mathcal{W}_{s+t,m}) \leq$$

$$\sum_{j=0}^{s-1} \sum_{i>j}^{s-1} \mathbb{P}(T^{-j}(A) \cap T^{-i}(A) \cap W_{s+t,m}).$$

Hence, using these last two computations we get:

$$\left| \mathbb{P}(\mathcal{W}_{0,s} \cap \mathcal{W}_{s+t,m}) - \mathbb{P}(\mathcal{W}_{0,m}) + \sum_{j=0}^{s-1} \mathbb{P}(A \cap \mathcal{W}_{s+t-j,m}) \right|$$

$$\leq \sum_{j=0}^{s-1} \sum_{i>j}^{s-1} \mathbb{P}(T^{-j}(A) \cap T^{-i}(A) \cap \mathcal{W}_{s+t,m})$$

$$\leq \sum_{j=0}^{s-1} \sum_{i>j}^{s-1} \mathbb{P}(T^{-j}(A) \cap T^{-i}(A)) \leq s \sum_{j=1}^{s-1} \mathbb{P}(A \cap T^{-j}(A)).$$

The result now follows directly from plugging the last estimate into (4.1.10). $\square$

**Proposition 4.1.15.** *Fix $A \in \mathcal{B}$ and $n \in \mathbb{N}$. Let $\ell, k, t \in \mathbb{N}$ be such that $\ell = \lfloor n/k \rfloor$ and $\ell \mathbb{P}(A) < 1$. We have:*

$$\left| \mathbb{P}(\mathcal{W}_{0,n}(A)) - (1 - \ell\,\mathbb{P}(A))^k \right| \leq 2kt\mathbb{P}(A) + 2n \sum_{j=1}^{\ell-1} \mathbb{P}(A \cap T^{-j}(A))$$

$$+ \sum_{i=0}^{k-1} \left| \ell\mathbb{P}(A)\mathbb{P}(\mathcal{W}_{0,i(\ell+t)}) - \sum_{j=0}^{\ell-1} \mathbb{P}(A \cap \mathcal{W}_{\ell+t-j,i(\ell+t)}) \right|.$$

*Proof.* The basic idea is to split the time interval $[0, n)$ into $k$ blocks of size $\lfloor n/k \rfloor$. Then, using Lemma 4.1.13 we add gaps of size $t$ between the blocks, and next we apply Lemma 4.1.14 recursively until we exhaust all the blocks.

Noting that $0 \leq k(\ell + t) - n \leq kt$ and using Lemma 4.1.13, with $s = n$, $t' = k(\ell + t) - n$, $m = 0$, and setting $\mathcal{W}_{i,0} := \mathcal{X}$, for all $i = 0, 1, 2, \dots$ as well as $\mathcal{W}_{i,n} = \mathcal{W}_{i,n}(A)$ for $n \in \mathbb{N}$, we have:

$$\left| \mathbb{P}(\mathcal{W}_{0,n}) - \mathbb{P}(\mathcal{W}_{0,k(\ell+t)}) \right| \leq kt\mathbb{P}(A). \tag{4.1.11}$$

Using Lemmas 4.1.13 and 4.1.14 we obtain

$$\left| \mathbb{P}(\mathcal{W}_{0,i(\ell+t)}) - (1 - \ell\mathbb{P}(A))\mathbb{P}(\mathcal{W}_{0,(i-1)(\ell+t)}) \right|$$

$$\leq \left| \mathbb{P}(\mathcal{W}_{0,i(\ell+t)}) - \mathbb{P}(\mathcal{W}_{0,\ell} \cap \mathcal{W}_{(\ell+t),(i-1)(\ell+t)}) \right|$$

$$+ \left| \mathbb{P}(\mathcal{W}_{0,\ell} \cap \mathcal{W}_{(\ell+t),(i-1)(\ell+t)}) - (1 - \ell\mathbb{P}(A))\mathbb{P}(\mathcal{W}_{0,(i-1)(\ell+t)}) \right|$$

$$\leq t\mathbb{P}(A) + \left| \ell\mathbb{P}(A)\mathbb{P}(\mathcal{W}_{0,(i-1)(\ell+t)}) - \sum_{j=0}^{\ell-1} \mathbb{P}(A \cap \mathcal{W}_{\ell+t-j,(i-1)(\ell+t)}) \right|$$

$$+ 2\ell \sum_{j=1}^{\ell-1} \mathbb{P}(A \cap T^{-j}(A)). \tag{4.1.12}$$

Let $\Upsilon_i := t\mathbb{P}(A) + \left| \ell\mathbb{P}(A)\mathbb{P}(\mathcal{W}_{0,i(\ell+t)}) - \sum_{j=0}^{\ell-1} \mathbb{P}(A \cap \mathcal{W}_{\ell+t-j,i(\ell+t)}) \right| + 2\ell \sum_{j=1}^{\ell-1} \mathbb{P}(A \cap T^{-j}(A))$. Since $\ell\mathbb{P}(A) < 1$, then it is clear that $|(1 - \ell\mathbb{P}(A))| < 1$. Also, note that $|\mathbb{P}(\mathcal{W}_{0,\ell+t}) - (1 - \ell\mathbb{P}(A))| \leq \Upsilon_0$. Now, we use (4.1.12) recursively to estimate $\left| \mathbb{P}(\mathcal{W}_{0,k(\ell+t)}) - (1 - \ell\mathbb{P}(A))^k \right|$. In fact, we have

$$\left| \mathbb{P}(\mathcal{W}_{0,k(\ell+t)}) - (1 - \ell\mathbb{P}(A))^k \right|$$

$$\leq \sum_{i=0}^{k-1} (1 - \ell\mathbb{P}(A))^{k-1-i} \left| \mathbb{P}(\mathcal{W}_{0,(i+1)(\ell+t)}) - (1 - \ell\mathbb{P}(A))\mathbb{P}(\mathcal{W}_{0,i(\ell+t)}) \right|$$

$$\leq \sum_{i=0}^{k-1} (1 - \ell\mathbb{P}(A))^{k-1-i}\, \Upsilon_i \leq \sum_{i=0}^{k-1} \Upsilon_i \tag{4.1.13}$$

The result follows now at once from (4.1.11) and (4.1.13). $\qquad\square$

We are now in a position to prove Theorem 4.1.4.

*Proof of Theorem* 4.1.4. Letting $A = A_n^{(q)}$, $\ell = \lfloor n/k_n \rfloor$, $k = k_n$ and $t = t_n$ on Proposition 4.1.15, we obtain

$$
\left| \mathbb{P}(\mathcal{W}_{0,n}(A_n^{(q)})) - \left( 1 - \left\lfloor \frac{n}{k_n} \right\rfloor \mathbb{P}(A_n^{(q)}) \right)^{k_n} \right| \leq 2k_n t_n \mathbb{P}(U_n)
$$

$$
+ 2n \sum_{j=1}^{\lfloor n/k_n \rfloor - 1} \mathbb{P}\left( A_n^{(q)} \cap T^{-j} A_n^{(q)} \right)
$$

$$
+ \sum_{i=0}^{k_n - 1} \left| \left\lfloor \frac{n}{k_n} \right\rfloor \mathbb{P}(A_n^{(q)}) \mathbb{P}\left( \mathcal{W}_{0,i(\ell_n + t_n)} \left( A_n^{(q)} \right) \right) \right.
$$

$$
\left. - \sum_{j=0}^{\lfloor n/k_n \rfloor - 1} \mathbb{P}\left( A_n^{(q)} \cap \mathcal{W}_{\ell_n + t_n - j, i(\ell_n + t_n)} \left( A_n^{(q)} \right) \right) \right|. \tag{4.1.14}
$$

Using condition Д$(u_n)$, we have that for the third term:

$$
\sum_{i=0}^{k_n - 1} \left| \left\lfloor \frac{n}{k_n} \right\rfloor \mathbb{P}(A_n^{(q)}) \mathbb{P}\left( \mathcal{W}_{0,i(\ell+t)} \left( A_n^{(q)} \right) \right) - \sum_{j=0}^{\lfloor n/k_n \rfloor - 1} \mathbb{P}\left( A_n^{(q)} \cap \mathcal{W}_{\ell+t-j, i(\ell+t)} \left( A_n^{(q)} \right) \right) \right|
$$

$$
\leq n\gamma(q, n, t_n). \tag{4.1.15}
$$

By (4.1.9), we have that there exists $C$ such that

$$
\left| \mathrm{e}^{-\left\lfloor \frac{n}{k_n} \right\rfloor k_n \mathbb{P}\left( A_n^{(q)} \right)} - \mathrm{e}^{-\theta_n \tau} \right|
$$

$$
\leq \mathrm{e}^{-\theta_n \tau} \left[ \left| \theta_n \tau - \left\lfloor \frac{n}{k_n} \right\rfloor k_n \mathbb{P}\left( A_n^{(q)} \right) \right| + o\left( \left| \theta_n \tau - \left\lfloor \frac{n}{k_n} \right\rfloor k_n \mathbb{P}\left( A_n^{(q)} \right) \right| \right) \right]
$$

$$
\leq C\mathrm{e}^{-\theta_n \tau} \left| \theta_n \tau - n\mathbb{P}\left( A_n^{(q)} \right) \right|
$$

$$
\leq C\mathrm{e}^{-\theta_n \tau} \left| \tau - n\mathbb{P}\left( U_n \right) \right|.
$$

Using (4.1.8) and (4.1.9), there exists $C' > 0$ such that

$$
\left| \left( 1 - \left\lfloor \frac{n}{k_n} \right\rfloor \mathbb{P}\left( A_n^{(q)} \right) \right)^{k_n} - \mathrm{e}^{-\left\lfloor \frac{n}{k_n} \right\rfloor k_n \mathbb{P}\left( A_n^{(q)} \right)} \right| =
$$

$$
\mathrm{e}^{-\left\lfloor \frac{n}{k_n} \right\rfloor k_n \mathbb{P}\left( A_n^{(q)} \right)} \left( \frac{\left( n\mathbb{P}\left( A_n^{(q)} \right) \right)^2}{2k_n} + o\left( \frac{1}{k_n} \right) \right)
$$

$$
\leq \mathrm{e}^{-\theta_n \tau} \left( \frac{\tau^2}{k_n} + C\frac{\tau^2}{k_n} \left| \tau - n\mathbb{P}\left( U_n \right) \right| + o\left( \frac{1}{k_n} \right) \right)
$$

$$
\leq C'\mathrm{e}^{-\theta_n \tau} \frac{\tau^2}{k_n}.
$$

Hence, there exists $C'' > 0$, independent of $\tau$, such that

$$\left| \left( 1 - \left\lfloor \frac{n}{k_n} \right\rfloor \mathbb{P}\left( A_n^{(q)} \right) \right)^{k_n} - e^{-\theta_n \tau} \right| \leq C'' e^{-\theta_n \tau} \left( |\tau - n\mathbb{P}\left( U_n \right)| + \frac{\tau^2}{k_n} \right). \quad (4.1.16)$$

Finally, by Proposition 4.1.12 we have

$$\left| \mathbb{P}(M_n \leq u_n) - \mathbb{P}\left( \mathcal{W}_{0,n}\left( A_n^{(q)} \right) \right) \right| \leq \sum_{j=1}^{q} \mathbb{P}\left( \mathcal{W}_{0,n}\left( A_n^{(q)} \right) \cap T^{-n+j}(U_n \setminus A_n^{(q)}) \right)$$

$$\leq q\mathbb{P}\left( U_n \setminus A_n^{(q)} \right). \quad (4.1.17)$$

Note that when $q = 0$ both sides of inequality (4.1.17) are 0.

The estimate in Theorem 4.1.4 follows from joining the estimates in (4.1.14), (4.1.15), (4.1.16) and (4.1.17). □

*Proof of Corollary* 4.1.5. By (4.1.9), we have that there exists $C > 0$ such that

$$\left| e^{-\theta_n \tau} - e^{-\theta \tau} \right| \leq e^{-\theta \tau} \left[ |\theta_n - \theta|\tau + o(|\theta_n - \theta|) \right]$$

$$\leq C e^{-\theta \tau} |\theta_n - \theta|\tau. \quad (4.1.18)$$

By Theorem 4.1.4 and (4.1.18), there exists $C' > 0$ such that

$$\left| \mathbb{P}(M_n \leq u_n) - e^{-\theta_n \tau} \right| \leq C' \Bigg( k_n t_n \frac{\tau}{n} + n\gamma(q, n, t_n) + n \sum_{j=1}^{\lfloor n/k_n \rfloor - 1} \mathbb{P}\left( A_n^{(q)} \cap T^{-j}\left( A_n^{(q)} \right) \right)$$

$$+ e^{-\theta \tau}(1 + C|\theta_n - \theta|\tau) \left( |\tau - n\mathbb{P}\left( U_n \right)| + \frac{\tau^2}{k_n} \right) + q\mathbb{P}\left( U_n \setminus A_n^{(q)} \right) \Bigg).$$

So, there exists $C'' > 0$ such that

$$\left| \mathbb{P}(M_n \leq u_n) - e^{-\theta_n \tau} \right| \leq C'' \Bigg( k_n t_n \frac{\tau}{n} + n\gamma(q, n, t_n) + n \sum_{j=1}^{\lfloor n/k_n \rfloor - 1} \mathbb{P}\left( A_n^{(q)} \cap T^{-j}\left( A_n^{(q)} \right) \right)$$

$$+ e^{-\theta \tau} \left( |\tau - n\mathbb{P}\left( U_n \right)| + \frac{\tau^2}{k_n} \right) + q\mathbb{P}\left( U_n \setminus A_n^{(q)} \right) \Bigg).$$

$$(4.1.19)$$

The result follows from (4.1.18) and (4.1.19).

□

## 4.2
## Extreme Values for Dynamically Defined Stochastic Processes

We next apply the above ideas and techniques to stochastic processes arising from dynamical systems. So, we start by defining such processes and afterwards discuss their extremal properties already furnished with the results described in the previous section, which were developed with such application in mind.

Take a system $(\mathcal{X}, \mathcal{B}, \mathbb{P}, f)$, where $\mathcal{X}$ is a Riemannian manifold, $\mathcal{B}$ is the Borel $\sigma$-algebra, $f : \mathcal{X} \to \mathcal{X}$ is a measurable map and $\mathbb{P}$ an $f$-invariant probability measure.

In Chapter 6, we will see applications to several types of dynamical systems. However, to give the reader some insight of applications we have in mind, we consider a very simple toy example that we will use throughout this section to illustrate the results and motivate possible applications. The example we will consider is often referred to as doubling map and is one of the simplest examples of a chaotic system.

**Example 4.2.1.** Consider the circle $\mathbb{S}^1$, in its additive representation, *i.e.*, we identify $\mathbb{S}^1 = [0, 1)$ where $[0, 1)$ represents the equivalence classes of the relation $\sim$ defined in $\mathbb{R}$ and given by $x \sim y$ if there exists some $n \in \mathbb{Z}$ such that $x = y + n$. Let $f : [0, 1) \to [0, 1)$, be such that $f(x) = 2x \mod 1$. See Figure 4.1. Please note that the doubling map is a special example of the family of Bernoulli shift maps

$$f(x) = qx \mod 1, \quad q \in \mathbb{N}, \quad q \geq 2 \tag{4.2.1}$$

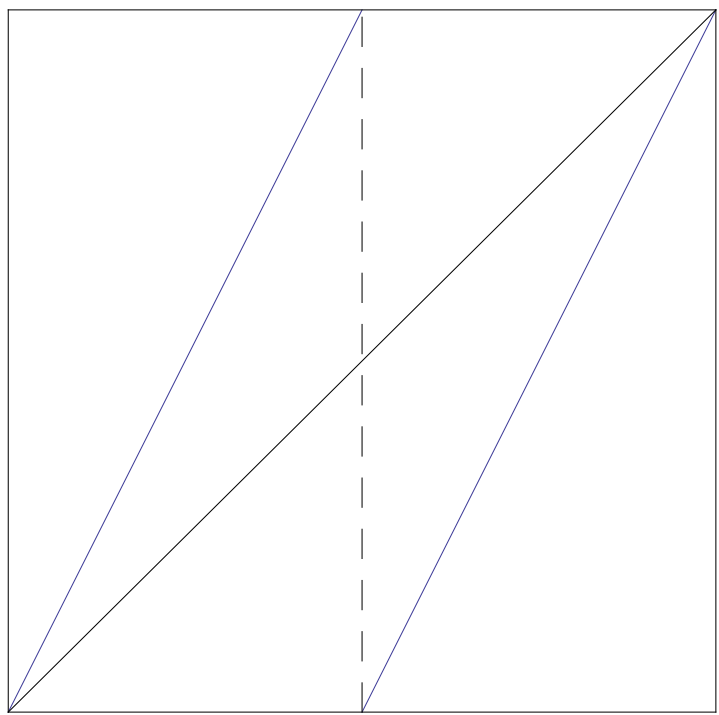

**Figure 4.1** Graphic of the Bernoulli shift (doubling) map depicted in blue and identity map in black

Observe that the map is differentiable in the whole $\mathbb{S}^1$ and for all $x \in [0, 1)$, we have $Df(x) = 2$.

It is easy to see that Lebesgue measure, that we will denote by $\mathrm{Leb}$, is an invariant probability measure. In fact, consider the interval $[a, b)$ with $a, b \in [0, 1)$. Note that $f^{-1}([a, b)) = [a/2, b/2) \cup [a/2 + 1/2, b/2 + 1/2)$ and $\mathrm{Leb}([a, b)) = b - a = \mathrm{Leb}(f^{-1}([a, b)))$.

The orbit of Leb-a.e. point $x$ is dense in $\mathbb{S}^1$. There exist countably many periodic points, which form a dense set of zero Lebesgue measure. In this case, the periodic points are easily identifiable because they correspond to the points with a periodic binary expansion, which means that they are among the rational points of $[0, 1)$.

In particular, we have that 0 is a fixed point ($f(0) = 0$); $1/3$ is a periodic point of prime period 2 ($f(1/3) = 2/3$ and $f(2/3) = 1/3$); and $1/5$ is a periodic point of prime period 4 ($f(1/5) = 2/5$, $f(2/5) = 4/5$, $f(4/5) = 3/5$ and $f(3/5) = 1/5$).

This system is both ergodic and mixing with respect to Leb.

Suppose that the time series $X_0, X_1, \ldots$ arises from such a system simply by evaluating a given observable $\varphi : \mathcal{X} \to \mathbb{R} \cup \{\pm\infty\}$ along the orbits of the system, or in other words, the time evolution given by successive iterations by $f$:

$$X_n = \varphi \circ f^n, \quad \text{for each } n \in \mathbb{N}. \tag{4.2.2}$$

Clearly, $X_0, X_1, \ldots$ defined in this way is not an independent sequence. However, $f$-invariance of $\mathbb{P}$ guarantees that this stochastic process is stationary.

We assume that the r.v. $\varphi : \mathcal{X} \to \mathbb{R} \cup \{\pm\infty\}$ achieves a global maximum at $\zeta \in \mathcal{X}$ (we allow $\varphi(\zeta) = +\infty$).

### 4.2.1
### Observables and Corresponding Extreme Value Laws

We assume that the observable $\varphi : \mathcal{X} \to \mathbb{R} \cup \{+\infty\}$ is of the form

$$\varphi(x) = g\big(\text{dist}(x, \zeta)\big) \tag{4.2.3}$$

where $\zeta$ is a chosen point in the phase space $\mathcal{X}$ and the function $g : [0, +\infty) \to \mathbb{R} \cup \{+\infty\}$ is such that 0 is a global maximum ($g(0)$ may be $+\infty$); $g$ is a strictly decreasing bijection $g : V \to W$ in a neighbourhood $V$ of 0; and has one of the following three types of behaviour:

Type $g_1$: there exists some strictly positive function $h : W \to \mathbb{R}$ such that for all $y \in \mathbb{R}$

$$\lim_{s \to g_1(0)} \frac{g_1^{-1}(s + yh(s))}{g_1^{-1}(s)} = e^{-y}; \tag{4.2.4}$$

Type $g_2$: $g_2(0) = +\infty$ and there exists $\beta > 0$ such that for all $y > 0$

$$\lim_{s \to +\infty} \frac{g_2^{-1}(sy)}{g_2^{-1}(s)} = y^{-\beta}; \tag{4.2.5}$$

Type $g_3$: $g_3(0) = D < +\infty$ and there exists $\gamma > 0$ such that for all $y > 0$

$$\lim_{s \to 0} \frac{g_3^{-1}(D - sy)}{g_3^{-1}(D - s)} = y^{\gamma}. \tag{4.2.6}$$

It may be shown that no non-degenerate limit applies if $\int_0^{g_1(0)} g_1^{-1}(s)ds$ is not finite. Hence, an appropriate choice of $h$ in the Type 1 case is given by

$$h(s) = \int_s^{g_1(0)} g_1^{-1}(t)dt/g_1^{-1}(s) \text{ for } s < g_1(0). \tag{4.2.7}$$

Examples of each one of the three types are as follows:

1) $g_1(x) = -\log x$ (in this case (4.2.4) is easily verified with $h \equiv 1$);
2) $g_2(x) = x^{-1/\alpha}$ for some $\alpha > 0$ (condition (4.2.5) is verified with $\beta = \alpha$);
3) $g_3(x) = D - x^{1/\alpha}$ for some $D \in \mathbb{R}$ and $\alpha > 0$ (condition (4.2.6) is verified with $\gamma = \alpha$).

The fact that the conditions on the shape of $g^{-1}$ imposed by (4.2.4), (4.2.5) and (4.2.6) correspond to the *sufficient* and *necessary* conditions (3.1.5), (3.1.6) and (3.1.7), respectively, on the tail of a distribution which guarantee a non-degenerate EVL in the i.i.d.setting, meaning that the only interesting cases for us are the ones where $g$ is of one of the three types above.

*Remark* 4.2.1. The choice of the observables in (4.2.3) implies that the shape of $g$ determines the type of extremal distribution we get. In particular, for observables of type $g_i$ we get an extremal law of type $\mathrm{e}^{-\tau_i}$, for $i = 1, 2, 3$. While the type of the extremal distribution is essentially determined by the shape of the observable, in the cases when types (2) and (3) apply, *i.e.*, the Fréchet and Weibull families of distributions, the exponents $\beta$ and $\gamma$ appearing in types (2) and (3), in Theorem 3.1.1, respectively, are also influenced by other quantities such as the EI and the local dimension of the stationary invariant measure $\mathbb{P}$. In particular, when such measure is absolutely continuous with respect to Lebesgue and its Radon-Nikodym derivative has a singularity at $\zeta$, then the order of the singularity also influences the value of $\alpha$.

*Remark* 4.2.2. If $\mathbb{P}$ is absolutely continuous with respect to the $d$-dimensional Lebesgue measure and the density of the invariant probability measure is continuous, we could write $\mathbb{P}(B_\eta(\zeta)) \sim \rho(\zeta)\mathrm{Leb}\,(B_\eta(\zeta)) \sim C\rho(\zeta)\eta^d$, where we assume that $\rho(\zeta) = \frac{d\mathbb{P}}{\mathrm{Leb}}(\zeta) > 0$ and $C > 0$ is a positive constant depending on the local metric used.

*Remark* 4.2.3. If $\mathbb{P}$ is supported on an attractor with fractal dimension so that the measure $\mathbb{P}$ has a local dimension at $\zeta$ given by [145]:

$$d(\zeta) := \lim_{\eta \to 0} \frac{\log(\mathbb{P}(B_\eta(\zeta))}{\log \eta} > 0, \tag{4.2.8}$$

then we can write $\log(\mathbb{P}(B_\eta(\zeta))) \sim d(\zeta) \log(\eta)$.We have that $d(\zeta)$ has constant value equal to $d$ if $\mathbb{P}$ follows the prescriptions given in Remark 4.2.2.

Following Remark 4.2.1, we can now relate the parameters $\beta$ and $\gamma$ appearing in types (2) and (3) in Theorem 3.1.1, which correspond to the shape parameter $\xi$ of the GEV ($\xi = 1/\beta$ for type (2) and $\xi = -1/\gamma$ for type (3)), to the local dimension $d(\zeta)$ of the measure.

Let us take $g_2(x) = x^{-1/\alpha}$, for some $\alpha > 0$, and $\varphi = g_2(\mathrm{dist}(x, \zeta))$ then let

$$Q_{s,y} := \frac{1 - F(sy)}{1 - F(s)} = \frac{\mathbb{P}(B_{(sy)^{-\alpha}}(\zeta))}{\mathbb{P}(B_{s^{-\alpha}}(\zeta))}.$$

From (4.2.8) we have that, for all $s$ sufficiently large, $\lim_{y \to \infty} \frac{\log(Q_{s,y})}{\log y} = -\alpha d(\zeta)$. Recalling that for Type (2), by (3.1.6), we should have that $\lim_{s \to \infty} \frac{\log(Q_{s,y})}{\log y} = -\beta$. Hence, we conclude that the only possible value for $\beta$ if the EVL holds is $\beta = \alpha d(\zeta)$, which we will use for numerical simulation purposes. Equivalently, we have:

$$\xi = \frac{1}{\alpha d(\zeta)}. \tag{4.2.9}$$

In the case, we take $g_3(x) = D - x^{1/\alpha}$, for some $\alpha > 0$, and $\varphi = g_3(\text{dist}(x, \zeta))$ then let

$$Q_{s,y} := \frac{1 - F(D - sy)}{1 - F(D - s)} = \frac{\mathbb{P}(B_{(sy)^\alpha}(\zeta))}{\mathbb{P}(B_{s^\alpha}(\zeta))}.$$

From (4.2.8) we have that, for all $s$ sufficiently small, $\lim_{y \to 0} \frac{\log(Q_{s,y})}{\log y} = \alpha d(\zeta)$. Recalling that for Type (3), by (3.1.7), we should have that $\lim_{s \to 0} \frac{\log(Q_{s,y})}{\log y} = \gamma$. This gives the following relation that we will use for numerical simulation purposes: $\gamma = \alpha d$ or, equivalently,

$$\xi = -\frac{1}{\alpha d(\zeta)}. \tag{4.2.10}$$

Note however that Eq. 4.2.8 is only necessary for the convergence of $Q_{s,y}$ in the case of type 2 and type 3 observables, though not sufficient for invariant measures with singular support. When the invariant probability measure $\mathbb{P}$ is not absolutely continuous with respect to the Lebesgue measure, sometimes, in order to simplify, we write, following [136], that

$$\varphi(x) = g\left(\mathbb{P}\big(B_{\text{dist}(x,\zeta)}(\zeta)\big)\right). \tag{4.2.11}$$

In here, since $\mathbb{P}$ is not an absolutely continuous invariant measure, the function $\hbar$ defined for small $\eta \geq 0$ and given by

$$\hbar(\eta) = \mathbb{P}(B_\eta(\zeta)) \tag{4.2.12}$$

is not absolutely continuous. However, we will require in the following that $\hbar$ is continuous in $\eta$. For example, if $\mathcal{X}$ is an interval and $\mathbb{P}$ a Borel probability with no atoms, *i.e.*, points with positive $\mathbb{P}$ measure, then $\hbar$ is continuous.

Throughout the text we assume that the observables and the measure are sufficiently regular so that $\hbar$ is continuous. We will refer to this assumption as condition $(R1)$. In some applications, we will consider that the point $\zeta$ is a repelling periodic point, which implies that condition $(R2)$ defined below holds.

(R1) for $u$ sufficiently close to $u_F := \varphi(\zeta)$, the event

$$U(u) = \{X_0 > u\} = \{x \in \mathcal{X} : \varphi(x) > u\}$$

corresponds to a topological ball centred at $\zeta$. Moreover, the quantity $\mathbb{P}(U(u))$, as a function of $u$, varies continuously on a neighbourhood of $u_F$.

(R2) If $\zeta \in \mathcal{X}$ is a repelling periodic point, of prime period[1] $p \in \mathbb{N}$, then we have that the periodicity of $\zeta$ implies that for all large $u$, $\{X_0 > u\} \cap f^{-p}(\{X_0 > u\}) \neq \emptyset$ and the fact that the prime period is $p$ implies that $\{X_0 > u\} \cap f^{-j}(\{X_0 > u\}) = \emptyset$ for all $j = 1, \ldots, p-1$. Moreover, the fact that $\zeta$ is repelling means that we have backward contraction implying that there exists $0 < \theta < 1$ so that $\bigcap_{j=0}^{i} f^{-jp}(X_0 > u)$ is another ball of smaller radius around $\zeta$ with

$$\mathbb{P}\left(\bigcap_{j=0}^{i} f^{-jp}(X_0 > u)\right) \sim (1 - \theta)^i \mathbb{P}(X_0 > u),$$

for all $u$ sufficiently large. Note that at repelling periodic points, if the measure is absolutely continuous with respect to Lebesgue, with a positive and sufficiently regular density, then the EI can be given by:

$$\theta = 1 - \frac{1}{|\det D(f^p)(\zeta)|}.$$

One of our applications is to equilibrium states, which we explain in Section 4.5. In some of these more general cases, although $g$ is invertible in a small neighbourhood of 0, the function $\hbar$ does not have to be. This means that the observable $\varphi$, as a function of the distance to $\zeta$, may not be invertible in any small neighbourhood of $\zeta$.

For that reason, we now set

$$\ell(\gamma) := \inf\{\eta > 0 : \mathbb{P}(B_\eta(\zeta)) = \gamma\}, \tag{4.2.13}$$

which is well defined for all small enough $\gamma \geq 0$, by the continuity of $\hbar$. Moreover, again by continuity of $\hbar$, we have

$$\mathbb{P}\left(B_{\ell(\gamma)}(\zeta)\right) = \gamma. \tag{4.2.14}$$

*Remark* 4.2.4. Observe that the choice of the observables in (4.2.11) and the assumption on $\mathbb{P}$ regarding the continuity of $\hbar$ guarantee that condition (R1) holds.

*Remark* 4.2.5. Observe that if at time $j \in \mathbb{N}$ we have an exceedance of level $u$ (sufficiently large), i.e., $X_j(x) > u$, then we have an entrance of the orbit of $x$ into the ball $B_{\ell(g^{-1}(u))}(\zeta)$ of radius $\ell(g^{-1}(u))$ around $\zeta$, at time $j$. This means that the behaviour of the tail of $F$, *i.e.*, the behaviour of $1 - F(u)$ as $u \to u_F$ is basically determined by $g^{-1}$ and by the local dimension of the measure at $\zeta$. The above conditions on $g^{-1}$ are just the translation in terms of the shape of $g^{-1}$, of the sufficient and necessary conditions on the tail of $F$ that appear in (3.1.5), (3.1.6) and (3.1.7).

In Sect. 4.2.2 and later in Sects. 6.4-6.8, we will study in detail the mathematical properties of the extremes of distance observables for various mathematical models of increasing level of complexity.

In Chap. 7 we will consider dynamical systems with random components.

---

[1] i.e., the *smallest* $n \in \mathbb{N}$ such that $f^n(\zeta) = \zeta$. Clearly $f^{ip}(\zeta) = \zeta$ for any $i \in \mathbb{N}$.

In Sect. 8.2.1 we will extend the geometrical point of view hinted at above. We will focus on relating the GPD parameters of the distance observables to the geometrical properties of the attractor of high dimensional statistical mechanical systems, assuming that they have Axiom A-like properties. In order to derive results of relevance for *generic* physical systems, we will sacrifice some rigour and use some heuristic arguments. Nonetheless, we will find closely related results.

### 4.2.2
### Extreme Value Laws for Uniformly Expanding Systems

In this section, we illustrate the application of the theory developed in the previous section to stochastic processes arising from specific systems, namely, uniformly expanding and piecewise expanding systems. Although, these are not as general as the non-uniformly expanding systems treated in Chapter 6, for these systems we can actually prove a dichotomy which basically states that either there exists an EI less than 1 at periodic repelling points or there exists an EI equal to 1 at every other point of continuity of the map. An example of such systems is the doubling map introduced in Example 4.2.1.

Up to our knowledge, the statement of this dichotomy appeared first in [49, Section 6], where it is proved for uniformly expanding systems in $S^1$ equipped with the Bernoulli measure and for the cylinder case. Moreover, in the introduction of [49], it is conjectured that this dichotomy should hold in much greater generality (both for more general systems and for the more general case of balls rather than cylinders). In [146], which appeared shortly after [49] on arXiv, the authors build up on the work of [147] and eventually obtain the dichotomy for balls and for conformal repellers. Then, in [141], making use of powerful spectral theory tools developed in [148], the dichotomy for balls is established for general systems such as those for which there exists as spectral gap for their respective Perron-Frobenius operator. In [142], the dichotomy for balls is obtained once again for the same type of systems considered in [141] but using as assumption the existence of decay of correlations against $L^1$ observables (see definition below). In the recent [149], the dichotomy for cylinders is established for mixing countable alphabet shifts, but also in the context of nonconventional ergodic sums.

We also mention that in the papers [150, 151, 152], the authors proved the existence of a limiting law for the cases when $\zeta$ is not a periodic point and we remark that this has been achieved for all such points.

Our basic assumption to prove conditions $\text{Д}_q(u_n), \text{Д}'_q(u_n)$ will be sufficiently fast decay of correlations, in some specific function spaces. Hence we define:

**Definition 4.2.6** (Decay of correlations). Let $\mathcal{C}_1, \mathcal{C}_2$ denote Banach spaces of real valued measurable functions defined on $\mathcal{X}$. We denote the *correlation* of non-zero functions $\phi \in \mathcal{C}_1$ and $\psi \in \mathcal{C}_2$ w.r.t. a measure $\mathbb{P}$ as

$$\text{Cor}_{\mathbb{P}}(\phi, \psi, n) := \frac{1}{\|\phi\|_{\mathcal{C}_1} \|\psi\|_{\mathcal{C}_2}} \left| \int \phi \, (\psi \circ f^n) \, \mathrm{d}\mathbb{P} - \int \phi \, \mathrm{d}\mathbb{P} \int \psi \, \mathrm{d}\mathbb{P} \right|.$$

We say that we have *decay of correlations*, w.r.t. the measure $\mathbb{P}$, for observables in $\mathcal{C}_1$ *against* observables in $\mathcal{C}_2$ if, for every $\phi \in \mathcal{C}_1$ with every $\psi \in \mathcal{C}_2$ we have

$$\mathrm{Cor}_{\mathbb{P}}(\phi, \psi, n) \to 0, \quad \text{as } n \to \infty.$$

We say that we have *decay of correlations against $L^1$ observables* whenever this holds for $\mathcal{C}_2 = L^1(\mathbb{P})$ with $\|\psi\|_{\mathcal{C}_2} = \|\psi\|_1 = \int |\psi| \, \mathrm{d}\mathbb{P}$.

We next state an abstract result giving general conditions to establish the dichotomy:

**Theorem 4.2.7** ([142]). *Consider a continuous dynamical system $(\mathcal{X}, \mathcal{B}, \mathbb{P}, f)$ for which there exists a Banach space $\mathcal{C}$ of real valued functions such that for all $\phi \in \mathcal{C}$ and $\psi \in L^1(\mathbb{P})$,*

$$Cor_{\mathbb{P}}(\phi, \psi, t) \le Ct^{-2}, \tag{4.2.15}$$

*for some $C > 0$. Let $X_0, X_1, \ldots$ be given by (4.2.2), where $\varphi$ achieves a global maximum at some $\zeta \in \mathcal{X}$ and condition (R1) holds. Assume $(u_n)_{n \in \mathbb{N}}$ is such that (2.2.2) holds.*

- *If $\zeta$ is a non periodic point and there exists some $C' > 0$ such that, for all $n \in \mathbb{N}$, we have $\|\mathbf{1}_{U(u_n)}\|_{\mathcal{C}} \le C'$, then conditions $\boxed{\square}_0(u_n)$ and $\boxed{\square}'_0(u_n)$ hold for $X_0, X_1, \ldots$ which means we have an EI equal to 1, i.e., we have an EVL with $H(\tau) = 1 - \mathrm{e}^{-\tau}$.*
- *If $\zeta$ is a periodic point of prime period $p$, if there exists some $C' > 0$ such that, for all $n \in \mathbb{N}$, we have $\|\mathbf{1}_{A_n^{(q)}}\|_{\mathcal{C}} \le C'$, where $q = p$ and the limit $\theta = \lim_{n \to \infty} \theta_n$ exists, where $\theta_n$ is as in (4.1.2), then conditions $\boxed{\square}_q(u_n)$ and $\boxed{\square}'_q(u_n)$ hold for $X_0, X_1, \ldots$ which means we have an EI equal to $\theta < 1$, i.e., we have an EVL with $H(\tau) = 1 - \mathrm{e}^{-\theta\tau}$.*

*Remark* 4.2.8. Observe that decay of correlations as in (4.2.15) against $L^1(\mathbb{P})$ observables is a very strong property. There are no obvious examples - except uniformly expanding maps - with decay of correlations against $L^1$. In fact, regardless of the rate (in this case $n^{-2}$), as long as it is summable, one can actually show that the system has exponential decay of correlations of Hölder observables against $L^\infty(\mathbb{P})$. See [153, Theorem B].

*Remark* 4.2.9. For simplicity, here, we will not deal with discontinuity points of the map $f$. However, that can be done and the existence of an EI less than 1 depends on the existence of some periodic behaviour. Otherwise, we also get an EI equal to 1. See [142, Section 3.3].

*Remark* 4.2.10. We observe that the second statement of Theorem 4.2.7 had already been established in [49]. The first statements of Theorem 4.2.7 and 4.3.5, which finally allowed to establish the dichotomy, were obtained in [142].

In Section 4.5 we will show some specific expanding and piecewise expanding systems for which we can verify the assumptions in Theorems 4.2.7 but to illustrate the essence of their content, we give here a straightforward application to the doubling map.

**Corollary 4.2.11.** *Consider the system* $f : \mathbb{S}^1 \to \mathbb{S}^1$ *given by* $f(x) = 2x \bmod 1$, *as in Example 4.2.1, equipped with the Lebesgue measure. Let* $X_0, X_1, \ldots$ *be given by (4.2.2), where* $\varphi$ *achieves a global maximum at some* $\zeta \in \mathbb{S}^1$ *(take for example* $\varphi(x) = 1 - |x - \zeta|$). *Then, if* $(u_n)_{n \in \mathbb{N}}$ *is such that (2.2.2) holds then*

- *If* $\zeta \in \mathbb{S}^1$ *is not periodic then the EI is equal to 1,* i.e., $\lim_{n \to \infty} \mathrm{Leb}(M_n \le u_n) = \mathrm{e}^{-\tau}$.
- *If* $\zeta \in \mathbb{S}^1$ *is periodic of prime period* $p \in \mathbb{N}$ *then the EI is equal to* $\theta = 1 - (1/2)^p$, i.e., $\lim_{n \to \infty} \mathrm{Leb}(M_n \le u_n) = \mathrm{e}^{-\theta\tau}$.

### 4.2.3
## Example 4.2.1 revisited

In order to illustrate the application of the above results to a specific dynamical system we consider the doubling map introduced in Example 4.2.1. To make things more concrete, we choose the observable

$$\varphi(x) = -\log(|x - 1/3|), \tag{4.2.16}$$

which archives a global maximum $(+\infty)$ at $\zeta = 1/3$ (see Figure 4.2).

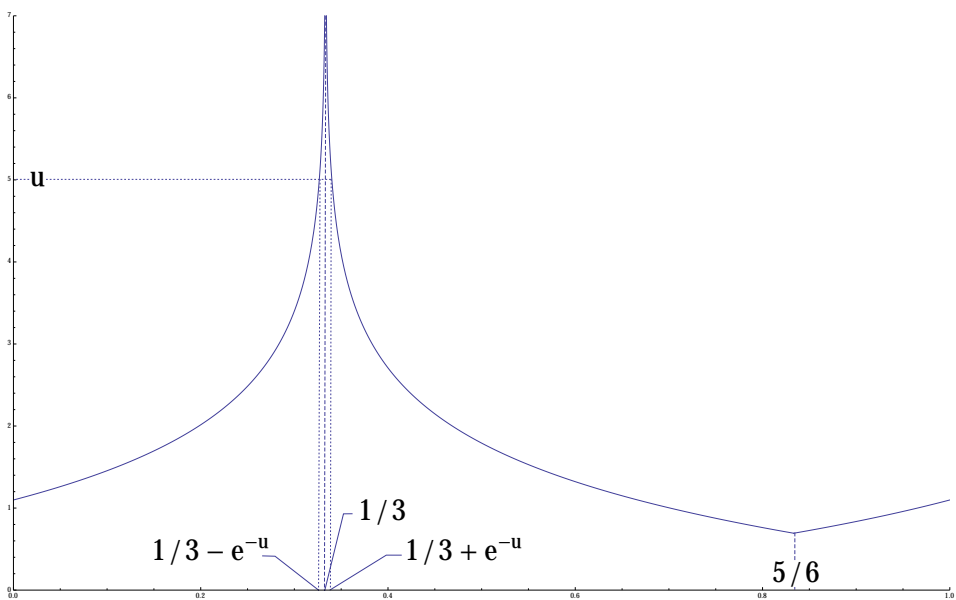

**Figure 4.2** Graphic of the observable $\varphi(x) = -\log(|x - 1/3|)$

We start by computing the distribution function of $X_0$. To do that we need to estimate the Lebesgue measure of the set $\{X_0 \le u\}$. By observing the Figure 4.2, we see that if $u \ge \log 3$ then $\{X_0 \le u\} = [0, 1/3 - \mathrm{e}^{-u}] \cup [1/3 + \mathrm{e}^{-u}, 1]$, hence $F(u) = \mathrm{Leb}(X_0 \le u) = 1 - 2\mathrm{e}^{-u}$. If $\log 2 \le u < \log 3$, then $\{X_0 \le u\} = [1/3 + \mathrm{e}^{-u}, 4/3 - \mathrm{e}^{-u}]$, hence $F(u) = \mathrm{Leb}(X_0 \le u) = 1 - 2\mathrm{e}^{-u}$. If $u < \log 2$ then $\{X_0 \le u\} = \emptyset$ and $F(u) = 0$. To recapitulate, we may write:

$$F(u) = \text{Leb}(X_0 \le u) = \begin{cases} 1 - 2e^{-u} & \text{if } u \ge \log 2 \\ 0 & \text{if } u < \log 2. \end{cases}$$

The normalising sequence is computed as in the i.i.d. case, see Eq. 2.2.2, *i.e.* we solve the equation: $n(1 - F(u_n)) = \tau$ for some $\tau \ge 0$. One easily gets that $u_n = -\log \tau + \log(2n)$. We want to write this in the form $u_n = y/a_n + b_n$, we take $y = -\log \tau \Leftrightarrow \tau = e^{-y}$, $a_n = 1$ and $b_n = \log(2_n)$. However, on the contrary to the i.i.d. case, we have clusters created by the periodicity of the point $1/3$.

The point $1/3$ is a periodic point of period $p = 2$, which is also repelling because $Df^p(1/3) = Df(2/3) \cdot Df(1/3) = 4$. Hence, condition (R2) is easily seen to hold with $\theta = 3/4$. In fact, for $n$ sufficiently large we have

$$U_n = \{X_0 > u_n\} = [1/3 - e^{-u_n}, 1/3 + e^{-u_n}]$$
$$U_n \cap f^{-1}(U_n) = U_n \cap ([1/6 - 1/2e^{-u_n}, 1/6 + 1/2e^{-u_n}]$$
$$\cup [2/3 - 1/2e^{-u_n}, 2/3 + 1/2e^{-u_n}]) = \emptyset$$

and

$$U_n \cap f^{-2}(U_n)$$
$$= U_n \cap \Bigg([1/12 - 1/4e^{-u_n}, 1/12 + 1/4e^{-u_n}] \cup [1/3 - 1/4e^{-u_n}, 1/3 + 1/4e^{-u_n}]$$
$$\cup [7/12 - 1/4e^{-u_n}, 7/12 + 1/4e^{-u_n}] \cup [5/6 - 1/4e^{-u_n}, 5/6 + 1/4e^{-u_n}]\Bigg)$$
$$= [1/3 - 1/4e^{-u_n}, 1/3 + 1/4e^{-u_n}].$$

The latter implies that $D'_0(u_n)$ cannot hold because $\lim_{n\to\infty} n\text{Leb}(X_0 > u_n, X_2 > u_n) = \tau/4$, which is larger than 0 for every positive $\tau$. Hence, we take $q = p = 2$ and set $A_n^{(q)} = U_n \setminus f^{-2}(U_n) = [1/3 - e^{-u_n}, 1/3 - 1/4e^{-u_n}] \cup [1/3 + 1/4e^{-u_n}, 1/3 + e^{-u_n}]$. Using Theorem 4.2.7, we can actually check that conditions $D_q(u_n)$ and $D'_q(u_n)$ hold and in that case we have that $M_n$ has an asymptotic distribution similar to the one obtained in the i.i.d. case, which is now affected by an Extremal Index equal to

$$\theta = \lim_{n\to\infty} \frac{\text{Leb}([1/3 - e^{-u_n}, 1/3 - 1/4e^{-u_n}] \cup [1/3 + 1/4e^{-u_n}, 1/3 + e^{-u_n}])}{\text{Leb}[1/3 - e^{-u_n}, 1/3 + e^{-u_n}]} = \frac{3}{4}.$$

Hence, we can write that for the stochastic process $X_0, X_1, \ldots$ defined in (4.2.2), with $\varphi$ given by (4.2.16):

$$\lim_{n\to\infty} \text{Leb}(M_n - \log(2n) \le y) = e^{-\frac{3}{4}e^{-y}}.$$

### 4.2.4
### Proof of the Dichotomy for Uniformly Expanding Maps

In this section we prove Theorems 4.2.7. We begin with the following notion:

**Definition 4.2.12.** For every $A \in \mathcal{B}$, we define the *first return time* to $A$, which we denote by $R(A)$, as the minimum of the return time function to $A$, *i.e.*,

$$R(A) = \min_{x \in A} r_A(x),$$

where $r_A(x)$ has been defined in in Eq. 2.3.1.

**Proposition 4.2.13.** *Consider a dynamical system $(\mathcal{X}, \mathcal{B}, \mathbb{P}, f)$ satisfying the assumptions of Theorems 4.2.7, in particular, for which there exists decay of correlations against $L^1$. For any point $\zeta$, assume that (R1) holds. Consider $X_0, X_1, \ldots$ defined as in (4.2.2), let $u_n$ be such that (2.2.2) holds. Then condition $\underline{D}'_q(u_n)$ holds for $X_0, X_1, \ldots$, with $q = 0$ if $\zeta$ is not periodic and with $q = p$ if $\zeta$ is periodic of prime period $p$.*

*Proof.* For every $n \in \mathbb{N}$, if $\zeta$ is not periodic take $q = 0$ if $\zeta$ is periodic of period $p$, let $q = p$. Also, set $R_n := R(A_n^{(q)})$.

By hypothesis, for all $n \in \mathbb{N}$ we have that $\mathbf{1}_{A_n^{(q)}} \in \mathcal{C}$ and $\|\mathbf{1}_{A_n^{(q)}}\|_{\mathcal{C}} \le C'$, for some $C' > 0$.

Taking $\phi = \mathbf{1}_{A_n^{(q)}}$ and $\psi = \mathbf{1}_{A_n^{(q)}}$ in (4.2.15) we get

$$\mathbb{P}\left(A_n^{(q)} \cap f^{-j}(A_n^{(q)})\right) \le (\mathbb{P}(A_n^{(q)}))^2 + C \left\|\mathbf{1}_{A_n^{(q)}}\right\|_{\mathcal{C}} \left\|\mathbf{1}_{A_n^{(q)}}\right\|_{L^1(\mathbb{P})} j^{-2}$$
$$\le (\mathbb{P}(A_n^{(q)}))^2 + C^* \mathbb{P}(A_n^{(q)}) j^{-2}, \qquad (4.2.17)$$

where $C^* = CC' > 0$. By considering the definition of $R_n$ and the estimate (4.2.17), and by recalling that $n\mathbb{P}(U_n) \to \tau$ as $n \to \infty$, it follows that there exists some constant $c > 0$ such that

$$n \sum_{j=1}^{\lfloor n/k_n \rfloor} \mathbb{P}(A_n^{(q)} \cap f^{-j}(A_n^{(q)})) = n \sum_{j=R_n}^{\lfloor n/k_n \rfloor} \mathbb{P}(A_n^{(q)} \cap f^{-j}(A_n^{(q)}))$$

$$\le n \lfloor \tfrac{n}{k_n} \rfloor \mathbb{P}(A_n^{(q)})^2 + n \, C^* \mathbb{P}(A_n^{(q)}) \sum_{j=R_n}^{\lfloor n/k_n \rfloor} j^{-2}$$

$$\le \frac{(n\mathbb{P}(A_n^{(q)}))^2}{k_n} + n \, C^* \mathbb{P}(A_n^{(q)}) \sum_{j=R_n}^{\infty} j^{-2} \le c \left( \frac{\tau^2}{k_n} + \tau \sum_{j=R_n}^{\infty} j^{-2} \right).$$

Since by hypothesis $k_n \to \infty$, as $n \to \infty$, the result is proved once we show that $R_n \to \infty$ as $n \to \infty$. We separate the proof of the latter in two cases.

Let $\zeta$ be a non periodic point, in which case, we have that $q = 0$ and $A_n^{(q)} = U_n$. Note that, by hypothesis, $f^i$ is continuous at $\zeta$, for all $i \in \mathbb{N}$. For some fixed $j$, we define $\epsilon = \min_{i=1,\ldots,j} \operatorname{dist}(f^i(\zeta), \zeta)$ Then, using the continuity of each $f^i$ at $\zeta$, for every $i = 1, \ldots, j$, there exists $\delta_i > 0$ such that $f^i(B_{\delta_i}(\zeta)) \subset B_{\epsilon/2}(f^i(\zeta))$. Now, let $U := \cap_{i=1}^j B_{\delta_i}(\zeta)$.

If we choose $N$ sufficiently large so that $U_n \subset U$ for all $n \geq N$ then, by definition of $\epsilon$, it is clear that $f^i(U_n) \cap U_n = \emptyset$ for all $i = 1, \ldots, j$, which implies that $R(U(u_n)) > j$. Since $j$ is arbitrary, the statement follows.

Regarding the case when $\zeta$ is periodic point of prime period $p$, by the Hartman-Grobman theorem there is a neighbourhood $V$ around $\zeta$ where $f^p$ is conjugate to its linear approximation given by the derivative at $\zeta$. Hence, for $n$ sufficiently large so that $U_n \subset V$, if a point starts in $A_n^{(q)}$ it takes a time $\beta_n$ to leave $V$, during which, it is guaranteed that it does not return to $U_n$. Moreover, since by condition (R1) and definition of $u_n$, we have that $U_n$ shrinks to $\zeta$ as $n \to \infty$, then $\beta_n \to \infty$ as $n \to \infty$. Since $R_n \geq \beta_n$ then the statement follows. □

We are left to prove condition $\mathring{\mathrm{A}}_q(u_n)$. As we have mentioned before, the main advantage of the theory developed in Section 4.1 is that the new condition $\mathring{\mathrm{A}}_q(u_n)$ is designed to follow easily by decay of correlations, in contrast to the original $D(u_n)$ condition of Leadbetter. In fact, we don't even need this strong type of decay of correlations against $L^1$ (see Section 4.4).

*Proof of Theorems* 4.2.7. Choose $\phi = \mathbf{1}_{A_n^{(q)}}$ and $\psi = \mathbf{1}_{\mathcal{W}_{0,\ell}(A_n^{(q)})}$. We have that we can take $\gamma(n, t) = CC't^{-2}$. Hence, condition $\mathring{\mathrm{A}}_q(u_n)$ is trivially satisfied for the sequence $(t_n)_n$ given by, *e.g.* $t_n = n^\alpha$ for $1/2 < \alpha < 1$. □

## 4.3
## Point Processes of Rare Events

If we enrich the process by considering multiple exceedances, we are led to studying point processes of rare events as a result of counting the number of exceedances within a certain time frame. For uniformly expanding systems, under the exact same assumptions just seen above, the REPP converges in distribution to a standard Poisson process, when no clustering is involved and to a compound Poisson process with intensity $\theta$ and a geometric multiplicity d.f., otherwise.

### 4.3.1
### Absence of Clustering

When condition $\mathring{\mathrm{A}}'_0(u_n)$ holds, there are no clusters and so we may benefit from a criterion, proposed by Kallenberg [127, Theorem 4.7], which applies only to simple point processes without multiple events. Accordingly, we merely need to adjust condition $\mathring{\mathrm{A}}_0(u_n)$ to this scenario of multiple exceedances in order to prove that the REPP converges in distribution to a standard Poisson process. We denote this adapted condition by:

**Condition** $(D_3(u_n))$**.** Let $A \in \mathcal{R}$ and $t \in \mathbb{N}$. We say that $D_3(u_n)$ holds for the sequence $X_0, X_1, \ldots$ if

$$|\mathbb{P}(\{X_0 > u_n\} \cap \{\mathcal{N}_{u_n}(A + t) = 0\}) - \mathbb{P}(\{X_0 > u_n\})\mathbb{P}(\mathcal{N}_{u_n}(A) = 0)| \leq \gamma(n, t),$$

where $\gamma(n, t)$ is nonincreasing in $t$ for each $n$ and $n\gamma(n, t_n) \to 0$ as $n \to \infty$ for some sequence $t_n = o(n)$. (The last equality means that $t_n / n \to 0$ as $n \to \infty$).

In [74, Theorem 5] it is proved a strengthening of [73, Theorem 1] which says the following:

**Theorem 4.3.1** ([74, Theorem 5])**.** *Let $X_1, X_2, \ldots$ be a stationary stochastic process for which conditions $D_3(u_n)$ and $\underline{D}'_0(u_n)$ hold for a sequence of levels $u_n$ such that (2.2.2) holds. Then the REPP $N_n$ defined in (3.3.1) is such that $N_n \xrightarrow{d} N$, as $n \to \infty$, where $N$ denotes a Poisson Process with intensity 1.*

### 4.3.2
### Presence of Clustering

Condition $\underline{D}'_0(u_n)$ prevents the existence of clusters of exceedances, which implies that the EVL is standard exponential $\bar{H}(\tau) = \mathrm{e}^{-\tau}$. When $\underline{D}'_0(u_n)$ fails, the clustering of exceedances is responsible for the appearance of a parameter $0 < \theta < 1$ in the EVL, called EI, and implies that, in this case, $\bar{H}(\tau) = \mathrm{e}^{-\theta\tau}$. In [49], the authors established a connection between the existence of an EI less than 1 and a periodic behavior. This has been later generalized for REPP in [140].

For the convergence of the REPP when there is clustering, one cannot use the aforementioned criterion of Kallenberg because the point processes are not simple anymore and possess multiple events. This means that a much deeper analysis must be done in order to obtain convergence of the REPP. This was carried out in [140] and we describe below the main results and conditions needed.

Let $\zeta$ be a periodic point of prime period $q$. Firstly, we consider the sequence $\left(U^{(\kappa)}(u_n)\right)_{\kappa \geq 0}$ of nested balls centred at $\zeta$ given by

$$U^{(0)}(u_n) = U(u_n) = U_n \quad \text{and} \quad U^{(\kappa)}(u_n) = T^{-q}(U^{(\kappa-1)}(u_n)) \cap U(u_n), \quad \forall \kappa \in \mathbb{N}. \text{(4.3.1)}$$

Then, for $i, \kappa, \ell, s \in \mathbb{N} \cup \{0\}$, we define the following events:

$$Q_{q,i}^{\kappa}(u_n) := T^{-i}\left(U^{(\kappa)}(u_n) - U^{(\kappa+1)}(u_n)\right). \tag{4.3.2}$$

These events correspond to $\kappa+1$ exceedances that appear separated by $q$ units of time and which are followed by a period of length larger than $q$, during which no exceedances occur. Note that $Q_{q,0}^0(u_n) = A_n^{(q)}$. Besides, $U_n = \bigcup_{\kappa=0}^{\infty} Q_{q,0}^\kappa(u_n) \bigcup\{\zeta\}$ which means that the ball centred at $\zeta$ which corresponds to $U_n$ can be decomposed into a sequence of disjoint annuli where $Q_{q,0}^0(u_n)$ are the most outward annuli and the inner annuli $Q_{q,0}^{\kappa+1}(u_n)$ are sent outward by $T^p$ onto the annuli $Q_{q,0}^\kappa(u_n)$, i.e., $T^q(Q_{q,0}^{\kappa+1}(u_n)) = Q_{q,0}^\kappa(u_n)$.

*Remark* 4.3.2. When $\varphi$ achieves a global maximum at a repelling periodic point $\zeta$, we are lead to the appearance of clusters of exceedances whose size depends on the severity of the first exceedance that begins the cluster. To be more precise, let $x \in \mathcal{X}$: if we have a first exceedance at time $i \in \mathbb{N}$, which means that $f^i(x)$ enters the ball $U(u)$, then we must have that $f^i(x) \in Q_{p,0}^\kappa(u)$ for some $\kappa \geq 0$, which we express by saying that the entrance at time $i$ had a *depth* $\kappa$. Notice that the deeper the entrance, the closer $f^i(x)$ got to $\zeta$ and the more severe is the exceedance. Now, observe that if $f^i(x) \in Q_{p,0}^\kappa(u)$ we must have $f^{i+p}(x) \in Q_{p,0}^{\kappa-1}(u), \ldots, f^{i+\kappa p}(x) \in Q_{p,0}^0(u)$ and $f^{i+(\kappa+1)p}(x) \notin U(u)$ which means that the size of the cluster initiated at time $i$ is exactly $\kappa + 1$ and ends with a visit to the outermost ring $Q_{p,0}^0(u)$, which plays the role of an escaping exit from $U(u)$. So the depth of the entrance in $U(u)$ determines the size of the cluster, and the deeper the entrance, the more severe is the corresponding exceedance and the longer the cluster.

We are now ready to state the adapted condition:

**Condition** $(D_q(u_n)^*)$. We say that $D_q(u_n)^*$ holds for the sequence $X_0, X_1, X_2, \ldots$ if for any integers $t, \kappa_1, \ldots, \kappa_\zeta, n$ and any $J = \cup_{i=2}^q I_j \in \mathcal{R}$ with $\inf\{x : x \in J\} \geq t$,

$$\left| \mathbb{P}\left( Q_{q,0}^{\kappa_1}(u_n) \cap \left( \cap_{j=2}^q \mathcal{N}_{u_n}(I_j) = \kappa_j \right) \right) - \mathbb{P}\left( Q_{q,0}^{\kappa_1}(u_n) \right) \mathbb{P}\left( \cap_{j=2}^q \mathcal{N}_{u_n}(I_j) = \kappa_j \right) \right|$$
$$\leq \gamma(n,t),$$

where for each $n$ we have that $\gamma(n,t)$ is nonincreasing in $t$ and $n\gamma(n,t_n) \to 0$ as $n \to \infty$, for some sequence $t_n = o(n)$.

In [140], for technical reasons only, the authors introduced a slight strengthening to $D_q'(u_n)$. This condition was denoted by $D_q'(u_n)^*$ and it requires that

$$\lim_{n \to \infty} n \sum_{j=1}^{[n/k_n]} \mathbb{P}(Q_{p,0}^{(0)}(u_n) \cap \{X_j > u_n\}) = 0, \tag{4.3.3}$$

which holds whenever condition $D_q'(u_n)$ does.

From the study developed in [140] and as noticed in [142, Appendix B], we can state the following result which applies to general stationary stochastic processes, whose full proof is appearing in the forthcoming paper [154].

**Theorem 4.3.3.** *Let $X_0, X_1, \ldots$ satisfy conditions $D_q(u_n)^*$, $D_q'(u_n)^*$ and*

$$\lim_{n \to \infty} \sum_{\kappa \geq 1} \mathbb{P}(U^{(\kappa)}(u_n)) = 0,$$

*where $(u_n)_{n \in \mathbb{N}}$ is such that (2.2.2) holds. Assume that the limit $\theta = \lim_{n \to \infty} \theta_n$ exists, where $\theta_n$ is as in (4.1.2) and moreover that for each $\kappa \in \mathbb{N}$, the following limit also exists*

$$\pi(\kappa) = \lim_{n \to \infty} \frac{\left( \mathbb{P}(Q_{q,0}^{\kappa-1}(u_n)) - \mathbb{P}(Q_{q,0}^{\kappa}(u_n)) \right)}{\mathbb{P}(Q_{q,0}^0(u_n))}. \tag{4.3.4}$$

*Then the REPP $N_n$ converges in distribution to a compound Poisson process with intensity $\theta$ and multiplicity d.f. $\pi$ given by (4.3.4).*

When $X_0, X_1, \ldots$ arise from a dynamical system and the point $\zeta$ is a repelling periodic point for which condition (R2) holds, then it follows immediately from (R2) that for each $\kappa$ the limit in (4.3.4) exists and equals $\pi(\kappa) = \theta(1-\theta)^{\kappa-1}$. Hence, as corollary we can state now the main result of [140].

**Corollary 4.3.4** ( [140, Theorem 1]). *Let $X_0, X_1, \ldots$ be given by (4.2.2), where $\varphi$ achieves a global maximum at the repelling periodic point $\zeta$, of prime period $p$, and conditions (R1) and (R2) hold. Let $(u_n)_{n \in \mathbb{N}}$ be a sequence satisfying (2.2.2). Assume that conditions $D_p(u_n)^*$, $D'_p(u_n)^*$ hold. Then the REPP $N_n$ converges in distribution to a compound Poisson process $N$ with intensity $\theta$ and multiplicity d.f. $\pi$ given by $\pi(\kappa) = \theta(1-\theta)^{\kappa-1}$, for every $\kappa \in \mathbb{N}_0$, where the extremal index $\theta$ is given by the expansion rate at $\zeta$ stated in (R2).*

### 4.3.3
### Dichotomy for Uniformly Expanding Systems for Point Processes

Making minor adjustments to the proof of Theorem 4.2.7, we can check the conditions of Theorem 4.3.1 and Corollary 4.3.4 to obtain the following stronger version of the dichotomy.

**Theorem 4.3.5** ([142]). *Under the same assumptions of Theorem 4.2.7,*

- *If $\zeta$ is a non periodic point and there exists some $C' > 0$ such that, for all $n \in \mathbb{N}$, we have $\|\mathbf{1}_{U(u_n)}\|_{\mathcal{C}} \le C'$, then conditions $D_3(u_n)$ and $\mathcal{A}'_0(u_n)$ hold for $X_0, X_1, \ldots$, which means that the REPP $N_n$ defined in (3.3.1) converges in distribution to a standard Poisson process as $n \to \infty$.*

- *If $\zeta$ is a periodic point of prime period $p$, at which condition (R2) holds and there exists some $C' > 0$ such that, for all $n \in \mathbb{N}$, we have $\|\mathbf{1}_{Q_p(u_n)}\|_{\mathcal{C}} \le C'$, then conditions $D^p(u_n)^*$ and $D'_p(u_n)^*$ hold for $X_0, X_1, \ldots$ which means that the REPP $N_n$ converges in distribution to a compound Poisson process $N$ with intensity $\theta$ and multiplicity d.f. $\pi$ given by $\pi(\kappa) = \theta(1-\theta)^{\kappa}$, for every $\kappa \in \mathbb{N}_0$, where the extremal index $\theta$ is given by the expansion rate at $\zeta$ stated in (R2).*

*Remark* 4.3.6. We remark that the second statement of the previous theorem had already been obtained in [140] but the final version of the dichotomy, which also includes the statement regarding non periodic points was only established in [142].

The proof of condition $D_3(u_n)$ follows after minor adjustments to the proof of $\mathcal{A}_q(u_n)$ in Section 4.2.4. Since condition $\mathcal{A}'_0(u_n)$ holds at every non-periodic point $\zeta$ (see Proposition 4.2.13), then for all such points $\zeta$ the corresponding REPP $N_n$ converges in distribution to a standard Poisson process.

The proof of condition $D_q(u_n)^*$ follows with some adjustments from the proof of $\mathcal{A}_q(u_n)$ in Section 4.2.4. Since condition $\mathcal{A}'_q(u_n)$ holds at every periodic point $\zeta$, as was shown in Proposition 4.2.13, then condition $D'_q(U_n)^*$ from [140] also holds for all such points $\zeta$. So, at every periodic point $\zeta$ of prime period $p$, the respective

REPP $N_n$ converges in distribution to a compound Poisson process with intensity $\theta$ and multiplicity d.f. $\pi$ given by $\pi(\kappa) = \theta(1-\theta)^\kappa$, for every $\kappa \in \mathbb{N}_0$, where the extremal index $\theta$ is given by the expansion rate at $\zeta$ stated in (R2).

*Remark* 4.3.7. We note that the pattern of appearance of clusters resulting from a periodic repelling point is quite rigid and obeys a very strict pattern. Let $\zeta$ be a repelling periodic point of period $p \in \mathbb{N}$, then a clustering is easily identified by a bulk of strictly decreasing exceedances which appear separated precisely by $p-1$ non-exceedances. This observation is quite helpful in order to identify clusters and to perform a data declustering in order to apply the POT method, for example, to make statistical inference on the tail of the distributions.

## 4.4
## Conditions $\Delta_q(u_n)$, $D_3(u_n)$, $D^p(u_n)^*$ and Decay of Correlations

In general terms the conditions $\Delta_q(u_n)$, $D_3(u_n)$, and $D^p(u_n)^*$ follow from sufficiently fast (e.g. polynomial) decay of correlations of the dynamical system. This is why these conditions can be considered being much more useful than Leadbetter's condition $D(u_n)$. Indeed $D(u_n)$ usually follows only from strong uniform mixing, like $\alpha$-mixing (see [155] for a definition), and even then applies only on certain subsequences,. This means that that EVLs can only be shown for cylinders. Conditions $\Delta_q(u_n)$, $D_3(u_n)$, and $D^p(u_n)^*$ follow from decay of correlations, which is much weaker and allows instead for deriving EVLs s for balls.

In fact, in order to prove $\Delta_q(u_n)$, $D_3(u_n)$, and $D^p(u_n)^*$, there is actually no need for such strong type of decay of correlations as such as against $L^1$, like in the assumptions of Theorems 4.2.7 or 4.3.5. It suffices to have decay of correlations against all $\psi$ in, *e.g.*, $L^\infty$.

Rates of decay of correlations are nowadays well known for many chaotic systems. Examples of these include hyperbolic or uniformly expanding systems as well as the non-hyperbolic or non-uniformly expanding admitting, *e.g.*, inducing schemes with a well behaved return time function. See Sect. 5.5. for some related discussion. In fact, in two remarkable papers L.-S..Young[156, 157] showed that the rates of decay of correlations of the original system are intimately connected with the recurrence rates of the respective induced map.

Just to give an idea of how simple it is to check $\Delta_q(u_n)$, $D_3(u_n)$, and $D^p(u_n)^*$ for systems with sufficiently fast decay of correlations, let us begin by defining the following Banach spaces:

Given a function $\psi : Y \to \mathbb{R}$ on an interval $Y$, the *variation* of $\psi$ is defined as

$$\text{Var}(\psi) := \sup \left\{ \sum_{i=0}^{n-1} |\psi(x_{i+1}) - \psi(x_i)| \right\},$$

where the supremum is taken over all finite ordered sequences $(x_i)_{i=0}^n \subset Y$.

We use the norm $\|\psi\|_{BV} = \sup |\psi| + \text{Var}(\psi)$, which makes $BV := \{\psi : Y \to \mathbb{R} : \|\psi\|_{BV} < \infty\}$ into a Banach space.

Now, let $\mathcal{X}$ be a compact subset of $\mathbb{R}^d$ and let $\psi : \mathcal{X} \to \mathbb{R}$. Given a Borel set $\Gamma \subset \mathcal{X}$, we define the oscillation of $\psi \in L^1(\text{Leb})$ over $\Gamma$ as

$$\text{osc}(\psi, \Gamma) := \underset{\Gamma}{\text{ess sup}} \, \psi - \underset{\Gamma}{\text{ess inf}} \, \psi.$$

It is easy to verify that $x \mapsto \text{osc}(\psi, B_\varepsilon(x))$ defines a measurable function (see [158, Proposition 3.1]). Given real numbers $0 < \alpha \leq 1$ and $\varepsilon_0 > 0$, we define $\alpha$-seminorm of $\psi$ as

$$|\psi|_\alpha = \sup_{0 < \varepsilon \leq \varepsilon_0} \varepsilon^{-\alpha} \int_{\mathbb{R}^d} \text{osc}(\psi, B_\varepsilon(x)) \, \text{dLeb}(x).$$

Let us consider the space of functions with bounded $\alpha$-seminorm $V_\alpha = \{\psi \in L^1(\text{Leb}) : |\psi|_\alpha < \infty\}$, and endow $V_\alpha$ with the norm $\|\cdot\|_\alpha = \|\cdot\|_{L^1(\text{Leb})} + |\cdot|_\alpha$ which makes it into a Banach space. We note that $V_\alpha$ is independent of the choice of $\varepsilon_0$.

In what follows we assume that, for every $\zeta \in \mathcal{X}$ condition (R1) holds, and, in case $\zeta$ is periodic of period $p$, then condition (R2) also holds.

**Proposition 4.4.1.** *Assume that for our system $f : \mathcal{X} \to \mathcal{X}$ we have decay of correlations for all $\phi \in \mathcal{C}$, where $\mathcal{C}$ is BV or $V_\alpha$, depending on whether $\mathcal{X}$ is a compact subset of $\mathbb{R}$ or $\mathbb{R}^d$ (with $d = 2, 3, \ldots$), and all $\psi \in L^\infty$ so that there exist $C$, independent of $\phi, \psi$, and a rate function $\varrho : \mathbb{N} \to \mathbb{R}$ such that*

$$\left| \int \phi \cdot (\psi \circ f^t) d\mathbb{P} - \int \phi d\mathbb{P} \int \psi d\mathbb{P} \right| \leq C \|\phi\|_\mathcal{C} \|\psi\|_\infty \varrho(t), \quad \forall t \in \mathbb{N}_0, \ (4.4.1)$$

*and $n\varrho(t_n) \to 0$, as $n \to \infty$ for some $t_n = o(n)$. Then conditions $\not\!\!\Delta_q(u_n)$, $D_3(u_n)$, and $D^p(u_n)^*$ hold.*

*Proof.* In what follows, for all $n \in \mathbb{N}$, let $A_n = U(u_n)$, in case $\zeta$ is non periodic, and $A_n = Q_p(u_n)$ if $\zeta$ is periodic of prime period $p$. Take $\phi = \mathbf{1}_{A_n}, \psi = \mathbf{1}_{\mathcal{W}_{0,\ell}(A_n)}$. Observe that whether $\|\cdot\|_\mathcal{C}$ is $\|\cdot\|_{BV}$ or $\|\cdot\|_\alpha$, we have that there exists some $C' > 0$ such that $\|\mathbf{1}_{A_n}\|_\mathcal{C} \leq C'$, for all $n \in \mathbb{N}$. Set $c = CC'$. Then (4.4.1) implies that, in case $\zeta$ is non periodic, condition $D_2(u_n)$ holds and, in case $\zeta$ is periodic of prime period $p$, condition $D^p(u_n)$ holds, where $\gamma(n, t) = \gamma(t) := c\varrho(t)$ and the sequence $t_n$ is such that $n\varrho(t_n) \to 0$, as $n \to \infty$. Note that to prove $D_3(u_n)$, we just have to take $\psi = \mathbf{1}_{\mathcal{N}(A)=0}$ and for $D^p(u_n)^*$, we would take $\phi = \mathbf{1}_{Q_{p,0}^{\kappa_1}(u_n)}$, $\psi = \mathbf{1}_{\left(\cap_{j=2}^\varsigma \mathcal{N}_{u_n}(I_j) = \kappa_j\right)}$ and the argument would follow as before. $\qquad \square$

Note that, in the proof above, it has been useful for checking $\not\!\!\Delta_q(u_n)$ that $\mathbf{1}_{A_n} \in \mathcal{C}$ and $\|\mathbf{1}_{A_n}\|_\mathcal{C} \leq C'$. However, these conditions can still be checked even when $\mathbf{1}_{A_n} \notin \mathcal{C}$. This is the case when $\mathcal{C}$ is the Banach space of Hölder observables, which is used, for example, to obtain decay of correlations for systems with Young towers. The idea is to use, as in [72, Lemma 3.3], an adequate Hölder approximation for $\mathbf{1}_{A_n}$.

Now, let $\mathcal{X}$ be a compact subset of $\mathbb{R}^d$ and let $\phi : \mathcal{X} \to \mathbb{R}$. Let $\mathcal{H}_\beta$ denote the space of Hölder continuous functions $\phi$ with exponent $\beta$ equipped with the norm

$\|\phi\|_{\mathcal{H}_\beta} = \|\phi\|_\infty + |\phi|_{\mathcal{H}_\beta}$, where

$$|\phi|_{\mathcal{H}_\beta} = \sup_{x \neq y} \frac{|\phi(x) - \phi(y)|}{|x - y|^\beta}.$$

**Proposition 4.4.2.** *Assume that $\mathcal{X}$ is a compact subset of $\mathbb{R}^d$ and $f : \mathcal{X} \to \mathcal{X}$ is a system with an absolutely continuous invariant measure $\mathbb{P}$, such that $\frac{d\mathbb{P}}{d\text{Leb}} \in L^{1+\epsilon}$. Assume, moreover, that the system has decay of correlations for all $\phi \in \mathcal{H}_\beta$ against any $\psi \in L^\infty$ so that there exists some $C > 0$ independent of $\phi, \psi$ and $t$, and a rate function $\varrho : \mathbb{N} \to \mathbb{R}$ such that*

$$\left| \int \phi \cdot (\psi \circ f^t) d\mathbb{P} - \int \phi d\mathbb{P} \int \psi d\mathbb{P} \right| \leq C \|\phi\|_{\mathcal{H}_\beta} \|\psi\|_\infty \varrho(t), \tag{4.4.2}$$

*and $n^{1+\beta(1+\max\{(0,\epsilon+1)/\epsilon-d\}+\delta)} \varrho(t_n) \to 0$, as $n \to \infty$ for some $\delta > 0$ and $t_n = o(n)$. Then conditions $D_2(u_n)$, $D_3(u_n)$, $D^p(u_n)$ and $D^p(u_n)^*$ hold.*

*Proof.* Since $\frac{d\mathbb{P}}{\text{Leb}} \in L^{1+\epsilon}$, by Hölder's inequality, this last fact implies that for all Borel sets $B$, there exists $C > 0$ such that

$$\mathbb{P}(B) \leq C(\text{Leb}(B))^\Theta, \tag{4.4.3}$$

where $\Theta = \epsilon/(1 + \epsilon)$.

For $\eta = \max\{\Theta^{-1} - d, 0\} + \delta > 0$ we build the Hölder approximation $\phi$ of $\mathbf{1}_{A_n}$, where $A_n$ is as in the proof of Proposition 4.4.1. Let $D_n := \{x \in A_n : \text{dist}(x, \overline{A_n^c}) \geq (\mathbb{P}(A_n))^{1+\eta}\}$, where $\bar{A}$ denotes the closure of $A$. Define $\phi : \mathcal{X} \to \mathbb{R}$ as

$$\phi(x) = \begin{cases} 0 & \text{if } x \in \overline{A_n^c} \\ \frac{\text{dist}(x, \overline{A_n^c})}{\text{dist}(x, \overline{A_n^c}) + \text{dist}(x, D_n)} & \text{if } x \in \overline{A_n^c}^c \setminus D_n \\ 1 & \text{if } x \in D_n \end{cases}.$$

Observe that $\phi$ is Hölder continuous with Hölder constant $(\mathbb{P}(A_n))^{-\beta(1+\eta)}$.

Now, we apply the decay of correlations to the Hölder continuous function $\phi$ against $\psi = \mathbf{1}_{\mathcal{W}_{0,\ell}(A_n)}$ to get

$$\left| \int \phi \cdot (\mathbf{1}_{\mathcal{W}_{0,\ell}(A_n)} \circ f^t) d\mathbb{P} - \int \phi d\mathbb{P} \int \mathbf{1}_{\mathcal{W}_{0,\ell}(A_n)} d\mathbb{P} \right| \leq C (\mathbb{P}(A_n))^{-\beta(1+\eta)} \varrho(t).$$

Observe that the support of $\mathbf{1}_{A_n} - \phi$ is contained in $A_n \setminus D_n$ whose Lebesgue measure is $O\left((\mathbb{P}(A_n))^{d+\eta}\right)$ and, using (4.4.3), we get that $\int \mathbf{1}_{A_n} - \phi d\mathbb{P} \leq O\left((\mathbb{P}(A_n))^{\Theta(d+\eta)}\right)$. It follows that

$$\begin{aligned} |\mathbb{P}(A_n \cap f^{-t}(\mathcal{W}_{0,\ell}(A_n))) &- \mathbb{P}(A_n)\mathbb{P}(\mathcal{W}_{0,\ell}(A_n))| \\ &\leq (\mathbb{P}(A_n))^{-\beta(1+\eta)} \varrho(t) + O\left((\mathbb{P}(A_n))^{\Theta(d+\eta)}\right). \end{aligned}$$

Hence, we take $\gamma(n, t) = O\left((\mathbb{P}(A_n))^{-\beta(1+\eta)}\right) \varrho(t) + O\left((\mathbb{P}(A_n))^{\Theta(d+\eta)}\right)$. Let $t_n$ be as in the hypothesis and recalling that $\mathbb{P}(A_n) \sim \theta\tau/n$, where $\theta < 1$ is given by (R2) if $\zeta$ is periodic and $\theta = 1$ otherwise, we have that $n\gamma(n, t_n) \leq O\left(n^{1+\beta(1+\eta)}\right) \varrho(t_n) + O\left(n^{-\delta}\right) \xrightarrow[n\to\infty]{} 0$. As before, to prove $D_3(u_n)$ the argument is the same except for the fact that we need to take $\psi = \mathbf{1}_{\mathcal{N}(A)=0}$. In order to prove $D^p(u_n)^*$ we just need to follow the proof as before and use a Hölder continuous approximation for $\mathbf{1}_{Q_{p,0}^{\kappa_1}(u_n)}$. The only extra difficulty is that we need an upper bound that works for all $\kappa_1 \in \mathbb{N}_0$. Now, taking $\psi = \mathbf{1}_{\left(\cap_{j=2}^{\varsigma} \mathcal{N}_{u_n}(I_j)=\kappa_j\right)}$, recalling that $\mathbb{P}(Q_{p,0}^{\kappa}) \sim \theta(1-\theta)^{\kappa} \mathbb{P}(X_0 > u_n)$ and following the same steps, we easily conclude that (4.4.2) leads to the following estimate

$$\left| \mathbb{P}\left(Q_{p,0}^{\kappa_1}(u_n) \cap \left(\cap_{j=2}^{\varsigma} \mathcal{N}_{u_n}(I_j)=\kappa_j\right)\right) - \mathbb{P}\left(Q_{p,0}^{\kappa_1}(u_n)\right) \mathbb{P}\left(\cap_{j=2}^{\varsigma} \mathcal{N}_{u_n}(I_j)=\kappa_j\right)\right|$$
$$\leq C\left(\left(\frac{n}{(1-\theta)^{\kappa_1}}\right)^{\beta(1+\eta)} \varrho(t) + \left(\frac{(1-\theta)^{\kappa_1}}{n}\right)^{\Theta(d+\eta)}\right),$$

for any $\eta > 0$, some $C > 0$ and where $d$ is the dimension of $\mathcal{X}$. Now, we have to be cautious because the first term in the right hand side explodes as $\kappa_1 \to \infty$. However, the trivial observation:

$$\left| \mathbb{P}\left(Q_{p,0}^{\kappa_1}(u_n) \cap \left(\cap_{j=2}^{\varsigma} \mathcal{N}_{u_n}(I_j)=\kappa_j\right)\right) - \mathbb{P}\left(Q_{p,0}^{\kappa_1}(u_n)\right) \mathbb{P}\left(\cap_{j=2}^{\varsigma} \mathcal{N}_{u_n}(I_j)=\kappa_j\right)\right|$$
$$\leq 2\mathbb{P}\left(Q_{p,0}^{\kappa_1}(u_n)\right),$$

allows us to set:

$$\gamma(n, t) = \min_{\kappa_1 \in \mathbb{N}_0} \left\{ 2\theta(1-\theta)^{\kappa_1} \mathbb{P}(X_0 > u_n), \right.$$
$$\left. C\left(\left(\frac{n}{(1-\theta)^{\kappa_1}}\right)^{\beta(1+\eta)} \varrho(t) + \left(\frac{(1-\theta)^{\kappa_1}}{n}\right)^{\Theta(d+\eta)}\right)\right\}.$$

Since by assumption, there is a sequence $(t_n)_{n\in\mathbb{N}}$ such that $n^{1+\beta(1+\eta)}\varrho(t_n) + n^{-\delta} \xrightarrow[n\to\infty]{} 0$, then $n\gamma(n, t_n) \to 0$ as $n \to \infty$, as required. $\qquad\square$

## 4.5
## Specific Dynamical Systems where the Dichotomy Applies

In Chapter 6, we compile a list of examples for which the existence of EVL has been proved. However, before we finish this chapter we give a class of maps for which we can prove the dichotomies stated in Theorem 4.2.7 and Corollary 4.3.4.

Let $f : \mathcal{X} \to \mathcal{X}$ be a measurable function as above. For a measurable potential $\phi : \mathcal{X} \to \mathbb{R}$, we define the *pressure* of $(\mathcal{X}, f, \phi)$ to be

$$P(\phi) := \sup_{\mathbb{P} \in \mathcal{M}_f} \left\{ h(\mathbb{P}) + \int \phi \, d\mathbb{P} : -\int \phi \, d\mathbb{P} < \infty \right\},$$

where $h(\mathbb{P})$ denotes the metric entropy of the measure $\mathbb{P}$, see [106] for details. If $\mathbb{P}$ is an invariant probability measure such that $h(\mathbb{P}_\phi) + \int \phi \, d\mathbb{P} = P(\phi)$, then we say that $\mathbb{P}$ is an *equilibrium state* for $(\mathcal{X}, f, \phi)$.

A measure $m$ is called a $\phi$-*conformal* measure if $m(\mathcal{X}) = 1$ and whenever for a Borel set $A$ holds that $f : A \to f(A)$ is a bijection, then $m(f(A)) = \int_A e^{-\phi} \, dm$. Therefore, setting

$$S_n\phi(x) := \phi(x) + \cdots + \phi \circ f^{n-1}(x),$$

if $f^n : A \to f^n(A)$ is a bijection then $m(f^n(A)) = \int_A e^{-S_n\phi} \, dm$.

Note that for example for a smooth interval map $f$, Lebesgue measure is $\phi$-conformal for $\phi(x) := -\log|Df(x)|$. Moreover, if for example $f$ is a topologically transitive quadratic interval map then as in Ledrappier [159], any absolutely continuous invariant measure $\mathbb{P}$ with $h(\mathbb{P}) > 0$ is an equilibrium state for $\phi$. This also holds for the even simpler case of piecewise smooth uniformly expanding maps, which we consider below. This is the case we principally consider in this subsection. For results on more general equilibrium states see [136].

### 4.5.1
### Rychlik Systems

The first class of examples to which we apply our results is the class of interval maps considered by Rychlik in [160], that is given by a triple $(Y, f, \phi)$, where $Y$ is an interval, $f$ a piecewise expanding interval map (possibly with countable discontinuity points) and $\phi$ a certain potential. This class includes, for example, piecewise $C^2$ uniformly expanding maps of the unit interval with the associated physical measures. We refer to [160] or to [49, Section 4.1] for details on the definition of such class and instead give the following list of examples of maps in such class:

- Given $m \in \{2, 3, \ldots\}$, let $f : x \mapsto mx \mod 1$ and $\phi \equiv -\log m$. Then $m_\phi = \mathbb{P}_\phi = \text{Leb}$.

- Let $f : x \mapsto 2x \mod 1$ and, for $\alpha \in (0,1)$, let

$$\phi(x) := \begin{cases} -\log \alpha & \text{if } x \in (0, 1/2) \\ -\log(1-\alpha) & \text{if } x \in (1/2, 1) \end{cases}$$

(and $\phi = -\infty$ elsewhere). Then $m_\phi = \mathbb{P}_\phi$ is the $(\alpha, 1-\alpha)$-Bernoulli measure on $[0, 1]$.

- Let $f : (0, 1] \to (0, 1]$ and $\phi : (-\infty, 0) \to \mathbb{R}$ be defined as $f(x) = 2^k(x - 2^{-k})$ and $\phi(x) := -k \log 2$ for $x \in (2^{-k}, 2^{-k+1}]$. Then $m_\phi = \mathbb{P}_\phi = \text{Leb}$.

In order to prove that the stated dichotomies hold for these systems, we basically need to show that these systems satisfy the conditions of Theorems 4.2.7 and 4.3.5.

In this setting, as in [160], there is a unique $f$-invariant probability measure $\mathbb{P}_\phi \ll m_\phi$ which is also an equilibrium state for $(Y, f, \phi)$ with a strictly positive density $\frac{d\mathbb{P}_\phi}{dm_\phi} \in BV$. Moreover, there exists exponential decay of correlations against $L^1(m_\phi)$, *i.e.*, there exist $C > 0$ and $\beta > 0$, such that for any $\upsilon \in BV$ and

$\psi \in L^1(m_\phi)$ we have

$$\left| \int \psi \circ f^n \cdot \upsilon \, d\mathbb{P}_\phi - \int \psi \, d\mathbb{P}_\phi \int \upsilon \, d\mathbb{P}_\phi \right| \leq C \|\upsilon\|_{BV} \|\psi\|_{L^1(m_\phi)} \, e^{-\beta n}.$$

Assume that $\zeta$ is such that $0 < \frac{d\mathbb{P}_\phi}{dm_\phi}(\zeta) < \infty$ and the observable $\varphi : \mathcal{X} \to \mathbb{R} \cup \{+\infty\}$ is of the form (4.2.3). The regularity of $\mathbb{P}_\phi$ and $\varphi$ guarantee that condition (R1) holds for every such $\zeta$. Besides, if $\zeta$ is a *repelling $p$-periodic point*, which means that $f^p(\zeta) = \zeta$, $f^p$ is differentiable at $\zeta$ and $0 < |\det D(f^{-p})(\zeta)| < 1$. As shown in [49, Theorem 5], we have that (R2) holds. Moreover, the EI is given by the formula $\theta = 1 - e^{S_p \phi(\zeta)}$.

Finally, since $U(u_n)$ is an interval and $Q_{p,0}^{\kappa_1}(u_n)$ is the union of two intervals, for all $\kappa_1$, we have that $\|\mathbf{1}_{U(u_n)}\|_{BV} \leq 3$ and $\|\mathbf{1}_{Q_{p,0}^{\kappa_1}(u_n)}\|_{BV} \leq 5$, which means all the assumptions of Theorems 4.2.7 and 4.3.5 hold.

### 4.5.2
### Piecewise Expanding Maps in Higher Dimensions

The second class of examples we consider here corresponds to a higher dimensional version of the piecewise expanding interval maps of the previous section. We refer to [158, Section 2] for precise definition of this class of maps and give a very particular example corresponding to a uniformly expanding map on the 2-dimensional torus:

- let $\mathbb{T}^2 = \mathbb{R}^2/\mathbb{Z}^2$ and consider the map $f : \mathbb{T}^2 \to \mathbb{T}^2$ defined by the action of a $2 \times 2$ matrix with integer entries and eigenvalues $\lambda_1, \lambda_2 > 1$.

According to [158, Theorem 5.1], there exists an absolutely continuous invariant measure $\mathbb{P}$. Moreover, in [158, Theorem 6.1], it is shown that, on the mixing components, $\mathbb{P}$ enjoys exponential decay of correlations against $L^1$ observables on $V_\alpha$. More precisely, if the map $f$ is as defined above, and if $\mathbb{P}$ is the mixing absolutely continuous invariant measure , then there exist constants $C < \infty$ and $\gamma < 1$ such that

$$\left| \int \psi \circ f^n \, h \, d\mathbb{P} \right| \leq C \|\psi\|_{L^1} \|h\|_\alpha \gamma^n, \; \forall h \in V_\alpha \text{ and } \forall \psi \in L^1, \text{ where } \int \psi \, d\mathbb{P} = 0.$$

Assume that the observable $\varphi : \mathcal{X} \to \mathbb{R} \cup \{+\infty\}$ is of the form (4.2.3). This guarantees that condition (R1) holds. If $\zeta$ is a *repelling $p$-periodic point*, which means that $f^p(\zeta) = \zeta$, $f^p$ is differentiable at $\zeta$ and $0 < |\det D(f^{-p})(\zeta)| < 1$. Then condition (R2) holds and the EI is equal to $\theta = 1 - |\det D(f^{-p})(\zeta)|$ (see [49, Theorem 3]). It is also easy to check that $\|\mathbf{1}_{U(u_n)}\|_\alpha$ and $\|\mathbf{1}_{Q_{p,0}^{\kappa_1}(u_n)}\|_\alpha$ are bounded by a positive constant, for all $\kappa_1$, which means that all conditions of Theorems 4.2.7 and 4.3.5 hold. $\mathbb{R}$

## 4.6
## Extreme Value Laws for Physical Observables

Let $(\mathcal{X}, \mathcal{B}, \mu, f)$ be a dynamical system, where $\mathcal{X}$ is a $d$-dimensional Riemannian manifold, $f : \mathcal{X} \to \mathcal{X}$ a measurable map and $\mu$ an $f$-invariant probability measure. Assume that there is a compact invariant set $\Omega \subset \mathcal{X}$ which supports the measure $\mu$. Specifically, our main interest is the situation where $\Omega$ is a strange attractor and $\mu$ is a Sinai-Ruelle-Bowen (SRB) measure.

By a strange attractor we mean a compact invariant set $\Omega$ attracting the orbits of a non-empty interior set of points and having a dense orbit with positive Lyapunov exponent (*sensitive dependence on initial condition*, see [70, 71]).

We recall that an SRB measure is an invariant measure of a dynamical system that is additionally characterised by having absolutely continuous conditional measure on unstable manifolds. As well known, it is possible to provide other equivalent definitions of the SRB measure, see [70, 71, 161].

Given an observable $\phi : \mathcal{X} \to \mathbb{R} \cup \{+\infty\}$, previous sections have considered the cases where $\phi$ has the form given in (4.2.3), in which case, the observable can be seen as a function depending directly on the distance to a specific point $\zeta$ chosen in the phase space, and most typically belonging to the attractor.

However, typical observable functions used in applications are not of this form, at least when $\text{dist}(\cdot, \cdot)$ is taken to be the ambient (usually Euclidean) metric. In certain geophysical applications, see [81], observables can take the form

$$A_E(x) = x^T E x, \quad A_W(x) = ||Wx||, \quad A_V(x) = Vx, \quad \text{respectively,} \quad (4.6.1)$$

where $x$ is a point in the phase space $\mathcal{X} = \mathbb{R}^d$, $|| \cdot ||$ denotes the Euclidean norm and $E \in \mathbb{R}^{d \times d}$, $W \in \mathbb{R}^{2 \times d}$, $V \in \mathbb{R}^{1 \times d}$ are matrices. When $\text{dist}(\cdot, \cdot)$ is the Euclidean metric then none of the observables in (4.6.1) has the form (4.2.3), assuming that the origin (which can be thought as representing a state of rest) does not belong to the attractor. In fact this situation is to be expected in many observables found in applications, including the atmospheric and oceanic models of [162]. Although observables such as (4.6.1) are usually unbounded in the system's phase space, the system attractor $\Omega$ is usually bounded due to the presence of dissipative processes in the models. Therefore, time series of such observables should be expected to have an upper bound and, hence, large values typically obey Weibull limit distributions.

In the discussion that follows we consider the invariant measure $\mu$ and the observable $A$ given, with $A$ of the form (4.6.1). If the attractor $\Omega$ is compact, then there exists a point $x_0 \in \Omega$ where the observable is maximized. An alternative approach to the problem could be to find a function $\psi : \mathbb{R}^+ \to \mathbb{R}$ and a metric in $\mathcal{X}$ such that the given observable $A$ can be rewritten in the form (4.2.3). In some particular cases $\psi$ and the metric can be made explicit, but in general finding this adapted construction may be just as difficult as working with the original observable given by the problem. Indeed the adapted metric would also depend on the geometry of the attractor.

Given an observable $A : \mathcal{X} \to \mathbb{R}$ and a threshold $u \in \mathbb{R}$, we define the level regions $L^+(u)$ (resp. level sets $L(u)$) as follows:

$$L^+(u) = \{x \in \mathcal{X} : A(x) \geq u\}, \qquad L(u) = \{x \in \mathcal{X} : A(x) = u\}. \quad (4.6.2)$$

We consider observables which achieve a finite maximum within $\Omega$, although the observable themselves could be unbounded in $\mathcal{X}$. We define

$$u_F = \sup_{x \in \Omega} A(x). \tag{4.6.3}$$

Since $\Omega$ is compact there exists (at least) one point $x_0 \in \Omega$ for which $A(x_0) = u_F$. We will assume that such an *extremal point* $x_0$ is unique. Given our focus on the Weibull case we consider sequences $u_n := u/a_n + b_n$, and let

$$\tau_n(u) := \lim_{n \to \infty} n\mu(L^+(u_n)), \tag{4.6.4}$$

and investigate whether the limit $\tau(u) := \lim_{n \to \infty} \tau_n(u)$ exists. If $L^+(u)$ is in the domain of attraction of a Type III distribution, then following [1], we can choose $b_n = u_F$, and we take $a_n \to \infty$. The precise form of $a_n$ depends on the regularity of $A$, and the regularity of the density of $\mu$ in the vicinity of the extremal point $\tilde{p}$ (if such a density exists). As before, we consider the process $M_n = \max(X_0, \ldots, X_{n-1})$ with $X_n = A \circ f^n$, and we investigate to what extent the following statement is true:

$$n\mu(\{x : A(x) \geq u_n\}) \to \tau(u) \quad \Leftrightarrow \quad \mu(\{M_n \leq u_n\}) \to e^{-\tau(u)}. \tag{4.6.5}$$

If $\tau(u) = u^\alpha$, then the process $M_n$ is described by a Type III extreme value distribution. The statement (4.6.5) and its analogues are known to hold for a wide class of dynamical systems, such as those governed by non-uniformly expanding maps, and systems with certain (non-)uniformly hyperbolic attractors, [163, 138]. These systems will be discussed in Chapter 6. However much of the theory for analysing extremes assumes that the level regions $L^+(u)$ introduced in (4.6.2) are described by balls in the ambient metric, and moreover that these balls are centred on points in $\Omega$ that are generic for $\mu$. These assumptions allow for the index $\alpha$ to be expressed in terms of local dimension formulae for measures.

When leading with more physical observables, we do not assume that the level sets are balls (in the Euclidean metric): for example we consider observables of the form $A(p) = A(x_1, \ldots, x_d) = \sum_i |x_i|^{a_i}$, where the level sets have cusps or are non-conformal. We also consider observables $A(p) = \sum_i c_i x_i$, for which the level sets are hyperplanes. For observables of these types (also compare with (4.6.1)) the standard machinery does not immediately apply. The first problem is to determine the sequence $u_n$ and the limit $\tau(u)$ defined in (4.6.4). Even if the measure $\mu$ is sufficiently regular then the sequence $u_n$ will depend on the geometry of the attractor close to where $A(p)$ achieves its maximum value on $\Omega$, in addition to depending on the form of $A$.

In Chapter 6 we illustrate the various geometrical scenarios that can arise using an hyperbolic toral automorphism as a simple example. When $\mu$ is a more general SRB measure, then even for uniformly hyperbolic systems (such as the solenoid map) it becomes a non-trivial problem to determine $u_n$ and $\tau(u)$. In fact for systems with general SRB measures we do no expect convergence of the left hand side of (4.6.4) to some limit function $\tau(u)$, at least for linear scaling sequences $u_n$.

Even if such a limit $\tau(u)$ exists, we must then check the two conditions $\not{D}_0(u_n)$, $\not{D}'_0(u_n)$. If the level sets have complicated geometry, or if the measure $\mu$ is supported on a fractal set then these conditions must be carefully checked. For uniformly hyperbolic systems, and for observables that are functions of balls in the ambient metric these conditions are checked in e.g. [138].

Later in Section 8.2.2, instead, we take a more heuristic point of view on physical observables. We sacrifice mathematical rigour and derive results using the GPD approach that we claim are correct - within the level of rigour necessary and suitable for physical investigations - for *generic* statistical mechanical systems.

# 5
# Hitting and Return time Statistics

## 5.1
## Introduction to Hitting and Return Time Statistics

Rare events in the context of dynamical systems have long been studied in the framework of Hitting Time Statistics (HTS)/Return Time Statistics (RTS) by, *e.g.*, Pitskel [164] Collet, Galves and Schmitt [165, 166], and Hirata [147]. Interestingly, there is a conceptual bridge between the limit laws described in this framework and the EVLs for distance observables. In this chapter, we first introduce the notions of HTS and RTS, and discuss a theorem which gives an explicit link between these classical concepts. Some of the results on HTS/RTS have long been presented and discussed in the literature, so that we will mostly just state results, explaining the main ideas, but omitting the proofs. Since the link between these two kinds of statistical laws and the EVL point of view is central to the development of the ideas of this book, we will give a full proof of the relevant theorem.

Since the HTS properties can be pulled back to the features of the corresponding RTS, we present some basic findings on specific systems by focussing on HTS only. We start by considering uniformly hyperbolic dynamical systems. We note that early results in this context were restricted in their application to particular sets shrinking down to $\zeta$ (cylinders rather than balls). Next we give more recent results for non-uniformly hyperbolic dynamical systems. One of the techniques here is inducing, where essentially a speeded up version of the dynamical system, which is uniformly hyperbolic, has the same recurrence properties as the original system. We close by giving a brief outline of systems which have been shown to have unusual recurrence laws in the context of Hitting Time Statistics. So, due to the link to the properties of extremes, this means that we also have unusual EVLs for distance observables.

Before delving into the more technical material contained in this chapter, we give some idea of the perspective of hitting times. Given a measure space $(\mathcal{X}, \mathcal{B})$, a measurable function $f : \mathcal{X} \to \mathcal{X}$ and a point $x \in \mathcal{X}$, we denote the iterates of $x$ by $f$ as $x, f(x), f^2(x), \ldots f^n(x), \ldots$. We can study the recurrence properties of this system, *i.e.*, understand how the iterates of *typical* $x$ make repeated visits to certain parts of the system, using the framework of rare events. Namely, we will fix a point $\zeta$, and a rare event will be an occurrence of a point $x$ having some iterate $f^n(x)$ very close to $\zeta$. One way to quantify *very close* is to assume that $\mathcal{X}$ is a metric space endowed with the Borel $\sigma$-algebra with metric denoted by $d : \mathcal{X} \times \mathcal{X} \to [0, \infty]$, and consider entries of iterates of typical $x$ to the $\varepsilon$-ball around $\zeta$ denoted $B_\varepsilon(\zeta)$. When $\varepsilon$ is very small, we would expect most points to take a long time before they hit the ball $B_\varepsilon(\zeta)$. This idea is discussed in detail below. We assume the existence of a probability measure $\mu$ on $\mathcal{X}$, *i.e.*, $\mu(\mathcal{X}) = 1$ which is, moreover, $f$-invariant, *i.e.* $\mu \circ f^{-1} = \mu$, and ergodic. We denote the set of such measures by $\mathcal{M}_f$. Specific members of $\mathcal{M}_f$ will be discussed later in this and the following chapters.

Note that in the dynamical systems context, a measure which in other parts of this book is written as $\mathbb{P}$ is often called $\mu$. Also, while iteration of dynamical systems has been discussed in previous chapters, we will go over this notion again here for the sake of clarity. *Repetita iuvant*.

### 5.1.1
### Definition of Hitting and Return Time Statistics

We start by taking a set $A \subset \mathcal{X}$ and we recall the definition of *first hitting time* given in Eq. 2.3.1 as $r_A : \mathcal{X} \to \mathbb{N} \cup \{\infty\}$:

$$r_A(x) := \inf\{n \in \mathbb{N} : f^n(x) \in A\}.$$

Note that if $x$ never enters $A$ under iteration by $f$ then $r_A(x) = \infty$. If we restrict our attention to points $x$ starting in $A$, the $r_A(x)$ is also called the *first return time* to $A$. For $A \in \mathcal{B}$, the relevant measure on these specific points is the *conditional measure* on $A$, where the measure of a set $B \in \mathcal{B}$ is defined as:

$$\mu_A(B) := \frac{\mu(B \cap A)}{\mu(A)}. \tag{5.1.1}$$

The first result is that if $A$ is measurable then the expected value of $r_A$ is the inverse of the measure of $A$ itself:

$$\int r_A \, d\mu_A = \frac{1}{\mu(A)}.$$

This is known as Kac's Theorem (see, *e.g.*, [167, Theorem 2.44]). Given this fact, in order to learn about the first return time function, it makes sense to normalise by the expectation value of $r_A$, *i.e.*, to study the distribution of $\mu(A) \cdot r_A$. It is common to approach this as the study of

$$\mu_A \left( \left\{ x \in A : r_A(x) > \frac{t}{\mu(A)} \right\} \right)$$

for all values $t \geq 0$. In order to study the asymptotic recurrence to a given point $\zeta \in \mathcal{X}$, one can take a family of sets $\{U(u)\}_u$ where $u$ is a real parameter and sets $U(u)$ shrink to $\zeta$ in some way as $u$ tends to some limiting parameter $u_F$.

#### 5.1.1.1  Hitting and Return Time Statistics to balls

First, we will restrict ourselves to the case that $\{U(u)\}_u$ is a family of balls centered in $\zeta$ and having a radius that tends to zero. In particular, we say that $(\mathcal{X}, f, \mu)$ has *Return Time Statistics* (RTS) $\tilde{G}$ if

$$\lim_{\delta \to 0} \mu_{B_\delta(\zeta)} \left( \left\{ x \in B_\delta(\zeta) : r_{B_\delta(\zeta)}(x) > \frac{t}{\mu(B_\delta(\zeta))} \right\} \right) = \tilde{G}(t)$$

for some non-degenerate function $\tilde{G} : [0, \infty] \to [0, 1]$. Returning to general measurable sets $A$, it turns out that $1/\mu(A)$ is also the relevant normalising factor for the first hitting time not restricted to $A$. So for the asymptotics of the first hitting time, say that $(\mathcal{X}, f, \mu)$ has *Hitting Time Statistics* (HTS) $G$ if

$$\lim_{\delta \to 0} \mu \left( \left\{ x \in \mathcal{X} : r_{B_\delta(\zeta)}(x) > \frac{t}{\mu(B_\delta(\zeta))} \right\} \right) = G(t)$$

for some non-degenerate function $G : [0, \infty] \to [0, 1]$.

In most cases we outline below, the functions $\tilde{G}$ and $G$ are of the form $G(t) = e^{-\theta t}$ for $\theta \in (0, 1]$. If the HTS (RTS) law above is $G(t) = e^{-t}$ ($\tilde{G}(t) = e^{-t}$), we say that we have *exponential HTS (RTS)* to balls around $\zeta$.

#### 5.1.1.2  Hitting and Return Time Statistics to Cylinders

Let $\mathcal{P}_0$ denote a partition of $\mathcal{X}$ and define the corresponding pullback partition $\mathcal{P}_n = \bigvee_{i=0}^{n-1} f^{-i}(\mathcal{P}_0)$, where $\vee$ denotes the join of partitions. We refer to the elements of the partition $\mathcal{P}_n$ as *cylinders of order $n$*. For every $\zeta \in \mathcal{X}$, we denote by $Z_n[\zeta]$ the cylinder of order $n$ that contains $\zeta$. For some $\zeta \in \mathcal{X}$ this cylinder may not be unique (this is the case when $\zeta$ belongs to the border of $\mathcal{P}_n$) , but we can make an arbitrary choice, so that $Z_n[\zeta]$ is well defined. In good cases, as in most cases where the system is uniformly expanding, $Z_n[\zeta] \to \{\zeta\}$ as $n \to \infty$. Note, however, that for our laws we generally only require that $\lim_{n \to \infty} \mu(Z_n[\zeta]) = 0$.

Now we let $\{U(u)\}_u$ be the (countable) set of cylinders $\{Z_n[\zeta]\}_n$. So we say that the system has HTS (RTS) $G$ ($\tilde{G}$) to cylinders at $\zeta$ if we have HTS (RTS) $G$ ($\tilde{G}$) when $U(u)$ is replaced by the cylinder $Z_n(\zeta)$, and the limit is taken as $n \to \infty$.

### 5.2
### HTS vs RTS and Possible Limit Laws

Tthe elegant paper [168] presented a complete description of the relationship between HTS and RTS. In order to describe the main aspects of the connection between HTS and RTS, we first need to clarify some notions of convergence of distributions and also set up our classes of possible distributions.

Given a sequence of d.f.'s $(F_n)_{n \in \mathbb{N}}$, we say that the sequence *converges weakly* to a function $F$ if $F$ is non-increasing right-continuous and satisfies $\lim_{n \to \infty} F_n(t) = F(t)$ at every point $t$ of continuity of $F$. If this is true, then the notation we use is $F_n \Rightarrow F$.

Set

$$\mathcal{F} := \left\{ \begin{array}{l} F : \mathbb{R} \to [0,1], \ F \equiv 0 \text{ on } (-\infty, 0], \ F \text{ is non-decreasing,} \\ \text{continuous and concave on } [0, +\infty), \ F(t) \leq t \text{ for } t \geq 0 \end{array} \right\}.$$

and

$$\tilde{\mathcal{F}} := \left\{ \begin{array}{l} \tilde{F} : \mathbb{R} \to [0,1], \ \tilde{F} \equiv 0 \text{ on } (-\infty, 0], \ \tilde{F} \text{ is non-decreasing,} \\ \text{right continuous and } \displaystyle\int_0^{+\infty} (1 - \tilde{F}(s)) \, ds \leq 1 \end{array} \right\}.$$

Given a dynamical system $(\mathcal{X}, f)$, for $U \subset \mathcal{X}$, we let

$$F_U(t) := \mu(\{\mu(U) r_U \leq t\})$$

and

$$\tilde{F}_U(t) := \frac{1}{\mu(U)} \mu(\{U \cap \{\mu(U) r_U \leq t\}\}).$$

The idea is that we can choose sequences $(U_n)_{n \in \mathbb{N}}$ such that $F_{U_n}$ will converge to a function $F \in \mathcal{F}$ and $\tilde{F}_{U_n}$ will converge to a function $\tilde{F} \in \tilde{\mathcal{F}}$. These limiting distributions are the tails of the HTS and RTS respectively.

**Theorem 5.2.1** ([168]). *Let $(\mathcal{X}, f, \mu)$ be an ergodic dynamical system and $\{U_n\}_{n \in \mathbb{N}}$ a sequence of positive measure measurable subsets of $\mathcal{X}$. Then $F_{U_n}$ converge weakly if and only if $\tilde{F}_{U_n}$ converge weakly. Moreover, given this convergence, to $F$ and $\tilde{F}$ respectively, we have the following relation:*

$$F(t) = \int_0^t (1 - \tilde{F}(s)) \, ds$$

*for $t \geq 0$, and we must have $F \in \mathcal{F}$ and $\tilde{F} \in \tilde{\mathcal{F}}$.*

The last part of this theorem, the assertion that $F \in \mathcal{F}$ and $\tilde{F} \in \tilde{\mathcal{F}}$, follows from [169, 170]. Furthermore, those papers show that for *any* $F \in \mathcal{F}$, there exists a sequence $U_n$ of sets with $\mu(U_n) \to 0$ such that $F_{U_n} \Rightarrow F$. We will discuss the possible limit laws further in Section 5.2 and Section 5.6. Our main interest is in the fixed point $1 - e^{-t}$ of the integral in Theorem 5.2.1, as well as the distribution $F(t) = 1 - e^{-\theta t}$ with $\theta \in (0, 1)$, which corresponds to $\tilde{F}(t) = 1 - \theta + \theta(1 - e^{-\theta t})$.

Since there is a clear relationship between HTS and RTS and since HTS is most directly linked to EVLs, after this subsection we will restrict ourselves to discussing HTS .

## 5.3
## The Link between Hitting Times and Extreme Values

Note that, by assumption (R1), we have that $f^{-1}(\{M_n \leq u\}) = \{r_{U(u)} > n\}$. This relationship was essentially implicit in the work of Collet in [72] and was mentioned in [171]. In [74], a complete relation between the existence of HTS/RTS and EVLs was rigorously proved for absolutely continuous (with respect to Lebesgue) invariant measure and natural observables depending on the distance to $\zeta$. Then, in [136], by adapting the observables, this relation was further developed to hold for general equilibrium states. Below we prove the link between HTS and EVL in the general case. Note that, in order to make the connection to EVL more explicit, we denote our measures by $\mathbb{P}$, again, as in the previous chapters.

We recall that the stochastic process $X_0, X_1, \ldots$ is defined by (4.2.2) where the observable $\varphi$ is defined by (4.2.11). We stress that $g$ stands for some $g_i$ with $i = 1, 2, 3$ defined by conditions (4.2.4), (4.2.5) and (4.2.6), respectively.

The theorems 5.3.1 and 5.3.2 below provide a *dictionary* for translating results obtained for HTS into EVLs and *vice versa*. Our first main result, which derives EVLs from HTS for balls, reads as follows.

**Theorem 5.3.1** ([136, Theorem 1]). *Let $(\mathcal{X}, \mathcal{B}, \mathbb{P}, f)$ be a dynamical system, $\zeta \in \mathcal{X}$ be in the support of $\mathbb{P}$ and assume that $\mathbb{P}$ is such that the function $\hbar$ defined on (4.2.12) is continuous.*

*If we have HTS $G$ to balls centred on $\zeta \in \mathcal{X}$, then we have an EVL $H = G$ for $M_n$ that applies to the observables (4.2.3) achieving a maximum at $\zeta$. Moreover, for all $i = 1, 2, 3$, if $u_n$ is chosen as linear normalising sequence, as in (3.1.1), then the shape $g_i$ for the observable corresponds to an extremal type law of the form $e^{-\tau_i}$, given in (4.6.4).*

*Proof of Theorem* 5.3.1. Set

$$
\begin{aligned}
& u_n = g_1\left(n^{-1}\right) + h\left(g_1\left(n^{-1}\right)\right) y, && \text{for } y \in \mathbb{R}, \text{ for type } g_1; \\
& u_n = g_2\left(n^{-1}\right) y, && \text{for } y > 0, \text{ for type } g_2; \\
& u_n = D - \left(D - g_3\left(n^{-1}\right)\right)(-y), && \text{for } y < 0, \text{ for type } g_3.
\end{aligned}
$$

For $n$ sufficiently large,

$$
\begin{aligned}
\{x : M_n(x) \leq u_n\} &= \bigcap_{j=0}^{n-1} \{x : X_j(x) \leq u_n\} = \bigcap_{j=0}^{n-1} \left\{x : g\left(\mathbb{P}\left(B_{\mathrm{dist}(f^j(x),\zeta)}(\zeta)\right)\right) \leq u_n\right\} \\
&= \bigcap_{j=0}^{n-1} \left\{x : \mathbb{P}\left(B_{\mathrm{dist}(f^j(x),\zeta)}(\zeta)\right) \geq g^{-1}(u_n)\right\}. \quad (5.3.1)
\end{aligned}
$$

Consequently, by (4.2.14),

$$
\begin{aligned}
\mathbb{P}(\{x : M_n(x) \leq u_n\}) &= \mathbb{P}\left(\bigcap_{j=0}^{n-1} \left\{x : \mathbb{P}\{B_{\mathrm{dist}(f^j(x),\zeta)}(\zeta)\} \geq \mathbb{P}\{B_{\ell(g^{-1}(u_n))}(\zeta)\}\right\}\right) \\
&= \mathbb{P}\left(\bigcap_{j=0}^{n-1} \left\{x : \mathrm{dist}(f^j(x),\zeta) \geq \ell(g^{-1}(u_n))\right\}\right) \\
&= \mathbb{P}\left(\left\{x : r_{B_{\ell(g^{-1}(u_n))}(\zeta)}(x) \geq n\right\}\right). \quad (5.3.2)
\end{aligned}
$$

Now, observe that (4.2.4), (4.2.5) and (4.2.6) imply

$$g_1^{-1}(u_n) = g_1^{-1}\left[g_1(n^{-1}) + p\left(g_1(n^{-1})\right)y\right] \sim g_1^{-1}\left[g_1(n^{-1})\right] e^{-y} = \frac{e^{-y}}{n};$$

$$g_2^{-1}(u_n) = g_2^{-1}\left[g_2(n^{-1})y\right] \sim g_2^{-1}\left[g_2(n^{-1})\right] y^{-\beta} = \frac{y^{-\beta}}{n};$$

$$g_3^{-1}(u_n) = g_3^{-1}\left[D - \left(D - g_3(n^{-1})\right)(-y)\right]$$

$$\sim g_3^{-1}\left[D - \left(D - g_3(n^{-1})\right](-y)^\gamma = \frac{(-y)^\gamma}{n}.$$

Thus, we may write

$$g^{-1}(u_n) \sim \frac{\tau(y)}{n},$$

meaning that

$$g_i^{-1}(u_n) \sim \frac{\tau_i(y)}{n}, \quad \forall i \in \{1, 2, 3\}$$

where $\tau_1(y) = e^{-y}$ for $y \in \mathbb{R}$, $\tau_2(y) = y^{-\beta}$ for $y > 0$, and $\tau_3(y) = (-y)^\gamma$ for $y < 0$.

Recalling (4.2.14), we have

$$\mathbb{P}\left(B_{\ell(g^{-1}(u_n))}(\zeta)\right) \sim \frac{\tau(y)}{n},$$

and so,

$$n \sim \frac{\tau(y)}{\mathbb{P}\left(B_{\ell(g^{-1}(u_n))}(\zeta)\right)}. \tag{5.3.3}$$

Now, we claim that using (5.3.2) and (5.3.3), we have

$$\lim_{n \to \infty} \mathbb{P}(M_n(x) \leq u_n) = \lim_{n \to \infty} \mathbb{P}\left(r_{B_{l(g^{-1}(u_n))}(\zeta)}(x) \geq \frac{\tau(y)}{\mathbb{P}\left(B_{\ell(g^{-1}(u_n))}(\zeta)\right)}\right) \tag{5.3.4}$$

$$= 1 - G(\tau(y)), \tag{5.3.5}$$

which gives the first part of the theorem.

To see that (5.3.4) holds, observe that by (5.3.2) and (5.3.3) we have

$$\left|\mathbb{P}(M_n \leq u_n) - \mathbb{P}\left(r_{B_{\ell(g^{-1}(u_n))}(\zeta)} \geq \frac{\tau(y)}{\mathbb{P}\left(B_{\ell(g^{-1}(u_n))}(\zeta)\right)}\right)\right|$$

$$= \left|\mathbb{P}\left(r_{B_{\ell(g^{-1}(u_n))}(\zeta)} \geq n\right) - \mathbb{P}\left(r_{B_{\ell(g^{-1}(u_n))}(\zeta)} \geq (1 + \varepsilon_n)n\right)\right|,$$

where $(\varepsilon_n)_{n \in \mathbb{N}}$ is such that $\varepsilon_n \to 0$ as $n \to \infty$. Since we have

$$\left\{r_{B_{\ell(g^{-1}(u_n))}(\zeta)} \geq m\right\} \setminus \left\{r_{B_{\ell(g^{-1}(u_n))}(\zeta)} \geq m + k\right\}$$

$$\subset \bigcup_{j=m}^{m+k-1} f^{-j}\left(B_{\ell(g^{-1}(u_n))}(\zeta)\right), \quad \forall m, k \in \mathbb{N}, \tag{5.3.6}$$

it follows by stationarity that

$$\left| \mathbb{P}\left( r_{B_{\ell(g^{-1}(u_n))}}(\zeta) \geq n \right) - \mathbb{P}\left( r_{B_{\ell(g^{-1}(u_n))}}(\zeta) \geq (1 + \varepsilon_n)n \right) \right|$$
$$\leq |\varepsilon_n| n \mathbb{P}\left( B_{\ell(g^{-1}(u_n))}(\zeta) \right) \sim |\varepsilon_n| \tau \to 0,$$

as $n \to \infty$, completing the proof of (5.3.4). $\qquad\square$

Now, we state a result in the other direction, *i.e.*, we show how to get HTS from EVLs for balls.

**Theorem 5.3.2** ([136, Theorem 2]). *Let $(\mathcal{X}, \mathcal{B}, \mathbb{P}, f)$ be a dynamical system, $\zeta \in \mathcal{X}$ be in the support of $\mathbb{P}$ and assume that $\mathbb{P}$ is such that the function $\hbar$ defined in (4.2.12) is continuous.*

*If we have an EVL $H$ for $M_n$ which applies to the observables (4.2.3) achieving a maximum at $\zeta \in \mathcal{X}$ then we have HTS $G = H$ to balls at $\zeta$.*

For the proof of Theorem 5.3.2, we require the following lemma. This is essentially contained in [1, Theorem 1.6.2]. See also [74, Lemma 2.1] where the lemma was proved for absolutely continuous invariant measures. We provide a proof in the general case for completeness.

**Lemma 5.3.3.** *Let $(\mathcal{X}, \mathcal{B}, \mathbb{P}, f)$ be a dynamical system, $\zeta \in \mathcal{X}$ and assume that $\mathbb{P}$ is such that the function $\hbar$ defined on (4.2.12) is continuous. Furthermore, let $\varphi$ be as in (4.2.3). Then, for each $y \in \mathbb{R}$, there exist sequences $(a_n)_{n \in \mathbb{N}}$ and $(b_n)_{n \in \mathbb{N}}$ so that the sequence $(u_n(y))_{n \in \mathbb{N}}$ as in (3.1.1) is such that*

$$n\mathbb{P}(\{x : \varphi(x) > u_n(y)\}) \xrightarrow[n \to \infty]{} \tau(y) \geq 0.$$

*Moreover, for every $t > 0$ there exists $y \in \mathbb{R}$ such that limit $\tau(y) = t$.*

*Proof.* We will prove the lemma in the case when $g$ is of type $g_2$. For the other two types of $g$ observables, the argument is the same, but with minor adjustments, see [1, Theorem 1.6.2].

First we show that we can always find a sequence $(\gamma_n)_{n \in \mathbb{N}}$ such that

$$n\mathbb{P}(X_0 > \gamma_n) \xrightarrow[n \to \infty]{} 1.$$

Take $\gamma_n := \inf\{y : \mathbb{P}(X_0 \leq y) \geq 1 - 1/n\}$, and let us show that it has the desired property. Note that $n\mathbb{P}(X_0 > \gamma_n) \leq 1$, which means that $\limsup_{n \to \infty} n\mathbb{P}(X_0 > \gamma_n) \leq 1$. Using (4.2.14) and (4.2.5), for any $z < 1$, we have

$$\liminf_{n \to \infty} \frac{\mathbb{P}(X_0 > \gamma_n)}{\mathbb{P}(X_0 > \gamma_n z)} = \liminf_{n \to \infty} \frac{\mathbb{P}(B_{\ell(g_2^{-1}(\gamma_n))}(\zeta))}{\mathbb{P}(B_{\ell(g_2^{-1}(z\gamma_n))}(\zeta))} = \liminf_{n \to \infty} \frac{g_2^{-1}(\gamma_n)}{g_2^{-1}(z\gamma_n)} = z^\beta,$$

where $\ell$ is the function defined in (4.2.13). Since, by definition of $\gamma_n$, for any $z < 1$, $n\mathbb{P}(X_0 > \gamma_n z) \geq 1$, letting $z \to 1$, it follows immediately that $\liminf_{n \to \infty} n\mathbb{P}(X_0 > \gamma_n) \geq 1$.

Let $u_n(y) = \gamma_n y$; thus means that, for all $n \in \mathbb{N}$, we take $a_n = \gamma_n^{-1}$ and $b_n = 0$ in (3.1.1). Hence, using (4.2.5), it follows that for all $y > 0$

$$
\begin{aligned}
n\mathbb{P}(X_0 > \gamma_n y) = n\mathbb{P}(B_{\ell(g_2^{-1}(\gamma_n y))}(\zeta)) &= n g_2^{-1}(\gamma_n y) \\
&\sim n y^{-\beta} g_2^{-1}(\gamma_n) = y^{-\beta} n\mathbb{P}(B_{\ell(g_2^{-1}(\gamma_n))}(\zeta)) \\
&= y^{-\beta} n\mathbb{P}(X_0 > \gamma_n) \xrightarrow[n\to\infty]{} y^{-\beta}.
\end{aligned}
$$

Hence, taking $y = t^{-1/\beta} > 0$ would suit our purposes. $\qquad\square$

*Proof of Theorem* 5.3.2. We assume that by hypothesis for every $y \in \mathbb{R}$ and some sequence $u_n = u_n(y)$ as in (3.1.1) such that $n\mathbb{P}(\{x : \varphi(x) > u_n(y)\}) \xrightarrow[n\to\infty]{} \tau(y)$, we have

$$
\lim_{n\to\infty} \mathbb{P}(\{x : M_n(x) \le u_n(y)\}) = \bar{H}(\tau(y)) = H(y).
$$

Observe that, by Khintchine's Theorem (see [1, Theorem 1.2.3]), up to linear scaling the normalising sequences are unique, which means that we may assume that they are the ones given by Lemma 5.3.3. Hence given $t > 0$, Lemma 5.3.3 implies that there exists $y \in \mathbb{R}$ such that

$$
n\mathbb{P}(\{x : \varphi(x) > u_n(y)\}) \xrightarrow[n\to\infty]{} t.
$$

Given $(\delta_n)_{n\in\mathbb{N}} \subset \mathbb{R}^+$ with $\delta_n \xrightarrow[n\to\infty]{} 0$, we define

$$
\kappa_n := \lfloor t/\mathbb{P}(B_{\delta_n}(\zeta)) \rfloor.
$$

We will prove

$$
g^{-1}(u_{\kappa_n}) \sim \mathbb{P}(B_{\delta_n}(\zeta)). \tag{5.3.7}
$$

If $n$ is sufficiently large, then

$$
\{x : \varphi(x) > u_n\} = \{x : g(\mathbb{P}(B_{\mathrm{dist}(x,\zeta)}(\zeta))) > u_n\} = \{x : \mathbb{P}(B_{\mathrm{dist}(x,\zeta)}(\zeta)) < g^{-1}(u_n)\}.
$$

By (4.2.14) and the definition of $\ell$ in (4.2.13) we obtain

$$
\begin{aligned}
\mathbb{P}(\{x : \varphi(x) > u_n\}) &= \mathbb{P}\left(\{x : \mathbb{P}(B_{\mathrm{dist}(x,\zeta)}(\zeta)) < g^{-1}(u_n)\}\right) \\
&= \mathbb{P}\left(\{x : \mathbb{P}(B_{\mathrm{dist}(x,\zeta)}(\zeta)) < \mathbb{P}(B_{\ell(g^{-1}(u_n))}(\zeta))\}\right) \\
&= \mathbb{P}\left(\{x : \mathrm{dist}(x,\zeta) < \ell(g^{-1}(u_n))\}\right) \\
&= \mathbb{P}\left(B_{\ell(g^{-1}(u_n))}(\zeta)\right).
\end{aligned}
$$

Hence, by assumption on the sequence $u_n$, we have $n\mathbb{P}\left(B_{\ell(g^{-1}(u_n))}(\zeta)\right) \xrightarrow[n\to\infty]{} \tau(y) = t$. As we know that $\mathbb{P}\left(B_{\ell(g^{-1}(u_n))}(\zeta)\right) = g^{-1}(u_n)$, we have $ng^{-1}(u_n) \xrightarrow[n\to\infty]{} t$. Thus, we can write $g^{-1}(u_n) \sim \frac{t}{n}$ and, after substituting $n$ by $\kappa_n$, we immediately obtain (5.3.7) by definition of $\kappa_n$.

Again, by the definition of $\ell$ in (4.2.13) and (4.2.14) we note that

$$
\begin{aligned}
\mathbb{P}(\{x : M_{\kappa_n}(x) \leq u_{\kappa_n}\}) &= \mathbb{P}\left(\bigcap_{j=0}^{\kappa_n-1}\{x : \mathbb{P}\{B_{\mathrm{dist}(f^j(x),\zeta)}(\zeta)\} \geq g^{-1}(u_{\kappa_n})\}\right) \\
&= \mathbb{P}\left(\bigcap_{j=0}^{\kappa_n-1}\left\{x : \mathbb{P}\{B_{\mathrm{dist}(f^j(x),\zeta)}(\zeta)\} \geq \mathbb{P}\{B_{\ell(g^{-1}(u_{\kappa_n}))}(\zeta)\}\right\}\right) \\
&= \mathbb{P}\left(\bigcap_{j=0}^{\kappa_n-1}\left\{x : \mathrm{dist}(f^j(x),\zeta) \geq \ell(g^{-1}(u_{\kappa_n}))\right\}\right) \\
&= \mathbb{P}\left(\left\{x : r_{B_{\ell(g^{-1}(u_{\kappa_n}))}(\zeta)}(x) \geq \kappa_n\right\}\right). \quad (5.3.8)
\end{aligned}
$$

We claim that

$$
\lim_{n\to\infty} \mathbb{P}\left(\left\{x : r_{B_{\delta_n}(\zeta)}(x) \geq \frac{t}{\mathbb{P}(B_{\delta_n}(\zeta))}\right\}\right) = \lim_{n\to\infty} \mathbb{P}(\{x : M_{\kappa_n}(x) \leq u_{\kappa_n}\}). (5.3.9)
$$

Using such a claim the first part of the theorem easily follows, since by hypothesis,

$$
\mathbb{P}\left(\{x : M_{\kappa_n}(x) \leq u_{\kappa_n}\}\right) \xrightarrow[n\to\infty]{} \bar{H}(\tau(y)) = \bar{H}(t).
$$

We need to show that (5.3.9) holds. First, we observe that

$$
\begin{aligned}
\mathbb{P}\left(r_{B_{\delta_n}(\zeta)} \geq \frac{t}{\mathbb{P}(B_{\delta_n}(\zeta))}\right) &= \mathbb{P}(M_{\kappa_n} \leq u_{\kappa_n}) + \left(\mathbb{P}\left(r_{B_{\delta_n}(\zeta)} \geq \kappa_n\right) - \mathbb{P}(M_{\kappa_n} \leq u_{\kappa_n})\right) \\
&\quad + \left(\mathbb{P}\left(r_{B_{\delta_n}(\zeta)} \geq \frac{t}{\mathbb{P}(B_{\delta_n}(\zeta))}\right) - \mathbb{P}\left(r_{B_{\delta_n}(\zeta)} \geq \kappa_n\right)\right).
\end{aligned}
$$

For the third term on the right, we note that, using the definition of $\kappa_n$ t, we have

$$
\begin{aligned}
&\left|\mathbb{P}\left(r_{B_{\delta_n}(\zeta)} \geq \kappa_n\right) - \mathbb{P}\left(r_{B_{\delta_n}(\zeta)} \geq \frac{t}{\mathbb{P}(B_{\delta_n}(\zeta))}\right)\right| \\
&= \left|\mathbb{P}\left(r_{B_{\delta_n}(\zeta)} \geq \kappa_n\right) - \mathbb{P}\left(r_{B_{\delta_n}(\zeta)} \geq (1+\varepsilon_n)\kappa_n\right)\right|,
\end{aligned}
$$

for some sequence $(\varepsilon_n)_{n\in\mathbb{N}}$ such that $\varepsilon_n \to 0$, as $n \to \infty$. Using (5.3.6) and stationarity, we derive that

$$
\left|\mathbb{P}\left(r_{B_{\delta_n}(\zeta)} \geq \kappa_n\right) - \mathbb{P}\left(r_{B_{\delta_n}(\zeta)} \geq (1+\varepsilon_n)\kappa_n\right)\right| \leq |\varepsilon_n|\kappa_n \mathbb{P}(B_{\delta_n}(\zeta)) \sim |\varepsilon_n|t \to 0,
$$

as $n \to \infty$.

As for the remaining term, using the definition of $\kappa_n$ and (5.3.8), we have

$$
\left| \mathbb{P}\left( \left\{ r_{B_{\delta_n}(\zeta)} \geq \kappa_n \right\} \right) - \mathbb{P}\left( \left\{ M_{\kappa_n} \leq u_{\kappa_n} \right\} \right) \right|
$$

$$
= \left| \mathbb{P}\left( \left\{ r_{B_{\delta_n}(\zeta)} \geq \kappa_n \right\} \right) - \mathbb{P}\left( \left\{ r_{B_{\ell(g^{-1}(u_{\kappa_n}))}(\zeta)} \geq \kappa_n \right\} \right) \right|
$$

$$
\leq \sum_{i=1}^{\kappa_n} \mathbb{P}\left( f^{-i}\left( B_{\delta_n}(\zeta) \triangle B_{\ell(g^{-1}(u_{\kappa_n}))}(\zeta) \right) \right)
$$

$$
= \kappa_n \mathbb{P}\left( B_{\delta_n}(\zeta) \triangle B_{\ell(g^{-1}(u_{\kappa_n}))}(\zeta) \right)
$$

$$
\sim \frac{t}{\mathbb{P}\left( B_{\delta_n}(\zeta) \right)} \left| \mathbb{P}\left( B_{\delta_n}(\zeta) \right) - \mathbb{P}\left( B_{\ell(g^{-1}(u_{\kappa_n}))}(\zeta) \right) \right|
$$

$$
= t \left| 1 - \frac{\mathbb{P}\left( B_{\ell(g^{-1}(u_{\kappa_n}))}(\zeta) \right)}{\mathbb{P}\left( B_{\delta_n}(\zeta) \right)} \right|,
$$

where $\triangle$ indicates the operation of symmetric difference, and, which, by (4.2.14) and (5.3.7), tends to 0 as $n \to \infty$; this ends the proof of (5.3.9). $\square$

In the rest of the chapter we discuss explicitly only HTS. The reader can then interpret the obtained results in terms of EVLs using the theorems 5.3.1 and 5.3.2 above.

## 5.4
## Uniformly Hyperbolic Systems

HTS were first applied in a dynamical context to study Markov chains. This was later extended to symbolic dynamical systems with good mixing properties, such as $\phi$-mixing. We refer the reader to the excellent reviews in [172, 173] for a more complete discussion of these results. In this chapter, we explore these ideas, focussing on (piecewise) smooth dynamical systems. We begin by considering the uniformly hyperbolic case.

Perhaps, the simplest case of uniformly hyperbolic systems is given by Markov shifts. Let $S := \{1, \ldots, n\}$ be a finite alphabet (note that we extend to an infinite alphabet later), and $A = (a_{i,j})$ be an $n \times n$ matrix of 0's and 1's, and consider the set of sequences

$$
\Sigma_A^+ := \{(x_0, x_1, \ldots) : a_{x_i, x_{i+1}} = 1 \text{ for each } i \in \mathbb{N}_0 \}.
$$

Any finite word $(x_0, x_1, \ldots, x_{k-1}) \in \Sigma^k$ of length $k$ with $a_{x_i, x_{i+1}} = 1$ for each $i = 0, \ldots, k-1$ is called *allowable*, and the set of all such words, where $k$ is any finite number is denoted $\Sigma_A^+$. Then the corresponding *one-sided subshift of finite type* (SFT) is the pair $(\Sigma_A^+, \sigma)$ where $\sigma$ is the left-shift map: $\sigma(x_0, x_1, \ldots) = (x_1, x_2, \ldots)$. If $A$ is irreducible[1], then we say that $(\Sigma_A^+, \sigma)$ is *topologically mixing*.

---

[1] A matrix is irreducible if if it is not similar via a permutation to a block upper triangular matrix that has more than one block of positive size.

A special case of this is the *one-sided full shift on $n$ symbols*, were $a_{i,j} = 1$ for all $i, j \in \{1, \ldots, n\}$.

An SFT has a canonical *cylinder* structure: for any allowable word $(x_0, x_1, \ldots, x_{k-1}) \in \Sigma_A^*$, the corresponding *$k$-cylinder* is the set

$$\{y = (y_0, y_1, \ldots) \in \Sigma_A^+ : y_i = x_i \text{ for } i = 0, \ldots, k-1\}.$$

Given $x \in \Sigma_A^+$, let $Z_k[x]$ denote the $k$-cylinder containing $x$. We denote the set of $k$-cylinders by $\mathcal{P}_k$. Hence, formally, $\mathcal{P}_0$ denotes the whole of $\Sigma_A^+$. Note that later we will use the same language and notation to describe dynamical cylinders in settings other than the purely symbolic one.

The cylinder structure induces a topology on $\Sigma_A^+$, which is the same as that induced by the metric $d : \Sigma_A^+ \times \Sigma_A^+ \to [0, 1]$, where $d(x, y) = 2^{-k}$, when $x$ and $y$ are in the same $k$-cylinder, but in different $(k+1)$-cylinders.

### 5.4.1
### Gibbs Measures

Given a function $\phi : \Sigma_A^+ \to \mathbb{R}$, the $k^{th}$ *Variation* is defined as

$$V_n(\phi) := \sup_{Z_k \in \mathcal{P}_k} \sup_{x, y \in Z_k} |\phi(x) - \phi(y)|.$$

We say that $\phi$ is Hölder if there exists $\alpha > 0$ such that for $k \in \mathbb{N}_0$, $V_k(\phi) = O(e^{-\alpha k})$.

Given an SFT $(\Sigma_A^+, \sigma)$ and a Hölder potential $\phi$, there exists a $\sigma$-invariant measure $\mu_\phi$ which has the following property: There exist $K \geq 1$ and $P \in \mathbb{R}$ where

$$\frac{1}{K} \leq \frac{\mu_\Phi(Z_k)}{e^{S_k \phi(x) - kP}} \leq K \qquad \text{for any } x \in Z_k,$$

where $S_k \phi(x)$ denotes the $k^{th}$ ergodic sum $\phi(x) + \phi(\sigma x) + \cdots + \phi(\sigma^{k-1} x)$. This is called the *Gibbs property*, and we call $P$ the *Gibbs constant*. Note that replacing $\phi$ by $\phi - P$ makes the Gibbs constant zero. See [174, 175, 176] for a proof of that $\mu_\phi$ exists and satisfies these properties.

### 5.4.2
### First HTS theorem

In particular, the Gibbs property means that $\mu_\phi$ has a strong form of mixing, which is used in [177] to prove the following statement, also proved in this context in [147]:

For $\mu_\phi$-a.e. $\zeta \in \Sigma_A^+$,

$$\lim_{\delta \to 0} \mu_{B_\delta(\zeta)} \left( \left\{ x \in B_\delta(\zeta) : r_{B_\delta(\zeta)}(x) > \frac{t}{\mu(B_\delta(\zeta))} \right\} \right) = e^{-t}.$$

We will formalise these results in a theorem using the more complete picture, given by the full treatment made in [146], which builds up on the work developed in [147] and also, to some extent, in [150]. The method there was a transfer operator approach.

**Theorem 5.4.1.** *Let $(\Sigma_A^+, \sigma)$ be a topologically transitive SFT with $\mu_\phi$ a $\sigma$-invariant Gibbs measure for a Hölder potential $\phi : \Sigma_A^+ \to \mathbb{R}$ with Gibbs constant $P = 0$. Then, for $\zeta \in \Sigma_A$,*

$$\lim_{\delta \to 0} \mu_{B_\delta(\zeta)} \left( \left\{ x \in B_\delta(\zeta) : r_{B_\delta(\zeta)}(x) > \frac{t}{\mu(B_\delta(\zeta))} \right\} \right) = e^{-\theta t},$$

*where*

$$\theta = \begin{cases} 1 - e^{S_k \phi(\zeta)} & \text{if } \zeta \text{ is periodic of prime period } k, \\ 1 & \text{otherwise} \end{cases}$$

Note that, in this case, the symbolic metric implies that balls are, in fact, like cylinders. Also note that a law for *all* points $\zeta$ in $\Sigma_A^+$, rather than just $\mu_\phi$-typical ones and periodic ones can also be found in [165, 178] along with subsequent works along that line such as [179, 180]. Nonetheless, the dichotomy between periodic and non-periodic behaviour, as well as the associated explicit expressions for $\theta$, are not presented in those contributions. We also refer to [172, 173] for a more complete picture of the early developments along that line, which include early estimates on the speed of convergence of the law of HTS.

Theorem 5.4.1 can be extended to a large class of uniformly hyperbolic dynamical systems via coding, as we explain below.

### 5.4.3
### Markov partitions

In some cases, a dynamical system $f : \mathcal{X} \to \mathcal{X}$ is equipped with a finite Markov partition. There are two natural ways to define this: the first being for non-invertible expansive systems, and the second being for invertible hyperbolic dynamical systems in two dimensions or more. Note that these definitions appear in a variety of different ways in the literature.

If $f : \mathcal{X} \to \mathcal{X}$ is non-invertible, let $\Lambda \subset \mathcal{X}$ be the set of points for which $f$ is defined for all time. Moreover, suppose that $\mathcal{P}_1$ is a set of maximal sets such that $f : Z \to f(Z)$ is a homeomorphism for each $Z \in \mathcal{P}_1$. Then $\mathcal{P}_1$ is a finite Markov partition for $f$ if

1) it has finitely many elements;
2) $\Lambda \subset \cup_{Z \in \mathcal{P}_1} Z$;
3) if $x \in \Lambda$ has $x \in Z \in \mathcal{P}_1$ and $f(x) \in Z' \in \mathcal{P}_1$, then $f(Z) \cap \Lambda \subset Z' \cap \Lambda$.
4) $\sup\{d(x, y) : f^i(x), f^i(y) \in Z_i \text{ for some } Z_i \in \mathcal{P}_1 \text{ with } i = 0, \ldots n - 1\} \to 0$ as $n \to \infty$.

In this case, we call $f$ a *Markov map*

This structure yields a coding since writing $\mathcal{P}_1 = \{Z^1, \cdots, Z^n\}$, for any $x \in \Lambda$, there is a code $x = (x_0, x_1, \ldots) \in \{1, \ldots, n\}^{\mathbb{N}_0}$ given by the rule $x_i = k$ if $f^i x \in Z^k$. If the closures of elements of $\mathcal{P}_1$ are pairwise disjoint, or if they do intersect then this intersection takes place outside $\Lambda$, then the coding is uniquely defined. Otherwise, there is ambiguity and the coding is ill-defined in some places. However, since this is usually a zero measure phenomenon, it is not a significant problem.

The Markov structure gives us a set of allowable transitions, namely, the transition from $Z^i$ to $Z^j$ is allowable if $Z^j \subset f(Z^i)$. This also defines a matrix $A$, where $a_{i,j} = 1$ if the transition from $Z^i$ to $Z^j$ is allowable and $a_{i,j} = 0$ otherwise. Thus, we obtain the SFT $(\Sigma_A^+, \sigma)$. It is now standard to look at these properties thorough the lens of the following commuting diagram:

$$
\begin{array}{ccc}
\Sigma_A^+ & \xrightarrow{\ \sigma\ } & \Sigma_A^+ \\
{\scriptstyle\pi}\downarrow & & \downarrow{\scriptstyle\pi} \\
\Lambda & \xrightarrow{\ f\ } & \Lambda
\end{array}
$$

Here, $\pi$ is the inverse of the coding above: we write it in this way so that the functions we deal with are well-defined. The continuous function $\pi : \Sigma_A^+ \to \Lambda$ is a *semi-conjugacy*. As above, if the coding map is uniquely defined, then $\pi$ is a homeomorphism, and is called a *conjugacy*.

We can now start to link the Gibbs measures on $\Sigma_A^+$ to those on $\Lambda$. First, we need to ensure smoothness of potentials, for which it is sufficient to assume that on each $Z \in \mathcal{P}_1$, if $x, y \in Z \cap \Lambda$ then $d(f(x), f(y)) > 1$. Note that this is *not* true for general hyperbolic system, so that we need to add such specific hypothesis. Then the following lemma is easily proved.

**Lemma 5.4.2.** *If $\phi : \cup_{Z \in \mathcal{P}_1} Z \to \mathbb{R}$ is $\alpha$-Hölder continuous, then there exists $\alpha' > 0$ such that the potential $\phi \circ \pi : \Sigma \to \mathbb{R}$ is $\alpha'$-Hölder continuous.*

Therefore, given such a potential $\phi : \cup_{Z \in \mathcal{P}_1} Z \to \mathbb{R}$, we can take the invariant Gibbs measure $\mu_{\Sigma_A^+, \phi}$ on $\Sigma_A^+$ for the associated symbolic version $\phi \circ \pi$ of $\phi$, then transport this to $\Lambda$ using $\pi$ to produce a measure $\mu_\phi = \mu_{\Sigma_A^+, \phi} \circ \pi^{-1}$. This is an $f$-invariant measure and which also has the Gibbs property.

Therefore, any HTS law for $(\Sigma_A^+, \sigma, \mu_{\Sigma_A^+, \phi})$ passes on to $(\Lambda, f, \mu_\phi)$. However, an important caveat is that the metrics on these two spaces usually do not match. In particular, not all balls in $\mathcal{X}$ are images of balls/cylinders in $\Sigma_A$. An approximation argument of balls by cylinders can sometimes be done, in particular for the following example.

Let $I := [0, 1]$ and let $f : \cup_{i=1}^n I_i \to I$ be a Markov map as above with Markov partition $\mathcal{P}_1 = \{I_1, \dots, I_n\}$ such that on each $Z \in \mathcal{P}_1$, $f$ is continuously differentiable and $\inf_{x \in Z} |Df(x)| > 1$. We call $(I, f)$ a *Markov interval map*.

**Theorem 5.4.3.** *Given a Markov interval map $(I, f)$, suppose that $\mu_\phi$ is a $\sigma$-invariant Gibbs measure for a Hölder potential $\phi : \cup_{Z \in \mathcal{P}_1} Z \to \mathbb{R}$ with Gibbs constant $P = 0$. Then, for $\zeta \in \Lambda$,*

$$
\lim_{\delta \to 0} \mu_\phi \left( \left\{ x \in I : r_{B_\delta(\zeta)}(x) > \frac{t}{\mu(B_\delta(\zeta))} \right\} \right) = e^{-\theta t},
$$

*where*

$$\theta = \begin{cases} 1 - e^{S_k \phi(\zeta)} & \textit{if } \zeta \textit{ is periodic of prime period } k, \\ 1 & \textit{otherwise} \end{cases}$$

This is Theorem 5.4.1, but adapted to balls using the results contained in [147, Section 5].

For non-Markov interval maps that are nevertheless uniformly expanding on a finite number of intervals, under certain conditions for the Gibbs-like measure, Paccaut proved the same result as above in [181].

### 5.4.4
### Two-sided Shifts

For a hyperbolic diffeomorphism, if there is a shift system that codes it, then the system is invertible. Hence, we have the motivation to extend the results above to invertible shift systems. As before, we start with an alphabet $S = \{1, \ldots, n\}$ and $A$ an $n \times n$ matrix of 0's and 1's. Then

$$\Sigma_A := \{(\ldots, x_{-1}, x_0, x_1, \ldots) : a_{x_i, x_{i+1}} = 1 \text{ for each } i \in \mathbb{Z}\}.$$

Therefore, the corresponding *two-sided subshift of finite type* is the pair $(\Sigma_A, \sigma)$ where $\sigma$ is the again the left-shift map

In this case the cylinder structure is given as follows. Given $a, b \in \mathbb{N}_0$, for any allowable word $(x_{-a}, x_{-a+1}, \ldots, x_{b-1}, x_b) \in \Sigma_A^*$, the corresponding *cylinder* is the set

$$\{y = (\ldots, y_{-1}, y_0, y_1 \ldots) \in \Sigma_A : y_i = x_i \text{ for } i = -a, \ldots, b\}.$$

If $a = b$, we call this a (two-sided) $k$-cylinder. Given $x \in \Sigma_A$, let $Z_k[x]$ denote the two-sided $k$-cylinder containing $x$.

It is a standard construction, see [174, Section 3] and [176, Lemma 1.6], to show that for any Gibbs measure on $(\Sigma_A, \sigma)$, which is defined naturally through two-sided $k$-cylinders, there is a Gibbs measure on $(\Sigma_A^+, \sigma)$ with precisely the same ergodic properties. In particular, this can be exploited as in [147] to prove a version of Theorem 5.4.3 for two-sided cylinders.

### 5.4.5
### Hyperbolic Diffeomorphisms

Suppose now that we have a diffeomorphism $f$ of a Riemannian manifold $\mathcal{X}$, with metric $d$, with the following hyperbolic structure on an an $f$-invariant set $\Lambda$. Set $\lambda^s < 0 < \lambda^u$ and suppose $C > 0$. For each $x \in \Lambda$ we require disks $W^s(x), W^u(x)$, called the *local stable manifold* and the *local unstable manifold* with the following properties:

1) $f(W^s(x)) \subset W^s(f(x))$ and $f^{-1}(W^u(x)) \subset W^u(f^{-1}(x))$;
2) For $n \in \mathbb{N}$,

$$d(f^n(x), f^n(y)) \leq C e^{\lambda^s n} d(x, y) \text{ for } y \in W^s(x),$$

$$d(f^{-n}(x), f^{-n}(y)) \leq C e^{-\lambda^u n} d(x, y) \text{ for } y \in W^u(x),$$

3) There exists $\varepsilon > 0$ such that for the sets $W^s_\varepsilon(x) := B_\varepsilon(x) \cap W^s(x)$ and $W^u_\varepsilon(x) := B_\varepsilon(x) \cap W^u(x)$, for each $x, y \in \Lambda$, the intersection $W^s_\varepsilon(x) \cap W^u_\varepsilon(y)$ consists of at most one point.

The numbers $\lambda^s$, $\lambda^u$ can be thought of as the exponential rates at which the systems contracts/expands on the stable/unstable directions. Note that in general, more structure is required, but we will only give a sketch of the fundamental properties. Given a set $A \subset \Lambda$, for $\varepsilon$ as above, let $W^{s,u}(x, A) := A \cap W^{s,u}_\varepsilon(x)$. This structure can give rise to a Markov partitions by 'rectangles': $R \subset \Lambda$ is a *rectangle* if

1) $R = \overline{\text{int}(R)}$;
2) $x, y \in R$ implies that $W^s(x, R) \cap W^u(y, R)$ is exactly one point.

Here closures and interiors are taken relative to $\Lambda$. Hence, in the context of a diffeomorphism $f : \mathcal{X} \to \mathcal{X}$ as above, with an invariant set $\Lambda$, a Markov partition for $f$ is a collection of rectangles $\{R_1, \ldots, R_n\}$ such that

1) $\Lambda = \cup_{i=1}^n R_i$;
2) If $x \in \text{int}(R_i)$ and $f(x) \in \text{int}(R_j)$, then

$$f(W^s(x, R_i)) \subset W^s(f(x), R_j) \text{ and } f(W^u(x, R_i)) \supset W^u(f(x), R_j).$$

3) If $x \in \text{int}(R_i) \cap f^{-1}(\text{int}(R_j))$, then

$$\text{int}(R_j) \cap f(W^u(x, \text{int}(R_i))) = W^u(f(x), \text{int}(R_j))$$

and

$$\text{int}(R_i) \cap f^{-1}(W^s(f(x), \text{int}(R_j))) = W^s(x, \text{int}(R_i)).$$

Under these conditions, see, *e.g.*, [182, Chapter 8], there is a topological conjugacy between $(\Lambda, f)$ and a SFT $(\Sigma_A, \sigma)$ where $A$ is determined by the third part of the definition above. This means that Theorem 5.4.1 holds, but for two-sided $k$-cylinders defined at points in $\Lambda$ rather than balls. The approximation results of balls by cylinders, performed in [147, 146] means that the result also extends to balls. The application given in [147] is to Axiom A diffeomorphisms, *e.g.*, the so-called Arnold cat map given by $\binom{x_{n+1}}{y_{n+1}} = \binom{2\ 1}{1\ 1} \binom{x_n}{y_n}$ on the 2-torus, since it is well-known that these have the requisite Markov partition, see, *e.g.* [176, 182]. We will present some results on this map in Chapters 6 and 9. The measures considered here are (two-sided) Gibbs measures on $\Lambda$.

The reader is encouraged to look into additional relevant investigations performed using different methods on the same type of systems [150, 151].

### 5.4.6
### Additional Uniformly Hyperbolic Examples

We close this section by making some brief comments on other uniformly hyperbolic dynamical systems for which the HTS/EVLs have been amenable to analysis. The first natural example is the class of tent maps. Given $a \in (1, 2]$, we define the *tent map* $T_a : [0, 1] \to [0, 1]$ by

$$T_a = \begin{cases} ax & \text{if } x \in [0, 1/2], \\ a(1-x) & \text{if } x \in (1/2, 1]. \end{cases}$$

This map has an absolutely continuous invariant measure $\mu_a$, and it is natural to ask whether we can prove that $\mu_a$-typical points have exponential HTS to balls around the point. This can be shown in various ways, one of which is described for more general interval maps in Sect. 5.5.3. since these maps fit into the framework of [160]. We note that, while this class of maps is uniformly expanding, not every map in the family has a finite Markov partition, so any proof via coding would be more problematic than the examples seen so far in this chapter.

The class of tent maps have many of the same topological characteristics as the logistic family we discuss in Section 5.5.3. Similarly, Hénon maps, which we describe below, share many topological properties with the *Lozi* family: a member of which is a homeomorphism of $\mathbb{R}^2$ with real parameters $a$ and $b$ defined by $f_{a,b}(x, y) = (1 + y - a|x|, bx)$. For a certain set of pairs $(a, b)$, this mapping has an invariant Sinai-Ruelle-Bowen (SRB) measure. It is shown in [138] that almost every point has exponential HTS to balls centered around the point. Again, there needs not be a Markov partition here.

We recall that an SRB measure is an invariant measure of a dynamical system that is additionally characterised by having absolutely continuous conditional measure on unstable manifolds. As well known, it is possible to provide other equivalent definitions of the SRB measure, see [70, 161].

Another class of uniformly hyperbolic systems are certain billiard systems. We refer the reader to [138] for results on these systems. One of the features of these systems is that they can be viewed as a suspension flow over a discrete dynamical system. In [139] it was shown how to go between the corresponding discrete system and its (semi-)flow: this was applied going from billiard maps to flows in [138].

### 5.5
### Non-uniformly Hyperbolic Systems

So far we have focussed in this chapter on systems which are uniformly hyperbolic, either uniformly expanding or with a hyperbolic splitting. Here we will focus on systems which only exhibit non-uniform hyperbolicity. In particular this means that these systems do not have nice Gibbs properties. (A more in-depth description of properties one might like good systems to have is given in Chapter 6.) These systems can be analysed directly, or by looking at a speeded-up version: by 'inducing' to produce a uniformly hyperbolic system. So before beginning to consider specific systems, we outline the inducing technique. In the specific examples given below, we will usually discuss absolutely continuous invariant measures, and their properties, but these ideas also extend to other classes of measures.

### 5.5.1
### Induced System

Given a dynamical system $f : \mathcal{X} \to \mathcal{X}$ with an $f$-invariant probability measure $\mu$, suppose that we are interested in the HTS to a point $\zeta \in \mathcal{X}$. If, *e.g.*, if our system is not uniformly hyperbolic, it can be useful to choose some subset $\hat{\mathcal{X}} \subset \mathcal{X}$ which also contains $\zeta$ and consider the *first return map* $\hat{f} : \hat{\mathcal{X}} \to \hat{\mathcal{X}}$, where $\hat{f} = f^{R_{\hat{\mathcal{X}}}}$ and $R_{\hat{\mathcal{X}}} : \hat{\mathcal{X}} \to \mathbb{N} \cup \{\infty\}$ is the *first return time* to $\hat{\mathcal{X}}$, *i.e.*, for $x \in \hat{\mathcal{X}}$,

$$R_{\hat{\mathcal{X}}}(x) = \inf\{n \in \mathbb{N} : f^n(x) \in \hat{\mathcal{X}}\}.$$

In particular, if $x$ never returns to $\hat{\mathcal{X}}$, then $R_{\hat{\mathcal{X}}}(x) = \infty$. We denote the *conditional measure* on $\hat{\mathcal{X}}$ by $\hat{\mu} = \mu_{\hat{\mathcal{X}}}$, so for $A \subset \mathcal{X}$, $\hat{\mu}(A) = \frac{\mu(A \cap \hat{\mathcal{X}})}{\mu(\hat{\mathcal{X}})}$ as in Eq. 5.1.1. By Kac's Theorem, $\hat{\mu}$ is $\hat{f}$-invariant. The basic idea here is that the HTS of this induced system can pass to the HTS of the original system. So if $(\hat{\mathcal{X}}, \hat{f}, \hat{\mu})$ is a system with well-understood HTS, *e.g.*, a uniformly hyperbolic system, then the HTS of the original system is also understood. We next give more details of this approach.

In order to discuss the HTS for the induced system, we need to consider the hitting time to a set $A$ with respect to $\hat{f}$, which we denote by $\hat{r}_A$. Suppose that $(A_n)_n$ is a sequence of sets shrinking down to $\zeta$. It is shown in [183] that for a distribution $H \in \mathcal{F}$,

$$\hat{\mu}(A_n)\hat{r}_{A_n} \stackrel{\hat{\mu}}{\Longrightarrow} H \quad \text{if and only if} \quad \mu(A_n)r_{A_n} \stackrel{\mu}{\Longrightarrow} H.$$

Note that the proof of this theorem for typical points $\zeta$ can be found in [184], while the proof for periodic points is given in [140]. The basic idea of the proof is that since $\hat{f}$ is a speeded up version of $f$, the hitting times should be speeded up too, by some average scaling factor related to the first returns to $\hat{\mathcal{X}}$. Indeed, by Kac's Theorem,

$$\lim_{j \to \infty} \frac{1}{j} \sum_{k=0}^{j-1} R_{\hat{\mathcal{X}}} = \frac{1}{\mu(\hat{\mathcal{X}})},$$

and because the important returns to $A_n$ take longer and longer as $n \to \infty$ (if $\zeta$ is periodic, then the periodic behaviour could be dealt with separately), $r_{\hat{A}_n}$ asymptotically behaves like $\frac{r_{A_n}}{\mu(\hat{\mathcal{X}})}$. Note also that by definition $\frac{\mu(A_n)}{\hat{\mu}(A_n)} = \mu(\hat{\mathcal{X}})$ if $A_n \subset \hat{\mathcal{X}}$.

### 5.5.2
### Intermittent Maps

Perhaps, the simplest non-uniformly hyperbolic systems are the map that are topologically equivalent the standard doubling map, $x \mapsto 2x \mod 1$ on $[0, 1)$, but are non-uniformly expanding at the fixed point. Such maps are usually referred to as Manneville-Pomeau maps [185]; though, we use the version considered in [186]. For $\alpha > 0$, set

$$g_\alpha(x) = \begin{cases} x(1 + 2^\alpha x^\alpha) & \text{if } x \in [0, 1/2), \\ 2x - 1 & \text{if } x \in [1/2, 1). \end{cases} \tag{5.5.1}$$

For all $\alpha > 0$, $g_\alpha$ has an acim $\mu_\alpha$, but this is only a probability measure, an absolutely continuous invariant measure, when $\alpha \in (0, 1)$. In none of these cases does $\mu_\alpha$ satisfy a Gibbs property. It is interesting to note that some of the earliest studies of HTS in this context were done for the case $\alpha \geq 1$, see [166, 187], although these were laws for cylinders where the partition $\mathcal{P}_1$ was not the canonical one $\mathcal{P}_1 = \{[0, 1/2), [1/2, 1)\}$, but rather one tailored to suit the specific problem. Namely, the partition elements were given by $A_n := [a_{n+1}, a_n)$ where $a_0 = 1$, $a_1 = 1/2$ and $g_\alpha(a_{n+1}) = a_n$ for $n \in \mathbb{N}$. Since then $g_\alpha(A_0) = [0, 1)$ and $g_\alpha(A_n) = A_{n-1}$ for $n \in \mathbb{N}$, the system is topologically the renewal shift: *i.e.*, the sequence space $\{(x_0, x_1, \ldots) : x_i \in \mathbb{N}$ follows $x_{i-1}$ if either $x_i = x_{i-1} - 1$ or $x_i = 1\}$, equipped with the left shift. Therefore, we refer to this as the *renewal partition*.

A significant step in the proof of HTS laws in general, and for non-uniformly hyperbolic systems in particular, comes from the framework devised in [188]. This gives useful tests on a system which, if satisfied, guarantees an exponential HTS law. One of the applications is to prove a HTS law for cylinders in the Manneville-Pomeau setting for $\mu_\alpha$ when $\alpha \in (0, 1)$. However, the cylinders come from the renewal partition.

In [184], these results have been extended to balls, as well as other kinds of sets, using the inducing technique described above in Section 5.5.1. In particular, $\hat{\mathcal{X}} = [1/2, 1)$ and $\hat{f}$ is the first return map of $g_\alpha$ to $\hat{\mathcal{X}}$. In order to prove that this induced system had exponential HTS to balls around almost-every point, they show that this induced map fits into the class of Rychlik maps defined in [160], which in particular have good mixing properties, such as the exponential mixing for BV against $L^1$. In concert with [188], this has allowed them to prove that for $\hat{\mu}$-a.e. $\zeta \in \hat{\mathcal{X}}$, there is exponential HTS to balls around $\zeta$. Then they have been able to pull this information back to the original system, proving the same result for $([0, 1), g_\alpha, \mu_\alpha)$. Later in [140], these ideas have been extended to periodic points, ultimately proving the following.

**Theorem 5.5.1.** *For $\alpha \in (0, 1)$, consider the Manneville-Pomeau map $([0, 1), g_\alpha, \mu_\alpha)$. Then for $\mu_\alpha$-a.e. $\zeta$,*

$$\lim_{\delta \to 0} \mu_\alpha \left( \left\{ x \in [0, 1) : r_{B_\delta(\zeta)}(x) > \frac{t}{\mu(B_\delta(\zeta))} \right\} \right) = e^{-t}.$$

*Moreover, if $\zeta$ is a periodic point of prime period $p$, then*

$$\lim_{\delta \to 0} \mu_\alpha \left( \left\{ x \in [0,1) : r_{B_\delta(\zeta)}(x) > \frac{t}{\mu(B_\delta(\zeta))} \right\} \right) = e^{-\theta t},$$

*for $\theta = 1 - \frac{1}{|Dg_\alpha^p(\zeta)|}$.*

### 5.5.3
### Interval Maps with Critical Points

The two key properties of the Manneville-Pomeau maps, in the context of HTS, are that despite being non-uniformly hyperbolic, they have a Markov structure, and that first return maps fit into the class given in [160]. In order to extend understanding of HTS beyond such settings, we consider the logistic family $f_\lambda : [0,1] \to [0,1]$ such that $f_\lambda(x) = \lambda x(1-x)$. We are most interested when the parameter $\lambda$ lies in $(3, 4]$. In that range of parameters there are maps $f_\lambda$ with a wide variety of behaviours. By convention, we denote the critical point $1/2$ by $c$. Notice that the only interesting dynamics is in the set $[f_\lambda^2(c), f_\lambda(c)]$, the *dynamical core*, since any other point, except 0 and 1 must map inside this interval eventually. Some of these maps are Markov (this is the case where the orbit of the critical point is finite), in which case the analysis is similar to that for the Manneville-Pomeau system. In some cases, the critical point gets attracted to an attracting cycle, *i.e.*, a cycle $x, f_\lambda(x), \ldots, f_\lambda^{n-1}(x)$ where $f_\lambda^n(x) = x$ and $|Df^n(x)| < 1$, the *Axiom A* case: here Lebesgue-a.e. point in $[0,1]$ is also attracted to that cycle and there is no absolutely continuous invariant measure. The set of $\lambda$ for which $f_\lambda$ is Axiom A is open and dense in $[0,4]$ ([189, 190]). We focus on the case of the positive Lebesgue measure set of parameters $\lambda$ such that there *is* an absolutely continuous invariant measure for $f_\lambda$ [191]. Note that while we specialise to the logistic family, the results discussed here extend to smooth multimodal maps and beyond.

If the critical point of $f_\lambda$ is non-recurrent, *i.e.*, there is a neighbourhood $U$ of $c$ such that the forward orbit of $c$ never enters $U$, then $f_\lambda$ is called a *Misiurewicz map*. Misiurewicz showed in [192] that such maps have an absolutely continuous invariant measure $\mu$. One of the nicest properties of these maps is that one can find an interval $X \subset [0,1]$ such that the first return map to $X$ is a full-branched mapping satisfying the conditions of [160], so using [184, 140] we can show an analogue of Theorem 5.5.1 for $(I, f_\lambda, \mu)$. Indeed, this type of result, for typical points appeared in [184].

The most interesting maps in the logistic family are the $f_\lambda$ whose critical orbit is dense in the dynamical core. In this case, first return maps to a given interval $X$ are not generally well-behaved. The best studied, and best-behaved, maps of this form are the *Collet-Eckmann* maps: there exist $C, \alpha > 0$ such that

$$|Df_\lambda^n(f(c))| \geq C e^{\alpha n}. \tag{CE}$$

It is shown in [193] that for such maps there exists an absolutely continuous invariant measure. In the language of HTS, Collet proved the following in [72] using a Young tower method.

**Theorem 5.5.2.** *Suppose that $\lambda$ is such that $f_\lambda$ satisfies (CE). Then for $\mu$ an absolutely continuous invariant measure, for $\mu$-a.e. $\zeta$,*

$$\lim_{\delta \to 0} \mu\left(\left\{x \in [0,1) : r_{B_\delta(\zeta)}(x) > \frac{t}{\mu(B_\delta(\zeta))}\right\}\right) = e^{-t}.$$

When the growth of $|Df_\lambda^n(f(c))|$ is subexponential, [194] uses a Hofbauer extension technique, see also [184], to prove an analogue of Theorem 5.5.2: for example, if $|Df_\lambda^n(f(c))|$ grows faster than $n^{1+\delta}$ for some $\delta > 0$ then the same result holds. This has been subsequently improved in [195], where it has been shown that all that is needed to derive the desired result is *only* the existence of an absolutely continuous invariant measure. Instead, properties such as growth along the critical orbit, and the related mixing behaviour do not play any role. The periodic case here has been studied in [140]. Unfortunately, in order to get good information on the density of the absolutely continuous invariant measure at the periodic points, it has been necessary to make assumptions on the map $f_\lambda$ stronger than (CE). As a result, the authors obtained the following result. Recall here that if $\zeta$ is a repelling periodic cycle of period $p$ for some $f_\lambda$, then for all $\lambda'$ close enough to $\lambda$, there is an equivalent point $\zeta_{\lambda'}$, the *hyperbolic continuation* of $\zeta$, which is a repelling periodic cycle of period $p$ for $f_{\lambda'}$.

**Theorem 5.5.3.** *Suppose that $\zeta$ is a periodic point of period $p$ for $f_4$. Then there exists a positive Lebesgue measure set $\Gamma_\zeta$ such that for the hyperbolic continuation $\zeta_\lambda$ of $\zeta$ for maps $f_\lambda$ where $\lambda \in \Gamma_\zeta$,*

$$\lim_{\delta \to 0} \mu\left(\left\{x \in [0,1) : r_{B_\delta(\zeta_\lambda)}(x) > \frac{t}{\mu(B_\delta(\zeta_\lambda))}\right\}\right) = e^{-\theta t}.$$

*Here $\mu$ is an absolutely continuous invariant measure and $\theta = 1 - \frac{1}{|Df_\lambda^p(\zeta_\lambda)|}$.*

Our catalogue of non-uniformly hyperbolic interval maps for which it is possible to prove HTS laws is not exhaustive. As an example, we do not discuss Lorenz-like interval maps considered in [139], where, moreover, these results have been extended to the relevant Lorenz flow.

### 5.5.4
### Higher Dimensional Examples of Non-uniform Hyperbolic Systems

To finish this section we make some short remarks, with very few details, on higher-dimensional non-uniformly hyperbolic systems. In [74] a non-uniformly hyperbolic system on the 2-torus was considered. This was non-uniformly expanding and it was possible to prove exponential HTS around almost every point w.r.t. the SRB measure.

In [150], a complete analogue of Theorem 5.4.3 was shown for partially hyperbolic diffeomorphisms w.r.t. the '$u$-Gibbs states' there. More recently, in [163], Hénon maps were considered (recall the simple model of Lozi maps above). Given real parameters $a, b \in \mathbb{R}$, the Hénon map is the homeomorphism $f_{a,b} : \mathbb{R}^2 \to \mathbb{R}^2$ defined by $f_{a,b}(x, y) = (1 - ax^2 + y, bx)$. In [196] it is shown that there is a positive Lebesgue measure set of parameters $(a, b) \in \mathbb{R}^2$ such that $f_{a,b}$ has an SRB measure $\mu_{a,b}$ with good mixing properties (one should think of $a$ being close to 2 and $b$ being small here). In [163] it is shown that for $\mu_{a,b}$-typical points $\zeta$, we have exponential HTS to balls around $\zeta$. We do not give any details at this stage, but we remark that, in contrast to the uniformly hyperbolic toral automorphisms, the structure of the stable and unstable manifolds here is very complicated due to the criticality at $(0, 0)$. There is certainly no simple Markov structure to work with here.

## 5.6
## Non-exponential Laws

We finish this chapter by commenting on some systems where the limiting laws are not of an exponential type. We have already mentioned that any law in $\mathcal{F}$, defined in Section 5.2 can be realised as the HTS for a certain sequence of sets. However, these sequences are arguably unnatural in general. However in [197], and later in [198], some very natural cases were considered, and we wish to give here a quick sense of these results.

Consider the circle rotation $R_\alpha : S^1 \to S^1$ with rotation number $\alpha$, induced from the transformation $x \mapsto x + \alpha$ on $\mathbb{R}$ where $S^1$ is identified with $\mathbb{R}/\mathbb{Z}$. Denote by $\mu$ the Haar measure on $S^1$. Given $\zeta \in S^1$, consider the interval/arc $J_\varepsilon := [\zeta, \zeta + \varepsilon] \subset S^1$. Note that is only interesting to consider $\alpha$ irrational here.

**Theorem 5.6.1.** *([198, Theorem 1]) For any irrational rotation number $\alpha$, the distribution $\mu(J_\varepsilon) r_{J_\varepsilon}$ does not converge as $\varepsilon \to 0$.*

This striking result is then followed by laws along special sequences of sets, which complement further such findings [197]. For example, for certain $\alpha$ (of unbounded type in the sense of the continued fraction expansion), one can choose a $c \in (0, 1)$ and an explicit sequence $\varepsilon_n \to 0$ as $n \to \infty$ such that the limiting distribution for $\mu(J_{\varepsilon_n}) r_{J_{\varepsilon_n}}$ is

$$F(t) = \begin{cases} t & \text{if } 0 \le t \le c, \\ c & \text{if } t > c. \end{cases}$$

One of the features of the system $(S^1, R_\alpha, \mu)$ is that it has zero entropy. For more work in this direction see [199, 200]. Conversely, in [201] it is shown that there is a zero-entropy system which, nevertheless, has exponential HTS along all sequences of cylinders at a typical point.

# 6
# Extreme Value Theory for Selected Dynamical Systems

## 6.1
## Rare Events and Dynamical Systems

As it can be seen in Section 5.3, the theory of HTS/RTS and the theory of EVL for dynamical systems are two sides of the same coin. This means that we can prove the existence of EVLs by proving HTS and the other way around.

The theory of HTS/RTS laws is now a well developed topic, applied first to cylinders and hyperbolic dynamics, and then extended to balls and also to non-uniformly hyperbolic systems. We refer to [172] and [202] for very nice reviews as well as many references on the subject. (See also [173], where the focus is more towards a finer analysis of uniformly hyperbolic systems.) Since the early papers [164, 147], several different approaches have been used to prove HTS/RTS: from the analysis of adapted Perron-Frobenius operators as in [147], the use of inducing schemes as in [184], to the relation between recurrence rates and dimension as explained in [202, Section 4].

For many mixing systems it is known that the HTS/RTS are standard exponential around almost every point. Among these systems we note the following: Markov chains [164], Axiom A diffeomorphisms [147], uniformly expanding maps of the interval [203], 1-dimensional non-uniformly expanding maps [188, 184, 194, 195], partially hyperbolic dynamical systems [150], toral automorphisms [151], higher dimensional non-uniformly hyperbolic systems (including Hénon maps) [163].

In most of the papers mentioned so far, exponential HTS and RTS have been proved for generic points, in the sense that there exists exponential HTS/RTS for almost all $\zeta$ in the phase space, with respect to the invariant measure. However, in [147] and [129], the authors consider specific points. Namely, they consider the cases when $\zeta$ is a periodic point and obtain the existence of an EI less than 1, although they did not stated in these terms because, at the time, the connection with EVL was not yet established. In fact, Galves and Schmitt [165] introduced a short correction factor $\lambda$ in order to get exponential HTS, that was then studied later in great detail by Abadi *et al.* [173, 204, 205, 206, 152, 180], and which, in case of being convergent, can actually be seen as the EI itself.

On the other hand, EVLs for the partial maximum of dynamically defined stochas-tic processes is a much recent topic and have been proved directly in the recent papers [72, 134, 74, 137, 138, 139, 49, 141, 140, 142, 207, 144, 143]. We highlight the pioneering work of Collet [72] for the innovative ideas introduced. The dynamical systems covered in these papers include non-uniformly hyperbolic 1-dimensional maps (in all of them), higher dimensional non-uniformly expanding maps in [74], suspension flows in [139], billiards and Lozi maps in [138].

The purpose of this chapter is to present the state of the art regarding the existence of EVL for specific dynamical systems. We compile a list of systems for which the existence of EVL have been proved and describe the main techniques used to obtain such results. Furthermore, we present some results elucidating the rate of convergence of BM to the asymptotic EVL.

## 6.2
## Introduction and Background on Extremes in Dynamical Systems

In this chapter we consider a dynamical system $(\mathcal{X}, \mathcal{B}, \mu, f)$, where $\mathcal{X} \subset \mathbb{R}^d$ is a compact Riemannian manifold, $\mathcal{B}$ is its Borel $\sigma$-algebra, $f : \mathcal{X} \to \mathcal{X}$ a measurable transformation, and $\mu$ is an $f$-invariant probability measure supported on $\mathcal{X}$. Given an observable $\phi : \mathcal{X} \to \mathbb{R}$ we consider the stationary stochastic process $X_1, X_2, \ldots$ defined as in (4.2.2) and its associated maximum process $M_n$ defined by (2.2.1) For the stochastic process defined in (4.2.2), much recent work has focused on the computation of the function $G$ appearing on the right hand side of (3.1.3), and showing that it exists and agrees, at least for sufficiently hyperbolic systems and for regular enough observations $\phi$ maximized at generic $\tilde{x}$, with that which would hold if $\{X_i\}$ were i.i.d. random variables. If $\tilde{x}$ is periodic we expect different behavior (for details see [49, 146, 147, 141]).

In this chapter we focus on a dynamical blocking argument approach that leads to an estimation of $\mu\{M_n \leq u_n\}$ in terms of the distribution types described above. This approach was first adopted in [72], and then applied more recently in [163, 138, 144, 139, 208, 137]. We begin this chapter by reviewing this blocking argument approach. In large part it is purely probabilistic and in applications we just require knowledge of mixing rates (decay of correlations), and quantitative recurrence statistics. We contrast the approach to the verification of conditions $D(u_n), D'(u_n)$ described in Section 3.2. Using this approach we then prove extreme value laws for observations maximized at generic points for a wide class of chaotic dynamical systems. We include a section on non-uniformly expanding systems [72, 139], a section on non-uniformly hyperbolic systems [138, 163, 208], and a section on partially hyperbolic systems [137]. The blocking argument approach also leads to estimates of rates on convergence to an EVL, and we describe this following [144].

Throughout we fix the following notations. For general sequences $x_n, y_n$, we say that $x_n \sim y_n$ if $x_n/y_n \to 1$ as $n \to \infty$. We say $x_n \approx y_n$ if there are real constants $c_1, c_2$ such that $c_1 \leq x_n/y_n \leq c_2$. For positive sequences we write $x_n = \mathcal{O}(y_n)$ if there is a constant $C > 0$ such that $x_n \leq Cy_n$, and we write $x_n = o(y_n)$ if $x_n/y_n \to 0$.

## 6.3
## The Blocking Argument for Non-uniformly Expanding Systems

We adapt the blocking argument used by Leadbetter in the classical setting described in Section 3.2.1 to a dynamical setting. This adjustments are consistent with the arguments in Section 4.1 and illustrate how the information about decay of correlations and quantitative recurrence of the system come to play and allow to show the existence of EVL for typical points in the absence of clusters. Moreover, these computations allow the reader to follow the successive estimates in order to obtain convergence rates.

We take $(f, \mathcal{X}, \mu)$ to be an ergodic dynamical system. Let $\tilde{g} : \mathbb{N} \to \mathbb{R}$ be a monotonically increasing function and let $E_n$ be a sequence of sets defined by:

$$E_n := \left\{ x \in \mathcal{X} : \operatorname{dist}(x, f^j(x)) \le \frac{1}{n}, \text{ for some } j \in [1, \tilde{g}(n)] \right\}. \qquad (6.3.1)$$

In [208], the sets $E_n$ are referred to as *recurrence sets*. If the time scale $\tilde{g}(n) = o(n)$, these sets contain points which have what we call 'fast' returns. We will see that the distributional rate of convergence to an EVL is partly determined by the rate at which $\mu(E_n)$ converges to zero.

**Assumptions on the invariant measure $\mu$.**
We will assume that the measure $\mu$ is absolutely continuous with respect to Lebesgue $m$. Some of our results require also that $\mu$ admits a density $\rho \in L^{1+\delta}(m)$ for some $\delta > 0$. For non-uniformly hyperbolic systems in dimension at least two, we require further assumptions on the regularity of $\mu$. We will discuss this further in Section 6.4.

**Dynamical assumptions on $(f, \mathcal{X}, \mu)$**
We list the following assumptions. The function $\tilde{g}(n)$ and the sets $E_n$ are from Eq. (6.3.1).

(H1) **(Decay of correlations)** There exists a monotonically decreasing sequence $\Theta(j) \to 0$ such that for all $\varphi_1$ Lipschitz continuous and all $\varphi_2 \in L^\infty$:

$$\left| \int \varphi_1 \cdot \varphi_2 \circ f^j d\mu - \int \varphi_1 d\mu \int \varphi_2 d\mu \right| \le \Theta(j) \|\varphi_1\|_{\text{Lip}} \|\varphi_2\|_{L^\infty},$$

where $\| \cdot \|_{\text{Lip}}$ denotes the Lipschitz norm. This is our decay assumption for non-uniformly expanding maps.

(H1s) **(Decay of correlations)** There exists a monotonically decreasing sequence $\Theta(j) \to 0$ such that for all Lipschitz $\varphi_1$ and $\varphi_2$:

$$\left| \int \varphi_1 \cdot \varphi_2 \circ f^j d\mu - \int \varphi_1 d\mu \int \varphi_2 d\mu \right| \le \Theta(j) \|\varphi_1\|_{\text{Lip}} \|\varphi_2\|_{\text{Lip}},$$

where $\| \cdot \|_{\text{Lip}}$ denotes the Lipschitz norm. This is our decay assumption for non-uniformly hyperbolic systems with a stable direction.

(H2a) **(Strong quantitative recurrence rates).** There exist numbers $\gamma, \alpha > 0$ such that

$$\tilde{g}(n) \sim n^{\gamma} \quad \implies \quad \mu(E_n) \leq \frac{C}{n^{\alpha}}. \tag{6.3.2}$$

Condition (H1) is an assumption on the rate of mixing in a suitable Banach space of functions. The correlation decay in the Banach spaces of Lipschitz versus $L^{\infty}$ has been shown to hold for our applications to non-uniformly expanding maps. In Section 4, where we detail results on hyperbolic systems with stable directions, we will assume a correlation decay in Lipschitz versus Lipschitz as given in (H1s). In the statement of the results we will make precise asymptotic statements about the rate of decay of $\Theta(j)$.

Condition (H2a) is a quantitative control on the recurrence statistics. For general non-uniformly expanding systems, checking this condition requires careful analysis. See for example [137, 139]. In Section 6.4 we show how to check this condition for uniformly expanding Markov maps and certain Markov intermittency maps. For systems having sub-exponential decay of correlations we specifically require (H2a) to hold in order to derive error rates. It is conjectured in [208] that (H2a) holds for a broad class of non-uniformly hyperbolic systems, and this is observed numerically (when analytic estimates are not available). For systems with exponential decay of correlations we can work with a weaker version of (H2a), which we label as (H2b). This is stated as follows:

(H2b) **(Weak quantitative recurrence rates).** For some $\gamma' > 1, \alpha > 0$:

$$\tilde{g}(n) \sim (\log n)^{\gamma'} \quad \implies \quad \mu(E_n) \leq \frac{C}{n^{\alpha}}. \tag{6.3.3}$$

This latter condition only requires control of the recurrence up to a slow time scale $\tilde{g}(n) \sim (\log n)^{\gamma'}$ and is easier to check analytically relative to (H2a), see for example [72, 138, 208]. To obtain convergence to an EVL (and to find corresponding error rates), it is sufficient to check this condition provided the system has exponential decay of correlations.

To establish convergence rates for certain non-uniformly hyperbolic systems (without absolutely continuous invariant measures) we will work with conditions analagous to (H2a)/(H2b). Such conditions will be discussed in Section 6.5.

### Assumption on the observable type

We will assume that the observable $\phi : \mathcal{X} \to \mathbb{R}$ is a distance observable, so that it can be written in the form $\phi(x) = \psi(\text{dist}(x, \tilde{x}))$, where $\psi : \mathbb{R}^{+} \to \mathbb{R}$ has a *unique* maximum at 0 and $\tilde{x}$ has a density $\rho(\tilde{x}) = \frac{d\mu}{dm}(x)$. Hence $\phi$ is maximized at the unique point $\tilde{x} \in \mathcal{X}$. In earlier chapters, the dependence of the distributional limit $G$ on the specific functional form of $\psi$ was considered. Within this section we will focus on the case where $\psi(y) = -\log y$, and thus we will be in the domain of attraction of a Type I distribution at $\mu$-a.e. $\tilde{x}$, provided the density of $\mu$ is absolutely continuous with respect to Lebesgue measure $m$. In the case of observable types having multiple maxima, distributional limit theorems can also be derived, see [139].

### Statement or Results

**Theorem 6.3.1.** *Suppose that $f : \mathcal{X} \to \mathcal{X}$ is a map with ergodic measure in $L^{1+\delta}(m)$ for some $\delta > 0$, and $\mu$ absolutely continuous with respect to $m$. We have the following cases.*

1) *Suppose that $\Theta(n) = O(\theta_0^n)$ for some $\theta_0 < 1$ and (H1), (H2b) hold. Then for all $\epsilon > 0$ and $\mu$-a.e. $\tilde{x} \in \mathcal{X}$ we have that*

$$\left| \mu\{M_n \leq u_n\} - \left(1 - \sqrt{n}\mu\{X > u_n\}\right)^{\sqrt{n}} \right| \leq C_1 \frac{(\log n)^{1+\epsilon}}{\sqrt{n}} + \frac{C_2}{n^{\alpha-\epsilon}}, \quad (6.3.4)$$

*where $C_1$, $C_2 > 0$ are constants independent of $n$, but dependent on $\tilde{x}$.*

2) *Suppose that $\Theta(n) = O(n^{-\zeta})$ for some $\zeta > 0$ and (H1), (H2a) hold. Then for all $\epsilon > 0$ and $\mu$-a.e. $\tilde{x} \in \mathcal{X}$ we have that*

$$\left| \mu\{M_n \leq u_n\} - \left(1 - \sqrt{n}\mu\{X > u_n\}\right)^{\sqrt{n}} \right| \leq C_1 n^{-\frac{1}{2}+\kappa} + C_2 n^{-\alpha+\kappa},$$

$$\text{with } \kappa = \epsilon + \frac{2(1+2\delta)}{\zeta\delta}. \quad (6.3.5)$$

*where $C_1$, $C_2 > 0$ are constants independent of $n$, but dependent on $\tilde{x}$.*

As we indicate in Section 6.8 the conclusions of Theorem 6.3.1 provide quantitative estimates on the convergence rates to an extreme value distribution. We state the following corollary which indicates that we do indeed get convergence to a distribution of given type. The type depends on the regularity of the measure and the precise form of the distance observable function $\phi(x)$. We state the following corollary in the case $\phi(x) = -\log \text{dist}(x, \tilde{x})$.

**Corollary 6.3.2.** *Suppose that $f : \mathcal{X} \to \mathcal{X}$ is a map with ergodic measure in $L^{1+\delta}(m)$ for some $\delta > 0$, and $\mu$ absolutely continuous with respect to $m$. Suppose that i) $\Theta(n) = O(\theta_0^n)$ for some $\theta_0 < 1$ and (H1), (H2b) hold, or ii) $\Theta(n) = O(n^{-\zeta})$ for some $\zeta > 0$ and (H1), (H2a) hold. Consider the observable function $\phi(x) = -\log \text{dist}(x, \tilde{x})$ and the associated process $M_n$. Then for $\mu$-a.e. $\tilde{x} \in \mathcal{X}$ we have:*

$$\lim_{n \to \infty} \mu\{M_n \leq u + \log n\} = \exp\{-2\rho(\tilde{x})e^{-u}\}. \quad (6.3.6)$$

### 6.3.1
### The Blocking Argument in One Dimension

To simplify our exposition we will suppose that our system is one-dimensional, the modifications for $d > 1$ are obvious. We suppose that there are scaling sequences $a_n, b_n$ such that Eq. (3.1.3) is valid. In particular for the distance observable type $\phi(x) = -\log(\text{dist}(x, \tilde{x}))$ in dimension one, we may take $u_n = v + \log n$, and

$$\lim_{n \to \infty} n\mu\{x \in \mathcal{X} : \phi(x) > u_n\} = \lim_{n \to \infty} n\mu\{\text{dist}(x, \tilde{x}) \leq \frac{e^{-v}}{n}\}$$
$$\to 2\rho(\tilde{x})e^{-v}. \quad (6.3.7)$$

### 6.3.2
### Quantification of the Error Rates

The following propositions give precise quantification of the error rates in terms of the assumptions on the correlation decay $\Theta(j)$, and the decay of $\mu(E_n)$. They will be used in the proof of Theorem 6.3.1. To state the propositions, we fix integers $p(n), q(n) > 0$ and let $n = pq + r$ with $0 \leq r < p$ (by Euclid's division algorithm). The blocking argument consists of writing $n = p(n)q(n)$ and between each of the $q$ gaps of length $p$ we take a gap of length $t = g(n)$. The decay of correlations over the gap of length $t = g(n)$ allows us to consider successive blocks as approximately independent. We suppose that $p, q \to \infty$ as $n \to \infty$. We let $u_n$ be the sequence with the property that $n\mu\{X_1 > u_n\} \to \tau(u)$, for some function $\tau(u)$, and $u_n = u/a_n + b_n$.

**Proposition 6.3.3.** *Suppose that $f : \mathcal{X} \to \mathcal{X}$ is ergodic with respect to a measure $\mu$ which has a density $\rho \in L^{1+\delta}(m)$ for some $\delta > 0$. Suppose that (H1) holds. Then for $\mu$-a.e. $\tilde{x} \in \mathcal{X}$, all $p, q$ such that $n = pq + r$, and $t < p$, we have*

$$|\mu\{M_n \leq u_n\} - (1 - p\mu\{X_1 > u_n\})^q| \leq \mathcal{E}_n, \qquad (6.3.8)$$

*where for any $\epsilon > 0$ and $\delta_1 = \delta/(1 + 2\delta)$:*

$$\mathcal{E}_n = \max\{qt, p\}\mu\{X_1 \geq u_n\} + qp^2(\mu\{X_1 > u_n\})^2$$
$$+ C_1 pq^2 \Theta(t)^{\delta_1 - \epsilon} + pq \sum_{j=2}^{t} \mu(X_1 > u_n, X_j > u_n). \qquad (6.3.9)$$

**Proposition 6.3.4.** *Suppose that $f : \mathcal{X} \to \mathcal{X}$ is ergodic with respect to a measure $\mu$ with density $\rho \in L^{1+\delta}(m)$ for some $\delta > 0$. Suppose that (H2a) or (H2b) hold, and suppose for given $\epsilon > 0$ that $g(n) = \tilde{g}(n)^{1-\epsilon}$. Then for all $\epsilon > 0$, and $\mu$-a.e. $\tilde{x} \in \mathcal{X}$:*

$$\sum_{j=2}^{g(n)} \mu(X_1 > u_n, X_j > u_n) \leq C(\tilde{x})\frac{g(n)}{n^{\tilde{\alpha}}}, \quad \text{with } \tilde{\alpha} = \alpha + 1 - \epsilon_1 \qquad (6.3.10)$$

*where the constant $C$ depends on $\tilde{x}$, and $\alpha$ is as defined in (H2a) or (H2b).*

#### 6.3.2.1  Proof of Proposition 6.3.3
To prove Proposition 6.3.3 we need the following result which is purely probabilistic, see [72, Proposition 3.2].

**Lemma 6.3.5.** *For any integers $t, r, m, k, p \geq 0$*

$$0 \leq \mu(M_r < u) - \mu(M_{r+k} < u) \leq k\mu(X_1 \geq u),$$

*and*

$$\left| \mu(M_{m+p+t} < u) - \mu(M_m < u) + \sum_{j=1}^{p} \mathbb{E}\left(\mathbf{1}_{\{X_1 \geq u\}} \mathbf{1}_{\{M_m \circ f^{p+t-j} < u\}}\right) \right|$$

$$\leq 2p \sum_{j=1}^{p} \mathbb{E}\left(\mathbf{1}_{\{X_1 \geq u\}} \mathbf{1}_{\{X_j < u\}}\right) + t\mu(X_1 > u).$$

*Proof of Proposition* 6.3.3. We will apply the estimates of Lemma 6.3.5. Note first of all that these estimates are expressed in terms as expectations of products $\mathbf{1}_{\{X_1 \geq u\}} \mathbf{1}_{\{X_j < u\}}$. To use decay of correlation estimates, we must approximate the indicator function $\mathbf{1}_{\{X > u_n\}}$ by a Lipschitz function. This is done by approximating the indicator function $\Phi := \mathbf{1}_{\{X_1 > u_n\}}$ by a Lipschitz continuous function $\Phi_B$, which is set equal to 1 inside a ball centered at $\tilde{x}$ of radius $\ell'_n := \psi^{-1}(u_n) - [\psi^{-1}(u_n)]^{1+\eta}$, for some $\eta > 0$, and decaying to 0 at a linear rate so that $\Phi_B$ vanishes on the boundary of $\{X_1 > u_n\}$. The Lipschitz norm of $\Phi_B$ is bounded by $[\psi^{-1}(u_n)]^{-(1+\eta)}$. For any measurable set $A \subset \mathbb{R}^d$ this leads to the estimate

$$|\mu(\mathbf{1}_{\{X > u_n\}} \cap f^{-t}(A)) - \mu(\mathbf{1}_{\{X > u_n\}})\mu(A)| \leq m(\mathbf{1}_{\{X > u_n\}})^{-\frac{(1+\eta)}{d}}\Theta(t)$$
$$+ \mathcal{O}(1)m(\mathbf{1}_{\{X > u_n\}})^{\theta(1+\eta)} \quad (6.3.11)$$

where $\Theta$ is the decay of Lipschitz functions agains the $L^\infty$-norm,

$$\left| \int \phi \cdot \psi \circ f^n d\mu - \int \phi d\mu \int \psi d\mu \right| \leq \Theta(n)\|\psi\|_{L^\infty}\|\phi\|_{\text{Lip}}.$$

This is obtained by approximating $\mathbf{1}_I$ with a piecewise-linear function.

Putting this together gives:

$$\left| \mu(M_n < u_n) - \mu(M_{q(p+t)} < u_n) \right| \leq \max\{qt, p\}\mu(X_1 \geq u_n),$$

and for $1 \leq \ell \leq q$ we have

$$|\mu(M_{\ell(p+t)} < u_n) - (1 - p\mu(\psi > u_n))\mu(M_{(\ell-1)(p+t)} < u_n)|$$

$$\leq \left| p\mu(X_1 \geq u_n)\mu(M_{(\ell-1)(p+t)} < u_n) - \sum_{j=1}^{p} \mathbb{E}\left(\chi_{\{\psi \circ f^j \geq u_n\}}\chi_{M_{(\ell-1)(p+t)} \circ f^{p+t} < u_n}\right) \right|$$

$$+ t\mu(X_1 \geq u_n) + 2p\sum_{j=1}^{p} \mu\left(\{X_1 \geq u_n\} \cap \{X_j \geq u_n\}\right).$$

This latter expression is bounded above by $\Gamma_n$ by the estimates above. $\qquad \square$

*Proof of Corollary* 6.3.2. We see that if $p\mu(\psi \geq u_n) < 2$ then an iterative argument shows that

$$|\mu(M_{q(p+t)} < u_n) - (1 - p\mu(\psi \geq u_n)^q| \leq q\Gamma_n$$

and so

$$|\mu(M_n < u_n) - (1 - p\mu(\psi \geq u_n)^q)| \leq q\Gamma_n + \max\{qt, p\}\mu(X_1 \geq u_n)$$

by the results above.

Hence if $q\Gamma_n \to 0$ and $qt\mu(\psi < u_n) \to 0$ then

$$\mu(\max\{X_1, \ldots, X_n\} \leq u_n) \to \exp(-\lim_n n\mu(X_1 \geq u_n)).$$

$\square$

We state and prove the following lemma. Combining this result with the blocking argument described earlier in Section 6.3.2 completes the proof of Proposition 6.3.3.

**Lemma 6.3.6.** *For any $g(n) < p$ we have that:*

$$\sum_{j=g(n)}^{p} \mu(X_0 > u_n, X_j > u_n) \leq pO(1)(\Theta(g(n)))^{\delta_1} + p(\mu\{X_1 > u_n\})^2,$$

$$\text{with } \delta_1 = \frac{\delta}{(1 + 2\delta)} - \epsilon, \tag{6.3.12}$$

*where $\epsilon > 0$ can be taken arbitrarily small, and the implied constant $O(1)$ depends on $\delta$.*

Under assumptions (H2a) or (H2b), and Propostion 6.3.4 a condition on the choice of $g(n)$ is that $g(n) = o(\tilde{g}(n))$. We now use decay of correlations to prove the second part of (6.3.12). We will write $\phi(x) = \psi(\text{dist}(x, \tilde{x}))$, and recall that we work with the explicit observable $\psi(y) = -\log y$, for $y > 0$. As before we approximate the indicator function $\Phi := 1_{\{X_1 > u_n\}}$ by a Lipschitz continuous function $\Phi_B$, which is set equal to 1 inside a ball centered at $\tilde{x}$ of radius $\ell'_n := \phi^{-1}(u_n) - [\psi^{-1}(u_n)]^{1+\eta}$, for some $\eta > 0$, and decaying to 0 at a linear rate so that $\Phi_B$ vanishes on the boundary of $\{X_1 > u_n\}$. The Lipschitz norm of $\Phi_B$ is bounded by $[\psi^{-1}(u_n)]^{-(1+\eta)}$.

We now take $j \in [g(n), n]$. We have the following triangle inequality:

$$\left| \int \Phi(\Phi \circ f^j)d\mu - \left(\int \Phi d\mu\right)^2 \right| \leq \left| \int \Phi_B(\Phi \circ f^j)d\mu - \int \Phi_B d\mu \int \Phi d\mu \right|$$

$$+ \left| \int (\Phi_B - \Phi)\Phi \circ f^j d\mu - \int (\Phi - \Phi_B)d\mu \int \Phi d\mu \right|,$$

and we estimate each term on the right hand side. By decay of correlations and for sufficiently large $n$:

$$\left| \int \Phi_B(\Phi \circ f^j)d\mu - \int \Phi_B d\mu \int \Phi d\mu \right| \leq \|\Phi_B\|_{\text{Lip}}\|\Phi\|_\infty \Theta(j)$$

$$= O\left(n^{1+\eta}\Theta(g(n))\right),$$

and

$$\left| \int (\Phi_B - \Phi)\Phi \circ f^j d\mu - \int (\Phi - \Phi_B) d\mu \int \Phi d\mu \right| \leq 2\|\Phi\|_\infty \mu(x : \Phi_B(x) \neq \Phi(x))$$
$$= O(n^{-\theta(1+\eta)}),$$

where we can take any $\theta < \delta/(1+\delta)$. Hence for each $j > g(n)$ we obtain:

$$\mu(X_1 > u_n, X_j > u_n) \leq C_1 (m\{X_1 \geq u_n\})^{-1-\eta} \Theta(j)$$
$$+ C_2 (m\{X_1 \geq u_n\})^{\theta(1+\eta)} + (\mu\{X_1 > u_n\})^2. \quad (6.3.13)$$

The constants $C_1, C_2$ depend on the regularity of $\rho(x)$ at $\tilde{x}$ and on the Lipschitz norm. The constant $\eta$ is arbitrary, and hence we can optimize the right hand side by varying $\eta$. If we let $x = (m\{X_1 \geq u_n\})^{1+\eta}$, then we can write the first two terms on the right hand side of inequality (6.3.13) as:

$$A(x) := C_1 x^{-1} \Theta(j) + C_2 x^\theta. \quad (6.3.14)$$

The minimizer of $A(x)$ is the value $x = O(1)\Theta(j)^{\frac{1}{\theta+1}}$, and leads to the bound:

$$\mu(X_1 > u_k, X_j > u_k) \leq O(1)(\Theta(j))^{\delta_1} + (\mu\{X_1 > u_n\})^2,$$
$$\text{with } \delta_1 = \frac{\delta}{(1+2\delta)} - \epsilon. \quad (6.3.15)$$

Here $\epsilon > 0$ can be made arbitrarily small. The implied constant depends only on $\delta$. This gives the required result.

### 6.3.2.2 Proof of Propositon 6.3.4: The maximal function argument

For a function $\varphi \in L^1(m)$ we define the Hardy–Littlewood maximal function

$$\mathcal{M}(x) := \sup_{a>0} \frac{1}{2a} \int_{x-a}^{x+a} \varphi(y) dm(y).$$

A theorem of Hardy and Littlewood [209], implies that

$$m(|\mathcal{M}(x)| > \lambda) \leq \frac{\|\varphi\|_1}{\lambda} \quad (6.3.16)$$

where $\| \cdot \|_1$ is the $L^1$ norm with respect to $m$. Recalling

$$E_n = \left\{ x : \text{dist}(x, f^j(x)) \leq \frac{1}{n} \text{ for some } j \leq \tilde{g}(n) \right\},$$

let $\rho(x)$ denote the density of $\mu$ with respect to $m$ and let $\mathcal{M}_n(x)$ denote the maximal function of $\varphi_n(x) := 1_{E_n}(x)\rho(x)$. For constants $a, b > 0$ to be fixed later consider sequences $\lambda_n = n^{-a}$ and $\alpha_n = \lfloor n^b \rfloor$. Inequality (6.3.16) gives

$$m(|\mathcal{M}_{\alpha_n}(x)| > \lambda_n) \leq \frac{\mu(E_{\alpha_n})}{\lambda_n} \leq \frac{1}{n^{\alpha b - a}}.$$

If $\alpha b - a > 1$ (first constraint required on $a$ and $b$), then the First Borel–Cantelli Lemma implies for $\mu$ a.e. $x$ there exists an $N := N(x)$ such that for all $n \geq N$ we have $|\mathcal{M}_{\alpha_n}(\tilde{x})| < \lambda_n$. Recall in the case of Theorem 6.3.1 that $\mu$ is absolutely continuous with respect to $m$. Hence, for all $n$ sufficiently large

$$
\begin{aligned}
\mu\big(\{x : \operatorname{dist}(x, \tilde{x}) < \alpha_n^{-1}\} \cap E_{\alpha_n}\big) &\leq \int_{\tilde{x} - \alpha_n^{-1}}^{\tilde{x} + \alpha_n^{-1}} \varphi_{\alpha_n}(y) dm(y) \\
&\leq 2\alpha_n^{-1} \mathcal{M}_{\alpha_n}(\tilde{x}) \\
&\leq 2\alpha_n^{-1} \lambda_n = O(n^{-a-b}).
\end{aligned}
\tag{6.3.17}
$$

Denote $A := \{X_1 > u_k, X_j > u_k\}$ with $2 \leq j \leq g(k)$, and $g(n) = \tilde{g}(n)^{(1-\epsilon)}$ for some $\epsilon > 0$. We assume that $\tilde{g}(n)$ has the representations given in either (H2a) or (H2b). For distance observables of the form $\phi(x) = \psi(\operatorname{dist}(x, \tilde{x}))$, we have $\psi^{-1}(u_k) \approx 1/k$. Hence there exists a $v > 0$ such that

$$
A \subset \left\{ x : \operatorname{dist}(\tilde{x}, x) \leq \frac{v}{k} \operatorname{dist}(\tilde{x}, f^j(x)) \leq \frac{v}{k} \text{ for some } j \leq g(k) \right\}.
$$

Given the sequence $\alpha_n$, let $k/(2v) \in [\alpha_n, \alpha_{n+1})$. Then (by monotonicity of $g(n)$),

$$
A \subset \left\{ x : \operatorname{dist}(\tilde{x}, x) \leq \frac{1}{2\alpha_n}, \ \operatorname{dist}(\tilde{x}, f^j(x)) \leq \frac{1}{2\alpha_n} \text{ for some } j \leq g((2v)\alpha_{n+1}) \right\}.
$$

Applying the triangle inequality $\operatorname{dist}(x, f^j(x)) \leq \operatorname{dist}(\tilde{x}, x) + \operatorname{dist}(\tilde{x}, f^j(x))$ gives

$$
A \subset \left\{ x : \operatorname{dist}(\tilde{x}, x) \leq \frac{1}{\alpha_n}, \ \operatorname{dist}(x, f^j(x)) \leq \frac{1}{\alpha_n} \text{ for some } j \leq g((2v)\alpha_{n+1}) \right\}.
$$

Since $\alpha_n = \lfloor n^b \rfloor$ we have that $\lim_{n \to \infty} \alpha_{n+1}/\alpha_n = 1$. By the growth properties of $g$ and $\tilde{g}$ (as given in Proposition 6.3.4), there exists $\kappa_v > 0$ and a sequence $c_n \to 0$ such that for all sufficiently large $\alpha_n$:

$$
g((2v)\alpha_{n+1}) \leq g(2(2v)\alpha_n) \leq c_n \tilde{g}((2(2v)\alpha_n) \leq c_n \kappa_v \tilde{g}(\alpha_n).
$$

Moreover, there exists $N$ such that $\forall n \geq N$ we have $c_n \kappa_v < 1$, and hence

$$
A \subset \left\{ x : \operatorname{dist}(\tilde{x}, x) \leq \frac{1}{\alpha_n}, \ \operatorname{dist}(x, f^j(x)) \leq \frac{1}{\alpha_n} \text{ for some } j \leq \tilde{g}(\alpha_n) \right\}.
$$

Applying inequality (6.3.17) gives

$$
\mu(X_1 > u_k, X_j > u_k) = O(k^{-1 - \frac{a}{b}}) \quad \text{for all } k > N,
$$

so that

$$
\sum_{j=1}^{g(k)} \mu(X_1 > u_k, X_j > u_k) = O(k^{-1 - \frac{a}{b}} g(k)) \quad \text{for all } k > N.
$$

To complete the proof of Proposition 6.3.4 we now choose optimal values of $a$ and $b$ (positive) to minimize $-1 - a/b$, subject to the constraint $\alpha b - a > 1$. This gives for all $k > N$:

$$\sum_{j=1}^{g(k)} \mu(X_1 > u_k, X_j > u_k) = O(k^{-1-\frac{a\alpha}{a+1}} g(k)) = O(k^{-1-\alpha+\epsilon} g(k)),$$

valid for all $\epsilon > 0$. Hence Eq. (6.3.10) is satisfied.

### 6.3.3
### Proof of Theorem 6.3.1

Before considering specific dynamical systems, we show how Theorem 6.3.1 follows from Propositions 6.3.3 and 6.3.4. The proof requires optimizing the choice of constants $q, p$ and the gap length $t < p$ which appear in the division algori theorem $n \sim q(n)p(n)$ of the blocking argument, as well taking into account the decay of correlations and regularity of $\mu$.

Given $\beta \in (0, 1)$ we will suppose first that $p \sim n^{1-\beta}$ and so $q \sim n^\beta$. The constant $r$ in $n = pq + r$ satisfies $r < p$ and hence $r = O(n^{1-\beta})$. Assuming (H1) and either (H2a) or (H2b), we immediately have from Propositions 6.3.3 and 6.3.4 that:

$$\mathcal{E}_n \leq C_1 \max\{n^{\beta-1}g(n), n^{-\beta}\} + C_2 n^{-\beta} + C_3 n^{2-\beta}\{\Theta(g(n))\}^{\delta_1-\epsilon} + C_4 \frac{g(n)}{n^{\alpha-\epsilon}}, (6.3.18)$$

where $\delta_1 = \delta/(1 + 2\delta)$, and $\epsilon > 0$ arbitrary. In this estimate we take $t = g(n)$ and used the fact that

$$\mu\{X_1 > u_n\} = \frac{\tau(u)}{n} + o\left(\frac{1}{n}\right), \tag{6.3.19}$$

Thus the constants $C_1$ and $C_2$ depend on the functional form of $\tau(u)$. In turn this behaves on the local behaviour of the invariant density at $\tilde{x}$ and on the functional form of the observable $\phi$. The constant $C_3$ depends on $\delta$ (and hence the regularity of $\mu$). The constant $C_4$ is from Proposition 6.3.4 and depends on the recurrence properties associated to $\tilde{x}$.

Let us now prove part 1 Theorem 6.3.1. In this case $\Theta(n) \leq O(\theta_0^n)$, and under assumption (H2), we can take $t = (\log n)^\kappa$ for some $\kappa > 1$. For this choice it follows that $\Theta(t) \to 0$ at a super-polynomial rate. By assumption the measure $\mu$ has density in $L^{1+\delta}$ for some $\delta > 0$, and hence for any choice of $p, q = o(n)$ we have

$$C_3 n^{2-\beta}\{\Theta(t)\}^{\frac{\delta}{1+2\delta}-\epsilon} = o(1/n).$$

Inspecting the first two terms on the right hand side of Eq. (6.3.9), gives an optimal choice $p \approx q \approx \sqrt{n}$. For this choice of $p$ and $q$, and noting that $n = pq + r$, we conclude the proof of part 1 of Theorem 6.3.1 via the following sequence of

estimates:

$$(1 - p\mu\{X_1 > u_n\})^q = \left(1 - \frac{(n-r)}{\sqrt{n}}\mu\{X_1 > u_n\}\right)^{\sqrt{n}},$$

$$= \left(1 - \sqrt{n}\mu\{X_1 > u_n\}\right)^{\sqrt{n}} + O\left(\frac{1}{\sqrt{n}}\right),$$

where in the last line we have used Proposition 6.8.1 in conjunction with the fact that $r\mu\{X_1 > u_n\} \leq O(1/\sqrt{n})$. The convergence is pointwise in $u$.

To prove part 2 of Theorem 6.3.1, we see from (H2a) and Proposition 6.3.4 that in the estimation of $\mathcal{E}_n$ we can take $t = n^\kappa$ for any $\kappa < \gamma$. We inspect each term of the error $\mathcal{E}_n$ in Eq. (6.3.9). Again we take $\beta = 1/2$. The contribution of the first error term on the right hand side of (6.3.9) gives a contribution of order $n^{-1/2+\kappa}$. The next significant error term is now the third right hand term of (6.3.9). Putting in $t = n^\kappa$ gives an error contribution:

$$C_3 n^{2-\beta}\{\Theta(n^\kappa)\}^{\delta_1} = C_3 n^{3/2} n^{-\zeta\kappa\delta_1}, \quad \text{with } \delta_1 = \frac{\delta}{1+2\delta} - \epsilon,$$

for any $\epsilon > 0$. In fact we require this error to be order $n^{-1/2}$, and hence this gives a bound on $\kappa$. We obtain for any $\epsilon_1 > 0$:

$$\kappa \geq \frac{2(1+2\delta)}{\delta\zeta} + \epsilon_1, \tag{6.3.20}$$

and thus we just take the minimal value of $\kappa$. This gives the required conclusion to part 2 of Theorem 6.3.1.

This gives the required conclusion to part 2 of Theorem 6.3.1.

## 6.4
## Non-uniformly Expanding Dynamical Systems

We consider examples of non-uniformly expanding dynamical systems that fit assumptions (H1) and (H2), and hence whose extreme statistics can be understood via Theorem 6.3.1. These systems will admit an invariant measure that is absolutely continuous with respect to the ambient Riemannian measure $m$. We will consider the following examples: uniformly expanding maps, quadratic maps [210], expanding Lorenz maps [211], and intermittency maps with subexponential mixing rates [156]. Extreme value laws have been proved for these systems on a case by case basis.

### 6.4.1
### Uniformly Expanding Maps

We derive explicit convergence rates for the tent map $f(x) = 1 - |1 - 2x|$ on $[0, 1]$. Let $E_n^{(j)} := \{x : \text{dist}(x, f^j(x)) < 1/n\}$, and suppose $I$ is a monotonicity sub-interval of $f^j$ and let $J = I \cap E_n^{(j)}$. Since $f^j(I) = [0, 1]$ and $f^j$ has slope $2^j$, it follows easily that $|J| = O(2^{-j}/n)$. Hence, summing over all such $J$, we have $\mu(E_n^{(j)}) = O(1/n)$ and hence $\mu(E_n) = O(\tilde{g}(n)/n)$. To optimize $\mathcal{E}_n$, we can take a functional form $\tilde{g}(n) = (\log n)^{\gamma'}$ for any $\gamma' > 1$. By exponential decay of correlations, $\Theta(g(n))$ tends to zero at a superpolynomial rate. Here $g(n) = (\log n)^\kappa$ for some $1 < \kappa < \gamma'$. Thus conditions (H1) and (H2b) are valid, in particular the latter for $\alpha > 1/2$. We summarize as follows:

**Proposition 6.4.1.** *Suppose $(f, [0, 1], \text{Leb})$ is the tent map. For the observation $\phi(x) = -\log|x - \tilde{x}|$, we have for Leb-a.e $\tilde{x} \in [0, 1]$ and all $\epsilon > 0$:*

$$\lim_{n \to \infty} \mu\{M_n \leq u + \log n\} = e^{-2e^{-u}} \tag{6.4.1}$$

*where $C(\tilde{x}) > 0$ is a uniform constant depending on $\tilde{x}$.*

*Remark* 6.4.2. The main technical step required to achieve Eq. (6.4.1) is in the estimation of $\mu(E_n)$. For a wide class of uniformly expanding Markov maps, such as those considered in [139, Section 3.3] we expect the same error estimate to apply.

### 6.4.2
### Non-uniformly Expanding Quadratic Maps

Consider the quadratic family $f(x) = a - x^2$ for $x \in [-2, 2]$ and parameter $a \simeq 2$. For a positive measure set of parameter values, it is known that $f$ admits an absolutely continuous invariant measure $\mu$ with density $\rho \in L^{1+\delta}$ for some $\delta > 0$. Moreover the system admits exponential decay of correlations, see [157]. Here $E_n$ had explicit representation:

$$E_n = \{x \in [0, 1] : \text{dist}(f^j x, x) < \frac{1}{n}, \text{ some } j \leq (\log n)^5\}. \tag{6.4.2}$$

It was shown in [72] that $\mu(E_n) \leq n^{-\alpha}$ for some $\alpha > 0$. Hence conditions (H1) and (H2b) are satisfied for this family, and we can take $g(n) = (\log n)^{1+\kappa}$ for some $\kappa > 0$.

**Proposition 6.4.3.** *Suppose $(f, [-2, 2], \mu)$ is the quadratic family with a a Benedicks-Carlesion parameter. For the observation $\phi(x) = -\log|x - \tilde{x}|$, we have for $\mu$-a.e $\tilde{x} \in [0, 1]$ and all $\epsilon > 0$:*

$$\lim_{n \to \infty} \mu\{M_n \leq u + \log n\} = e^{-2\rho(\tilde{x})e^{-u}}, \tag{6.4.3}$$

*where $C(\tilde{x}) > 0$ is a uniform constant independent of $n$, but dependent on $\tilde{x}$.*

We remark that for observables maximized at specific points $\tilde{x} \in \mathcal{X}$, such as values of the critical orbit, EVLs where established in [134]. The proof used specific combinatorics of the critical orbit to verify $D(u_n)$ and $D'(u_n)$.

### 6.4.3
### One-dimensional Lorenz Maps

In this section we consider extreme statistics for a family of uniformly expanding maps with singularities. In [138] the family of Lorenz maps where considered, but the methods easily extend to other expanding maps with discontinuities or singularities. The Lorenz family of maps $f : [-1, 1] \setminus \{0\} \to [-1, 1]$, are defined as follows:

(L1) There exists $C > 0$ and $\lambda > 1$ such that for all $x \in I$ and $n > 0$, $|(f^n)'(x)| > C\lambda^n$.

(L2) There exists $\beta \in (0, 1)$ such that $f'x = |x|^{\beta-1} g(x)$ where $g \in C^{\beta\epsilon}(X)$, $g > 0$.

(L3) $f$ is locally eventually onto. *i.e.* For all intervals $J \subset X$, there exists $k = k(J) > 0$ such that $f^k(J) = X$.

These maps can be modelled by a Young tower with exponential tails. Hence there exists an absolutely continuous invariant measure $\mu$ and $(f, X, \mu)$ has exponential decay of correlations in the space of Hölder continuous functions. It is known that the the density of $\mu$ is of bounded variation. We will return again to the Lorenz family in Section 6.11.1, but in the context of non-uniformly expanding systems we have the following result:

**Proposition 6.4.4.** *Suppose $(f, [0, 1], \mu)$ is a Lorenz map satisying (L1), (L2) and (L3). Consider the observable $\phi(x) = -\log |x - \tilde{x}|$. Then for all $\epsilon > 0$ and $\mu$-a.e. $\tilde{x} \in \mathcal{X}$ we have that*

$$\lim_{n \to \infty} \mu\{M_n \leq u_n\} = e^{-2\rho(\tilde{x})e^{-u}}, \tag{6.4.4}$$

*where $C(\tilde{x}) > 0$ is a constant independent of $n$.*

### 6.4.4
### Non-uniformly Expanding Intermittency Maps.

Consider the following interval map defined for $b > 0$:

$$f(x) = \begin{cases} x(1 + (2x)^b) & \text{for } 0 \leq x < \frac{1}{2}, \\ 2x - 1 & \text{for } \frac{1}{2} \leq x \leq 1. \end{cases} \tag{6.4.5}$$

This map is non-uniformly expanding, and it has a neutral fixed point at $x = 0$. It was introduced in [186] as a simple model of intermittency and is sometimes called the Liverani-Saussol-Vaienti map. For $b \in (0, 1)$, the map admits an absolutely continuous invariant measure, where the density lies in $L^{1+\delta}$ for any $\delta < 1/b - 1$. The system has polynomial decay of correlations: $\Theta(n) = O\left(n^{1-1/b}\right)$, see [156]. We have the following result:

**Proposition 6.4.5.** *Suppose $(f, [0, 1], \mu)$ is the intermittent system (6.4.5) with $b < 1/20$ and consider the observable $\phi(x) = -\log |x - \tilde{x}|$. Then for all $\epsilon > 0$ and $\mu$-a.e. $\tilde{x} \in \mathcal{X}$ we have that*

$$\lim_{n \to \infty} \mu\{M_n \leq u_n\} = e^{-2\rho(\tilde{x})e^{-u}}, \tag{6.4.6}$$

*where $C(\tilde{x}) > 0$ is a constant independent of $n$.*

*Proof.* The following lemma, a version proved in [139] will be of use to us: Define $\mathcal{E}_n(\epsilon) := \{x : \text{dist}(x, f^n(x)) < \epsilon\}$.

**Lemma 6.4.6.** *Suppose $f$ is the interval map given by Eq. (6.4.5). There exists a uniform constant $C > 0$, such that $\forall n \geq 0$, and $\forall \epsilon > 0$ we have $m(\mathcal{E}_n(\epsilon)) \leq C\sqrt{\epsilon}$.*

*Proof.* The following proof is specific to the intermittent family. In [139], an extended result is proved for general non-uniformly expanding maps with convex slope. We prove Lemma 6.4.6 as follows. Let $\mathcal{P} = \{[0, 1/2], [1/2, 1]\}$, and $\mathcal{P}_n = \bigvee_{k=0}^{n-1} f^{-k}\mathcal{P}$. Here $\bigvee$ denotes the join of partitions. Elements of $\mathcal{P}$ consist of intervals of the form $J_k = [a_k, a_{k+1}]$ where $f^n : J_k \to [0, 1]$. In particular $f^n \mid J_k$ is a function with increasing derivative, with slope $(f^n)' \geq 1$. We take $f^n(a_k) = 0$, $f^n(a_{k+1}) = 1$, where it is understood on $J_k$, that $\lim_{x \to a_{k+1}} f^n(x) = 1$, while on $J_{k+1}$ we have $\lim_{x \to a_{k+1}} f^n(x) = 0$.

Suppose that there are points $x_k^\pm \in (a_k, a_{k+1})$ such that $f^n(x_k^\pm) = x_k^\pm \pm \epsilon$. If these points exist then they are unique, and in particular $\mathcal{E}_n(\epsilon) \cap J_k = [x_k^-, x_k^+]$. Let $x \geq x_k^-$ in $J_k$. Since $(f^n)'$ is increasing on $J_k$,

$$(f^n)'(x) \geq (f^n)'(x_k^-) \geq \frac{\int_{[a_k, x_k^-]} (f^n)(t)'\, dt}{x_k^- - a_k} = \frac{x_k^- - \epsilon}{x_k^- - a_k},$$

and therefore

$$(f^n)'(x) - 1 \geq \frac{a_k - \epsilon}{x_k^- - a_k} \geq \frac{a_k - \epsilon}{m(J_k)}.$$

By integrating each term, we obtain:

$$2\epsilon = \int_{[x_k^-, x_k^+]} \left[ (f^k)'(t) - 1 \right]\, dt \geq (x_k^+ - x_k^-) \frac{(a_k - \epsilon)}{m(J_k)},$$

and thus

$$m(\mathcal{E}_n(\epsilon) \cap J_k) \leq \frac{2\epsilon}{(a_k - \epsilon)} m(J_k). \tag{6.4.7}$$

This estimate is useful provided $(a_k - \epsilon)$ is not small. Given $\eta > 0$, let $k' = \sup\{k : a_k < \eta\}$, and let $W(\eta) = \cup_{k+1 \leq k'} J_k$. Then

$$m(\mathcal{E}_n(\epsilon)) = m(\mathcal{E}_n(\epsilon) \cap W(\eta)) + m(\mathcal{E}_n(\epsilon) \cap J_{k'}) + m(\mathcal{E}_n(\epsilon) \cap (W(\eta) \cup J_{k'})^c)$$

$$\leq m(W(\eta)) + m(\mathcal{E}_n(\epsilon) \cap J_{k'}) + \frac{2\epsilon}{\eta - \epsilon} m((W(\eta) \cup J_{k'})^c).$$

We estimate these sets in two different ways depending on whether $n$ is large or small. Following [139], the optimal choice turns out to be $\eta = \sqrt{\epsilon}$ and so we work with this value. This would give the desired result modulo the middle term (especially when $n$ is small). If $\text{diam}(\mathcal{P}_n) < \sqrt{\epsilon}$. Then for all $J_k \subset W(\sqrt{\epsilon})$ we have $|J_k| < \sqrt{\epsilon}$, and moreover $|J_{k'}| < \sqrt{\epsilon}$. Hence,

$$m(\mathcal{E}_n(\epsilon)) \leq C\sqrt{\epsilon}.$$

Note that the analysis above assumes that we can always solve for such a $x^{\pm}$. If we cannot, then the estimates are actually improved in the sense that either $\mathcal{E}_n(\epsilon) \cap J_k = \emptyset$ or $J_k$ partially crosses $\mathcal{E}_n(\epsilon)$ and the required measure would be smaller than that computed in Eq. (6.4.7).

Suppose now that $n$ is small with $\operatorname{diam}(\mathcal{P}_n) > \sqrt{\epsilon}$. By the explicit form of the intermittency map, the largest element of $\mathcal{P}_n$ is in fact $J_0 = [a_0, a_1]$, with $a_0 = 0$. Hence $a_1 > \sqrt{\epsilon}$. In this case $x^- = 0$, and on $J_0$ we have $f^n(x) - x \geq 2^b x^{b+1}$. Therefore $\mathcal{E}_n \cap J_0 \subset \{x : 2^b x^{b+1} < \epsilon\}$, and so $x^+ \leq O(\epsilon^{1/(1+b)}) = o(\sqrt{\epsilon})$. Therefore for $n$ small, we have

$$m(\mathcal{E}_n(\epsilon)) = m(\mathcal{E}_n(\epsilon) \cap J_0) + m(\mathcal{E}_n(\epsilon) \cap J_0^c)$$
$$\leq C\sqrt{\epsilon} + \frac{2\epsilon}{\sqrt{\epsilon} - \epsilon} m(J_0^c) \leq C\sqrt{\epsilon}.$$

Hence the conclusion of the Lemma follows. $\qquad \square$

It follows that $m(E_n) \leq C\tilde{g}(n) n^{-1/2}$, and hence by Hölder's inequality

$$\mu(E_n) \leq C(\tilde{g}(n) n^{-1/2})^{1-b}.$$

Therefore if we take $\tilde{g}(n) = n^{\gamma}$ for any $\gamma \in (0, 1/2)$, $\mu(E_n)$ tends to zero, and so (H2a) applies. In order to control the error term $\mathcal{E}_n$ in Proposition 6.3.3 we require that $\Theta(g(n)) \to 0$ sufficiently fast, for some $g(n) = n^{\kappa} = o(\tilde{g}(n))$. The latter choice being made also so that Proposition 6.3.4 applies with any $\alpha \leq (1/2 - \gamma)(1 - b)$. We can follow again the proof of part 2 of Theorem 6.3.1, and work out the minimal choice of $\kappa$. Recalling that $\delta < 1/b - 1$, and $\Theta(n) = O\left(n^{1-1/b}\right)$ we obtain using Eq. (6.3.20) the following bound on $\kappa$ (valid for all $\epsilon_1 > 0$):

$$\gamma > \kappa \geq \frac{2b(2-b)}{(1-b)^2} + \epsilon_1.$$

Notice for $b$ close to one, this bound is of little utility since we require $\kappa < 1/2$. However for $b < 1/20$, we can take $\kappa \geq 5b + \epsilon_1$. From Proposition 6.3.4 we also have the additional error term:

$$\sum_{j=2}^{g(n)} \mu(X_1 > u_n, X_j > u_n) \leq \frac{C(\tilde{x})}{n^{\tilde{\alpha}}}, \quad \text{with } \tilde{\alpha} = 1 + \alpha - \gamma + \epsilon \qquad (6.4.8)$$

with $\epsilon > 0$ arbitrary and $\alpha = (1/2 - \gamma)(1 - b)$. By choice of $\kappa$ and hence $\gamma$, this latter term gives a contribution of order $n^{-1/2 + \tilde{b}}$ with $\tilde{b} = 10b + \epsilon$.

$\qquad \square$

## 6.5
## Non-uniformly Hyperbolic Systems

In this section we discuss extreme laws and corresponding convergence rates for certain hyperbolic and non-uniformly hyperbolic dynamical systems. Examples include the Anosov Cat map, the Lozi map, hyperbolic billiard systems, and systems with rank one attractors such as the Hénon family. Extreme statistics for these systems have been investigated in [163, 138, 81, 208]. For these systems, rate estimates are possible in the spirit of Theorem 6.3.1 though there are complications in choosing the correct scaling $u_n := u/a_n + b_n$ as the invariant measure may fluctuate wildly. To establish rates of convergence (or even just convergence) to a limit $G(u)$ along the sequence $\mu\{M_n \leq u_n\}$ for a linear scaling $u/a_n + b_n$ information on the regularity of $\mu$ is required. In general such estimates on the regularity of $\mu$ are not established for non-uniformly hyperbolic systems. Ideally we would like to know the scaling $u_n$ for broad classes of physical systems without fine knowledge of the invariant measure, which will typically have a fractal structure.

We suppose $(f, \mathcal{X}, \mu)$ is a non-uniformly hyperbolic dynamical system, and that the invariant measure $\mu$ is a physical or SRB measure supported on the $f$-invariant set $\mathcal{X} \subset \mathbb{R}^n$. In most of our examples $n = 2$. In the late 1990's L.S-Young [157, 156] developed a method of inducing, subsequently called a Young tower construction, to determine the rate of decay of correlations for non-uniformly hyperbolic dynamical systems. This method revolutionized the field of dynamics and in particular solved a longstanding problem by proving that Sinai dispersing billiards have exponential decay of correlations for Hölder observables.

We will assume that $(f, \mathcal{X}, \mu)$ is modelled by a Young tower. In dimension one this is not an assumption as it has been shown that a system with a positive Lyapunov exponent and an absolutely continuous invariant measure has a Young tower [153].

Most of the dynamical systems we consider have exponential decay of correlations for Hölder observables, but in the statement of our results we will allow also for sub-exponential decay of correlations. We will assume that condition (H1s) as defined in Section 6.3 holds. For non-uniformly hyperbolic systems, we introduce the following *short return time* (SRT) conditions in place of conditions (H2a) and (H2b). Let $B(x, r)$ be the ball of radius $r$ centered at $x \in \mathcal{X}$.

(SRT1) Suppose $\{B(\tilde{x}, u_n)\}$ is a sequence of balls centered at $\tilde{x} \in \mathcal{X}$ with $\limsup n\mu(B(\tilde{x}, u_n)) < \infty$. Then there exists $\alpha > 0$, and $\gamma > 0$ such that for $\mu$-a.e. $\tilde{x} \in \mathcal{X}$:

$$\mu\left(B(\tilde{x}, u_n) \cap f^{-k} B(\tilde{x}, u_n)\right) \leq \mu\left((B(\tilde{x}, u_n))^{1+\alpha}\right) \tag{6.5.1}$$

for all $k = 1, \ldots \tilde{g}(n)$, with $\tilde{g}(n) = n^\gamma$.

(SRT2) As in (SRT1) but $\tilde{g}(n) = (\log n)^{\gamma'}$ for some $\gamma' > 1$.

The (SRT) conditions are easier to verify that (H2a) and (H2b) for systems where the SRB measure $\mu$ is not absolutely continuous with respect to Lebesgue measure. For example the Hénon family of maps possess an SRB measure that is not absolutely continuous with respect to two-dimensional Lebesgue measure. For systems with absolutely continuous invariant measures conditions (H2a) and (H2b) actually imply the (SRT) conditions for generic points $\tilde{x} \in \mathcal{X}$. This follows from the proof of Proposition 6.3.4. We make further comment on how (H2a) and (H2b) relate to the (SRT) conditions in Section 6.5.1.

In our analysis of extremes for non-uniformly hyperbolic systems, we require some control on the regularity of $\mu$. Recall that the pointwise local dimension of $\mu$ is given by:

$$d(x) := \lim_{r \to 0} \frac{\log \mu(B(x, r))}{\log r}, \qquad (6.5.2)$$

whenever this limit exists. For the examples we consider the local dimension of $\mu$ exists and is constant $\mu$-a.e. In addition we need to control the regularity of $\mu$ on certain shrinking annuli. We assume (H3):

(H3) **(Regularity of $\mu$ on shrinking annuli).** For all $\delta > 1$ and $\mu$-a.e. $x \in \mathcal{X}$, there exists $\sigma > 0$ such that

$$|\mu(B(x, r + r^\delta)) - \mu(B(x, r))| \leq C r^{\sigma \delta}. \qquad (6.5.3)$$

The constant $C$ and $\sigma$ depend on $x$ but are independent of $\delta$.

Condition (H3) (and versions thereof) have been stated and verified in [163, 138, 208] for the systems we consider here.

**Theorem 6.5.1.** *Suppose that $f : \mathcal{X} \to \mathcal{X}$ is a non-uniformly hyperbolic system with ergodic SRB measure $\mu$, and the the local dimension $d$ exists $\mu$-a.e. Assume that condition (H3) holds.*

1) *If we have exponential decay of correlations, i.e. $\Theta(n) = O(\theta_0^n)$ with $\theta_0 < 1$ and conditions (H1s) and (SRT2) hold, then for all $\epsilon > 0$ and $\mu$-a.e. $\tilde{x} \in \mathcal{X}$,*

$$\left| \mu\{M_n \leq u_n\} - G_{\sqrt{n}}(u) \right| \leq C_1 \frac{(\log n)^{1+\epsilon}}{\sqrt{n}} + C_2 \frac{(\log n)^{1+\epsilon}}{n^\alpha} \qquad (6.5.4)$$

*where $C_1$, $C_2 > 0$ are constants independent of $n$, but dependent on $\tilde{x}$.*

2) *If we have polynomial decay of correlations i.e. $\Theta(n) = O(n^{-\zeta})$ for some $\zeta > 1$ and conditions (H1) and (SRT1) hold, then for all $\epsilon > 0$ and $\mu$-a.e. $\tilde{x} \in \mathcal{X}$*

$$\left| \mu\{M_n \leq u_n\} - G_{\sqrt{n}}(u) \right| \leq C_1 n^{-\frac{1}{2} + \kappa} + C_2 n^{-\alpha + \kappa}, \quad \text{with } \kappa = \epsilon + \frac{C_\sigma}{\zeta}. \quad (6.5.5)$$

*where $C_1$ is independent of $n$. The constant $C_\sigma > 0$ is independent of $\zeta$, but depends on $\sigma$ in (H3).*

The proof of this theorem follows in large part from the proof of Propositions 6.3.3 and 6.3.4, where in the latter we use the (SRT) conditions in place of conditions (H2a) and (H2b). Within [163, 138, 208], a direct blocking argument approach is used to verify $D'(u_n), D(u_n)$ for hyperbolic systems. We point out in Section 6.5.1, how the proofs of these propositions are modified in the hyperbolic setting. Lozi maps, hyperbolic billiardss and the Hénon family satisfy conditions (H1s), (SRT) and (H3) with exponential decay of $\Theta(n)$ (see [163, 138, 208]). The constant $\alpha$ will depend on the mixing rates, the regularity of the invariant measure, and on geometrical properties of the map. In the present situation we do not attempt to optimize an estimate for $\alpha$, though with work it could be done for the systems we have mentioned. To show how to do so we shall estimate $\alpha$ for the Arnold cat map (discussed below), which is perhaps the simplest invertible uniformly hyperbolic system. This map has both a stable and unstable foliation. For higher dimensional hyperbolic systems with polynomial decay of correlations, less is known about convergence to EVL (*i.e.* examples that satisfy (H1s), (SRT) and (H3)), but some progress in this direction is made in [212].

For non-uniformly hyperbolic systems the estimates (6.5.4) and (6.5.5) are the best convergence rates that our techniques allow, at least for regular observables $\phi : \mathcal{X} \to \mathbb{R}$. The main difficulty is in the control of the fluctuations of $\mu\{\phi(x) > u_n\}$ as $n \to \infty$. From the definition of local dimension the function $r \mapsto \mu\{B(x, r)\}$ need not be regularly varying as $r \to 0$, and hence even for smooth observables such as $\phi(x) = -\log(\text{dist}(x, \tilde{x}))$, the sequence $\tau_n(u) = n\mu\{\phi(x) > u + \log n\}$ may fluctuate wildly as $n \to \infty$. Thus, in the statement of the (SRT) conditions we choose $u_n$ so that $\limsup \tau_n(u) < \infty$. However, along other (non-linear) scalings $u_n(u)$ of $u$ control on the rate of convergence of $\tau_n(u)$ to $\tau(u)$ might be achievable, though determining these scalings is impracticable. Alternatively, for observables $\phi : \mathcal{X} \to \mathbb{R}$ tailored to the measure $\mu$, so that $\mu\{\phi(x) > u_n\}$ is regularly varying in $u$, convergence rates can again be achieved. Such observables are considered in [136].

## 6.5.1
## Proof of Theorem 6.5.1

The blocking argument described in Section 6.3.2 is purely probabilistic. For dynamical systems, we require an extra hypothesis (H3) on the regularity of the measure $\mu$. For non-uniformly expanding systems, we used the fact that the density of $\mu$ belonged to $L^{1+\delta}$ for some $\delta > 0$, which implies condition (H3).

We point out the main modifications required over and above the details presented in the proofs of Propositions 6.3.3 and 6.3.4 in order to prove Theorem 6.5.1. For non-uniformly hyperbolic systems, we assume (SRT1) and (SRT2) instead of (H2a) and (H2b).

The scheme of proof of Proposition 6.3.4, under the assumptions (H2a) and (H2b) encounters a problem when using the maximal function approach, as the Borel-Cantelli argument used to establish Eq. (6.3.17) is no longer valid. It required $\mu$ to be absolutely continuous with respect to $m$ and this is no longer the case. If we assume intend the (SRT) conditions, then Eq. (6.3.10) in Proposition 6.3.4 becomes:

$$\sum_{j=2}^{g(n)} \mu(X_1 > u_n, X_j > u_n) = \sum_{j=1}^{g(n)-1} \mu\left(B(\tilde{x}, \psi^{-1}(u_n)) \cap f^{-j}B(\tilde{x}, \psi^{-1}(u_n))\right),$$

$$\leq \sum_{j=1}^{g(n)} \mu\left((B(\tilde{x}, \psi^{-1}(u_n)))^{1+\alpha}\right),$$

$$\leq g(n)n^{-1-\alpha},$$

where we recall that $\phi(x) = \psi(\text{dist}(x, \tilde{x}))$, and $u_n$ is chosen in such a way that $\limsup n\mu\{B(\tilde{x}, \psi^{-1}(u_n))\} < \infty$. This gives us the corresponding error term contributions as stated in Theorem 6.5.1.

Under the assumption of (H1s) and (H3), the modification of proof of Proposition 6.3.3 in the hyperbolic case again follows the proofs of [138, 208]. In particular condition (H3) comes into the Lipschitz approximation of $1_{\{X_1 > u_n\}}$, to estimate $\|\Phi(x) - \Phi_B(x)\|_1$. The result is that Eq. (6.3.15) is modified to

$$\mu(X_1 > u_k, X_j > u_k) \leq O(1)(\Theta(j))^{\delta'} + (\mu\{X_1 > u_n\})^2, \tag{6.5.6}$$

for some constant $\delta' > 0$. A lower bound on the constant $\delta'$ can be achieved by using the proofs within [138, Lemmas 3.1, 3.2]. In particular the bound depends on $\sigma$ (from (H3)), and also on the rates of contraction along stable manifolds. Theorem 6.5.1 now follows.

*Remark* 6.5.2. To check the (SRT) conditions in specific examples, arguments based upon Hardy Littlewood maximal functions are still applicable, see [138, 208]. However, the approach used to prove versions of Propostion 6.3.4 must be modified. Instead the problem is reduced to studying the quantity:

$$\mathcal{M}(x) := \sup_{r > 0} \frac{1}{m_\gamma(B(x, r))} \int_{B(x, r)} \varphi(y) dm_\gamma(y),$$

where $m_\gamma$ is the Riemannian measure on the local unstable manifold $\gamma$ centered at $x$, and $\varphi_n(x) := 1_{E_n}(x)\rho_\gamma(x)$ where $\rho_\gamma$ is the conditional density of $\mu$ on $\gamma$ (which is absolutely continuous with respect to $m_\gamma$.) Conditions analogous to conditions (H2a) and (H2b) are then verified on local unstable manifolds. On an example driven basis arguments based on the regularity of unstable foliations, and on the geometrical properties of the system can sometimes be used to establish the (SRT) conditions.

## 6.6
## Hyperbolic Dynamical Systems

Hyperbolic toral automorphisms are canonical models of invertible uniformly hyperbolic dynamics. We will consider in detail one of the best studied, the Arnold cat map, mainly from the point of view of establishing the rate of convergence to EVLs for observables maximized at generic points.

### 6.6.1
### Arnold Cat Map

The Arnold cat map $(f, \mathbb{T}^2, \mu)$ is given by

$$f(x, y) = (2x + y, x + y) \mod \mathbb{T}^2. \tag{6.6.1}$$

where $\mathbb{T}^2 = \mathbb{R}^2/\mathbb{Z}^2$ is the two-torus and $\mu$ is two-dimensional Lebesgue measure or Haar measure. We are able to establish the following convergence rates:

**Proposition 6.6.1.** *Suppose $(f, \mathbb{T}^2, \mu)$ is the Arnold cat map defined by Eq. 6.6.1. If $\phi(x) = -\log(dist(x, \tilde{x}))$, then for $\mu$-a.e $\tilde{x} \in [0, 1]$ and all $\epsilon > 0$:*

$$\left| \mu\{M_n \leq u_n\} - (1 - \sqrt{n}\mu\{X_1 > u_n\})^{\sqrt{n}} \right| \leq \frac{C}{n^{\frac{1}{2} - \epsilon}}, \tag{6.6.2}$$

*Here $C$ is a constant depending on $\epsilon$ and $\tilde{x}$. Hence*

$$\lim_{n \to \infty} \mu\{M_n \leq (u + \log n)/2\} = e^{-\pi e^{-u}}, \tag{6.6.3}$$

*where $C$ is independent of $n$, but dependent on $\tilde{x}$.*

*Proof.* This map is uniformly hyperbolic, admits a finite Markov partition, and has expansion estimates on unstable manifolds. Thus condition (H1s) holds, with $\Theta(j) \leq O(\theta_0^j)$ for some $\theta_0 < 1$. Since the invariant SRB measure is Lebesgue measure it is more natural to check condition (H2b) (rather than (SRT2)) and apply the proof of Proposition 6.3.4 to establish the error bound.

Consider the set $E_n$ as defined in Eq. (6.3.1), with $\tilde{g}(n) = (\log n)^{\gamma'}$ for some $\gamma' > 1$. Using the Markov structure of the Arnold cat map we can repeat the method of Section 6.4.1 to estimate $\mu(E_n)$. Since the conditional invariant measures on unstable manifolds are uniformly equivalent to Lebesgue measure we obtain the estimate

$$\mu(E_n) \leq \frac{(\log n)^{\gamma'}}{n}.$$

We then follow the proof of Proposition 6.3.4 step by step. The only technical modification being that we apply a Hardy Littlewood maximal operator on functions $\varphi : \mathbb{R}^2 \to \mathbb{R}$, with $\varphi \in L^1(m)$, (and $m$ is two dimensional Lebesgue measure). See also [139] for analogous higher dimensional arguments. We then obtain the bound

$$\sum_{j=1}^{g(k)} \mu(X_1 > u_k, X_j > u_k) = O(k^{-3/2 + \epsilon} g(k)),$$

valid for all $\epsilon > 0$. We now take $g(k) = (\log k)^{1 + \epsilon'}$ for some $\epsilon' < \gamma'$, and following the proof of Proposition 6.3.3 this latter bound is the main contribution to the error term. Due to exponential decay of correlations the remaining error terms are of higher order. Condition (H3) clearly holds since $\mu$ is Lebesgue measure. Hence for any $\epsilon > 0$:

$$\left| \mu\{M_n \leq u_n\} - (1 - \sqrt{n}\mu\{X_1 > u_n\})^{\sqrt{n}} \right| \leq \frac{C}{n^{\frac{1}{2} - \epsilon}}, \tag{6.6.4}$$

and the result follows. $\qquad \square$

### 6.6.2
### Lozi-like Maps

The Lozi map $f$ is a homeomorphism of $\mathbb{R}^2$ given by

$$f_{a,b}(x, y) = (1 + y - a|x|, bx) \tag{6.6.5}$$

where $a$ and $b$ are parameters. A graphic representation of the attractor of this system is given in Fig. 6.1 for $a = 1.7$ and $b = 0.2$. It has been studied as a model of chaotic dynamics intermediate in complexity (or difficulty) between Axiom A diffeomorphisms and Hénon diffeomorphisms [213, 214, 215]. The derivative is discontinuous on the $y$-axis and this leads to arbitrarily short local unstable manifolds. Misuiurewicz [213] proved that there exists a positive measure set $G$ of parameters such that if $(a, b) \in G$ the map $f$ is hyperbolic. If $(a, b) \in G$, then $f_{a,b}$ has invariant stable and unstable directions (where the derivative is defined) and the angle between them is bounded below by $\pi/5$. We will restrict our attention to maps with parameters in the set $G$.

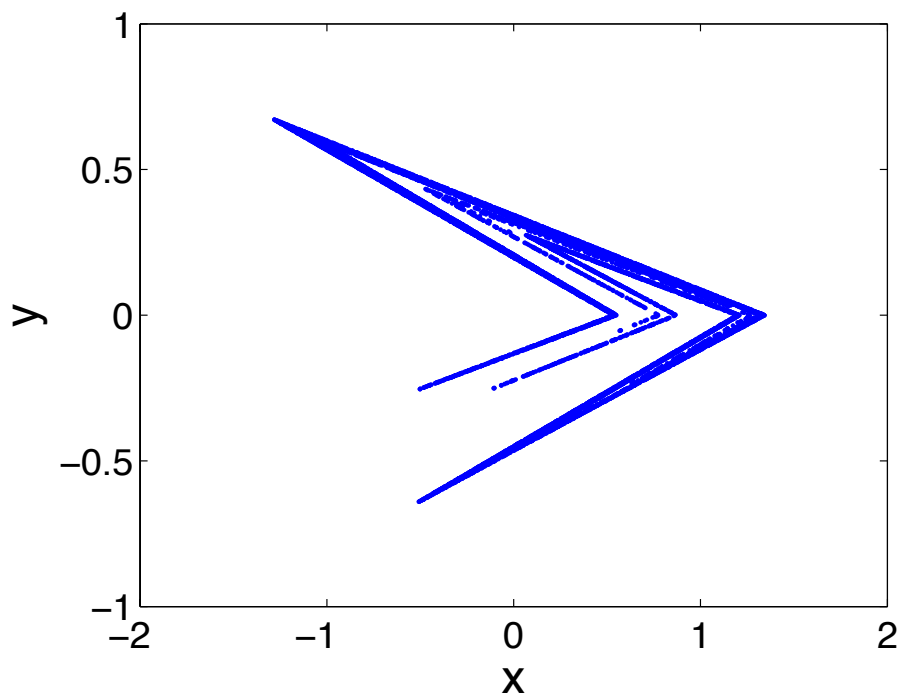

**Figure 6.1** Numerical approximation to the attractor of the Lozi map given in Eq. 6.6.5 for $a = 1.7$ and $b = 0.5$.

The tangent derivatives, where defined, satisfy uniform expansion estimates [214] in that there exists $\lambda > 1$ such that $|Df^n v| \geq \lambda^n v$ for all $v \in E^u$, the unstable direction and correspondingly for $E^s$, the stable direction. $f_{a,b}$ has an invariant ergodic probability measure $\mu$ [214] which is absolutely continuous with respect to the one-dimensional Lebesgue measure along local unstable curves. One reason for restricting to maps $f_{a,b}$, $(a, b) \in G$ is that for such maps Collet and Levy have shown that for $\mu$ almost every point on the attractor the Hausdorff dimension of $\mu$ exists and is constant [214]. Hence these maps fall within the class of systems we have axiomatized. The existence of a pointwise local dimension $d$ which is $\mu$ a.e. constant implies that for almost every $x$ in the attractor, the dimension constant $d$ in the definition of $u_n$ is the same. We will use a sequence of scaling constants $u_n(x_0, v)$ defined for a generic point $x_0$ by the requirement that $n\mu(B(x_0, e^{-u_n}) \to$

$e^{-v}$.

The Lozi map $f_{a,b}$ with $(a, b) \in G$, $b$ sufficiently small admits such a Tower with exponential tails [156]. Hence the Lozi maps we consider satisfy exponential decay of correlations for Hölder continuous observations.

**Proposition 6.6.2.** *Let $f_{a,b} : M \to M$ be a Lozi map with $(a, b) \in G$ with $b$ sufficiently small. Then for $\mu$ a.e. $x_0$ the stochastic process defined by $X_n(x) = -\log(d(x_0, f^n x))$ satisfies a Type I extreme value law in the sense that $\lim_{n \to \infty} \mu(M_n \leq u_n(v)) = e^{-e^{-v}}$.*

*Remark* 6.6.3. We don't have useful convergence rates in this setting.

We do not know the precise scaling constants $u_n$, but for all $\epsilon > 0$, $\lim_{n \to \infty} \mu(M_n \leq (1 - \epsilon)(\log n + v)/d)) \leq e^{-e^{-v}} \leq \lim_{n \to \infty} \mu(M_n \leq (1 + \epsilon)(\log n + v)/d))$ which provides an estimate of the correct sequence $u_n$.

### 6.6.3
### Sinai Dispersing Billiards

Let $\Gamma = \{\Gamma_i, i = 1 : k\}$ be a family of pairwise disjoint, simply connected $C^3$ curves with strictly positive curvature on the two-dimensional torus $\mathbb{T}^2$. We consider the billiard map $f : \mathbb{T}^2 \setminus \Gamma \to \mathbb{T}^2 \setminus \Gamma$. We assume the finite horizon condition, namely, that the number of successive tangential collisions of the freely moving particle with the convex scatterers is bounded above. More precisely we take $M := \Gamma \times [-\pi/2, \pi/2]$ and we will let $f$ be the Poincaré map that gives the position and instantaneous velocity after collision, according to the rule angle of incidence equals angle of reflection. The billiard map preserves a measure $\mu$ equivalent to 2-dimensional Lebesgue measure $m$ with density $\rho(x_0) = \frac{d\mu}{dm}(x_0)$. We refer to the book by Chernov and Markarian [216] for a comprehensive exposition of the ergodic properties of hyperbolic billiards.

**Proposition 6.6.4.** *Let $f : M \to M$ be a Sinai dispersing billiard map. Then for $\mu$ a.e. $x_0$ the stochastic process defined by $X_n(x) = -\log(d(x_0, f^n x))$ satisfies a Type I extreme value law in the sense that $\lim_{n \to \infty} \mu(M_n \leq (v + \log n + \log(\rho(x_0)))/2) = e^{-e^{-v}}$.*

### 6.6.4
### Hénon maps

Hénon maps $f$ are diffeomorphisms of $\mathbb{R}^2$ given by

$$f_{a,b}(x, y) = (1 + y - ax^2, bx), \tag{6.6.6}$$

for parameters $a, b \in \mathbb{R}$. A representation of the attractor of this map is given in Fig. 6.2 for $a = 1.4$ and $b = 0.3$. Unlike the Lozi maps, these maps are non-uniformly hyperbolic. It has been shown [196] that there is a positive measure set of parameters ($a \simeq 2$ and $|b| \ll 1$) for which the invariant set is a strange attractor (containing a dense orbit with positive Lyapunov exponent). For the same parameters, the Hénon map admits an SRB measure $\mu$, and can be modeled by a Young tower with exponential decay of correlations, see [217, 156]. Thus condition (H1) holds. The verification of conditions (H2b) and (H3) rely on technical calculations using the Besicovitch Covering Lemma, see [163]. The main technical obstacles in the Hénon map being the non-invariance of the stable and unstable foliations.

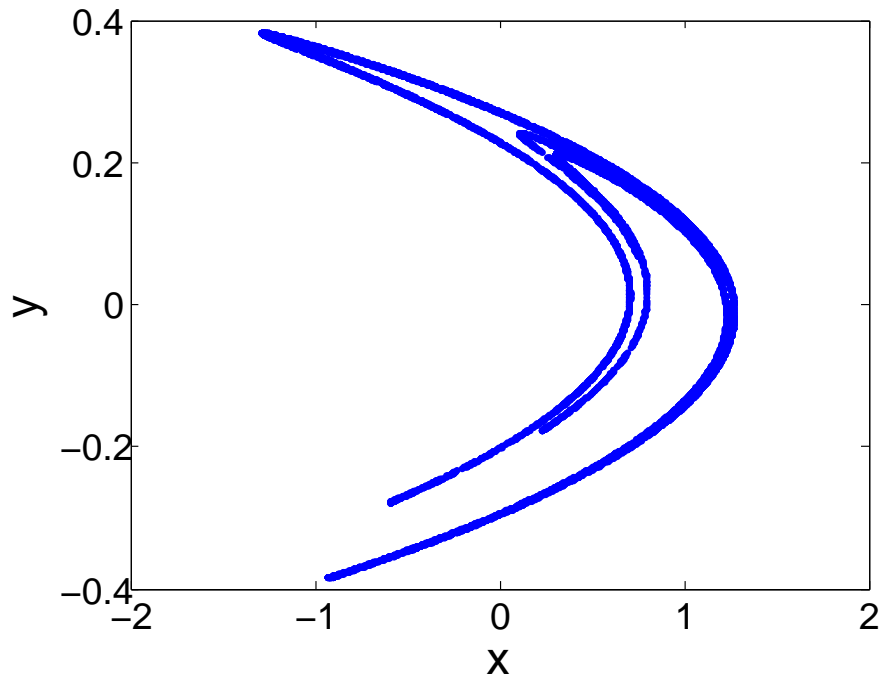

**Figure 6.2** Numerical approximation to the attractor of the Hénon map given in Eq. 6.6.6 for $a = 1.4$ and $b = 0.3$.

**Proposition 6.6.5.** *Let $f_{a,b} : M \to M$ be a Hénon map with $(a, b) \in \Omega$. Then for $\mu$ a.e. $x_0$ the stochastic process defined by $X_n(x) = -\log(d(x_0, f^n x))$ satisfies a Type I extreme value law in the sense that $\lim_{n \to \infty} \mu(M_n \leq u_n(v)) = e^{-e^{-v}}$*

*Remark 6.6.6.* We do not know the precise scaling constants $u_n$, but for all $\epsilon > 0$, $\lim_{n \to \infty} \mu(M_n \leq (1 - \epsilon)(\log n + v)/d) \leq e^{-e^{-v}} \leq \lim_{n \to \infty} \mu(M_n \leq (1 + \epsilon)(\log n + v)/d))$ which provides an estimate of the correct sequence $u_n$.

## 6.7
## Skew-product Extensions of Dynamical Systems.

In this section we review the results of [137] concerning extreme value laws for skew product extensions, in particular compact group extensions of hyperbolic systems. Suppose that $Y$ is a compact, connected, $M$-dimensional manifold with metric $d_Y$ and $X$ is a compact $N$-dimensional manifold with metric $d_X$. We let $D = M + N$ and define a metric on $X \times Y$ by

$$d((x_1, \theta_1), (x_2, \theta_2)) = \sqrt{d_X(x_1, x_2)^2 + d_Y(\theta_1, \theta_2)^2}.$$

We denote the Lebesgue measure on X by $m_X$, the Lebesgue measure on $Y$ by $m_Y$ and the product measure on $X \times Y$ Y by $m = m_X \times m_Y$.

If $f : X \to X$ is a measurable transformation and $u : X \times Y \to Y$ a measurable function, then we may define $g$, the $Y$-skew extension of $f$ by $u$, via: $g : X \times Y \to X \times Y$:

$$g(x, \theta) = (f(x), u(x, \theta)).$$

We assume further that $T : X \to X$ has an ergodic invariant measure $\mu_X$, and $f$ preserves an invariant probability measure $\mu$ with density in $h_\mu \in L^p(m)$. Given $(\tilde{x}, \tilde{\theta})$, we consider the observable function $\phi(x, \theta)$ with representation

$$\phi(x, \theta) = \psi \mathrm{dist}((x, \theta), (\tilde{x}, \tilde{\theta})),$$

where $\psi : \mathbb{R}^+ \to \mathbb{R}$ takes it maximum value at 0.

**Theorem 6.7.1.** *Assume that $h_\mu \in L^{1+\delta}$ for some $\delta > 0$, and*

*1) There exist constants $C_1 > 0$, $\beta > 0$ and an increasing function $l(n) \approx n^{\gamma D}$ (with $0 < \gamma < 1$) such that if:*

$$E_n^X := \left\{ x \in X : d_X(f^j x, x) < \frac{1}{n}, \text{ some } g \leq l(n) \right\},$$

*then $\mu_X(E_n^X) < Cn^{-\beta}$.*

*2) Condition (H1) holds for $(g, \mu)$, and there exists $\zeta > 0$ such that $\Theta(n) \leq n^{-\zeta}$.*

*Suppose $\phi(x) = -\log \mathrm{dist}((x, \theta), (\tilde{x}, \tilde{\theta}))$, then for $\mu$-a.e $(\tilde{x}, \tilde{\theta})$ and for every $u \in \mathbb{R}$:*

$$\lim_{n \to \infty} \mu(M_n < u + \log n) = \exp\{-h_\mu(\tilde{x}, \tilde{\theta})e^{-Du}\} \qquad (6.7.1)$$

The method of proof is simple [137] and relies on the observation that if points don't have short returns in the base system they also don't have short returns in the skew-product system. In [137] applications to the extreme value theory of compact group extensions of hyperbolic systems and other partially hyperbolic dynamical systems are considered.

## 6.8
## On the Rate of Convergence to an Extreme Value Distribution

In this section we focus on the speed of convergence of $\mu\{M_n \leq u_n\}$ to $G(u)$. In the case of i.i.d random variables $\hat{X}_i$, with probability distribution function $F(u) :=$ $P(\hat{X}_i \le u)$ the rate of convergence depends on the normalization sequences $a_n, b_n$, and the functional form of $F(u) := P(\hat{X}_i \le u)$ as $u \to u_F := \max \hat{X}_i$, see [1]. For illustration in the i.i.d case, uniform estimates on rates of convergence can be obtained for the exponential distribution and the Gaussian distribution. For the exponential distribution $F(u) = 1 - e^{-u}$, it is shown in [218] that:

$$\sup_u \left| P\left\{ M_n \le u + \frac{\log n}{\lambda} \right\} - e^{-e^{-u}} \right| = \sup_u \left| \left(1 - \frac{e^{-u}}{n}\right)^n - e^{-e^{-u}} \right|$$

$$\le \frac{1}{n}\left(1 + \frac{2}{n}\right)e^{-2}. \tag{6.8.1}$$

However, in general the convergence rate can be quite slow, and this is evident for the Gaussian distribution [219], where it is established that:

$$\frac{C_1}{\log n} \le \sup_u \left| \{\Phi(u/a_n + b_n)\}^n - e^{-e^{-u}} \right| \le \frac{C_2}{\log n}. \tag{6.8.2}$$

Here $\Phi(u)$ is the standard Gaussian distribution function, $C_1, C_2$ are uniform constants, and $a_n, b_n$ satisfy:

$$a_n = b_n, \; 2\pi b_n^2 e^{b_n^2} = n^2.$$

In this case, the choice of constants is optimal for the error rate.

In the setting of dynamical systems, we consider the corresponding quantity $\mu\{x : \phi(x) < u\}$ (which corresponds to $F(u)$), and study the behaviour of this measure as $u \to \max \phi$. We take the scaling sequence $u_n = u/a_n + b_n$ (as we would in the i.i.d case), and study the errors involved in approximating the limit distibution $G(u)$ by $\mu\{M_n \le u_n\}$. For a given sequence $u_n$, we define functions $\tau_n(u)$ and $G_n(u)$ by

$$\tau_n(u) = n\mu\{\phi(x) \ge u_n\}, \quad G_n(u) = \left(1 - \frac{\tau_n(u)}{n}\right)^n.$$

The function $\tau_n(u) \to \tau(u)$ uniformly for all $u$ lying in a compact subset of $\mathbb{R}$.

The problem is to estimate $|\mu\{M_n \le u_n\} - G(u)|$ (as $n \to \infty$, for $u$ lying in compact subsets of $\mathbb{R}$). However, we cannot immediately express $\mu\{M_n \le u_n\}$ in terms of $G_n(u)$. However, the blocking arguments used in proving Theorem 6.3.1 give a partial answer towards estimating the error. Indeed as we have shown in Theorem 6.3.1, $\mu\{M_n \le u_n\}$ can be approximated by $G_{n^\beta}(u)$ (for some $\beta \in (0,1]$), up to an error of order $O(n^{-\alpha})$ for some constant $\alpha > 0$. Here

$$G_{n^\beta}(u) = \left(1 - \frac{\tau_n(u)}{n^\beta}\right)^{n^\beta},$$

and the optimal choice of $\beta$ turns out to be $1/2$ for most applications. Bounds on the constant $\alpha$ depend on the recurrence statistics and rates of decay of correlations. In applications, we achieve $\alpha$ close to $1/2$. Estimating the error between $G_{n^\beta}(u)$ and $G(u)$ depends on the regularity of the observable $\phi$ and on the regularity of the measure $\mu$. We now describe the methods involved. We have the following elementary results, see [1].

**Proposition 6.8.1.** *For all $n \geq 2$:*

$$\Delta_{n^\beta} := \left| \left( 1 - \frac{\tau_n(u)}{n^\beta} \right)^{n^\beta} - e^{-\tau_n(u)} \right| \leq \frac{(\tau_n(u))^2 e^{-\tau_n(u)}}{2(n^\beta - 1)} \leq \frac{0.3}{(n^\beta - 1)}. \quad (6.8.3)$$

*Suppose that $\tau_n$ and $\tau$ satisfy $|\tau_n - \tau| \leq \log 2$. Then for some $\theta \in (0,1)$:*

$$\Delta'_n := |e^{-\tau_n(u)} - e^{-\tau(u)}| \leq e^{-\tau(u)} \{ |\tau(u) - \tau_n(u)| + \theta(\tau(u) - \tau_n(u))^2 \}. \quad (6.8.4)$$

We remark that the bound for $\Delta_{n^\beta}$ is uniform in $u$. This gives the corollary:

**Corollary 6.8.2.** *Suppose $(f, \mathcal{X}, \mu)$ is an ergodic dynamical system, then:*

$$|\mu\{M_n \leq u_n\} - G(u)| \leq |\mu\{M_n \leq u_n\} - G_{\sqrt{n}}(u)| + \Delta_{\sqrt{n}} + \Delta'_n. \quad (6.8.5)$$

Theorem 6.3.1 is used to estimate the first term on the right hand side of Eq. (6.8.5). The estimation of $\Delta_n$ and $\Delta'_n$ require knowledge of the explicit representation of $\phi(x)$. We will consider four explicit forms of $\phi(x)$, and then comment on the general cases. For these explicit forms, we will give explicit bounds on $\Delta_n$ and $\Delta'_n$ under the additional assumption that the density at $\tilde{x}$ is Lipschitz. For the class of non-uniformly expanding systems under consideration, this is a reasonable assumption. Consider the following representations $\phi = \phi_i : \mathcal{X} \to \mathbb{R}$ (for $i =$1,2 and 3) defined by:

$$\phi_1(x) = -\log(\text{dist}(x, \tilde{x})), \ \phi_2(x) = \text{dist}(x, \tilde{x})^{-\alpha}, \ \phi_3(x) = C - \text{dist}(x, \tilde{x})^\alpha, \quad (6.8.6)$$

for $\alpha > 0$ and $C \in \mathbb{R}$. We also consider a fourth observable:

$$\phi_4(x) = \psi_4(\text{dist}(x, \tilde{x})) : \quad \psi_4^{-1}(y) = -(y \log y)^{-1}. \quad (6.8.7)$$

In the cases of Eq. (6.8.6) the scaling sequences $u_n = u/a_n + b_n$ can be made explicit with $n\mu\{\phi > u_n\} \sim \tau(u)$, and we have:

$$\mu\{\phi_1(x) > u + \log n\} \sim 2\rho(\tilde{x})e^{-u}/n,$$
$$\mu\{\phi_2(x) > un^\alpha\} \sim 2\rho(\tilde{x})u^{-1/\alpha}/n,$$
$$\mu\{\phi_3(x) \geq C - u/n^\alpha\} \sim 2\rho(\tilde{x})u^{1/\alpha}/n.$$

Since the density is Lipschitz, the higher order terms in the above set of asymptotics are all $O(1/n^2)$, where the implied constant depends on the Lipschitz norm. Hence by Proposition 6.8.1, $\Delta_{\sqrt{n}}$ is bounded by $O(1/\sqrt{n})$ and $\Delta'_n$ is bounded by $O(1/n)$. If instead the invariant density has Hölder exponent $\beta \in (0,1)$, the error $\Delta'_n$ is bounded by $O(1/n^\beta)$. For observations that are general (regularly varying) functions of $\text{dist}(x, \tilde{x})$, the scaling sequences $a_n$ and $b_n$ cannot always be made explicit and this leads to weaker estimates on the error bounds $\Delta_n$ and $\Delta'_n$. Indeed, consider the observable $\phi_4(x)$, and for simplicity take $\mu$ to be Lebesgue measure. This observable is regularly varying with index $-1$, and we are therefore in the domain of attraction of a Type II distribution. We take $b_n = 0$ and estimate $a_n$ so that $n\mu\{\phi(x) > u_n\} \to \tau(u) = u^{-1}$. In this case $\phi(x) = \psi_4(\text{dist}(x, \tilde{x}))$ with $\psi_4^{-1}(y) = -(y \log y)^{-1}$. Using asymptotic inversion the orginal function $\psi(y)$ takes the form $\psi(y) \sim -y^{-1}(\log y)$ as $y \to 0$. Using the particular form of $\psi(y)$ we try $a_n = \log n/n$, and obtain:

$$
\begin{aligned}
|\tau_n(u) - \tau(u)| &= \left| n\text{Leb}\{\text{dist}(x, \tilde{x}) \leq \psi^{-1}(un \log n)\} - 2u^{-1} \right| \\
&= 2u^{-1}\left| \log n \, (\log \log n - \log u - \log n)^{-1} - 1 \right| \\
&= 2u^{-1}\frac{\log \log n}{\log n} + O\left(\frac{(\log \log n)^2}{(\log n)^2}\right),
\end{aligned}
$$

and hence for fixed $u$, this error is of order $1/(\log n)^{1-\epsilon}$ for any $\epsilon > 0$. It is feasible other sequences $a_n$ give rise to smaller errors, but significant improvement is not expected. For the Gaussian distribution, only bounds of order $1/(\log n)$ can be achieved, see [1].

### 6.8.1
### Error Rates for Specific Dynamical Systems

Based on the results of Theorem 6.3.1 and Corollary 6.8.2 we can obtain error rates for the systems studied in Section 6.4. For simplicity we consider the distance observable $\phi(x) = -\log|x - \tilde{x}|$, and hence we obtain the convergence rate associated to a Type I distribution. Further details can be found in in [144], but we give the main steps below

**The Tent Map.** Convergence to EVL for the tent map was proved in Section 6.4.1. We have the following result concerning the rate of convergence to an EVL.

**Proposition 6.8.3.** *Suppose* $(f, [0, 1], \text{Leb})$ *is the tent map. For the observation* $\phi(x) = -\log|x - \tilde{x}|$, *we have for* Leb-*a.e* $\tilde{x} \in [0, 1]$ *and all* $\epsilon > 0$:

$$
|\mu\{M_n \leq u + \log n\} - e^{-2e^{-u}}| \leq C\frac{(\log n)^{1+\epsilon}}{\sqrt{n}}, \tag{6.8.8}
$$

*where* $C(\tilde{x}) > 0$ *is a uniform constant depending on* $\tilde{x}$.

To analyse $\Delta_n, \Delta_n'$ from the function form $\phi(x) = -\log|x - \tilde{x}|$, we have in this case equality $\mu\{X_1 > u + \log n\} = 2e^{-u}/n$ (since $\mu$ is Lebesgue measure). Hence $\tau(u) = 2e^{-u}$, and Proposition 6.8.1 implies that:

$$
\Delta_{\sqrt{n}} + \Delta_n' \leq Cn^{-1/2},
$$

and so the result follows.

**Intermittent Maps.** We consider the map intermittent map defined by Eq. (6.4.5) in Section 6.4.4. We have the following result

**Proposition 6.8.4.** *Suppose* $(f, [0, 1], \mu)$ *is the intermittent system* (6.4.5) *with* $b < 1/20$ *and consider the observable* $\phi(x) = -\log|x - \tilde{x}|$. *Then for all* $\epsilon > 0$ *and* $\mu$-*a.e.* $\tilde{x} \in \mathcal{X}$ *we have that*

$$\left| \mu\{M_n \leq u_n\} - e^{-2\rho(\tilde{x})e^{-u}} \right| \leq Cn^{-\frac{1}{2}+\tilde{b}}, \quad with \ \tilde{b} = \epsilon + 10b. \tag{6.8.9}$$

*where* $C(\tilde{x}) > 0$ *is a constant independent of* $n$.

The proof of this proposition has two steps. In the first step, equation (6.3.5) of Theorem 6.3.1 applies, and its right hand side is estimated using properties of the regularity of the invariant density, and knowledge of the rate of decay of correlations, each of which can be expressed in terms of the parameter $b$. Moreover Lemma 6.4.6 can be used to estimate the measure of the recurrence sets $E_n$.

In the second step, we analyse, $\Delta_{\sqrt{n}}$ and $\Delta'_n$ (recalling a block size $q \sim \sqrt{n}$). For $\tilde{x} \neq 0$, we can apply Lebesgue differentiation to deduce that

$$\mu\{\phi(x) > u + \log n\} \sim 2\rho(\tilde{x})e^{-u}/n.$$

However we need higher order regularity information on $\mu$ to deduce bounds on $\Delta_n, \Delta'_n$. In [186], the density of $\rho$ is in fact locally Lipschitz $\mu$-a.e. The fact that $\rho \in L^{1+\delta}$ is proved from an analysis of the singularity in $\rho(x)$ at $x = 0$. The density scales as $O(x^{-b})$ as $x \to 0$. Hence for $\tilde{x} \neq 0$, we have

$$\mu\{\phi(x) > u + \log n\} = 2\rho(\tilde{x})e^{-u}/n + O(n^{-2}),$$

and this turns out to be of higher order error relative to the estimate in Eq. (6.3.5).

**The quadratic map** We consider the quadratic map defined by Eq. (6.4.3) in Section 6.4.2. We have the following result:

**Proposition 6.8.5.** *Suppose* $(f, [-2, 2], \mu)$ *is the quadratic family with* $a$ *is a Misiurewicz parameter. For the observation* $\phi(x) = -\log|x - \tilde{x}|$, *we have for* $\mu$-*a.e* $\tilde{x} \in [0, 1]$ *and all* $\epsilon > 0$:

$$|\mu\{M_n \leq u + \log n\} - e^{-2\rho(\tilde{x})e^{-u}}| \leq C_1 \frac{(\log n)^{1+\epsilon}}{\sqrt{n}} + \frac{C_2}{n^{\alpha-\epsilon}}. \tag{6.8.10}$$

*where* $C_1, C_2 > 0$ *are uniform constants independent of* $n$, *but dependent on* $\tilde{x}$.

This is proved in [144]. To get quantitative error estimates we require the critical orbit to satisfy a Misiurewicz condition, see [220]. Such a condition is required to ensure the invariant density is sufficiently regular. If $\alpha > 1/2$ then the error bound in Eq. (6.4.3) is improved to $O(n^{-1/2+\epsilon})$ (for any $\epsilon > 0$). It is feasible to extract bounds on the constant $\alpha$ using the proof of [72, Proposition 2.4]. However such bounds will in turn depend on the (exponential) rate of decay of correlations.

**Hyperbolic systems**  For general hyperbolic systems, the error term in Theorem 6.5.1 is the best that can be achieved. For systems with SRB measures, the function $\tau_n(u) = n\mu\{\phi(x) > u_n\}$ need not converge to a limit as $n \to \infty$. An example of a (uniformly) hyperbolic system for which error rates can be achieved is in the case of the Arnold cat map $(f, \mathbb{T}^2, \mu)$ as discussed in Section 6.6.1. We have the following result:

**Proposition 6.8.6.**  *Suppose $(f, \mathbb{T}^2, \mu)$ is the Arnold cat map defined by Eq. (6.6.1). For the observation $\phi(x) = -\log(dist(x, \tilde{x}))$, we have for $\mu$-a.e $\tilde{x} \in [0, 1]$ and all $\epsilon > 0$:*

$$|\mu\{M_n \leq (u + \log n)/2\} - e^{-\pi e^{-u}}| \leq \frac{C(\tilde{x})}{n^{\frac{1}{2} - \epsilon}}, \qquad (6.8.11)$$

*where $C$ is independent of $n$, but dependent on $\tilde{x}$.*

The proof follows that of Proposition 6.6.1 but we keep track of the decay rate of $\mathrm{Leb}(E_n)$. The method of estimating $\Delta_n$ and $\Delta_n'$ is identical to that for the tent maps. The result follows.

## 6.9
## Extreme Value Theory for Deterministic Flows

A common model of a continuous time dynamical system $(\phi_t, M)$ (where $\phi_t$ is a flow on a manifold $M$) is a suspension flow. One way this model arises is to take a Poincaré section $X$, where $X$ is a codimension-one submanifold tangent to the flow direction, and study the first return map $f : X \to X$. More precisely if we define the first return time $h : X \to \mathbb{R}$ by $h(x) := \inf\{t > 0 : \phi_t(x) \in X\}$ then $f : X \to X$ is given by $f(x) := \phi_{h(x)}(x)$.

We then define the suspension space

$$X^h = \{(x, u) \in X \times \mathbb{R} \mid 0 \leq u \leq h(x)\}/\sim, \qquad (x, h(x)) \sim (f(x), 0)$$

The flow $\phi_t : M \to M$ is modeled by the suspension (semi) flow $f_s : X^h \to X^h$,

$$f_s(x, u) = \begin{cases} (x, u + s) & \text{if } u + s < h(x); \\ (f(x), 0) & \text{if } u + s = h(x). \end{cases}$$

The projection $\pi(x, t) = \phi_t(x)$ semi-conjugates the original flow and the semi-flow, and $\phi_t \circ \pi(x, u) = \pi f_t(x, u)$. For more details of this construction see [221]

In general the study of the ergodic and mixing properties of flows is much more delicate than that of discrete time transformations. One advantage of the suspension flow model is that it (sometimes) enables ergodic and statistical properties of the underlying base map to be 'lifted' to the suspension flow. A systematic study in this direction is given in [222].

The results we discuss in this section concerning extreme value theory for suspension flows of chaotic base maps are mainly found in [139].

These results may also be seen as an extension of certain theorems in extreme value theory where instead of sampling at a 'random' time $h(x)$ a discrete sampling at constant times is considered(see [1, Theorem 13.3.2] and related results in [1, Chapter 13]).

To use a suspension flow to understand statistical properties we need to relate invariant measures of the base map $f : X \to X$ to invariant measures for the suspension flow. In fact, we may always lift an invariant measure $\mu$ for $f : X \to X$ to an invariant measure $\mu^h$ for the suspension flow $f_s : X^h \to X^h$.

Assume that $T : X \to X$ preserves the probability measure $\mu$. We suppose that $h \in L^1(\mu)$ is a positive roof function.

The $f$ invariant probability measure $\mu$ lifts to a flow-invariant probability measure $\mu^h$ on $X^h$ given by $d\mu^h = d\mu \times du/\bar{h}$ and $\bar{h} = \int_X h \, d\mu$, where $du$ is Lebesgue measure.

Given a (measurable) observation $\tilde{\phi} : M \to \mathbb{R}$ we may lift it to $\phi : X^h \to \mathbb{R}$ by defining $\phi(x, u) = \phi \circ \pi(x, u)$. By the semi-conjugacy that $\pi$ provides statistical properties of the stationary process $\tilde{\phi}_t$ sampled from $(M, \mu)$ are the same as those of $\phi_t$ sampled from $(X^h, \mu^h)$.

Define $\Phi : X \to \mathbb{R}$ by

$$\Phi(x) := \max\{\phi(f_s(x)) \mid 0 \le s < h(x)\}. \tag{6.9.1}$$

Denote

$$\phi_t(p) := \max\{\phi(f_s(p)) \mid 0 \le s < t\}$$
$$\Phi_N(x) := \max\{\Phi(f^k(x)) \mid 0 \le k < N\}. \tag{6.9.2}$$

In this setting [139]:

**Theorem 6.9.1** (discrete or continuous suspensions). *Assume that $T : (X, \mu) \to (X, \mu)$ is ergodic and $h \in L^1(\mu)$. Suppose also that the normalizing constants $a_n > 0$ and $b_n$ satisfy:*

$$\lim_{\epsilon \to 0} \limsup_{n \to \infty} a_n |b_{[n+\epsilon n]} - b_n| = 0, \tag{6.9.3}$$

$$\lim_{\epsilon \to 0} \limsup_{n \to \infty} \left| 1 - \frac{a_{[n+\epsilon n]}}{a_n} \right| = 0. \tag{6.9.4}$$

*Then, if $G$ is a non-degenerate distribution,*

$$a_N(\Phi_N - b_N) \to G \implies a_{\lfloor T/\bar{h} \rfloor}(\phi_T - b_{\lfloor T/\bar{h} \rfloor}) \to G. \tag{6.9.5}$$

*Furthermore if $h^{-1} \in L^1(\mu)$ then*

$$a_N(\Phi_N - b_N) \to G \iff a_{\lfloor T/\bar{h} \rfloor}(\phi_T - b_{\lfloor T/\bar{h} \rfloor}) \to G. \tag{6.9.6}$$

Equations (6.9.3) and (6.9.4) are satisfied for most of the linear scalings discussed in this book. For example, they are satisfied in the Type I case where $b_n = \log n$ and $a_n = 1$, in the Type II and Type III scenarios when $a_n$ is regularly varying and $b_n$ is constant. Thus the assumptions on $a_n, b_n$ hold for logarithmic singularities, power singularities and power function maxima. In the flow case the constants $a_n, b_n$ are not determined by the requirement $\mu^h(\phi > \frac{v}{a_n} + b_n) = O(\frac{1}{n})$.

The proof of this theorem follows this route. First we show that the convergence is mixing in distribution [223].

A sequence of random variables $S_n : X \to \mathbb{R}$ on a probability space $(X, \mu)$ *converges mixing in distribution* to $G$ if for each $A \subset X$ with positive measure, $S_n|_A \to G$ with respect to the conditional measure $\mu_A(B) := \mu(B \cap A)/\mu(A)$ on $A$. This notion was introduced by Renyi [224].

A useful criterion to establish mixing in distribution was given by Eagleson [223, Thm. 6]: if $X_1, X_2, \ldots$ have a trivial invariant $\sigma$-field, (which is a consequence of ergodicity), $T_n(X_1, \ldots, X_n) \to T$ for a sequence of statistics $T_n$, and $T_n(X_1, \ldots, X_n) - T_n(X_{k+1}, \ldots, X_{k+n}) \to 0$ for each $k$, then $T_n(X_1, \ldots, X_n) \to T$ is mixing in distribution.

Following proving mixing in distribution we establish the distributional convergence of $\phi_T(p)$ by obtaining bounds on the distributional convergence of $\Phi_N(x)$, where $N = N(x, T)$ is a function of time (governed by measure $\mu$) which satisfies the strong law of large numbers.

Let $(X, \mu, f)$ be a m.p.t. on a probability space and $\Phi : X \to \mathbb{R}$ an (a.e. finite) random variable.

Denote $\Phi_{n,m}(x) := \max\{\Phi \circ f^k(x) \mid n \leq k < m\}$, for $n, m \geq 0$. Hence $\Phi_n = \Phi_{0,n}$.

**Lemma 6.9.2.** *If $f$ is ergodic then:*

(a) $\Phi_n \to \text{ess} - \sup(\Phi)$ *a.e.*

(b) $a_n(\Phi_{0,n} - \Phi_{k,n+\ell}) \to 0$ *as $n \to \infty$, for each $k \geq 0$ and $\ell$. Therefore, by [223, Thm. 6], if $a_n(\Phi_n - b_n) \to G$, then the convergence to $G$ is mixing in distribution. Moreover, for each $k \geq 0$ and $\ell$, $a_n(\Phi_{k,n+\ell} - b_n) \to G$ mixing as well.*

(c) *Let $h \in L^1(\mu)$ be a roof function and denote $\hat{\Phi}_n := \Phi_n \circ \pi^h$. If $a_n(\hat{\Phi}_n - b_n) \to G$ on $(X^h, \mu^h)$, then the convergence to $G$ is mixing in distribution .*

*Remark 6.9.3.* Note that if $a_n(\Phi_n - b_n) \to G$, then (b) implies that $a_n(\Phi_{n+k} - b_n) \to G$ for any fixed $k$. Hence, by Khintchine's Theorem (see, e.g., [1, Theorem 1.2.3]), $a_{n+k}/a_n \to 1$ and $a_{n+k}(b_n - b_{n+k}) \to 0$ as $n \to \infty$. Conditions (6.9.3) and (6.9.4) are stronger versions of this.

*Proof.* Part (a) is a straightforward consequence of the Birkhoff ergodic theorem.

For (b) consider the case when $\ell \geq 0$, the other being similar. Then

$$\Phi_{0,n} - \Phi_{0,n+\ell} \leq \Phi_{0,n} - \Phi_{k,n+\ell} \leq \Phi_{0,n+\ell} - \Phi_{k,n+\ell}$$

Ergodicity implies a.e. convergence to zero if $\text{esssup}(\Phi)$ is finite. For the case in which $\text{esssup}(\Phi) = \infty$, we show that the left side (which is non-positive) and right side (which is non-negative) in the above inequalities converge in probability to zero. Let $\epsilon > 0$ and pick $a < \text{esssup}(\Phi)$ such that $\mu(\Phi_{0,k} \geq a) \leq \epsilon$ (this is possible because $\Phi$ is a.e. finite and $\operatorname{esssup}(\Phi) = \infty$). Then

$$\mu(a_n(\Phi_{0,n+\ell} - \Phi_{k,n+\ell}) \geq \epsilon) = \mu(\Phi_{0,k} \geq \Phi_{k,n+\ell} + \epsilon/a_n)$$
$$\leq \mu(\Phi_{0,k} \geq a) + \mu(\Phi_{0,k} < a, \Phi_{0,k} \geq \Phi_{k,n+\ell} + \epsilon/a_n)$$
$$\leq \mu(\Phi_{0,k} \geq a) + \mu(\Phi_{k,n+\ell} < a - \epsilon/a_n)$$
$$\leq \mu(\Phi_{0,k} \geq a) + \mu(\Phi_{k,n+\ell} < a)$$

which converges to $\mu(\Phi_{0,k} \geq a)$ by (a). The other side is dealt with similarly, using the stationarity of the process:

$$\mu(a_n(\Phi_{0,n+\ell} - \Phi_{0,n}) \geq \epsilon) \leq \mu(\Phi_{n,n+\ell} \geq a) + \mu(\Phi_{0,n} < a - \epsilon/a_n)$$
$$\leq \mu(\Phi_{n,n+\ell} \geq a) + \mu(\Phi_{0,n} < a).$$

Using ergodicity and the result of Eagleson [223, Thm. 6] the other statements in (b) follow.

Part (c) follows since $\hat{\Phi}_{1,n+1} - \hat{\Phi}_{0,n} \to 0$ on $(X^h, \mu^h)$ because $\Phi_{1,n+1} - \Phi_{0,n} \to 0$ on $(X, \mu)$ and $h \in L^1(\mu)$. Furthermore the invariant $\sigma$-field of $\{(\Phi \circ f^k) \circ \pi^h\}_k$ is the pull-back through $\pi^h$ of the invariant $\sigma$-field on $X$, hence it is still trivial, hence all the conditions for Eagleson's Theorem hold [223, Thm. 6] □

### 6.9.1
### Lifting to $X^h$

Mixing in distribution convergence allows to relate the extreme value laws for observations on $X$ to observations on $X^h$.

**Lemma 6.9.4.** *Let $f : (X, \mu) \to (X, \mu)$ be ergodic and $X^h$ the suspension space with a roof function $h \in L^1(\mu)$. Let $\Phi : X \to \mathbb{R}$ be an observation.*

*Set $\Phi_N(x) = \max_{k \leq N-1} \Phi(f^k(x))$ and define $\hat{\Phi}_N : X^h \to \mathbb{R}$ by $\hat{\Phi}_N(x, u) = \Phi_N \circ \pi^h(x, u) = \Phi_N(x)$.*

(a) *If $a_N(\Phi_N - b_N) \to_d G$ on $X$, then $a_N(\hat{\Phi}_N - b_N) \to_d G$ on $X^h$.*
(b) *If $1/h \in L^1(\mu)$ and $a_N(\hat{\Phi}_N - n_N) \to_d G$ on $X^h$, then $a_N(\Phi_N - b_N) \to_d G$ on $X$.*

*Proof.* By a result of Renyi [225], mixing convergence in distribution on $(X, \mu)$ implies convergence in distribution on $(X, \nu)$ whenever $\nu$ is absolutely continuous with respect to $\mu$.

For (a), take $\nu$ to be the probability measure on $X$ give by $d\mu := h/\overline{h}d\mu$. Then

$$\int_X \exp\{ita_N(\Phi_N - b_N)\} h/\overline{h} \, d\mu \to \mathbb{E}(e^{itG}), \qquad \forall t \in \mathbb{R}$$

because the convergence $a_N(\Phi_N - b_N) \to_d G$ is mixing in distribution. However,

$$\int_{X^h} \exp\{ita_N(\hat{\Phi}_N - b_N)\} \, d\mu^h = \frac{1}{\overline{h}} \int_X \int_0^{h(x)} \exp\{ita_N(\hat{\Phi}_N - b_N)\} \, du \, d\mu$$
$$= \int_X \exp\{ita_N(\Phi_N - b_N)\} h/\overline{h} \, d\mu.$$

Hence $\int_{X^h} \exp\{it a_N(\hat{\Phi}_N - b_N)\} \, d\mu^h \to \mathbb{E}(e^{itG})$, and the result follows.

For (b), note that the probability measure $d\nu^h := \overline{h}/h \, d\mu^h$ is an absolutely continuous probability with respect to $d\mu^h$, so we may repeat this argument.

$\square$

### 6.9.2
### The Normalization Constants

**Lemma 6.9.5.** *Assume that $g_n : X \to \mathbb{Z}$ are measurable functions such that $g_n(x)/n \to 0$ a.e. Let $S_n$ be an increasing sequence of random variables on $X$. If conditions* (6.9.3) *and* (6.9.4) *hold then*

$$a_n(S_n - b_n) \to G \iff a_n(S_{n+g_n(x)}(x) - b_n) \to G.$$

*Proof.* If we define $X_n = a_n(S_n - b_n)$ and $Y_n = a_n(S_{n+g_n} - b_n)$ then the lemma is a straightforward consequence of the inequality

$$X_n \leq Y_n \leq \frac{a_n}{a_{n+\epsilon n}} X_{a_{n+\epsilon n}} + a_n(b_{n+\epsilon n} - b_n) \tag{6.9.7}$$

valid on the set $\left|\frac{g_n}{n}\right| < \epsilon$.

$\square$

As a consequence:

*Remark 6.9.6*(a) If $\liminf_{n\to\infty} \mu(a_n(S_{n+g_n(x)}(x) - b_n) \leq v) \geq G(v)$ at each continuity point $v$ of $G$, then $\liminf_{m\to\infty} \mu(a_m(S_m - b_m) \leq v) \geq G(v)$ at each continuity point $v$ of $G$.

(b) If $\limsup_{n\to\infty} \mu(a_n(S_{n+g_n(x)}(x) - b_n) \leq v) \leq G(v)$ at each continuity point $v$ of $G$, then $\limsup_{m\to\infty} \mu(a_m(S_m - b_m) \leq v) \leq G(v)$ at each continuity point $v$ of $G$.

### 6.9.3
### The Lap Number

Now we relate the sequence of return times to the base $X$, namely the sum $h(x) + h(fx) + \ldots + h(f^n x)$ to the flow time $t$ via the strong law of large numbers. By the strong law of large numbers,

$$h_N = N\overline{h} + o(N) \text{ a.e.} \quad \text{as } N \to \infty \tag{6.9.8}$$

where $h_N(x) = h(x) + h(f(x)) + \ldots + h(f^{N-1}(x))$.

Given a time $T \geq 0$ define the *lap number* $N(x, T)$ by

$$h_{N(x,T)}(x) \leq T < h_{N(x,T)+1}(x). \tag{6.9.9}$$

The strong law of large numbers also implies that

$$\lim_{T\to\infty} N(x, T) = \infty \text{ a.e.} \tag{6.9.10}$$

and thus

$$\lim_{T\to\infty} \frac{T}{N(x, T)} = \overline{h} \text{ a.e.} \tag{6.9.11}$$

**Sublemma 6.9.7.** *For $\mu$ almost all $x \in X$ we have*

$$\lim_{T \to \infty} \frac{N(x, T + h(x))}{N(x, T)} = 1$$

*Proof.* Let $Z_a = \{x \in X : h(x) \leq a\}$, and given $\epsilon > 0$ and $T > 0$ let

$$X_{\epsilon, T} = \left\{ x \in X : \left| \frac{N(x, t + h(x))}{N(x, t)} - 1 \right| \geq \epsilon \text{ for some } t \geq T \right\}.$$

Then we have

$$\mu(X_{\epsilon, T}) \leq \mu\left( x \in X : \left| \frac{N(x, t + h(x))}{N(x, t)} - 1 \right| \geq \epsilon \text{ for some } t \geq T, \ h(x) \in [0, a] \right) + \mu(X \setminus Z_a).$$

Now for given $a > 0$ we have that $\mu$-a.e.

$$\frac{t}{N(x, t)} = \overline{h} + o(1), \quad \frac{t + a}{N(x, t + a)} = \overline{h} + o(1) \text{ as } t \to \infty,$$

and therefore $N(x, t + a)/N(x, t) \to 1$ almost surely as $t \to \infty$. Hence by taking $a$ arbitrarily large and then $T \to \infty$ it follows that $\mu(X_{\epsilon, T}) \to 0$. The result follows. $\qquad \square$

### 6.9.4
### Proof of Theorem 6.9.1.

The main observation is that for $(x, u) \in X^h$ with $0 \leq u < h(x)$

$$\Phi_{1, N(x, T)}(x) \leq \phi_T(x, u) \leq \Phi_{N(x, T + h(x)) + 1}(x) \tag{6.9.12}$$

(recall that $\Phi_{n, m}(x) := \max\{\Phi \circ f^k(x) \mid n \leq k < m\}$). Indeed, $\Phi_N(x) = \phi_{h_N(x)}(x, 0)$ for $x \in X$ and thus, taking into account the identifications of $X^h$:

$$\Phi_{1, N(x, T)}(x) = \max\{\phi(x, t) \mid h(x) \leq t < h_{N(x, T)}(x)\}$$
$$\phi_T(x, u) = \max\{\phi(x, t) \mid u \leq t < u + T\}$$
$$\Phi_{N(x, T + h(x)) + 1}(x) = \max\{\phi(x, t) \mid 0 \leq t < h_{N(x, T + h(x)) + 1}(x)\}.$$

The definition (6.9.9) of the lap number gives

$$h_{N(x, T)}(x) \leq T, \quad u + T < T + h(x) < h_{N(x, T + h(x)) + 1}(x),$$

and (6.9.12) follows.

We will also use that

$$\frac{N(x, T)}{\lfloor T/\overline{h} \rfloor} \to 1, \quad \frac{N(x, T + h(x)) + 1}{\lfloor T/\overline{h} \rfloor} \to 1 \text{ a.e. on } X \text{ (and hence on } X^h\text{)}, \tag{6.9.13}$$

which follow from (6.9.11) and Sublemma 6.9.7.

We first prove the implication

$$a_n(\Phi_n - b_n) \to G \text{ on } X \implies a_{T/\overline{h}}(\phi_T - b_{T/\overline{h}}) \to G \text{ on } X^h. \tag{6.9.14}$$

By Lemma 6.9.4(a), the left hand side of (6.9.14) implies

$$a_{\lfloor T/\bar{h}\rfloor}(\hat{\Phi}_{\lfloor T/\bar{h}\rfloor} - b_{\lfloor T/\bar{h}\rfloor}) \to G \text{ on } X^h, \tag{6.9.15}$$

where $\lfloor r \rfloor$ denotes the largest integer not exceeding $r$. By Lemma 6.9.2(b)

$$a_{\lfloor T/\bar{h}\rfloor}(\Phi_{1,\lfloor T/\bar{h}\rfloor} - \Phi_{\lfloor T/\bar{h}\rfloor}) \to 0 \text{ on } X$$

and, because $h \in L^1(\mu)$, this convergence in probability also holds on $X^h$ if we extend the function $\Phi_{1,N}$ to $\hat{\Phi}_{1,N} := \Phi_{1,N} \circ \pi^h$. Together with (6.9.15) this implies

$$a_{\lfloor T/\bar{h}\rfloor}(\hat{\Phi}_{1,\lfloor T/\bar{h}\rfloor} - b_{\lfloor T/\bar{h}\rfloor}) \to G \text{ on } X^h. \tag{6.9.16}$$

By Lemma 6.9.5, (6.9.13), (6.9.16) and (6.9.15) imply

$$a_{\lfloor T/\bar{h}\rfloor}(\hat{\Phi}_{1,N(x,T)}(x,u) - b_{\lfloor T/\bar{h}\rfloor}) \to G,$$

$$a_{\lfloor T/\bar{h}\rfloor}(\hat{\Phi}_{N(x,T+h(x))+1}(x,u) - b_{\lfloor T/\bar{h}\rfloor}) \to G \quad \text{on } X^h.$$

Use (6.9.12) to obtain

$$a_{\lfloor T/\bar{h}\rfloor}(\phi_T - b_{\lfloor T/\bar{h}\rfloor}) \to G \text{ on } X^h,$$

We now show the converse,

$$a_{\lfloor T/\bar{h}\rfloor}(\phi_T - b_{\lfloor T/\bar{h}\rfloor}) \to G \text{ on } X^h \implies a_n(\Phi_n - b_n) \to G \text{ on } X. \tag{6.9.17}$$

Denote by $\Omega \subset \mathbb{R}$ the continuity points of $G$.

Since (6.9.12) implies

$$\mu^h(a_{\lfloor T/\bar{h}\rfloor}(\hat{\Phi}_{1,N(x,T)}(x,u) - b_{\lfloor T/\bar{h}\rfloor}) \le v) \ge \mu^h(a_{\lfloor T/\bar{h}\rfloor}(\phi_T(x,u) - b_{\lfloor T/\bar{h}\rfloor}) \le v)$$

$$\mu^h(a_{\lfloor T/\bar{h}\rfloor}(\hat{\Phi}_{N(x,T+h(x))+1}(x,u) - b_{\lfloor T/\bar{h}\rfloor}) \le v) \le \mu^h(a_{\lfloor T/\bar{h}\rfloor}(\phi_T(x,u) - b_{\lfloor T/\bar{h}\rfloor}) \le v)$$

we obtain from the left hand side of (6.9.17) that

$$\liminf_{T\to\infty} \mu^h(a_{\lfloor T/\bar{h}\rfloor}(\hat{\Phi}_{1,N(x,T)}(x,u) - b_{\lfloor T/\bar{h}\rfloor}) \le v) \ge G(v),$$

$$\limsup_{T\to\infty} \mu^h(a_{\lfloor T/\bar{h}\rfloor}(\hat{\Phi}_{N(x,T+h(x))+1}(x,u) - b_{\lfloor T/\bar{h}\rfloor}) \le v) \le G(v), \qquad v \in \Omega.$$

By Remark 6.9.6 and (6.9.13), we conclude that

$$\liminf_{T\to\infty} \mu^h(a_{\lfloor T/\bar{h}\rfloor}(\hat{\Phi}_{1,\lfloor T/\bar{h}\rfloor} - b_{\lfloor T/\bar{h}\rfloor}) \le v) \ge G(v), \tag{6.9.18}$$

$$\limsup_{T\to\infty} \mu^h(a_{\lfloor T/\bar{h}\rfloor}(\hat{\Phi}_{\lfloor T/\bar{h}\rfloor} - b_{\lfloor T/\bar{h}\rfloor}) \le v) \le G(v), \qquad v \in \Omega. \tag{6.9.19}$$

Use that, by Lemma 6.9.2(b), $a_{\lfloor T/\bar{h}\rfloor}(\hat{\Phi}_{\lfloor T/\bar{h}\rfloor} - \hat{\Phi}_{1,\lfloor T/\bar{h}\rfloor}) \to 0$ on $X^h$ to deduce from the first relation above that

$$\liminf_{T\to\infty} \mu^h(a_{\lfloor T/\bar{h}\rfloor}(\hat{\Phi}_{\lfloor T/\bar{h}\rfloor} - b_{\lfloor T/\bar{h}\rfloor}) \le v) \ge G(v), \qquad v \in \Omega. \tag{6.9.20}$$

From (6.9.19) and (6.9.20) it follows that

$$a_{\lfloor T/\bar{h}\rfloor}(\hat{\Phi}_{\lfloor T/\bar{h}\rfloor} - b_{\lfloor T/\bar{h}\rfloor}) \to G \text{ on } X^h.$$

and Lemma 6.9.4(b) completes the proof. $\qquad\square$

## 6.10
## Physical Observables and Extreme Value Theory

In this section we discuss the extent to which we can obtain convergence to an EVL for physical observables. The reader is advised that such an issue is addressed also in Chap. 8 taking a more heuristic point of view aimed at providing results for high dimensional statistical mechanical systems.

While in Chap. 8 a POT point if view is taken and results are derived relating GPD parameters and geometrical property of the attractor, in this section we focus on the BM approach and aim at deriving the GEV description of extremes.

We start with Arnold's Cat Map which is a uniformly hyperbolic system with Lebesgue measure as the ergodic invariant measure. For this system, and for general observables we can obtain a functional form for $\tau(u)$ which is regularly varying. However, for systems with general SRB measures only bounds on $\tau(u)$ exist, at least for linear scaling sequences $u_n$. We discuss this situation for the solenoid map, the Lozi map, the Hénon map, and the Lorenz '63 model. In the following, the constant $\alpha$ will refer to the constant of regular variation of $\tau(u)$. For a generalized extreme value distribution, the tail index $\xi$ is precisely $-1/\alpha$.

### 6.10.1
### Arnold Cat Map

Let $\mathbb{T}^2 = \mathbb{R}^2 \mod 1$ be the 2 dimensional torus and let's consider the Arnold cat map given in Eq. 6.6.1. This system is Anosov and it has Lebesgue measure $\mu$ on the torus $\mathbb{T}^2$ as an invariant measure. With this example we want to study the role of the observable in determining extreme value laws. For this purpose we will consider $f$ as a map of $\mathbb{R}^2$ having the square $[0,1)^2$ as the invariant set. In other words, $\mathcal{X} = \mathbb{R}^2$ and $\Lambda = [0,1)^2$, hence $\Lambda$ is not an *attractor*, strictly speaking. The advantage is that this allows us to take functions of $\mathbb{R}^2$ as observables, rather than functions of $\mathbb{T}^2$. In this way, we can construct observables which are maximised at points in the interior or in the complement of $\Lambda$ and whose level sets have different shapes.

The main point of this section is that the value of the tail index is determined by the interaction between the shape of level sets (4.6.2) of the observable and the shape of the support of the invariant measure. To illustrate our ideas, and without attempting to cover all possible cases, we consider the following two observables $\phi_\gamma, \phi_{a,b} : \mathbb{R}^2 \to \mathbb{R}$

$$\phi_\gamma(x,y) = 1 - \text{dist}(p, \tilde{p})^\gamma, \quad \text{with } p = (x,y) \in \mathbb{R}^2. \tag{6.10.1}$$

$$\phi_{a,b}(x,y) = 1 - |x - \tilde{x}|^a - |y - \tilde{y}|^b, \tag{6.10.2}$$

where, given our focus on the Weibull case, we require $a, b, \gamma > 0$. Both observables are maximised at a point $\tilde{p} = (\tilde{x}, \tilde{y}) \in \mathbb{R}^2$. When $\tilde{p}$ is in the interior of $\Lambda$, observable (6.10.1) has the form so far analysed in the mathematical literature about extremes in dynamical systems, but we will also consider the case $\tilde{p} \notin \Lambda$. Observable (6.10.2) has been chosen to illustrate the effect of the shape of the level sets:

the level regions of (6.10.2) are not (Euclidean) balls unless $a = b = 2$, in which case (6.10.2) can be written as (6.10.1) for $\gamma = 2$.

The level regions $L^+(u)$ as defined in (4.6.2) are always balls for observable (6.10.1). However three (main) different situations occur, depending on the location of the point $\tilde{p}$ relative to the support of the invariant measure. We have the following result.

**Theorem 6.10.1.** *Consider the process $M_n = \max(X_1, \ldots, X_n)$ with $X_n = \phi_\gamma \circ f^{n-1}$, where $f$ is the map 6.6.1 and $\phi_\gamma : \mathbb{R}^2 \to \mathbb{R}$ is the observable in* (6.10.1). *Then for $\mu$-a.e. $\tilde{p} = (\tilde{x}, \tilde{y}) \in \mathbb{R}^2$, statement 4.6.5 holds with $\tau(u) = u^\alpha$, and we have the following cases:*

$$\alpha = 2\gamma^{-1}, \text{ for } \tilde{p} \in \Lambda; \tag{6.10.3}$$

$$\alpha = 3/2, \text{ for } \tilde{p} \notin \Lambda, \text{ with either } \tilde{y} \in (0, 1) \text{ or } \tilde{x} \in (0, 1); \tag{6.10.4}$$

$$\alpha = 2, \text{ for } \tilde{p} \notin \Lambda, \text{ with both } \tilde{x}, \tilde{y} \notin [0, 1]; \tag{6.10.5}$$

For observable (6.10.2) the shape of the level sets $L(u)$ depends on $a$ and $b$. For example, $L(u)$ has a convex elliptic-like shape when both $a, b > 1$, or an asteroid-like shape when both $a, b < 1$. Clearly various possibilities arise, depending on the geometry of the level sets, on whether the point $\tilde{p}$ is in the interior of $\Lambda$ and on the local geometry of $\Lambda$ near the extremal point $\hat{p} = (\hat{x}, \hat{y})$ with minimum distance from $\tilde{p}$.

**Theorem 6.10.2.** *Consider the process $M_n = \max(X_1, \ldots, X_n)$ with $X_n = \phi_\gamma \circ f^{n-1}$, where $f$ is the map (6.6.1) and $\phi_\gamma : \mathbb{R}^2 \to \mathbb{R}$ is the observable in* (6.10.2). *Then for $\mu$-a.e. $\tilde{p} = (\tilde{x}, \tilde{y}) \in \mathbb{R}^2$, statement 4.6.5 holds with $\tau(u) = u^\alpha$, and*

$$\alpha = \frac{1}{a} + \frac{1}{b} \quad \text{for } \tilde{p} \in \Lambda. \tag{6.10.6}$$

To prove Theorems 6.10.1-6.10.2 the main step is to determine the explicit sequence $u_n$ and functional form of $\tau(u)$ as defined in (4.6.4).

The verification of $\amalg_0(u_n), \amalg_0'(u_n)$ follows from the techniques of [138, 139] for this class of observables. Since the observable geometry is non-standard, we discuss briefly the idea of proof at the end of this subsection and point out the limitations. The main proof of 6.10.1 is contained in Lemmas 6.10.3, 6.10.4 and 6.10.5.

The proof of 6.10.2 is given in Lemma 6.10.6.

**Lemma 6.10.3.** *Suppose $\tilde{p} \in \text{int}(\Lambda) = (0, 1)^2$ and $\phi$ takes the form of (6.10.1), then $\alpha = 2/\gamma$.*

*Proof.* If $\tilde{p}$ is an interior point of $\Lambda$ then we see that

$$n\mu\{\phi(x, y) \geq u_n\} = n\mu\{d((x, y), \tilde{p}) \leq (1 - u_n)^{1/\gamma}\}$$
$$= C_\mu n (1 - u_n)^{2/\gamma}.$$

Thus the correct scaling laws are $a_n = C_\mu' n^{\gamma/2}, b_n = 1$ and $\tau(u) = (-u)^{2/\gamma}$. Here $C_\mu, C_\mu'$ are uniform constants. $\qquad \square$

**Lemma 6.10.4.** *Suppose $\tilde{p} \notin \overline{\Lambda} = [0,1]^2$ and $\phi$ takes the form of (6.10.1). If $\tilde{y} \in (0,1)$ or $\tilde{x} \in (0,1)$ then $\alpha = 3/2$.*

*Proof.* If $\tilde{p} \notin \overline{\Lambda}$ then there will exist a unique extremal point $\tilde{p} = (\tilde{x}, \tilde{y}) \in \overline{\Lambda}$ where $\phi(p)$ achieves its supremum with value $\tilde{u}$ as in (4.6.3). Since $\tilde{y} \in (0,1)$ or $\tilde{x} \in (0,1)$ then this point $\tilde{p}$ will not be a vertex of $\partial \Lambda$. The scaling $u_n$ will be chosen to so that

$$n\mu\{\phi(x,y) \geq u_n\} = n\mu\{p = (x,y) \in \Lambda : d(p,\tilde{p}) \leq (1-u_n)^{1/\gamma}\} \to \tau(u). \quad (6.10.7)$$

The middle term is no longer $\mathcal{O}(n(1-u_n)^{1/\gamma})$ since the level region that intersects $\Lambda$ is not a ball. We first of all set $u_n = u/a_n + \tilde{u}$ so that

$$\mu\{\phi(x,y) \geq u_n\} = \mu\left\{p = (x,y) \in \Lambda : (1-\tilde{u})^{1/\gamma} \leq d(p,\tilde{p}) \leq \left(1 - \tilde{u} - \frac{u}{a_n}\right)^{1/\gamma}\right\}. \quad (6.10.8)$$

To choose $a_n$ we first note that the level set $L(\tilde{u}^{1/\gamma})$ as defined in (4.6.2) is a circle that is tangent to $\partial \Lambda$ (since the extremal point $\tilde{p}$ is not a vertex). However the level set $L((\tilde{u} - \frac{u}{a_n})^{1/\gamma}))$ crosses $\partial \Lambda$ transversely (and is concentric to $L(\tilde{u}^{1/\gamma})$). To estimate the measure in Eq. (6.10.8), we suppose without loss of generality that $\tilde{y} \in (0,1)$ and $\tilde{x} > 1$. Hence

$$\tilde{u} = \sup_{(x,y) \in \Lambda} \{1 - d((x,y),(\tilde{x},\tilde{y}))^\gamma\} = 1 - |1 - \tilde{x}|^\gamma$$

Thus the measure (6.10.8) is of the order $\Delta x \Delta y$, where

$$\Delta x = (1 - \tilde{u} - u/a_n)^{1/\gamma} - (1 - \tilde{u})^{1/\gamma},$$

and

$$\Delta y = 2\left\{(1 - \tilde{u} - u/a_n)^{2/\gamma} - (1 - \tilde{u})^{2/\gamma}\right\}^{1/2}.$$

The former expression was obtained by solving $\phi(x,\tilde{y}) = u_n$ for $x$ and taking the difference of this (smaller) root with $\tilde{x}$, while the latter result for $\Delta y$ was obtained by solving $\phi(1,y) = u_n$ for $y$ and then taking the difference between the roots. If $u/a_n$ is sufficiently small, then a Taylor expansion implies that

$$\mu\{\phi(x,y) \geq u_n\} \approx \Delta x \Delta y \approx (u/a_n)^{3/2}, \quad (6.10.9)$$

where $\approx$ means *equal to* up to a uniform multiplication constant. Setting $a_n = n^{2/3}$ implies that $\tau(u) = \mathcal{O}((-u)^{3/2})$. $\qquad \square$

**Lemma 6.10.5.** *Suppose $\tilde{p} \notin \Lambda$ and $\phi$ takes the form of (6.10.1). If both $\tilde{x}, \tilde{y} \notin [0,1]$ then $\xi = -1/2$.*

*Proof.* Without loss of generality we consider $\tilde{p} = (\tilde{x}, \tilde{y})$ defined so that $\tilde{x} = 1 + \lambda \cos\theta$, $\tilde{y} = 1 + \lambda \sin\theta$ for $\lambda > 0$ and $\theta \in (0, \pi/2)$. For such values of $(\tilde{x}, \tilde{y})$, the corner point $(1,1) \in \partial \Lambda$ will always maximise $\phi$. The proof is identical to Lemma 6.10.4 except that the level sets are not tangent to $\partial \Lambda$ at $(1, 1)$. Setting $u_n = u/a_n + \tilde{u}$ and arguing as in the proof of Lemma 6.10.4, we obtain:

$$\tilde{u} = \sup_{(x,y) \in \Lambda} \{1 - d((x,y),(\tilde{x},\tilde{y}))^\gamma\} = 1 - ((\tilde{x}-1)^2 - (\tilde{y}-1)^2)^{\gamma/2},$$

$$\Delta x = 1 - \tilde{x} + \{(1 - \tilde{u} - u/a_n)^{2/\gamma} - (1 - \tilde{y})^2\}^{1/2},$$

$$\Delta y = 1 - \tilde{y} + \{(1 - \tilde{u} - u/a_n)^{2/\gamma} - (1 - \tilde{x})^2\}^{1/2}.$$

Again, if $u/a_n$ is sufficiently small, then a Taylor expansion implies that

$$\mu\{\phi(x,y) \geq u_n\} \approx \Delta x \Delta y \approx (u/a_n)^2, \tag{6.10.10}$$

and hence setting $a_n = n^{1/2}$ implies that $\tau(u) = \mathcal{O}((-u)^2)$. ☐

This concludes the proof of Theorem 6.10.1. For the proof of Theorem 6.10.2 we have the following lemma.

**Lemma 6.10.6.** *Suppose that $\tilde{p} \in \mathrm{int}(\Lambda) = (0,1)^2$ and $\phi$ takes the form of (6.10.2). Then for $u \lesssim 1$ we have that $\mathrm{Leb}(L(u)) = C(1-u)^{\frac{1}{a}+\frac{1}{b}}$ for some $C_0 \leq C \leq 4$ where $C_0 > 0$.*

*Proof.* Let $u = 1 - \varepsilon$. For $\varepsilon$ sufficiently small the level region can be written as

$$L^+(u) = \{(x,y) \in \mathrm{int}(\Lambda) : |x|^a + |y|^b \leq \varepsilon\}. \tag{6.10.11}$$

The area of this set is bounded from above by the area of a rectangle of sides $2\varepsilon^{1/a}$ and $2\varepsilon^{1/b}$. Also, for any $q \in (0,1)$, the area of the set is bounded from below by that of a rectangle of sides $2q^{1/a}\varepsilon^{1/a}$ and $2(1-q)^{1/b}\varepsilon^{1/b}$, so we can choose $C_0 = 4q^{1/a}(1-q)^{1/b}$. ☐

Hence if $(\tilde{x}, \tilde{y}) \in \mathrm{int}(\Lambda)$, we see that (for uniform $C_\mu > 0$),

$$n\mu\{p = (x,y) : \phi(p) \geq u_n\} \to C_\mu(-u)^{\frac{1}{a}+\frac{1}{b}} \quad \text{with} \quad a_n = n^{\frac{ab}{a+b}}, b_n = 1. (6.10.12)$$

We now explain how to check $\mathcal{D}_0'(u_n)$ and $\mathcal{D}_0(u_n)$ for the observables given in (6.10.2). The other scenarios are similar. The methods used in [138] are primarily geared towards observables that are expressed as functions of distance. In our situation the observables are not given explicitly in this form, but they do have a bounded geometry in the sense that the level set $\{\phi(p) = \epsilon\}$ can be circumscribed by a ball of radius $\epsilon^{d'}$ with $d' = \max\{a^{-1}, b^{-1}\}$. This fact is useful when checking the $\mathcal{D}_0'(u_n)$ and $\mathcal{D}_0(u_n)$ conditions.

More specifically, to check $\mathcal{D}_0(u_n)$ following [138] it suffices to show that for (fixed) $r > 0$, $\mu\{r \leq \phi(p) \leq r + \delta\}$ is bounded by a power of $\delta$ as $\delta \to 0$. By a simple integration calculation, this estimate holds for the observables (6.10.2). Another property required is that the system has exponential decay of correlations. This is property holds by uniform hyperbolicity.

Checking $D'(u_n)$ is generally harder, but the main estimate, see [138, 139, 72] involves a control of the measure of the set $\{p : d(p, f^j(p)) \leq a_n\}$ for some specific sequence $a_n \to 0$, usually power law in $n$ with $j = o(n)$. The aim is then to show that the measure of this set goes to zero at a rate bounded by a power law in $n$. When the observables are functions of distance, then a natural choice for $a_n$ is $n^{-1/d}$. This follows by observing that (by choice of $u_n$)

$$\{\phi(p) \geq u_n\} \subset \{p : d(p, \tilde{p}) \leq \frac{1}{n^{1/d}}\}. \tag{6.10.13}$$

When the observables have the form given in (6.10.2) then a similar statement holds, but for the right hand set in (6.10.13) we instead have $d(p, \tilde{p}) = \mathcal{O}(n^{-\gamma})$ for some $\gamma = \gamma(a, b)$. This relationship is sufficient to allow the methods, such as the *maximal function technique* utilized in [72] to be applied to this class of observables.

However we remark that the methods have limitations, and at present the arguments do not directly extend to more exotic observables, such as

$$\phi(x, y) = \exp\left\{-\frac{1}{|x - \tilde{x}|}\right\} + |y - \tilde{y}|^c. \tag{6.10.14}$$

This observable has the property that the Lebesgue measure of the level set $\{\phi(p) = \epsilon\}$ is $\mathcal{O}\left(\epsilon^{1/c}/\log(1/\epsilon)\right)$, but any circumscribing ball must have radius at least $\mathcal{O}\left(1/\log(1/\epsilon)\right)$. Thus, if an extreme value law is to be proved for this observable then the methods of [138, 139] would need to be adapted to situations where the $a_n$ have sub-polynomial asymptotics.

## 6.10.2
### Uniformly HyperbolicHyperbolic Attractors: the Solenoid Map

Consider the solid torus as the product of $\mathbb{T} = \mathbb{R}/\mathbb{Z}$ times the unit disc in the complex plane $\mathbb{D}_R = \{z \in \mathcal{C} \,|\, |z| < R\}$, for some $R$ with $0 < R < 1$. Then the solenoid map is defined as follows:

$$\begin{aligned} f_\lambda : \quad \mathbb{T} \times \mathbb{D}_R \quad &\rightarrow \quad \mathbb{T} \times \mathbb{D}_R \\ (\psi, w) \quad &\mapsto \quad \left(2\psi, \lambda w + K e^{i2\pi\psi}\right). \end{aligned} \tag{6.10.15}$$

In order to have the map well defined we need $K + \lambda R < R$ and $\lambda R < K$. For our purposes it is convenient to have the torus embedded in $\mathbb{R}^3$. Consider Cartesian coordinates $(x, y, z) \in \mathbb{R}^3$ and define corresponding cylindrical coordinates $r, \psi, z$ by $x = r\cos(\psi)$ and $y = r\sin(\psi)$. Then the torus of width $R$ can be identified with the set $D = \{(r-1)^2 + z^2 \leq R^2\}$ for $R < 1$. The torus $\mathbb{T} \times \mathbb{D}_R$ (with coordinates $(\psi, u + iv)$) can be identified with $D$ taking $r = 1 + u$ and $z = v$. We thus obtain an embedded solenoid map

$$g_\lambda : D \to D, \quad g_\lambda(\psi, r, z) = (2\psi, 1 + K\cos(\psi) + \lambda(r-1), K\sin(\psi) + \lambda z). \tag{6.10.16}$$

The solenoid attractor is defined as the attracting set of the map $g_\lambda$:

$$\Lambda = \bigcap_{j \geq 1} g_\lambda^j(D).$$

For $\lambda < \frac{1}{2}$ we have

$$\dim_H(\Lambda) = 1 + \frac{\log 2}{\log \lambda^{-1}}, \tag{6.10.17}$$

where $\dim_H$ denotes the Hausdorff dimension [226]. We consider the following observables $\phi_\gamma, \phi_{abcd} : \mathbb{R}^3 \to \mathbb{R}$:

$$\phi_\gamma(x,y,z) = 1 - \mathrm{dist}(p,\tilde{p})^\gamma, \quad \text{with } p = (x,y,z) \in \mathbb{R}^3, \tag{6.10.18}$$

$$\phi_{abcd}(x,y,z) = ax + by + cz + d, \tag{6.10.19}$$

Observable (6.10.18) is maximised at a point $\tilde{p} \in \mathbb{R}^3$, whereas (6.10.19) is unbounded in the phase space $\mathbb{R}^3$ (except for the trivial choice $a = b = c = 0$).

For the case of observable (6.10.18) and when $\tilde{p} \in \Lambda$, results on convergence of $M_n$ to an extreme value distribution where discussed in Section 6.6, and indeed only bounds on the convergence could be achieved since $\mu$ does not admit a density with respect to Lebesgue measure on $\mathbb{R}^3$.

More interesting considerations arise for the observable (6.10.19), and similarly for observable (6.10.18) when $\tilde{p} \notin \Lambda$. In the following we focus on observable (6.10.19). The case for obsevable (6.10.18) turn out to be similar, since geometrically the level sets are smooth curves, *i.e.* their tangent vectors vary continuously. As a simple case, consider first the degenerate solenoid with $\lambda = 0$ and take a planar observable $\phi := ax + by + d$, thus reducing the problem to the $(x,y)$-plane. In this case we have the trivial dimension formula $\dim_H(\Lambda) = 1$ since $\Lambda$ is a circle. However, for computing the tail index we lose a factor of $1/2$ due to the geometry of the level set. Indeed, level sets are straight lines within the $(x,y)$-plane, and at the extremal point $\tilde{p} = (\tilde{x}, \tilde{y})$ the critical level set $L(\tilde{u})$ is tangent to $\Lambda$. Since the tangency is quadratic, we find that

$$\mu(L^+(\tilde{u} - \epsilon)) = m_{\gamma^u}\{\gamma^u(\tilde{p}) \cap L^+(\tilde{u} - \epsilon)\} = \mathcal{O}(\sqrt{\epsilon}). \tag{6.10.20}$$

Here $\gamma^u(\tilde{p})$ is the unstable manifold through $\tilde{p}$ (*i.e.* it is the unit circle), and $m_{\gamma^u}$ is the one-dimensional conditional (Lebesgue) measure on $\gamma^u(\tilde{p})$. Hence $\tau(u) = u^\alpha$ with

$$\alpha = \dim_H(\Lambda) - \frac{1}{2} = \frac{1}{2}.$$

The mechanism described above is similar to that for the cat map, leading to formula (6.10.4): indeed, there we have $\dim_H(\Lambda) = 2$, yielding the value $3/2$ for the constant $\alpha$.

For $\lambda > 0$, the attractor has more complicated geometry and is locally the product of a Cantor set with an interval [226]. Planar cross sections that intersect $\Lambda$ transversely form a Cantor set of dimension $\dim_H(\Lambda) - 1 = -\log 2/\log \lambda$. To calculate $\mu(L^+(\tilde{u} - \epsilon))$ we would like to repeat the calculation above using Eq. (6.10.20), but now the set of unstable leaves that intersect $L^+(\tilde{u} - \epsilon)$ form a Cantor set (for each $\epsilon > 0$). The extremal point $\tilde{p}$ where $\phi(p)$ attains its maximum on $\Lambda$ forms a *tip* of $\Lambda$ relative to $L(\tilde{u})$. Such a tip corresponds to a point on $\tilde{p} \in \Lambda$ whose unstable segment $\gamma^u(\tilde{p})$ is tangent to $L(\tilde{u})$ at $\tilde{p}$, and moreover normal to $\nabla\phi(\tilde{p})$ at $\tilde{p}$. Given $\epsilon > 0$, we (typically) expect to find a Cantor set of values $t \in [0, \epsilon]$ for which the level sets $L(\tilde{u} - t)$ are tangent to some unstable segment $\gamma^u \subset \Lambda$. For other values of $t$ these level sets cross the attractor transversely. Given (fixed) $\epsilon_0 > 0$ we can define the *tip set* $\Gamma \equiv \Gamma(\epsilon_0) \subset \Lambda$ as follows: let $T_p\gamma^u(p)$ be the tangent space to $\gamma^u$ at $p$. Then we define

$$\Gamma = \{p \in L^+(\tilde{u} - \epsilon_0) \cap \Lambda : T_p\gamma^u(p) \cdot \nabla\phi(p) = 0\}. \tag{6.10.21}$$

This tip set plays a role in proving the following result, which in turn provides us with information on the form of the tail index $\xi$.

**Proposition 6.10.7.** *Suppose that $g_\lambda$ is the map (6.10.16) and $\phi = \phi_{abcd}$. Define $\tau(\epsilon) = \mu\{\phi(p) \geq \tilde{u} - \epsilon\}$. If $\dim_H(\Gamma) < 1$, then modulo a zero measure set of values $(a, b, c, d)$,*

$$\lim_{\epsilon \to 0} \frac{\log \tau(\epsilon)}{\log \epsilon} = 1/2 + \dim_H(\Gamma).$$

We give a proof below. Based on this proposition it is possible to bound the limiting behaviour of $\mu\{M_n \leq u_n\}$ relative to the function $e^{-\tau(u)}$, with $\tau(u) = u^\alpha$, and

$$\alpha = \dim_H(\Lambda) - \frac{1}{2} = \frac{1}{2} + \frac{\log 2}{\log \lambda^{-1}}. \tag{6.10.22}$$

We outline the main technical steps required. Firstly, conditions $D_0(u_n)$ and $D_0'(u_n)$ should be checked. We believe that this should follow from [138], however the proof would be non-standard due to the level set geometry. Secondly, we claim that $\dim_H(\Gamma) = \dim_H(\Lambda) - 1$. The proof of this requires estimates on the regularity of the holonomy map between stable disks, such as Lipschitz regularity.

*Proof of Proposition* 6.10.7. For each $\epsilon < \epsilon_0$, consider the set $\Gamma(\epsilon) \subset \Gamma \cap L^+(\tilde{u} - \epsilon)$. Then for each $p \in \Gamma(\epsilon)$, there exists $t < \epsilon$ such that $\gamma^u(p)$ is tangent to $L(\tilde{u} - t)$. If the observable $\phi$ takes the form of (6.10.19), then by the same calculation as (6.10.20) we obtain

$$m_{\gamma^u}\{\gamma^u(p) \cap L^+(\tilde{u} - \epsilon)\} = \mathcal{O}(\sqrt{\epsilon - t}). \tag{6.10.23}$$

Thus to compute $\mu(L^+(\tilde{u} - \epsilon))$, we integrate (6.10.23) over all relevant $t < \epsilon$ using the measure $\mu_\Gamma$, which is the measure $\mu$ conditioned on $\Gamma$. Provided $\dim_H(\Gamma) < 1$, the projection of $\Gamma$ onto the line in the direction of $\nabla\phi$ is also a Cantor set of the same dimension for typical (full volume measure) $(a, b, c, d)$, see [226]. Thus the set of values $t$ corresponding to when $L(\tilde{u} - t)$ is tangent to $\Gamma$ form a Cantor set of dimension $\dim_H(\Gamma)$. If $\pi$ is the projection from $\Gamma$ onto a line in the direction of $\nabla\phi$, then the projected measure $\pi_*\mu_\Gamma$ has local dimension $\dim_H(\Gamma)$ for typical $(a, b, c, d)$. We have

$$\mu(L^+(\tilde{u} - \epsilon)) = \int_0^\epsilon \int_{L^+(\tilde{u} - \epsilon)} dm_{\gamma^u} d\mu_\Gamma. \tag{6.10.24}$$

To estimate this integral we bound it above via the inequality $m_{\gamma^u}(\gamma^u \cap L^+(\tilde{u} - \epsilon)) \leq C\sqrt{\epsilon}$, and bound it below using the fact that for $t > \epsilon/2$, $m_{\gamma^u}(\gamma^u \cap L^+(\tilde{u} - \epsilon)) \geq C\sqrt{\epsilon}$. Here $C > 0$ is a uniform constant. Putting this together we obtain for typical $(a, b, c, d)$

$$\mu(L^+(\tilde{u} - \epsilon)) = \int_0^\epsilon \int_{L^+(\tilde{u}-\epsilon)} dm_{\gamma^u} d\mu_\Gamma = \sqrt{\epsilon} \cdot \epsilon^{\dim_H(\Gamma) + \delta(\epsilon)}$$
$$= \epsilon^{1/2 + \dim_H(\Gamma) + \delta(\epsilon)}, \qquad (6.10.25)$$

$$\square$$

where $\delta(\epsilon) \to 0$ as $\epsilon \to 0$. The constant $\delta(\epsilon)$ comes from the definition of dimension of an SRB measure, namely that for $\mu$-a.e. $x \in \Lambda$, $\log \mu(B(x, \epsilon))/\log \epsilon \to \dim_H(\mu)$.

## 6.11
## Non-uniformly Hyperbolic Examples: the Hénon and Lozi maps

We here consider the Hénon map and the Lozi map given in Eqs. 6.6.6 and 6.6.5, respectively, and the observables

$$\phi_\gamma(x, y) = 1 - \text{dist}(p, \tilde{p})^\gamma, \quad \text{with } p = (x, y) \in \mathbb{R}^2, \qquad (6.11.1)$$

$$\phi_\theta(x, y) = x \cos(2\pi\theta) + y \sin(2\pi\theta), \qquad (6.11.2)$$

where $\gamma > 0$ and $\theta \in [0, 2\pi]$ are parameters and $\tilde{p}$ is a point in $\mathbb{R}^2$. Following the discussion for the solenoid map, and defining as before $\tau(\epsilon) = \mu\{\phi(p) \geq \tilde{u} - \epsilon\}$, we could conjecture that for $\phi = \phi_\theta$, or for $\phi = \phi_\gamma$ with $\tilde{p} \notin \Lambda$:

$$\lim_{\epsilon \to 0} \frac{\log \tau(\epsilon)}{\log \epsilon} = \dim_H(\Lambda) - \frac{1}{2}. \qquad (6.11.3)$$

For the Hénon map under observable (6.10.1), and in view of the results of a recent paper [227], it is expected that formula (6.11.3) holds for so-called *Benedicks-Carleson parameter values* [134]. Such parameter values, however, are obtained by a perturbative argument near $(a, b) = (2, 0)$, where the bound on the smallness of $b$ is not explicit. Moreover, the parameter exclusion methods used to define the Benedicks-Carleson parameter values are not constructive. For these reasons, it is not possible to say whether Benedicks-Carleson behaviour is also attained at the classical parameter values $(a, b) = (1.4, 0.3)$.

For planar observables, we again study the *tip set* $\Gamma \subset \Lambda$ as defined for the Solenoid map, namely, for fixed $\epsilon_0 > 0$ and $p = (x, y)$, let

$$\Gamma = \{p \in L^+(\tilde{u} - \epsilon_0) \cap \Lambda : T_p \gamma^u(p) \cdot \nabla\phi(p) = 0\}, \qquad (6.11.4)$$

and for each $\epsilon < \epsilon_0$, consider the set $\Gamma(\epsilon) \subset \Gamma \cap L^+(\tilde{u}-\epsilon)$. Then for each $p \in \Gamma(\epsilon)$, there exists $t < \epsilon$ such that $\gamma^u(p)$ is tangent to $L^+(\tilde{u} - t)$. For the planar observable $\phi$ we would expect to obtain (as with the solenoid):

$$m_{\gamma^u}\{\gamma^u(p) \cap L^+(\tilde{u} - \epsilon)\} = \mathcal{O}(\sqrt{\epsilon - t}), \qquad (6.11.5)$$

where $m_{\gamma^u}$ is the conditional (Lebesgue) measure on the one-dimensional unstable manifold. However in this calculation we have assumed that the tangency between $\gamma^u(p)$ and $L(\tilde{u}-t)$ is quadratic, and that the unstable segment is sufficiently long so as to cross $L(\tilde{u}-\epsilon)$ from end to end. For the Hénon map both of these conditions can fail. In particular, the Hénon attractor admits a critical set of folds that correspond to points where the attractor curvature is large. More precisely the critical set is formed by homoclinic tangency points between stable and unstable manifolds. This set has zero measure, but it is dense in the attractor. Furthermore the attractor has complicated geometry, where local stable/unstable manifolds can fold back and forth upon themselves. However, the regions that correspond to these folds (of high curvature) occupy a set of small measure. See [228] and references therein for a more detailed discussion.

As discussed for the solenoid, for linear scaling sequences $u_n$ we expect not to have convergence to an EVL. However similar to [138], it is possible to obtain bounds on limiting behaviour of $\mu\{M_n \leq u_n\}$ relative to the function $e^{-\tau(u)}$, with $\tau(u) = u^\alpha$. To estimate $\alpha$ we conjecture to have the following formula:

$$\lim_{\epsilon \to 0} \frac{\log \tau(\epsilon)}{\log \epsilon} = \dim_H(\mu) - \frac{1}{2} \tag{6.11.6}$$

where $\mu$ is the SRB measure for the Hénon map (at Benedicks-Carleson parameters). This would follow from the estimate:

$$\mu(L^+(\tilde{u}-\epsilon)) = \int_0^\epsilon \int_{L^+(\tilde{u}-t)} dm_{\gamma^u} d\mu_\Gamma = \sqrt{\epsilon} \cdot \epsilon^{\dim_H(\Gamma)+\delta(\epsilon)}, \tag{6.11.7}$$

where the factor of $\sqrt{\epsilon}$ comes from Eq. (6.10.23) and $\delta(\epsilon) \to 0$ as $\epsilon \to 0$. To obtain Eq. (6.11.6), we would need to show that $\dim_H(\Gamma) = \dim_H(\mu)-1$. This is perhaps harder to verify and it will depend on the regularity of the holonomy map taken along unstable leaves. Finally we would project this set onto a line in the direction of $\nabla\phi(p)$, and typically the projection would preserve the dimension.

## 6.11.1
## Extreme Value Statistics for the Lorenz '63 Model

In this section we examine what results are known about recurrence and extremes for the classical Lorenz '63 model [105]:

$$\dot{x} = \sigma(y - x),$$
$$\dot{y} = x(\rho - z) - y, \tag{6.11.8}$$
$$\dot{z} = xy - \beta z.$$

These equations provide the time evolution of three modes of the temperature and streamfunction obtained by performing a severe truncation of the full dynamics of the Rayleigh-Bènard convection problem in two dimensions [229]. Here $\sigma$ is the Prandtl and $\rho$ the Rayleigh number, while $\beta$ is a geometrical factor taking into account the aspect ratio of the modes in physical space. We set $\sigma = 10$, $\beta = 8/3$ and $\rho =$ 28, which is a fairly common choice in the vast literature on the Lorenz model, see e.g. [211, 230]. A graphical representation of the Lorenz attractor for this choice of the parameters' value is given in Fig. 6.3. See also [231, 232] for a discussion of the properties of generalised Lorenz-like models. Considering the Lorenz '63 model (6.11.8) we recall some geometrical facts of the Poincaré map to $z = $ constant sections. Given the planar sections $\Sigma = \{(x, y, 1) : |x|, |y| \leq 1\}$, and $\Sigma' = \{(1, y, z) : |y|, |z| \leq 1\}$, the map $P : \Sigma \to \Sigma$ decomposes as $P = P_2 \circ P_1$, where $P_1 : \Sigma \to \Sigma'$ and $P_2 : \Sigma' \to \Sigma$. To describe the form of $P$, let $\beta = |\lambda_s|/\lambda_u$, $\beta' = |\lambda_{ss}|/\lambda_u$, where $\lambda_s$, $\lambda_{ss}$ and $\lambda_u$ are the eigenvalues of the linearised Lorenz flow at the origin, with $\lambda_s = -8/3$, $\lambda_{ss} = -22.83$ and $\lambda_u = 11.83$ for our choice of parameters. Then it can be shown that $P_1(x, y, 1) = (1, x^{\beta'} y, x^{\beta})$, and $P_2$ is a diffeomorphism. Thus the rectangle $\Sigma^+ = \{(x, y, 1) : x > 0, |y| \leq 1\}$ gets mapped into a region $P_1(\Sigma^+)$ with a cusp at $y = 0$. The cusp boundary can be represented as the graph $|y| = z^{\beta'/\beta} \approx z^8$. The flow has a strong stable foliation, and we form the quotient space $\widehat{\Sigma} = \Sigma/\sim$ by defining an equivalence relation $p \sim q$ if $p \in \gamma^s(q)$, for a stable leave $\gamma^s$. Hence the map $P : \Sigma \to \Sigma$ can be reduced to a uniformly expanding one-dimensional map $f : \widehat{\Sigma} \to \widehat{\Sigma}$, with a derivative singularity at $x = 0$. Here $\widehat{\Sigma}$ identified with $[-1, 1]$, and $f'(x) \approx |x|^{\beta-1}$ near $x = 0$. The extreme statistics associated to such a map was discussed in Section 6.4.3.

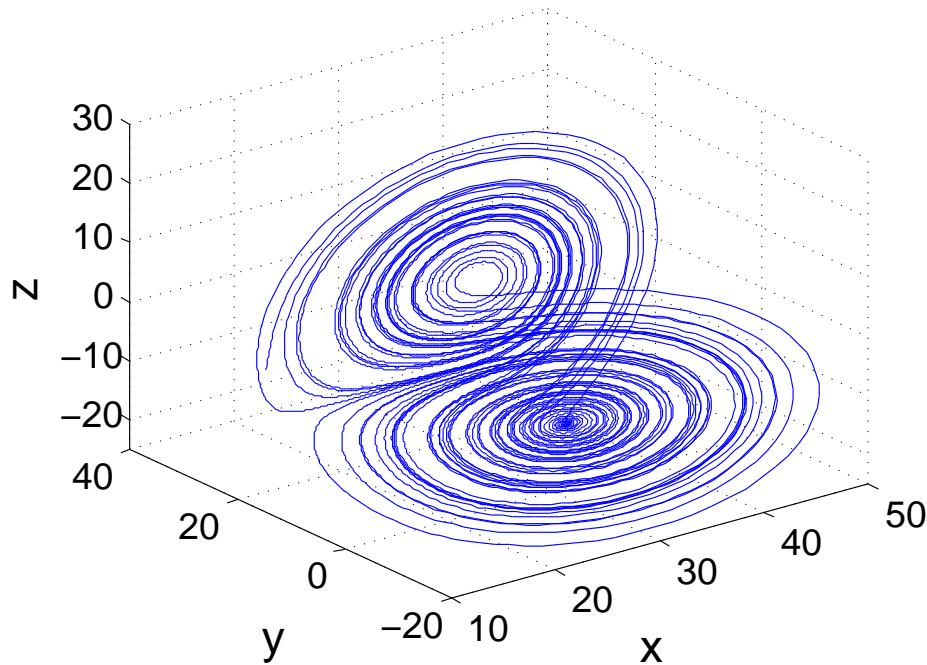

**Figure 6.3** Numerical approximation to the attractor of the Lorenz '63 model given in Eq. 6.11.8, where the classical values for the parameters parameters $\sigma = 10$, $\beta = 8/3$ and $\rho = 28$ are used.

The Lorenz flow admits an SRB measure $\mu$ which can be written as $\mu = \mu_P \times \text{Leb}$ (up to a normalisation constant). The measure $\mu_P$ is the SRB measure associated to the Poincaré map $P$, and the local dimension is defined (and constant) $\mu$-a.e., see [233]. Using the existence of the stable foliation, and the SRB property of $\mu$, we can write $\mu_P$ as the (local) product $\mu_{\gamma^u} \times \mu_{\gamma^s}$ where $\mu_{\gamma^u}$ is the conditional measure on unstable manifolds, and $\mu_{\gamma^s}$ is the conditional measure on stable manifolds. We can identify each measure $\mu_{\gamma^u}$ (via a holonomy map) with that of the invariant measure $\mu_f$ associated to $f$. The measure $\mu_f$ is absolutely continuous with respect to Lebesgue measure, but it has zero density at the endpoints of $\widehat{\Sigma}$, that is

$$\mu_f([1-\epsilon, 1]) \approx \epsilon^{1/\beta} \approx \epsilon^{4.4} \quad \text{as } \epsilon \to 0. \tag{6.11.9}$$

Results on extreme statistics for the two dimensional map $P : \Sigma \to \Sigma$ and the Lorenz model were established recently in [234]. The proof required verification of conditions (H1s), (SRT1), and (H3) as described in Section 6.5. Using our discussion in Section 6.5 the result of [234] can be expressed as follows:

**Proposition 6.11.1.** *Let $P : \Sigma \to \Sigma$ be the Poincaré return map associated to the Lorenz model for the classical parameters, and let $M_n$ denote the associated maximal process. Then for $\mu_P$ a.e. $\tilde{x} \in \Sigma$*

$$\left| \mu\{M_n \le u_n\} - G_{\sqrt{n}}(u) \right| \le C_1 \frac{(\log n)^{1+\epsilon}}{\sqrt{n}} + C_2 \frac{(\log n)^{1+\epsilon}}{n^\alpha} \tag{6.11.10}$$

*where $C_1$, $C_2 > 0$ are constants independent of $n$, but dependent on $\tilde{x}$.*

*Remark* 6.11.2. The situation here is similar to that for the Lozi map, and the Hénon map in the sense that we cannot precise estimates on convergence rates due to the irregularity of the measure. We do not know the precise scaling constants $u_n$. For the observable $\phi(x) = -\log \text{dist}(x, \tilde{x})$ we have that for all $\epsilon > 0$, $\lim_{n\to\infty} \mu(M_n \le (1-\epsilon)(\log n + v)/d)) \le e^{-e^{-v}} \le \lim_{n\to\infty} \mu(M_n \le (1+\epsilon)(\log n + v)/d))$ which provides an estimate of the correct sequence $u_n$.

Using methods of Section 6.9, extreme statistics for the corresponding Lorenz model can be deduced. This is discussed also in [234] and the following is achieved.

**Proposition 6.11.3.** *Assume that $P : \Sigma \to \Sigma$ is the Lorenz map and $f_t$ is the corresponding Lorenz flow. Assume there is a sequence $u_n := u_n(u)$ such that*

$$n\mu\{\Phi_0 \le u_n\} \to \tau(u),$$

*and suppose this sequence gives rise to normalizing constants $a_n > 0$ and $b_n$ satisfying:*

$$\lim_{\epsilon \to 0} \limsup_{n\to\infty} \, a_n |b_{[n+\epsilon n]} - b_n| = 0, \tag{6.11.11}$$

$$\lim_{\epsilon \to 0} \limsup_{n\to\infty} \left| 1 - \frac{a_{[n+\epsilon n]}}{a_n} \right| = 0, \tag{6.11.12}$$

*and*

$$a_n(\Phi_n - b_n) \to e^{-\tau(u)}. \tag{6.11.13}$$

*Then*

$$a_{\lfloor T/\bar{h} \rfloor}(\phi_T - b_{\lfloor T/\bar{h} \rfloor}) \to e^{-\tau(u)}. \tag{6.11.14}$$

*Remark* 6.11.4. The above result assumes the existence of a linear scaling sequence $u_n$. In general the observable $\Phi_0$ would not vary regularly if $n\mu\{\Phi_0 \leq u_n\} \to \tau(u)$. If $\Phi_0$ is regularly varying, then a similar statement to Remark 6.11.2 applies.

To end, we note that for physical observations corresponding extreme value laws depend on the geometry of the level sets, and also on the geometry of the attractor. For a planar observable, e.g. $\phi(x) = ax + by + cz + d$, it was conjectured in [81] that

$$\lim_{\epsilon \to 0} \frac{\log \tau(\epsilon)}{\log \epsilon} = \frac{1}{\beta} + \frac{1}{2} + \tilde{d}_s. \tag{6.11.15}$$

The constant $\tilde{d}_s$ comes from the dimension of $\mu_s$ which is (numerically) seen to be small due to the strong stable foliation. In contrast with the solenoid and Hénon maps, the tail index associated to this planar observable comes from an estimate of the measure of $\mu(L^+(\tilde{u} - \epsilon))$ which we assume scales as the product of the three factors: $\sqrt{\epsilon} \cdot \epsilon^{d_u} \cdot \epsilon^{8d_s}$. Here the factor $\sqrt{\epsilon}$ comes from the measure $\mu$ conditioned on $\Lambda \cap L^+(\tilde{u} - \epsilon)$ in the (central)-flow direction, while the factor $\epsilon^{d_u}$ comes from the $\mu_P$-measure conditioned on unstable manifolds that terminate at the cusp. In a generic case we would expect $d_u = 1$. However, since we are near the cusp (namely near the boundary $\partial\hat{\Sigma}$) we have $d_u = 1/\beta = 4.4$ due to the zero in the density of $\mu_f$, see (6.11.9). Finally we have a contributing factor $\epsilon^{8d_s}$ that comes from the the strength of the cusp at $P(\partial\Sigma)$, with $d_s$ the local dimension of $\mu_{\gamma^s}$.

# 7
# Extreme Value Theory for Randomly Perturbed Dynamical Systems

## 7.1
## Introduction

The aim of this chapter is to show that EVT can be defined and developed for dynamical systems that are randomly perturbed. Two kinds of perturbations will be considered: *random transformations* and the *observational noise*. Random transformations will also be split in two categories: *additive noise* and *randomly applied stochastic perturbations*. At this point we would like to stress that in the general theory of randomly perturbed dynamical systems one could consider perturbations other than the previous two. Our choice is motivated by the fact that our main result for the extreme values in presence of noise could be relatively easily shown with those assumptions, but it is also clear from the proof where possible generalizations could occur and we will quote a few of them. Let us notice that other authors basically used additive noise when they studied statistical properties of random dynamical systems [235, 236, 237], for instance. Random transformations with additive noise will be investigated with two methods: the first is the usual probabilistic approach which is at the base of all the previous mathematical formalisations presented in this book: we will briefly refer to it as the *probabilistic approach*. The second method is a spectral approach using perturbation theory of the Perron-Frobenius operator [238]: we will refer to it as the *spectral approach*. The probabilistic approach is in some sense more flexible since it is adapted to higher dimensions and, what is more important, it allows to get extreme value statistics even if the original unperturbed system does not obey it; we will show that on rotations and pure contractions on the interval. The spectral technique has the advantage to condensate in a very limited amount of assumptions what is really needed to prove the existence of limiting distribution and even to characterize the extremal index. The disadvantage is that the unperturbed map should have strong mixing properties, in particular the transfer operator must be quasi-compact on suitable functional spaces.

We will then describe two other kinds of noises, the *observational noise* and the *randomly applied stochastic perturbation* which will allow us respectively to detect the fractal properties of the invariant measure whenever the distribution of the maximum converges toward the Gumbel law and to get such an asymptotic distribution even for piecewise contracting maps. We will finally spend a few words about the extreme value statistics for *sequential dynamical systems*, where the stationarity of the random process is broken due to the non-autonomous concatenation of (close) maps.

## 7.2
## Random Transformations *via* the Probabilistic Approach: Additive Noise

The random transformations approach consists in perturbing a given deterministic dynamical system [1]. The latter will be given, in the following considerations, by a map $f$ acting either on a compact subset of $\mathbb{R}^D$ or on the $D$-dimensional torus. In both cases we will denote it with $\mathcal{M}$, and we denote with $\mathcal{B}$ the related Borel $\sigma$-algebra. If moreover the set is the one dimensional torus, we will denote it with $\mathcal{S}^1$; still for tori we will let $\mathrm{dist}(\cdot, \cdot)$ be the induced usual quotient metric on $\mathcal{M}$ and Leb a normalised volume form on the Borel sets of $\mathcal{M}$ that we call Lebesgue measure. As in the previous Chapters, we denote the ball of radius $\varepsilon > 0$ around $x \in \mathcal{M}$ as $B_\varepsilon(x) := \{y \in \mathcal{M} : \mathrm{dist}(x, y) < \varepsilon\}$.

Consider now a sequence of i.i.d. random variables $W_1, W_2, \dots$ taking values in the set $\Omega$ and with common distribution given by $\theta$; we associate to each $\omega \in \Omega$ a measurable map $f_\omega : \mathcal{M} \to \mathcal{M}$. Let $\Omega^{\mathbb{N}}$ denote the space of realisations of such process and $\theta^{\mathbb{N}}$ the product measure defined on its Borel subsets. Given a point $x \in \mathcal{M}$ and the realisation of the stochastic process $\underline{\omega} = (\omega_1, \omega_2, \dots) \in \Omega^{\mathbb{N}}$, we define the random orbit of $x$ as $x, f_{\underline{\omega}}(x), f_{\underline{\omega}}^2(x), \dots$ where, the evolution of $x$, up to time $n \in \mathbb{N}$, is obtained by the concatenation of the respective randomly perturbed maps in the following way:

$$f_{\underline{\omega}}^n(x) = f_{\omega_n} \circ f_{\omega_{n-1}} \circ \cdots \circ f_{\omega_1}(x), \tag{7.2.1}$$

with $f_{\underline{\omega}}^0$ being the identity map on $\mathcal{M}$.

As a very special case, but important for the next sections, we construct the *additive noise* by taking, for some $\varepsilon > 0$, $\Omega = B_\varepsilon(0)$ and $\theta := \theta_\varepsilon$ be a probability measure defined on the Borel subsets of $B_\varepsilon(0)$, such that

$$\theta_\varepsilon = g_\varepsilon \mathrm{Leb} \quad \int g_\varepsilon \, d\mathrm{Leb} = 1, \text{and} \quad 0 < \underline{g_\varepsilon} \leq g_\varepsilon \leq \overline{g_\varepsilon} < \infty. \tag{7.2.2}$$

where $\underline{g}_\epsilon$ and $\overline{g}_\epsilon$ are respectively a strictly positive lower bound and a finite upper bound for the density $g_\epsilon$. For each $\omega \in B_\varepsilon(0)$, we define the additive perturbation of $f$ that we denote by $f_\omega$ as the map $f_\omega : \mathcal{M} \to \mathcal{M}$, given by

$$f_\omega(x) = f(x) + \omega \bmod d \tag{7.2.3}$$

*Remark* 7.2.1. In the previous formula we implicitly considered $\mathcal{M}$ as an $D$-dimensional torus and this explains the mod-$d$ operation. We could consider any compact subset $\mathcal{M}$ of $\mathbb{R}^D$ by simply requiring that all the images of $\mathcal{M}$ are still in $\mathcal{M}$ and this of course imposes additional assumptions on the choice of the unperturbed map. Most of the following results will be stated for maps on the torus.

---

[1] We defer to [142] for a more detailed introduction to randomly perturbed dynamical systems and for exhaustive references

The ext definition introduces a notion that plays the role of invariance in the deterministic setting.

**Definition 7.2.2.** We say that the probability measure $\mu$ on the Borel subsets of $\mathcal{M}$ is stationary if

$$\iint \phi(f_\omega(x)) \, d\mu(x) \, d\theta(\omega) = \int \phi(x) \, d\mu(x),$$

for every $\phi : \mathcal{M} \to \mathbb{R}$ integrable w.r.t. $\mu$.

The previous equality could also be written as

$$\int \mathcal{U}\phi \, d\mu = \int \phi \, d\mu$$

where the operator $\mathcal{U} : L^\infty(\text{Leb}) \to L^\infty(\text{Leb})$, is defined as $(\mathcal{U}\phi)(x) = \int_\Omega \phi(f_\omega(x)) \, d\theta$ and it is called the *random evolution operator*. We define the *random Perron-Frobenius operator* as the linear operator $\mathcal{P} : L^1(\text{Leb}) \to L^1(\text{Leb})$ acting by duality as

$$\int \mathcal{P}\psi \cdot \phi \, d\text{Leb} = \int \mathcal{U}\phi \cdot \psi \, d\text{Leb} \tag{7.2.4}$$

where $\psi \in L^1$ and $\phi \in L^\infty$ and from now on these spaces will be be intended w.r.t. the Lebesgue measure, when the latter will be not explicitly mentioned.

It is immediate from this definition to get another useful representation of this operator, namely for $\psi \in L^1$:

$$(\mathcal{P}\psi)(x) = \int_\Omega (\mathcal{P}_\omega \psi)(x) \, d\theta(\omega),$$

where $\mathcal{P}$ is the Perron-Frobenius operator associated to $f_\omega$.

We recall that the stationary measure $\mu$ is absolutely continuous w.r.t. the Lebesgue measure and with density $h$ if and only if such a density is a fixed point of the random Perron-Frobenius operator: $\mathcal{P}h = h$

We can give a deterministic representation of this random setting using the following skew product transformation:

$$S : \begin{array}{ccc} \mathcal{M} \times \Omega & \longrightarrow & \mathcal{M} \times \Omega \\ (x, \underline{\omega}) & \longmapsto & (f_{\omega_1}, \sigma(\underline{\omega})), \end{array} \tag{7.2.5}$$

where $\sigma : \Omega \to \Omega$ is the one-sided shift $\sigma(\underline{\omega}) = \sigma(\omega_1, \omega_2, \dots) = (\omega_2, \omega_3, \dots)$. We remark that $\mu$ is stationary if and only if the product measure $\mu \times \theta^{\mathbb{N}}$ is an $S$-invariant measure.

Hence, the random evolution can fit the original deterministic model by taking the product space $\mathcal{X} = \mathcal{M} \times \Omega$, with the corresponding product Borel $\sigma$-algebra $\mathcal{B}$, where the product measure $\mathbb{P} = \mu \times \theta^{\mathbb{N}}$ is defined. The system is then given by the skew product map $T = S$. We recall here the notion of decay of correlation given in Definition 4.2.6 and point out that the correlations taken with respect to the (product) measure $\mathbb{P}$ are called *annealed*, to distinguish them from the *quenched* correlations, where the expectations are taken with respect to the stationary (or the ambient) measure and with a fixed realisation.

In the random setting, we will be mostly interested in Banach spaces of functions that do not depend on $\underline{\omega} \in \Omega$, hence, we assume that $\phi, \psi$ are actually functions defined on $\mathcal{M}$ and the correlation between these two observables can be written more simply as

$$\mathrm{Cor}_{\mathbb{P}}(\phi, \psi, n) := \frac{1}{\|\phi\|_{\mathcal{C}_1}\|\psi\|_{\mathcal{C}_2}} \left| \int \left( \int \psi \circ f_{\underline{\omega}}^n \, d\theta^{\mathbb{N}} \right) \phi \, d\mu - \int \phi \, d\mu \int \psi \, d\mu \right|$$

$$= \frac{1}{\|\phi\|_{\mathcal{C}_1}\|\psi\|_{\mathcal{C}_2}} \left| \int \mathcal{U}^n \psi \cdot \phi \, d\mu - \int \phi \, d\mu \int \psi \, d\mu \right| \quad (7.2.6)$$

where $(\mathcal{U}^n \psi)(x) = \int \cdots \int \psi(f_{\omega_n} \circ \cdots \circ f_{\omega_1} x) \, d\theta(\omega_n) \ldots d\theta(\omega_1) = \int \psi \circ f_{\underline{\omega}}^n(x) \, d\theta^{\mathbb{N}}$.

Note that when $\mu$ is absolutely continuous with respect to Leb and the respective Radon-Nikodym derivatives are bounded above and below by positive constants, then $L^1(\mathrm{Leb}) = L^1(\mu)$.

As in the deterministic case, the goal is to study the existence of an EVL for the partial maximum of observations made along the time evolution of the system. To be more precise consider the time series $X_0, X_1, X_2, \ldots$ arising from such a system simply by evaluating a given observable $\phi : \mathcal{M} \to \mathbb{R} \cup \{+\infty\}$ along the random orbits of the system:

$$X_n(x, \underline{\omega}) := \phi \circ f_{\underline{\omega}}^n(x), \quad \text{for each } n \in \mathbb{N}. \quad (7.2.7)$$

Clearly, $X_0, X_1, \ldots$ defined in this way is not an independent sequence, nevertheless the stationarity of $\mu$ guarantee that the stochastic process (7.2.7) ruled out by $\mathbb{P}$ is itself stationary.

We assume that the real-valued function $\phi : \mathcal{M} \to \mathbb{R} \cup \{\pm\infty\}$ achieves a global maximum at $\zeta \in \mathcal{M}$ (we allow $\phi(\zeta) = +\infty$). As before, we also assume that $\phi$ and $\mathbb{P}$ are sufficiently regular so that for $u$ sufficiently close to $u_F := \phi(\zeta)$, the event

$$U(u) = \{X_0 > u\} = \{x \in \mathcal{M} : \phi(x) > u\} \quad (7.2.8)$$

corresponds to a topological ball centred at $\zeta$; explicit examples of these prescriptions will be presented in the following. We assume that (R1) given in Chapter 4 holds.

Let $M_n$ be defined as in 2.2.1. As before, for every sequence $(u_n)_{n\in\mathbb{N}}$ satisfying (2.2.2) we define:

$$U_n := \{X_0 > u_n\}. \quad (7.2.9)$$

When $X_0, X_1, X_2, \ldots$ are not independent, the standard EVL law still applies under some conditions on the dependence structure; these conditions are the same as in the deterministic case since they are expressed solely in terms of the process $((X_k)_{k \geq 1}, \mathbb{P})$ and of the event $U_n$: we call them again $Ð_0(u_n)$ and $Ð'_0(u_n)$ and we refer to Chapter 4 for the definitions. We stress that Corollary 4.1.5 says that, if conditions $Ð_0(u_n)$ and $Ð'_0(u_n)$ hold for $X_0, X_1, \ldots$, then there exists an EVL for $M_n$ and the EVL reads $H(\tau) = 1 - e^{-\tau}$.

## 7.2.1
## Main Results

We anticipated above that by perturbing a very regular system like rotations or contractions, one could prove the existence of an EVL for some observable. One more advantage relies in the possibility of finding the scaling sequence $(u_n)_{n \in \mathbb{N}}$ in the case where it is affine $u_n = \frac{y}{a_n} + b_n$. In particular, this is achieved when the distance observable $\phi$ is of the form: $\phi(x) = -\log(\text{dist}(x, \zeta))$, being $\zeta$ a given point in $\mathcal{M}$, and for a particular choice of the probability measure $\theta$.

Before stating this result, we should recall that by using Theorem 1.7.13 in [1], one disposes of a sufficient condition to guarantee the existence of the limit (2.2.2) for $0 < \tau < \infty$. Such a condition requires that $\frac{1-F(x)}{1-F(x-)} \to 1$, as $x \to u_F$, where $F$ is again the distribution function of $X_0$, the term $F(x-)$ in the denominator denotes the left limit of $F$ at $x$ and $u_F$ was defined above as $\sup\{x; F(x) < 1\}$. For the observable just introduced $u_F = \infty$ and if the probability $\mathbb{P}$ satisfies (R1) we have that $F$ is continuous at $\zeta$ and therefore the above ratio goes to 1. We said above that in the context of random transformations the observables do not depend on the noise, therefore the probability $(\mathbb{P} > u_n)$ simply reduces to $(\mu_\varepsilon > u_n)$; therefore it will be sufficient to have a stationary measure which is not atomic at $\zeta$ to verify (R1). Of course this general result will not allow to compute explicitly the scaling coefficients $a_n$ and $b_n$.

*Remark* 7.2.3. From now on in this chapter we index all the random quantities, in particular the stationary measure $\mu$, the distribution $\theta$ on $\Omega$, the evolution operator $\mathcal{U}$, and the random Perron-Frobenius operator $\mathcal{P}$, with the lower index $\varepsilon$, where $\varepsilon$ represents the magnitude of the noise.

We start by recalling the following result from [239]:

**Proposition 7.2.4.** *Let us consider the dynamical systems $(\mathcal{M}, \mathcal{B}, f)$, where, as usual $\mathcal{M}$ is a compact set in $\mathbb{R}^D$ or it is the $D$-dimensional torus. We perturb it with additive noise admitting the stationary measure $\mu_\varepsilon$. We consider the associated process $X_n(x, \overline{\omega}) := -\log(|f_{\underline{\omega}}^n(x) - \zeta|)$ endowed with the probability $\mathbb{P} = \mu_\varepsilon \times \theta_\varepsilon^{\mathbb{N}}$, and suppose moreover that $\theta_\varepsilon$ is the Lebesgue measure measure on $B_\varepsilon(0)$. Then the linear sequence $u_n := u/a_n + b_n$ given in (2.2.2), verifies:*

$$a_n = d; \quad b_n = \frac{1}{D} \log\left(\frac{n \, K_d \, \mu_\varepsilon(B_\varepsilon(\zeta))}{(2\varepsilon)^D}\right),$$

*where $K_D$ is the volume on the unit hypercube in $\mathbb{R}^D$.*

We observe that the proof of the previous Proposition holds because $\theta_\varepsilon$ behaves exactly as the Lebesgue measure. This fact has another interesting consequence whenever we consider the additive noise. In fact by the very definition of the stationary measure, we have that for any Borel set $A$: $\mu_\varepsilon(A) = \int \theta_\varepsilon(\omega; f(x) + \omega \in A) d\mu_\varepsilon$. But $\theta_\varepsilon(\omega; f(x) + \omega \in A) \leq \frac{\text{Leb}(A)}{(2\varepsilon)^d}$, which implies that *the stationary measure is absolutely continuous with respect to the Lebesgue measure on $\mathcal{M}$.*

We now come to the proof of conditions $\not\!\!\Delta_0(u_n)$ and $\not\!\!\Delta_0'(u_n)$. We first collect a few facts and definitions:

- **(Y1)** We call *random transformations* the collections of random orbits given in (7.2.1 ), defined on the topological space $\mathcal{M}$ with the stationary measure $\mu_\varepsilon$ associated to the distribution $\theta_\varepsilon$. We remind that additive noise is a special case of it. We consider on the measurable space $(\mathcal{M}, \mathcal{B})$ two Banach spaces of observable, where $\mathcal{C}_2 = L^1(\text{Leb})$ and $\mathcal{C}_1$ are defined as the adapted spaces for which the decay of correlations applies.
- **(Y2)** We suppose that the decay actually holds and that the characteristic functions of measurable sets belong to $\mathcal{C}_1$.
- **(Y3)** We consider the process (7.2.7) $X_n := \phi \circ f_{\underline{\omega}}^n$, with the probability $\mathbb{P} = \mu_\varepsilon \times \theta_\varepsilon^{\mathbb{N}}$, where the observable $\phi$ obeys condition (R1) and it defines the level sets $U_n := \{\phi > u_n\}$ (see 7.2.9)), where $u_n$ is a sequence of real numbers. We finally suppose that condition (2.2.2) is satisfied.

From [142], we have

**Proposition 7.2.5.** *If a random transformation **(Y1)** enjoys property **(Y3)**, condition $(\not\!\!\Delta_0(u_n))$ is satisfied for the random process $X_n$ defined in **(Y3)**, provided there exists a constant $C'$ such that $||\mathbf{1}_{U_n}||_{\mathcal{C}_1} \leq C'$.*

*Proof. Sketch.* The proof is a simple computation of the correlation integral in Definition (4.2.6), by choosing the right observable. We should actually take [2]

$$\phi(x) = \mathbf{1}_{\{X_0 > u_n\}} = \mathbf{1}_{\{\varphi(x) > u_n\}},$$
$$\psi(x) = \int \mathbf{1}_{\{\varphi(x),\, \varphi \circ f_{\tilde{\omega}_1}(x),\, \dots,\, \varphi \circ f_{\underline{\tilde{\omega}}}^{\ell-1}(x) \leq u_n\}} \, \mathrm{d}\theta_\varepsilon^{\ell-1}(\tilde{\underline{\omega}}).$$

which, after substitution in the random evolution operator $\mathcal{U}_\varepsilon$ yields

$$\int \phi(x) \left(\mathcal{U}_\varepsilon^t \psi\right)(x) \, \mathrm{d}\mu_\varepsilon$$
$$= \int \mu_\varepsilon \Big(\varphi(x) > u_n, \varphi \circ f_{\underline{\omega}}^t(x) \leq u_n, \,\dots,\, \varphi \circ f_{\underline{\omega}}^{t+\ell-1}(x) \leq u_n\Big) \mathrm{d}\theta_\varepsilon^{\mathbb{N}}(\underline{\omega}).$$

---

[2] The integer $\ell$ in the definition of $\psi$ is motivated by the definition of the condition $\not\!\!\Delta_0(u_n)$, which actually requires to control the following correlation for all $\ell$, $t$ and $n$

$$|\mathbb{P}\left(X_0 > u_n \cap \max\{X_t, \dots, X_{t+\ell-1} \leq u_n\}\right) - \mathbb{P}(X_0 > u_n)\mathbb{P}(M_\ell \leq u_n)| \leq \gamma(n, t),$$

where $\gamma(n, t)$ is decreasing in $t$ for each $n$ and $n\gamma(n, t_n) \to 0$ when $n \to \infty$ for some sequence $t_n = o(n)$.

The result then follows by observing that

$$\int \phi(x)\mathrm{d}\mu_\varepsilon = \mu_\varepsilon(X_0(x) > u_n) = \int \mu_\varepsilon(X_0(x) > u_n)\,\mathrm{d}\theta_\varepsilon^{\mathbb{N}}(\underline{\omega})$$

$$\int \psi(x)\mathrm{d}\mu_\varepsilon = \int \mu_\varepsilon\Big(\varphi(x) \leq u_n, \varphi \circ f_{\omega_1}(x) \leq u_n, \ldots, \varphi \circ f_{\underline{\omega}}^{\ell-1}(x) \leq u_n\Big)\mathrm{d}\theta_\varepsilon^{\mathbb{N}}(\underline{\omega}).$$

$\square$

If we restrict ourselves to additive noise, it is possible to find systems that verify the previous proposition. Here we present three of such systems, for which we successively investigate property $\mathcal{A}_0'(u_n)$. We would like to warn the reader that the regularity conditions of the maps we introduce could be sometimes weakened; we defer to the quoted references for that; our aim here is to focus on the qualitative features of such maps.

- **S1**: *one-dimensional uniformly expanding maps*. We have already emphasized above that there will not be a substantial difference between maps defined on the unit circle and on the unit interval, if in the latter situation the image of the unit interval will be mapped strictly into itself when we add the noise to the unperturbed map. With this precision, we will refer in the following to maps $f$ defined on the unit interval $\mathcal{M}$. We will therefore suppose that

1) $f$ is locally injective on the open intervals $A_k$, $k = 1, \ldots, m$, that give a (finite) partition of the unit interval $\mathcal{M}$ up to zero measure sets.
2) $f$ is $C^2$ on each $A_k$ and has a $C^2$ extension to the boundaries. Moreover there exist $\Lambda > 1$, $C_1 < \infty$, such that $\inf_{x \in M}|Df(x)| \geq \Lambda$ and $\sup_{x \in M}\left|\frac{D^2 fx}{Df(x)}\right| \leq C_1$. This assumption ensures that the map $f$ admits many absolutely continuous invariant measure. We will therefore suppose that
3) $f$ admits a unique absolutely continuous invariant measure which is mixing. The Banach spaces invoked at point (**Y1**) above will be respectively $\mathcal{C}_2 = L^1(\text{Leb})$ and $\mathcal{C}_1$ : *the space of bounded variation (BV) functions defined in Section* 4.4.

- **S2**: *multidimensional uniformly expanding maps* These maps have been extensively investigated in [158, 129, 142, 153, 240, 241] and we defer to those papers for more details. We consider it a particular case corresponding to smooth boundaries. Let $\mathcal{M}$ be a compact subset of $\mathbb{R}^N$ which is the closure of its non-empty interior. We take a map $f : \mathcal{M} \to \mathcal{M}$ and let $\mathcal{A} = \{A_i\}_{i=1}^m$ be a finite family of disjoint open sets such that the Lebesgue measure of $\mathcal{M} \setminus \bigcup_i A_i$ is zero, and there exist open sets $\tilde{A}_i \supset \overline{A_i}$ and $C^{1+\alpha}$ maps $f_i : \tilde{A}_i \to \mathbb{R}^N$, for some real number $0 < \alpha \leq 1$ and some sufficiently small real number $\varepsilon_1 > 0$ such that

1) $f_i(\tilde{A}_i) \supset B_{\varepsilon_1}(f(A_i))$ for each $i$, where $B_\varepsilon(V)$ denotes a neighborhood of size $\varepsilon$ of the set $V$. The maps $f_i$ are the local extensions of $f$ to the $\tilde{A}_i$.
2) there exists a constant $C_1$ so that for each $i$ and $x, y \in f(A_i)$ with $\text{dist}(x, y) \leq \varepsilon_1$,

$$|\det Df_i^{-1}(x) - \det Df_i^{-1}(y)| \leq C_1|\det Df_i^{-1}(x)|\text{dist}(x, y)^\alpha;$$

3) there exists $s = s(f) < 1$ such that $\forall x, y \in f(\tilde{A}_i)$ with $\mathrm{dist}(x,y) \leq \varepsilon_1$, we have

$$\mathrm{dist}(f_i^{-1}x, f_i^{-1}y) \leq s\,\mathrm{dist}(x,y);$$

4) each $\partial A_i$ is a codimension-one embedded compact piecewise $C^1$ submanifold and

$$s^\alpha + \frac{4s}{1-s} Z(f) \frac{\gamma_{N-1}}{\gamma_N} < 1, \tag{7.2.10}$$

where $Z(f) = \sup_x \sum_i \#\{$smooth pieces intersecting $\partial A_i$ containing $x\}$ and $\gamma_N$ is the volume of the unit ball in $\mathbb{R}^N$.

5) The Banach space $\mathcal{C}_2$ will be again $L^1(\mathrm{Leb})$, while $\mathcal{C}_2$ will be the space of quasi-Hölder functions defined in Section 4.4. According to [242, Theorem 5.1], the assumptions (1)-(4) imply the existence of (many) acip's $\mu$ which will be in $V_\alpha$. Since these measures have densities which are convex combinations of the eigenfunctions of the dual of the transfer operator which is quasi-compact on $V_\alpha$, we could find a finite number of mixing measures for some power of $f$. In the following we will assume without restriction that:

6) There is only one mixing measure which is invariant with respect to $f$.

One of the main results in [142] is

**Theorem 7.2.6.** *Let us consider the class of transformations described in the items* **S1** *and* **S2** *and perturbed with additive noise with probability density equivalent to the Lebesgue measure on the compact set $\Omega$. We suppose moreover that properties* **Y1, Y2, Y3** *are satisfied, with $\mathrm{Cor}_\mathbb{P}(\phi, \psi, t) \leq Ct^{-2}$, for all $\phi \in \mathcal{C}_1$ and $\psi \in \mathcal{C}_2$. Then the process $X_n = \phi \circ f_{\underline{\omega}}^n$ satisfies $\underline{\underline{D}}_0(u_n)$ which together condition $\underline{\underline{D}}_0'(u_n)$ established in Proposition 7.2.5 imply that we have an EVL for $M_n$ such that $\overline{H}(\tau) = e^{-\tau}$.* [3]

*Proof. Sketch* We remind that we have to estimate the quantity

$$\lim_{n \to \infty} n \sum_{j=1}^{\lfloor n/k_n \rfloor} \mathbb{P}(X_0 > u_n, X_j > u_n), \tag{7.2.11}$$

where $(k_n)_{n \in \mathbb{N}}$ satisfies $k_n \to \infty$ and $k_n t_n = o(n)$. By remembering that the event $X_j > u_n$ equals, for our class of observable, the event $f_{\underline{\omega}}^{-j} U_n$, where $U_n = \{X_0 > u_n\}$, we decompose (7.2.11) as

$$n \sum_{j=1}^{\lfloor n/k_n \rfloor} \mathbb{P}(U_n \cap f_{\underline{\omega}}^{-j}(U_n)) \leq n \sum_{j=\alpha_n}^{\lfloor n/k_n \rfloor} \mathbb{P}\Big(\big\{(x,\underline{\omega}) : x \in U_n,\ f_{\underline{\omega}}^j(x) \in U_n\big\}\Big)$$

---

[3] Actually, the previous assumptions imply an exponential decay of correlations when the perturbation is not too large since the random Perron-Frobenius operator remains quasi-compact, see also Remark 7.3.2

$$+n \sum_{j=1}^{\lfloor n/k_n \rfloor} \mathbb{P}\Big( \big\{ (x, \underline{\omega}) : x \in U_n, R^{\underline{\omega}}(U_n) \leq \alpha_n \big\} \Big) := I + II. \tag{7.2.12}$$

where, for a given set $A$, we have defined $R^{\underline{\omega}}(A) := \min_{x \in A} \min_{j \in \mathbb{N}} \{ f_{\underline{\omega}}^j(x) \in A \}$, and $(\alpha_n)$ is some sequence such that $\alpha_n \to \infty$ and $\alpha_n = o(\log k_n)$. The first term takes into account the long returns and it can be handled by decay of correlations. We should stress at this point that decay against $L^1$ functions plays here an essential role. We in fact have for suitable constants $C$ and $C^*$ :

$$\mathbb{P}\big( \{ (x, \underline{\omega}) : x \in U_n, f_{\underline{\omega}}^j(x) \in U_n \} \big) \leq (\mu_\varepsilon(U_n))^2 + C \left\| \mathbf{1}_{U_n} \right\|_{\mathcal{C}_1} \left\| \mathbf{1}_{U_n} \right\|_{L^1(\mu_\varepsilon)} j^{-2}$$
$$\leq (\mu_\varepsilon(U_n))^2 + C^* \mu_\varepsilon(U_n) j^{-2}, \tag{7.2.13}$$

where $\mathcal{C}_1$ is the Banach space adapted to $L^1$. We should also stress here another useful fact: the Banach norms are actually computed on characteristic functions and this will turn out to be crucial in a forthcoming situation. By returning to (7.2.12) and by summing over $j$, we see that the first term on the right hand side will go to zero since by assumption $\mu_\epsilon(U_n) \sim \tau/n$ and for the same reason the second goes to 0 after having summed the rest of the convergent series in $j^{-2}$. The second term $II$ deals with short returns and it is usually the difficult one to compute because it uses the topological features of the map. In particular since our maps are defined on (quotient of) Banach vector spaces, and by using the translation invariance of the Lebesgue measure (we denote here with $g_\varepsilon$ the density of $\theta_\varepsilon$), we have for another suitable constant which takes into account also the finite expectation of $g_\varepsilon$ :

$$\theta_\varepsilon^{\mathbb{N}}\big( \{ \underline{\omega} : R^{\underline{\omega}}(U_n) \leq \alpha_n \} \big) \leq \sum_{j=1}^{\alpha_n} \int \theta_\varepsilon \Big( \big\{ \omega_j : f\big(f_{\underline{\omega}}^{j-1}(\zeta)\big) + \omega_j \in B_{2\eta^j |U_n|}(\zeta) \big\} \Big) \mathrm{d}\theta_\varepsilon^{\mathbb{N}}$$
$$= \sum_{j=1}^{\alpha_n} \int \theta_\varepsilon \Big( \big\{ \omega_j : \omega_j \in B_{2\eta^j |U_n|}(\zeta) - f\big(f_{\underline{\omega}}^{j-1}(\zeta)\big) \big\} \Big) \mathrm{d}\theta_\varepsilon^{\mathbb{N}}$$
$$= \sum_{j=1}^{\alpha_n} \iint_{B_{2\eta^j |U_n|}(\zeta) - f(f_{\underline{\omega}}^{j-1}(\zeta))} g_\varepsilon(x) \mathrm{d}\mathrm{Leb} \, \mathrm{d}\theta_\varepsilon^{\mathbb{N}}$$
$$\leq \sum_{j=1}^{\alpha_n} \overline{g_\varepsilon} \mathrm{Leb}\big( B_{2\eta^j |U_n|}(\zeta) \big) \leq C \mathrm{Leb}(U_n) \frac{\eta}{\eta - 1} \eta^{\alpha_n}.$$

The quantity $\eta > 1$ is a Lipschitz constant for the map $f$ which could be easily defined locally if the map is piecewise continuous on a finite number of domains. By summing over $j$ and remembering the scaling of $k_n$, $\alpha_n$ and $\mu_\varepsilon(U_n)$, it is easy to check that this second term converges to zero too. $\qquad \square$

We said above that in order to establish the EVL, it will be enough to get decay of correlations for characteristic functions. This turns out to be relevant in the rather unexpected situation of random perturbations of rotations (including rational ones). Let us therefore introduce our third class of maps:

- **S3**: *rotations of the unit circle.* We consider on the unit circle $\mathcal{M}$ the map $f(x) = x + \alpha \bmod 1$, $\alpha \in \mathbb{R}$, and we perturb it with random noise in this way: $f_{\varepsilon\omega} = x + \alpha + \varepsilon\omega \bmod 1$, where $\omega$ is a random variable uniformly distributed over the interval $[-1, 1]$.

Notice that we slightly modified the previous notations, especially to emphasize the role of the strength of the perturbation $\varepsilon$. Therefore $\theta_\varepsilon$, the distribution of the noise, will now be the normalized Lebesgue measure over $[-1, 1]$ and the stationary measure $\mu_\varepsilon$ is now simply the Lebesgue measure $(dx)$ over the unit interval. Finally, we continue to call $\mathcal{U}_\varepsilon$ the random evolution operator. We observe that Proposition 7.2.4 applies to our case and therefore the scaling (2.2.2) is guaranteed: we consider again the process generated by the observations introduced in (7.2.7) with the probability $\mathbb{P}$ given by assumption (**Y3**) above. The first interesting result is the decay of correlations for the evolution of measurable sets, which can be proved using the Fourier transform as in [79]:

**Proposition 7.2.7.** *Let us take the observable $\phi = \mathbf{1}_A$ and $\psi = \mathbf{1}_B$, where $B = \cup_{l=1}^{\ell} B_\ell$, for some $\ell \in \mathbb{N}$ and $A, B_1, \ldots, B_\ell \subset \mathcal{M}$ are connected intervals, then*

$$C_{j,\varepsilon} := \left| \int_0^1 \mathcal{U}_\varepsilon^j(\psi)\phi\, dx - \int_0^1 \psi\, dx \int_0^1 \phi\, dx \right| \leq 4e^{-j\varepsilon^2 \log(2\pi)}, \qquad (7.2.14)$$

*as long as $\varepsilon^2 < 1 - \log 2/\log(2\pi)$.*

This will enable us to prove immediately condition $\not{D}_0(u_n)$ and assumption $\not{D}_0'(u_n)$ is proved similarly to Theorem 7.2.6. We could summarise by establishing one of the main results in [79]:

**Theorem 7.2.8.** *Let us consider rotations of the unit circle perturbed with additive noise with uniform distribution. Then the process $X_0, X_1, \ldots$, given by (7.2.7) verifies the scaling of Proposition 7.2.4 and the distribution of the maximum $M_n$ will converge to an EVL with $\overline{H}(\tau) = e^{-\tau}$.*

This result address the question of whether it is possible to distinguish numerically the real nature of the underlying (unperturbed) system when we look at the extremal statistics of the randomly perturbed data. To answer this question which is of practical interest, we used in [79] the fact that the numerical round off is comparable to a random noise on the last precision digit [243]. This observation allowed us to claim that, for systems featuring periodic or quasi periodic motions (as the rotations on the circle), such a noise level is not sufficient for producing detectable changes regarding the observed extremal behavior. In fact the simulations produced in [79] clearly showed that EVLs are obtained when considering small but finite noise amplitudes only when very long trajectories are considered. The quality of the fitting improves when larger bins are considered: this is in agreement with the idea that we should get EVL for infinitely small noises in the limit of infinitely long samples. In our case, EVLs are obtained only for $\varepsilon > 10^{-4}$, which is still considerably larger than the noise introduced by round-off resulting from double precision, as the round-off procedure is equivalent to the addition to the exact map of a random noise of order $10^{-7}$ [243, 244].

## 7.3
## Random Transformations *via* the Spectral Approach

We introduced above the quantity $R^{\underline{\omega}}(A)$ as the first random return of the set $A$ to itself; it is equivalently defined as as the infimum over $x \in A$ of the first random hitting time defined by

$$r_A^{\frac{\omega}{}}(x) := \min\{j \in \mathbb{N} : f_{\underline{\omega}}^j(x) \in A\} \qquad (7.3.1)$$

As we will see in the next section this quantity is intimately related to the EVL statistics; for our actual pourposes it will be enough to observe that if we consider again the distribution of the maximum $\{M_m \leq u_m\}$ (we keep here the notations of the previous section), where $U_m = \{\phi > u_m\}$ is the event given by a topological ball shrinking to the point $\zeta$, then

$$(\mu_\varepsilon \times \theta_\varepsilon^{\mathbb{N}})((x, \underline{\omega}) : r_m^{\frac{\omega}{}}(x) > m) = (\mu_\varepsilon \times \theta_\varepsilon^{\mathbb{N}})(M_m \leq u_m) \qquad (7.3.2)$$

Let us write the measure on the left-hand side of (7.3.2) in terms of integrals: it is given by

$$\iint\limits_{\{r_m^{\frac{\omega}{}} > m\}} \mathrm{d}(\mu_\varepsilon \times \theta_\varepsilon^{\mathbb{N}})$$
$$= \iint h_\varepsilon \mathbf{1}_{U_m^c}(x) \mathbf{1}_{U_m^c}(f_{\omega_1} x) \cdots \mathbf{1}_{U_m^c}(f_{\omega_{m-1}} \circ \cdots \circ f_{\omega_1} x) \, \mathrm{d}\mathrm{Leb} \, \mathrm{d}\theta_\varepsilon^{\mathbb{N}} \qquad (7.3.3)$$

which is in turn equal to

$$\int_M \widetilde{\mathcal{P}}_{\varepsilon,m}^m h_\varepsilon(x) \, \mathrm{d}\mathrm{Leb} \qquad (7.3.4)$$

where we have now defined

$$\widetilde{\mathcal{P}}_{\varepsilon,m}\psi(x) := \mathcal{P}_\varepsilon(\mathbf{1}_{U_m^c}\psi)(x), \qquad (7.3.5)$$

Let us note that the operator $\widetilde{\mathcal{P}}_{\varepsilon,m}$ depends on $m$ via the set $U_m$, and not on $\varepsilon$ which is kept fixed and that $\widetilde{\mathcal{P}}_{\varepsilon,m}$ "reduces" to $\mathcal{P}_\varepsilon$ as $m \to \infty$. It is therefore tempting to consider $\widetilde{\mathcal{P}}_{\varepsilon,m}$ as a small perturbation of $\mathcal{P}_\varepsilon$, the random Perron-Frobenius operator defined in (7.2.4), when $m$ is large and to see if the spectral properties of $\mathcal{P}_\varepsilon$ could give information on the behavior of $\widetilde{\mathcal{P}}_{\varepsilon,m}$. This approach has been successfully proposed by Keller [141] for deterministic systems where the Perron-Frobenius operator is quasi-compact on suitable adapted spaces. A few properties should be verified by those systems. These properties are referred to as *Rare event Perron-Frobenius operators* (REPFO) by Keller and we summarise them below. Although those assumptions are stated in a general way, it turns out that they can be successfully checked when one of the adapted spaces is $L^1(\mathrm{Leb})$. This is equivalent to ask exponential decay of correlations on the mixing components against $L^1(\mathrm{Leb})$ observables. We remind that such a decay reveals to be particularly useful in the usual probabilistic approach to EVL in order to control the short returns in condition $\not{\Pi}'_0(u_n)$. As a byproduct of the *REPFO* technique, one gets an explicit computation of the EI. Next, we illustrate Keller's method in the random setting, and in particular for random maps defined on the unit interval or on the unit circle. We choose as the couple of adapted spaces $L^1$(Leb) and the space BV of bounded variation functions introduced in Section 4.4. We assume:

**(A1)** There are constants $A > 0, B > 0, D > 0$ such that:

$$\forall m \geq 1, \forall \psi \in \mathcal{V}, \forall n \in \mathbb{N} : ||\widetilde{\mathcal{P}}^n_{\varepsilon,m}\psi||_1 \leq D||\psi||_1 \qquad (7.3.6)$$

$$\exists \alpha \in (0,1), \forall m \geq 1, \forall \psi \in \mathcal{V}, \forall n \in \mathbb{N} : ||\widetilde{\mathcal{P}}^n_{\varepsilon,m}\psi||_{BV} \leq A\alpha^n ||\psi||_{BV} + B||\psi||_1 \qquad (7.3.7)$$

**(A2)** The unperturbed operator $\mathcal{P}_\varepsilon$ verifies the mixing condition

$$\mathcal{P}_\varepsilon = h_\varepsilon \otimes \text{Leb} + Q_\varepsilon,$$

where $h_\varepsilon$ is the density of the stationary measure and $Q_\varepsilon$ has spectral radius less than 1.

**(A3)** $\exists C > 0$ such that

$$\eta_m := \sup_{\|\psi\|_{BV} \leq 1} \left| \int (\mathcal{P}_\varepsilon - \widetilde{\mathcal{P}}_{\varepsilon,m})\psi \, \text{dLeb} \right| \to 0, \text{ as } m \to \infty \qquad (7.3.8)$$

**(A4)** and

$$\eta_m \left\| (\mathcal{P}_\varepsilon - \widetilde{\mathcal{P}}_{\varepsilon,m})\varphi_0 \right\|_{BV} \leq C \, \Delta_{\varepsilon,m} \qquad (7.3.9)$$

where $\Delta_{\varepsilon,m} = \mu_\varepsilon(U_m) = \text{Leb}((\mathcal{P}_\varepsilon - \widetilde{\mathcal{P}}_{\varepsilon,m})h_\varepsilon)$.

In [142] the following results was obtained

**Proposition 7.3.1.** *The maps introduced in **(S1)** perturbed with additive noise verify Assumptions **(A1)** to **(A2)** provided we add to the three items quoted in the former **(S1)** and the additional requirement:*
*(4) The density $h_\varepsilon$ of the stationary measure is bounded from below* Leb-*a.e.* .

*Remark* 7.3.2. It should be pointed out that in order to prove the mixing condition **(A2)** we use the fact that the original unperturbed map, with Perron-Frobenius operator $\mathcal{P}$, is mixing and the fact that the $L^1$-norm of the difference $\mathcal{P} - \mathcal{P}_\varepsilon$ is of order $\varepsilon$ for the additive noise.

We now return to Keller's theory; the above assumptions ensure that the operator $\widetilde{\mathcal{P}}_{\varepsilon,m}$ verifies the spectral equations

$$\widetilde{\mathcal{P}}_{\varepsilon,m}\varphi_{\varepsilon,m} = \lambda_{\varepsilon,m}\varphi_{\varepsilon,m}; \; \nu_{\varepsilon,m}\widetilde{\mathcal{P}}_{\varepsilon,m} = \lambda_{\varepsilon,m}\nu_{\varepsilon,m}, \; \lambda^{-1}_{\varepsilon,m}\widetilde{\mathcal{P}}_{\varepsilon,m} = \varphi_{\varepsilon,m}\otimes\nu_{\varepsilon,m} + Q_{\varepsilon,m}$$

where $\sup_{m\in\mathbb{N}}||Q^n_{\varepsilon,m}||_{BV}$ decays exponentially as $n \to \infty$ and the eigenvalues $\lambda_{\varepsilon,m}$ obey $1 - \lambda_{\varepsilon,m} = \Delta_{\varepsilon,m}\vartheta_\varepsilon(1 + o(1))$, where in turns $\vartheta_\varepsilon$ verifies $\lim_{\varepsilon\to 0}\frac{1-\lambda_{\varepsilon,m}}{\Delta_{\varepsilon,m}} = \vartheta_\varepsilon := 1 - \sum_{k=0}^{\infty} q_{k,\varepsilon}$, with

$$q_{k,\varepsilon} := \lim_{m\to\infty} q_{k,\varepsilon,m} := \lim_{m\to\infty}\frac{\text{Leb}((\mathcal{P}_\varepsilon - \widetilde{\mathcal{P}}_{\varepsilon,m})\widetilde{\mathcal{P}}^k_{\varepsilon,m}(\mathcal{P}_\varepsilon - \widetilde{\mathcal{P}}_{\varepsilon,m})(h_\varepsilon))}{\Delta_{\varepsilon,m}} \quad (7.3.10)$$

Armed with these asymptotic expansions of the eigenvalues $\lambda_{\varepsilon,m}$ for $m$ large, we have

$$(\mu_\varepsilon \times \theta_\varepsilon^{\mathbb{N}})(M_m \leq u_m) = \int_M \widetilde{\mathcal{P}}_{\varepsilon,m}^m h_\varepsilon(x) \, \mathrm{dLeb}$$

$$= \lambda_{\varepsilon,m}^m \int h_\varepsilon \, \mathrm{d}\nu_{\varepsilon,m} + \lambda_{\varepsilon,m}^m \int Q_{\varepsilon,m}^m h_\varepsilon \, \mathrm{dLeb}$$

$$= e^{-(\vartheta_\varepsilon m \, \mu_\varepsilon(U_m) + m \, o(\mu_\varepsilon(U_m)))} \int h_\varepsilon \, \mathrm{d}\nu_{\varepsilon,m} + \mathcal{O}(\lambda_{\varepsilon,m}^m \left\| Q_{\varepsilon,m}^m \right\|_{BV})$$

Remember that we are under the assumption that $m \, (\mu_\varepsilon \times \theta_\varepsilon^{\mathbb{N}})(\phi > u_m) = m \, \mu_\varepsilon(\phi > u_m) = m \, \mu_\varepsilon(U_m) \to \tau$, when $m \to \infty$; moreover it follows from the theory of [148] that $\int h_\varepsilon \, \mathrm{d}\nu_{\varepsilon,m} \to \int h_\varepsilon \, \mathrm{dLeb} = 1$, as $m$ goes to infinity. In conclusion we get

$$(\mu_\varepsilon \times \theta_\varepsilon^{\mathbb{N}})(M_m \leq u_m) = e^{-\tau \vartheta_\varepsilon}(1 + o(1))$$

in the limit, as $m \to \infty$. We remind that in the deterministic case the exponent multiplying $\tau$ was defined as the EI at the point $\zeta$, see Definition 3.2.4. It is interesting to observe that in the random setting the extremal index will always be the same ($= 1$), for all the points $\zeta$. This is the content of the next proposition drawn from [142], whose proof is reminiscent of that of Theorem 7.2.6.

**Proposition 7.3.3.** *Let us suppose that $f$ verifies the assumptions of Proposition 7.3.1, but with the density $h_\varepsilon$ not necessarily bounded away from zero. Then for each $k$,*

$$\lim_{m\to\infty} q_{k,\varepsilon,m} = 0,$$

*i.e., the limit in the definition of $q_{k,\varepsilon}$ in (7.3.10) exists and equals zero. Also the EI verifies $\vartheta_\varepsilon = 1 - \sum_{k=0}^{\infty} q_{k,\varepsilon} = 1$ and this is independent of the point $\zeta$ serving as centre of the ball $U_m$.*

Another interesting application of this technique is that it provides an explicit formula for HTS with respect to the (annealed) probability.

**Proposition 7.3.4** ([142])**.** *Let us suppose that $f$ verifies the assumptions of Proposition 7.3.1; then there exists a constant $C > 0$ such that for all $m$ big enough there exists $\xi_m > 0$ s.t. for all $t > 0$*

$$\left| (\mu_\varepsilon \times \theta_\varepsilon^{\mathbb{N}}) \left\{ r_{U_m}^\omega > \frac{t}{\xi_m \, \mu_\varepsilon(U_m)} \right\} - e^{-t} \right| \leq C \delta_m (t \vee 1) e^{-t}$$

*where $\delta_m = O(\eta_m \log \eta_m)$, and*

$$\eta_m = \sup \left\{ \left| \int_{U_m} \psi \, \mathrm{dLeb} \right|; \|\psi\|_{BV} \leq 1 \right\}.$$

*Moreover, $\xi_m$ goes to $\vartheta_\varepsilon$ as $m \to \infty$, and $\eta_m$ goes to zero being bounded by $Leb(U_m)$.*

This result is the random counterpart of the relation between the EVL and the HTS established in Section 5.3. We get here basically the same result based on the identification (7.3.2). We would like to point out that the convergence towards the exponential law and for *any* $\zeta$ could be obtained as well using the pure probabilistic approach of Section 7.2, see in [142, Corollary E]. The spectral approach strengthens the latter, since it provides the error in the convergence to the exponential law. Instead the probabilistic approach allows us to get a further statistical property, by considering the distribution of the number of exceedances (or hits to $U(u_n)$) during a suitable re-scaled time period. Consider the REPP defined in (3.3.1). In [142] the following result was proved:

**Proposition 7.3.5.** *Let us suppose that $f$ verifies the assumptions of Proposition 7.3.1; then for the stochastic process defined in (7.2.7) the REPP $N_n$ converges in distribution to $N$ for $n \to \infty$, where $N$ denotes a Poisson Process with intensity 1.*

*Remark* 7.3.6. The paper [245] gets similar results for the annealed distribution of the first hitting and return times defined with respect to the probability $\mathbb{P}$ introduced above. The author used super-polynomial decay of correlations against $L^\infty$ observables instead of $L^1$ functions. The perturbations used are different from the additive noise.

## 7.4
## Random Transformations *via* the Probabilistic Approach: Randomly Applied Stochastic Perturbations

We have seen in a previous section how a very regular systems like rotations could enjoy EVL whenever it is perturbed with additive noise. We now introduce another kind of noise which will allows us to get EVL even for perturbed piecewise contacting maps (PCM). The noise we refer to was introduced by Lasota and Mackey (see [246], for instance) and correspond to *randomly applied stochastic perturbations* (RASP). They consist in operating a random reset of the initial condition of the original dynamical system $(\mathcal{M}, f)$[4], at each failure of a Bernoulli random variable: if $(x_n)_{n\in\mathbb{N}}$ denotes the successive states of such a random dynamical systems, then at each time $n \in \mathbb{N}$ we have $x_{n+1} = f(x_n)$ with probability $(1 - \epsilon)$ and $x_{n+1} = \pi_n$ with probability $\epsilon$, where $\pi_n$ is the realization of a random variable with value in $X$. This kind of perturbation corresponds to the family $(f_\omega)_{\omega \in \Omega_\varepsilon}$ of random transformations defined by

$$f_\omega(x) = \eta f(x) + (1 - \eta)\pi \qquad \forall x \in X, \tag{7.4.1}$$

where $\omega = (\eta, \pi)$ is a random vector with value in $\Omega_\varepsilon = \{0, 1\} \times X$. The two components $\eta$ and $\pi$ of $\omega$ are independent and $\eta$ is a Bernoulli variable with the probability of being 0 equal to $\varepsilon$, while $\pi$ is a random variable that we will suppose Lebesgue-uniformly distributed on $X$. The joint distribution $\theta_\varepsilon$ of these two components is the product of the Bernoulli measure with weights $(1 - \varepsilon, \varepsilon)$ and the uniform measure on $X$. We will consider again the process (7.2.7) $X_n := \varphi \circ f_{\underline{\omega}}^n$, with the probability $\mathbb{P} = \mu_\varepsilon \times \theta_\varepsilon^{\mathbb{N}}$, where the observable $\phi$ defines the level sets $U_n := \{\phi > u_n\}$ (see 7.2.9)), where $u_n$ is a sequence of real numbers. We now apply this kind of noise to contracting maps on the unit interval; we present here a few already established results [239], further generalizations are under investigations and they will be simply quoted in the following [247].

---

4) We are using the notation introduced in the previous section, where the maps $f$ act on a compact subset $\mathcal{M}$ of $\mathbb{R}^D$ or on the $D$-dimensional torus.

Let us therefore take the map $f$ defined on the unit interval $I = [0, 1]$ by

$$S(x) = \alpha x, \quad \alpha \in (0, 1).$$

We perturb it according to eq. (7.4.1) giving raise to the family of random maps defined for each $n \geq 1$ by

$$f_{\omega_n}(x) = \eta_n S(x) + (1 - \eta_n)\pi_n, \quad \forall x \in I,$$

where $\omega_n = (\eta_n, \pi_n)$.

The first interesting result is that there is only one stationary measure $\mu_\varepsilon$ which is absolutely continuous with respect to the Lebesgue measure on the unit interval and the density $h_\varepsilon$ reads [246]:

$$h_\varepsilon(x) = \varepsilon \sum_{k=0}^{p-1} \frac{(1 - \varepsilon)^k}{\alpha^k} \quad \forall x \in (\alpha^p, \alpha^{p-1}], \ p \geq 1.$$

Notice that the density is bounded for $(1 - \varepsilon) < \alpha$. We can also compute explicitly the scaling coefficients of the affine normalization $u_n$:

**Proposition 7.4.1** ([247])**.** *Let $\tau > 0$, $y = -\ln(\tau)$ and $u_n = \frac{y}{a_n} + b_n$ with*

$$a_n = 1 \quad and \quad b_n = \log\left(2n\varepsilon \sum_{k=0}^{p-1} \frac{(1 - \varepsilon)^k}{\alpha^k}\right) \qquad \forall n \in \mathbb{N},$$

*If:*
*(i) $z \neq 0$ on the interval but not in the countably many discontinuity points of $h_\varepsilon$, namely $z \notin \cup_{j \in \mathbb{N}}\{\alpha^j\}$;*
*(ii) $p \geq 1$ such that $z \in (\alpha^p, \alpha^{p-1})$ and $n$ is large enough such that the ball $B(z, e^{-u_n}) \subset (\alpha^p, \alpha^{p-1})$, then:*

$$n\mathbb{P}(X_0 > u_n) = \tau$$

The next step will be to check the conditions $\text{Ð}_0(u_n)$ and $\text{Ð}_0'(u_n)$. To verify condition $\text{Ð}_0(u_n)$ we need to show that for specific observables $\phi \in L^\infty$ and $\psi \in L^1$ the correlations

$$Cor(\phi, \psi, n) := \left| \int \mathcal{U}_\varepsilon^n(\phi(x))\psi(x)d\mu_\varepsilon - \int \phi(x)d\mu_\varepsilon \int \psi(x)d\mu_\varepsilon \right|$$

$$\left| \int \int \phi(f_{\underline{\omega}}^n(x))\psi(x)d\mu_\varepsilon d\theta_\varepsilon^{\mathbb{N}} - \int \phi(x)d\mu_\varepsilon \int \psi(x)d\mu_\varepsilon \right| \tag{7.4.2}$$

decay sufficiently fast with $n$; $\mathcal{U}$ denotes again the random evolution operator.

**Proposition 7.4.2** ([247]). *If $\phi \in L^\infty$ and $\psi \in L^1 \cap L^\infty$, then*

$$Cor(\phi, \psi, n) \leq 2(1-\varepsilon)^n ||\phi||_{L^\infty} ||\psi h_\varepsilon||_{L^1}.$$

We notice that condition $\amalg_0(u_n)$ requires that $\psi$ and $\phi$ are characteristic functions of measurable sets. Condition $\amalg_0'(u_n)$ is needed to control short returns in the ball around $z$. We can prove that it holds for $z \neq 0$ and the proof is basically based on the fact that the image of the ball does not intersect the ball itself for large $n$. Instead, whenever $z = 0$ an extremal index appears in the limiting law for the distribution of the maxima.

**Proposition 7.4.3** ([239, 247]). *For the map $S(x) = \alpha x$, $\alpha \in (0, 1)$ perturbed with the noise (7.4.1), and by considering the observable $\phi(x) = -\log(|x - z|)$, conditions $\amalg_0(u_n)$ and $\amalg_0'(u_n)$ hold.*
*If $z = 0$ we have the existence of an extremal index less than 1, which is given by $\varepsilon$.*

The previous result could be generalised in the following direction. We first define the system

**Definition 7.4.4. RASP maps** Let us take $X$ a compact subset of $\mathbb{R}^D$ which is the closure of its interior $X = \overline{\text{int}(X)}$, equipped with some metric $d$ and with the normalized Lebesgue measure Leb; suppose $\{X_i\}_{i=1}^N$ is a collection of $N$ disjoint open subsets of $X$ such that $X = \bigcup_{i=1}^N \overline{X_i}$ and $\text{Leb}(\Delta) = 0$, where $\Delta := X \setminus \bigcup_{i=1}^N X_i$ is the *singular set*. We will consider non singular maps $f : X \to X$ which are injective and such that $f|_{X \setminus \Delta}$ is a $C^1$-diffeomorphism. We say that $f$ is continuous in a point $w$ of the boundary of $X$ if there is an open ball of radius $\kappa$, $B_\kappa(w)$ such that $f$ is continuous on $B_\kappa(w) \cap \text{int}(X)$ and it can be extended continuously on $B(w, \kappa)$.

In the examples we have treated, the set $\Delta$ will be composed by points where $f$ is discontinuous and by the points of the boundary of $X$. We would like to stress that a point of that boundary is not necessarily a discontinuity point of $f$.

The density of the stationary measure has the expression

$$h_\varepsilon(x) = \varepsilon \sum_{k=0}^\infty (1-\varepsilon)^k J_k(x) \mathbf{1}_{\Lambda_k}(x) \quad \forall x \in X. \tag{7.4.3}$$

where $\{\Lambda_k\}_{k\in\mathbb{N}}$ be the sequence of sets defined by $\Lambda_0 := X$ and $\Lambda_{k+1} := f(\Lambda_k \setminus \Delta)$ for all $k \geq 1$. We denote with $J_k(x) := \prod_{l=1}^k |\det(f'(f^{-l}(x))|^{-1}$ for all $k \geq 1$ and $J_0(x) := \mathbf{1}(x)$.

We then decompose our study of the existence of EVL in two parts. The first part concerns the point of $\Lambda_1 \setminus \tilde{\Lambda}$, where

$$\tilde{\Lambda} := \bigcap_{k \in \mathbb{N}} \overline{\Lambda}_k. \tag{7.4.4}$$

The set $\tilde{\Lambda}$ is called the *global attractor* of $f$ and contains all the point of $\Lambda$. In $\Lambda_1 \setminus \tilde{\Lambda}$ we can verify the conditions $Ð_0(u_n)$ and $Ð'_0(u_n)$ to show the existence of EVL still for the observable $\phi(\cdot) = -\log d(\cdot, \zeta), \zeta \in X$, but in $\tilde{\Lambda}$ condition $Ð'_0(u_n)$ is not verified. For this reason in the second part, we introduce other conditions which allow us to show the existence of an EVL but with an extremal index different from 1; this is achieved by assuming an additional property, namely, by choosing $f$ as a piecewise contracting map. It has been shown in [248] that for any piecewise contracting map defined on a compact space, if the global attractor $\tilde{\Lambda}$ does not intersect the set of the discontinuities, then it is composed of a finite number of periodic orbits. Moreover, this condition, for injective maps as ours, is generic in the $C^0$ topology [249].

**Theorem 7.4.5** ([247]). *Suppose $f : X \to X$ satisfies Definition 7.4.4 and consider the random perturbations defined in (7.4.1) with $\mathbb{P} = \mu_\varepsilon \times \theta_\varepsilon^{\mathbb{N}}$, where $\mu_\varepsilon$ is the stationary probability measure. Then, for any $z \in \Lambda_1 \setminus \tilde{\Lambda}$, the sequence $\{M_n\}_{n \in \mathbb{N}}$ of the maxima of the process defined for every $n \in \mathbb{N}$ by $Y_n(x, \underline{\omega}) := -\log d(f_{\underline{\omega}}^n(x), z)$ admits the Gumbel's law as extreme values distribution.*

*Suppose now that $f$ is piecewise contracting, that is, there exists $\alpha \in (0, 1)$ such that for every $i \in \{1, \ldots, N\}$ we have*

$$d(f(x), f(y)) \le \alpha d(x, y) \quad \forall x, y \in X_i.$$

*Suppose that the set $\tilde{\Lambda}$ does not contain any discontinuity point of $f$, and let $z \in \tilde{\Lambda}$. Let $Y_0, Y_1, \ldots$ be the stochastic process given by $Y_n(x, \underline{\omega}) := -\log d(f_{\underline{\omega}}^n(x), \mathcal{O}(z)$, where $\mathcal{O}(z)$ is the orbit of $z$. Let $(u_n)_{n \in \mathbb{N}}$ be a sequence such that $\lim_{n \to \infty} n\mathbb{P}(X_0 > u_n) = \tau > 0$. Then we have*

$$\lim_{n \to \infty} \mathbb{P}(M_n \le u_n) = e^{-\varepsilon \tau}.$$

We notice that the existence of the limit $\lim_{n \to \infty} n\mathbb{P}(X_0 > u_n) = \tau > 0$ could be proved again with general arguments, since $\mathbb{P}$ satisfies (R1), but the scaling coefficients $a_n$ and $b_n$ could also be computed explicitly even on the global attractor provided the density of the stationary measure is essentially bounded.

Typical examples of piecewise contracting maps for which the attractor is generically a finite set of periodic orbits are piecewise affine. That is, if for each $i \in \{1, \ldots, N\}$ the restriction $f_i := f|_{X_i}$ of the map $f : X \to X$ of Definition 7.4.4 to a piece $X_i$ has the following form:

$$f_i(x) = A_i x + c_i \qquad \forall x \in X_i,$$

where $A_i : \mathbb{R}^D \to \mathbb{R}^D$ is a linear contraction and $c_i \in X$.
The simplest example of a piecewise contracting map of this class is given in dimension 1 on the interval $[0, 1]$ by $f(x) = ax + c$ with $a, c \in (0, 1)$ (here $X_1 = (0, (1-c)/a)$, $X_2 = ((1-c)/a, 1)$ and $\Delta = \{0, 1, (1-c)/a\}$). For almost all the values of the parameters $a$ and $c$, the global attractor $\widetilde{\Lambda}$ is composed of a unique periodic orbit (whose period depends on the specific values of the parameters), but for the remaining set of parameters the global attractor is a Extremal Index - EI set supporting a minimal dynamics. We refer to [248] for a detailed description of the asymptotic dynamics of this map. We just quoted above the case with $c = 0$ and 0 the (unique) fixed point.

In the unit square $X = [0, 1]^2$, another interesting example is given by a map with four pieces $X_1 := (0, T_1) \times (0, T_2)$, $X_2 := (T_1, 1) \times (0, T_2)$, $X_3 := (0, T_1) \times (T_2, 1)$ and $X_4 := (T_1, 1) \times (T_2, 1)$, singular set $\Delta$ given by the two segments $\{x = T_1\}$, $\{y = T_2\}$ and the boundary of the unit square, and restricted maps $f_1(x, y) := a(x, y) + (1-a)(1, 0)$, $f_2(x, y) := a(x, y) + (1-a)(1, 1)$, $f_3(x, y) := a(x, y)$ and $f_4(x, y) := a(x, y) + (1-a)(0, 1)$. Once again, for almost all values of the parameter $a \in [0, 1)$, $T_1$ and $T_2$, the global attractor is composed of a finite number of periodic orbits. However, contrarily to the 1 dimensional case, the number of periodic orbits increases with the parameter $a$ and tends to infinity when $a$ goes to 1. This example belongs to a larger class of higher dimensional piecewise affine contracting maps which are models for genetic regulatory networks and are studied in [250]. We would like to conclude this section with two comments:

1) We stress that the first part of Theorem 7.4.5 applies to any map verifying the assumptions in Definition 7.4.4, not necessarily to PCM, provided we take the target point $z$ outside the global attractor. As an example we could quote the baker's transformation defined on the unit square, see [247] for a detailed description of it.

2) Piecewise contacting maps could also be perturbed with additive noise. In [239], a numerical analysis was performed. In this situation it is difficult to get rigorous results since we do not dispose of useful tools like the Fourier series technique for rotations. Nevertheless we could show a convergence toward a Gumbel's law together with the presence of an EI less than 1 for some periodic target points $z$. It is remarkable that for PCM the two kinds of perturbations considered, RASP and additive noise, do not smoothen out completely the periodicity features of the unperturbed map, as this is reflected in the persistency of an extremal index less than one.

## 7.5
## Observational Noise

A different type of perturbation is given by the *observational noise*; it affects the observations of the orbits of a dynamical systems $(\mathcal{M}, f, \mu)$, where $\mu$ is $f$-invariant, but does not affect the dynamics itself (here again $\mathcal{M}$ is a compact subset of $\mathbb{R}^D$ or on the $D$-dimensional torus). Precisely, it consists in replacing the orbit $(f^n(x))_{n \in \mathbb{N}}$ of a point $x \in \mathcal{M}$, by the sequence $(y_n)_{n \in \mathbb{N}}$ defined by

$$y_n := y_n(x, \overline{\pi}) := f^n(x) + \varepsilon \pi_n \quad \forall n \in \mathbb{N}, \tag{7.5.1}$$

where $\varepsilon > 0$ and $\pi_n = \Pi_n(\overline{\pi})$, $n \in \mathbb{N}$, where $\Pi_n$ projects on the $n$-th component of $\overline{\pi} = (\pi_0, \pi_1, \ldots, \pi_n, \ldots)$, and these components are i.i.d. random vectors, which take values in the hypercube $\Omega_1$ of $\mathbb{R}^D$ centered at 0 and of side 2, and with common distribution $\theta$, which we choose absolutely continuous with density $\rho \in L^\infty_{\text{Leb}}$, namely, $d\theta(\pi) = \rho(\pi)d\text{Leb}(\pi)$, with $\int_{\Omega_1} \rho(\pi)d\text{Leb}(\pi) = 1$[5]).

An orbit perturbed with observational noise mimics the behavior of an instrumental recorded time series. Instruments characteristics, defined as precision and accuracy, act both by truncating and randomly displacing the real values of a measured observable. It has also been shown in [251] that in the computation of some statistical quantities, the dynamical noise corresponding to the random transformations described above, could be considered as an observational noise with the Cauchy distribution.

We now take a function $\phi$ as in (7.2.7), and we define the stochastic process

$$X_0(x, \overline{\pi}) = \phi(x + \varepsilon \Pi_0 \overline{\pi}), \, X_1(x, \overline{\pi}) = \phi(f(x) + \varepsilon \Pi_1 \overline{\pi}), \ldots, X_n(x, \overline{\pi}) = \phi(f^n(x) + \varepsilon \Pi_n \overline{\pi})$$

It is easy to see that this process is stationary with respect to the product measure $\mathbb{P} = \mu \times \theta^{\mathbb{N}}$, defined on the product space $\mathcal{M} \times \Omega_1^{\mathbb{N}}$ with the product $\sigma$-algebra.

If the observable $\phi$ is chosen as $\phi(x) = -\log(\text{dist}(x, \zeta))$, being dist a distance on the metric space $\mathcal{M}$, then we have the analogous of Proposition 7.2.4, namely:

**Proposition 7.5.1** ([65]). *Let us consider the dynamical systems $(\mathcal{M}, f, \mu)$; we perturb it with observational noise and we consider the associated process $X_n(x, \overline{\pi}) := -\log(||f^n x + \varepsilon \pi_n) - \zeta||)$ endowed with the probability $\mathbb{P} = \mu \times \theta^{\mathbb{N}}$. We suppose moreover that $\theta$ is the Lebesgue measure measure on $S$. Then the linear sequence $u_n := y/a_n + b_n$ defined in (2.2.2) verifies*

$$a_n = d \quad and \quad b_n = \frac{1}{D} \log\left(\frac{n \, K_D \, \mu(B_\varepsilon(\zeta))}{(2\varepsilon)^D}\right). \tag{7.5.2}$$

where $K_D$ is the volume on the unit hypercube in $\mathbb{R}^D$ and $\zeta$ belongs to the support of the measure $\mu$.

The next step will be to provide class of maps $f$ for which conditions $\boxed{\Pi}_0(u_n)$ and $\boxed{\Pi}'_0(u_n)$ are satisfied; these maps are very similar to those given in the items **S1** and **S2** above. We call this class:

• **SO**. They are defined on the torus $\mathcal{M} = \mathbb{T}^D$ with the norm $||\cdot||$ and satisfy:

1) There exists a finite partition (mod-0) of $\mathcal{M}$ into open sets $Y_j$, $j = 1, \cdots, p$, namely $\mathcal{M} = \cup_{j=1}^p \overline{Y_j}$, such that $T$ has a Lipschitz extension on the closure of each $Y_j$ with a uniform and strictly larger than 1 Lipschitz constant $\eta$, $||T(x) - T(y)|| \leq \eta ||x - y||$, $\forall x, y \in \overline{Y_j}$, $j = 1, \cdots, p$.

2) We will suppose that $f$ preserve a Borel probability measure $\nu$ which is also mixing with decay of correlations given by

$$\left| \int f \circ f^m h d\nu - \int f d\nu \int h d\nu \right| \leq C ||h||_{\mathcal{B}} ||f||_1 m^{-2} \tag{7.5.3}$$

---

[5]) Each $\pi$ is a vector with $d$ components; all these components are independent and distributed withnoise common density $\rho'$; the product of such marginals $\rho'$'s gives $\rho$.

where the constant $C$ depends only on the map $T$, $||\cdot||_1$ denotes the $L_\nu^1$ norm with respect to $\nu$ and finally $\mathcal{B}$ is a Banach space included in $L_{\mathrm{Leb}}^\infty$. We will also need that $\nu$ be equivalent to Leb with density in $L_{\mathrm{Leb}}^\infty$.

**Proposition 7.5.2** ( [65]). *Let us suppose that our dynamical systems verifies Assumption* **SO** *and it is perturbed with* observational noise *defined above. Then conditions* $\bar{\Pi}_0(u_n)$ *and* $\bar{\Pi}'_0(u_n)$ *hold for the observable* $\phi$.

An interesting application of the two previous Propositions is obtained when the distribution of the additive noise is exactly the Lebesgue measure. For the given magnitude of the noise $\varepsilon$ we could take as first approximation to the invariant measure of the ball $\nu(B_\varepsilon(\zeta,\ ) \approx \varepsilon^{d(\zeta)}$, where $d(\zeta)$ is the local dimension defined in Eq. 4.2.8; see a more detailed discussion of such a scaling in the Chapter 8, when noise is, instead, absent. Then the linear scaling parameter $b_n$ of the EVL is expected to behave as

$$b_n \sim \frac{1}{D} \, \log(n\epsilon^{d(\zeta)-D}). \qquad (7.5.4)$$

where $D$ is the ambient space dimension. Therefore we have an useful technique to detect the local dimension of the measure, with a *finite* resolution given by the strength of the noise. This will also allow us to compute directly the distribution of the maxima with the affine linearization given by the explicit expression of the $u_n$. This would be particularly useful whenever the invariant measure is singular and therefore the GEV distribution, see Section 1.1, does not admit a probability density function: see also [77, Section 3.1] for a detailed discussion on this point. Of course the rigorous applicability of these arguments is ensured for the moment for systems verifying the two previous propositions; nevertheless we investigated the scaling (7.5.4) for more general systems like one-dimensional repellers with a cantorian structure and strange attractors. In both cases the dimension $D$ computed with formula (7.5.4) is in good agreement with the value of the Kaplan-Yorke dimension. We defer to Chapter 9 for a detailed description of these numerical experiments and also to other applications of the observational noise, for instance the possibility to discriminate between highly recurrent and sporadic points. This is based again on formula (7.5.4) and the simple observation that if the ball around $\zeta$ is visited with less frequency, the local density is of lower order with respect to $\varepsilon$, which means that one should go to higher values of $n$ in order to have a reliable statistics, namely a good convergence for $b_n$. As we argued in [65] the main advantage of studying recurrence properties in this way over applying other techniques is due to the built-in test of convergence of this method: even for a point rarely recurrent there will be a time scale $n$ such that the fit converges. In [80] we used this technique to define rigorous recurrences in long temperature records collected at several weather stations in Europe.

## 7.6
## Non-stationarity – the Sequential Case

We now discuss briefly the case when the process $X_0, X_1, \ldots$ is not stationary: we stress that stationarity plays an important role to establish extreme value statistics. On the other hand, Hüsler [252] gave a non-stationary version of EVL in the case of non-stationary sequences without identical marginal distributions. In the domain of non-autonomous dynamical systems it is natural to relax the stationary assumption. In a forthcoming paper [253] we will present an approach to this problem with applications to the class of sequential maps and we will now briefly mention these results. The starting point is to provide a generalisation of Hüsler's results in a pure probabilistic setting and with an adjustment of the dependence conditions $\not\!\!\Delta_q(u_n)$ and $\not\!\!\Delta'_q(u_n)$ introduced in Chapter 4.

By using the notations introduced in the previous sections, we state that our main goal is now to determine the limiting law of

$$P_n = \mathbb{P}(X_0 \leq u_{n,0}, X_1 \leq u_{n,1}, \ldots, X_{n-1} \leq u_{n,n-1})$$

as $n \to \infty$, where $\{u_{n,i}, i \leq n - 1, n \geq 1\}$ is considered a real-valued boundary. Define $x_{0,i} = \sup\{x : F_i(x) < 1\}$ and let $F_i(x_{0,i}-) = 1$ for all $i$ and assume that

$$\bar{F}_{\max} := \max\{\bar{F}_i(u_{n,i}), i \leq n - 1\} \to 0 \text{ as } n \to \infty, \tag{7.6.1}$$

which is equivalent to

$$u_{n,i} \to x_{0,i} \text{ as } n \to \infty, \text{ uniformly in } i,$$

where $\bar{F}_i(x) = 1 - F_i(x)$ for all $i$.

Let us denote $F_n^* := \sum_{i=0}^{n-1} \bar{F}_i(u_{n,i}),$. and assume that there is $\tau > 0$ such that

$$F_n^* := \sum_{i=0}^{n-1} \bar{F}_i(u_{n,i}) \to \tau. \tag{7.6.2}$$

The main result in the first part of the forthcoming paper [253] is to show that for $X_0, X_1, \ldots$ a non-stationary stochastic process such that (7.6.1) and (7.6.2) hold for some $\tau > 0$,then as long as some specially adapted conditions $\not\!\!\Delta_q(u_n)$ and $\not\!\!\Delta'_q(u_n)$ hold then

$$\lim_{n \to \infty} P_n = e^{-\theta\tau},$$

where $\theta$ is also defined as before with necessary adjustments to this non-stationary setting.

Let us give now an example where the previous theorem applies; since we are interested in connection with dynamical systems, we could get a situation which fits with the previous considerations by taking, as before, a concatenation of maps but which are now chosen without any distribution for the $\omega_k$, namely we take the maps $(f_{\omega_k})_{k \geq 0}$ in some set and compose them. Let us consider for simplicity a distance observable of the form (7.2.7) with $\phi(\cdot) = -\log(\text{dist}(\zeta, \cdot))$, and take the very particular dynamical systems given by $\beta$-transformations. Let us call $T(x) = T_\beta(x) = \beta x$ mod 1, the original unperturbed $\beta$-transformation and take the other transformations of the same kind $x \to \beta_k x$ mod 1, with $\beta_k \geq 1 + a, \forall k \geq 1$, where $a$ is a given positive number and moreover $|\beta_n - \beta| \leq \frac{1}{n^\iota}$, with $\iota > 1$. Notice that the process $X_k(\cdot) = \phi(f_{\omega_k} \circ \cdots \circ f_{\omega_1}(\cdot))$, equipped with the probability $\mathbb{P}$ given by the Lebesgue measure Leb is not necessarily stationary nor independent.

By using the theory of sequential $\beta$-transformations developed in [254], we can apply the generalisation of Extreme Value Theory to non-stationary sequences obtained in the first part of [253] and actually verify the adapted conditions $\underline{\amalg}_q(u_n)$ and $\underline{\amalg}'_q(u_n)$. We can also obtain that, when $\zeta$ is periodic, explicit expressions for the EI. It is interesting to notice that even in this sequential setting the extremal indices could be smoothed out. Let us consider for instance $\beta = 5/2$ and $T = T_\beta = 5/2x$ mod 1. Let $\zeta = 2/3$. Note that $T(2/3) = 2/3$. Consider a sequence $\beta_n = 5/2 + \varepsilon_n$, with $\varepsilon_n = n^{-\alpha}$, where $\alpha < 1$. Note that $1/n = o(\varepsilon_n)$. It is possible to show that at $\zeta = 2/3$ although for the unperturbed system $T$ shows an EI equal to $1 - 2/5 = 3/5$, for the sequential systems chosen as above the EI is equal to 1.

# 8
# A Statistical Mechanical Point of View

## 8.1
## Choosing a Mathematical Framework

In this chapter, we want to present some results of practical significance for studying the extremes of complex systems featuring chaotic behavior. We take a different point of view with respect to the previous chapters of the book, as we will sacrifice some mathematical rigour and try to address the properties of *typical* - to be understood below - physically relevant systems, having the statistical mechanical perspective of envisioning high-dimensional dynamical systems.

We focus our attention on Axiom A systems [70], which are a special class of dynamical systems possessing a Sinai-Ruelle-Bowen (SRB) invariant measure [161] and featuring uniform hyperbolicity in the attracting set. Such invariant measure coincides with the Kolmogorov's physical measure, *i.e.* it is robust against infinitesimal stochastic perturbations. Another important property of Axiom A systems is that it is possible to develop a response theory for computing the change in the statistical properties of any observable due to small perturbations to the flow [83, 84]. Such a response theory has recently been the subject of intense theoretical [255, 256], algorithmic [257, 258, 259], and numerical investigations [260, 261, 262, 263] and is gaining prominence especially for geophysical fluid dynamical applications. Moreover, the response theory seems to provide powerful tools for studying multiscale systems and deriving parametrizations of the impact of the fast variables on dynamics of the slow variables [264, 265].

Finally, most importantly, Axiom A systems are a good paradigm to describe high dimensional systems, which can be described by statistical mechanics. While the dynamics of natural or artificial systems are most often definitely not Axiom A, Axiom A systems can serve as good *effective* and tractable models already revealing features of non-equilibrium statistical mechanics, that is features which should hold in a properly interpreted manner for a much larger and more realistic class of systems. In other words, we subscribe to the so-called chaotic hypothesis, which is somewhat the equivalent in the non-equilibrium framework of the classic ergodic hypothesis for equilibrium dynamics [266]. Moreover, as discussed in [263], when we perform numerical simulations we effectively assume implicitly that the system under inves-tigation is in large parts like an Axiom A or very similar to an Axiom A system. In order to avoid misconception, this does not mean that the story ends with Axiom A systems, quite the contrary, it just shows the enormous relevance to extend the theory of chaotic dynamical system further and further, which is a mathematical challenge under active investigation [267, 268].

In this chapter we concentrate on Axiom A systems, which seems, in light of the previous discussion, a good mathematical framework to provide insights and results useful for a large spectrum of applications in statistical mechanics and physics in general.

The strong technology developed for Axiom A systems is instrumental in the derivation of various results on the relationship of parameters describing the extremes to the dynamical and geometrical properties of the system, and will allow addressing the problem of the sensitivity of extremes to small perturbations of the system in a general conceptual way. The dependence of the properties of extremes of parametric modulations of the underlying dynamics is an issue of relevant theoretical as well as applicative interest. The practical interest stems from the fact that it is relevant to be able to control or predict variations in extreme events due to small perturbations to the dynamics to quantify for example model errors. The theoretical interest comes from the fact that when considering extremes, universal parametric probability distributions can be defined, as opposed to the case of the bulk statistical properties. Because of this, we may hope to reconstruct the parameters descriptive of the extremes from simple moments of the distributions, express these in terms of observables of the system, and use the Ruelle response theory for expressing rigorously the sensitivity of extremes to small perturbations to the dynamics. Our construction will be based on the POT method and will lead to deriving explicitly the parameters of the corresponding GPD, but our results can be pulled back to the equivalent GEV formulation.

## 8.2
## Generalized Pareto Distributions for Observables of Dynamical Systems

Let us consider a continuous-time Axiom A dynamical system $\dot{x} = G(x)$ on a compact manifold $\mathcal{X} \subset \mathbb{R}^D$ (phase space), where $x(t) = f^t(x_{in})$, with $x(t = 0) = x_{in} \in \mathcal{X}$ initial condition and $f^t$ evolution operator, is defined for all $t \in \mathbb{R}_{\geq 0}$. Let us define $\Omega$ as the attracting invariant set of the dynamical system, so that $\mu$ is the associated SRB measure supported in $\Omega = \text{supp}(\mu)$. We consider two different classes of function mapping $\mathcal{X}$ to $\mathbb{R}$, the so-called *distance observables* and the so-called *physical observables*, which we have already encountered in Sect. 4.2.1 and Sect. 4.6, respectively, and briefly recapitulated below.

## 8.2.1
## Distance Observables

Distance observables can be expressed as functions $g : \mathcal{X} \to \mathbb{R} \cup \{+\infty\}$ written as $g(r)$, with $r = \text{dist}(x(t), \zeta) \geq 0$, where $\zeta \in \Omega$ is a reference point belonging to the attractor. In Section 4.2.1 we have already introduced the three observables $g_i$, $i = 1, 2, 3$:

$$g_1(r) = -\log(r) \tag{8.2.1}$$

$$g_2(r) = r^{-1/\alpha}, \alpha > 0 \tag{8.2.2}$$

$$g_3(r) = -r^{1/\alpha}, \alpha > 0 \tag{8.2.3}$$

We analyze exceedances of $g$ (chosen among the $g_i$'s, $i = 1, 2, 3$, given above) above a certain threshold $T$. Due to the invertibility of the function $g$, the threshold $T$ is in one to one correspondence to a radius $r^*$, namely $T = g(r^*)$. An above-threshold event happens every time the distance between the orbit of the dynamical system and $\zeta$ is smaller than $r^*$. See Figure 8.1 for clarification. In order to address the problem of extremes, we have to consider balls of small radii. Therefore, we denote the exceedances above $T$ by $z = g(r) - T$. That is, if at time $t$ the dynamical systems is at $x(t)$, then we have an exceedance $g(\text{dist}(x(t), \zeta)) - T$, if this expression is non-negative. The number of exceedances above $z + T$ relative to the number of exceedances above $T$ up to time $t$ can be written as:

$$\frac{\int_0^t \Theta(g(\text{dist}(x(s), \zeta)) - T - z) ds}{\int_0^t \Theta(g(\text{dist}(x(s), \zeta)) - T) ds} \tag{8.2.4}$$

or discrete version of this expression, if discrete dynamics is considered. Hence, by ergodicity of the system we can express this ratio for a large observation window in time by the ergodic measure $\mu$. If we choose the starting point with respect to the invariant measure $\mu$ then the dynamics become a stationary process and if we denote its law by $\mathbb{P}$ we can express the ratio as

$$\mathbb{P}(g(\text{dist}(x, \zeta)) < z + T | g(\text{dist}(x, \zeta)) < T) = \frac{\mathbb{P}(g(\text{dist}(x, \zeta)) < z + T)}{\mathbb{P}(g(\text{dist}(x, \zeta)) r < T)}, \tag{8.2.5}$$

In terms of the invariant measure $\mu$ of the system, we have that the probability $H_{g,T}(z)$ of observing an exceedance of at least $z$ above $T$ given that an exceedance above $T$ has occured is given by:

$$H_{g,T}(z) \equiv \frac{\mu(B_{g^{-1}(z+T)}(\zeta))}{\mu(B_{g^{-1}(T)}(\zeta))}. \tag{8.2.6}$$

Obviously, the value of the previous expression is 1 if $z = 0$. In agreement with the conditions given on $g$, the expression contained in Eq. 8.2.6 monotonically decreases with $z$ and vanishes when $z = z^{max} = g^{max} - T$. Note that the corresponding d.f. is given by $F_{g,T}(z) = 1 - H_{g,T}(z)$, so that, using the convention introduced in Chapter 2, we write $\bar{F}_{g,T}(z) = H_{g,T}(z)$. Equation 8.2.6 clarifies that we have translated the computation of the probability of above-threshold events into a geometrical problem.

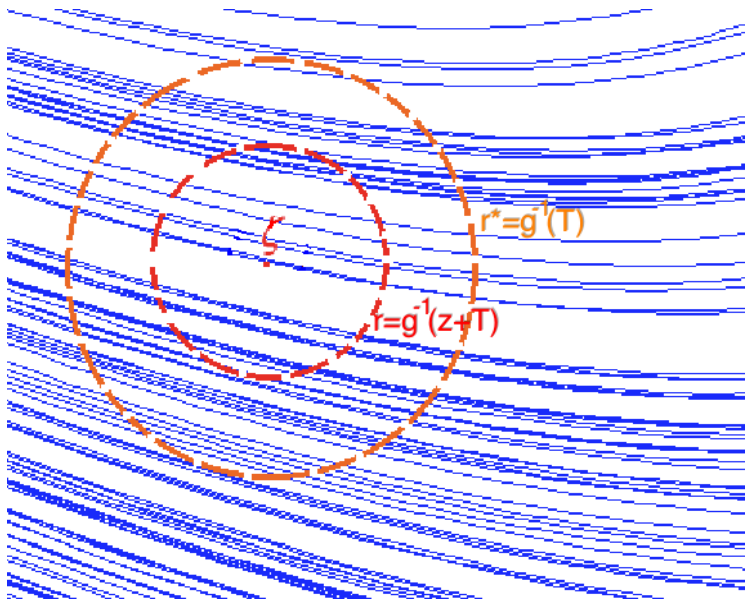

**Figure 8.1** If we consider a distance observable of the form $g = g(r)$, where $r$ is the distance from $\zeta$ and $g$ is monotonically decreasing with a maximum for $r = 0$, events above the threshold $T$ are given by close returns of the orbit near $\zeta$, at a distance $r$ from $\zeta$ smaller than $r^* = g^{-1}(T)$. The conditional probability of exceedances is constructed by taking the ratio of the mass of the attractor contained inside a sphere of radius $r$ divided by the mass of the attractor contained inside a sphere of radius $r^*$. For more details, cf. the main text.

For Axiom A systems one has that the local dimension around $\zeta$ introduced in Eq. 4.2.8 and given by $d(\zeta) = \lim_{r\to 0}(\log(\mu(B_r(\zeta)))\log(r))$ is such that $d(\zeta) = d_H$ almost everywhere on the attractor [269, 270], where $d_H$ is the Hausdorff dimension, and we additionally have that $d_H = d_q, \forall q \in \mathbb{N}$, where $d_q$ are the generalised Renyi dimensions [71]. Moreover, we follow the conjecture that $d_H = d_{KY}$, where $d_{KY}$ is the Kaplan-Yorke dimension [71] (note that this conjecture is not essential for most of the considerations in this chapter):

$$d_{KY} = n + \frac{\sum_{k=1}^{n} \lambda_k}{|\lambda_{n+1}|}, \qquad (8.2.7)$$

where the $\lambda_j$'s are the Lyapunov exponents of the systems, ordered from the largest to the smallest, $n$ is such that $\sum_{k=1}^{n} \lambda_k$ is positive and $\sum_{k=1}^{n+1} \lambda_k$ is negative.

Following [269], we can also write $d_H = d_u + d_n + d_s$, where $d_s$, $d_u$ and $d_n$ are the dimensions of the attractor $\Omega$ restricted to the stable, unstable and neutral directions, respectively, at the point $x = x_0$. We have that $d_u$ is equal to the number of positive Lyapunov exponents $\lambda_j^+$, $d_u = \#(\{\lambda_j > 0, j = 1, \ldots, d\})$, the dimension $d_n$ is unitary for Axiom A systems, while $d_s$ is given by $d_s = d_H - d_u - d_n$. Note that, if $d_H = d_{KY}$, it follows that $\{d_s\} = d_s - \lfloor d_s \rfloor = \sum_{k=1}^{n} \lambda_k/|\lambda_{n+1}|$, because the last term gives a positive contribution smaller than 1, and $d_u$ and $d_n$ are both integer numbers.

In order to proceed with the derivation of an EVL, the asymptotics provided in Eq. 4.2.8 is not sufficient to guarantee the convergence of the expression in (8.2.6). One needs to assume that

$$\mu(B_r(\zeta)) = f_\zeta(r) r^{d_H}, \tag{8.2.8}$$

where $f_\zeta(r)$ is a slowly varying function of $r$ as $r \to 0$ possibly depending on $\zeta$, *i.e.*, for small enough $s$, we have that $\lim_{r \to 0} \frac{f_\zeta(sr)}{f_\zeta(r)} = 1$. Two functions $f$ and $g$ as such that $f \sim g$ for $x \to y$ if $\lim_{x \to y} f(x)/g(x) = 1$. An equivalent way to express that $f$ is slowly varying is to say that for small enough $s$ one has that $f(st) \sim f(t)$ at $t \to 0$. Hence, if $f \sim g$ and if $f$ is slowly varying, then also $g$ is slowly varying.

Inserting the expression 8.2.8 in Eq. 8.2.6 we obtain the following expression for the tail probability of exceedance:

$$\bar{F}_{g,T}(z) \sim \left( \frac{g^{-1}(z+T)}{g^{-1}(T)} \right)^{d_H}, \tag{8.2.9}$$

where the slowly varying terms coming from $f_\zeta(r)$ cancel out, as a result of the property of $f$ mentioned above. By replacing $g$ with the specific observables we are considering, we obtain explicitly the corresponding distribution for the extremes.

By choosing an observable of the form given by either $g_1(r) = -\log(x)$, $g_2(r) = x^{-1/\alpha}$, or $g_3(x) = C - x^{1/\alpha}$, we derive as EVL a member of the GPD family. In fact, substituting $g_1^{-1}(y) = \exp(-y)$, $g_2^{-1}(y) = y^{-\alpha}$, and $g_3^{-1}(y) = (C-y)^\alpha$ into Eq. 8.2.9 and comparing with Eq. 3.1.8, we derive that

$$F_{g_i,T}(z) = GPD_\xi \left( \frac{z}{\sigma} \right),$$

where the parameters for the different cases are:

- $g_1$-type observable:

$$\sigma = \frac{1}{d_H}, \qquad \xi = 0; \tag{8.2.10}$$

- $g_2$-type observable:

$$\sigma = \frac{T}{\alpha d_H}, \qquad \xi = \frac{1}{\alpha d_H}; \tag{8.2.11}$$

- $g_3$-type observable:

$$\sigma = \frac{C-T}{\alpha d_H}, \qquad \xi = -\frac{1}{\alpha d_H}. \tag{8.2.12}$$

The previous expressions show that there is a simple algebraic link between the parameters of the GPD and the Hausdorff dimension of the attractor, and suggest multiple ways to extract it from different outputs of statistical inferences procedures on data. Assuming Axiom A dynamics allows to derive global results from the recurrence properties of almost every point of the attractor. One needs to note that assuming

Axiom A property is relevant only for relating local measurements to global properties. Additionally, we remind that in the case of Axiom A systems the conditions $\mathit{Д}_0(u_n)$ and $\mathit{Д}'_0(u_n)$ discussed in Chapter 4 are automatically obeyed for the considered observables $g_i$'s, $i = 1, 2, 3$. This is fundamentally why the results presented in this section agree with what obtained in Sect. 4.2.1 using the GEV approach, at least when one does not consider as $\zeta$ one of the period points of the attractor.

Let's now for step out of the comfortable Axiom A world. In this case, $d(\zeta)$ is not necessarily constant a.e. in $\Omega$, but the formulas in Eqs. 8.2.8-8.2.12 are still valid, *modulo* substituting $d_H$ with $d(\zeta)$. Hence the analysis of extremes of $g_i(dist(x, \zeta))$, $i = 1, 2, 3$ provides information on the local fine structure of the attractor, and, in particular, on the local dimension $d(\zeta)$ around the point $\zeta$. Note that nowhere in the derivation we need to introduce conditions on the decay of correlations of the observables, because our construction is purely geometrical, and no dynamics is involved (except assuming ergodicity). Interestingly, the GPD approach can be used also when studying regular systems or systems featuring slow decay of correlations for the considered observables [78].

Instead, as discussed in the previous Chapters, deriving corresponding EVLs following the GEV approach requires assuming $\mathit{Д}_0(u_n)$ and $\mathit{Д}'_0(u_n)$ conditions for the observable $g(dist(x, \zeta))$ [77]. If such conditions apply, the GPD and GEV points of view on the extremes are, indeed, equivalent.

## 8.2.2
## Physical Observables

Physical observables have been introduced in Sect. 4.6 and some examples of relevance in, *e.g.*, fluid dynamics are given in Eq. 4.6.1. The physical observables can be written as functions $A : \Omega \to \mathbb{R}$ whose maximum restricted to the support of $\mu$ is unique, so that is there exists a unique $x_0 \in \Omega$ such that $A_{\max} \equiv \max(A)|_\Omega = A(x_0)$. Moreover, we assume that $x_0$ is not a critical point, so that $\nabla A|_{x=x_0} \neq 0$, where the gradient is taken in $\mathbb{R}^D$. Therefore, we have that the neutral manifold and the unstable manifold have to be tangent to the manifold $\{A(x) = A_{\max}\} = \{x \in \mathcal{X} : A(x) = A_{\max}\}$ in $x = x_0$. We also have that the intersection between the manifolds $\{x \in \mathcal{X} : A(x) = \tilde{A}\}$ and $\Omega$ is the empty set if $\tilde{A} > A_{max}$. We define as $\tilde{\Sigma}^T_{A_{\max}}$ the subset of $\mathbb{R}^D$ included between the manifolds $\{x \in \mathcal{X} : A(x) = A_{\max}\}$ and $\{x \in \mathcal{X} : A(x) = T\}$.

Furthermore, we define as $\Sigma^T_{A_{\max}}$ the subset of $\mathbb{R}^D$ included between the hyperplane $\beta_{max}$ tangent to the manifold $A(x) = A_{\max}$ in $x = x_0$, that is $\{x \in \mathcal{X} : (x - x_0)\hat{n} = 0\}$, where $\hat{n} = \nabla A|_{x=x_0}/|\nabla A|^2_{x=x_0}$ is the normal vector of length $|\nabla A|^{-1}_{x=x_0}|$, and the hyperplane $\beta_T$, which is obtained by applying the translation given by the vector $(T - A_{\max})\hat{n}$ to the hyperplane $\beta_{\max}$, in other words all points with the distance $(T-A_{\max})$ from $\beta_{\max}$. Hence, $\Sigma^T_{A_{\max}} = \{x \in \mathcal{X} : T - A_{\max} \leq (x - x_0)\hat{n} \leq 0\}$.

As $T \to A_{\max}$, which is the limit we are interested in, we have that $\tilde{\Omega}^T_{A_{\max}} = \Omega \cap \tilde{\Sigma}^T_{A_{\max}}$ and $\Omega^T_{A_{\max}} = \Omega \cap \Sigma^T_{A_{\max}}$ become undistinguishable up to leading

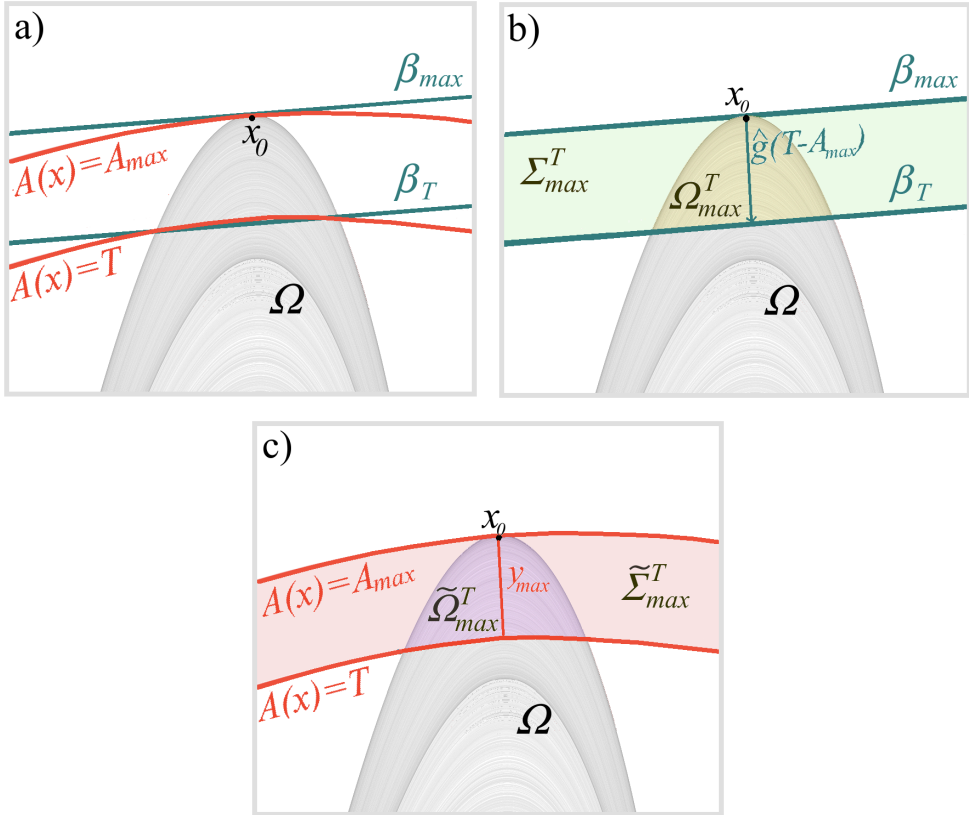

**Figure 8.2** A low-dimensional cartoon of the geometrical construction used for deriving the EVL for exceedances above the threshold $T$ for the observable $A(x)$ such that $\max(A)|_\Omega = A_{max}$ is realized for $x = x_0$. a) The manifolds $A(x) = A_{max}$ and $A(x) = T$ are depicted, together with the attracting invariant set $\Omega$ and the two hyperplanes $\beta_{max}$ and $\beta_T$. $\beta_{max}$ is tangent to $A(x) = A_{max}$ in $x_0$ and $\beta_T$ is obtained from $\beta_{max}$ via translation along $(T - A_{max})\hat{g}$. b) The hyperplanes $\beta_{max}$ and $\beta_T$ delimit the region $\Sigma_{max}^T$. Its intersection with $\Omega$ is $\Omega_{max}^T$. c) The manifolds $A(x) = A_{max}$ and $A(x) = T$ delimit the region $\bar{\Sigma}_{max}^T$. Its intersection with $\Omega$ is $\bar{\Omega}_{max}^T$. As $T \to A_{max}$, we have that $\Omega_{max}^T \to \bar{\Omega}_{max}^T$.

order. See Fig. 8.2 for a depiction of this geometrical construction.

More in general, we denote as $\Omega_V^U$, with $V > U > T$, the intersection between $\Omega$ and the subset of $\mathbb{R}^D$ included between the hyperplane $\beta_U$ and $\beta_V$, where $\beta_X$, $X = U, V$, is obtained from $\beta_{max}$ by applying to it the translation given by the vector $(X - A_{\max})\hat{n}$, that is $\Omega_V^U = \{x \in \mathcal{X} : U - A_{\max} \le (x - x_0)\hat{n} \le V - A_{\max}\}$.

It is now clear that we observe an exceedance of the observable $A(x)$ above $T$ each time the systems visits a point belonging to $\Omega_{A_{max}}^T$. In more intuitive terms, and taking the linear approximation described above, an exceedance is realized each time the system visit a point $x \in \Omega$ whose distance $\mathrm{dist}(x, \beta_{max})$ from the hyperplane $\beta_{max}$ is smaller than $y_{max} = (A_{max} - T)/|\nabla A|_{x=x_0}|$.

We define the exceedances $Z$ above $T$ as all points $x \in \Omega$ with distance $y = \mathrm{dist}(x, \beta_T) = Z/|\nabla A|_{x=x_0}|$ from $\beta_T$, and the maximum exceedance $A_{\max} - T$ corresponds to a distance $(A_{\max} - T)/|\nabla A|_{x=x_0}|$ between $x_0$ and $\beta_T$.

As $P(z > Z|z > 0) = P(z > Z)/P(z > 0)$, the probability $\bar{F}_T(Z)$ of observing an exceedance of at least $Z$ conditionally that an exceedance occurs is given by:

$$\bar{F}_T(Z) \equiv \frac{\mu(\Omega_{A_{\max}}^{T+Z})}{\mu(\Omega_{A_{\max}}^{T})}, \tag{8.2.13}$$

where we have used the same arguments as in the case of distance observables based on the ergodicity of the system. Obviously, the value of the previous expression is 1 if $Z = 0$. The expression contained in Eq. (8.2.13) monotonically decreases with $Z$ (as $\Omega_{A_{\max}}^{T+Z_2} \subset \Omega_{A_{\max}}^{T+Z_1}$ if $Z_1 < Z_2$) and vanishes when $Z = A_{\max} - T$.

Let us now estimate $\mu(\Omega_{A_{\max}}^{T})$ as a function of $y_{\max}$ in the case of generic quadratic tangency between the hyperplane $A(x) = A_{\max}$ and the unstable manifold in $x = x_0$.

### 8.2.3
### Derivation of the Generalised Pareto Distribution Parameters for the Extremes of a Physical Observable

We now wish to understand how to express the numerator and denominator of Eq. 8.2.13 as a function of $T$, $Z$, and $A_{max}$.

Following [44], we build upon the construction proposed by [81] and discussed thoroughly in Sect. 6.10. We derive the result by considering the following heuristic argument. Near $x_0$, the attractor could be seen as the cartesian product of a multi-dimensional paraboloid (of dimension $d_u + d_n$) times a fractal set of dimension $d_s$ immersed in $\mathbb{R}^{D-d_u-d_n}$. Note that this excludes, *e.g.*, conservative chaotic systems, whose attractor has the same dimension of the phase space, and systems that can be decomposed into a conservative part and a purely contractive part, whose attractor also has integer dimension. The mass of the cut of the paraboloid contained in $\Omega_{A_{\max}}^{Z+T}$ with the stable and the neutral direction is $\sim h_{u,x_0}(r)r^{d_u+d_n}$, where $r$ is the distance from the minimum and $h_{u,x_0}(r)$ is a slowly varying function of $r$ as $r \to 0$ for an SRB-measure $\mu$. Instead, the cut of the paraboloid with the unstable directions has asymptotically a mass of corresponding fractal set $\sim h_{s,x_0}(l)l^{d_s}$, where $l$ is the distance along the cartesian projection and $h_{s,x_0}(l)$ is a slowly varying function of $l$ as $l \to 0$, possibly depending on $x_0$. The existence of such a function is an additional requirement, cf. the discussion in the case of the distance variable.

In our case, $l = \gamma y_{\max}$ and $r = \kappa\sqrt{y_{\max}}$, which results from the functional form of the paraboloid, see also Fig. 8.2. The first asymptotic follows from the fact that $\beta_{\max}$ has to be tangential to the stable and the neutral direction in order that $x_0$ can be a maximum of $A$ in $\Omega$. The second asymptotic follows, if the direction $\hat{n}$ is not exactly pointing in a direction particular for the fractal structure in the stable direction, which should hold generically. Hence, we obtain that $\mu(\Omega_{A_{\max}}^{T}) \sim \tilde{h}_{x_0}(y_{\max})y_{\max}^{\delta}$, where

$$\delta = d_s + (d_u + d_n)/2, \tag{8.2.14}$$

and $\tilde{h}_{x_0}(y_{\max}) = h_{s,x_0}(\gamma y_{\max})$ is also a slowly varying function of its argument.

This construction can be made more formal by considering the disintegration of the SRB measure $\mu$ along the stable and unstable directions of the flow [71]. If the generic orientation of $\hat{n}$ is not generic as in our previous hypothesis, then we have to attribute a factor $1/2$ also to some of the stable directions. In general, the following should hold:

$$(d_s + d_u + d_n)/2 \leq \delta \leq d_s + (d_u + d_n)/2. \tag{8.2.15}$$

As a side note, we emphasize that in the case of more general tangencies between the unstable and neutral manifold, and the manifold $A(x) = A_{max}$, we expect that $\delta = d_s + \sum_{j=1} d_u(2j)/(2j) + d_n(2j)/(2j)$, where $d_u(2j)$ ($d_n(2j)$) gives the number of directions along the unstable (neutral) manifold where a tangency of order $2j$ with the manifold $A(x) = A_{max}$ is found. We obviously have that $\sum_{j=1} d_u(2j) = d_u$.

We continue our discussion considering the case of generic tangency. Following the same argument as above, we have that

$$\mu(\Omega_{A_{max}}^{T+Z}) \sim \tilde{h}_{x_0} \left( \frac{A_{max} - T - Z}{|\nabla A_{x=x_0}|} \right) \left( \frac{A_{max} - T - Z}{|\nabla A_{x=x_0}|} \right)^{\delta}.$$

We define

$$\tilde{z} = 1 - \frac{Z}{A_{max} - T}$$

and obtain $\mu(\Omega_{A_{max}}^{T+Z}) \sim \tilde{h}_{x_0}(\tilde{z} y_{max}) \tilde{z}^{\delta} y_{max}^{\delta}$, where $y_{max} = \frac{A_{max} - T}{|\nabla A_{x=x_0}|}$. Using the definition of slowly varying function and considering Eq. 8.2.13, we derive that in the limit $T \to A_{max}$:

$$\bar{F}_T(Z) \sim \left( 1 - \frac{Z}{A_{max} - T} \right)^{\delta}. \tag{8.2.16}$$

Note that the corresponding d.f. is given by $F_T(Z) = 1 - \bar{F}_T(Z)$. Comparing Eqs. 3.1.8 and 8.2.16, one obtains that $F_T(Z)$ belongs to the GPD family, so that

$$F_T(Z) = GPD_{\xi} \left( \frac{Z}{\sigma} \right),$$

and that the GPD parameters can be expressed as follows:

$$\xi = -1/\delta \tag{8.2.17}$$

$$\sigma = (A_{max} - T)/\delta. \tag{8.2.18}$$

These results generalise what discussed in Sect. 6.10.2 for two specific examples of uniformly hyperbolic systems using the GEV approach, see also [81]. This comes to relatively little surprise, given the properties of Axiom A systems.

It is important to remark that Eq. 8.2.16 has been obtained in the limit of $T \to A_{max}$, and under the assumption that $\mu(B_r(x_0))$ is a regularly varying function of degree D as $r \to 0$. When considering a finite range $A_{max} - T$, one should expect deviations of the empirical distributions of extremes of $A$ from what prescribed in Eq. 8.2.16, which are intrinsic to the fractal nature of the measure. See also discussion and Fig. 1 in [78]. When finite ranges for $A$ are considered, one expects that, in some averaged sense, Eq. 8.2.16 fits well the distributions of extremes of $A$ and Eqs. 8.2.17 and 8.2.18 give the value of the two relevant parameters of the GPD, analogously to the idea that the number of points of the attractor at distance smaller than a small but finite $r$ from the point $x_0$ scales approximately, on the average as $r^D$.

### 8.2.4
### Comments

It is clear that the $\xi_A$ parameter is always negative (so that the distribution of extremes is upper limited), reflecting the fact that the observable is smooth and the attractor is a compact set. Equation 8.2.17 provides a very valuable information, as it shows that the shape parameter $\xi$ of the GPD does not depend on the considered observable so long it fulfils the general conditions given above, but only on the dimensions of the stable and of the unstable manifold. Moreover, the shape parameter is always negative, which is hardly surprising as we are considering compacts attractor and a well-behaved observable, whose values on the attractor have an upper bound. Note that for Axiom A systems, $d_s$ and $d_u$ are constant almost everywhere in the attractor $\Omega$, so that the information gathered for $x = x_0$ has a global value. Therefore, the expression for $\xi_A$ is universal, in the sense that we can gather fundamental properties of the dynamical system by looking at the shape parameters of the extremes of a generic observables with the properties described above. Measuring $\xi_A$ allows us to provide an upper and lower bound for $d_{KY}$ and *vice versa*. Note also that $\delta$ can be used to put upper and lower bounds to the Kaplan-Yorke dimension of the system (which is also independent of the assumption that the direction $\hat{n}$ is generic), as

$$d_{KY} = d_s + d_u + d_n > d_s + (d_u + d_n)/2 = \delta = -1/\xi_A$$

and

$$d_{KY} = d_s + d_u + d_n < 2d_s + d_u + d_n = 2\delta = -2/\xi_A,$$

so that

$$-1/\xi_A < d_{KY} < -2/\xi_A.$$

In particular, we have that $\xi_A$ is small and negative if and only if the Kaplan-Yorke dimension of the attractor is large. If we consider a chaotic system with a high dimensional attractor (*e.g.* in the case of an extensive chaotic system with many degrees of freedom), we derive that $\xi_A \approx 0$. This may well explain why in a multitude of applications in natural sciences the special $\xi = 0$ member of the GPD family often gives a good first guess of the statistics of observed extremes [72]. Alternatively, this result suggests that if we perform a statistical analysis of the extremes using the POT method for a high-dimensional chaotic system and obtain, as a result of the statistical inference, a shape parameter $\xi_A \ll 0$ or $\xi_A \geq 0$, we should conclude that our sample is not yet suited for an EVT statistical fit. This may depend on the fact that we have selected an insufficiently stringent value for $T$. Obviously, choosing higher values for $T$ implies that we need to have longer time series of the observable under investigations.

Furthermore, the expression we obtain for $\sigma$ leads to interpreting it as a scale parameter. We derive, as anticipated, that $\sigma > 0$ and we observe that it is proportional to the actual range of values considered as extremes of the observable $A(x)$, by incorporating the difference between the absolute maximum of the observable $A_{max}$ and the selected threshold $T$. Therefore, if we consider as observable $A_1(x) = \alpha A(x)$, with $\alpha > 0$, and take as threshold for $A_1(x)$ the value $\alpha T$, we have that $\xi_{A_1} = \xi_A$ and $\sigma_{A_1} = \alpha \sigma_A$. In physical terms, $\sigma$ changes if we change the unit of measure of the observable, whereas $\xi$ does not. More generally, we can make the following construction. Let's define $\min(A)|_\Omega = A_{min}$. If we select an observable $A_2(x) = \Phi(A(x))$, with $\max(\Phi)|_{[A_{min}, A_{max}]} = \Phi(A_{max})$, $\Phi$ differentiable and $d\Phi(y)/dy$ positive in a sufficiently wide neighbourhood around $y = A(x_{max})$ so to ensure monotonicity of $A_2(x)$ near $x = x_{max}$, we get $\xi_{A_2} = \xi_A$ and $\sigma_{A_2} = \gamma \sigma_A$, where $\gamma = d\Phi(y)/dy|_{y=A(x_{max})}$. If one fits 8.2.16 to data, $\sigma$ will give an estimate for the absolute maximal extreme, even when it has not been observed yet $A_{max} = T + \delta \sigma$.

## 8.2.5
## Partial Dimensions along the Stable and Unstable Directions of the Flow

It is worth considering the following strategy of investigation of the local properties of the invariant measure near $x = x_0$, where $A(x_0) = A_{max}$. By performing statistical inference of the extremes of $A$ we can deduce as a result of the data fitting the best estimate of $\xi_A = 1/\delta$. If, following [78], we select as observable, e.g. $B(x) = C - (\text{dist}(x, x_0))^{1/\alpha}$, $\alpha > 0$, we have that the extremes of the observable $B$ feature as shape parameter $\xi_B = -1/(\alpha D) = -1/(\alpha d_{KY})$ and scale parameter $\sigma_B = (C - \tau)/(\alpha D) = (C - \tau)/(\alpha d_{KY})$ [78], where $C$ is a constant and $\tau$ is the chosen threshold.

We can then derive:

$$\frac{2}{\xi_A} - \frac{2}{\alpha \xi_B} = d_u + d_n, \tag{8.2.19}$$

$$\frac{\beta}{\xi_B} - \frac{2}{\xi_A} = d_s, \tag{8.2.20}$$

where, as discussed above, we can take $d_n = 1$. Therefore, using rather general classes of observables, we are able to deduce the partial dimensions along the stable and unstable manifolds, just by looking at the index of extremes related to $x = x_0$. Note that, more generally, $d_u$ and $d_s$ can be deduced from the knowledge of any pair of values $(\xi_A, \xi_B)$, $(\sigma_A, \xi_B)$, $(\xi_A, \sigma_B)$, and $(\sigma_A, \sigma_B)$.

This provides further support to the idea that extremes can be used as excellent diagnostic indicators for the detailed dynamical properties of a system. The message we like to conclude is that one can construct observables whose large fluctua-tions give precise information on the dynamics. While considering various sorts of anisotropic scalings of the neighborhood of a point of the attractor allows to derive its partial dimensions [271], the specific result we obtain here is that using an arbitrary *physical* observable and studying its extremes, we automatically select a special, non ellipsoidal neighborhood, which splits automatically the stable and the unstable part of the dynamics and identifies its dimensions.

### 8.2.6
### Expressing the shape parameter in terms of the GPD moments and of the invariant measure of the system

We consider the physical observable $A$. We denote by

$$
\begin{aligned}
f_{GPD}(z; \xi_A, \sigma_A) &= \frac{\mathrm{d}}{\mathrm{d}z}\left(GPD_{\xi_A}\left(\frac{z}{\sigma_A}\right)\right) \\
&= \frac{1}{\sigma_A}\left(1 + \frac{\xi_A z}{\sigma_A}\right)^{-1/\xi_A - 1}.
\end{aligned}
\tag{8.2.21}
$$

We can express its moments as follows:

$$
\int_{-\sigma_A/\xi_A}^{0} dz\, z\, f_{GPD}(z; \xi_A, \sigma_A) = \frac{\sigma_A}{1 - \xi_A} = M_1
\tag{8.2.22}
$$

$$
\int_{-\sigma_A/\xi_A}^{0} dz\, z^2\, f_{GPD}(z; \xi_A, \sigma_A) = \frac{2\sigma_A^2}{(1 - \xi_A)(1 - 2\xi_A)} = M_2
\tag{8.2.23}
$$

$$
\cdots
$$

$$
\int_{-\sigma_A/\xi_A}^{0} dz\, z^n\, f_{GPD}(z; \xi_A, \sigma_A) = \frac{n!\sigma_A^n}{\Pi_{k=0}^{n}(1 - k\xi_A)} = M_n
\tag{8.2.24}
$$

where convergence is obtained for all moments because the shape parameter $\xi_A$ is negative. Using the expression of the first two moments of the distribution, it is easy to derive that

$$
\xi_A = \frac{1}{2}\left(1 - \frac{M_1^2}{M_2 - M_1^2}\right) = \frac{1}{2}\left(1 - \frac{1}{\mathrm{ind}_A}\right)
\tag{8.2.25}
$$

and

$$
\sigma_A = \frac{M_1 M_2}{2(M_2 - M_1^2)}
\tag{8.2.26}
$$

where we indicate explicitly that we refer to the observable $A$ and we have introduced the index of dispersion $\mathrm{ind}_A$, the ratio between the variance and the squared first moment of the considered stochastic variable.

One can express $\xi_A$ and $\sigma_A$ as a function of higher moments as well. In general, one obtains (considering $n \geq 2$):

$$
\xi_A = \frac{1}{n(n-1)}\left(n - 1 - \frac{M_{n-1}^2}{M_{n-2}M_n - M_{n-1}^2}\right)
\tag{8.2.27}
$$

and

$$\sigma_A = \frac{1}{n(n-1)} \frac{M_{n-1}M_n}{M_{n-2}M_n - M_{n-1}^2}. \tag{8.2.28}$$

We now connect the previous formulas to the properties of the invariant measure of the dynamical system. As we know, the GPD is the exact asymptotic model for the extremes of the observable $A$, so that we can express the results in terms of the conditional invariant measure as follows:

$$M_n^T = \frac{\langle \tilde{A}_n^T \rangle}{\langle \tilde{A}_0^T \rangle} \tag{8.2.29}$$

where $\tilde{A}_n^T(x) = \Theta(A(x) - T)(A(x) - T)^n$ and

$$\langle \tilde{A}_n^T \rangle = \int \mu(dx)\Theta(A(x) - T)(A(x) - T)^n, \tag{8.2.30}$$

with $\Theta$ being the usual Heaviside distribution. Using the definition of the first moments of the distributions, we obtain the following expression for the shape and dispersion parameters, respectively:

$$\xi_A^T = \frac{1}{2}\left(1 - \frac{(\langle \tilde{A}_1^T \rangle)^2}{\langle \tilde{A}_0^T \rangle\langle \tilde{A}_2^T \rangle - (\langle \tilde{A}_1^T \rangle)^2}\right), \tag{8.2.31}$$

and

$$\sigma_A^T = \frac{1}{2}\frac{\langle \tilde{A}_1^T \rangle\langle \tilde{A}_2^T \rangle}{\langle \tilde{A}_2^T \rangle\langle \tilde{A}_0^T \rangle - \langle \tilde{A}_1^T \rangle^2}, \tag{8.2.32}$$

where these results are exact in the limit for $T \to A_{max}$. The more general expressions obtained using Eqs. 8.2.27 and 8.2.28 read as follows:

$$\xi_A^T = \frac{1}{n(n-1)}\left(n - 1 - \frac{(\langle \tilde{A}_{n-1}^T \rangle)^2}{\langle \tilde{A}_{n-2}^T \rangle\langle \tilde{A}_n^T \rangle - (\langle \tilde{A}_{n-1}^T \rangle)^2}\right), \tag{8.2.33}$$

and

$$\sigma_A^T = \frac{1}{n(n-1)}\frac{\langle \tilde{A}_{n-1}^T \rangle\langle \tilde{A}_n^T \rangle}{\langle \tilde{A}_{n-2}^T \rangle\langle \tilde{A}_n^T \rangle - \langle \tilde{A}_{n-1}^T \rangle^2}, \tag{8.2.34}$$

so that, if we consider $n \geq 3$, the normalization factor $\langle \tilde{A}_0^T \rangle = \mu(A(x) \geq T)$ is not present in the previous formulas and, the higher the value of $n$, the smoother the functions defining the $\tilde{A}^T$'s observables.

As a check, it is useful to verify that the right hand side of Eq. 8.2.31 gives the same general result as given in Eq. 8.2.17. By definition we have for $T \to A_{\max}$:

$$\mu(\Omega_{A_{\max}}^T) = \langle \tilde{A}_0^T \rangle = \int \mu(dx)\Theta(A(x) - T)$$

$$\sim \tilde{h}_{x_0}(A_{\max} - T)(A_{\max} - T)^{\delta^\epsilon}. \tag{8.2.35}$$

Using the fundamental theorem of calculus, it is possible to derive that:

$$\langle \tilde{A}_n^T \rangle = \int_T^{A_{max}} dz \, n(z-T)^{n-1} \langle \tilde{A}_0^z \rangle_0. \tag{8.2.36}$$

we obtain:

$$\langle \tilde{A}_1^T \rangle \sim \frac{\tilde{h}_{x_0}(A_{max}-T)}{(\delta+1)}(A_{max}-T)^{\delta+1} \tag{8.2.37}$$

and

$$\langle \tilde{A}_2^T \rangle \sim \frac{2\tilde{h}_{x_0}(A_{max}-T)}{(\delta+1)(\delta+2)}(A_{max}-T)^{\delta+2}. \tag{8.2.38}$$

By plugging these expression into Eq. 8.2.31, we indeed obtain $\xi = -1/\delta$, which, as we expect, agrees with Eq. 8.2.17. We also note that it is possible to generalize the results for Eqs. 8.2.33 and 8.2.34 for arbitrary $n$.

Moreover, note that the representation of the parameters via moments as we presented in this subsection can be replicated step by step for the distance observables $B(x, x_0) = C - \text{dist}(x, x_0)^{1/alpha}$ discussed above. We obtain:

$$\xi_B^T = \frac{1}{n(n-1)} \left( n - 1 - \frac{(\langle \tilde{B}_{n-1}^T \rangle)^2}{\langle \tilde{B}_{n-2}^T \rangle \langle \tilde{B}_n^T \rangle - (\langle \tilde{B}_{n-1}^T \rangle)^2} \right), \tag{8.2.39}$$

and

$$\sigma_B^T = \frac{1}{n(n-1)} \frac{\langle \tilde{B}_{n-1}^T \rangle \langle \tilde{B}_n^T \rangle}{\langle \tilde{B}_{n-2}^T \rangle \langle \tilde{A}_n^T \rangle - \langle \tilde{B}_{n-1}^T \rangle^2}, \tag{8.2.40}$$

where the quantities $\langle \tilde{B}_k^T \rangle_0$, $k \geq 0$, and $\langle \tilde{B}_0^T \rangle_0 = 1$ are constructed analogously to how described in Eq. 8.2.29 .

We wish to remark that Eqs. 8.2.31-8.2.32 and Eqs. 8.2.39-8.2.40 combined could in fact provide a very viable method for estimating the GPD parameters from data, since moments estimators are in general more stable than maximal likelihood methods, and one can extract the value of both $d_u$ and $d_s$ using Eqs. 8.2.19-8.2.20.

## 8.3
## Impacts of Perturbations: Response Theory for Extremes

We wish to present some ideas on how to use response theory and the specific expressions given in Eqs 8.2.17-8.2.18 to derive a response theory for extremes of physical and distance observables in Axiom A dynamical systems. Let's assume that we alter the Axiom A dynamical system under consideration as $\dot{x} = G(x) \rightarrow \dot{x} = G(x) + \epsilon X(x)$, where $\epsilon$ is a small parameter and $X(x)$ is a smooth vector field, so that the evolution operator, that is the flow, is altered as $f^t \rightarrow f_\epsilon^t$ and the invariant measure is altered as $\mu \rightarrow \mu_\epsilon$. Ruelle's response theory allows to express the change in the expectation value of a general measurable observable $\Psi(x)$ as a perturbative series as

$$\langle \Psi \rangle^\epsilon = \langle \Psi \rangle_0 + \sum_{j=1}^{\infty} \epsilon^j \langle \Psi \rangle_0^{(j)},$$

with $j$ indicating the order of perturbative expansion, where

$$\langle \Psi \rangle^\epsilon = \int \mu_\epsilon(dx) \Psi(x)$$

is the expectation value of $\Psi$ over the perturbed invariant measure, that is the invariant measure with respect to $G + \epsilon X$, and

$$\langle \Psi \rangle_0 = \int \mu(dx) \Psi(x) \tag{8.3.1}$$

is the unperturbed expectation value of $\Psi$. The term corresponding to the perturbative order of expansion $j$ is given by $\langle \Psi \rangle_0^{(j)}$, which can be expressed in terms of the time-integral of a suitably defined Green function [83]:

$$\langle \Psi \rangle_0^{(j)} = \int_{-\infty}^{\infty} d\tau_1 \ldots \int_{-\infty}^{\infty} d\tau_n G_\Psi^{(n)}(\tau_1, \ldots, \tau_n). \tag{8.3.2}$$

The integration kernel $G_\Psi^{(n)}(\tau_1, \ldots, \tau_n)$ is the $n^{th}$ order Green function, which can be written as:

$$G_\Psi^{(n)}(\tau_1, \ldots, \tau_n) = \langle \Theta(\tau_1) \ldots \Theta(\tau_n - \tau_{n-1}) \Lambda \Pi(\tau_n - \tau_{n-1}) \ldots \Lambda \Pi(\tau_1) \Psi(x) \rangle_0, \tag{8.3.3}$$

where $\Lambda(\bullet) = X(\cdot) \cdot \nabla(\bullet)$ describes the impact of the perturbation field and $\Pi(\sigma)$ is the unperturbed time evolution operator such that $\Pi(\sigma)F(x) = F(x(\sigma))$. The Green function obeys two fundamental properties

- its variables are time-ordered: if $j > k$, $\tau_j < \tau_k \rightarrow G_\Psi^{(n)}(\tau_1, \ldots, \tau_n) = 0$;
- the function is causal: $\tau_1 < 0 \rightarrow G_\Psi^{(n)}(\tau_1, \ldots, \tau_n) = 0$.

It is important to stress that many authors suggest that, for all practical purposes, the validity of the response theory extends well beyond the (well understood but somewhat limited) mathematical world of Axiom A systems if one consider reasonable physical systems and reasonable observables. See discussions in, *e.g.*, [263] and references therein. This is closely related to the chaos hypothesis [266].

Limiting our attention to the linear case we have:

$$\langle \Psi \rangle_0^{(1)} = \int_{-\infty}^{+\infty} d\tau_1 G_\Psi^{(1)}(\tau_1), \tag{8.3.4}$$

where the first order Green function can be expressed as follows:

$$G_\Psi^{(1)}(\tau_1) = \langle \Theta(\tau_1) X(x) \cdot \nabla \Psi(x(\tau_1)) \rangle_0. \tag{8.3.5}$$

Recall that in general:

$$\frac{d^n \langle \Psi \rangle^\epsilon}{d\epsilon^n} \bigg|_{\epsilon=0} = n! \langle \Psi \rangle_0^{(n)},$$

which we use, in particular, for the $n = 1$ case.

### 8.3.1
## Sensitivity of the Shape Parameter as Determined by the Changes in the Moments

We wish to propose a linear response formula for the parameter $\xi_A$ using Eq. 8.2.31. We start by considering that in Eq. 8.2.31 the shape parameter is expressed for every $T < A_{max}$ as a function of actual observables of the system. Unfortunately, in order to apply Ruelle's response theory, we need the observables to be smooth, which is in contrast with the presence of the $\Theta$ in the definition of the terms $\langle \tilde{A}_j^T \rangle^\epsilon$. Nonetheless, replacing the $\Theta$'s with a smooth approximation $\Theta_S$, the Ruelle response theory can be rigorously applied. Considering a sequence of approximating $\Theta_S^m$ such that the measure of the support of $\Theta - \Theta_S^m$ is smaller than $\delta_m = (A_{max} - T)/m$, it is reasonable to expect that, as $\delta_m \to 0$, the effect of the smoothing becomes negligible, because a smaller and smaller portion of the extremes is affected, and the response of the smoothed observable approaches that of $\langle \tilde{A}_j^T \rangle^\epsilon$. Therefore, we derive formally retaining the $\Theta$ in the definition of the $\langle \tilde{A}_j^T \rangle^\epsilon$ for every $T < A_{max}$:

$$\frac{d\xi_A^{T,\epsilon}}{d\epsilon}\bigg|_{\epsilon=0} = -\frac{1}{n(n-1)} \frac{d}{d\epsilon} \left\{ \frac{(\langle \tilde{A}_{n-1}^T \rangle^\epsilon)^2}{\langle \tilde{A}_{n-2}^T \rangle^\epsilon \langle \tilde{A}_n^T \rangle^\epsilon - (\langle \tilde{A}_{n-1}^T \rangle^\epsilon)^2} \right\}\bigg|_{\epsilon=0} \quad (8.3.6)$$

and

$$\frac{d\sigma_A^{T,\epsilon}}{d\epsilon}\bigg|_{\epsilon=0} = \frac{1}{n(n-1)} \frac{d}{d\epsilon} \left\{ \frac{\langle \tilde{A}_{n-1}^T \rangle^\epsilon \langle \tilde{A}_n^T \rangle^\epsilon}{\langle \tilde{A}_{n-2}^T \rangle^\epsilon \langle \tilde{A}_n^T \rangle^\epsilon - (\langle \tilde{A}_{n-1}^T \rangle^\epsilon)^2} \right\}\bigg|_{\epsilon=0}, \quad (8.3.7)$$

where $n \geq 2$. By expanding the derivative in Eq. 8.3.6, the previous expression can be decomposed in various contributions containing exclusively $\langle \tilde{A}_{n-2}^T \rangle_0$, $\langle \tilde{A}_{n-1}^T \rangle_0$, $\langle \tilde{A}_n^T \rangle_0$ and $\langle \tilde{A}_{n-2}^T \rangle_0^{(1)}$, $\langle \tilde{A}_{n-1}^T \rangle_0^{(1)}$, $\langle \tilde{A}_n^T \rangle_0^{(1)}$, while higher order terms are not included.

We obtain:

$$\frac{d\xi_A^{T,\epsilon}}{d\epsilon}\bigg|_{\epsilon=0} = \alpha_{\xi,T,n-2,n}^{(1)} \frac{d\langle \tilde{A}_{n-2}^T \rangle}{d\epsilon}\bigg|_{\epsilon=0} + \alpha_{\xi,T,n-1,n}^{(1)} \frac{d\langle \tilde{A}_{n-1}^T \rangle}{d\epsilon}\bigg|_{\epsilon=0} + \alpha_{\xi,T,n,n}^{(1)} \frac{d\langle \tilde{A}_n^T \rangle}{d\epsilon}\bigg|_{\epsilon=0},$$
$$(8.3.8)$$

$$= \frac{1-(n-1)\xi_A^T}{n-1} \left( \frac{1}{\langle \tilde{A}_{n-2}^T \rangle} \frac{d\langle \tilde{A}_{n-2}^T \rangle}{d\epsilon}\bigg|_{\epsilon=0} - \frac{2}{\langle \tilde{A}_{n-1}^T \rangle} \frac{d\langle \tilde{A}_{n-1}^T \rangle}{d\epsilon}\bigg|_{\epsilon=0} + \frac{1}{\langle \tilde{A}_n^T \rangle} \frac{d\langle \tilde{A}_n^T \rangle}{d\epsilon}\bigg|_{\epsilon=0} \right),$$
$$(8.3.9)$$

where

$$\alpha_{\xi,T,n-2,n}^{(1)} = \frac{1}{n(n-1)} \frac{\left(\langle \tilde{A}_{n-1}^T \rangle_0\right)^2 \langle \tilde{A}_n^T \rangle_0}{\left(\langle \tilde{A}_n^T \rangle_0 \langle \tilde{A}_{n-2}^T \rangle_0 - (\langle \tilde{A}_{n-1}^T \rangle_0)^2\right)^2},$$

$$\alpha_{\xi,T,n-1,n}^{(1)} = -\frac{2}{n(n-1)} \frac{\langle \tilde{A}_{n-2}^T \rangle_0 \langle \tilde{A}_{n-1}^T \rangle_0 \langle \tilde{A}_n^T \rangle_0}{\left(\langle \tilde{A}_n^T \rangle_0 \langle \tilde{A}_{n-2}^T \rangle_0 - (\langle \tilde{A}_{n-1}^T \rangle_0)^2\right)^2},$$

$$\alpha_{\xi,T,n,n}^{(1)} = \frac{1}{n(n-1)} \frac{\left(\langle \tilde{A}_{n-1}^T \rangle_0\right)^2 \langle \tilde{A}_{n-2}^T \rangle_0}{\left(\langle \tilde{A}_n^T \rangle_0 \langle \tilde{A}_{n-2}^T \rangle_0 - (\langle \tilde{A}_{n-1}^T \rangle_0)^2\right)^2}, \qquad (8.3.10)$$

and

$$\frac{\mathrm{d}\sigma_A^{T,\epsilon}}{\mathrm{d}\epsilon}\bigg|_{\epsilon=0} = \alpha_{\sigma,T,n-2,n}^{(1)} \frac{\mathrm{d}\langle \tilde{A}_{n-2}^T \rangle}{\mathrm{d}\epsilon}\bigg|_{\epsilon=0} + \alpha_{\sigma,T,n-1,n}^{(1)} \frac{\mathrm{d}\langle \tilde{A}_{n-1}^T \rangle}{\mathrm{d}\epsilon}\bigg|_{\epsilon=0} + \alpha_{\sigma,T,n,n}^{(1)} \frac{\mathrm{d}\langle \tilde{A}_n^T \rangle}{\mathrm{d}\epsilon}\bigg|_{\epsilon=0},$$
$$(8.3.11)$$

$$= \frac{\sigma_A^T}{(n-1)(1-n\xi_A^T)} \left( \frac{n - n(n-1)\xi_A^T}{\langle \tilde{A}_{n-2}^T \rangle} \frac{\mathrm{d}\langle \tilde{A}_{n-2}^T \rangle}{\mathrm{d}\epsilon}\bigg|_{\epsilon=0}, \qquad (8.3.12) \right.$$

$$\left. -\frac{3n-1-3n(n-1)\xi_A^T}{\langle \tilde{A}_{n-1}^T \rangle} \frac{\mathrm{d}\langle \tilde{A}_{n-1}^T \rangle}{\mathrm{d}\epsilon}\bigg|_{\epsilon=0} + \frac{2n-1-2(n-1)n\xi_A^T}{\langle \tilde{A}_n^T \rangle} \frac{\mathrm{d}\langle \tilde{A}_n^T \rangle}{\mathrm{d}\epsilon}\bigg|_{\epsilon=0} \right),$$

where

$$\alpha_{\sigma,T,n-2,n}^{(1)} = -\frac{1}{n(n-1)} \frac{\langle \tilde{A}_{n-1}^T \rangle_0 \left(\langle \tilde{A}_n^T \rangle_0\right)^2}{\left(\langle \tilde{A}_n^T \rangle_0 \langle \tilde{A}_{n-2}^T \rangle_0 - (\langle \tilde{A}_{n-1}^T \rangle_0)^2\right)^2},$$

$$\alpha_{\sigma,T,n-1,n}^{(1)} = \frac{1}{n(n-1)} \frac{\langle \tilde{A}_n^T \rangle_0 \left(\langle \tilde{A}_n^T \rangle_0 \langle \tilde{A}_{n-2}^T \rangle_0 + \left(\langle \tilde{A}_{n-1}^T \rangle_0\right)^2\right)}{\left(\langle \tilde{A}_n^T \rangle_0 \langle \tilde{A}_{n-2}^T \rangle_0 - (\langle \tilde{A}_{n-1}^T \rangle_0)^2\right)^2},$$

$$\alpha_{\sigma,T,n,n}^{(1)} = -\frac{1}{n(n-1)} \frac{\left(\langle \tilde{A}_{n-1}^T \rangle_0\right)^3}{\left(\langle \tilde{A}_n^T \rangle_0 \langle \tilde{A}_{n-2}^T \rangle_0 - (\langle \tilde{A}_{n-1}^T \rangle_0)^2\right)^2}. \qquad (8.3.13)$$

We wish to remark the special relevance of the observable $\langle \tilde{A}_0^T \rangle^\epsilon$, which is the normalizing factor in Eq. 8.2.29, and, in practice, measures the fraction of above-$T$-threshold events. Therefore, once $T$ is chosen, the sensitivity of $\langle \tilde{A}_0^T \rangle^\epsilon$ with respect to $\epsilon$ informs on whether the $\epsilon$-perturbation to the vector flow leads to an increase or decrease in the number of extremes. We obtain formally:

$$
\begin{aligned}
\frac{d\langle A_0^T\rangle^\epsilon}{d\epsilon}\bigg|_{\epsilon=0} &= \langle A_0^T\rangle_0^{(1)} = \int d\tau\, G_{A_0^T}^{(1)}(\tau) = \int d\tau\, \Theta(\tau)\,\langle X_k(x)\partial_k\Theta\left(A(f^\tau(x)-T)\right)\rangle_0 \\
&= \int d\tau\, \Theta(\tau)\,\langle X_k(x)\partial_k A(f^\tau(x))\delta\left(A(f^\tau(x)-T)\right)\rangle_0 \\
&= \int d\tau\, \Theta(\tau)\langle X_k(x)\partial_k f_i^\tau(x)\partial_y A(y)\bigg|_{y=f^\tau(x)}\delta\left(A(f^\tau(x))-T\right)\rangle_0,
\end{aligned}
$$
$$(8.3.14)$$

where $\delta$ is the derivative of the $\Theta$ function, with all the caveats discussed above, and $f^\tau(x)$ is the position of the dynamics $x(\tau)$ at time $\tau$ if $x(0) = x$. The formula can be interpreted as follows. In the last formula, $\partial_k f_i^\tau$ is the adjoint matrix of the tangent linear of the unperturbed flow. At each instant $\tau$ we consider, in the unperturbed system, all the trajectories starting in the infinite past from points distributed according to the invariant measure such that the observable $A$ at time zero has value equal to $T$. For each of these trajectories, we can measure whether the presence of the perturbation field $X(x)$ at time $-\tau$ would lead to a decrease or increase in $A$ at time zero. Summing over all trajectories, we get whether there is a net positive or negative change in the above threshold events at time zero. We integrate over all times $\tau$ at which the perturbation can impact and get the final result. Considering the geometrical construction given in Fig. 8.2, the previous formula can also be approximated as follows, because $f^\tau(x)$ has to be near to the maxima $x_0$:

$$
\frac{d\langle A_0^T\rangle^\epsilon}{d\epsilon}\bigg|_{\epsilon=0} \approx \int d\tau\, \Theta(\tau)\langle X_k(x)\partial_k f_i^\tau(x)\partial_i A(x_0)\delta\left(A(f^\tau(x)-T)\right)\rangle_0 \quad (8.3.15)
$$

Note that, when considering $\langle A_n^T\rangle^\epsilon$ for $n \geq 1$, we obtain:

$$
\begin{aligned}
\frac{d\langle A_n^T\rangle^\epsilon}{d\epsilon}\bigg|_{\epsilon=0} &= \langle A_n^T\rangle_0^{(1)} = \int d\tau\, G_{A_n^T}^{(1)}(\tau) \\
&= \int d\tau\, \Theta(\tau)\,\langle X_k(x)\partial_k\left[A(x(\tau))-T\right]^n\Theta\left(A(x(\tau))-T\right)]\rangle_0 \\
&= n\int d\tau\, \Theta(\tau)\langle X_k(x)\partial_k A(x(\tau))\left(A(x(\tau))-T\right)^{n-1}\Theta\left(A(x(\tau)-T)\right)\rangle_0,
\end{aligned}
$$
$$(8.3.16)$$

using that the derivative of $x^n\Theta(x)$ is $nx^{n-1}\Theta(x)$ for $n \geq 1$. Therefore, Eqs. 8.3.6-8.3.16 provide recipes for computing the sensitivity of $\xi_A^{T,\epsilon}$ and $\sigma_A^{T,\epsilon}$ at $\epsilon = 0$ for any case of practical interest. In fact, $A_{max} - T$ is indeed finite, because in order to collect experimental data or process the output of numerical simulations we need to select a threshold which is high enough for discriminating true extremes and low enough for allowing a sufficient number of samples to be collected for robust data processing. Note that all statistical procedures used in estimating GPD parameters from data are actually based on finding a reasonable value for $T$ such that both conditions described above apply by testing in an appropriate sense that parameters' estimates do not vary appreciably when changing $T$ [40, 3].

We wish to underline that apparently potential problems emerge when taking the limit in Eqs. 8.3.6-8.3.7 for higher and higher values of $T$. It is indeed not clear at this stage whether

$$\lim_{T \to A_{max}} \frac{d\xi_A^{T,\epsilon}}{d\epsilon}\bigg|_{\epsilon=0} = \lim_{T \to A_{max}} \lim_{\epsilon \to 0} \frac{\xi_A^{T,\epsilon} - \xi_A^{T,0}}{\epsilon} \tag{8.3.17}$$

exists, because we cannot apply the smoothing argument presented above in the limit of vanishing $A_{max} - T$. More important, it is not clear whether such limit is equal to

$$\lim_{\epsilon \to 0} \lim_{T \to A_{max}} \frac{\xi_A^{T,\epsilon} - \xi_A^{T,0}}{\epsilon}, \tag{8.3.18}$$

which seems at least as well suited for describing the change of the shape observable given in Eq. 8.2.31 due to an $\epsilon$-perturbation in the dynamics and gives the link to the universality of the asymptotic distribution. Obviously, if the two limits given in Eqs. 8.3.17 and 8.3.18 exist and are equal, then a rigorous response theory for $\xi_A$ can be established. Same applies when considering the properties of $\sigma_A^{T,\epsilon}$.

The same derivation and discussion can be repeated for the $B$ observables introduced above and we can derive the corresponding formulas for $d\xi_B^{T,\epsilon}/d\epsilon|_{\epsilon=0}$ and $d\sigma_B^{T,\epsilon}/d\epsilon|_{\epsilon=0}$, where the relevant limit for $T$ is $T \to C$.

Let us try to give a more intuitive interpretation to the results given above. Consider Eq. 8.2.25 and assume that, indeed, $\xi_A^T$ is differentiable with respect to $\epsilon$. We have:

$$\frac{d\xi_A^\epsilon}{d\epsilon}\bigg|_{\epsilon=0} = -\frac{1}{2}\frac{d}{d\epsilon}\left\{\frac{1}{\text{ind}_A^\epsilon}\right\}\bigg|_{\epsilon=0} = \frac{1}{2\text{ind}_A^2}\frac{d}{d\epsilon}\left\{\text{ind}_A^\epsilon\right\}\bigg|_{\epsilon=0}, \tag{8.3.19}$$

which implies that the sensitivity of the shape parameter is half of the opposite of the sensitivity of the inverse of the index of dispersion $\text{ind}_A$. Alternatively

$$\frac{1}{1 - 2\xi_A}\frac{d(1 - 2\xi_A^\epsilon)}{d\epsilon}\bigg|_{\epsilon=0} = -\frac{1}{\text{ind}_A}\frac{d\,\text{ind}_A^\epsilon}{d\epsilon}\bigg|_{\epsilon=0}, \tag{8.3.20}$$

which implies that relative sensivity of $1 - 2\xi_A^T$ and $\text{ind}_A$ coincides. Therefore, a positive sensitivity of the index of dispersion (larger relative variability of the extremes of the observable $A$ with positive values of $\epsilon$) implies a larger value (closer to 0) of $\xi_A$, and so the possibility that larger and larger extremes are realized. The analogous interpretation applies for the $B$ observables.

### 8.3.2
### Sensitivity of the shape parameter as determined by the modification of the geometry

In the previous subsection we have discerned that the Ruelle response theory supports the idea that the shape parameters descriptive of the extremes of both the physical observables $A$ and the distance observables $B$ change with a certain degree of regularity when considering $\epsilon$−perturbations to the dynamics.

In this subsection, we wish consider the sensitivity of extremes with respect to perturbation from another angle, *i.e.* investigating the relationship between the shape parameters $\xi_A$ and $\xi_B$ and the partial dimension of the attractor along the stable, neutral and unstable manifolds of the underlying dynamical system, see Eqs. 8.2.19-8.2.20. As long as the $\epsilon$-perturbation is small, the modified dynamical system belongs to the Axiom A family, indeed we have structural stability, so that the results presented above apply. Therefore, we can write in more general terms:

$$\xi_A^\epsilon = -1/\delta^\epsilon = -1/(d_s^\epsilon + d_u^\epsilon/2 + d_n^\epsilon/2) \tag{8.3.21}$$

$$\xi_B^\epsilon = -1/(\alpha d_{KY}^\epsilon) = -1/(\alpha(d_s^\epsilon + d_u^\epsilon + d_n^\epsilon)). \tag{8.3.22}$$

In the following, we introduce somewhat heuristically derivatives with respect to $\epsilon$ of quantities for which we do not know *a priori* that they are differentiable. The main point we want to make is that if $\xi_A$ and $\xi_B$ are differentiable with respect to $\epsilon$, then various quantities describing the structure of the attractor have to be differentiable. Therefore, the existence of the limits given in Eqs. 8.3.17 and 8.3.18 (and their equivalent for the $B$ observables) would have far-reaching consequences. We will discuss the obtained results at the end of the calculations. Another *caveat* we should mention is that Eqs. 8.3.21-8.3.22 are in general true almost anywhere, so that we may have to interpret the derivatives in this section in some suitable weak form.

It seems relevant to add the additional hypothesis of strong transversality for the unperturbed flow, which is equivalent to invoking structural stability [71]. We take such pragmatic point of view and proceed assuming that derivatives with respect to $\epsilon$ are well defined. Linearizing the dependence of $\xi_A$ on $\epsilon$ around $\epsilon = 0$ in Eq. 8.3.21, we obtain:

$$\frac{d\xi_A^\epsilon}{d\epsilon}\bigg|_{\epsilon=0} = \left\{\frac{d\xi_A^\epsilon}{d(d_s^\epsilon)}\frac{d(d_s^\epsilon)}{d\epsilon}\right\}\bigg|_{\epsilon=0} + \left\{\frac{d\xi_A^\epsilon}{d(d_u^\epsilon)}\frac{d(d_u^\epsilon)}{d\epsilon}\right\}\bigg|_{\epsilon=0} + \left\{\frac{d\xi_A^\epsilon}{d(d_n^\epsilon)}\frac{d(d_n^\epsilon)}{d\epsilon}\right\}\bigg|_{\epsilon=0}. \tag{8.3.23}$$

We have that $d(d_u^\epsilon)/d\epsilon|_{\epsilon=0} = d(d_n^\epsilon)/d\epsilon|_{\epsilon=0} = 0$, as, thanks to structural stability, small perturbations do not alter the qualitative properties of the dynamics, and cannot change in a step-wise way the integer dimension of the expanding or neutral directions. Hence we have for the local dimension $D$ that $d(d_s^\epsilon)/d\epsilon|_{\epsilon=0} = d(D^\epsilon)/d\epsilon|_{\epsilon=0}$ which implies that

$$\frac{d\xi_A^\epsilon}{d\epsilon}\bigg|_{\epsilon=0} = \left\{\frac{1}{(d_s^\epsilon + d_u^\epsilon/2 + d_n^\epsilon/2)^2}\frac{d(D^\epsilon)}{d\epsilon}\right\}\bigg|_{\epsilon=0} = \left\{\frac{1}{(\xi_A^\epsilon)^2}\frac{d(D^\epsilon)}{d\epsilon}\right\}\bigg|_{\epsilon=0}. \tag{8.3.24}$$

This implies that the shape parameter $\xi$ increases, thus attaining a value closer to zero ($\xi_A$ is always negative), when the perturbation increases the Kaplan-Yorke dimension of the attractor. This corresponds to the case when the perturbation favours *forcing* over *dissipation*. This matches quite well, at least qualitatively, with the discussion following Eq. 8.3.19.

We have that, when considering a distance observable of the form $B(x) = -\text{dist}(x, x_0)^{1/\alpha}$, along the same steps described above one gets the following result:

$$\left.\frac{\mathrm{d}\xi_B^\epsilon}{\mathrm{d}\epsilon}\right|_{\epsilon=0} = \left\{\frac{1}{\alpha(D^\epsilon)^2}\frac{\mathrm{d}(D^\epsilon)}{\mathrm{d}\epsilon}\right\}\bigg|_{\epsilon=0} ; \qquad (8.3.25)$$

such result can be easily generalized by considering the class of observables described in [78].

Combining Eq. 8.3.6 with Eq. 8.3.24, and the derivative with respect to $\epsilon$ of Eq. 8.2.39 with Eq. 8.3.25, we can derive two expressions for the sensitivity of the unstable dimension and of the Kaplan Yorke dimension at $\epsilon = 0$:

$$\left.\frac{\mathrm{d}(d_{KY}^\epsilon)}{\mathrm{d}\epsilon}\right|_{\epsilon=0} = -\left\{\frac{(d_s^\epsilon + d_u^\epsilon/2 + d_n^\epsilon/2)^2}{2}\right\}\bigg|_{\epsilon=0} \left\{\frac{\mathrm{d}}{\mathrm{d}\epsilon}\frac{((\langle\tilde{A}_1^T\rangle^\epsilon)^2}{\langle\tilde{A}_2^T\rangle^\epsilon\langle\tilde{A}_0^T\rangle^\epsilon - ((\langle\tilde{A}_1^T\rangle^\epsilon)^2}\right\}\bigg|_{\epsilon=0}$$
$$(8.3.26)$$

$$= -\left\{\frac{\alpha d^\epsilon_{KY}{}^2}{2}\right\}\bigg|_{\epsilon=0} \left\{\frac{\mathrm{d}}{\mathrm{d}\epsilon}\frac{((\langle\tilde{B}_1^T\rangle^\epsilon)^2}{\langle\tilde{B}_2^T\rangle^\epsilon\langle\tilde{B}_0^T\rangle^\epsilon - ((\langle\tilde{B}_1^T\rangle^\epsilon)^2}\right\}\bigg|_{\epsilon=0},$$
$$(8.3.27)$$

where, rigorously, we have to work in the limit $T \to A_{max}$ in Eq. 8.3.26 and $T \to 0$ in Eq. 8.3.27, and we have assumed that the Hausdorff dimension coincides with the Kaplan Yorke dimension.

The previous results imply that if one of $\xi_A$, $\xi_B$ or the Kaplan-Yorke dimension of the underlying Axiom A system change smoothly with $\epsilon-$perturbations to the dynamics, so do the other two quantities. The analogues formulas using Eqs. 8.2.33 and 8.2.39 are differentiable for $n \geq 3$ for fixed $T$. In order to make these considerations more rigorous one should identify in which sense one may interchange the limits in $T$ and $\epsilon$ in Eqs. 8.3.17 and 8.3.18. This may suggest ways to study the regularity of the Kaplan-Yorke dimension by resorting to the analysis of the regularity of a much more tractable expressions involving moments of given observables only.

This result provides useful insight also not considering the problematic limits discussed above. Taking a more qualitative point of view, this suggests that, when considering small perturbation in the dynamics of chaotic systems behaving like Axiom A systems, there is a link between the presence (or lack) of differentiability with respect to $\epsilon$ of $\xi_A$, $\xi_B$ and $d_{KY}$, so that either all of them or none of them is differentiable with respect to $\epsilon$. Specifically, we obtain that if the perturbation tends to increase the dimensionality of the attractor (thus, in physical terms, favoring forcing over dissipation), the value of $\xi$ becomes closer to zero, so that the occurrence of very large extreme events becomes more likely. The system, in this case, has more freedom to perform large fluctuations.

Taking a more pragmatic point of view, these results provide at least a rationale for the well-known fact that in moderate to high-dimensional strongly chaotic systems the Kaplan-Yorke dimension (and, actually, all the Lyapunov exponents) change smoothly with the intensity of the perturbating vector field [272] using simplified yet relevant fluid dynamical model. A detailed investigation of the apparent regularity for all practical purposes of the Lyapunov exponents with respect to small perturbations in the dynamics of intermediate complexity to high-dimensional models has been presented in [273]. We also wish to remark that if these regularity hypotheses were not satisfied, the very widespread (and practically successful) procedure of parametric tuning of high-dimensional models of natural, engineered or social phenomena would be absolutely hopeless, and delicate numerical procedures such as those involved in data assimilation of high-dimensional dynamical systems would lack any sort of robustness, contrary to the accumulated experience.

## 8.4
## Remarks on the Geometry and the Symmetries of the Problem

A specification is needed in the case of physical observables. We have here considered the case where the observable $A$ has a unique maximum restricted to $\Omega$ in $x = x_0 \in \Omega$. If $\Omega$ and $A$ share some symmetries, $x_0$ is not unique, and instead there is a set of points $x_0$'s belonging to $\Omega$, finite or infinite, depending of the kind of symmetries involved, where $A$ reaches its maximum value restricted to $\Omega$. Let's consider the relevant case where $A$ and $\Omega$ share a discrete symmetry, so that $\chi_0$, the set of the maximal point $x_0$'s, has finite cardinality. The results discussed here for the extremes of $A$ will be the direct analogues, and all arguments are valid in the appropriate sense, because we can perform an equivalent geometrical construction as in Fig. 8.2 for each element of $\chi_0$. When we consider an $\epsilon$-perturbation to the dynamics which respects the discrete symmetry, it is clear that all the results of the response presented here apply. Finally, one can deduce that if the considered perturbation, instead, breaks the discrete symmetry, the results presented here will still be valid as the break of the degeneracy will make sure that only one of the $x_0$'s (or a subset of $\chi_0$, if the corresponding perturbed vector flow obeys to a subgroup of the original symmetry group) still accounts for the extreme events of $A$.

We need to remark that such results have been derived using some intuitive geometrical construction and assuming generic relations between the direction of the gradient of $A$ at $x = x_0$ and the stable directions. It is possible to devise a special pair of Axiom A systems and observables such that the strange attractors do not fulfill such generic conditions. One can easily construct a situation where the gradient of $A$ is orthogonal also to stable manifold by immersing the attractor in a higher dimensional space and taking observables defined on such a space. In this case, the factor $1/2$ appearing in Eq. 8.2.14 will affect also some of the stable dimensions. Nonetheless, the result that for high-dimensional systems the distribution of extremes is indistinguishable from the Gumbel as the shape parameters tends to zero from below is of general validity. We believe that for a typical combinations of Axiom A systems and an observable given by a function allow for the generic conditions to be obeyed. We still need to understand how to frame precisely such a concept of genericity, which obviously differs from the traditional one, which focuses either on the observables, or on the systems. This should be the subject of theoretical investigation and accurate numerical testing.

A final remark: let's assume, instead, that the gradient of $A$ is vanishing in $x_0$ and that, at leading order near $x_0$, $A(x) \sim A_{max} + (x - x_0)^T H (x - x_0)$, where $H$ is a negative definite symmetric matrix the scalar product. The geometrical considerations will in this case imply, that apart from a linear change in the coordinates and rescaling, the statistical properties of the extremes of $A(x)$ will match those of $B(x) = C - dist(x, x_0)^2$.

# 9

# Extremes as Dynamical and Geometrical Indicators

In the previous chapters it has been shown that two different approaches have been devised for the study of extreme events, the BM method, which leads to using GEV (see Eq. 3.1.4) as statistical models, and the POT method, where instead the statistical model of reference is the GPD (see Eq. 3.1.8). As detailed in Chapter 3, if suitable conditions of weak short and long time scale correlations are met, the two approaches are equivalent *in the asymptotic limit*, so that, even if the selection procedure of the extremes is different, the information we derive is equivalent. For finite datasets, it is crucial to investigate which approach provides more reliable results and if differences between the use of the two methods arise when studying the extremes of the $g_i$'s observables.

We recall that analysing extremes poses the challenge of requiring extremely long time series, because the first step in the procedure of statistical inference relies of discarding the vast majority of the entries, and a sufficient number of *true* potential extremes have to be left for the fitting procedures to converge and have low uncertainty. In this chapter we wish to explore systematically how to get convergence of the empirical d.f. of BM and POT to the corresponding EVLs. The theory presented in the previous chapters is correct under the assumptions of an infinite sample of maxima (minima) each extracted among an infinite number of observations posing a question about the usefulness of the theorems when one has to deal with finite samples. In order to address this problem, we will construct the experiments as a sort of *numerical proofs* of the theorems introduced in the previous chapters, *i.e.*, define a procedure able to reproduce step-by-step the theoretical description and the convergence issues. One-dimensional and two-dimensional maps, are analyzed throughout the chapter, and results are provided both for deterministic and stochastic systems. The reader finds below also an extended review of the methods used to perform statistical inference of EVLs.

As an overall goal, we want to show here how extremes can be used for deriving crucial properties on the geometry of the system of mixing systems, allowing to derive precise information on the local and global properties of the attractor, and for studying qualitatively and quantitatively their dynamical properties, including the possibility of separating regular from mixing chaotic systems.

## 9.1
## The Block Maxima Approach

In Chapter 3 we have introduced the intimate relationship between the BM procedure of selection of extremes of time series and GEV statistical model, specifying the relevance of the decorrelation properties $D$ and $D'$.

In Chapter 4, we have then adapted this point of view for the study of observables of dynamical systems, introducing a new set of conditions $\amalg_q$ and $\amalg'_q$, which can be more easily related to standard properties of mixing the underlying dynamics. Finally, in Chapter 6, we have investigated some dynamical systems and derived explicit results for EVLs of specific observables.

Most of our theoretical results have been obtained considering the observables $g_i(dist(f^t x, \zeta))$, $i = 1, 2, 3$, described in Sect. 4.2.1. These are functions of the distance $r$ between a point in the orbit and a chosen point $\zeta$ belong to the attractor of the system. The basic requirements imposed on the $g_i(r)$ for achieving convergence towards the classical EVLs are the existence of a global maximum at $r = 0$ and a monotonic decrease with the argument. These properties imply that, by looking at maxima of $g_i(dist(x, \zeta))$, we are indeed looking at minima of $dist(x, \zeta)$, regardless the chosen metrics. This idea is schematically depicted in Fig. 9.1.

In the examples treated in this chapter, the BM approach boils down to to dividing the time series of the observable $g_i$, $i = 1, 2, 3$ of length $NK$ into bins of equal length $K$ and to selecting the maximum value in each of them, so that $N$ potential extremes are selected [3]. The left plot depicts a hypothetical series of $g_i(dist(f^t x, \zeta)$ at different values of time $t$, and indicates three different bin lengths $m_A, m_B, m_C$. The right picture is a 2-D representation of the balls centred on $\zeta$ inside which we are more likely to sample extremes, with smalller balls corresponding to longer bins. By choosing the BM approach we do not have a clear control of an effective radius of the sampled balls, because the largest value selected in a bin can easily be smaller than, *e.g.* the second largest of another bin. This situation is extremely common if significant clusters of extremes is present.

We remark that the BM approach is widely used in climatological and financial applications since it represents a very natural way of looking at extremes sampled at fixed time intervals. As an example, when dealing with a long time series of meteorological observations containing daily values - *e.g.* temperature records - one year is a natural and common choice for the bin length. The BM method helps us, in this case, in constructing a statistical model for the yearly maxima of temperature. Of course, one should consider whether such a choice for the bin length is sensible, given the total length of the time series and its rate of decay of correlations.

In general, the *a priori* knowledge of the asymptotic EVLs is usually precluded in many practical cases, so that we must proceed heuristically, taking often into consideration external constraints on the availability of data. In this sense, the use of the $g_i$ observables in chaotic systems can help us in understanding the issues related to finiteness and correlation of the data samples, because we are able to predict the exact EVLs., by using the results contained in Chapters 3, 4, and 6.

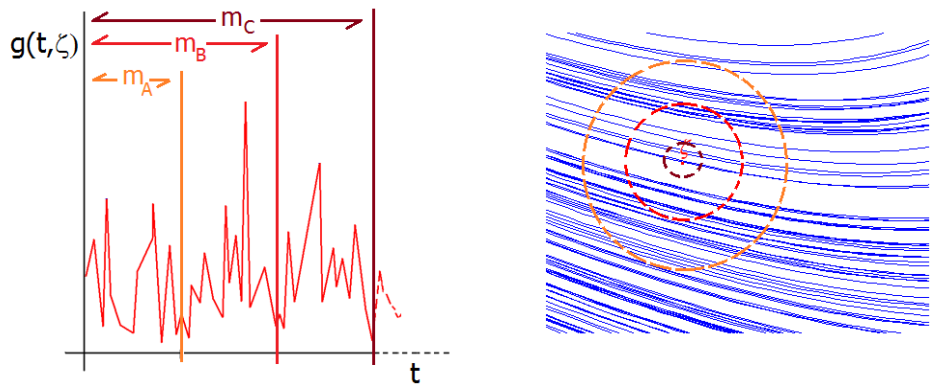

**Figure 9.1** A schematic representation of the BM approach for the observable $g_i$. Increasing the bin length increases the probability of sampling extremes in a smaller ball centred on $\zeta$.

### 9.1.1
### Extreme Value Laws and the Geometry of the Attractor

We would like to provide further details on the relationship between EVLs and the geometry of the attractor $\Omega$ of the dynamical system generating the stochastic process, thus completing some of the calculations presented in Chapter 4 and 6. We recapitulate the necessary ingredients we need in order to proceed further:

- Assumption 1: we consider observables of the form $g_i(dist(x, \zeta))$, $i = 1, 2, 3$ defined in Sect. 4.2.1, and $\zeta \in \Omega$;
- Assumption 2: the dynamical system is such that the time series of the observables obey the conditions $Ɗ_0$ and $Ɗ'_0$;
- Assumption 3: the measure of a ball of radius $r$ centred on the point $\zeta \in \Omega$ is a continuous function of $r$ the radius for almost all the center points $\zeta$; this can be guaranteed simply by requiring that the measure has no atoms.;
- Assumption 4: we can define for almost all the center points $\zeta$ the local dimension following Eq. 4.2.8, so that $d(\zeta) = \lim_{r \to 0} \frac{\log \mu(B_r(\zeta))}{\log r}$.

In Sect. 4.2.1 we have shown how to derive the shape parameters for each of the three corresponding families of observables $g_i(dist(x, \zeta))$, $i = 1, 2, 3$ and have discussed their relationship with the local dimension of the attractor. In some cases, *e.g.* if the system is Axiom A (see Chap. 8), we have seen that distance observables with reference point $\zeta \in \Omega$ can provide global information on the properties of the measure.

The missing ingredients in the picture presented so far are the normalising sequences $a_n$ and $b_n$ (where $n$ refers to the length of the bins used according to the BM approach) first introduced in Eq. 3.1.2. Such sequences are key to avoiding the derivation of degenerate d.f. for extremes, as thoroughly discussed in Theorem 3.1.1. What we want to do here is to find explicit expressions (or at least asymptotic estimates) for $a_n$ and $b_n$ and show their link with the geometrical properties of the attractor.

### 9.1.2
### Computation of the Normalizing Sequences

We proceed by focusing, for illustrative purposes, on the simpler yet relevant case of mixing dynamical systems possessing an absolutely continuous invariant measure $\mu(x)$ supported on the attractor $\Omega$ such that the Radon-Nikodym derivative with respect to the Lebesgue measure - $\rho(x) = \mathrm{d}\mu(x)/\mathrm{d}\,\mathrm{Leb}$, in the given notation - is constant. In this case, $d(\zeta)$ is constant everywhere in the attractor and is equal to the dimension $d$ of the phase space of the system. Low dimensional examples of such are the uniformly hyperbolic maps $f : [0, 1) \to [0, 1)$ of the form $f(x) = qx \bmod 1$ $q \in \mathbb{N}, q \geq 2$ (see Example 4.2.1 for the case $q = 2$), and the algebraic automorphisms of the torus better known as the Arnold cat map, already discussed in Sect. 6.6.1. These systems are studied later in the chapter.

In [46] one can find details on how to adapt the results presented below to the case of systems whose invariant measure is absolutely continuous with respect to Lebesgue, but the Radon-Nikodym derivative is *not* constant.

Constructing the normalising sequences $a_n$ and $b_n$ for distance observables $g$ in the case of systems whose invariant measure is not absolutely continuous with respect to Lebesgue is indeed more challenging. Some results and conjectures are presented in the later Sect. 9.4, where we also discuss the numerical investigations of extremes in various systems possessing a fractal invariant measure.

We go back to the simple setting discussed above. Starting from the definitions and constructions provided by Gnedenko we derive the exact expression for the normalising sequences $a_n$ and $b_n$ which constitute the backbone of EVT. Constructing such sequences have a great practical significance because when performing a GEV fit of BM taken over bins of length $n$, we have that the best fit $GEV_\xi((x - \mu)/\sigma)$ is such that $a_n = 1/\sigma$ and $b_n = \mu$, where, clearly, $\mu = \mu(n)$ and $\sigma = \sigma(n)$. As opposed to $\xi$, whose estimate, when we are in the asymptotic regime, does not depend on $n$, taking maxima over longer and longer bins does in general affect the obtained estimates of $\sigma$ and $\mu$ *also when the asymptotic regime is realized*.

This fact can be seen as a linear change of variable: the variable $y = a_n(x - b_n)$ has a GEV distribution $GEV_\xi(y)$ (that is an EV one parameter distribution with $a_n$ and $b_n$ normalising sequences) while $x$ is GEV distributed according to $GEV_\xi(a_n(x - b_n))$. This allows relating such somewhat abstract sequences introduced by Gnedenko to the output of a numerical inference procedure.

**Case 1: $g_1(x) = -\log(dist(x, \zeta))$.** Using the definition of the d.f. $F(u)$ and of its complement $\bar{F}(u) = 1 - F(u)$, we have

$$
\begin{aligned}
1 - F(u) = \bar{F}(u) &= 1 - \mu(g(\mathrm{dist}(x, \zeta)) \leq u) \\
&= 1 - \mu(-\log(\mathrm{dist}(x, \zeta)) \leq u) \\
&= 1 - \mu(\mathrm{dist}(x, \zeta) \geq e^{-u}).
\end{aligned}
\tag{9.1.1}
$$

Then, for maps with constant density measure, we can write:

$$1 - F(u) = \bar{F}(u) = \mu(B_{e^{-u}}(\zeta)) = \Omega_d e^{-ud} \tag{9.1.2}$$

where $d$ is the dimension of the space and $\Omega_d$ is a constant. In order to use Gnedenko's corollary it is necessary to calculate $u_F = \sup\{u; F(u) < 1\}$, which, in this case, is infinite. Following the Gnedenko approach outlined in Sect. 3.1.1, we can construct a function $h(t)$ as follows:

$$h(t) = \frac{\int_t^{u_F} \mathrm{d}u(1 - F(u))}{1 - F(t)} = \frac{\int_t^{\infty} \mathrm{d}u\, e^{-ud}}{e^{-td}} = \frac{1}{d} \int_{td}^{\infty} \frac{e^{-v}}{e^{-td}} \mathrm{d}v = \frac{1}{d}. \tag{9.1.3}$$

Hence, according to the proof of Gnedenko theorem given in [1], we can study both $a_n$ and $b_n$ or $\gamma_n$ convergence as:

$$\lim_{n \to \infty} m(1 - F\{\gamma_n + xh(\gamma_n)\}) = e^{-x},$$

so that, using Eq. 9.1.2, we obtain:

$$\lim_{n \to \infty} m\Omega_d e^{-d(\gamma_n + xh(\gamma_n))} = \lim_{n \to \infty} m\Omega_d e^{-d(\gamma_n + x/d)} = e^{-x}. \tag{9.1.4}$$

As a result, we obtain that

$$\gamma_n \simeq \frac{\log(n\Omega_d)}{d}$$

so that, since $a_n = 1/h(\gamma_n)$ and $b_n = \gamma_n$, we have:

$$a_n = d = 1/\sigma(n) \qquad b_n = \frac{1}{d}\log(n) + \frac{\log(\Omega_d)}{d} = \mu(n). \tag{9.1.5}$$

From Sect. 4.2.1 we know that $\xi = 0$.

**Case 2: $g_2(x) = (dist(x, \zeta))^{-1/\alpha}$**   We can proceed as previously done for the observable $g_1$. We derive:

$$
\begin{aligned}
1 - F(u) = \bar{F}(u) &= 1 - \mu(\mathrm{dist}(x, \zeta)^{-1/\alpha} \leq u) \\
&= 1 - \mu(\mathrm{dist}(x, \zeta) \geq u^{-\alpha}) \\
&= \mu(B_{u^{-\alpha}}(\zeta)) = \Omega_d u^{-\alpha d}
\end{aligned} \tag{9.1.6}
$$

Also in this case we have that $u_F = +\infty$. We derive that:

$$\gamma_n = F^{-1}(1 - 1/n) = (n\Omega_d)^{1/(\alpha d)} = 1/a_n = \sigma(n) \tag{9.1.7}$$

Gnedenko's theorem suggests that $b_n = 0 = \mu$; another possible choice of the normalising sequences has been proposed by [274] as follows:

$$b_n = c \cdot n^{-\xi} = c \cdot n^{-1/(d\alpha)} = \mu, \tag{9.1.8}$$

where $c \in \mathbb{R}$ is a constant and we know from Sect. 4.2.1 that $\xi = 1/(\alpha d)$.

**Case 3: $g_3(x) = C - dist(x, \zeta))^{1/\alpha}$.** Eventually, we compute $a_n$ and $b_n$ for the $g_3$ observable class:

$$\begin{aligned}
1 - F(u) &= 1 - \mu(C - \text{dist}(x, \zeta)^{1/\alpha} \leq u) \\
&= 1 - \mu(\text{dist}(x, \zeta) \geq (C - u)^{\alpha}) \\
&= \mu(B_{(C-u)^{\alpha}}(\zeta)) = \Omega_d(C - u)^{\alpha d}
\end{aligned} \tag{9.1.9}$$

in this case $u_F = C$.

$$\gamma_n = F^{-1}(1 - 1/n) = C - (n\Omega_d)^{-1/(\alpha d)} \tag{9.1.10}$$

For type 3 distribution:

$$a_n = (u_F - \gamma_n)^{-1} = (n\Omega_d)^{1/(\alpha d)} = 1/\sigma, \quad b_n = u_F = \mu. \tag{9.1.11}$$

*Remark* 9.1.1. The previous expressions given in Eqs. 9.1.5, 9.1.8, and 9.1.11 provide important constraints on what our statistical inference procedure should give as values for $\sigma$ and $\mu$ when we change the length $n$ of bins considered according to the BM method. Additionally, such expressions suggest multiple (and *independent*) ways to derive two important properties of the attractor, namely the dimension $d$ and the constant density $\Omega_d$, from the values of $\sigma$ and $\mu$, which can estimated by fitting the empirical d.f. of BM with the GEV model,

### 9.1.3
### Inference Procedures for the Block Maxima Approach

We present here some methods used for the inference of the GEV distribution parameters describing the extremes of a time series. Our goal is to provide the basic information needed to understand how the methods used may be adapted to the observables extracted from the orbits of dynamical systems. We study the sequence of maxima obtained by subdividing the available data $X_j$, $1 \leq j \leq s \gg 1$ into $k \gg 1$ bins of equal size $n \gg 1$, and extract the maximum $M_j$ from each bin, $1 \leq j \leq k$. We point the reader to the problem that given a time series of length $s$, many choices are possible for constructing the bins and selecting their maxima, boiling down to different combinations of values of $n$ $k$ such that $s = nk$. The selection of suitable candidates for true extremes requires choosing large values for $n$, while the need for achieving a robust statistical fit requires choosing large values for $k$, so that a compromise needs to be found in all practical circumstances and the robustness of the results of the inference procedure must be accurately tested [3]. The results discussed in remark 9.1.1 provide prescriptions on what we should obtain for the best estimates of $\mu$, $\sigma$, and $\xi$ as a function of $n$.

Several approaches have been proposed for estimating the parameters of the GEV distribution. We present here a brief summary and refer instead to [112] and to [3] for comprehensive descriptions.

We first focus on the Maximum Likelihood Estimation (MLE), which is most commonly implemented in the software packages currently used for the analysis of extremes as default choice, and then on the L-moment procedure - based on the computation of specific moments of the distribution of maxima and versatile enough to be implemented also by non expert programmers. The MLE procedure relies on maximizing the log likelihood function:

$$l(\mu, \sigma, \xi) = \prod_{j=1}^{k} \log(f_{GEV}(M_j; \mu, \sigma, \xi)) \tag{9.1.12}$$

where

$$f_{GEV}(x; \mu, \sigma, \xi) = \frac{\mathrm{d}}{\mathrm{d}x} GEV_\xi \left( \frac{x - \mu}{\sigma} \right).$$

Using Eq. 3.1.4, the log likelihood function can be rewritten as:

$$-m \log(\sigma) - \left(1 + \frac{1}{\xi}\right) \sum_{j=1}^{k} \left\{ \log \left[ 1 + \xi \left( \frac{M_j - \mu}{\sigma} \right) \right] - \left[ 1 + \xi \left( \frac{M_j - \mu}{\sigma} \right) \right]^{-\frac{1}{\xi}} \right\} \tag{9.1.13}$$

if $\xi \neq 0$, and as:

$$-m \log(\sigma) - \sum_{j=1}^{k} \left\{ \left( \frac{M_j - \mu}{\sigma} \right) - \exp \left[ - \left( \frac{M_j - \mu}{\sigma} \right) \right] \right\} \tag{9.1.14}$$

if $\xi = 0$.

We can obtain a profile likelihood of $\mu, \xi$ or $\sigma$ by setting the other two parameters to their maximum likelihood estimates $\tilde{\mu}, \tilde{\xi}, \tilde{\sigma}$ in Eqs. 9.1.13 or 9.1.14. For example, to compute the profile likelihood for the parameter $\xi$, we can construct the graph:

$$(x, y) = (\xi, l(\tilde{\mu}, \tilde{\sigma}, \xi)). \tag{9.1.15}$$

The intersections of the horizontal line with the profile likelihood graph allows for estimating a confidence interval:

$$y = \tilde{\xi} - 0.5 q_{0.95}, \tag{9.1.16}$$

where $q_{0.95}$ is the 95% quantile of the $\chi^2$ distribution with 1 degree of freedom.

In the source code provided along with the book (see Appendix A), the functions used for the GEV model inference via MLE give as output the 95% confidence intervals for the estimates of the parameters.

Whenever the probability density function is not absolutely continuous, the maximum likelihood estimation may fail as the minimization procedure based on the derivatives of the d.f. is not well defined when the density presents singularities causing unexpected divergences of the parameters. In these cases is better to rely on a L-moment estimation based upon the computation of integrals rather than upon the solution of a variational problem. The L-moment are analogous to ordinary moments, but are computed from linear combinations of the data values, arranged in an increasing order. For a random variable $X_i$, the $r$-th population L-moment is:

$$\lambda_r = r^{-1} \sum_{p=0}^{r-1} (-1)^p \binom{r-1}{p} \mathrm{E}[X_{r-p:r}], \tag{9.1.17}$$

where $X_{p:q}$ denotes the $p^{th}$ order statistics - $p^{th}$ largest value - in an independent sample of size $q$ from the distribution of $X$ and $\mathrm{E}[\,]$ denotes the expected value. In particular, the first four population L-moment are

$$\lambda_1 = \mathrm{E}[X] \tag{9.1.18}$$

$$\lambda_2 = (\mathrm{E}[X_{2:2}] - \mathrm{E}[X_{1:2}])/2 \tag{9.1.19}$$

$$\lambda_3 = (\mathrm{E}[X_{3:3}] - 2\mathrm{E}[X_{2:3}] + \mathrm{E}[X_{1:3}])/3 \tag{9.1.20}$$

$$\lambda_4 = (\mathrm{E}[X_{4:4}] - 3\mathrm{E}[X_{3:4}] + 3\mathrm{E}[X_{2:4}] - \mathrm{E}[X_{1:4}])/4. \tag{9.1.21}$$

$\lambda_1$ is conventionally referred to as *mean, L-mean* or *L-location* and $\lambda_2$ as *L-scale*. It has been shown that asymptotic approximations to sampling distributions are better for L-moment than for ordinary moments [114]. Moreover, the estimation provided in the i.i.d. case and the associated uncertainties are comparable to the MLE method. Whenever we have used this inference, we have computed confidence intervals using dispersion indicators of an ensemble of realizations of a particular map. In fact, when the data are correlated, the usual bootstrap procedure based on a reshuffling of the sample fails in giving reliable parameters uncertainty estimation due to the loss of information about the dependency structure of the data.

## 9.2
## The Peaks Over Threshold Approach

As mentioned before in the book, the POT method consists in choosing a threshold value $T$ and considering as extremes all the observations *above* the threshold. It is widely used in hydrological applications, since rivers and lakes clearly present banks whose height can be taken as a threshold. In the same way, in financial applications is common to use the POT approach and set critical values (associated to specific risk scenarios) as thresholds.

We have discussed in detail in Sect. 8.2.1 how to construct and interpret the POT method for distance observables. We remind the reader that, when considering the family of observables $g_i(dist(x, \zeta))$, $i = 1, 2, 3$, the selection of the threshold $T$ corresponds exactly to a specific choice of a radius $r^*$ such that all extremes belong to the intersection of the balls $B_{g^{-1}_i(T)}(\zeta)$ with the attractor of the system. Therefore, when studying extreme events, we are sampling the invariant measure in the vicinity of $\zeta$ with the desired precision as defined by the value of $r^*$. In Fig. 9.2, this is schematically depicted by showing how three different threshold values $T_A, T_B, T_C$ correspond to three specific radii.

It is immediate to see that by choosing the POT approach instead of the BM approach we lose some information about the dynamics of the system, because we restrict ourselves to a purely geometrical criterion for selecting extremes, based on belonging to the ball of radius $r^*$ centred on $\zeta$.

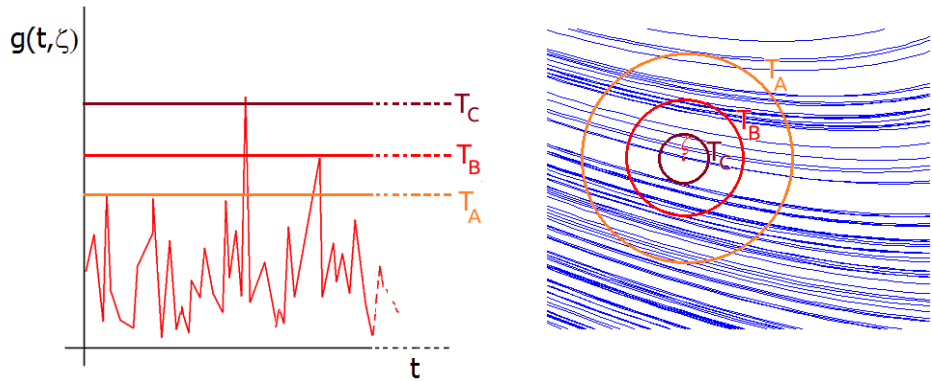

**Figure 9.2** A schematic representation of the POT approach for the observable $g_i$.

#### 9.2.0.1 Inference Procedures for the Peaks Over Threshold Approach

Similarly to the BM method, we construct a MLE of the parameters of the GPD by maximizing the log likelihood function:

$$l(\mu, \sigma, \xi) = \prod_{j=1}^{k} \log(f_{GPD}(y_j; \xi, \sigma)) \qquad (9.2.1)$$

where $f_{GPD}(x; \xi, \sigma)$ is defined in Eq. 8.2.21. Using Eq. 3.1.8, we obtain:

$$\log l(\sigma, \xi) = -n \log \sigma - (1 - \xi)\Sigma_{j=1}^{k} y_j, \quad y_j = -\xi^{-1}\log(1 - \xi(X_j - T)/\sigma), \quad (9.2.2)$$

where we assume that $k$ observations $X_j$ are beyond the threshold $T$. Giving a well-defined meaning to the maxima of the log likelihood function is not always entirely trivial. The previous expression may be made arbitrary large by taking $\xi > 1$ and $\sigma/\xi$ close to $\max(X_i)$. The maximum likelihood estimators are taken to be $\bar{\xi}$ and $\bar{\sigma}$, which yield a local maximum of $\log l$ [275]. For $k < 1/2$ a local maximum is well defined and the asymptotic variance for the estimators can be written as:

$$L \cdot Var \left[ \begin{array}{c} \bar{\sigma} \\ \bar{\xi} \end{array} \right] \sim \left[ \begin{array}{cc} 2\sigma^2(1-\xi) & \sigma(1-\xi) \\ \sigma(1-\xi) & (1-\xi)^2 \end{array} \right], \quad \xi < 1/2. \qquad (9.2.3)$$

The explicit expression for $\xi > 1/2$ is more problematic, but we can avoid to consider this range of parameter values by making specific choices for the exponent appearing in the definition of the $g_2$ functions.

The inference based on the L-moment procedure follows exactly the set-up already described for the GEV distribution. As we have already said for the MLE, if $\xi > 1/2$, the inference becomes problematic as the L-scale $\lambda_2$ is not defined anymore. The exact expression of the parameters as a function of the L-moment are presented in [114].

## 9.3
## Numerical Experiments: Maps Having Lebesgue Invariant Measure

We begin our journey on numerical investigations of the relationship between extremes and geometry of the attractor of the underlying system by focusing on maps possessing invariant measures proportional to Lebesgue. We mostly focus on the BM method and related GEV results, and then compare the findings with what obtained using the POT method.

Recalling remark 9.1.1, and keeping in mind that the total length of the time series is defined as $s = n \cdot k$ and is constant, we can express the asymptotic values of the GEV normalising sequences and, correspondingly, of the best GEV estimates for $\sigma$, $\mu$, and $\xi$ as a function of $k$ for the three types of $g$ distance observables discussed in Sect. 9.1.1 as follows:

- For $g_1$-type observables:

$$\sigma = \frac{1}{d} \qquad \mu = C_1 + \frac{1}{d}\log(n) = C_2 - \frac{1}{d}\log(k) \qquad \xi = 0, \qquad (9.3.1)$$

  where $C_1$ and $C_2$ is a positive constant.

- For $g_2$-type observables:

$$\sigma \propto n^{1/(\alpha d)} \propto k^{-1/(\alpha d)} \qquad \mu \propto n^{1/(\alpha d)} \propto k^{-1/(\alpha d)} \qquad \xi = \frac{1}{\alpha d} \quad (9.3.2)$$

- For $g_3$-type observables:

$$\sigma \propto n^{-1/(\alpha d)} \propto k^{1/(\alpha d)} \qquad \mu = C \qquad \xi = -\frac{1}{\alpha d} \qquad (9.3.3)$$

where $x \propto y$ indicates that $x/y$ is a constant.

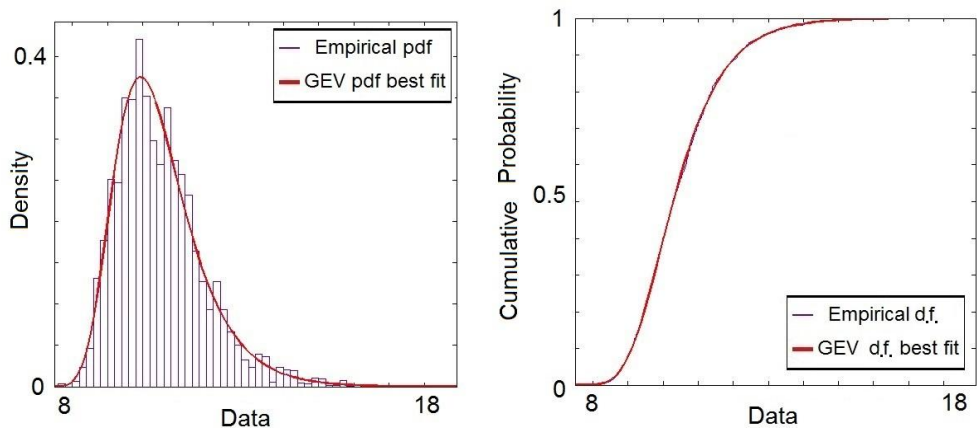

**Figure 9.3** Left: $g_1$ observable empirical histogram and fitted GEV pdf. Right: $g_1$ observable empirical cdf and fitted GEV cdf. Logistic map, $n = 10^4$, $m = 10^4$

In order to provide a numerical test of our results, we consider a one-dimensional and a two dimensional map. The one dimensional map considered here is a Bernoulli shift map defined on $\mathbb{S}^1$ and already introduced in the general form Eq. 4.2.1:

$$f(x) = qx \mod 1, \quad q = 3 \tag{9.3.4}$$

Note that one should avoid studying numerically the doubling map obtained setting $q = 2$ in Eq. 4.2.1, because, as a result of the binary representation of numbers in computers, all orbits converge to 0 when performing long trajectories, regardless of the initial condition. The considered two dimensional map is the Arnold's cat map defined on the 2-torus $\mathbb{T}^2$, which we have already presented in Eq. 6.6.1:

$$f(x, y) = (2x + y, x + y) \mod \mathbb{T}^2 \tag{9.3.5}$$

An extensive description of the properties of these maps can be found in [276] and [277]. In this book, we have discussed the Bernoulli shift map in example 4.2.1 (albeit for the case $q = 2$), and the Arnold cat map in Sect. 6.6.1.

For each map we run a long simulation starting from a given initial condition $x_0$. From the trajectory we compute the sequence of observables $g_1, g_2, g_3$ with respect to a point $\zeta$ randomly chosen according to the invariant measure. We then divide the time series of length $s$ into $k$ bins each containing $n = s/k$ observations. Then, we test the degree of agreement between the empirical distribution of the maxima and the GEV distribution $GEV_\xi((x - \mu)/\sigma)$ using the MLE method.

All the numerical analysis contained in this work has been performed using MAT-LAB© Statistics Toolbox function `gevfit`. These functions return MLE of the parameters for the GEV distribution giving 95% confidence intervals for estimates. As in every fitting procedure, it is necessary to test the a posteriori goodness of fit. We anticipate that in every case considered, fitted distributions passed, with maximum confidence interval, the Kolmogorov-Smirnov test described in [278]. For illustration purposes, we present in figure 9.3 an empirical pdf and d.f. with the corresponding best GEV fits.

Once $s$ is set to a given value (in our case $s = 10^7$), the numerical simulations allow us to explore two limiting cases of great interest in applications where the statistical inference is intrinsically problematic:

1) $k$ is small ($n$ is large), so that we extract only few maxima, each corresponding to a very good candidate for being a *truly* extreme event.
2) $n$ is small ($k$ is large), so that we extract many maxima but most of those will not be, in fact, *soft* extremes.

In case 1), we have only few data - of high quality - to fit our statistical models, whereas in case 2) we have many data but the sampling may be polluted by the inclusion of data not really corresponding to genuine extreme events. In general, we have that in order to obtain a reliable fit for a distribution with $p$ parameters we need $10^p$ independent data [45] so that we expect that the fit procedure gives reliable results for $k \geq k_{min} \sim 10^3$. As the value of $n$ determines to which extent the extracted bin maximum is representative of an extreme, below a certain value $n_{min}$ our selection procedure will be unavoidably misleading. We have no obvious theoretical argument to define the value of $n_{min}$. We expect to obtain good fits throughout the parametric region where both constraints on $n$ and $k$ are satisfied. Clearly, our flexibility in choosing *good* pairs $(n, k)$ increases with larger values of $s$, *i.e.*, when we have longer time series at hand.

For $g_1$-type observables, the behaviour of the three parameters against $k$ is presented in figure 9.4. According to Eq. 9.3.1, we expect to find $\xi = 0$. For relatively small values of $k$ the sample is too small to ensure a good convergence to analytical $\xi$ and confidence intervals are wide. On the other hand we see deviations from the expected value as $n < 10^3$, *i.e.*, when $k > 10^4$. For the scale parameter a similar behavior is achieved and deviations from expected theoretical values $\sigma = 1/2$ for the Arnold cat map and $\sigma = 1$ for the Bernoulli shift map are found when $k < 10^3$ or $n < 10^3$. The location parameter $\mu$ shows a logarithm decay with $k$, as expected from Eq. 9.3.1.

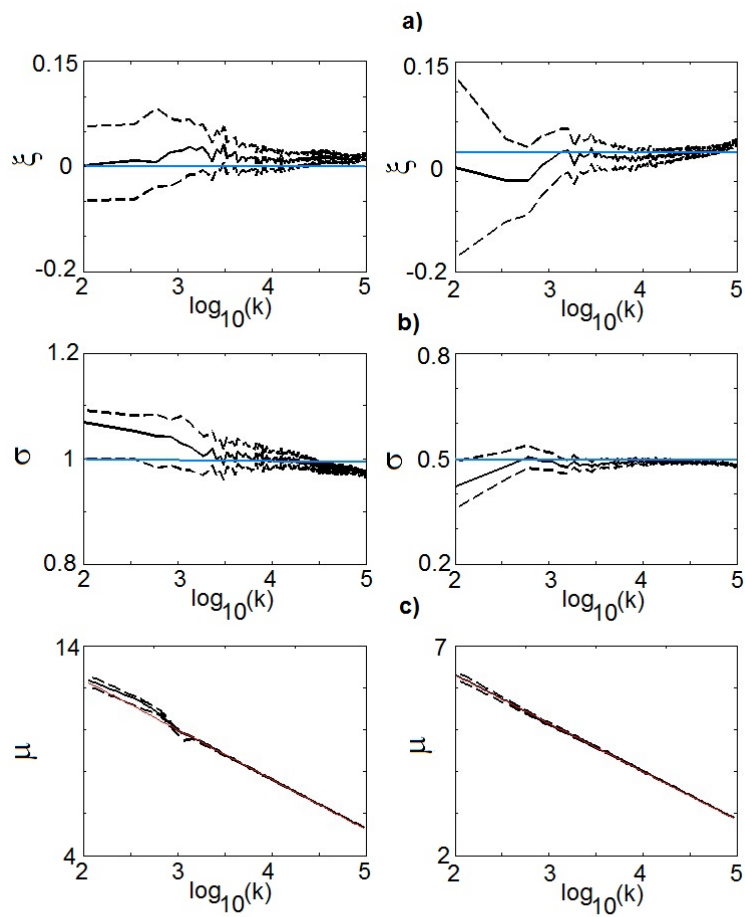

**Figure 9.4** $g_1$ observable, **a)** $\xi$ VS $\log_{10}(k)$; **b)** $\sigma$ VS $\log_{10}(k)$; **c)** $\mu$ VS $\log_{10}(k)$. Right: Bernoulli shift map. Left: Arnold cat map. Dotted lines represent the computed confidence intervals, gray lines represent linear fits and theoretical values.

A linear fit of $\mu$ against $\log_{10}(k)$ is shown with a gray line in figure 9.4. The estimated angular coefficients $K^*$ agree very accurately with the theoretically predicted value of $1/d$ given in Eq. 9.3.1: for the Bernoulli shift map we obtain $|K^*| = 1.001 \pm 0.001$, while for Arnold cat map $|K^*| = 0.489 \pm 0.001$.

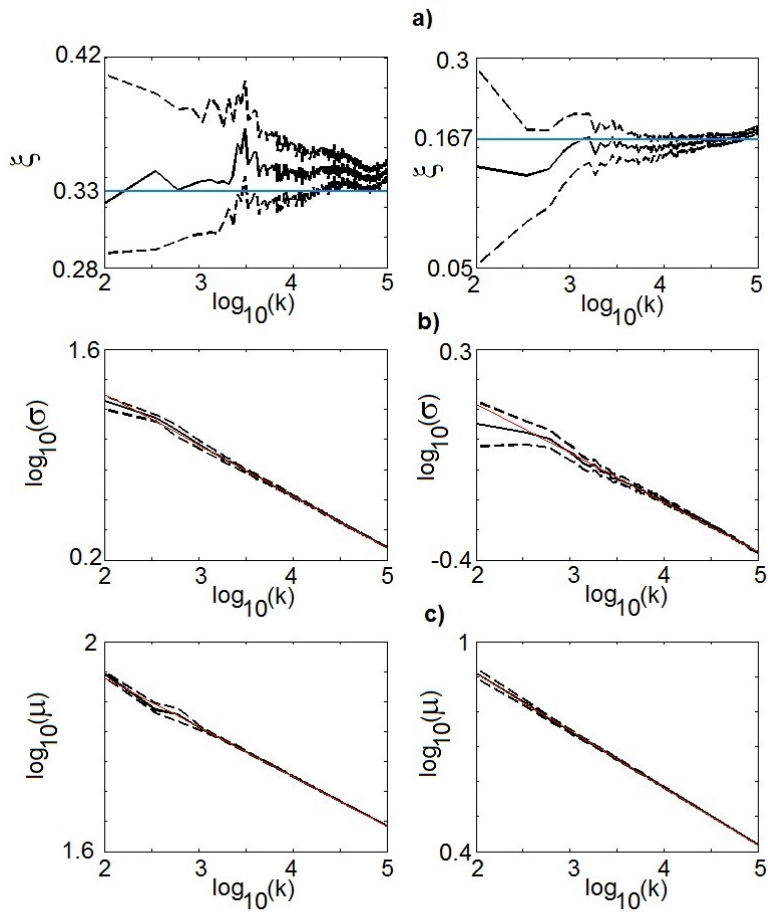

**Figure 9.5** $g_2$ observable. **a)** $\xi$ VS $\log_{10}(k)$; **b)** $\log_{10}(\sigma)$ VS $\log_{10}(k)$; **c)** $\log_{10}(\mu)$ VS $\log_{10}(k)$. Right: Bernoulli shift map. Left: Arnold cat map. Dotted lines represent computed confidence interval, gray lines represent linear fits and theoretical values.

We find that $\xi$ values have best matching with theoretical ones with reliable confidence interval when both $n > 10^3$ and $k > 10^3$.

In the following, without loss of generality, we discuss the main findings obtained for $g_2$-type and $g_3$-type observables taking $\alpha = 3$ (see Sect. 9.1.1).

For a $g_2$ observable function we expect to have $\xi = 1/3$ for the Bernoulli shift map and $\xi = 1/6$ for the Arnold cat map. In both cases - see Fig. 9.5a) - the best match is obtained for $n > 10^3$ and $k > 10^3$. We can also check that $\mu$ and $\sigma$ parameters follow a power law as a function of $k$ as described in Eq. 9.3.2. In the log-log plot in Figure 9.5b), 9.5c), we can see a very clear linear behavior. For the Bernoulli shift map, we obtain $|K^*| = 0.330 \pm 0.001$ for $\mu$ series, $|K^*| = 0.341 \pm 0.001$ for $\sigma$, in good agreement with theoretical value of $1/3$. For the Arnold cat map we expect to find $K^* = 1/6$, from the experimental data we obtain $|K^*| = 0.163 \pm 0.001$ for $\mu$ and $|K^*| = 0.164 \pm 0.001$ for $\sigma$.

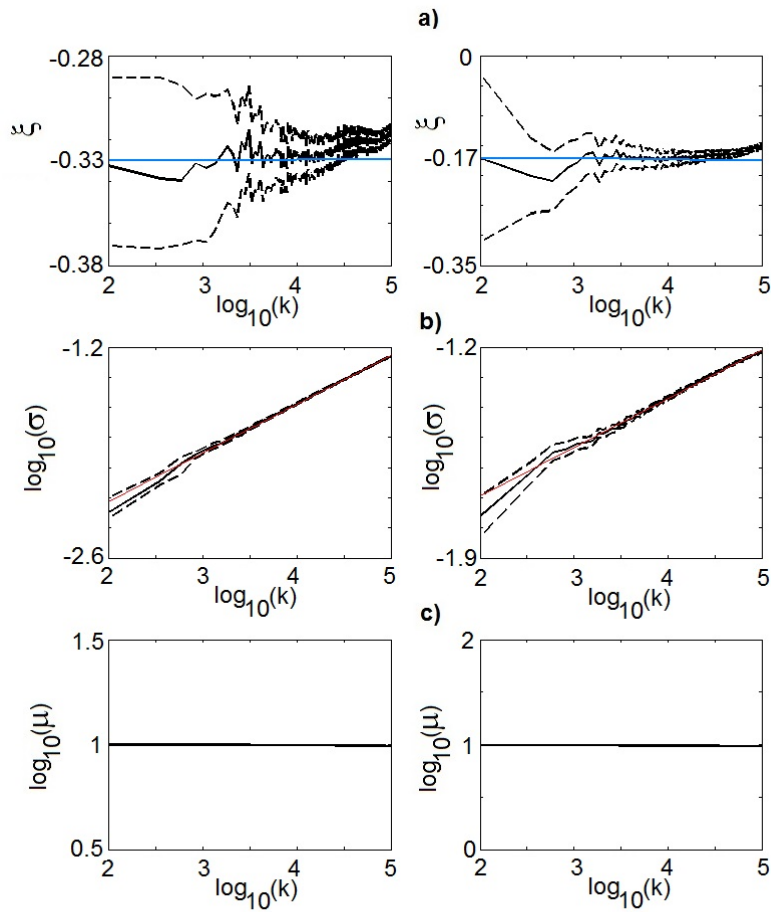

**Figure 9.6** $g_3$ observable, **a)** $\xi$ VS $\log_{10}(k)$; **b)** $\log_{10}(\sigma)$ VS $\log_{10}(k)$; **c)** $\log_{10}(\mu)$ VS $\log_{10}(k)$. Right: Bernoulli shift map. Left: Arnold cat map. Dotted lines represent computed confidence interval, gray lines represent linear fits and theoretical values.

Finally, we consider $g_3$-type observables, where we choose $C = 10$ and we expect to have $\xi = -1/3$ for the Bernoulli shift map and $\xi = -1/6$ for the Arnold cat map. In both cases, as above the best match is obtained for $n > 10^3$ and $k > 10^3$ - see Fig. 9.6a).

The results presented in Eq. 9.3.3 suggest that in the asymptotic regime, $\mu$ is constant, while $\sigma$ grows following a power law of $k$. As in the $g_2$ case we expect that the value of the corresponding slope of the graph in the log-log plot is $|K^*| = 1/(\alpha d)$. The numerical results are shown in Figs. 9.6b), 9.6c) and are consistent with the theoretical ones, since we find $|K^*| = 0.323 \pm 0.006$ for Bernoulli shift map and $|K^*| = 0.162 \pm 0.006$ for the Arnold cat map.

### 9.3.1
### Maximum Likelihood vs L-moment Estimators

We want now to compare the performance of two main methods of performing GEV statistical inference for BM described in Sec. 9.1.3. In Fig. 9.7 we show the estimates for the shape parameters $\xi$ of the two maps given in Eqs 4.2.1-**??** obtained using the MLE method (which has been used to produce the figures above in this chapter) and the L-moment method. For both methods, we display (solid line) the mean value of $\xi$ obtained by averaging the results of an ensemble of 30 realizations started from random initial conditions, whereas the shaded regions represent cover one standard deviation of the mean.

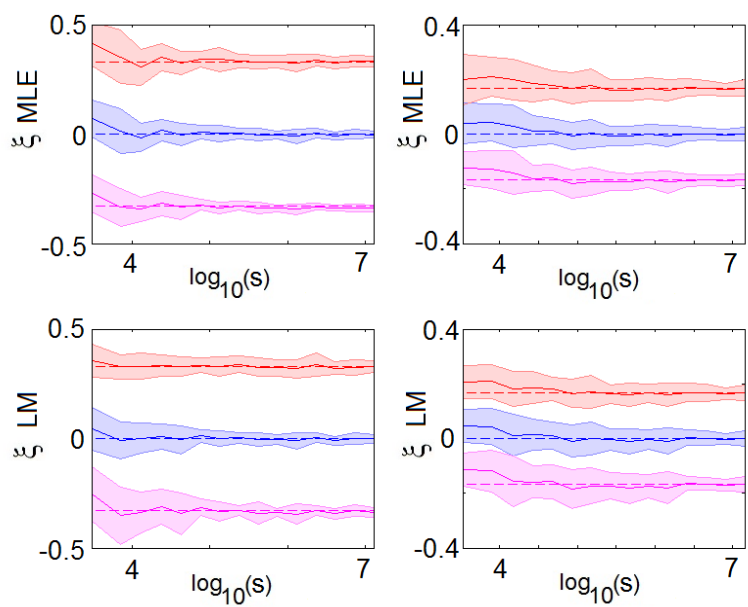

**Figure 9.7** Shape parameters of the GEV distribution $\xi$ VS total length of the series $s = NK$. Solid lines represent the mean values among 30 realizations, the shaded regions represent the standard deviation of the mean. The dotted lines are the theoretical expected values; blue: $g_1$, red: $g_2$, magenta: $g_3$. Top: MLE , bottom: L-moment. Left: Bernoulli shift map, right: Arnold cat map.

Maxima of $g_i$ are extracted every $n$ observations over a simulation containing $s = nk$ points, where $k = 10^3$, so that we consider simulations of varying length, as opposed to above, where the length $s$ is kept fixed. The value of $k$ has been chosen, following what discussed before, in order to follow the idea that in order to obtain a reliable fit for a distribution with $p$ parameters we need $10^p$ independent data [45].

The colors code refer to the three different observables: blue is used for $g_1 = -\log(dist(x, \zeta))$, red for $g_2 = dist(x, \zeta)^{-3}$ and magenta for $g_3 = -(dist(x, \zeta))^3$. The observables are computed with respect to a point $\zeta$ randomly chosen on the attractor. The dotted lines represent the theoretical expected values discussed above.

The analysis shows that in this case the two statistical inference methods have com-parable convergence speed and precision.We first observe that the value of the low threshold $n_{min}$ such that the selection procedure leads to a biased parameters estimation depends on the properties of the map. For both maps the confidence interval of the estimates for the values of $\xi$ includes the theoretical values already for $n = 10$, which corresponds to $s = 10^4$. Note that, if we increase the number of maxima we consider in order to decrease the uncertainty on the inferred parameters, we in fact discover that our best fit is not compatible with the theory; compare the rightmost range of Figs 9.4a), 9.5a), 9.6a), where the case $s = 10^7$, $k = 10^5$ and $n = 10^2$ are considered.

## 9.3.2
## Block Maxima vs Peaks Over Threshold Methods

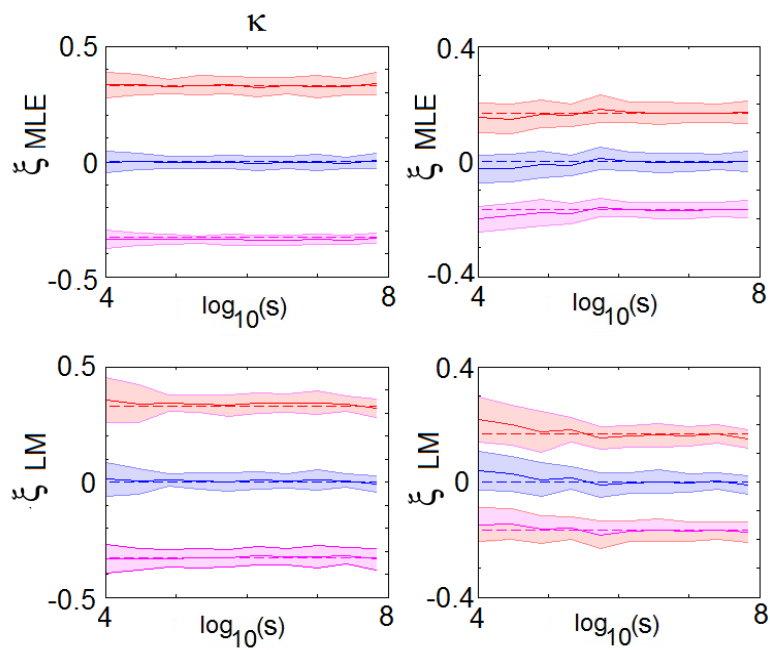

**Figure 9.8** Shape parameters of the GPD $\xi$ VS total length of the series $s = NK$. Solid lines represent the mean values among 30 realizations, the shaded regions the standard deviation of the mean. The dotted lines represent the values predicted by the theory. blue: $g_1$, red: $g_2$, magenta: $g_3$. Top: MLE , bottom: L-moment. Left: Bernoulli shift map, right: Arnold cArnold cat mapat map.

We repeat the analysis described in the previous section by using the POT approach. As in the previous case, we refer to [45] for the choice of the number of maxima to be extracted, so that we choose $k = 10^3$. The total length of the series $s$ is gradually increased through the experiments, which are devised similarly to what described in the previous section for the BM method. The results are shown in Fig. 9.8. The agreement between theory and experiments is very good - and with rather narrow confidence intervals as compared to what reported in Fig. 9.7 - already for very short time series $NK = 10^4$, due to the particular properties of the measure considered. In fact, for uniform densities, $d(\zeta) = 1$ no matter the radius of the ball considered. The POT approach allows for achieving a better convergence for the two dimensional case with respect to the BM approach. This can be explained by considering that in the particular case of measure having a constant density, the choice of the radius does not really affect the EVLs as the same scaling properties hold ubiquitously on the invariant measure. We remark that, as in the previous case, no particular differences can be observed when the inference procedure is performed via the MLE or the L-moment method.

## 9.4
## Chaotic Maps with Singular Invariant Measures

In Sect. 4.2.1 we have shown that it is indeed possible to construct EVLs for the distance observables $g_i$'s, $i = 1, 2, 3$ also in the case of mixing systems possessing a fractal invariant measure without atoms supported on the attractor $\Omega$. This case is extremely relevant for a variety of applications, and, in particular, when constructing models for nonequilibrium statistical mechanical systems [70, 71].

In particular, if one can define the local dimension $d(\zeta)$ for a point $\zeta \in \Omega$, then the value of the shape parameters $\xi$ of the GEV distribution describing the BM extremes of the observable can be derived as follows:

- Observable $g_1(dist(x, \zeta)) = -\log(dist(x, \zeta)) \rightarrow \xi = 0$;
- Observable $g_2(dist(x, \zeta)) = -(dist(x, \zeta))^{-1/\alpha} \rightarrow \xi = 1/(\alpha d(\zeta))$;
- Observable $g_3(dist(x, \zeta)) = C - (dist(x, \zeta))^{1/\alpha} \rightarrow \xi = -1/(\alpha d(\zeta))$.

These results are the natural extension of what discussed in Sect. 9.1.1 in the case of systems possessing absolutely continuous invariant measure.

We wish to remind that, as discussed in Chapter 8, deriving EVLs descriptive of the asymptotic behaviour of POT requires assuming that the mass of the ball of radius r entered in $\zeta$ can be approximated as $\mu(B_r(\zeta)) \sim f(r)r^{d(\zeta)}$, where $f$ is a slowly varying function [44]. Correspondingly, in the BM approach, in the case of singular invariant measures it is not possible to derive the properties of the normalising sequences $a_n$ and $b_n$ with the same level of precision as in the case of systems possessing Lebesgue as an invariant measures. We present below some related results.

### 9.4.1
### Nomalizing Sequences

**Case 1: $g_1(x) = -\log(dist(x, \zeta))$.** Substituting Eq. 4.2.3 into Eq. 4.2.8 we obtain that:

$$\begin{aligned}
1 - F(u) &= 1 - \mu(g(\text{dist}(x, \zeta)) \leq u) \\
&= 1 - \mu(-\log(\text{dist}(x, \zeta)) \leq u) \\
&= \mu(\text{dist}(x, \zeta) < e^{-u}) = \mu(B_{e^{-u}}(\zeta)) \qquad (9.4.1)
\end{aligned}$$

where $u_F = \infty$. Following Sect. 3.1.1, we have that according to Corollary 1.6.3 in [1], for type 1 the following results hold: $a_n = [h(\gamma_n)]^{-1}$ and $b_n = \gamma_n = F^{-1}(1 - \frac{1}{n})$, where $h$ is defined in Eq. 4.2.7. We now show how to get the limiting value of $\gamma_n$; a similar proof will hold for type II and III.

**Proposition 9.4.1.** *Let us suppose that our system verifies Assumptions 1, 2, 3, 4 given in Sect. 9.1.1 and let us consider the observable $g_1$; then:*

$$\lim_{n \to \infty} \frac{\log n}{\gamma_n} = d(\zeta).$$

*Proof* Using the definition of $g_1$, we have that $1 - F(\gamma_n) = \mu(B_{e^{-\gamma_n}}(\zeta)) = \frac{1}{n}$. Since the measure is not atomic and varies continuously with the radius, we have necessarily that $\gamma_n \to \infty$ when $m \to \infty$. Now we set $\delta > 0$ and small enough; there is $n_{\delta,\zeta}$ depending on $\delta$ and on $\zeta$, such that for any $n \geq n_{\delta,\zeta}$ we have

$$-\delta\gamma_n \leq \log \mu(B_{e^{-\gamma_n}}(\zeta)) + d(\zeta)\gamma_n \leq \delta\gamma_n. \tag{9.4.2}$$

Since $\log n - d(\zeta)\gamma_n = -[\log \mu(B_{e^{-\gamma_n}}(\zeta)) + d(\zeta)\gamma_n]$ and by using the bounds in the previous Eq. 9.4.2 we immediately have

$$-\delta\gamma_n \leq \log n - d(\zeta)\gamma_n \leq \delta\gamma_n,$$

which proves the Proposition.

It should be clear that the previous proposition will not give us the value of $\gamma_n = b_n$. We have instead a rigorous limiting behavior:

$$\gamma_n = b_n \sim \frac{1}{d(\zeta)} \log n$$

The values for finite $n$ could be obtained if one knew the functional dependence of $\mu(B_r(\zeta))$ on the radius $r$ and the center $\zeta$. The same reason prevents us to get a rigorous limiting behavior for $a_n = [h(\gamma_n)]^{-1}$. The only rigorous statement we can state is that $h(\gamma_n) = o(\gamma_n)$. This follows by adapting the previous proof of the proposition to another result (see [1]) which says that for type I observables one has $\lim_{n \to \infty} n(1 - F\{\gamma_n + xh(\gamma_n)\}) = e^{-x}$, for all real $x$: choosing $x = 1$ gives us the previous domination result. In the following and again for numerical purposes we will assume the validity of the following approximate formula

$$a_n = [h(\gamma_n)]^{-1} \sim \frac{1}{d(\zeta)},$$

which follows by replacing in Eq. 9.4.1 the simple scailng law $\mu(B_r(\zeta)) \sim r^d(\zeta)$ (whose relevance has been already critically discussed in Chap. 8) for $r$ small.

**Case 2:** $g_2(x) = dist(x, \zeta)^{-1/\alpha}$. In this case we have

$$
\begin{aligned}
1 - F(u) &= 1 - \mu(\mathrm{dist}(x, \zeta)^{-1/\alpha} \leq u) \\
&= 1 - \mu(\mathrm{dist}(x, \zeta) \geq u^{-\alpha}) \\
&= \mu(B_{u^{-\alpha}}(\zeta))
\end{aligned}
\tag{9.4.3}
$$

and $x_F = +\infty$. Since $b_m = 0$ we have only to compute $a_n$ which is the reciprocal of $\gamma_n$ which is in turn defined by $\gamma_n = F^{-1}(1 - 1/n)$. By adapting Proposition 9.4.1 we immediately get that

$$\lim_{n \to \infty} \frac{\log n}{\log \gamma_n} = \alpha d(\zeta)$$

so that we derive the approximate relation $a_m \sim n^{-1/(\alpha d(\zeta))}$.

**Case 3: $g_3(x) = C - dist(x, \zeta)^{1/\alpha}$.** We have:

$$
\begin{aligned}
1 - F(u) &= 1 - \mu(C - \text{dist}(x, \zeta)^{1/\alpha} \leq u) \\
&= \mu(B_{(C-u)^\alpha}(\zeta))
\end{aligned}
\tag{9.4.4}
$$

In this case $x_F = C < \infty$ and $a_n = (C - \gamma_n)^{-1}$; $b_n = C$. The previous proposition immediately shows that $\lim_{n \to \infty} \frac{\log n}{-\alpha \log(C - \gamma_n)} = d(\zeta)$, which gives the asymptotic scaling $\gamma_n \sim C - \frac{1}{n^{\frac{1}{\alpha d(\zeta)}}}$; $a_n \sim n^{\frac{1}{\alpha d(\zeta)}}$; $b_n = C$.

Concluding, we derive the following prescriptions for the estimates of the parameters $\xi$, $\sigma$, and $\mu$ for given values of $n$ or $k$. They mirror exactly what was given in Eqs. 9.3.1-9.3.3, except that we have specific reference to the local dimension $d(\zeta)$.

- For $g_1$-type observables:

$$\sigma = \frac{1}{d(\zeta)} \qquad \mu = C_1 + \frac{1}{d(\zeta)} \log(n) = C_2 - \frac{1}{d(\zeta)} \log(k) \qquad \xi = 0, \tag{9.4.5}$$

where $C_1$ and $C_2$ are positive constants.

- For $g_2$-type observables:

$$\sigma \propto n^{1/(\alpha d(\zeta))} \propto k^{-1/(\alpha d(\zeta))} \qquad \mu \propto n^{1/(\alpha d(\zeta))} \propto k^{-1/(\alpha d(\zeta))} \qquad \xi = \frac{1}{\alpha d(\zeta)}. \tag{9.4.6}$$

- For $g_3$-type observables:

$$\sigma \propto n^{-1/(\alpha d(\zeta))} \propto k^{1/(\alpha d(\zeta))} \qquad \mu = C \qquad \xi = -\frac{1}{\alpha d(\zeta)}. \tag{9.4.7}$$

## 9.4.2
## Numerical Experiments

In general, singular measures resulting from chaotic dynamics correspond to observing at least a direction of contraction along which the relevant invariant measure appears as a supported on a Cantor set [71]. Hence, the empirical d.f. of the distance (with respect to $\zeta$ in the attractor) observables feature non trivial behaviour (lack of smoothness) corresponding to range of distance from $\zeta$ where lacunae appear in the attractor. Such lack of smoothness indeed interferes with the optimisation procedure performed in the MLE, as the smoothness properties of the cost function are jeopardized. Conversely, the L-moment procedure allows for overcoming these difficulties since the normalization procedure carried out with this method consists in dividing each quantity computed via L-moment by a function of the total number of data, so that no derivatives of the d.f. are involved.

### 9.4.2.1 A First Example: the Cantor Set

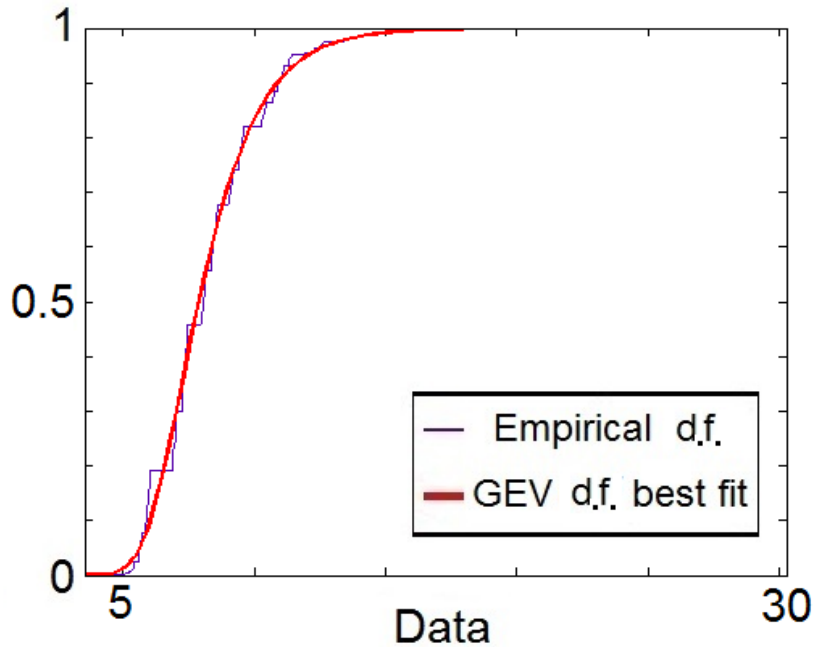

**Figure 9.9** Empirical (blue) and fitted (red) d.f. for the extremes of the distance observable $g_1$ for the IFS that generates a Cantor set. Reference point $\zeta = 0.775$,

In order to illustrate the additional challenges related to finding a satisfactory statistical model for extremes of distance observables for systems possessing a singular invariant measure, we first briefly analyse the so-called middle one third Cantor set, which can be obtained by the Iterating Function System (IFS) $\{f_1, f_2\}$ defined as:

$$
\begin{aligned}
f_1(x) &= x/3 & p < 0.5 \\
f_2(x) &= (x+2)/3 & p \geq 0.5,
\end{aligned} \qquad (9.4.8)
$$

where $x \in [0,1]$ and we set $p$ as a random variable extracted at each time step from a uniform distribution supported in $[0,1]$ so that, at each time step, we have the same probability of iterating $f_1(x)$ or $f_2(x)$. The support of the measure is a simple fractal, *i.e.* all the points $\zeta$ have the same non integer local dimension $d(\zeta) = d_H$, which also agrees with the value of all of Renyi's dimensions $d_{(q)}$'s, $\forall q \in \mathbb{N}$. Such dimension can be directly computed and its value is $d_1 = \log(2)/\log(3)$ (see *e.g.* [279]). We consider the usual observables $g_i$, $i = 1, 2, 3$, and we have that the conditions $Ð_0(u_n)$ and $Ð'_0(u_n)$ apply for this system. First of all, we have analyzed the empirical d.f. $F(u)$ of the extremes for $g_1$ observable. An example of the histogram and the corresponding fit to the GEV model for the observable $g_1$ is shown in Fig. 9.9. The histogram is obtained by iterating the map in Eq. 9.4.8 for $s = 5 \cdot 10^7$ iterations, at the point $\zeta \simeq 0.775$, by choosing $n = 5000, m = 1000$. The fit is produced with the MLE procedure. The empirical d.f. contains plateaux

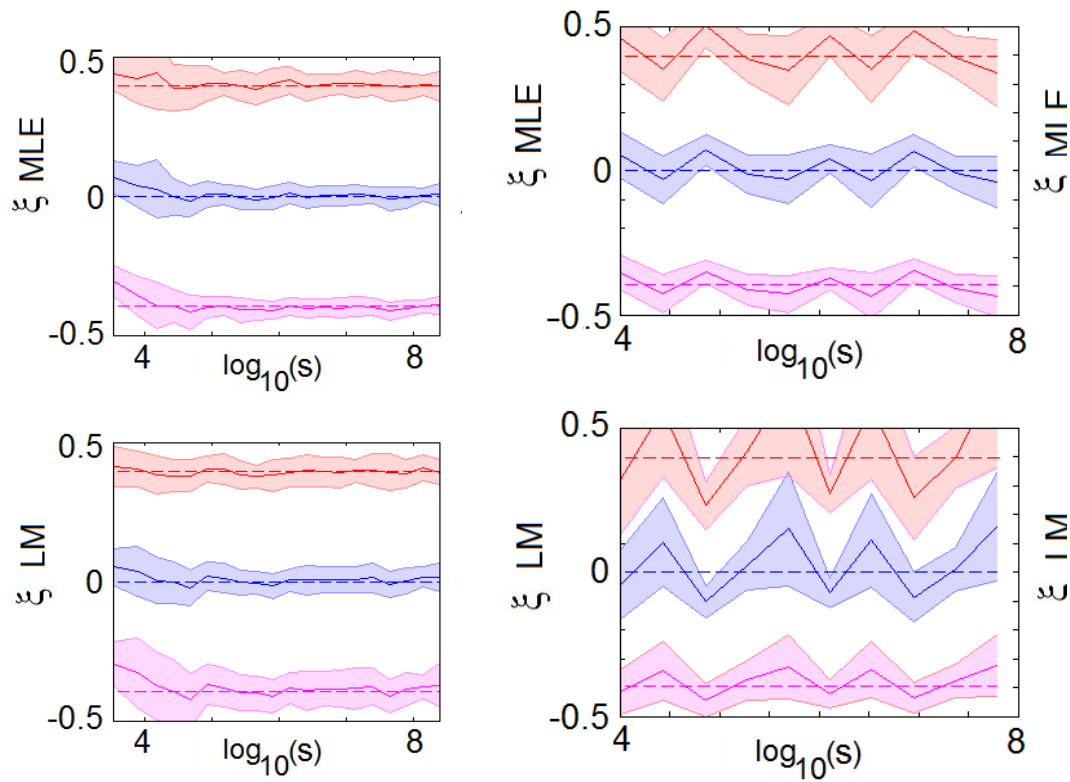

**Figure 9.10** IFS generating a Cantor set. Shape parameter $\xi$ of the GEV (left) and of the GPD (right) VS the total length of the series $s = nk$, $k = 10^3$ for the observables: $g_1$ (blue), $g_2$ (red), $g_3$ (magenta). Solid lines represent the mean values computed over 30 realizations, the shaded regions include one standard deviation of the mean. The dotted lines are the theoretical expected values. Top: MLE , bottom: L-moment.

which correspond to non accessible distances in correspondence of the holes of the Cantor set.

The numerical experiments on the IFS generating the Cantor set follow exactly the setting already described for the absolutely continuous case and they are reported in Fig. 9.10. The convergence towards the theoretical parameter is not as good both for the MLE and the L-moment procedure as in case where the underlying measure is absolutely continuous, compare with Fig. 9.7. Nonetheless, what we find is encouraging in terms of supporting the use of GEV methods also in the case where the invariant measure is singular. Apparently, there is not much difference in using the L-moment or the MLE for estimating $\xi$. However, we observed that many fits obtained via MLE returned unreliable uncertainties range as the minimization procedure failed. This is due to the fact that the MLE procedure works on continuous densities with a well defined maximum here not unequivocally detectable due to the jump occurring in the d.f. Therefore in these cases we recommend the use of the L-moment procedure which requires only the computation of integral expressions still well defined for singular continuous d.f. Further details on the BM estimates of $\xi$ for such a system are given in Table 9.1.

The experiments have been repeated with the POT approach and the related results are also presented in Fig.9.10. We have observed more serious issues in finding agreement between the numerical estimates and what the theory suggests. We find that in many cases the estimates of $\xi$ vary in a non-monotonic (actually, oscillating) manner as longer series of size $s$ are considered, while no convergence is found to the theoretical value even for very high values of $s$. This effect can be explained by observing that the discontinuities of the measure make the fit to the GPD model very sensible to the value of the density around the chosen threshold. The fit to this distribution is clearly less stable than with the BM method, because in this latter case the GEV density is constrained to go to zero at $\pm\infty$, so that the fitting procedure is less perturbed by the presence of holes in the measure.

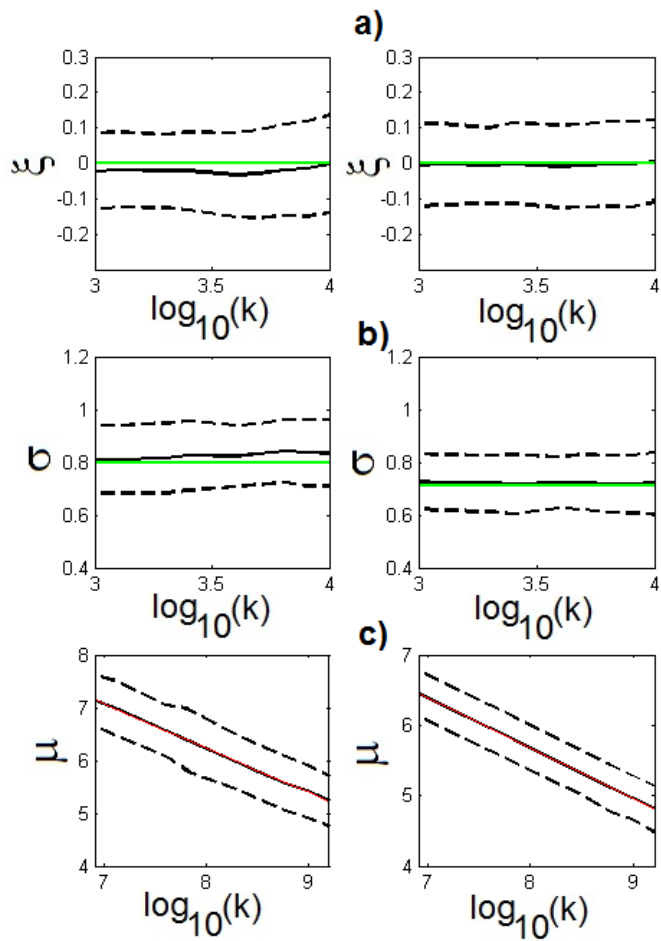

**Figure 9.11** $g_1$ observable. **a)** $\xi'$ VS $\log_{10}(k)$; **b)** $\sigma$ VS $\log_{10}(k)$; **c)** $\mu$ VS $\log(k)$. Left: Hénon map, Right: Lozi map. Dotted lines represent one standard deviation, red lines represent a linear fit, green lines are theoretical values.

#### 9.4.2.2 The Hénon and Lozi Maps

We want to study the extremes of distance observables for the Lozi and Hénon maps, which have already been introduced in Sect. 6.6.2 and Sect. 6.6.4, respectively. We briefly recapitulate here some of their properties. The Hénon map, already presented in Eq. 6.6.6, is defined as:

$$f_{a,b}(x, y) = (y + 1 - ax^2, bx) \tag{9.4.9}$$

while in the Lozi map, previously presented in Eq. 6.6.5, the term $ax^2$ is substituted with $a|x|$:

$$f_{a,b}(x, y) = (y + 1 - a|x|, bx) \tag{9.4.10}$$

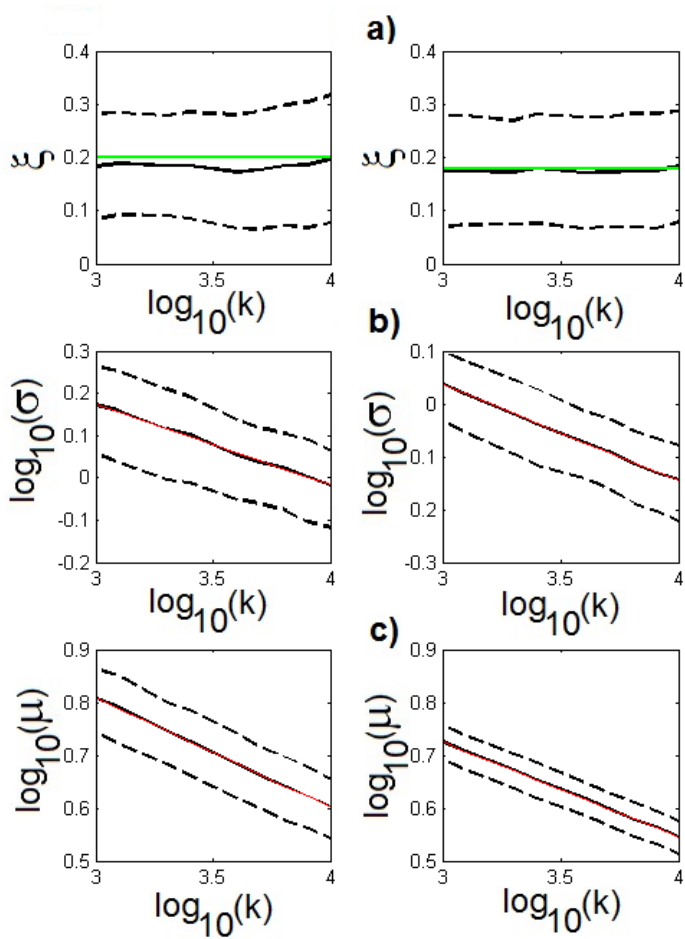

**Figure 9.12** $g_2$ observable **a)** $\xi'$ VS $\log_{10}(k)$; **b)** $\log_{10}(\sigma)$ VS $\log_{10}(k)$; fc) $\log_{10}(\mu)$ VS $\log_{10}(k)$. Left: Hénon map, Right: Lozi map. Dotted lines represent one standard deviation„ red lines represent a linear fit, green lines are theoretical values.

In the numerical experiments presented here, we consider the classical set of parameters $a = 1.4$, $b = 0.3$ for the Hénon map and $a = 1.7$ and $b = 0.5$ for the Lozi map. [280] proved the existence of the SRB measure for the Lozi map, whereas for the Hénon map no such rigorous proof exists, even if convincing numerical results suggest its existence [281]. Note that [196] proved the existence of an SRB measure for the Hénon map with a different set of parameters. Using the classical Young results [280] which make use of the Lyapunov exponents, we obtain an exact result for the information dimension $d_1 = \int \mathrm{d}\mu(\zeta)d(\zeta)$ for the Lozi attractor:

$$d_1 \simeq 1.40419$$

Instead, in the case of the Hénon attractor, we consider the numerical estimate provided by [282]:

$$d_1 = 1.25826 \pm 0.00006$$

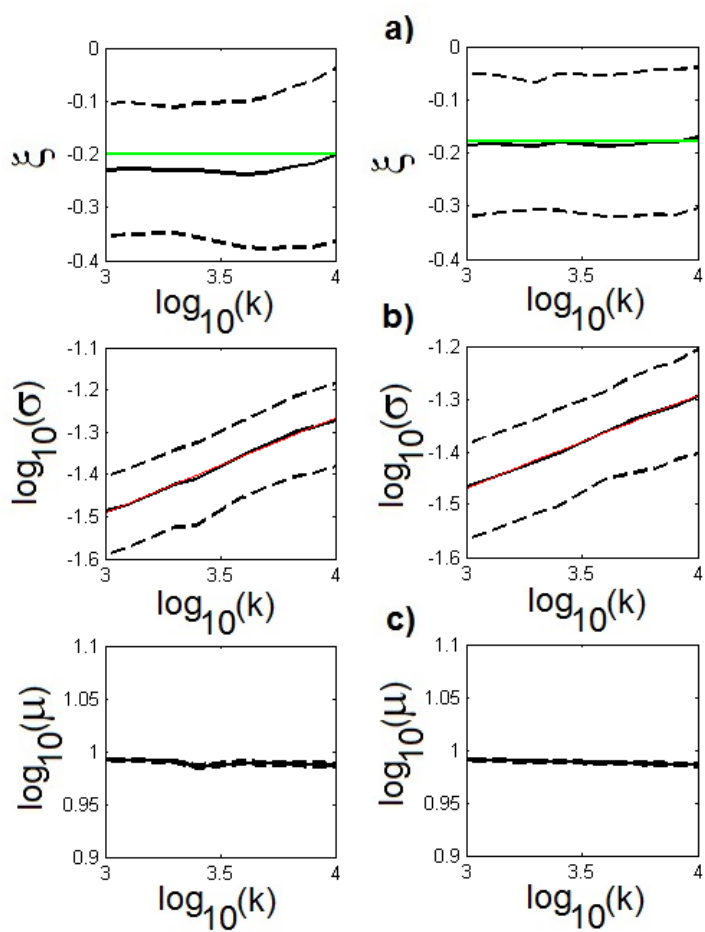

**Figure 9.13** $g_3$ observable. **a)** $\xi$ VS $\log_{10}(k)$; **b)** $\log_{10}(\sigma)$ VS $\log_{10}(k)$; **c)** $\log_{10}(\mu)$ VS $\log_{10}(k)$. Left: Hénon map, Right: Lozi map. Dotted lines represent one standard deviation, red lines represent a linear fit, and green lines are theoretical values.

| $d_1$ | **Baker** | **Hénon** | **Lozi** | **Cantor** |
|---|---|---|---|---|
| Theor. | 1.4357 | 1.2582 | 1.4042 | $\log(2)/|log(3) \sim 0.6309$ |
| $\sigma(g_1)$ | $1.43 \pm 0.03$ | $1.21 \pm 0.02$ | $1.39 \pm 0.02$ | $0.635 \pm 0.005$ |
| $\mu(g_1)$ | $1.48 \pm 0.03$ | $1.23 \pm 0.02$ | $1.40 \pm 0.01$ | $0.64 \pm 0.01$ |
| $\xi(g_1)$ | $1.41 \pm 0.02$ | $1.24 \pm 0.02$ | $1.41 \pm 0.01$ | $0.63 \pm 0.01$ |
| $\sigma(g_2)$ | $1.39 \pm 0.04$ | $1.35 \pm 0.07$ | $1.38 \pm 0.02$ | $0.63 \pm 0.01$ |
| $\mu(g_2)$ | $1.47 \pm 0.02$ | $1.24 \pm 0.01$ | $1.40 \pm 0.01$ | $0.64 \pm 0.01$ |
| $\xi(g_3)$ | $1.45 \pm 0.02$ | $1.28 \pm 0.02$ | $1.4 \pm 0.01$ | $0.64 \pm 0.01$ |
| $\sigma(g_3)$ | $1.56 \pm 0.08$ | $1.15 \pm 0.07$ | $1.42 \pm 0.01$ | $0.64 \pm 0.01$ |

**Table 9.1** Estimate of the information dimension $d_1$ obtained by averaging over $p = 1000$ ensemble members the estimates of $d(\zeta)$ computed by taking the logarithm of Eqs. 9.4.5-9.4.7 and computing the angular coefficient $\xi$ of a linear fit of data; for Baker, Hénon, and Lozi, maps, and the IFS generating the Cantor set. The Baker map is not discussed in this book; see [77] for details.

Following the considerations given in the previous section, the fit of the BM to the GEV distribution is performed using the L-moment methods. We average our results over $q = 1000$ chosen $\zeta$ reference points chosen on the attractor of the system according to the invariant measure. The simulation includes $s = nk = 10^7$ time steps and perform our statistical analysis considering different values of $n$ and $k$. We label each statistical analysis with the corresponding $\zeta$ and compute $d_1 = \int \mathrm{d}\mu(\zeta)d(\zeta)$ as $d_1 \sim 1/p \sum_{j=1}^{q} d(\zeta_j)$.

Our estimates for the parameters $\xi$, $\sigma$, and $\mu$ for the three observables $g_i$, $i = 1, 2, 3$ are presented in Figs. 9.11-9.13. The plots on the left-hand side refer to the Hénon map, while on the right-hand side the results refer to the Lozi map. When considering $\xi$, the numerical results are in agreement with the theoretical estimates. Nevertheless, the parameters distribution have a rather range spread which indicates a slower convergence towards the expected values in respect to what is observed for the IFS case. The estimates for $\sigma$ and $\mu$ obey the predictions given in Eqs- 9.4.5-9.4.7, even if substantial uncertainties persist.

We note that there are several ways to derive $d(\zeta)$ from the estimates of the GEV parameters using Eqs. 9.4.5-9.4.7. We focus on the expressions of $\sigma$ and $\mu$ for the $g_1$ observable, $\xi$, $\mu$ and $\sigma$ for the $g_2$ observable, the expression of $\xi$ and $\sigma$ for the $g_3$ observable. We derive for the estimates of such parameters the corresponding estimates of $d(\zeta)$, and then average over the set of $q = 1000$ chosen $\zeta$'s.

Results are presented in Table 9.1, where we present for each considered observable and for each map the best estimate for $d_1$ and its uncertainty as measured by twice the standard deviation. We also include, for reference the results obtained for the IFS generating the Cantor set discussed before and for the Baker map [69]. A more thorough discussion of these findings is provided in [77]. We have satisfactory agreement with the theory. Nonetheless, substantial uncertainty persists: the relatively slow convergence for these maps may be related to the difficulties experienced computing the dimension with all box-counting methods, as shown in [281, 282]. In that case it has been proved that the number of points that are required to cover a fixed fraction of the support of the attractor diverges faster than the number of boxes itself for this kind of non uniform attractor. In our case the situation is similar since we consider balls around the initial condition $\zeta$.

The best result for the dimension is achieved using the parameters provided by $g_1$ observable, since the logarithm modulation of the distance takes into proper account real extrema while weighting less extremely large extreme events which, within a finite datasets context, might appear as outliers spoiling the statistics.

## 9.5
## Analysis of the Distance and Physical Observables for the Hénon map

In this section we want to bring together the analysis of extremes of distance and physical observables for the Hénon map given in Eq. 9.4.9, taking the point of view of POT, so that we will mostly refer to results contained in Chap. 8, even if results from Chap. 6 are also relevant.

We consider two sets of parameter values for which chaotic behavior is observed, $a = 1.4$ $b = 0.3$ and $a = 1.2$, $b = 0.3$. In the first case, the largest Lyapunov exponent $\lambda_1 \sim 0.416$ and the Kaplan-Yorke dimension is estimated as $d_{KY} = 1 + \lambda_1/|\lambda_2| = 1 + \lambda_1/|\log(b) - \lambda_1| \sim 1.26$, where $d_u = 1$ and $d_s = \lambda_1/|\log(b) - \lambda_1| \sim 0.26$. In the second case, the largest Lyapunov exponent $\lambda_1 \sim 0.305$ and the Kaplan-Yorke dimension is estimated as $d_{KY} = 1 + \lambda_1/|\lambda_2| = 1 + \lambda_1/|\log(b) - \lambda_1| \sim 1.20$.

Since the local dimension $d(\zeta)$ is not constant on the attractor (see discussion in Sect. 9.4.2.2 for the first choice of the parameters), we have that these systems are not exact dimensional. Therefore, the considered pairs of values of $a$ and $b$ do not belong to the Benedicks-Carleson set of parameters, which, instead, lead to an invariant SRB measure for the system. Moreover, this sheds doubts on the validity on the formula given Eq. 6.11.3 (obtained using the BM point of view) and its high-dimensional version (obtained using a POT-based approach) given in Eqs. 8.2.15 and 8.2.17 for estimating the shape parameter for the family of physical observables given in Eq. 6.11.2 by setting $\theta = 0$.

We proceed as follows for both pairs of parameters $a = 1.4$ $b = 0.3$ and $a = 1.2$, $b = 0.3$. The initial conditions are selected in the basin of attraction of the strange attractors. We perform long integrations ($s = 10^{10}$ iterations) and select the maximum value of $A(x, y) = A(\vec{x}) = x$, which can be obtained from the family of observables given in Eq. 6.11.2 by setting $\theta = 0$. Note that in this section we use he bold font $\vec{x}$ to refer to the point $(x, y)$. We denote such maximum as $A_{max}$, and define as $\vec{x_0}$ the unique point belonging to the attractor such that $A(\vec{x_0}) = A_{max}$. We then construct the observable $B(\vec{x}) = -dist(\vec{x}, \vec{x_0})$, which measures the distance between the orbit and the point $\vec{x_0}$. As discussed in Chapter 8 and detailed in [78], the asymptotic properties of the extremes (maxima) of the $B$ observable allow to derive easily the local dimension $d(\vec{x_0})$. We then repeat the investigation using, instead, the observable $A(\vec{x}) = -x$. In all the analyses presented below, we have chosen extremely high thresholds $T$ for studying the statistical properties of the extremes of the $A$ and $B$ observables, in such a way to include only about a fraction of about $10^{-5}$ or less of the total number of points of the orbit. All the results are insensitive to choice of $T$, which suggests that we are well into the asymptotic regime.

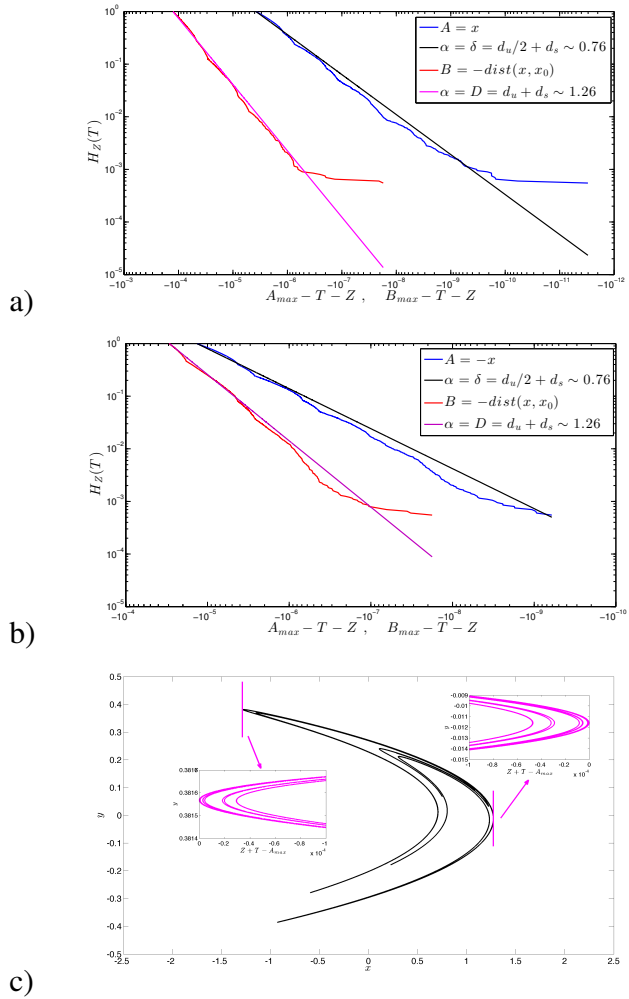

**Figure 9.14** Results of numerical simulations performed on thnon mape Hénon map with parameters' value $a = 1.4$ and $b = 0.3$. a) Blue curve: empirical $\bar{F}_T(Z)$ for the observable $A = x$, with $A_{max} = A(\vec{x_0}) \sim 1.2730$. Black line: power law behavior deduced from the theory. Red curve: empirical $\bar{F}_T(Z)$ for the observable $B = -dist(\vec{x}, \vec{x_0})$, with $B_{max} = 0$. Magenta line: power law deduced from the theory. b) Same as a), for the observable $A = -x$, with $A_{max} \sim 1.2847$ and $B_{max} = 0$. c) Approximation to the attractor with blow-ups of the portions of the invariant measure corresponding to the extremes of the $A$ observables ($\Omega^T_{A_{max}}$ regions); the vertical lines indicate the thresholds. In both inserts, we consider $A_{max} - T = 10^{-4}$. See also Fig. 8.2.

The results obtained for the Hénon system featuring $a = 1.4$ $b = 0.3$ are shown in Fig. 9.14, where we present the complementary cumulative distribution of excesses $\bar{F}_T(Z)$ (see Eq. 8.2.16) for $A(\vec{x}) = x$ ($A(\vec{x}) = -x$) and for the corresponding $B(\vec{x}) = -dist(\vec{x}, \vec{x_0})$ in panel a) (panel b)). The empirical values of $\bar{F}_T(Z)$ for the $A$ and $B$ observables are shown by the blue and red curves, respectively, and the power law behavior $\bar{F}_T(Z) = (1 - Z/(A_{max} - T))^{\alpha}$ given by the theory (assuming Axiom A properties) are shown by the black and magenta lines, respectively.

The error bars on the empirical $\bar{F}_T(Z)$ (estimated by varying the initial conditions of the simulation) are for almost all values of $Z$ so small that they cannot be graphically reproduced. Instead, the flat region obtained for very low values of $\bar{F}_T(Z)$ results from the finiteness of the sampling and gives the baseline uncertainty. Note that the straight lines are obtained out of the theoretical predictions, without any procedure of optimization or of fit, so that no uncertainties are involved. The empirical and theoretical distributions obey the same normalization.

We first observe that the local dimension in the vicinity of both $\vec{x}_0$'s is extremely close to the $d_{KY} \sim 1.26$, as $\bar{F}_T(Z)$ scales to a very good approximation with an exponent $\alpha \sim d_{KY}$; compare the red curves and the magenta lines. Note that, considering that the local dimension has rather large variations across the attractor of the Hénon system, such a correspondence was not intentionally pursued. However, these are favorable circumstances to check the theory. We find that the distributions $\bar{F}_T(Z)$ for the observables $A(\vec{x}) = x$ and $A(\vec{x}) = -x$ also obey accurately the power law scaling with exponent $\alpha \sim \delta = d_u/2 + d_s \sim 0.76$ given in Eq. 6.11.3 and Eqs. 8.2.15 and 8.2.17, compare the blue curves and the black lines. In panel c) we present a simple description of the geometry of the problem, by showing an approximation to the map's attractor with blow-ups of the portions of the invariant measure corresponding to the extremes of the $A$ observables (the regions $\Omega_{max}^T$ introduced in Fig. 8.2). Even if the geometrical properties of the regions of the attractor around the two $x_0$'s seem indeed different, when zooming in, the two $\Omega_{max}^T$ regions look similar. The presence of many parabolas-like smooth curves stacked according to what looks qualitatively like a Cantor set fits with the comments and calculations given in Chaps. 6 and 8.

In Fig. 9.15 we report the corresponding results obtained for the Hénon system featuring $a = 1.2\, b = 0.3$. By looking at the empirical $\bar{F}_T(Z)$ of the $B$ observables, we note that also in this case the local dimension is close to the value of $d_{KY} \sim 1.20$ for both extremal points $\vec{x_0}$'s (compare the red curves and the magenta lines in panels a and b). The agreement between the predicted value of the power law scaling for the $\bar{F}_T(Z)$ of the $A$ observables is not as good as in the case reported in Fig. 9.14. The predicted scaling exponent $\delta = d_u/2 + d_s \sim 0.70$ seems to overestimate the very large extremes. Nonetheless, a power law scaling is apparent for the empirical $\bar{F}_T(Z)$. Note that the bias between the theoretical and empirical scalings is of the same sign for both the $A$ and $B$ observables, suggesting that also for the $A$ observables part of the disagreement is due to the discrepancy between the local dimension and the Kaplan-Yorke dimension (there is a shift in the values of the slopes). Also here, panel c) provides an approximate representation of the attractor of the system, and, in particular of the $\Omega_{max}^T$ regions: by comparing it with panel c) of Fig. 9.15, and considering that they contain the same number of points, one can intuitively grasp that the local dimension is lower in this case.

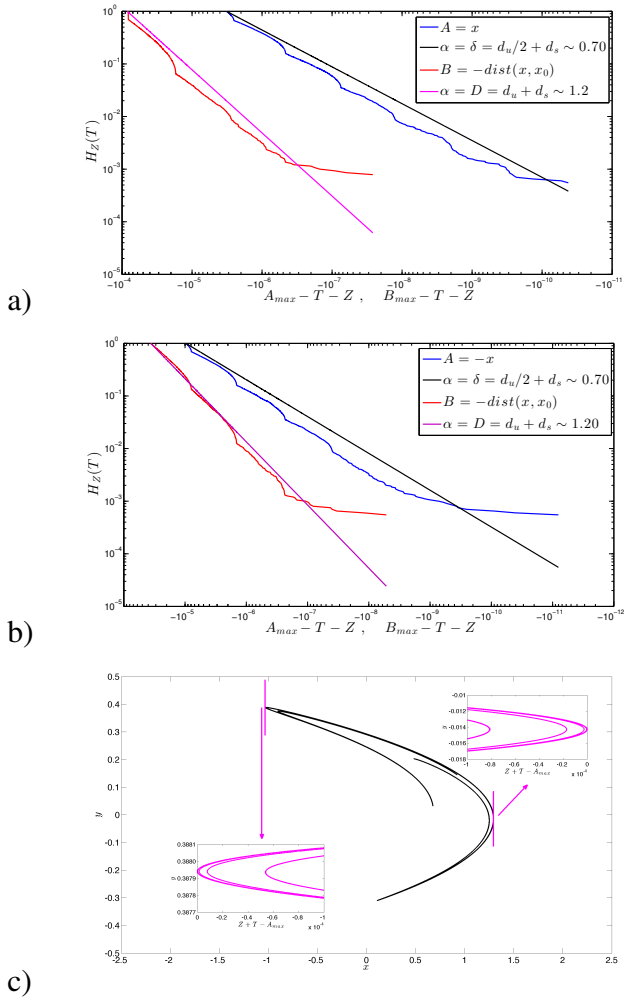

**Figure 9.15** Same as in Fig. 9.14, but for parameters' value $a = 1.2$ and $b = 0.3$. In this case in a) $A_{max} \sim 1.2950$ and $B_{max} = 0$, and in b) $A_{max} \sim 1.0328$ and $B_{max} = 0$.

## 9.5.1
## Remarks

We would like to emphasize that in panels a) and b) for Figs. 9.14 and 9.15, we observe deviations of the empirically obtained $\bar{F}_T(Z)$ from the power law behaviour, in the form of fluctuations above and below a straight line in a $\log$-$\log$ plot (this is quite clear in Fig. 9.15). As discussed in Sect. 9.4.2, the presence of such modulations across scales result from the fact that gaps are present along the stable manifold containing $\vec{x_0}$, with a Cantor set-like structure. See the inserts in Figs. 9.14c) and 9.15c), where the stable manifold (not shown) is, as opposed to the unstable manifold, not orthogonal to the gradient of $A$ (the $x$ direction, in this case).

Therefore, when we integrate the density of states along the direction of the gradient of the $A$ observable starting in $\vec{x_0}$ in order to obtain $\mu(\Omega_{A_{max}}^{T+Z})$ and $\mu(\Omega_{A_{max}}^{T})$ (see Fig. 8.2), we get a factor $(A_{max} - T - Z)^{d_u/2}$ ($d_u = 1$) coming from the (local) paraboloidal form of the unstable manifold discussed in Sec. 8.2.3, times a devil's staircase which can, *on the average*, be approximated by the power law $(A_{max} - T - Z)^{d_s}$. The same geometric arguments apply when considering integrations along the spherical shells centred on $\vec{x}_0$ for constructing the extreme value laws for $B$ observables. See also the discussion in the context of the GEV approximation to the BM statistics related to Fig. 9.9.

Overall, such results suggest that it is indeed promising to use the combined statistical properties of the extremes of physical and distance observables for determining the geometry of the attractor in terms of its partial dimensions along the stable and the unstable manifold.

## 9.6
## Extremes as Dynamical Indicators

In the previous sections we have shown how to use extremes for deducing important properties of the geometry of the attractor supporting the invariant measure. In other terms extremes act as a sort of microscopes allowing us to focus on the fine structure of the dynamics of the system. In the last example, we have shown how considering distance and physical observables at the same time makes it possible to gather specific information on the directions of contraction and of expansion in the attractor. In this section we would like to present some results which can clarify how observing extremes we can learn about the qualitative properties of the dynamics of the underlying system, and, in particular, detecting changeovers between regular and chaotic behaviours as a function of a control parameter.

In the context of dynamical systems theory, a large number of tools known as *indicators of stability* have been developed for determining whether a system is regular or chaotic. Quantities like Lyapunov exponents [283, 70, 284, 285] and the indicators related to the RTS [286, 287, 288, 289] have been used for a long time for such a task. Nevertheless, in the recent past, the need for computing stability properties with faster algorithms and for systems with many degrees of freedom resulted in a renewed interest in the technique and different dynamical indicators have been introduced. The Smaller Alignment Index (SALI) described in [290] and [291], the Generalized Alignment Index (GALI), introduced in [292] and the Mean Exponential Growth factor of Nearby Orbits (MEGNO) discussed in [293], [294] are suitable to analyze the properties of a single orbit. They are based on the divergence of nearby trajectories and require in principle the knowledge of the exact dynamics. Another class of indicators is based on the properties of the error due to the numerical round off and has been discussed in [244]. In this case, the focus is on illustrating the dynamical properties of a system by computing the divergence between two trajectories where we choose the same initial condition but different numerical precision in the numerical integration. The so called *Reversibility Error*, which measures the distance between a certain initial condition and the end point of a trajectory iterated forward and backward for the same number of time steps, gives basically the same information.

The purpose of all these indicators is to provide a *global* information on the struc-ture of the physical measure of the system. Few of them are explicitly designed to sample local properties. Moreover, even when this is the case, the level of *zoom* on the physical measure cannot be manually set or changed. The purpose of this chapter is to describe the geometrical and dynamical indicators derived by the extreme value analysis of recurrences. These indicators are naturally intended to provide a local information around the point chosen to sample recurrences. They will therefore act as a magnifying glass on the physical measure providing different levels of detail according to the specific problems to address.

### 9.6.1
### The Standard Map: Peaks Over Threshold Analysis

The standard map [295] is an area-preserving chaotic map defined on the bidimensional torus, and it is one of the most widely-studied examples of dynamical chaos in physics. The corresponding mechanical system is usually called a kicked rotator. The maps is defined as follows:

$$f(x,y) = \left( y - \frac{K}{2\pi} \sin(2\pi x), x + y + 1 \right) \mod \mathbb{T}^2 \tag{9.6.1}$$

The dynamics of the map given in Eq. (9.6.1) can be regular or chaotic. For $K \ll 1$ the motion follows quasi periodic orbits for all initial conditions, whereas if $K \gg 1$ the motion turns to be chaotic and irregular. An interesting behavior is achieved when $K \sim 1$: in this case we have coexistence of regular and chaotic motions depending on the chosen initial conditions [296].

We perform for various values of $K$ ranging from $K = 10^{-4}$ up to $K = 10^2$ an ensemble of 200 simulations, each characterised by a different initial condition $\zeta$ randomly taken on the bidimensional torus, and we compute for each orbit the observables $g_i$, $i = 1, 2, 3$ discussed above. In particular, without loss of generality, we choose the useful setting $g_i = g_i(dist(f^t\zeta, \zeta)), i = 1, 2, 3$, so that we study the recurrence properties of the orbits in the neighbourhood of the initial condition.

In each case, the map is iterated until obtaining a statistics consisting $10^4$ exceedances, where the threshold $T = 7 \cdot 10^{-3}$ and $\alpha = 3$. One can check that all the results are indeed robust with respect to the choice of the threshold and of the value of $\alpha$. For each orbit, we fit the statistics of the $10^4$ exceedances values of the observables to a GPD distribution, using a MLE [297] implemented in the MATLAB© function *gpdfit* [298].

The results are shown in Fig. 8.2 for the inferred values of $\xi$ and $\sigma$ and should be compared with Eqs. (8.2.10)-(8.2.12). When $K \ll 1$, we obtain that the estimates of $\xi$ and $\sigma$ are compatible with a dimension $d(\zeta) = 1$ for all the initial conditions: we have that the ensemble spread is negligible. Similarly, for $K \gg 1$, the estimates for $\xi$ and $\sigma$ agree remarkably well with having a local dimension $d(\zeta) = 2$ for all the initial conditions. In the transition regime, which occurs for $K \simeq 1$, the ensemble spread is much higher, because the scaling properties of the measure are different among the various initial conditions. As expected, the ensemble averages of the parameters change monotonically from the value pertaining to the regular regime to

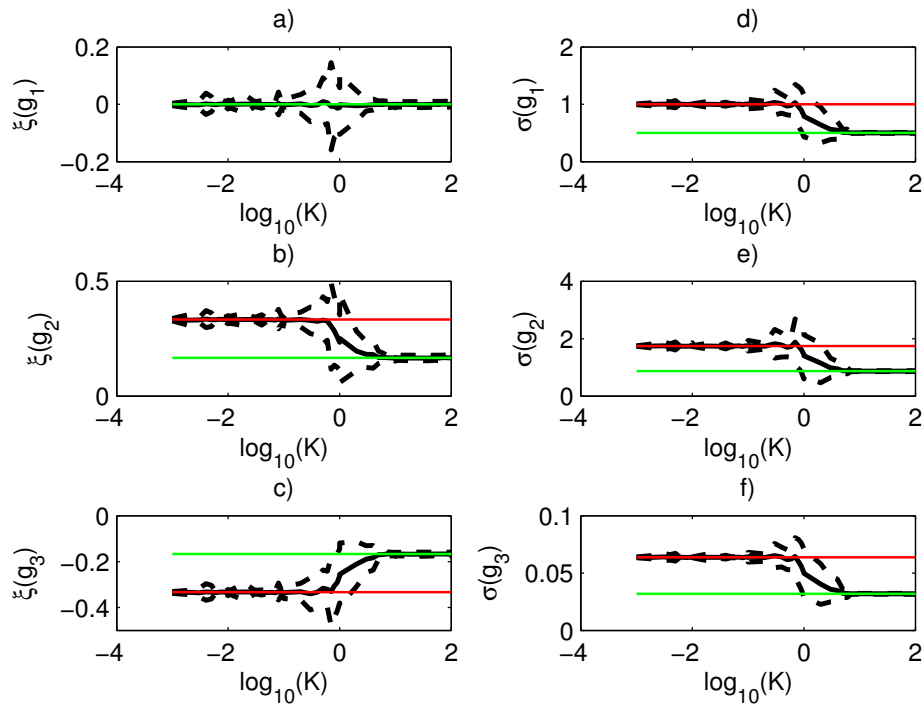

**Figure 9.16** GPD parameters for an ensemble of 200 initial conditions, Standard map, different values of $K$. **a)** $g(\xi_1)$ VS $K$, **b)** $g(\xi_2)$ VS $K$, **c)** $g(\xi_3)$ VS $K$, **d)** $g(\sigma_1)$ VS $K$, **e)** $g(\sigma_2)$ VS $K$, **f)** $g(\sigma_3)$ VS $K$. Black solid lines: averaged value. Black dotted lines: uncertainty evaluated as one standard deviation of the ensemble. Green lines: theoretical expected value for regular orbits. Red lines: theoretical expected value for chaotic orbits

that pertaining to the chaotic regime with increasing values of $K$.

Basically, this measures the fact that the so-called regular islands shrink with $K$. Note that in the case of the observable $g_1$, the estimate of the $\xi$ is robust in all regimes, even if, as expected, in the transition between low and high values of $K$ the ensemble spread is larger.

### 9.6.2
### The Standard Map: Block Maxima Analysis

As it should be by now clear to the reader, it is especially interesting to compare the results of analysing extremes of observables generated by dynamical systems via the BM and the POT methods when the applicability of the theory developed in this book is borderline. We have above shown that the transition from regular to chaotic dynamics of the standard map as a function of $K$ is well captured by looking at extremes using the POT approach. For both very large and very low values of $K$ the POT method works well, because its applicability is oblivious to whether the dynamics is chaotic or regular [78]. Instead, problems emerge in the transitional range $K \sim 1$.

Instead, the BM approach works only if the dynamics is sufficiently mixing. Therefore, we expect that a different picture of the transition between chaotic and regular motion appears when trying to construct GEV models for extremes across a vast range of values of $K$. We then take an ensemble of 500 initial conditions centered around $(x_0, y_0) = (0.305, 0.7340)$ in a small subset of the two-dimensional torus. Each orbit comprises of $s = 10^6$ iterations, so that we select $k = 1000$ maxima over bins of length $n = 1000$ of the observables $g_i = g_i(dist(f^t\zeta, \zeta))$, $i = 1, 2, 3$ as above. We then compute the best fit of the GEV parameters for each realization, and then averaged them over the different initial conditions. We consider a range of $K$ spanning from $K = 10^{-4}$ up to $K = 10^2$.

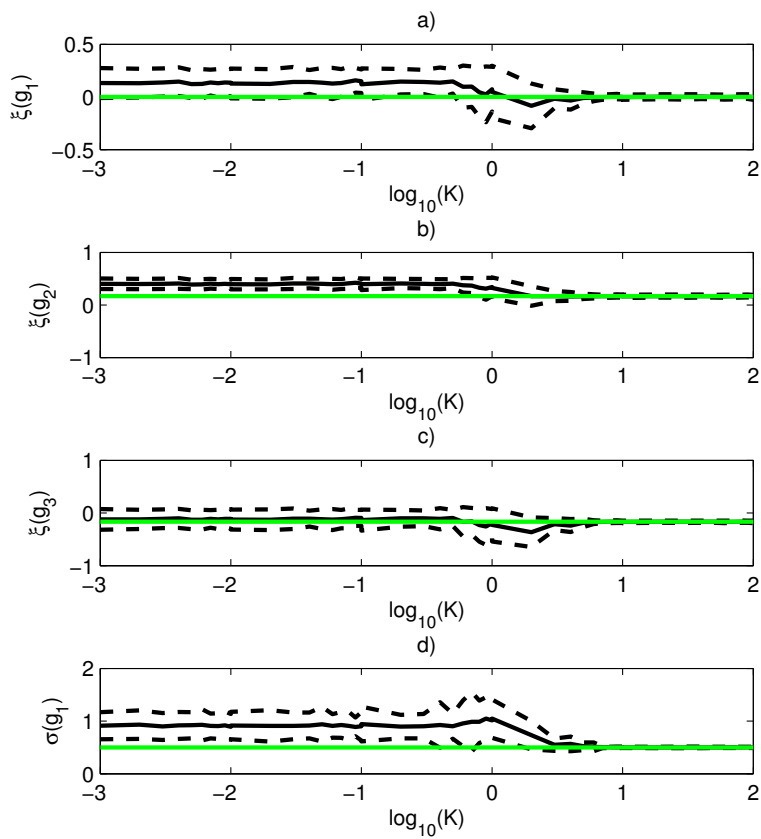

**Figure 9.17** Standard map: GEV parameters averaged over 500 different initial condition centred on $x_0 = 0.305, y_0 = 0.7340$ VS $K$. a) $\xi(g_1)$, b) $\xi(g_2)$, c) $\xi(g_3)$, d) $\sigma(g_1)$.

In order to study the changeover from $K \ll 1$ to $K \gg 1$, we choose as indicators the shape parameters for the three type observables $\xi(g_1), \xi(g_2), \xi(g_3)$ and the scale parameter for the type 1 observable $\sigma(g_1)$. Looking at Eqs. 9.3.1-9.3.3, these seem the best suited because they are related to the local dimension of the attractor $d(\zeta)$ but do not have a dependence of $m$, therefore, once we are in the asymptotic regime, the results are independent on the number of observations in each bin. The results are presented in Fig. 9.17. In this example we set $\alpha = 3$ for the $g_2$ and $g_3$ observables. For each parameter, the averaged value is represented with a solid line whereas the dotted lines represent one standard deviations of the ensemble. It is clear that for $K \gg 1$ the parameters converge towards the theoretical values predicted by the theory (we are, in fact, in a regime of mixing dynamics), whereas for $K \ll 1$, where we have regular motions, the GEV fits simply fail, as can be assessed by a Kolmogorov-Smirnov test [278]. As flag for this, we have that the fitted parameters have a very large spread, which is more than five times larger with respect to the case of the chaotic counterpart. The results are virtually unchanged if we change the initial conditions and the value of $\alpha$. The next step is to investigate extensively on the attractor the properties of the extremes of the distance observables.

### 9.6.2.1 Using Extremes to Separate Islands of Regular Dynamics from the Sea of Chaos

We want now to give evidence that a BM analysis of extremes of the observables $g_i = g_i(dist(f^t\zeta, \zeta))$, $i = 1, 2, 3$ can provide a great deal of information on the local (in the attractor) predictability of the system. We first introduce two standard methods used for this purpose, the analysis of the *divergence of orbits* due to the numerical round-off and the *reversibility error*. We briefly summarize here some definitions and suggest the reader to look into [244] for further clarifications.

**Divergence of the Orbits.** The arithmetic operations performed on a computer are unavoidably affected by round-off, which cause error to be propagated each time an operation is performed. Round-off algebraic procedures are hardware dependent, as detailed in [243]. Suppose we are given a map $f^t(x)$ then we will indicate with $f_*^t(x)$ the correspondent numerical map both evaluated at the $t$-th iteration. We define the divergence of orbits as:

$$\Delta_t(x) = \text{dist}(f_S^t(x), f_D^t(x)), \tag{9.6.2}$$

where $f_S^t$ and $f_D^t$ stand for single and double precision iterations, respectively, and dist is a suitable metrics.

**Reversibility Error.** If the map is invertible we can also define the reversibility error as

$$R_t(x) = \text{dist}(f_*^{-t} \circ f_*^t(x), x) \tag{9.6.3}$$

which is nonvanishing since the numerical inverse $f_*^{-1}$ of the map is not exactly the inverse of $f_*$ namely $f_*^{-1} \circ f_*(x) \neq x$. The reversibility error is easier to compute than the divergence of orbits (if we know explicitly the inverse map) and provide a comparable information. Both quantities grow on the average linearly if $f$ is a regular map together, while they grow exponentially if $f$ has a positive Lyapunov exponent. When computing $R_t(x)$ we set $f_* = f_S$ in order to be able to compare it with $\Delta_t(x)$.

Let us now choose the value of $K = 6.5$, which corresponds to a regime where we begin to have a good agreement between the inferred GEV parameters and what predicted by the theory in the presence of suitable mixing conditions, see Fig. 9.17. We want to show that we can capture the main properties of the standard map with the indicators presented above, by considering $500 \times 500 = 2.5 \times 10^5$ initial conditions $\zeta$ uniformly distributed on the two-dimensional torus. We stick to the previously mentioned values for the number of iterations $s = 10^6$, the number of bins $k = 10^3$, and the length of the bins $n = 10^3$ as well as to the value of $\alpha = 3$ for the observables $g_2$ and $g_3$.

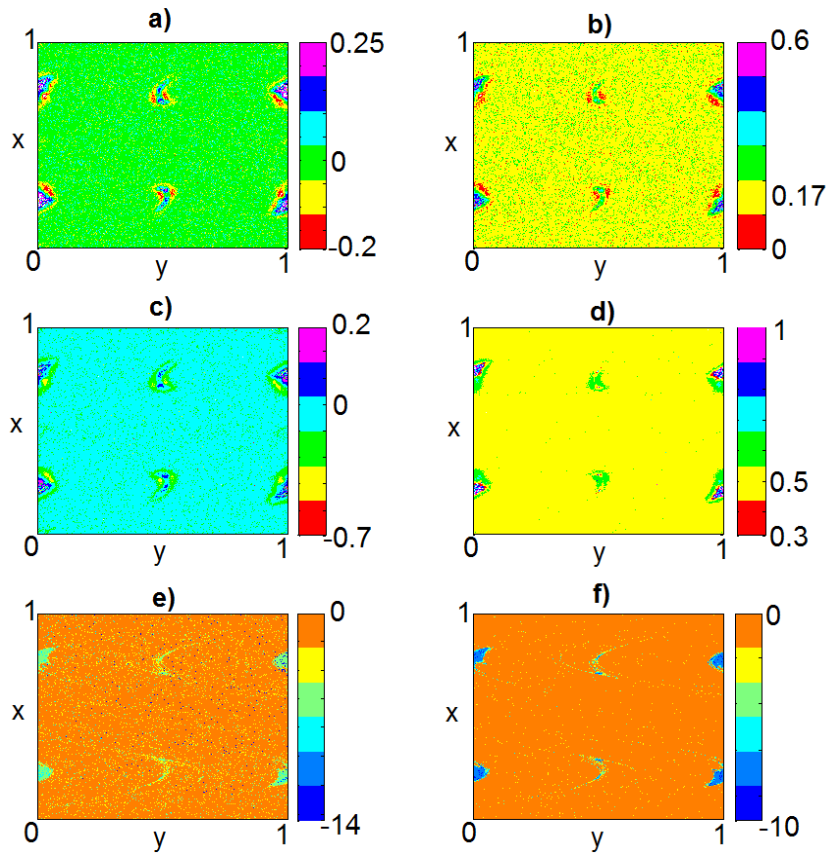

**Figure 9.18** Predictability of the standard map ($K = 6.5$) for $2.5 \times 10^5$ uniformly distributed initial conditions. Each simulation runs for $s = 10^6$ time steps. **a)** $\xi(g_1)$, **b)** $\xi(g_2)$, **c)** $\xi(g_3)$, **d)** $\sigma(g_1)$: We consider $k = 10^3$ maxima, each taken each over bins of length $n = 10^3$. **e)** $\log_{10}(R_{t=100})$: Reversibility error computed over 100 time steps. **f)** $\log_{10}(\Delta_{t=100})$: Divergence of trajectories computed over 100 time steps.

Results are shown in Fig. 9.18 where we present for each $\zeta$ the four parameters of GEV distribution (top and middle panels) for the extremes of $g_i(f^t \zeta, \zeta)$, $i = 1, 2, 3$, alongside the reversibility error and divergence of orbits in logarithm scale (lower panel). The number of iterations for the round off indicators is $t = 100$. The latter indicators provide a clear-cut separation between the so-called (small) *islands of regularity* and the so-called (widespread) *sea of chaos*, where they have very low and very high values, respectively. The attractor of the standard map is extremely non-homogeneous: in the islands of regularity, the system is highly predictable, as opposed to the the sea of chaos. Of course, if we consider higher (lower) values of $K$, the portion of the attractor occupied by the islands of regularity shrinks (increases).

Such basic dynamical structure is well highlighted by all the indicators based on GEV distribution. In all cases, we have that the indicators give values compatible with a mixing dynamics taking place on an attractor with local dimension $d(\zeta) = 2$ *only* in the chaotic regions. Instead, in the small regular islands we observe significant deviations from the expected values. Additional details can be found in [76].

## 9.7
## Extreme Value Laws for Stochastically Perturbed Systems

In Chapter 7 has been devoted to developing a mathematical framework for studying EVLs in randomly perturbed systems, focusing on two types of stochastic perturbation: additive noise and observational noise. The problem of understanding the impact of adding noise in a system or of observing with finite accuracy its state on the properties of its extremes is challenging and of general interest. In particular, the theory we have developed has practical relevance in a wide range of applications such as the analysis of the role of truncation errors for instrument with low accuracy, the statistics of points visited sporadically in the analysis of recurrence of time series, and the possibility of computing attractor dimension by using Eq. 7.5.4, which may serve as an alternative with respect to other well established techniques. See discussion in [65].

Before delving into the analysis of the extremes of few specific numerical models, let us stress an important point. Adding (sufficiently fast decorrelating) noise of amplitude $\varepsilon$ in a system either in the dynamics or in the process of observing its state tends to make it obey the decorrelation conditions $Ɗ_0$, $Ɗ_0'$ and their variants, thus allowing for an easier use of the probabilistic framework of EVT. The other side of the coin is that the presence of noise tends to mask the fine structure of the underlying deterministic dynamics below a certain scale defined by $\varepsilon$.

We can find a clear example of these effects by close inspection to the proof of Proposition 7.5.1, where the parameter $\varepsilon$ appears in the denominator of one factor in the r.h.s. of the term (II) at the end of the proof. This means that the convergence gets better when $\varepsilon$ is large, which is not surprising since a large value of the perturbation implies a more stochastic independence of the process. On the other hand, if we want to use the form of the linear scaling parameter $b_n$ to catch the local properties of the invariant measure $\mu$, we need small values of $\epsilon$. Hence, A careful compromise on the value of $\varepsilon$ between these the two regime of ephemeral and very strong noise is necessary when we pursue effective numerical investigations.

### 9.7.1

### Additive Noise

We will analyse the effects of the additive noise for two one dimensional maps: the rotation map introduced in Sec. 5.6 and defined as $f(x) = x + \alpha \mod 1$, $\alpha \in \mathbb{R}$ and the chaotic Bernoulli shift map $f(x) = qx \mod 1$, $q \in \mathbb{N}$, $q \geq 2$, introduced in Sect. 4.2. These toy models are well know to have relevance for studying also higher dimensional dynamical systems.

The rotation map has a regular dynamics and the statistics of its extremes does *not* conform to the EVLs, as discussed in Sect. 5.6. Instead, in Sect. 7.2.1 we have shown that EVLs are found as correct asymptotic model of the extremes of distance observables for the same map if we add random perturbations in the form of additive noise. We show below that this effect is practically detectable only if the intensity of the noise $\varepsilon$ is high enough.

The Bernoulli shift map features chaotic dynamics and the presence of a random perturbation has no effect on the convergence properties of the extremes of the distance observables to the EVLs nor in the estimate of the corresponding parameters. The only exception comes from the study of the recurrence of periodic points: in the deterministic case, the extremes obey modified EVLs with $EI < 1$ (see Sect. 4.2.2)[1], while the so-called dichotomy between periodic points and the rest of the attractor is washed out as soon as noise in added into the system, see Sect. 7.3.

#### 9.7.1.1 Rotations

We first discuss the properties of the extremes of distance observables for the stochastically perturbed rotation map $f_{\varepsilon\omega} = x + \alpha + \varepsilon\omega \mod 1$, where $\omega$ is a random variable uniformly distributed over the interval $[-1, 1]$. The results are displayed in Fig. 9.19 where the green lines correspond to experiments where we have chosen a bin length $n = 10^4$, whereas the blue lines refer to experiments where the chosen bin length is $n = 10^3$. The red lines indicate the values of the parameters predicted by the theory for a one-dimensional map satisfying the mixing conditions $Ɖ_0$ and $Ɖ'_0$ (which, again, are *not* obeyed by the unperturbed map). We set $\epsilon = 10^{-p}$ to analyse the role of the perturbations on scales ranging from values smaller than those typical for the numerical noise up to $\mathcal{O}(1)$.

The solid lines display the values obtained by averaging over the 500 realisations of the stochastic process, while the error bars indicate the standard deviation of the sample. Finally, with the dotted lines we indicated the experiments where less than 70% of the 500 realisations produce a statistics of extremes such that the empirical d.f. passes successfully the non-parametric Kolmogorov-Smirnov test [278] when compared to the best GEV fit.

Even though the mathematical findings presented in Sect. 7.2.1 guarantee the existence of EVL for the rotations perturbed with an arbitrarily weak noise (in practical terms, *e.g.*, comparable with the round-off), the simulations clearly show that EVLs are obtained when considering small but finite noise amplitudes only when very long trajectories are considered. The quality of the fit improves when larger bins are considered (compare blue and green lines in Fig.9.19). This is in agreement with the idea that we should get EVL for infinitely small noises in the limit of infinitely long samples. In our case, EVLs are obtained only for $\epsilon > 10^{-4}$, which is still considerably larger than the noise introduced by round-off resulting from double precision, as the round-off procedure is equivalent to the addition to the exact map of a random noise of order $10^{-7}$ [243, 244].

---

1) We recall that $\theta = \theta(z) = 1 - |\det D(f^{-p}(z)|$, where $z$ is a periodic point of prime period $p$, see [49, Theorem 3].

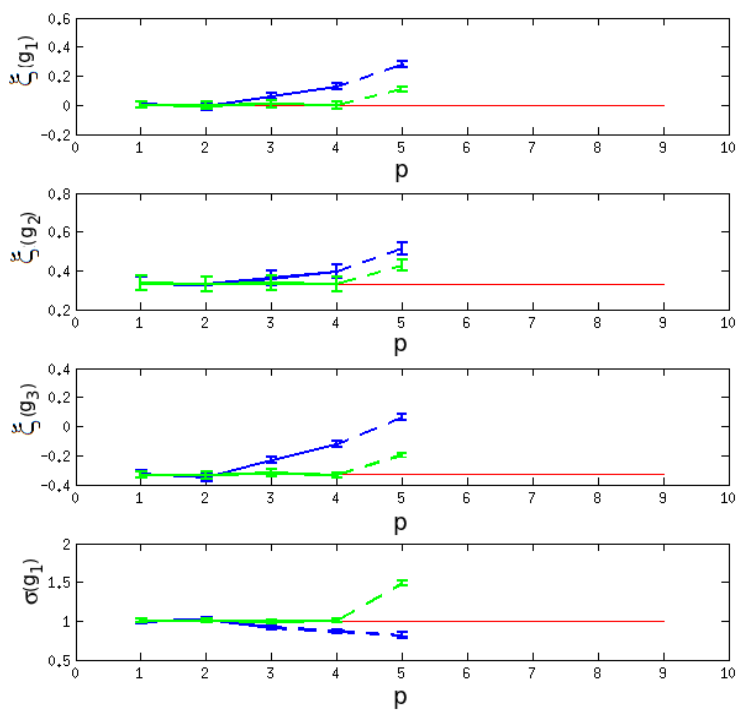

**Figure 9.19** GEV parameters VS intensity of the noise $\epsilon = 10^{-p}$ for the circle rotations perturbed map. Blue: $n = 10^3$, $m = 10^3$, Green: $n = 10^4$, $m = 10^3$. Red lines: expected values. $z \simeq 0.7371$. From the top to the bottom: $\xi(g_1), \xi(g_2), \xi(g_3), \sigma(g_1)$.

This suggests that is relatively hard to get rid of the properties of the underlying deterministic dynamics just by adding some noise of unspecified strength and considering generically long time series: the emergence of the smoothing due to the stochastic perturbations is indeed non-trivial when considering very local properties of the invariant measure as we do here.

### 9.7.1.2 Bernoulli Shift Map
In this set of numerical experiments, we want to study the impact of random perturbation on the statistics of extremes for the distance observables for the Bernoulli shift map mentioned above. Therefore, we consider the following map:

$$f_{\varepsilon\xi}(x) = 3x + \varepsilon\omega \bmod 1 \qquad (9.7.1)$$

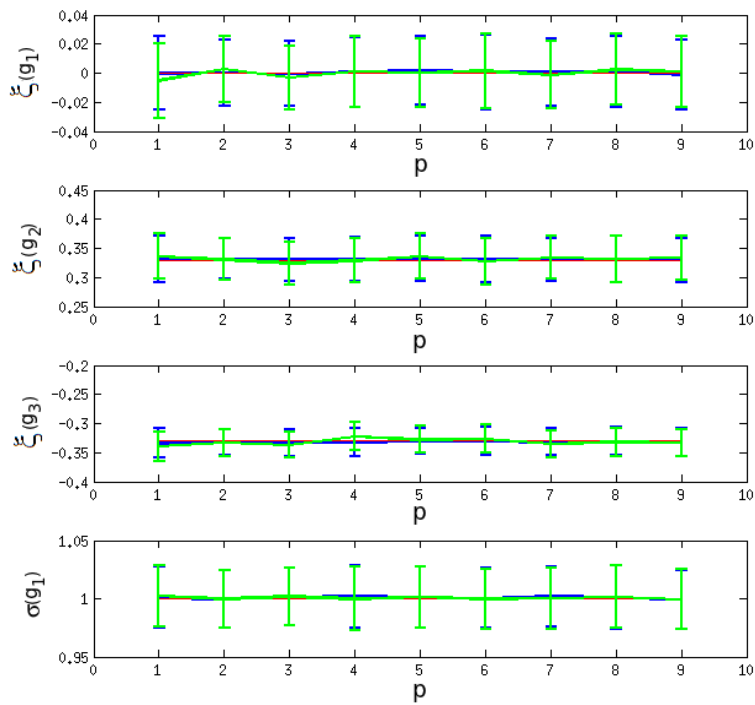

**Figure 9.20** GEV parameters VS intensity of the noise $\epsilon = 10^{-p}$ for the ternary shift perturbed map. Blue: $n = 10^3$, $m = 10^3$, Green: $n = 10^4$, $m = 10^3$. Red lines: expected values. $z \simeq 0.7371$. From the top to the bottom: $\xi(g_1), \xi(g_2), \xi(g_3), \sigma(g_1)$.

where $\omega$ is a stochastic variable with uniform distribution in $[-1, 1]$. As discussed in Chap. 7, the stationary measure for such a map is the uniform Lebesgue measure on the unit interval independently of the value of $\varepsilon$.

As a first check, we want to verify that when considering a non-periodic point $\zeta$, no difference should emerge between the deterministic and the randomly perturbed system. Results are shown in Fig. 9.20. It is clear that the stochastic perturbations do not introduce any changes in the type of statistical behaviour observed for a non-periodic point $z = 0.7371$ and no differences are encountered, even when the number of observations in each bin is increased. This is compatible with the idea that the intrinsic chaoticity of the map overcomes the effect of the stochastic perturbations. Summarising, extremes do not help us in this case to distinguish the effect of intrinsic chaos and the effect of adding external noise.

As discussed above and thoroughly studied in Chap. 7, the EVLs of distance observables for systems like the Bernoulli shift map are modified whenever $z$ is a periodic point of prime period $p$, so that the limit law reads as $e^{-\theta\tau}$ instead of the usual $e^{-\tau}$ realised for the remaining points of the attractor, where $0 < \theta < 1$ is the EI. Hence, we study what is the level of noise such that we indeed observe the disappearance of such a dichotomy, see Sect. 7.3.

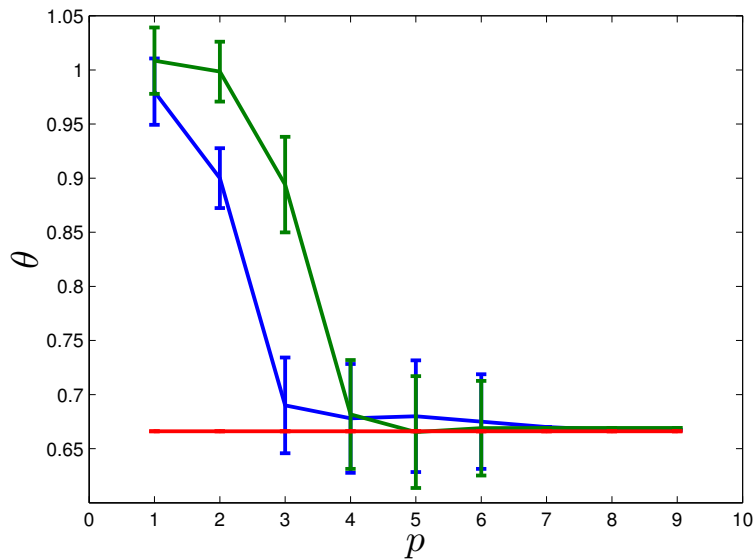

**Figure 9.21** Extremal index $\theta$ VS intensity of the noise $\epsilon = 10^{-p}$ for the ternary shift perturbed map. Blue: $n = 10^3$, $m = 10^3$, Green: $n = 10^4$, $m = 10^3$. Red line: theoretical $\theta$ for $z = 0.5$ of the unperturbed map.

In Fig. 9.21 we present the results for the EI obtained taking as reference point $\zeta$ the periodic point $z = 1/2$ of prime period 1, for which $\theta = 2/3$. As discussed in [79], we cannot use the usual fitting procedure for the GEV, since it always renormalises in such a way that the EI seems to be one. Instead, in order to observe extremal indices different from one, we have to fit the series of minimum distances to the exponential distribution by normalising a priori the data.

The results clearly show that we are able to recognise the perturbed dynamics as the extremal index goes to one when $\epsilon$ increases. Interestingly, the separation from the value expected in the purely deterministic case is observed only for relatively intense noise. Moreover, when longer time series are considered (green experiment), the stochastic nature of the map becomes evident also for weaker perturbations. Finally, it is clear that the numerical noise (corresponding to $\epsilon \simeq 10^{-7}$) is definitely not sufficiently strong for having a notable impact on the statistics of the deterministic system.

### 9.7.2
### Observational Noise

As a final numerical exercise, we wish to propose a simple example aimed at clarifying one of the possible effects of having observational noise on the statistical properties of extremes. We can construct a mathematical model for studying observational noise on time series following what proposed in Sect. 7.5. We consider a one-dimensional map describing the time evolution of the quantity of interest for the observations. The basic idea for simulating observational noise, as explained in Eq. 7.5.1, boils down to considering as actual observation the sequence $(y_n)_{n \in \mathbb{N}}$ obtained by perturbing the deterministic orbit $(f^n(x))_{n \in \mathbb{N}}$ of a point $x \in \mathcal{M}$ at each time instant with uncorrelated noise with a compact support and controlling the intensity of the noise with a parameter $\varepsilon$. It is important to keep in mind that the noise does not impact the orbit (it has no dynamical effect).

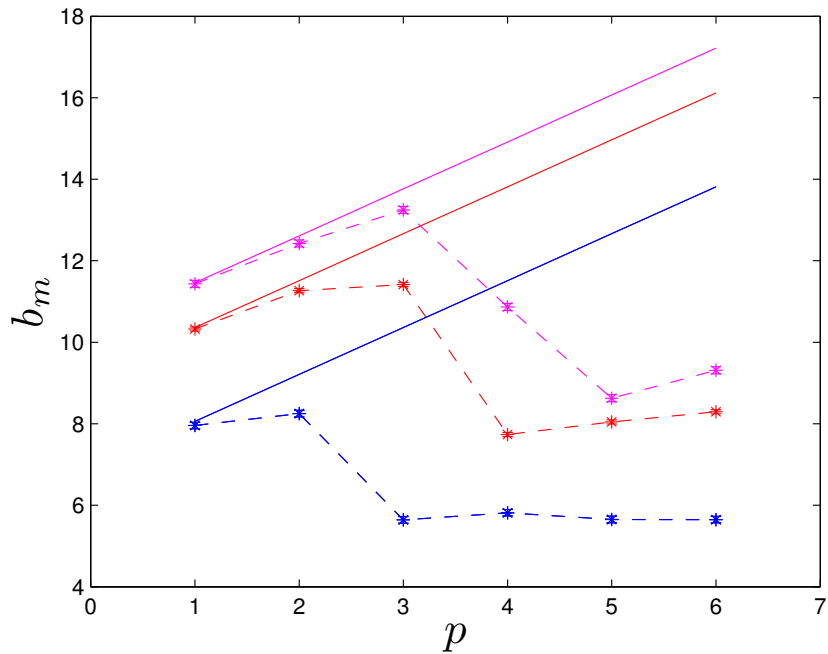

**Figure 9.22** Normalizing sequence $b_m$ vs intensity of the noise in terms of $p$ (we recall that $\epsilon = 10^{-p}$) for the Manneville-Pomeau map (Eq. 5.5.1). We recall that the dashed error-bars display the average of $b_m$ over 30 realizations and the standard deviation of the sample. Solid lines the theoretical values. The blue, red and magenta curves respectively refers to $m = 1000, 10000, 30000$, $z = 0$. $n = 1000$ for all the cases considered.

It is especially interesting to consider to study the extremes of the function $\phi(x) = -\log(y, \zeta)$, where $y$ is constructed as above and $\zeta$ is a certain chosen value for our observation. Large values of $\phi$ correspond to having values of the observations very close to $\zeta$. For given time series, we can say that a certain value $\zeta$ is highly recurrent if there are many occurrences of high values of function $\phi$, and highly sporadic if the opposite is true. Going back to what discussed in Chap. 1, if we make the identification between rare and extreme, we have that *highly sporadic values of our observable correspond to extreme events*. This is a rather different point of view with respect to what we have most typically proposed in this book, even if it is also based on EVT.

In Eq. 7.5.2, we find that the normalising sequence $b_n$ for the extremes of the function $\phi$ depends on the target point $\zeta$ *via* the local density of the invariant measure in a ball whose radius is given by the error $\varepsilon$. Therefore, if the point $\zeta$ is visited with less frequency, so that its the local density is small, one needs to go to higher values of $n$ in order to have a reliable statistics.

Here we want to test that the order of $n$ needed to get convergence to the asymptotic $b_n$ is lower for a highly recurrent point then for a sporadic one. As highly recurrent point we choose $z = 0$ for the Pomeau-Manneville map given in Eq. 5.5.1. The experiment consists in computing 30 realizations of the maps perturbed with observational noise. Again, we fit the maxima of the observable $\phi$ to the GEV distribution by using the $L$-moments procedure and compare the experimentally obtained values for $b_n$ (equal to the fitted value of $\mu$) to the theoretical ones stated in Proposition 7.2.4.

We report here the results for three different bin lengths $m = 1000, 10000, 30000$ in Fig. 9.22 for the Pomeau Manneville map. The figures show how $b_n$ varies as a function of the noise $\epsilon = 10^{-p}$, in terms of $p$. We observe convergence towards the theoretical values (solid lines) for high values of $\epsilon$ (low values of $p$), whereas in the limit of weak noise one must increase the bin lengths to get convergence.

The main result to be highlighted here is the better convergence of highly recurrent points with respect to the ones visited sporadically. This important property can be used to study time series recurrences and identify extremes as the points visited rarely for which the convergence towards the asymptotic parameters is bad. So, we end up with the paradoxical idea that a way to define extremes in a time series is to look for those values such that the GEV fit of a given function ($\phi$) does not converge given a time scale given by $n$. So, we can attach the quality of being extreme to an observation given how often it recurs (its rarity) on the time scale $n$.

The main advantage of studying recurrence properties in this way over applying other techniques is due to the built-in test of convergence of this method: even for a point rarely recurrent there will be a time scale $\bar{n}$ such that the fit converges. For smaller $m$, we can therefore classify such a $\zeta$ as a *sporadically recurrent* point of the orbit as explained in [80]. There, we show how to use this property to define rigorous recurrences in long temperature records collected at several weather stations. This application is discussed in detail Chap. 10.

In another study [299], we suggested a quantitative way to discriminate between highly recurrent points and sporadic points of the dynamics in a rather algorithmic way. Basically, one can assess the minimum bin length $m$ such that the fit to one of the EVL converges, that is the value of $m$ such that a sporadic point becomes a normally recurring one. For highly recurring points this typical value of $m$ is of order $10^3$, whereas for quasiperiodic dynamics it can be larger than $m = 10^9$.

# 10
# Extremes as Physical Probes

The previous chapter was devoted to showing examples of how the features the EVLs for distance observables of a dynamical system can inform on the mathematical properties of the underlying dynamics. In this chapter, we want to provide examples of how EVT for dynamical systems can be helpful in gaining understanding on the physical properties of the system we are studying. Building upon some the results described in Chapter 4 and Chapter 9 and the theory relevant for the extremes of distance observables, we present a new method for proposing an alternative definition of extreme in a temporal record of data, defined how as a rare recurrence for the specific time scale of reference. We provide an example where such a method is applied successfully for studying extremes in European surface temperature records. Extremes of physical observables are instead shown to provide a new tool for studying the properties of critical transitions in complex system and predicting the critical value of the control parameter determining the occurrence of tipping points. We first present an example where EVT helps understanding the multi stability properties of the plane Couette flow in the regimes around the transition to turbulence. We then apply the method to a well-known low dimensional stochastic dynamical system developed as toy model for studying the transition to turbulence.

## 10.1
## Surface Temperature Extremes

We consider times series of spatially gridded daily mean surface temperature taken from the European Climate Assessment and Dataset project database. One can find additional information on the data and download them at `http://eca.knmi.nl`. At each grid point, consider two different time series: the series of daily mean temperatures $T_k$ and the series of daily mean temperature anomalies $T_k^a$, obtained by subtracting from $T_k$ the best estimate of the seasonal cycle $S_k$. An example of $T_k$ and $T_k^a$ with relative histograms is presented in Fig. 10.1. The basic assumption we are taking is that by subtracting the seasonal cycle we isolate in $T_k^a$ the chaotic and stationary component to the variability of the time series. We assume that the influence of actual meteorological disturbances are responsible for defining the chaotic behaviour of the time series. The hypothesis of stationarity can, instead, directly checked by applying, *e.g.*, a so-called *KPSS* test [300].

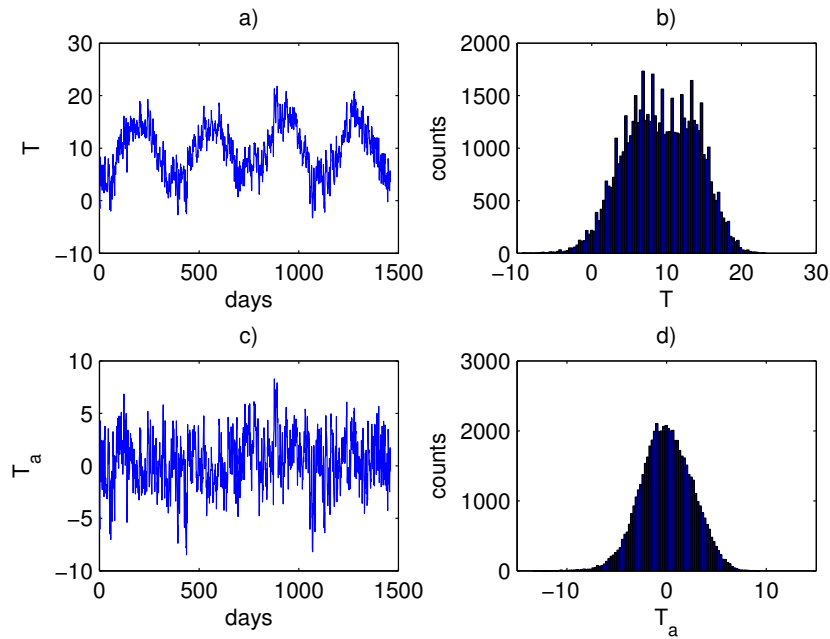

**Figure 10.1** Top: An example of temperature series (a) and its histogram (b). Bottom: the same for a temperature anomalies series (c) and the correspondent histogram (d). All the plots refer to Armagh (UK) weather station.

We report the analysis of three stations chosen for their distinct climate features. We first treat the case of the temperature records $T_k$.

- Armagh (UK) is situated in Northern Ireland, whose climate is influenced by the Atlantic Ocean with moderate temperature ranges; see Fig. 10.1.
- Milan (Italy) is situated in the Po Valley, so that its climate is influenced by the victim of the Mediterranean sea and that, at the same time, has relatively strong continental features, due to the proximity to the Alps and the Apennines mountain ranges.
- Vienna (Austria), which lies in a more central region of Europe and has a distinctly continental climate .

These three locations possess rather extensive (length of 161, 246 and 156 years, respectively) continuous daily records of temperature [301].

### 10.1.1
### Normal, Rare and Extreme Recurrences

We show how the link between EVT and HTS can be used to study the recurrence properties of chaotic time series and defining in a rigorous way what it is means that an even is *rare* given a specific reference time scale. A key aspect is the investiga-tion of the properties of if and how (fast) the statistics of BM converges to the GEV statistical model.

Our analysis is inspired by a methodology developed for Lyapunov exponents [284]. We introduce the following algorithm aimed at defining which are the rarely recurrent values of a time series:

1) Take a given time series $X_j$, $j = 1, 2, \ldots, s$,
2) Choose an index value $j$ and the corresponding element of the time series $X_j$; this constitutes the our reference point $\zeta$ (see Sect. 4.2.1);
3) Compute the series $Y_k^j = -\log(\text{dist}(|X_k - X_j|)$, $k = 1, 2, \ldots, s$, $j = 1, 2, \ldots, s$.
4) Divide the series $Y_k^j$ $k = 1, \ldots, s$ into $K$ bins each containing $N$ data ($NK = s$) and extract the BM $M_p^j$, $p = 1, \ldots, K$ for each reference value $\zeta = T_j$, $j = 1, \ldots, N$.
5) Perform a GEV fit of the empirical d.f., perform a Lilliefors test [302] against the hypothesis of Gumbel law to check whether the fit has succeeded or has failed. The rationale for this is the following: by construction, $\zeta$ should be zero, and if the best fit is not compatible with this hypothesis, we have to conclude conclude that the BM selection procedure does not choose true extremes, the reason being that $N$ is too small.

- If the fit is satisfactory, one repeats the experiment for shorter and shorter bin lengths $N$ and finds the smallest $N_{min}$ such that, GEV fit of the BM $M_j$, $j = 1, \ldots, s/N_{min}$ converges.
- If the fit fails, one repeats the experiment for longer and longer bin lengths $N$ until the value $N_{min}$ is reached, when the BM eventually selects only good candidates for extremes and the GEV fit is successful.

$N_{min}$ defines the longest time scale over which the value $X_j$ can be considered as occurring rarely. We propose here the following viewpoint for defining a rare record. Given a value $\bar{N}$, we can define which of the values of the time series $X_j$ $j = 1, 2, \ldots, s$ are rare according to time scales equal to $\bar{N}$ or shorter. In order to achieve this, we exploit the results given in [46, 77].

## 10.1.2
### Analysis of the Temperature Records

Before starting our analysis, we have to consider the nontrivial problem of addressing the fact that such time series are truncated at $p = 3$ digits because of the way the recording has been performed. As discussed by [80], a blind application of the methodology presented above would lead to a divergent fit as the recurrence distribution will appear as a collection of Dirac delta functions. In order to solve this issue, we need to transform our time series into a time series extracted from a fictitious instrument able to have very high precision (ideally, *continuous*) readings of the temperature data. We can solve the problem by adding to each reading $T_k$ a random number extracted from uniform distribution, so that the observations are not altered (at the level of precision of the real instrument), but continuity is recovered. We then redefine perform our analyses on the noisy version of the original time series. This idea is based on the mathematical results discussed in Chapter 7 and first presented in [65].

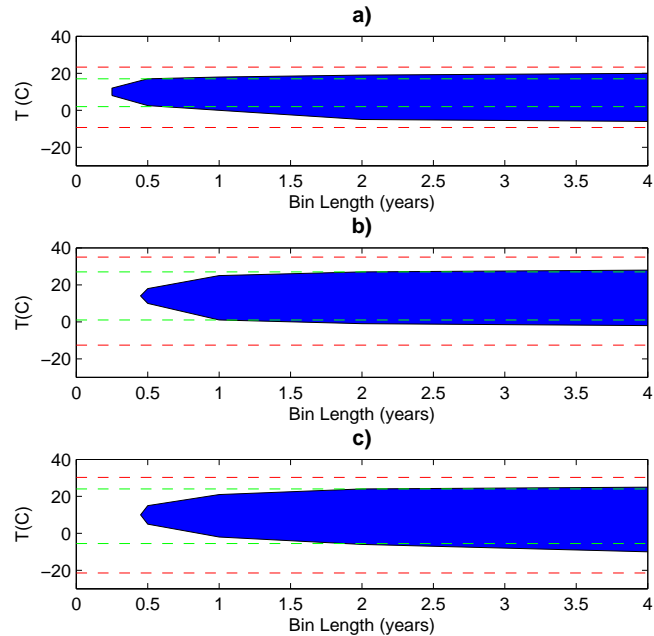

**Figure 10.2** Range of temperatures with convergent Gumbel fit (blue) for different bin lengths. The blue area indicates ranges of normal variability in the time scale given by the bin length. Red dotted lines: absolute extremes of the temperature series for Armagh (a), Milan (b) and Vienna (c). Green dotted lines: thresholds detected with the classical GPD approach.

We follow the procedure described in Sect. 10.1.1 and perform the GEV fits using an MLE technique, as discussed in the previous Chapter. We are then able to find for a given value $N_{min}$ whether a certain $\zeta$ between the absolute recorded extremes is or is not rarely recurrent. The results are presented in Fig 10.2 for the stations located in Armagh (a), Milan (b), and Vienna (c). The experiments have been repeated for different bin lengths between 3 months and 4 years. In Fig 10.2, the blue area represents the range of recurrent values of the reference temperature $\zeta$, obtained as those whose corresponding extremes of the $Y$ distance observables can be fitted successfully by a Gumbel distribution. Therefore, we can say that for each time time scale, the blue range defines what we propose as a rigorous definition of normal variability with respect to the time scale defined by the bin length. On the same time scale observing a fluctuation which goes beyond the blue region constitutes, instead, a genuine extreme event. These are all located outside in the white region and correspond to unsuccessful Gumbel fits of the time series of $M_p^j$.

Note that for bin lengths shorter than six months, the Gumbel fit fails at any reference temperature, suggesting that below the bin length is too short for observing proper recurrences near any $\zeta$. The only exception is registered at Armagh where, due to the limited seasonal temperature excursions, the convergence is achieved for $8\,\mathrm{C} < \zeta < 12\,\mathrm{C}$, already for a 3 months bin length. In general, one observes that, when the bin length is increased, the temperature range accepted as including non uncommon climate fluctuations (normal variability) increases. We observe that in Milan one could define $0\,\mathrm{C}$ as an extreme temperature with respect to a bin length of half a year but not when considering a bin length of 4 years. On such time scale, $0\,\mathrm{C}$ is indeed part of the normal variability. It is interesting to check whether the extremes found with this recurrence analysis (*extremes as rare events*) are related to the minimal thresholds one would derive using a POT of the actual temperature records (*extremes as large events*). These values have computed by fitting GPD distributions and they are represented in Fig. 10.2 by green dotted lines. The threshold values estimated with the POT approach are similar to the one obtained with our method, although they are remarkably different at Armagh station: the GPD method identifies as extremes temperatures values classified as normal by the recurrences method for bin longer than one year. The POT approach does not discriminate between extreme events belonging to the same cluster - temperatures beyond the threshold in consecutive days. In other words, this method does not carry information about the correlation structure of the time series: any random resorting of the data would produce exactly the same threshold values. This is why thresholds appear lower than the one detected by our approach.

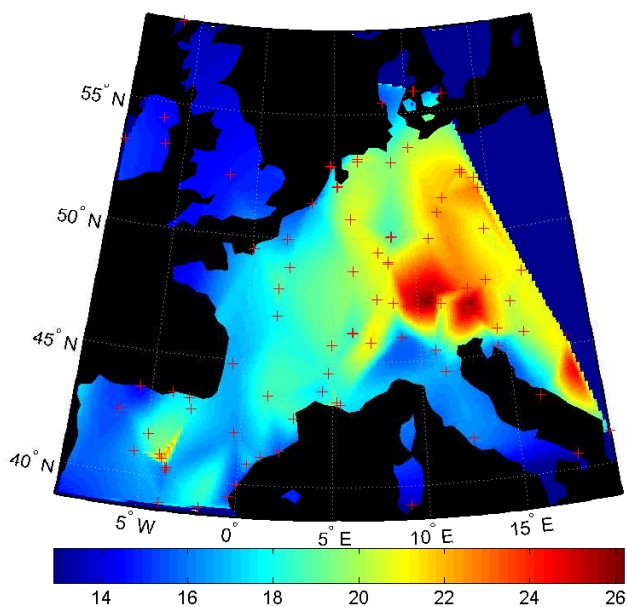

**Figure 10.3** Map of the range of *admissible temperature excursions* for the European region, obtained by considering the interval of temperature anomalies $\zeta$ such that the fit passes the Lilliefors test. The red crosses represent the location of the stations used for the analysis. The straight line near the right border represents the limit of the data-set.

We have repeated the same analysis reported in Fig. 10.2 for the series of the anomalies $T_a^k$. The main advantage of using temperature anomalies consists in the possibility of comparing the climatology of different locations. Let us consider a one year bin length: at Armagh, in a temperate, marine climate , only anomalies up to $\pm 6$ C can be considered as normal annual variability according to the described method. In Vienna, the continental climate permits large temperature excursions so that we find that anomalies up to $\pm 10$ C are normal part of the annual variability. For Milan, the range is reduced to $\pm 7$ C, an intermediate situation between Armagh and Vienna, in agreement with climatic features influenced both by the both the Mediterranean sea and by the location in the Po valley. Given a time scale, we can plot the corresponding *normal* value of the range of variability of $T_a$ at all the European locations for which at least 60 years of daily data are available. Hence, we can construct a map of Europe providing novel climatological information on variability at a given time scale. Results obtained by selecting a bin range of 1 $y$ are reported in Fig. 10.3. Different climatic regions are well highlighted: the British Isles, Brittany, Italy and the coastal areas of the Iberic peninsula have a milder climate with a significantly lower range of admissible temperature excursions. Very large fluctuations of temperature are instead amicable as normal part of the yearly internal variability in continental Europe and in mountainous areas. Note that near the central Iberic peninsula the field is strongly influenced by a specific station, Navacerrada, which is situated at 1800 m.s.l. and features extremely large annual variability, where normal anomalies of up to 24 C are observed.

## 10.2
## Dynamical Properties of Physical Observables: Extremes at Tipping Points

Let us consider a dynamical system, either deterministic or noisy, controlled by some parameter $\lambda$ which, when decreased below some value $\lambda_{\mathrm{crit}}$, drives it through a *critical transition*, which leads to a qualitative change in the dynamics. Here the word *critical* has the meaning it takes, say, in environmental sciences, where the expression *tipping point* is also used. In dynamical systems theory, one would speak of a *saddlenode* bifurcation or some appropriate generalization of it, namely a *crisis* [303].[1]

In high dimensional systems, critical transitions can be easily understood (as well as portrayed in a graph) by studying how specific physical observables (*e.g.* energy) experience abrupt changes when the system crosses its tipping point. Hence, EVT enters the picture as we are interested in studying the nature of the large fluctuations of such observables in the vicinity of a crisis.

---

1) In statistical physics, the present situation would correspond to a (discontinuous) *first-order* phase transition, *e.g.* a liquid-gas transition, as opposed to a (continuous) *second-order* phase transition, *e.g.* a para-ferromagnetic transition, studied within the framework of *critical phenomena*, where the word 'critical' thus gets a different meaning through the definition of *critical exponents* and related *universality classes* [304].

As discussed in Chapter 8, on general grounds we may expect that in Axiom A-like systems (in the sense of the chaotic hypothesis [266]) physical observables have bounded fluctuations and that their extremes follow Weibull distributions [81, 44]. The closer we are to a crisis, the more likely is for the system to explore regions of the phase space close to the saddle, so that there is an increasing probability that the physical observable will have anomalous values and feature (rare) very large fluctuations, *much* larger than the extreme fluctuations observed in the system far away from the crisis.

The idea we wish to explore here is that as we get close to the crisis, *e.g.* increasing the system's control parameter $\lambda$ towards the critical value $\lambda_{\mathrm{crit}}$ which defines the bifurcation, the shape parameter $\xi$ describing the extremes of the physical observable becomes larger and larger and crosses the zero value exactly when the tipping point is reached. Therefore, studying how $\xi$ depends on $\lambda$ might provide useful information for reconstructing how far we are from the the critical value $\lambda_{\mathrm{crit}}$.

A traditional approach in studying the fluctuations of physical observables near tipping points is based on the idea that near the crisis fluctuations of greater amplitude will be observed towards the state the system is *doomed* to fall into as $\lambda$ takes the value $\lambda_{\mathrm{crit}}$. The idea is then to use anomalous values of the skewness of the probability distribution of the observable, which measures its asymmetry, as an early warning indicator of a tipping point [305]. Such a method has good potential if the probability distribution of the observable is approximately symmetric for $\lambda \ll \lambda_{\mathrm{crit}}$, but may fail if the distribution is already skewed.

However, we can adapt this idea to the analysis of the extremes. In particular, we show below how comparing statistics the negative and positive extremes of an observable may be extremely effective for locating the critical transition. We will use as especially instructive example the analysis of the extreme fluctuations of turbulent kinetic energy in the plane Couette flow, refering to results contained in [75].

## 10.2.1
## Extremes of Energy for the Plane Couette Flow

The plane Couette flow can be described as resulting from the shearing of a viscous fluid in the space between two parallel plates in relative motion. The plates, at a distance $2h$, translate in opposite directions at a speed $U_w$ and the flow results from the viscous drag acting on the fluid with kinematic viscosity $\nu$. The nature of the flow regime, either laminar or turbulent, is controlled by a single parameter, the Reynolds number $R = U_w h/\nu$, which plays the role of the control parameter $\lambda$ mention above; see Fig. 10.4 for a schematic representation of the flow.

The laminar flow depends linearly on the coordinate normal to the plates and is known to remain stable against infinitesimal perturbations for all values of $R$, while turbulent flow is instead observed under usual conditions when $R$ is sufficiently large, typically of order 400–500, when increasing $R$ without particular care. As $R$ is decreased from high values for which the flow is turbulent, a particular regime appears at about $R_{\mathrm{t}} \approx 415$ where turbulence intensity is modulated in space [306]. When the experimental setup is sufficiently wide, a pattern made of oblique bands, alter-natively laminar and turbulent, becomes conspicuous. Bands have a pretty well defined wavelength and make a specific angle with the streamwise direction. As $R$ is further decreased, they break down and leave room to the laminar base flow below $R_g \approx 325$. Experiments show that the streamwise period[2] $\Lambda_x$ of the band pattern is roughly constant ($\Lambda_x \simeq 110h$) while the spawise period $\Lambda_z$ increases from about $55h$ close to $R_t$ to about $85h$ as $R$ decreases and approaches $R_g$ [306].

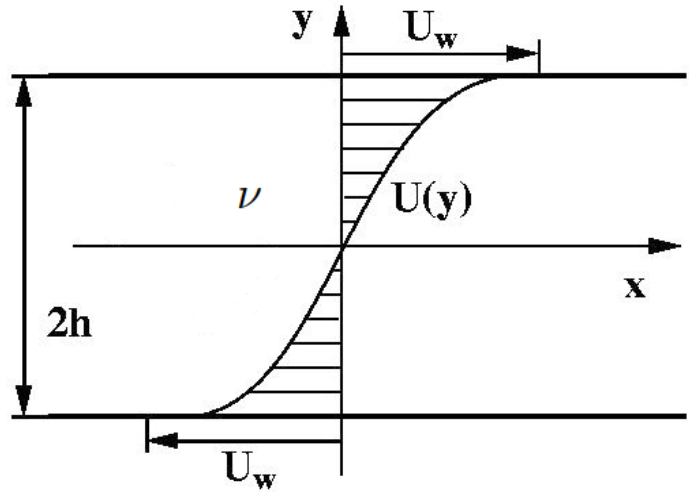

**Figure 10.4** Schematic representation of the experimental setting of a plane Couette flow. A three-dimensional fluid with kinetic viscosity $\nu$ lies between two plates of infinite extension along $y$ and $z$ (not shown). The two plates are separated by a distance $2h$ along $y$ and translate in opposite direction along $x$ at a constant speed $U_w$, so that a time dependent flow $\mathbf{v}(\mathbf{x}, \mathbf{y}, \mathbf{z}, \mathbf{t})$ is established (the time and $z$-average of the $x$ velocity profile $U(y)$ is portrayed in the figure).

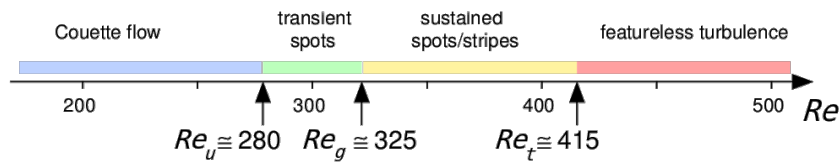

**Figure 10.5** Dynamical regimes of the plane Couette flow. See [307].

Whereas the turbulence self-sustainment process in wall-bounded flows is well understood [308], the mechanisms explaining band formation are still somewhat mysterious. The transition displays a large amount of hysteresis. A similar situation is to be found in several other flow configurations, circular Couette flow, the Couette flow sheared by coaxial cylinders rotating in opposite directions, plane channel, the flow between two plates driven by a pressure gradient, as well as in Poiseuille flow in a circular tube. See [309, 307] for a comprehensive review of these phenomena.

---

2) By convention, the streamwise direction is along $x$, the normal to the moving plates defines the $y$ direction, and $z$ denotes the spanwise direction (not shown in Fig. 10.4).

### 10.2.1.1   Conditions of the Numerical Experiment

The transition to turbulence in plane Couette flow has been studied numerically by a number of authors. System sizes required to observe the oblique band regime in Navier–Stokes DNSs are numerically quite demanding [310]. In order to reduce the computational load, Barkley and Tuckerman performed their computations in a cleverly chosen narrow but inclined domain [311]. The drawback is however to freeze the orientation beforehand, forbidding any angle or orientation fluctuation. A recent work has shown that another way to decrease computer requirements was to accept some under-resolution of the space dependence, especially in the wall-normal direction $y$ [312]. All qualitative features of the transitional range are indeed well reproduced in such a procedure, including orientations fluctuations. Quantitatively, the price to pay appears to be a systematic downward shift of the $[R_g, R_t]$ interval as the resolution is decreased.

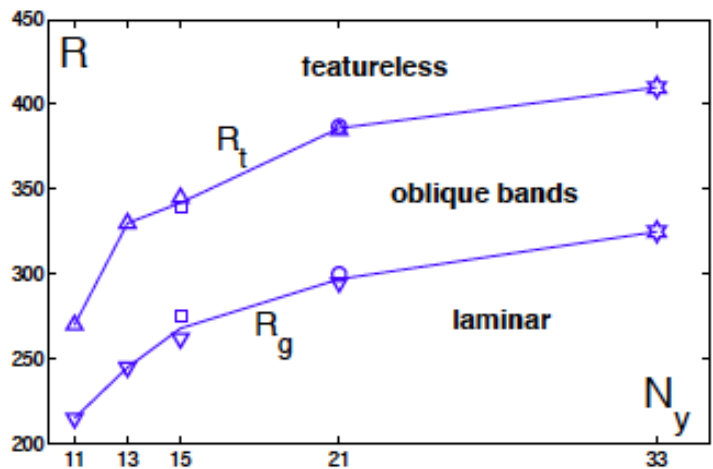

**Figure 10.6** Effect of changing resolution on the dynamical regimes of the plane Couette flow. A downward shift of the $[R_g, R_t]$ interval is found as the resolution is decreased [312].

The results shown here are obtained using simulation performed in a domain of constant size able to contain one pattern wavelength in each direction, i.e. $(L_x, L_z) \equiv (\Lambda_x \times \Lambda_z)$, with $L_x = 108$ and $L_z = 64$, using The open-source software CHANNELFLOW [313]. This size seems well adapted to the central part of the transitional domain, *i.e.* slightly too wide for $R \approx R_t$ and slightly too narrow at $R \approx R_g$, with mild consequence on the effective value of these thresholds, as guessed from a Ginzburg–Landau approach to this pattern forming problem [314]. Such finite-size effects [315] also account for the intermittent reentrance of featureless turbulence.

CHANNELFLOW is a Fourier–Chebyshev–Fourier pseudo-spectral code is dedicated to the numerical simulation of flow between parallel plates with periodic in-plane boundary conditions. In the wall-normal direction (see Note 2), the spatial resolution is a function of the number $N_y$ of Chebyshev polynomials used. The in-plane resolution depends on the numbers $(N_x, N_z)$ of collocation points used in the evaluation of the nonlinear terms. From the 3/2 rule applied to remove aliasing, this corresponds to solutions evaluated in Fourier space using $\frac{2}{3} N_{x,z}$ modes, or equivalently to effective space steps $\delta_{x,z}^{\text{eff}} = \frac{3}{2} L_{x,z}/N_{x,z}$. Numerical computations have been performed using three different resolutions: low ($N_y = 15$, $N_x = L_x$, $N_z = 3L_z$), medium ($N_y = 21$, $N_x = 2L_x$, $N_z = 6L_z$), and high ($N_y = 27$, $N_x = 3L_x$, $N_z = 6L_z$) for which we expect $[R_{\text{g}}, R_{\text{t}}] \approx [275, 350]$, $[300, 380]$, and $[325, 405]$, respectively; see Fig. 6 in [312].

### 10.2.1.2 The EVT Analysis of Turbulent Energy near the Critical Transition

In this section we show results about the changes in the extreme value distributions of quantity $E_{\text{t}}$ defined as the mean-square of the perturbation velocity $\tilde{\mathbf{v}}$, the difference between the full velocity field $\mathbf{v}$ and the base flow velocity $\mathbf{v}_{\text{b}} = yU_w/h\,\mathbf{e}_x$, where $\mathbf{e}_x$ is the unit vector in the $x$ direction. Physically speaking, apart from a factor $\frac{1}{2}$, this is the kinetic energy contained in the perturbation, which is zero in the case the laminar flow is obtained. Here we focus on the determination of $R_{\text{g}}$ using extremes as sketched above.

For each value of the Reynolds number, very long simulations are performed and, once the time series of $E_{\text{t}}$ has reached a stationary state, maxima (minima) are extracted in bins of fixed block length as described in the previous section. We then fit the maxima (minima after sign change) to the GEV distribution by using, in this case, a MLE method, as discussed in the previous Chapter. As now clear, the choice of the bin length $m$ is crucial: in the asymptotic regime the value of the shape parameter should be independent of $m$. We have tested that, within the confidence intervals, this happens for $m > 1000$.

The shape parameter $\xi$ is next analyzed as a function of the Reynolds number for both the positive and negative extremes of the energy. A first intuition on how the method should work comes from looking at the data series and the histograms shown in figure 10.7, see caption for details. The series in red refers to a value of Reynolds inside the band regime ($R = 300$), with fluctuations exploring a limited interval. The series in blue, with $R$ is fixed just above $R_{\text{g}}$ ($R = 277$), illustrates a clear tendency to intermittently visit states with very low values of the energy. These events, spotted in the green ovals, crucially contribute to a shift towards Fréchet laws since the fit to the GEV returns a Weibull EVL when removing them from the histogram.

### 10.2.1.3 High Resolution

Let us start with the localization of the global stability threshold $R_{\text{g}}$ in simulations performed at high resolution, namely $N_x = 216, N_y = 27, N_z = 384$. Results are shown in Fig. 10.8 (left column) for the shape parameter (upper panel), to be compared to the two common early warnings indicators based on the bulk statistics: the skewness (middle panel) and the variance (lower panel). When approaching $R = 322$ the shape parameter for the distribution of minima changes its sign, whereas for the maxima it remains negative in agreement with what was stated in the previous

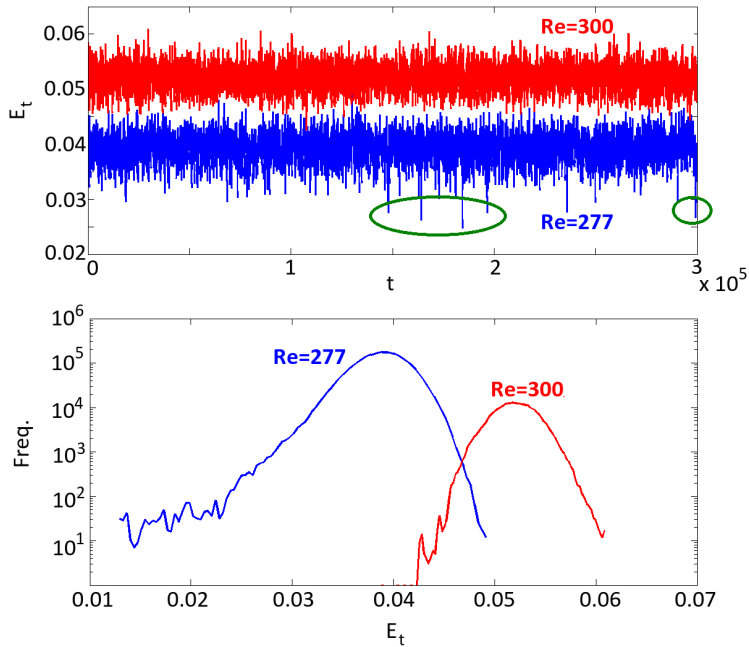

**Figure 10.7** Perturbation Energy $E_t$ for two simulations at low resolution. Upper panel: $E_t$ as a function of time. Lower panel: Histograms of $E_t$ in log-linear scale. Green circles indicate extremely rare events.

sections.

It is however evident that these results need confirmation since a limited set of Reynolds numbers has been studied and a single slightly positive value of $\xi$ has been obtained for $R = 322$, with error bars so large that the significance of the result is rather limited. The variance and the skewness of the time series follow what is expected from the statistics of global observables at a tipping point, namely a monotonic trend towards larger values. Since $E_t$ visits lower energy states, the skewness becomes more negative so that only the distribution of minima is affected. However, as noticed previously, no definite threshold value $R_g$ can be inferred from the consideration of the variance and skewness curves.

### 10.2.1.4 Medium Resolution

Using very high resolution is computationally expensive: the cumulated amount of CPU time required to produce series of length $s = 2.5 \cdot 10^5$ time units was beyond $10^5$ CPU hours, making it practically impossible to obtain much longer series with the available resources. In order to support our results, we have exploited the fact that downgrading the resolution preserves the qualitative features of the transition, up to a shift of transitional range, as shown in Fig. 10.6 and discussed in [312]. At medium resolution ($N_x = 216$, $N_y = 21$, $N_z = 384$), all time series have been stopped at $s = 2 \cdot 10^5$ time units. In these conditions, band breakdown was never observed for $R > 306$. The results shown in Fig. 10.8 (center column) confirm those at high resolution, with a slightly more pronounced change of sign of the shape parameter at $R = 306$ but point out the need of more and much longer series around the global stability threshold.

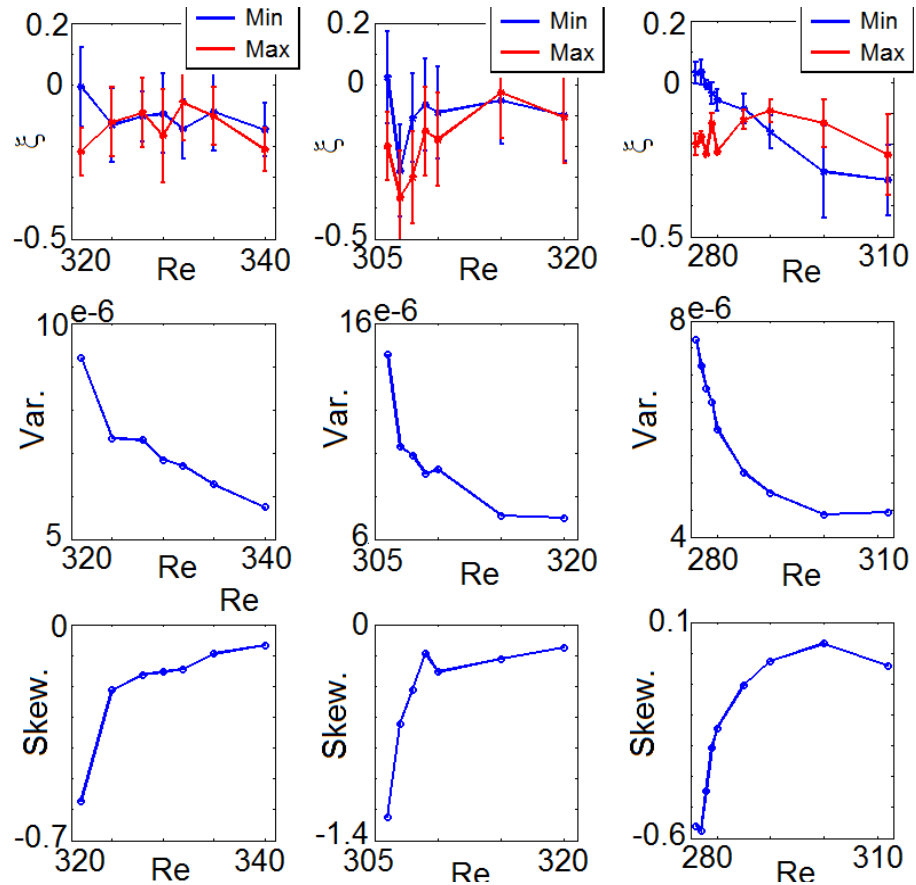

**Figure 10.8** Tipping point indicators for plane Couette flow as functions of $R$. Upper panel: Shape parameter $\xi$; red: maxima, blue: minima; error bars represent 95% confidence intervals, $m = 1000$. Center and bottom panels: Variance and skewness of the full series, respectively. Left: High resolution. Center: Medium resolution. Right: Low resolution.

### 10.2.1.5 Low Resolution

Further downgrading the resolution to $N_x = 108, N_y = 15, N_z = 192$ allowed us to produce series lasting nearly one order of magnitude longer than above, up to $2 \cdot 10^6$ time units. As a matter of fact, by collecting a greater statistics of maxima, the uncertainty on the estimation of the shape parameter could be greatly reduced, as shown in Fig. 10.8 (right column). In view of our proposal of trying to define $R_\mathrm{g}$ using extreme value statistics, the results at low resolution look much more convincing than those produced at higher resolutions since a clear monotonic variation of the shape parameter for minima is now observed upon decreasing $R$. As soon as $R_\mathrm{g} \leq 278$, the shape parameter $\xi$ describing the minima of the turbulent energy changes sign crossing the zero value.

The estimates of $R_\mathrm{g}$ determined here for the different resolutions studied are not much different from those given in [312] obtained by inspection of individual cases without any systematic criterion and using much shorter time series.

## 10.2.2
## Extremes for a Toy Model of Turbulence

In order to give more robustness to the link between critical transitions and change in the sign of the shape parameter $\xi$ describing the extremes of a suitably chosen physical observable, we should investigate a computationally cheaper model where it is easier to achieve a high-quality statistical inference. A good candidate for testing the identification of the global stability threshold using methods based on GEV parameters is a slightly modified version the model originally introduced in [316]:

$$dX/dt = -(\mu + u\xi(t))X + Y^2, \qquad dY/dt = -\nu Y + X - XY. \quad (10.2.1)$$

Here $X$ and $Y$ may be related to the amplitudes involved in the self-sustaining process of turbulence. Parameters $\mu$ and $\nu$ are damping coefficients accounting for viscous effects and assumed to vary as $1/R$. Non-linearities preserve the energy $E = \frac{1}{2}(X^2 + Y^2)$ in the same way as the advection term of the Navier–Stokes equations. Noise is here introduced in a multiplicative way *via* the term $u\xi(t)$, where $\xi(t)$ is a white noise and $u$ its amplitude, as proposed by Barkley [317]. A saddle-node bifurcation takes place at $\mu\nu = \frac{1}{4}$. The trivial solution $X = Y = 0$, corresponding to laminar flow, competes with two nontrivial solutions on the interval $\mu\nu = \left[0, \frac{1}{4}\right]$, the nontrivial solution being assimilated to turbulent flow. Unlike the additive noise considered in [318], the multiplicative noise taken here does not affect the trivial state and can be understood as a fluctuating turbulent-like contribution to effective viscous effects. Whenever the system undergoes a transition towards the laminar state, the simulation is restarted from the stable nontrivial fixed point.

Whereas for plane Couette flow only one simulation could be performed at each Reynolds number because of computational limitations, here we can easily produce ensembles of realizations for a given set of parameters $(\nu, \mu, u)$, extract corresponding GEV shape parameters, and average them over the realizations. We focus on the extremes of the energy $E$, which provides an efficient way for distinguishing the turbulent state from the laminar. In view of locating the critical transition, we want to relate the change of sign of the GEV parameter $\xi$ with the fact that the probability of transition between the two regimes becomes significant.

A remark is needed. Note that, as opposed the case of the plane Couette flow, the presence of noise allows for (exceedingly rare) transitions from the laminar to the turbulent state also far from the critical condition $\mu\nu = \frac{1}{4}$. This results into the fact that, in principle, the extremes of the energy fluctuations are *always* Gumbel distributed, even with a noise of extremely small amplitude. Nonetheless, the convergence of the empirical data to the Gumbel is exceedingly slow, so that that at any practical purpose, when finite (even if relatively long) time series are analysed, the best fit will be a Weibull distribution. Instead, when we get close to the bifurcation *and* the noise is strong enough, the convergence to the Gumbel distribution is very rapid and practically detectable.

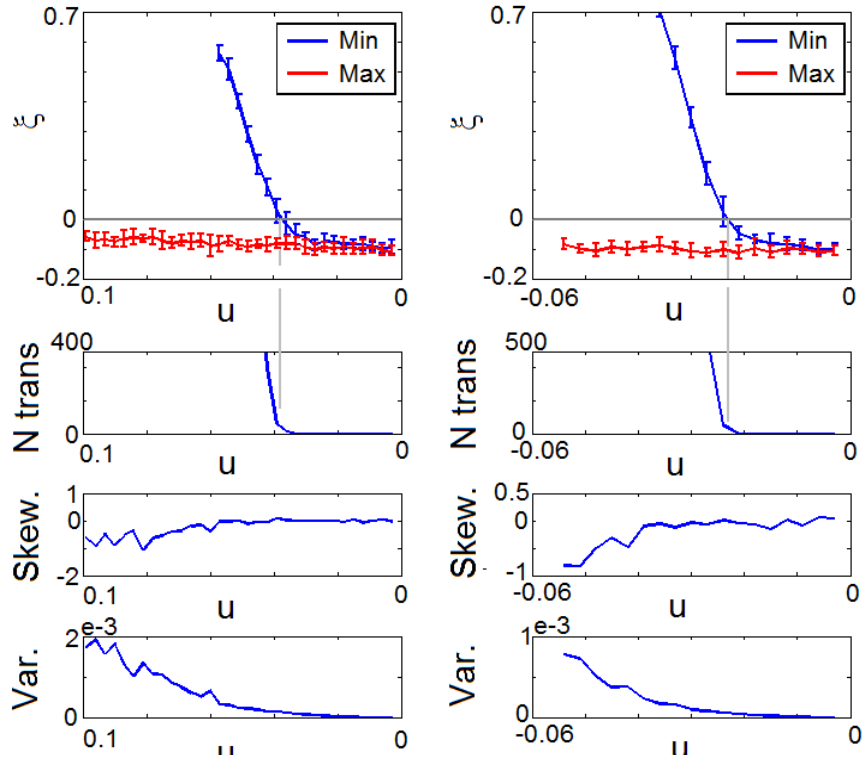

**Figure 10.9** Variation of different indicators of critical transitions as functions of the noise intensity for model (10.2.1) with $\mu = 1$. Left: $\nu = 0.2475$. Right: $\nu = 0.2487$. Top row: Averaged shape parameter $\xi$; red: maxima, blue: minima; error bars represent the standard deviation over the ensemble of 30 independent realizations. Second row: Number of transitions observed (see text). Vertical lines between the two top panels point to the critical value $u_c$ for which $\xi = 0$. Third and fourth row, respectively, averaged variance and skewness. Note that $u$ increases to the left and not to the right as usual.

The results of two different set of simulations are shown in Figure 10.9. Here, the control parameter $\lambda$ is the intensity of the noise $u$ whereas $\mu$ and $\nu$ are kept fixed. The left plots refer to the case $dt = 0.01$, $\mu = 1$, $\nu = 0.2487$, $n = 10^3$, $m = 10^6$, whereas the right ones refer to $\nu = 0.2475$ with the other parameters left unchanged. For each values of $u$, ensembles of 30 realizations has been prepared. The upper panels show the variation with $u$ of the shape parameter averaged over the realizations with error bars corresponding to the standard deviation over each ensemble. The plots in the second row display the number $N_{\mathrm{tr}}$ of times the system has undergone a critical transition from turbulent to laminar during the simulation performed at given $u$.

In both sets of simulations, the shape parameters vary similarly to the case of the plane Couette flow. For the distribution of extreme negative fluctuations, $\xi$ crosses zero from below when $u$ reaches a specific critical value $u_c$ (vertical lines between the first and second rows in Fig. 10.9), which depends on the other parameters $\mu$ and $\nu$. Instead, when looking at the distribution of the maxima of the turbulent energy, so features are observed, as expected.

We are confident in saying that $\xi = 0$ acts as a flag suggesting the presence of critical transitions because when $\xi$ crosses 0 and $u$ crosses the value $u_c$ the number of observed transitions $N_{\mathrm{tr}}$ increases substantially, as can be seen by comparising the two upper rows of Fig. 10.9.

The plots in the third and fourth row of Fig. 10.9 refer, respectively, to the variance and the skewness of the bulk statistics of $E$. In contrast, these indicators, while displaying the expected trends, do not show any specific feature or obvious flag allowing us to locate the threshold.

## 10.3
## Concluding Remarks

In this chapter we introduced some EVT-based tools of potential utility for studying physical systems. We have shown that investigating the properties of the invariant measure via the recurrence approach can be used to link the definition of normal and rare events to specific time scales of the dynamics. Such a characterization is a key problem for the mitigation of geophysical extremes. Return times of extreme events are of crucial relevance for stakeholders in order to define suitable mitigation and adaptation plans, *e.g.* in the case of the construction of specific infrastructures. The statistical approach for the computation of return times involves the implicit assumption of the existence of a long enough time scale such that if we estimate the GEV or GPD shape parameter over correspondingly long (or longer) time series, we get stable and robust results. Hence, one can compute return times for events with an infinitesimal probability. This can be misleading for geophysical extreme events as in some cases different dynamical phenomena are relevant at different time scales. Only a dynamical based approach, as the one presented here, prevents from the computation of biased return times. The use of all the sampled dynamical time scales available (so that finer and finer features of the attractor are taken into account) allow for the construction of robust return times for observed phenomena and avoid estimates over longer time scales whose dynamics remains, by definition, unknown.

We have also shown the power of EVT in detecting phase transitions. With respect to the commonly used indicators of the proximity of critical transitions, studying the changeover between different EVLs (in particular, from Weibull to Gumbel and then Frechét distributions when) provide a precise threshold for the warning, linked to the sign change of of the shape parameter of the GEV probability distribution. This is a further motivation to extend these techniques to geophysical extremes. Interestingly, the high sensitivity of extremes in the vicinity of tipping points suggests that observing unusual return times for extreme weather events may provide a flag signalling at an an early stage qualitative changes of atmospheric circulation, which could hardly be predicted by looking at the bulk of the statistics; see also [319].

# 11
# Conclusions

## 11.1
## Main Concepts of This Book

This book originates from the desire to develop a common framework for many closely related mathematical results and ideas linking the theory of extreme events with the theory of dynamical systems, and to show the potential of these concepts for studying complex systems in physical, engineering and social sciences. Such a connection is made possible by the observation that dynamical systems generate infinitely many random process. One can construct a random process by considering the time evolution of an observable, defined as a suitably well-behaved function of the phase space of the dynamical system. Hence, the properties of the random process are defined by a) features of the dynamics of the underlying system; and b) the specific choice of the observable one is considering. It seems useful at this stage to recapitulate some of the main concepts discussed in the book.

Using a dynamical system as generator of stochastic processes is far from trivial when one has the goal of establishing extreme value laws (EVLs). The main difficulty emerges from the fact that extreme value theory (EVT) was originally developed with the goal of finding the distribution function of the maximum of a set of $N$ independent and identically distributed random variables [38], and later extended to the case where such random variables feature a weak - with a specific technical definition - correlation [1]. If such conditions are met, the statistical properties of extremes defined as *block maxima* (BM) are asymptotically equivalent those of extremes taken as *Peaks Over Threshold* (POT), where choosing a very high threshold corresponds, conceptually, to taking the maxima over very large blocks.

By equivalent we mean the following. The statistics of BM is asymptotically described by the Generalised Extreme Value (GEV) distribution, while the statistics of POT is described by the Generalised Pareto Distribution (GPD). The parameters of the two distributions are in a one-to-one correspondence, and, in particular, the *shape parameter*, which determines the qualitative properties of the extremes (finite vs infinite) is the same. In other terms, if we perform a POT analysis, we are able to derive the information we would have obtained using a BM approach, and *vice versa*. It is useful to note that this is the case even if the BM and POT approaches lead to selecting *different candidates for the extremes*. In practical terms, the GPD method is well-known to the more efficient that the BM method when finite time series are considered, the basic reason being that more information is retained in the procedure [3].

The presence of long-term correlations in the stochastic variables lead to the unwanted result that the black maxima are correlated, and a different point of view must be taken on the problem, see, *e.g.*, [54].

The presence of strong short-term correlations has a less serious yet relevant impact on EVLs. When considering the POT method, in this case we find clusters of extreme events, whose occurrence can be modelled as a compound Poisson process, where a large time scale separates the occurrence of a cluster and the detailed structure of the extreme events occurring within each cluster [140]. Extremes selected as BM are, instead, not affected by such short-term correlations. In this case, fitting BM using a GEV model and POTs using a GPD model for the same time series would lead to different estimates for the the shape parameter. The equivalence between GPD and GEV is reestablished by introducing the *extremal index* (EI) [128], which is in most cases (see [320] for a counterexample) equal to the inverse of the average size of the clusters, and to the ratio between the shape parameter estimated using the POT and the BM methods. Declustering techniques have been introduced in order to address these issue in the analysis of data [130, 3].

Obviously, time series of observables of dynamical system are correlated, so that the conditions on the correlations of stochastic processes mentioned above pull back to the properties of the dynamics of the system and to the observable whose evolution generates the random process [73]. Having a dynamics where mixing is sufficiently strong is the most favourable setting for constructing a dynamical theory of extremes. Fortunately, this is, *e.g.*, the case of Axiom A systems [70, 71], which play a central role as useful mathematical models for high dimensional physically relevant systems [266]. This is extremely promising for constructing a mathematical framework suitable for linking extremes and dynamics in physical system, thus going beyond the usual context of statistical inference.

We are left with the important task of choosing suitable observables. A crucial step relies on studying extremes of monotonically decreasing functions of the distance between the orbit of a system and a point on its attractor, such that the (finite or infinite) maximum is obtained when such distance is zero. This framework, proposed by [72], permits a powerful connection between the recurrence properties of a system around a chosen point of its attractor and the possibility of establishing EVLs for the corresponding *distance* observables. One finds that the existence of exponential hitting or return time statistics leads naturally to being able to derive, where the shape parameter of the GPD or GEV distributions is determined by the functional form of the observable.

Interestingly, if one takes the reference point to be an unstable periodic point, and so an islet of regularity amidst chaos, clusters appears as a result of the periodicity of the point, and the extremal index - inverse of the average cluster size - can be related to how unstable the point is (the stronger the instability of the point, the weaker the clustering) [49]. At practical level, this result allows one to create an approximate dynamical framework for interpreting the reason behind clustering of extremes in observed time series [80].

An equally important aspect is that selecting extremes amounts to looking at the neighbourhood of the reference point with microscope, whose lens has the properties defined by the functional form of the observable. The corresponding parameters of the EVLs are directly related to the local geometrical properties. In particular, the shape and location parameters are directly linked to the local dimension of the attractor around the reference point [77, 78]. This implies that by studying the extremes of distance observables, we can learn about the local geometry of the system. Numerical experiments confirm the feasibility of this approach and indicate how to extract important information from the rate of convergence to the EVLs [46, 77]. If the system is exact dimensional, as in the case of the Axiom A systems [269], the local dimension is the same for almost every point of the attractor, and agrees with the Hausdorff dimension, so that we are able to deduce global geometrical properties [44]. Another interesting property is that since the presence of suitable mixing conditions lead to EVLs for the distance observables, one has that if the extremes of distance observables do *not* obey EVLs, then the underling dynamics in *not*, roughly speaking, chaotic. In other terms, extremes can also provide *qualitative* information on the properties of the underlying dynamics [76]. An interesting example of link between extremes and qualitative properties of the dynamics of a system is the possibility of using EVLs of suitable observables for studying tipping points in complex system, along the lines presented in [75] for the case of plane Couette flow. In this case, one can observe that when the approach to the critical transition between turbulent and laminar states obtained by reducing the Reynolds number is accompanied by a clear changes in the probability distribution of the tail of the energy fluctuations. In particular, the tipping point comes together with a changeover from Weibull to Frechét distributions in the extremes (minima, in this case) of the energy.

Distance observables are useful only for a limited class of practical problems, *i.e.*, when recurrences are relevant. Nonetheless, this approach allows also for a rather unexpected shift in the point of view. If one considers a time series, one can see extremes as rarely recurrent values. This construction can be made rigorous by choosing for each value of the time series distance observables of the form described above (*i.e.* such that they have a maximum when the reference value is realized), and investigating the time needed to achieve a good fit of the observables to the EVLs. If for a value of the time series such convergence is not obtained, one can conclude that recurrences around its value are too few, and one can conclude that we have found a true extreme on the time scale given by the length of the time series. An interesting application of this idea on climatological time series is proposed in [80].

In many other cases, one is interested in observables - like energy in a fluid flow - which have a different functional form. These are called *physical observables* [81], and occurrence of an extreme can be geometrically related to visits of the orbit of the evolving point in the phase space to specific regions of the attractor [44]. Making some assumptions on the geometry of the attractor and of the isosurfaces of the observables, one obtains universal results constraining the value of the shape parameter to be negative and to be related in a simple way to the dimension of the attractor on the stable and unstable directions. Therefore, a link can be drawn between the statistics of extremes and fundamental dynamical properties of the underlying system. An important aspect is that when we consider a high dimensional chaotic system, the shape parameter is close to zero, which is suggestive of a very general properties of observables of statistical mechanical systems. Additionally, making some hypotheses, one can adapt Ruelle's response theory [83, 84] in such a way to write explicit formulas for the sensitivity of the parameters descriptive of the EVLs to perturbations to the dynamics [44].

An especially important challenge is understanding the impact on the properties of the extremes resulting from introducing (suitably defined) stochastic forcing to deterministic dynamical systems. This is an extremely interesting mathematical question, as it gives rise to questions such as: under which conditions on the invariant measure and the correlations do the recurrence properties of a stochastically perturbed system converge, in the limit of noise intensity going to zero, to those of the deterministic system without noise? We have approached the problem using different perturbative approaches [141, 142]. The basic idea is that if the deterministic system is uniformly hyperbolic, like in the case of contracting maps, noise has little impact, for the basic reason that the deterministic system is already *mixing enough*. Moreover, the presence of noise washes out the structure of the unstable periodic orbits, which are very important in the deterministic case, so that clusters are absent. The relatively small impact of adding noise to a *sufficiently chaotic* dynamical system can also be seen from the fact that, in general, Ruelle's response theory suggests that centred noise has only a second order effect on the statistical properties of Axiom A dynamical systems [321]. Nonetheless, noise seems to facilitate relating the geometrical properties of the deterministic attractor to the statistics of extremes of the distance observables by removing some of the singular features of the SRB measure [65].

Instead, if we add noise to dynamical systems featuring regular dynamics, we introduce the kind of decay of correlations needed for establishing suitable EVLs. Nonetheless, in numerical simulations convergence to EVLs is very slow if the noise is very weak, so that for practical purposes a relatively strong noise is needed to create the necessary mixing in the trajectories of the flow [79]. Moreover, when considering distance observables, the underlying regular dynamics persists in the sense that clustering of extremes is found. Adding noise also allows for EVLs in other systems, for example contracting maps, whose deterministic version does not feature EVLs [239].

Many challenges remain open when trying to link extremes, recurrence properties, and dynamics of (complex) systems. Realising such a program would greatly improve our ability to predict extreme events, with many important impacts on a variety of fields of science and technology. With the goal of introducing some fascinating future lines of investigations, in the next sections we give to some thoughts and preliminary results regarding the following problems:

- how to study extremes of coarse-grained observables, and how to treat extremes in models of multiscale systems, where variability is present on a large range of scales, and parametrization of unresolved processes are a necessary aspect of the modelling exercise;
- how to treat extremes in the context of non-stationary (deterministic and stochastic) dynamics, taking into account the time-dependent measure induced by the explicit time-dependence of the dynamics;
- how to look at systems whose attractor can be roughly described as the union of the quasi-disjoint parts, the transitions between which are rare and erratic;
- how to study more effectively clusters of extremes, and how to investigate recurrence properties of extremes of functions where the maximizing set is not just a point but a more complex geometrical set.
- how to approach the problem of studying the underpinning properties of spatial and temporal correlations of extreme events in spatially extended random fields.

## 11.2
## Extremes, Coarse Graining, and Parametrizations

In many cases, one is interested in studying spatially extended systems (like a fluid or the climate system), and it is relevant to relate local extremes to spatially averaged ones, in order to capture spatial coherence of such events. The properties of local climate extremes depend critically on the spatial resolution of the numerical model, and spatial averaging operations lead to nontrivial changes in the statistics [27]. Moreover, extended large fluctuations are significant because they can require non-trivial feedbacks and compensating behaviours in other parts of the domain.

In many practical problems, extremes are indeed spatially coherent, and the relationship between such coherence properties and those of the typical fluctuations are not obvious. In a different yet complementary direction, one might be interested in studying long-lasting large fluctuations, where persistence can be key to causing serious impacts. A natural way to try to approach such problems starting from high-resolution (in space and in time) fields is to perform coarse graining, and study the resulting spatial and/or time averaged fields. More specifically, one may be interested in studying the properties of something reading like:

$$Y_{j,k} = \frac{1}{(L+1)(M+s1)} \sum_{l=-L/2}^{L/2} \sum_{m=-M/2}^{M/2} X_{j+l,k+m}, \qquad (11.2.1)$$

where the $X$'s are the dynamically generated stochastic variables of the high-resolution discrete spatio-temporal field, where, *e.g.*, the first index refers to time and the second index to space, while the $Y$'s are the coarse-grained variables. The statistics of the $X$'s might or might not be identical. The perfect setting for such investigations - both for space and temporal averages, when correlations are weak enough, and L and M are large enough - is large deviations theory [322], which provides powerful tools for deriving many important results in, *e.g.*, statistical mechanics [323, 324].

Large deviations are well known to be able to provide information also on the typical behaviour of statistical mechanical systems, similarly to the case of extremes.

Large deviation and extremes are both *large*, so it might seem straightforward to relate them, despite their different definitions. In fact, things are more complex than this, because extremes of averages of variables address the true tail of the distribution of the resulting stochastic variables, whereas large deviations are aimed at describing the most likely among the very large and very unlikely fluctuations of the averaged variables. Therefore, extremes are *more extreme* than large deviations. It is not currently clear whether EVT or large deviations theory is the best approach for studying spatially and/or temporally extended *intense* events, in the sense of providing results useful for applications such as assessing risk. A detailed investigation in this direction is of great scientific and practical urgency.

The issue of coarse graining emerges also in a different context, *i.e.* when we need to model numerically a multiscale system. In simple terms, a multiscale system can be written as:

$$\dot{X} = f_X(X) + g_X(X, Y) \tag{11.2.2}$$

$$\dot{Y} = f_Y(Y) + g_Y(X, Y) \tag{11.2.3}$$

where $X$ are the slow variables and $Y$ are the fast variables, $f$ represent the autonomous dynamics, and $g$ represent the couplings between the fast and slow variables. When trying to construct a numerical model for such a system, it is unavoidable to resort to parametrizations, *i.e.* to construct methods for accounting, at least approximately, for the impact of fast processes related to $Y$ occurring at small spatial scales on the slow variables $X$, often describing large scale features of primary interest. The goal is to construct an effective model for the slow variables *only*. If one can safely assume a vast time-scale separation between the slow and the fast variables, it is possible to use averaging and homogenization methods [325, 326], which have found extensive applications in geophysical fluid dynamics . The basic results are that one can approximately treat the effect of fast modes on the slow dynamics by adding suitably defined deterministic and stochastic (with white spectrum) forcing terms in the evolution equations of the slow variables [327, 328]. While the use of deterministic, *mean field* parametrizations for processes like convection, which cannot yet be captured by the relatively coarse grids of most weather and climate models is a standard practise in geophysical fluid dynamical modelling [329, 330, 24], currently weather and climate modelling centres are moving in the direction of introducing stochastic parametrizations [331, 332], as mounting evidences suggest that they are more effective than usual deterministic methods [333, 334].

It is indeed not clear to what extent such methods, which aim at being able to describe the typical behaviour of the slow variables, perform in terms of providing a good representation of extreme events. Recently Wouters and Lucarini [264, 265] proposed a method for constructing parametrisations based upon the Ruelle response theory, which bypasses the problem of assuming a vast scale separation between the slow and fast variables. We note that the parametrizations derived according to such theory include deterministic and stochastic (with non-white spectrum) correction to the autonomous dynamics of the slow variables, plus a memory term, which introduces non-markovian properties to the dynamics. A mathematically rigorous system-atisation of these results has been recently proposed in [335, 336]. Interestingly, the method proposed in [264, 265] allows for constructing parametrizations which are good for any observable of the slow system. Since in [44] it is shown that the statistics of extremes can reconstucted by calculating some moments of the above-threshold events, we are led to think that using this approach could in fact deal effectively also with such non-typical fluctuations of the system.

## 11.3
## Extremes of Non-Autonomous Dynamical Systems

In many situations of practical interest, we need to estimate the probability of occurrence of extremes in non-autonomous dynamical systems. Time-dependence can be related to the presence of natural periodic phenomena, such as in the case of the seasonal cycle when looking at hot or cold extremes of temperatures (or of energy consumption of heating/cooling buildings), or of slow modulations to the parameters of the system, as in the eponymous case of climate change. We will present two different yet related approaches for studying this problem.

Let us consider a continuous-time dynamical system $\dot{x} = G(x, t)$ on a compact manifold $\mathcal{Y} \subset \mathbb{R}^d$, where $x(t) = \phi(t, t_0)x(t_0)$, with $x(t = t_0) = x_{in} \in \mathcal{Y}$ initial condition and $\phi(t, t_0)$ is defined for all $t \geq t_0$ with $\phi(s, s) = \mathbf{1}$.

The two-time evolution operator $\phi$ generates a two-parameter semi-group. In the autonomous case, the evolution operator generates a one-parameter semigroup, because of time translational invariance, so that $\phi(t, s) = \phi(t - s) \, \forall t \geq s$. In the non-autonomous case, in other terms, there is an *absolute clock*. We want to consider forced and dissipative systems such that with probability one initial conditions in the infinite past are attracted at time $t$ towards $\mathcal{A}(t)$, a time-dependent family of geometrical sets. In more formal terms, we say a family of objects $\cup_{t \in \mathbb{R}} \mathcal{A}(t)$ in the finite-dimensional, complete metric phase space $\mathcal{Y}$ is a pullback attractor for the system $\dot{x} = G(x, t)$ if the following conditions are obeyed:

- $\forall t$, $\mathcal{A}(t)$ is a compact subset of $\mathcal{Y}$ which is covariant with the dynamics, *i.e.* $\phi(s, t)\mathcal{A}(t) = \mathcal{A}(s), s \geq t$.
- $\forall t \lim_{t_0 \to -\infty} d_{\mathcal{Y}}(\phi(t, t_0)B, \mathcal{A}(t)) = 0$ for *a.e.* measurable set $B \subset \mathcal{Y}$.

where $d_{\mathcal{Y}}(P, Q)$ is the Hausdorff semi-distance between the $P \subset \mathcal{Y}$ and $Q \subset \mathcal{Y}$. We have that $d_{\mathcal{Y}}(P, Q) = \sup_{x \in P} d_{\mathcal{Y}}(x, Q)$, with $d_{\mathcal{Y}}(x, Q) = \inf_{y \in Q} d_{\mathcal{Y}}(x, y)$. We have that, in general, $d_{\mathcal{Y}}(P, Q) \neq d_{\mathcal{Y}}(Q, P)$ and $d_{\mathcal{Y}}P, Q = 0 \Rightarrow P \subset Q$. In some cases, the geometrical set $\mathcal{A}(t)$ support useful measures $\mu(\mathrm{d}x)$. These can be obtained as evolution at time $t$ through the Ruelle-Perron-Frobenius operator [337] of the Lebesgue measure supported on $B$ in the infinite past, as from the conditions above. Taking the point of view of the *chaotic hypothesis*, we assume that when considering sufficiently high-dimensional, chaotic and dissipative systems, at all practical levels - *i.e.* when one considers macroscopic observables - the corresponding measure $\mu_t(\mathrm{d}x)$ is of the SRB type. This amounts to the fact that we can construct at all times $t$ a meaningful (time-dependent) physics for the system. Obviously, in the autonomous case, and under suitable conditions - *e.g.* in the case of of Axiom A system - $\mathcal{A}(t) = \Omega$ is the attractor of the system (where the $t-$ dependence is dropped), which supports the SRB invariant measure $\mu(\mathrm{d}x)$ discussed above.

In practical terms, when we want to construct the statistical properties of a numerical model describing a non-autonomous forced and dissipative system, we often follow - sometimes inadvertently - a protocol that mirrors precisely the definitions given above: we start many simulations in the distant past with initial conditions chosen according to an a-priori distribution. After a sufficiently long time, related to the slowest time scale of the system, at each instant the statistical properties of the ensemble of simulations do not not depend anymore on the choice of the initial conditions.

A prominent example of this procedure is given by simulations of past and historical climate conditions performed by the modeling group working in [329, 330, 24], where time-dependent changes in the climate forcings due to changes in greenhouse gases, volcanic eruptions, changes in the solar irradiance and other astronomical effects are taken into account in defining the radiative forcing to the system. Note that future climate projections are *always* performed using as initial conditions the final states of simulations of historical climate conditions, so that the covariance properties of the $\mathcal{A}(t)$ set is maintained.

It is tempting to extend to the time-dependent case the results presented in a) Sec. 8.2 for distance observables and in b) Sec. 8.2.2 for the physical observables. One needs to add as a *caveat* in case a) that it is far from obvious that the reference point $\zeta$ is contained in $\mathcal{A}(t)$ for all $t$.

In the case b), the probability of occurrence of an above-T-threshold event for a given observable $A(x)$ at a time $t$ is the fraction of members of ensembles initialized in the infinite past that finds itself in the region of $\mathcal{A}(t)$ such that $A(x) > T$, having measure $\mu_t(\{x \in \mathcal{Y} : A(x) > T\})$. it is not clear whether one could define a high-enough threshold $T$ such that $\forall t \, \mu_t(\mathbf{1}(A(x) > T))$ is small (so that we are considering genuine extreme events) but is non-vanishing (*i.e.* $\max A(x)|_{\mathcal{A}(t)} > T$), so that one may need to suitably define a time-dependent threshold $T(t)$. As a straightforward example, one may consider the fact that - disregarding the effect of global warming - any feasible definition of year-round extremes (cold and hot) of air surface temperatures in a given location requires considering thresholds reflecting the seasonal cycle. Taking into account such caveats, one can study the extremes of distance and physical observables described above by replacing the invariant measure $\mu(\mathrm{d}x)$ used in Secs. 8.2 and 8.2.2 with its time dependent version $\mu_t(\mathrm{d}x)$ discussed above.

Note that in this case we cannot substitute time and ensemble averages, because time-dependence makes individual trajectories in principle useless for inferring statistical properties. This results in a potentially enormous computational burden when considering high-dimensional non-autonomous complex dynamical systems; see also [338] for a clear discussion of these issues in a geophysical context. Additionally, this clarifies that in this case the BM approach is not feasible, because each block would contain, strictly speaking, information which cannot be merged together with other blocks, as stationarity is lost.

Such a fundamental problem can be eased in the case one considers a system where explicit time dependence in the evolution equations is due to very slow - in comparison to the time it takes for the system to describe the statistical properties of an observable of interest - modulations of some parameters. Taking the so-called *adiabatic approximation*, one assumes that the change in the properties of the extremes is so slow in time that analyzing an individual trajectory is sufficient for capturing such $t-$dependence, with time being introduced as a covariate in the statistical inference procedure [72, 56]. One needs to underline that giving a scientific meaning to such a - common - assumption is possible only in an intuitive, heuristic fashion. In [93] it is shown that when different time scales are present in the forcing, such picture breaks down, because complex effects due to the specific protocol of the forcing appear in the statistics of extremes.

We want now to take advantage of the response theory framework introduced in Sec. 8.3 for addressing the problem of treating extremes in a non-autonomous setting and provide some explicit results. Let's now assume that we can write

$$\dot{x} = G(x, t) = G(x) + \epsilon X(x, t) \qquad (11.3.1)$$

where $|\epsilon X(x, t)| \ll |G(x)| \, \forall t \in \mathbb{R}$ and $\forall x \in \mathcal{Y}$. Under appropriate mild regularity conditions, it is always possible to perform a Schauder decomposition [339] that is $X(x, t) = \sum_{k=1}^{\infty} X_k(x) T_k(t)$. Since we will consider only linear response properties, we restrict our analysis without loss of generality to the case where $G(x, t) = G(x) + \epsilon X(x) T(t)$ and all formula can be extended by linearity to the general case. Following [83, 84], it is possible to extend the response theory in such a way to include the effect of the presence of time-modulations in the forcing. In general, the corrections which are expectation value of an observable $\Psi$ depends explicitly on time. One can construct corrections to the time-independent invariant measure $\mu(\mathrm{d}x)$ of the unforced system introduced in Eq. 8.3.1 so that

$$\langle \Psi \rangle^{\epsilon}(t) = \langle \Psi \rangle_0 + \sum_{j=1}^{\infty} \epsilon^j \langle \Psi \rangle_0^{(j)}(t),$$

where $\langle \Psi \rangle_0^{(j)}(t)$ can be expressed as time-convolution of the Green function given in Eq. 8.3.3 with the time modulation $T(t)$:

$$\langle \Psi \rangle_0^{(j)}(t) = \int_{-\infty}^{\infty} d\tau_1 \ldots \int_{-\infty}^{\infty} d\tau_n G_{\Psi}^{(j)}(\tau_1, \ldots, \tau_n) T(t-\tau_1) \ldots T(t-\tau_j). \quad (11.3.2)$$

Accordingly, we can construct the first order correction to the observables $A_n^T$ introduced in Eq. 8.2.29 as

$$\langle A_n^T \rangle_0^{(1)}(t) = \int d\tau G_{A_n^T}^{(1)}(\tau) T(t-\tau), \qquad (11.3.3)$$

where the Green functions have been introduced in Eqs. 8.3.14 and 8.3.16. Subsequently, we can construct the first order time-dependent corrections to the shape and scale parameters describing the extremes of the $A$ observables following Eqs. 8.3.8 and 8.3.11 as follows:

$$\xi_A^{(1)}(t) = \int d\tau \, G_{\xi_A}^{(1)}(\tau) T(t-\tau) \tag{11.3.4}$$

where, using Eq. 8.3.10, we have

$$G_{\xi_A}^{(1)}(\tau) = \alpha_{\xi,T,n-2,n}^{(1)} G_{A_{n-2}^T}^{(1)}(\tau) + \alpha_{\xi,T,n-1,n}^{(1)} G_{A_{n-1}^T}^{(1)}(\tau) + \alpha_{\xi,T,n,n}^{(1)} G_{A_n^T}^{(1)}(\tau) \tag{11.3.5}$$

and

$$\sigma_A^{(1)}(t) = \int d\tau \, G_{\sigma_A}^{(1)}(\tau) T(t-\tau) \tag{11.3.6}$$

where, using Eq. 8.3.13, we derive

$$G_{\sigma_A}^{(1)}(\tau) = \alpha_{\sigma,T,n-2,n}^{(1)} G_{A_{n-2}^T}^{(1)}(\tau) + \alpha_{\sigma,T,n-1,n}^{(1)} G_{A_{n-1}^T}^{(1)}(\tau) + \alpha_{\sigma,T,n,n}^{(1)} G_{A_n^T}^{(1)}(\tau). \tag{11.3.7}$$

See a discussion of these results in a geophysical context in [263, 102, 340]. These results provide some explicit formulas for studying rigorously extremes in a time-dependent setting.

We would like to point out some recent results [92, 93, 59, 94] which seem to suggest that a mathematically sound treatment of time-dependent extremes in numerical models is within reach. The authors studied a time dependent modification of the so-called Lorenz '84 [341] minimal model of the mid-latitude atmospheric circulation:

$$\dot{x} = -y^2 - z^2 - ax - aF(t) \tag{11.3.8}$$
$$\dot{y} = xy - bxz - y + 1 \tag{11.3.9}$$
$$\dot{z} = xz + bxy - z \tag{11.3.10}$$

where one unit of time corresponds to $5 \, d = 1/73 \, y$, and we use the classical values for the parameters $a = 1/4$ and $b = 4$. Additionally the forcing is defined as $F(t) = F_0(t) + A\sin(\omega t)$, with $F_0 = 9.5$, $\omega = 2\pi/73$ (seasonal cycle), and $F_0(t) = 9.5$ if $t \leq 100 \, y$, while $F_0(t) = 9.5 - 2/(7300)(t - 7300)$ is $t > 7300$ (monotonically decreasing linear ramp starting after $100 \, y$). In this model $X$ represents schematically the intensity of the westerly winds, while $Y$ and $Z$ are related to two modes of meridional heat transport, $F$ is the (baroclinic) forcing, which is larger in winter than in summer, so that a slow decreasing ramp and the seasonal cycle, described by $\omega$, modulate $F$. After running a large ensemble of initial conditions, the authors have been able to construct the pullback attractor, whose supported measure $\mu_t$, after transients have died out, is affected by the ramp and by the periodic forcing, as suggested by the response formulas 11.3.2. Let's consider the observable $A(x, y, z) = z$ for in the year centered around $t = 250y$. The light grey points in Fig. 11.1a) and b) represent extremes in winter and summer conditions, respectively. The structure of the two attractors is clearly dissimilar and the statistics should not be merged.

(a)

(b)

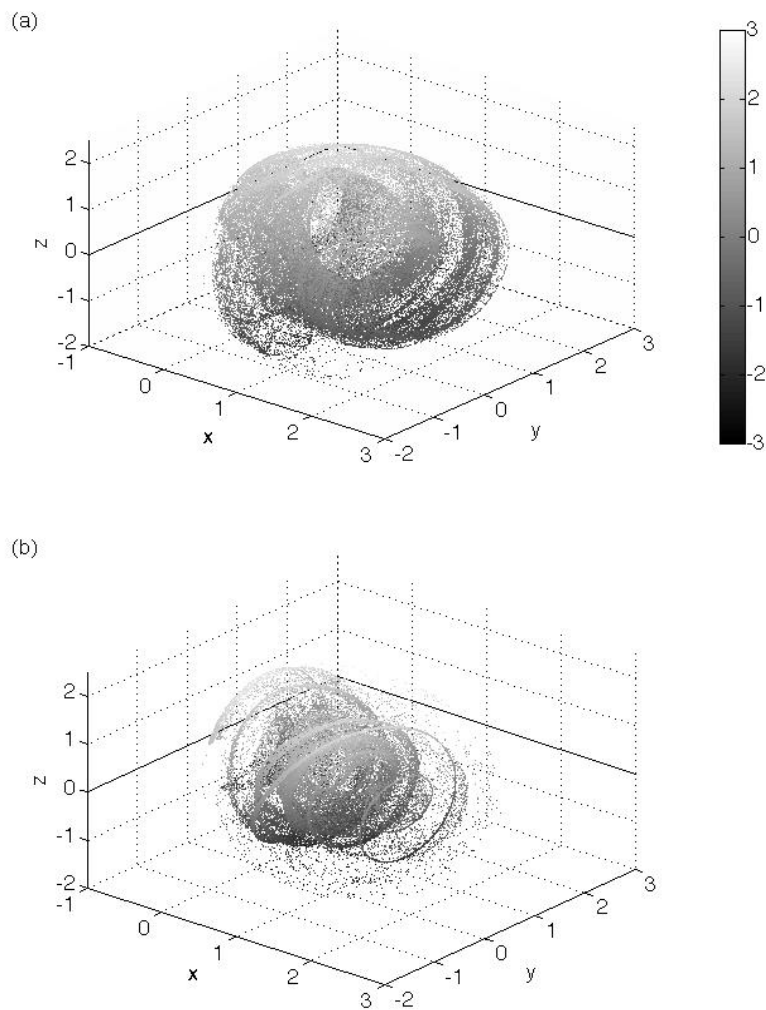

**Figure 11.1** Pullback attractor of the system 11.3.8 for a) winter and b) summer conditions, defined by the phase of the $\sin$ in the definition of $F$, see Eq. 11.3.8. The grey tone bar refers to the values of $Z$. The figure is a courtesy of T. Bódai. See [94] for details.

### 11.3.1
### A note on Randomly Perturbed Dynamical Systems

When considering a random time-dependent forcing, the pullback attractor introduced above becomes a *random pullback attractor*. While the overall construction is similar to the deterministic case, some additional details need to be provided. In particular, roughly speaking, one needs to construct an explicit time parametrization of the probability space defining the randomness by time, and follow the individual random trajectory. This is, in fact, the basic framework which led to the introduction of the very concept of pullback attractor. We then consider a probability space $(\Gamma, \mathcal{F}, \mathbb{P})$, $\mathcal{F}$ is a $\sigma$-algebra of the measurable subsets of $\Omega$, and $\mathbb{P}$ is a probability measure. We consider a measurable map $\theta_t : \Omega \to \Omega$ such that the classic semigroup property is obeyed, that is $\theta_s \circ \theta_t = \theta_{s+t}$ and $\theta_0 = \mathbf{1}$, and such that the probability measure is conserved - a stationarity conditions $\theta_t \mathbb{P} = \mathbb{P}$. The map $\theta$ gives a protocol of invertible transformations of $\Omega$ allowing for tracking the noise. A random evolution on $\mathcal{Y}$ is given by an evolution operator $\phi(t, \omega)x$, where $\omega$ describes the random influence, with the so-called *cocycle* property, that is $\psi(t+s, \omega) = \phi(t, \theta_s \omega) \circ \phi(s, \omega)$. We now construct the extended space $\mathcal{P} \times \mathcal{Y}$ and introduce the evolution operator $\Phi_t : (\omega, x) \to (\theta_t \omega, \phi(t, \omega)x)$, bringing together coherently the evolution of the point $x \in \mathcal{Y}$ and the evolution of the random variable $\omega \in \Gamma$.

Using this definition, we have that $\Phi_t$ is a semi-group, that is $\Phi_{t+s} = \Phi_t \circ \Phi_s$.

In Chapter 7 the following type of random systems is considered: let $T_t$ be a sequence of identically distributed random variables taking values in the measurable space $\Gamma$, where $t$ is a continuous index such that $\text{Cor}_{\mathbb{P}}(T, T, |t-s|)$ vanishes sufficiently fast as $|t-s| \to \infty$.

We associate to each $\gamma \in \Gamma$ a measurable function $\epsilon X(x)_{\gamma_t} : \mathcal{Y} \to \mathcal{Y}$, by $\epsilon X(x)_{\gamma_t} = \epsilon T_{\gamma_t} X(x)$, where $\epsilon$ is a small constant $\in \mathbb{R}$, $X : \mathcal{Y} \to \mathcal{Y}$. Using these definitions, we study:

$$\dot{x} = G(x, t) = G(x) + \epsilon T_{\gamma_t} X(x) \qquad (11.3.11)$$

as a Langevin equation, where we take the Stratonovich convention.

In [89] attention is focussed on the trajectories of the individual ensemble members driven by the same sequence of random variables $T_t$ and constructed according to the time-dependent measure. While mathematically sound, this approach seems a bit unnatural with respect to describing reasonable experimental conditions because we are assuming to be able to reproduce the same sequence of random events for all initial conditions considered.

When studying the statistical properties of the system, it seems more natural to take expectation values over the stochastic variables defining the random forcing. Along these lines, in [321] it is shown that it is possible to study the impact of weak random perturbations to Axiom A dynamical systems $\dot{x} = G(x)$ using the response theory approach. An interesting and general result is that if the stochastic process is centred, *i.e.* the expectation value of $T_{\gamma_t}$ vanishes, the linear response ($\propto \epsilon$) vanishes for any observable, whereas the second order response ($\propto \epsilon^2$) gives the leading order of perturbation. In practice, we expect that adding a moderate noise to an Axiom A (or Axiom A-equivalent) chaotic dynamical system will not alter significantly the parameters $\xi$ and $\sigma$ describing the EVT of both physical and distance observables, because $\xi$ and $\sigma$ can be expressed in terms of suitable observables, see Eqs. 8.2.33-8.2.34. At a more theoretical level, this agrees with the main findings of Chapter 7 referring to mixing systems, even if the point of view on the problem is indeed different. The relatively weak impact of noise on the statistics of extremes for sufficiently chaotic systems has also been observed in [92], while in [59] it is found, in fundamental agreement with the response theory approach, that the impact on the statistics of extremes is enhanced when the characteristic time scale of the noise is close to the characteristic internal time scale of the unperturbed system.

Both in the case of deterministically and random driven system, the challenge is to extend the encouraging results discussed in [92, 93, 59, 94] to high dimensional chaotic systems, *e.g.* when studying climate models data, and find effective ways to beat the curse of dimensionality, while keeping a sound mathematical approach and physical significance in the analysis of the results.

## 11.4
## Quasi-disconnected Attractors

We briefly wish to mention here some additional features which may appear and be extremely relevant at finite time in practical cases, where the dynamics can deviate from Axiom A-like when certain time scales are considered. Let's assume that one can, to first approximation, partition the (unique) attractor of a chaotic dynamical system into, say, two disconnected pieces, so that the system has two time scales, a short one related to the transitive dynamics within each of the two pieces, and a long one corresponding to intermittent transitions from one to the other piece. This corresponds to the scenario of quasi-intransitive dynamics proposed by Lorenz [342] and relevant for many problems of climate dynamics. Note that such setting is similar to that proposed in [75] when studying tipping points, where a parameter is varied until the actual catastrophic change in the attractor structure takes place. We are here considering the case where we are close to the tipping point, but the connection between two regions of the attractor is only ephemeral

In this case, if we observe the system for a time scale intermediate between the short and the long time scale, the properties of the extremes will depend only on the properties of the visited portion of the attractor and we will observe a Weibull distribution, as discussed here, as the dynamics is effectively an Axiom A system in one of the pieces. When our observation time nears the long time scale, we might observe extraordinary large events, corresponding to excursions directed towards the other piece of the attractor, until an irreversible (on the short time scale) transition takes place.

Such extraordinary events will not fit the Weibull law found on smaller time scales, because they result from properties of the attractor, which have not been sampled yet. Therefore, on these intermediate scales, the results proposed here will not be valid. Instead, one may interpret such *extraordinary* events as Dragon Kings [343], which will manifest as outliers spoiling the Weibull statistics and pushing the statistics of extremes towards an (apparently) unphysical Frechét distribution. Observing extremes over even longer time scales, so that the orbit visits many times both parts of the attractor, we shall recover a Weibull law, which reflects the global properties of the attractor. In a system with these properties, small perturbations to the dynamics might impact substantially the long time scale discussed above, resulting in a high sensitivity of the statistics of extremes when a fixed time window of observation is considered.

## 11.5
## Clusters and Recurrence of Extremes

The convergence of Rare Event Point Processes (REPP), proved in [74, 140, 142] can be used to obtain relevant information such as the expected time between the occurrence of catastrophic events, or the intensity of clustering, which ultimately are crucial for assessing risk. However, from a actuarial point of view, not only the frequency of rare undesirable events is relevant for the evaluation of risk associated to certain phenomenon. In fact, insurance companies are also interested on the severity and impact of aggregate damage. This motivates widening the scope of the analysis from counting only tthe number of exceedances (within a cluster and/or withiin a given period of time), to, *e.g.*, computing their aggregate impact by adding their magnitudes, thus studying the so-called Area Over Threshold (AOT). Developing a theory that addresses the convergence and the properties of AOT point processes for time series arising from chaotic dynamical systems is a main goal for future perspectives. The potential of application in actuarial science of such results is huge since they could be used to model and estimate the money losses by claim payments made by insurance companies.

As already hinted at when looking at physical observables, it seems extremely promising to study the extremes of observables whose maximizing sets are not limited to a single point $\zeta$. In trying to study more general sets, the first step would be to consider a finite or countable number of maxima. Some work has already been done in that direction in [139] for typical points chosen independently. Hence, a natural question would be analyze the problem of multiple correlated maxima. One of the interesting consequences of having such general maximising sets is the possibility of having some new recurrence properties for extremes effect associated to the fact that an orbit of the dynamical systems may enter and exit multiple times the neighbourhood of the maximal points, as in the case where various points of the maximising sets belong to the same orbit (or are shadowed by it), with ensuing clustering or short recurrence of extremes. A suitable interplay of structure of the orbits and of shape of the maximising set of the observable might emulate a periodic behaviour leading to clusters of exceedances. One example of that situation is observed in neuroscience, on the records of the activity of a neuron, where the bursts or accumulated spikes may appear in clusters.

The unfolding of such more complex options for clustering of extremes provides a fascinating area of study between EVT and dynamical systems, as one could find new and interesting ways to relate occurrence of extremes and properties of the dynamics. Moreover, these short recurrence mechanisms have a great potential as a source of examples for studying new stochastic processes in the classical EVT context, and, on the other side, for serving as a model for several practical problems in natural, engineering, and social sciences. We believe that through this more complex framework we can also recreate different clustering patterns. So far, the only pattern for clusters of extremes we can construct and explain in detail arises dynamically from the presence of repelling periodic points, as discussed in Chapter 4, and consists of a sequence of decreasing exceedances observed at specific times, as imposed by the period of the repellor. We plan to make an exhaustive classification of the clustering patterns and of the respective geometric/recurrence properties that cause them, in order to be able to categorize clusters and their dynamical causes. In the case of preventing structural failures, the anticipated study of the type of the clustering patterns of a certain natural phenomenon is of crucial importance, on one side to predict the likelihood of such failure and on, another side, to help designing the material and structures to stand stronger against the natural causes they have to face.

## 11.6
## Towards Spatial Extremes: Coupled Map Lattice Models

Coupled map lattice (CML) are often used to study at a qualitative level the chaotic dynamics of spatially extended systems [344]. This includes the dynamics of spatiotemporal chaos where the number of effective degrees of freedom diverges as the size of the system increases. The CML approach has been proven to be effective in studying systems of population dynamics, chemical reactions, convection, fluid flow and biological networks [345]. The idea is that CML could be provide a good framework for extending the application of EVT to infinite-dimensional systems, like those described by partial differential equations.

As Keller suggested in [141], CML given the possibility of account for strong and localised coupling mechanisms able to determine the possibility of observing a broader variety of situations compatible with or producing rare and extreme events. Several techniques presented in this book could be applied to study CML. In particular, we see as a potentially promising development the idea of using the spectral theoretic approach based on the transfer operator and its perturbations on suitable spaces of distributions on the lattice.

It will be fascinating to study whether GEV and/or GPD-like EVLs emerge when studying the extremes of observables closely related to and able to account for how different (how distant, in a suitable defined metrics) two spatially separated lattice points are at a given time. This is a change from the point of view taken in this book, where we consider *temporal* recurrence. In this latter case, we have specified in Chaps. 6 and 8, and used in Chaps. 9 and 10, that the parameters of the GEV and GP distributions are directly linked to the local dimension around the reference point we wish to study the recurrences of. Analyzing, instead, *spatial* recurrence as proposed here, one could investigate whether the physical properties of the systems (cluster size, percolation coefficients [346]) can be related to the EVLs parameters. Finally one could combine together space and time dimensions in a bivariate extreme value theory, thus providing a solid basis for a comprehensive understanding of the spatial and temporal correlation of extremes events in spatially extended random fields.

# A
# Codes

In this appendix we present a series of Extremal Index - EI functions useful to perform the extreme events analysis described in the previous chapters. This selection is not a replacement of an extreme value package for the estimation of the parameter of the GPD or the GEV distribution, but they instead provide a compact and complete suite to compute the relevant quantities introduced for the extreme value analysis of dynamical systems.

The first two functions `extremal_FerroSegers` and `extremal_Sueveges` compute the extremal index $\theta$ introduced in Section 3.3.3 with two different methods. Both these functions are necessary to run the main function `Gumbel_analysis` which computes, for a univariate time series $Y$, EI - the extremal index - and the parameters of the GEV distribution for a set of points $\zeta$ automatically chosen in $[\min(Y) \max(Y)]$. The function is called `Gumbel_analysis` because it assumes that the observable chosen is of type $g_1$ (see Section 4.2.1), but it can be easily modified for observables $g_2$ or $g_3$. The function `Gumbel_analysis2D` works as the previous one but it accepts as input two univariate time series $Y_1$ and $Y_2$, measured simultaneously.

The last script is an example of application for the Bernoulli Shift Map introduced in Section 10.2. It performs the analysis via the previous scripts and plot the results with respect to the points $\zeta$ selected. By changing the definition of the map is possible to study all one-dimensional dynamical systems.

## A.1
## Extremal index

The following two functions compute the extremal index introduced in Chapter 3 with the methods of Ferro-Segers [131] and Süveges [347]. They take as input a univariate time series $Y$ and a quantile $p$. The output is the value $\theta$ of EI.

```
function [theta]=extremal_FerroSegers(Y,p)

% This function computes the extremal index theta by using the
% method proposed by Ferro-Segers (Ferro, C. A. T., and
% J. Segers (2003), Inference for clusters of extremes,
% J. R. Stat. Soc., Ser. B, 65, 545-556.).
```

```
% INPUTS:
% -Y: a vector containing a univariate time series
% -p: a quantile value
% OUTPUTS:
% -theta: the estimate of the extrimal index.

% Extract the threshold u corresponding to the quantile p
 u=quantile(Y, p);

% Compute the exceedances
 Si=find(Y>u);

% Compute the cluster lengths
 Ti=diff(Si);

% Compute the total number of clusters
 N=length(Ti);

% Use the Ferro-Segers formula to extract theta
 theta=2.*(sum(Ti-1)).^2./(N.*sum((Ti-1).*(Ti-2)));
 end

%%%%%%%%%%%%%%
%%%%%%%%%%%%%%
%%%%%%%%%%%%%%
%%%%%%%%%%%%%%
%%%%%%%%%%%%%

function [theta]=extremal_Sueveges(Y,p)

% This function computes the extremal index theta by using the
% method proposed % by  Sueveges, M. (2007).
% Likelihood estimation of the extremal index. Extremes, 10, 41-55.

% INPUTS:
% -Y: a vector containing a univariate time series
% -p: a quantile value
% OUTPUTS:
% -theta: the estimate of the extrimal index.

% Extract the threshold u corresponding to the quantile p
 u=quantile(Y, p);
 q=1-p;

% Compute the exceedances
 Li=find(Y>u);

% Compute the cluster lengths
 Ti=diff(Li);
 Si=Ti-1;
```

```
% Compute the total number of clusters longer than one
Nc=length(find(Si>0));

% Compute the total number of clusters
N=length(Ti);

% Use the Sueveges formula to extract theta
theta=(sum(q.*Si)+N+Nc- sqrt( (sum(q.*Si) +N+Nc).^2 ...
-8*Nc*sum(q.*Si))  )./(2*sum(q.*Si));

 end
```

## A.2
## Recurrences - Extreme Value Analysis

The following two functions perform an extreme value analysis for recurrences presented in Chapter 9 of a univariate time series $Y$ (or bivariate $Y_1, Y_2$ for the second function). $Y$ can be a orbit of a dynamical systems or a series of observations sampled at intervals $\Delta t$. the method selects as much $\zeta$-points between $\min(Y)$ and $\max(Y)$ as much specified by the variable $num\_fich$, and compute the observable $g_1(i,t) = -\log(dist(Y(t), \zeta(i)))$ returning – for each $i$ – the L-moment of the GEV shape parameter $\xi(i)$, the scale parameter $\sigma(i)$ and the location parameter $\mu(i)$. It also computes the extremal index EI with the methods of Sueveges $\theta_1(i)$ and Ferro-Segers $\theta_2(i)$, according the quantile specified as $p$.

```
function [theta1 theta2 Zeta Csi Sigma Mu]=Gumbel_analysis(Y, num_fich, p)

% This function performs the  extreme value analysis for recurrences
% of a univariate time series Y.  The method selects  as much Zeta
% points between min(Y) and max(Y) as much specified by the variable
% num_fich, and compute the observable g1(i,t)=-log(dist(Y(t), Zeta(i)))
% returning for each of them the MLE of the GEV shape parameter Csi,
% the scale parameter Sigma and the location parameter Mu.
% It also computes the extremal index with the methods of Sueveges
% (theta1)% and Ferro-Segers (theta2),
% according the quantile specified as p.

% INPUTS:
% -Y: a vector containing a univariate time series
% -num_fich: the number of points desired as output
% -p: a quantile value

% OUTPUTS:
% -theta1: a vector containing the estimates of the extremal index
%     by using the  Sueveges methodology
% -theta2: a vector containing the estimates of the extremal index
%     by using the Ferro-Segers methodology
% -Zeta: vector of reference points for the series Y
% -Csi: a vector containing the estimates of the shape parameter
```

```
%           for the GEV distribution via MLE techniques
% -Sigma: a vector containing the estimates of the scale parameter
%           for the GEV distribution via MLE techniques
% -Mu: a vector containing the estimates of the location parameter
%           for the GEV distribution via MLE techniques

% Initial check and error messages
 if num_fich <3
    warning('Number of points must be = or greater than 3')
 end
 if p < 0.5
    warning('ple too small to obtain reliable estimates on theta')
 end
 if length(Y) < 10000
    warning('Series too short to perform a reliable analysis')
 end

% Definition of the points zeta, as maxima and minima
 Zeta=ones(1,num_fich);
 for i=1:num_fich
 Zeta(i)=min(Y)+ (i)*(abs(max(Y)-min(Y)))/num_fich;

% Definition of the bin length m and the number of maxima n.
% It can be changed according the needs. Here the number of data
% is taken as the square root of the total length of the series and the
% bin length is adjusted accordingly.
 s=length(Y);
 n=fix(sqrt(s)/2);
 m=fix(s/n)-2;
    clear dmin
    init=1;
    logdista=-log(abs(Y-Zeta(i)));

% Computation of the extremal index according to the two methods
    [theta1(i)]=extremal_Sueveges(logdista,p);
    [theta2(i)]=extremal_FerroSegers(logdista,p);

% Loop to extract the observable g1
 for j=1:n-init
    ddmin(j)=min((abs(Y(((j+init-1)*m+1):(j+init)*m)-Zeta(i)))));
    dmin(j)=-log(ddmin(j));
 end

% Eclusions outliers up to nout, here nout=1. Can be changed
% according to your needs
 nout=1;
 dminS=sort(dmin);
 clear dmin
 dmin=dminS(nout+1:end-nout);

% Fit the GEV distribution via the Matlab function gevfit
```

```
   [tpar tpari]=gevfit(dmin);
    Csi(i)=tpar(1);
    Sigma(i)=tpar(2);
    Mu(i)=tpar(3);
 end
 end
```

The function `Gumbel_analysis2D` is similar to the previous function, but accepts as input two time series $Y_1$ and $Y_2$. It assumed that the two sets of measurements are taken simultaneously.

```
 function  [theta1 theta2 Zeta1 Zeta2 Csi Sigma Mu ...
 ]=Gumbel_analysis2D(Y1, Y2, num_fich, m,  p)
% This function performs the  extreme value analysis for recurrences
% of a bivariate time series Y1, Y2.  The method selects as much Zeta
% points between min(Y) and max(Y) as much specified by the variable
% num_fich, and compute the observable
% g1(i,t)=-log(dist(Y(t), Zeta(i)))
% returning for each of them the MLE of the GEV shape parameter Csi,
% the scale parameter Sigma and the location parameter Mu.
% It also computes the extremal index with the methods of Sueveges
% (theta1)% and Ferro-Segers (theta2),
% according the quantile specified as p.

% INPUTS:
% -Y1: a vector containing a univariate time series
% -Y2: a vector containing a univariate time series
% -num_fich: the number of points desired as output
% -p: a quantile value

% OUTPUTS:
% -theta1: a matrix containing the estimates of the extremal index
%     by using the  Sueveges methodology
% -theta2: a matrix containing the estimates of the extremal index
%     by using the Ferro-Segers methodology
% -Zeta1: vector of reference points for the series Y1
% -Zeta2: vector of reference points for the series Y2
% -Csi: a matrix containing the estimates of the shape parameter
%            for the GEV distribution via MLE techniques
% -Sigma: a matrix containing the estimates of the scale parameter
%            for the GEV distribution via MLE techniques
% -Mu: a matrix containing the estimates of the location parameter
%            for the GEV distribution via MLE techniques

% Initial check and error messages
 if num_fich <3
    warning('Number of points must be = or greater than 3')
 end
 if p < 0.5
    warning('ple too small to obtain reliable estimates on theta')
 end
 if length(Y1) < 10000
```

```
    warning('Series too short to perform a reliable analysis')
end

% Definition of the points zeta, as maxima and minima
Zeta1=ones(1,num_fich);
Zeta2=ones(1,num_fich);
for i=1:num_fich
for j=1:num_fich
Zeta1(i)=min(Y1)+ (i)*(abs(max(Y1)-min(Y1)))/num_fich;
Zeta2(j)=min(Y2)+ (j)*(abs(max(Y2)-min(Y2)))/num_fich;

% Definition of the bin length m and the number of maxima n.
% It can be changed according the needs. Here the number of data
% is taken as the square root of the total length of the series and
% bin length is adjusted accordingly.
s=length(Y1);
n=fix(s./m)-1;
    clear dmin
    init=1;
    logdista=-log(sqrt( (1./std(Y1)).*(Y1-Zeta1(i)).^2 +...
    (1./std(Y2)).*(Y2-Zeta2(j)).^2));

% Computation of the extremal index according to the two methods
[theta1(i,j)]=extremal_Sueveges(logdista,p);
[theta2(i,j)]=extremal_FerroSegers(logdista,p);

% Loop to extract the observable g1
for k=1:n
    ddmin(k)=...
    min(sqrt( (1./std(Y1)).*(Y1(((k+init-1)*m+1):(k+init)*m)-...
    Zeta1(i)).^2 +...
    (1./std(Y2)).*(Y2(((k+init-1)*m+1):(k+init)*m)-Zeta2(j)).^2));
    dmin(k)=-log(ddmin(k));
end

% Eclusions outliers up to nout, here nout=1. Can be changed
% according to your needs
nout=1;
dminS=sort(dmin);
clear dmin
dmin=dminS(nout+1:end-nout);

% Fit the GEV distribution via the Matlab function gevfit
    [tpar tpari]=gevfit(dmin);
     Csi(i,j)=tpar(1);
     Sigma(i,j)=tpar(2);
     Mu(i,j)=tpar(3);

end
end
end
```

## A.3
## Sample Program

This sample program uses the previous functions to compute the recurrences - extreme value analysis for the Bernoulli Shift map $3x \bmod 1$, and can be easily adapted to the study one dimensional dynamical systems.

```
% This sample program compute the recurrences for the Bernoulli Shift
% map 3x mod 1 and plot the results of the parameter estimation.

% I: Compute of the orbit of 100000 iteration for the Bernoulli map 3x mod 1
 Npoints=100000;
% x stores the orbit
 x= ones(1, Npoints);

% Initial condition (irrational or rational  gives different results)
% here we set to irrational as a random number between 0 and 1
 x(1)=rand;

% Computation of the map
 for i=2:Npoints
    x(i)=mod(3*x(i-1),1);
 end

% II: Gumbel analysis
 num_fich=100;
 p=0.95;

 [theta1 theta2 Zeta Csi Sigma Mu]=Gumbel_analysis(x, num_fich, p);

%III: figures
figure(101)

 subplot(2,2,1)
 plot(Zeta,theta1)
 xlabel('\zeta')
 ylabel('\theta')

 subplot(2,2,2)
 plot(Zeta,Csi)
 xlabel('\zeta')
 ylabel('\csi')

 subplot(2,2,3)
 plot(Zeta,Csi)
 xlabel('\zeta')
 ylabel('\sigma')

 subplot(2,2,4)
 plot(Zeta,Csi)
 xlabel('\zeta')
 ylabel('\mu')
```

# Index

# References

**1** Leadbetter, M.R., Lindgren, G., and Rootzén, H. (1983) *Extremes and related properties of random sequences and processes*, Springer Series in Statistics, Springer-Verlag, New York.

**2** Embrechts, P., Klüppelberg, C., and Mikosch, T. (1999) *Modelling Extremal Events for Insurance and Finance*, Springer-Verlag, New York.

**3** Coles, S. (2001) *An introduction to statistical modeling of extreme values*, Springer-Verlag, New York.

**4** Ghil, M., Yiou, P., Hallegatte, S., Malamud, B.D., Naveau, P., Soloviev, A., Friederichs, P., Keilis-Borok, V., Kondrashov, D., Kossobokov, V., Mestre, O., Nicolis, C., Rust, H.W., Shebalin, P., Vrac, M., Witt, A., and Zaliapin, I. (2011) Extreme events: dynamics, statistics and prediction. *Nonlinear Processes in Geophysics*, **18** (3), 295–350, doi:10.5194/npg-18-295-2011. URL `http://www.nonlin-processes-geophys.net/18/295/2011/`.

**5** Katz, R., Parlange, M., and Naveau, P. (2002) Statistics of extremes in hydrology. *Advances in Water Resources*, **25**, 1287–1304.

**6** Pisarenko, V. and Sornette, D. (2003) Characterization of the frequency of extreme earthquake events by the generalized pareto distribution. *Pure and Applied Geophysics*, **160** (12), 2343–2364.

**7** Malevergne, Y., Pisarenko, V., and Sornette, D. (2006) On the power of generalized extreme value (gev) and generalized pareto distribution (gpd) estimators for empirical distributions of stock returns. *Applied Financial Economics*, **16** (3), 271–289.

**8** Smith, A. and Katz, R. (2013) Us billion-dollar weather and climate disasters: data sources, trends, accuracy and biases. *Natural Hazards*, **67**, 387–410.

**9** Gbegbelegbe, S., Chung, U., Shiferaw, B., Msangi, S., and Tesfaye, K. (2014) Quantifying the impact of weather extremes on global food security: A spatial bio-economic approach. *Weather and Climate Extremes*, **4** (0), 96 – 108, doi:http://dx.doi.org/10.1016/j.wace.2014.05.005. URL `http://www.sciencedirect.com/science/article/pii/S2212094714000474`.

**10** Ganguly, A.R., Kodra, E.A., Agrawal, A., Banerjee, A., Boriah, S., Chatterjee, S., Chatterjee, S., Choudhary, A., Das, D., Faghmous, J., Ganguli, P., Ghosh, S., Hayhoe, K., Hays, C., Hendrix, W., Fu, Q., Kawale, J., Kumar, D., Kumar, V., Liao, W., Liess, S., Mawalagedara, R., Mithal, V., Oglesby, R., Salvi, K., Snyder, P.K., Steinhaeuser, K., Wang, D., and Wuebbles, D. (2014) Toward enhanced understanding and projections of climate extremes using physics-guided data mining techniques. *Nonlinear Processes in Geophysics*, **21** (4), 777–795, doi:10.5194/npg-21-777-2014. URL `http://www.nonlin-processes-geophys.net/21/777/2014/`.

**11** Smolka, A. (2006) Natural disasters and the challenge of extreme events: risk management from an insurance perspective. *Philosophical Transactions of the Royal Society of London A: Mathematical, Physical and Engineering Sciences*, **364** (1845), 2147–2165, doi:10.1098/rsta.2006.1818.

**12** Pelinovski, E. and Kharif, C. (2008) *Extreme Ocean Waves*, Springer. URL https://books.google.com/books?id=N7TUHnxKLBYC.

**13** Kharif, C., Pelinovski, E., and Slunyaev, A. (2009) *Rogue Waves in the Ocean*, Springer.

**14** Didenkulova, I. and Pelinovsky, E. (2011) Rogue waves in nonlinear hyperbolic systems (shallow-water framework). *Nonlinearity*, **24** (3), R1. URL http://stacks.iop.org/0951-7715/24/i=3/a=R01.

**15** Berz, G. (2005) Windstorm and storm surges in Europe: loss trends and possible counter-actions from the viewpoint of an international reinsurer. *Philosophical Transactions of the Royal Society of London A: Mathematical, Physical and Engineering Sciences*, **363** (1831), 1431–1440, doi:10.1098/rsta.2005.1577.

**16** Ulbrich, W., Leckebusch, G., and Donat, M.G. (2013) Windstorms, the most costly natural hazard in europe, in *Natural Disasters and Adaptation to Climate Change* (eds S. Boulter, J. Palutikof, D. Karoly, and D. Guitart), Cambridge University Press, pp. 109–120.

**17** Robine, J.M., Cheung, S.L.K., Roy, S.L., Oyen, H.V., Griffiths, C., Michel, J.P., and Herrmann, F.R. (2008) Death toll exceeded 70,000 in europe during the summer of 2003. *Comptes Rendus Biologies*, **331** (2), 171 – 178, doi:http://dx.doi.org/10.1016/j.crvi.2007.12.001. URL http://www.sciencedirect.com/science/article/pii/S1631069107003770.

**18** Battisti, D.S. and Naylor, R.L. (2009) Historical warnings of future food insecurity with unprecedented seasonal heat. *Science*, **323** (5911), 240–244, doi:10.1126/science.1164363. URL http://www.sciencemag.org/content/323/5911/240.abstract.

**19** Bennett, J.E., Blangiardo, M., Fecht, D., Elliott, P., and Ezzati, M. (2014) Vulnerability to the mortality effects of warm temperature in the districts of england and wales. *Nature Clim. Change*, **4** (4), 269–273. URL http://dx.doi.org/10.1038/nclimate2123.

**20** Barriopedro, D., Fischer, E., Luterbacher, J., Trigo, R., and Garcia-Herrera, R. (2011) The hot summer of 2010: Redrawing the temperature record map of europe. *Science*, **332** (6026), 220–224, doi:10.1126/science.1201224. URL http://www.sciencemag.org/content/332/6026/220.abstract.

**21** Luterbacher, J., Dietrich, D., Xoplaki, E., Grosjean, M., and Wanner, H. (2004) European seasonal and annual temperature variability, trends, and extremes since 1500. *Science*, **303** (5663), 1499–1503, doi:10.1126/science.1093877. URL http://www.sciencemag.org/content/303/5663/1499.abstract.

**22** Coumou, D. and Rahmstorf, S. (2012) A decade of weather extremes. *Nature Clim. Change*, **2** (7), 491–496. URL http://dx.doi.org/10.1038/nclimate1452.

**23** Intergovernmental Panel on Climate Change [Eds: C.B. Field et al.] (2012) *Managing the Risks of Extreme Events and Disasters to Advance Climate Change Adaptation. A Special Report of Working Groups I and II of the Intergovernmental Panel on Climate Change*, Cambridge University Press, Cambridge, UK, and New York, USA.

**24** Intergovernmental Panel on Climate Change [Eds.: T. Stocker et al.] (2014) *Climate Change 2013 - The Physical Science Basis IPCC Working Group I Contribution to AR5*, Cambridge University Press, Cambridge, doi:10.1017/cbo9781107415324. URL http://dx.doi.org/10.1017/cbo9781107415324.

**25** Davison, A.C., Padoan, S.A., and Ribatet, M. (2012) Statistical modeling of spatial extremes. *Statistical Science*, **27** (2),

161–186, doi:10.1214/11-STS376. URL `http://dx.doi.org/10.1214/11-STS376`.

26 Cooley, D., Cisewski, J., Erhardt, R.J., Jeon, S., Mannshardt, E., Omolo, B.O., and Sun, Y. (2012) A survey of spatial extremes: Measuring spatial dependence and modeling spatial effects. *REVSTAT*, **10** (1). URL `http://www.ine.pt/revstat/pdf/rs120106.pdf`.

27 Volosciuk, C., Maraun, D., Semenov, V.A., and Park, W. (2014) Extreme precipitation in an atmosphere general circulation model: Impact of horizontal and vertical model resolutions. *Journal of Climate*, **28** (3), 1184–1205, doi:10.1175/JCLI-D-14-00337.1. URL `http://dx.doi.org/10.1175/JCLI-D-14-00337.1`.

28 Adejuwon, J., Leary, N., Barros, V., Burton, I., Kulkarni, J., and Lasco, R. (2012) *Climate Change and Adaptation*, Earthscan: Climate, Taylor & Francis. URL `http://books.google.com/books?id=qW2n7P54jkEC`.

29 Gioli, G., Khan, T., Bisht, S., and Scheffran, J. (2014) Migration as an Adaptation Strategy and its Gendered Implications: A Case Study From the Upper Indus Basin. *Mountain Research and Development*, **34** (3), 255–265, doi:10.1659/MRD-JOURNAL-D-13-00089.1. URL `http://dx.doi.org/10.1659/MRD-JOURNAL-D-13-00089.1`.

30 Gioli, G., Khan, T., and Scheffran, J. (2014) Climatic and environmental change in the Karakoram: making sense of community perceptions and adaptation strategies. *Regional Environmental Change*, **14** (3), 1151–1162, doi:10.1007/s10113-013-0550-3. URL `http://dx.doi.org/10.1007/s10113-013-0550-3`.

31 Kennett, D.J., Breitenbach, S.F.M., Aquino, V.V., Asmerom, Y., Awe, J., Baldini, J.U., Bartlein, P., Culleton, B.J., Ebert, C., Jazwa, C., Macri, M.J., Marwan, N., Polyak, V., Prufer, K.M., Ridley, H.E., Sodemann, H., Winterhalder, B., and Haug, G.H. (2012) Development and Disintegration of Maya Political Systems in Response to Climate

Change. *Science*, **338** (6108), 788–791, doi:10.1126/science.1226299. URL `http://www.sciencemag.org/content/338/6108/788.abstract`.

32 Fang, J.Q. and Liu, G. (1992) Relationship between climatic change and the nomadic southward migrations in eastern asia during historical times. *Climatic Change*, **22**.

33 Hvistendahl, M. (2012) Roots of empire. *Science*, **337** (6102), 1596–1599, doi:10.1126/science.337.6102.1596. URL `http://www.sciencemag.org/content/337/6102/1596.short`.

34 Jorion, P. (1996) *Value at Risk: The New Benchmark for Controlling Market Risk*, Irwin, Chicago.

35 Christoffersen, P.F., Diebold, F.X., and Schuermann, T. (1998) Horizon Problems and Extreme Events in Financial Risk Management, *Center for Financial Institutions Working Papers 98-16*, Wharton School Center for Financial Institutions, University of Pennsylvania. URL `http://ideas.repec.org/p/wop/pennin/98-16.html`.

36 Embrechts, P. (2000) Extreme value theory: Potential and limitations as an integrated risk management tool. *Derivatives Use, Trading & Regulation*, **6**, 449–456.

37 Gencay, R. and Selcuk, F. (2004) Extreme value theory and value-at-risk: Relative performance in emerging markets. *International Journal of Forecasting*, **20** (2), 287 – 303, doi:http://dx.doi.org/10.1016/j.ijforecast.2003.09.005. URL `http://www.sciencedirect.com/science/article/pii/S0169207003001031`, forecasting Economic and Financial Time Series Using Nonlinear Methods.

38 Gnedenko, B. (1943) Sur la distribution limite du terme maximum d'une série aléatoire. *Ann. of Math. (2)*, **44**, 423–453.

39 Kharin, V., Zwiers, F., and Zhang, X. (2005) Intercomparison of near surface temperature and precipitation extremes in AMIP-2 simulations, reanalyses and observations. *Journal of Climate*, **18**, 5201–5223.

40 Pickands III, J. (1975) Statistical inference

using extreme order statistics. *the Annals of Statistics*, pp. 119–131.

41 Balkema, A. and De Haan, L. (1974) Residual life time at great age. *The Annals of Probability*, pp. 792–804.

42 Simiu, E. and Heckert, N. (1996) Extreme wind distribution tails: a'peaks over threshold' approach. *Journal of Structural Engineering*, **122** (5), 539–547.

43 Deidda, R. (2010) A multiple threshold method for fitting the generalized pareto distribution to rainfall time series. *Hydrology and Earth System Sciences*, **14** (12), 2559–2575, doi:10.5194/hess-14-2559-2010. URL http://www.hydrol-earth-syst-sci.net/14/2559/2010/.

44 Lucarini, V., Kuna, T., Faranda, D., and Wouters, J. (2014) Towards a general theory of extremes for observables of chaotic dynamical systems. *Journal of Statistical Physics*, **154**.

45 Eckmann, J. and Ruelle, D. (1992) Fundamental limitations for estimating dimensions and lyapunov exponents in dynamical systems. *Physica B*, **56**, 185–187.

46 Faranda, D., Lucarini, V., Turchetti, G., and Vaienti, S. (2011) Numerical convergence of the block-maxima approach to the generalized extreme value distribution. *J. Stat. Phys.*, **145**, 1156–1180, doi:10.1007/s10955-011-0234-7. URL http://dx.doi.org/10.1007/s10955-011-0234-7.

47 Serinaldi, F. and Kilsby, C.G. (2014) Rainfall extremes: Toward reconciliation after the battle of distributions. *Water Resources Research*, **50** (1), 336–352, doi:10.1002/2013WR014211. URL http://www.ncbi.nlm.nih.gov/pmc/articles/PMC4016761/.

48 Vannitsem, S. (2007) Statistical properties of the temperature maxima in an intermediate order Quasi-Geostrophic model. *Tellus A*, **59** (1), 80–95.

49 Freitas, A.C.M., Freitas, J.M., and Todd, M. (2012) The extremal index, hitting time statistics and periodicity. *Adv. Math.*, **231** (5), 2626 – 2665, doi:10.1016/j.aim.2012.07.029. URL

http://www.sciencedirect.com/science/article/pii/S0001870812002939.

50 Beirlant, J., Goegebeur, Y., Segers, J., Teugels, J., De Waal, D., and Ferro, C. (2006) *Statistics of Extremes: Theory and Applications*, Wiley Series in Probability and Statistics, Wiley. URL https://books.google.com/books?id=jqmRwfG6aloC.

51 Robinson, P.J. (2001) On the definition of a heat wave. *Journal of Applied Meteorology*, **40** (4), 762–775, doi:10.1175/1520-0450(2001)040<0762: OTDOAH>2.0.CO;2. URL http://dx.doi.org/10.1175/1520-0450(2001)040<0762:OTDOAH>2.0.CO;2.

52 Shiferaw, B., Tesfaye, K., Kassie, M., Abate, T., Prasanna, B., and Menkir, A. (2014) Managing vulnerability to drought and enhancing livelihood resilience in sub-saharan africa: Technological, institutional and policy options. *Weather and Climate Extremes*, **3** (0), 67 – 79, doi:http://dx.doi.org/10.1016/j.wace.2014.04.004. URL http://www.sciencedirect.com/science/article/pii/S2212094714000280, high Level Meeting on National Drought Policy.

53 Altmann, E.G. and Kantz, H. (2005) Recurrence time analysis, long-term correlations, and extreme events. *Phys. Rev. E*, **71**, 056 106, doi:10.1103/PhysRevE.71.056106. URL http://link.aps.org/doi/10.1103/PhysRevE.71.056106.

54 Bunde, A., Eichner, J., Kantelhardt, J., and Havlin, S. (2007) Statistics of return intervals and extreme events in long-term correlated time series, in *Nonlinear Dynamics in Geosciences*, Springer New York, pp. 339–367, doi:10.1007/978-0-387-34918-3_19. URL http://dx.doi.org/10.1007/978-0-387-34918-3_19.

55 Nogaj, M., Parey, S., and Dacunha-Castelle, D. (2007) Non-stationary extreme models and a climatic application. *Nonlinear Processes in Geophysics*, **14** (3), 305–316,

doi:10.5194/npg-14-305-2007. URL `http://www.nonlin-processes-geophys.net/14/305/2007/`.

56 Felici, M., Lucarini, V., Speranza, A., and Vitolo, R. (2007) Extreme value statistics of the total energy in an intermediate-complexity model of the midlatitude atmospheric jet. part ii: Trend detection and assessment. *Journal of the Atmospheric Sciences*, **64** (7), 2159–2175, doi:10.1175/JAS4043.1. URL `http://dx.doi.org/10.1175/JAS4043.1`.

57 Smith, R. (1989) Extreme value analysis of environmental time series: an application to trend detection in ground-level ozone. *Statistical Science*, **4** (4), 367–377.

58 de Haan, L., Tank, A.K., and Neves, C. (2014) On tail trend detection: modeling relative risk. *Extremes*, pp. 1–38, doi:10.1007/s10687-014-0207-8. URL `http://dx.doi.org/10.1007/s10687-014-0207-8`.

59 Bódai, T., Károlyi, G., and Tél, T. (2013) Driving a conceptual model climate by different processes: Snapshot attractors and extreme events. *Phys. Rev. E*, **87**, 022 822, doi:10.1103/PhysRevE.87.022822. URL `http://link.aps.org/doi/10.1103/PhysRevE.87.022822`.

60 Goswami, B.N., Venugopal, V., Sengupta, D., Madhusoodanan, M.S., and Xavier, P.K. (2006) Increasing trend of extreme rain events over india in a warming environment. *Science*, **314** (5804), 1442–1445, doi:10.1126/science.1132027. URL `http://www.sciencemag.org/content/314/5804/1442.abstract`.

61 Franzke, C. (2013) A novel method to test for significant trends in extreme values in serially dependent time series. *Geophysical Research Letters*, **40** (7), 1391–1395, doi:10.1002/grl.50301. URL `http://dx.doi.org/10.1002/grl.50301`.

62 Maraun, D. (2013) When will trends in european mean and heavy daily precipitation emerge? *Environmental Research Letters*, **8** (1), 014 004. URL `http://stacks.iop.org/1748-9326/8/i=1/a=014004`.

63 Asadieh, B. and Krakauer, N.Y. (2015) Global trends in extreme precipitation: climate models versus observations. *Hydrology and Earth System Sciences*, **19** (2), 877–891, doi:10.5194/hess-19-877-2015. URL `http://www.hydrol-earth-syst-sci.net/19/877/2015/`.

64 Madsen, H., Lawrence, D., Lang, M., Martinkova, M., and Kjeldsen, T. (2014) Review of trend analysis and climate change projections of extreme precipitation and floods in europe. *Journal of Hydrology*, **519**, 3634 – 3650, doi:http://dx.doi.org/10.1016/j.jhydrol.2014.11.003. URL `http://www.sciencedirect.com/science/article/pii/S0022169414008889`.

65 Faranda, D. and Vaienti, S. (2014) Extreme value laws for dynamical systems under observational noise. *Physica D: Nonlinear Phenomena*, **280-281**, 86 – 94, doi:http://dx.doi.org/10.1016/j.physd.2014.04.011. URL `http://www.sciencedirect.com/science/article/pii/S0167278914000906`.

66 Deidda, R. and Puliga, M. (2006) Sensitivity of goodness-of-fit statistics to rainfall data rounding off. *Physics and Chemistry of the Earth, Parts A/B/C*, **31** (18), 1240 – 1251, doi:http://dx.doi.org/10.1016/j.pce.2006.04.041. URL `http://www.sciencedirect.com/science/article/pii/S1474706506002889`, time Series Analysis in Hydrology.

67 Deidda, R. (2007) An efficient rounding-off rule estimator: Application to daily rainfall time series. *Water Resources Research*, **43**, W12 405, doi:10.1029/2006WR005409.

68 Deidda, R. and Puliga, M. (2009) Performances of some parameter estimators of the generalized pareto distribution over rounded-off samples.

*Physics and Chemistry of the Earth*, **34**, doi:10.1016/j.pce.2008.12.002.

**69** Katok, A. and Hasselblatt, B. (1997) *Introduction to the modern theory of dynamical systems*, Cambridge Univ. Press.

**70** Eckmann, J.P. and Ruelle, D. (1985) Ergodic theory of chaos and strange attractors. *Rev. Mod. Phys.*, **57**, 617–656, doi:10.1103/RevModPhys.57.617. URL `http://link.aps.org/doi/10.1103/RevModPhys.57.617.`

**71** Ruelle, D. (1989) *Chaotic Evolution and Strange Attractors*, Cambridge University Press, Cambridge.

**72** Collet, P. (2001) Statistics of closest return for some non-uniformly hyperbolic systems. *Ergodic Theory Dynam. Systems*, **21** (2), 401–420, doi:10.1017/S0143385701001201. URL `http://dx.doi.org/10.1017/S0143385701001201.`

**73** Freitas, A.C.M. and Freitas, J.M. (2008) On the link between dependence and independence in extreme value theory for dynamical systems. *Statist. Probab. Lett.*, **78** (9), 1088–1093, doi:10.1016/j.spl.2007.11.002. URL `http://dx.doi.org/10.1016/j.spl.2007.11.002.`

**74** Freitas, A.C.M., Freitas, J.M., and Todd, M. (2010) Hitting time statistics and extreme value theory. *Probab. Theory Related Fields*, **147** (3), 675–710, doi:10.1007/s00440-009-0221-y. URL `http://dx.doi.org/10.1007/s00440-009-0221-y.`

**75** Faranda, D., Lucarini, V., Manneville, P., and Wouters, J. (2014) On using extreme values to detect global stability thresholds in multi-stable systems: The case of transitional plane couette flow. *Chaos, Solitons & Fractals*, **64**, 26 – 35, doi:http://dx.doi.org/10.1016/j.chaos.2014.01.008. URL `http://www.sciencedirect.com/science/article/pii/S0960077914000162`, nonequilibrium Statistical Mechanics: Fluctuations and Response.

**76** Faranda, D., Lucarini, V., Turchetti, G., and Vaienti, S. (2012) Generalized extreme value distribution parameters as dynamical indicators of stability. *International Journal of Bifurcation and Chaos*, **22** (11), 1250 276, doi:10.1142/S0218127412502768. URL `http://www.worldscientific.com/doi/abs/10.1142/S0218127412502768.`

**77** Lucarini, V., Faranda, D., Turchetti, G., and Vaienti, S. (2012) Extreme value distribution for singular measures. *Chaos*, **22** (2), 023 135, doi:http://dx.doi.org/10.1063/1.4718935.

**78** Lucarini, V., Faranda, D., and Wouters, J. (2012) Universal behavior of extreme value statistics for selected observables of dynamical systems. *J. Stat. Phys.*, **147** (1), 63–73.

**79** Faranda, D., Freitas, J.M., Lucarini, V., Turchetti, G., and Vaienti, S. (2013) Extreme value statistics for dynamical systems with noise. *Nonlinearity*, **26** (9), 2597. URL `http://stacks.iop.org/0951-7715/26/i=9/a=2597.`

**80** Faranda, D. and Vaienti, S. (2013) A recurrence-based technique for detecting genuine extremes in instrumental temperature records. *Geophysical Research Letters*, **40** (21), 5782–5786, doi:10.1002/2013GL057811. URL `http://dx.doi.org/10.1002/2013GL057811.`

**81** Holland, M.P., Vitolo, R., Rabassa, P., Sterk, A.E., and Broer, H.W. (2012) Extreme value laws in dynamical systems under physical observables. *Physica D: Nonlinear Phenomena*, **241** (5), 497 – 513, doi:10.1016/j.physd.2011.11.005. URL `http://www.sciencedirect.com/science/article/pii/S0167278911003095.`

**82** Felici, M., Lucarini, V., Speranza, A., and Vitolo, R. (2007) Extreme value statistics of the total energy in an intermediate-complexity model of the midlatitude atmospheric jet. part i: Stationary case. *Journal of the Atmospheric Sciences*, **64** (7), 2137–2158, doi:10.1175/JAS3895.1. URL `http://dx.doi.org/10.1175/JAS3895.1.`

**83** Ruelle, D. (1998) Nonequilibrium statistical mechanics near equilibrium:

computing higher-order terms. *Nonlinearity*, **11** (1), 5–18.

84 Ruelle, D. (2009) A review of linear response theory for general differentiable dynamical systems. *Nonlinearity*, **22** (4), 855–870.

85 Towler, E., Rajagopalan, B., Gilleland, E., Summers, R.S., Yates, D., and Katz, R.W. (2010) Modeling hydrologic and water quality extremes in a changing climate: A statistical approach based on extreme value theory. *Water Resources Research*, **46** (11), n/a–n/a, doi:10.1029/2009WR008876. URL `http://dx.doi.org/10.1029/2009WR008876`.

86 Gilleland, E. and Katz, R.W. (2011) A new software to analyze how extremes change over time. *EOS*, **92**, 13–14.

87 Cheng, L., AghaKouchak, A., Gilleland, E., and Katz, R.W. (2014) Non-stationary extreme value analysis in a changing climate. *Climatic Change*, **127** (2), 353–369, doi:10.1007/s10584-014-1254-5. URL `http://dx.doi.org/10.1007/s10584-014-1254-5`.

88 Kloeden, P. and Rasmussen, M. (2011) *Nonautonomous Dynamical Systems*, Mathematical surveys and monographs, American Mathematical Society. URL `https://books.google.com/books?id=ByCCAwAAQBAJ`.

89 Chekroun, M.D., Simonnet, E., and Ghil, M. (2011) Stochastic climate dynamics: Random attractors and time-dependent invariant measures. *Physica D: Nonlinear Phenomena*, **240** (21), 1685 – 1700, doi:http://dx.doi.org/10.1016/j.physd.2011.06.005. URL `http://www.sciencedirect.com/science/article/pii/S016727891100145X`.

90 Carvalho, A.N., Langa, J., and Robinson, J.C. (2013) The pullback attractor, in *Attractors for infinite-dimensional non-autonomous dynamical systems*, *Applied Mathematical Sciences*, vol. 182, Springer New York, pp. 3–22, doi:10.1007/978-1-4614-4581-4_1. URL `http://dx.doi.org/10.1007/978-1-4614-4581-4_1`.

91 Lai, Y. and Tél, T. (2011) *Transient Chaos: Complex Dynamics on Finite Time Scales*, Applied Mathematical Sciences, Springer. URL `https://books.google.com/books?id=rl00oCTsMzEC`.

92 Bódai, T., Károlyi, G., and Tél, T. (2011) A chaotically driven model climate: extreme events and snapshot attractors. *Nonlinear Processes in Geophysics*, **18** (5), 573–580, doi:10.5194/npg-18-573-2011. URL `http://www.nonlin-processes-geophys.net/18/573/2011/`.

93 Bódai, T. and Tél, T. (2012) Annual variability in a conceptual climate model: Snapshot attractors, hysteresis in extreme events, and climate sensitivity. *Chaos: An Interdisciplinary Journal of Nonlinear Science*, **22** (2), 023110, doi:http://dx.doi.org/10.1063/1.3697984. URL `http://scitation.aip.org/content/aip/journal/chaos/22/2/10.1063/1.3697984`.

94 Drótos, G., Bódai, T., and Tél, T. (2015) Probabilistic concepts in a changing climate: a snapshot attractor picture. *Journal of Climate*, doi:10.1175/JCLI-D-14-00459.1. URL `http://dx.doi.org/10.1175/JCLI-D-14-00459.1`.

95 Sterk, A.E., Holland, M.P., Rabassa, P., Broer, H.W., and Vitolo, R. (2012) Predictability of extreme values in geophysical models. *Nonlinear Processes in Geophysics*, **19**, 529–539, doi:10.5194/npg-19-529-2012.

96 Bódai, T. (2014) Predictability of threshold exceedances in dynamical systems. *ArXiv e-prints*.

97 Kalnay, E. (2003) *Atmospheric Modeling, Data Assimilation and Predictability*, Cambridge University Press. URL `http://books.google.com/books?id=Uqc7zC7NULMC`.

98 Kantz, H., Altmann, E., Hallerberg, S., Holstein, D., and Riegert, A. (2006) Dynamical interpretation of extreme events: predictability and predictions. *Extreme events in nature and society*, pp. 69–93.

99 Hallerberg, S. and Kantz, H. (2008) Influence of the event magnitude on the predictability of an extreme event. *Physical Review E*, **77** (1), 11 108.

100 Leadbetter, M.R. and Rootzén, H. (1998) On extreme values in stationary random fields, in *Stochastic processes and related topics*, Birkhäuser Boston, Boston, MA, Trends Math., pp. 275–285.

101 Pelloni, B., Cullen, M., and Lucarini, V. (2013) Mathematics of the Fluid Earth. *Mathematics Today*, **49** (2), 80–83.

102 Lucarini, V., Blender, R., Herbert, C., Ragone, F., Pascale, S., and Wouters, J. (2014) Mathematical and physical ideas for climate science. *Reviews of Geophysics*, **52** (4), 809–859, doi:10.1002/2013RG000446. URL http://dx.doi.org/10.1002/2013RG000446.

103 Kingman, J.F.C. and Taylor, S.J. (1966) *Introduction to measure and probability*, Cambridge University Press, London.

104 Billingsley, P. (1995) *Probability and measure*, Wiley Series in Probability and Mathematical Statistics, John Wiley & Sons Inc., New York, 3rd edn.. A Wiley-Interscience Publication.

105 Lorenz, E.N. (1963) Deterministic nonperiodic flow. *J. Atmos. Sci.*, **20**, 130–141.

106 Walters, P. (1982) *An introduction to ergodic theory*, Graduate Texts in Mathematics, vol. 79, Springer-Verlag, New York.

107 Galambos, J. (1978) *The asymptotic theory of extreme order statistics*, John Wiley & Sons, New York-Chichester-Brisbane. Wiley Series in Probability and Mathematical Statistics.

108 Resnick, S.I. (1987) *Extreme values, regular variation, and point processes*, Applied Probability. A Series of the Applied Probability Trust, vol. 4, Springer-Verlag, New York, doi:10.1007/978-0-387-75953-1. URL http://dx.doi.org/10.1007/978-0-387-75953-1.

109 Beirlant, J., Goegebeur, Y., Teugels, J., and Segers, J. (2004) *Statistics of extremes*, Wiley Series in Probability and Statistics, John Wiley & Sons, Ltd., Chichester, doi:10.1002/0470012382.

URL http://dx.doi.org/10.1002/0470012382, theory and applications, With contributions from Daniel De Waal and Chris Ferro.

110 Falk, M., Hüsler, J., and Reiss, R.D. (2011) *Laws of small numbers: extremes and rare events*, Birkhäuser/Springer Basel AG, Basel, extended edn., doi:10.1007/978-3-0348-0009-9. URL http://dx.doi.org/10.1007/978-3-0348-0009-9.

111 Fisher, A. and Tippett, L. (1928) Limiting forms of the frequency distribution of the largest or smallest member of a sample. *Math. Proc. Cambridge Philos. Soc.*, **24** (2), 180–190, doi:http://dx.doi.org/10.1017/S0305004100015681.

112 Hosking, J., Walls, J., and Wood, E. (1985) Estimation of the generalized extreme-value distribution by the method of probability-weighted moments. *Technometrics*, **27**, 251–261.

113 Diebolt, J., Guillou, A., Naveau, P., and Ribereau, P. (2008) Improving probability-weighted moment methods for the generalized extreme value distribution. *REVSTAT–Statistical Journal*, **6** (1), 33–50.

114 Hosking, J. (1990) L-moments: analysis and estimation of distributions using linear combinations of order statistics. *Journal of the Royal Statistical Society. Series B (Methodological)*, **52** (1), 105–124.

115 Scarrott, C. and MacDonald, A. (2012) A review of extreme value threshold estimation and uncertainty quantification. *REVSTAT–Statistical Journal*, **10** (1), 33–50.

116 Hill, B. (1975) A simple general approach to inference about the tail of a distribution. *The Annals of Statistics*, **3** (5), 1163–1174.

117 Gomes, M., de Haan, L., and Rodrigues, L. (2008) Tail index estimation for heavy-tailed models: accommodation of bias in weighted log-excesses. *J. R. Stat. Soc. Ser. B Stat. Methodol.*, **10**, 31–52.

118 de Haan, L. and Rootzen, H. (1993) On the estimation of high quantiles. *J. Statist. Plann. Inference*, **35**, 1–13.

119 Beirlant, J., Dierckx, G., Goegebeur, Y., and Matthys, G. (1999) Tail index estimation and an exponential regression

model. *Extremes*, **2**, 177–200.

**120** Feuerverger, A. and Hall, P. (1999) Estimating a tail exponent by modelling departure from a pareto distribution. *Ann. Stat.*, **27**, 760–781.

**121** Leadbetter, M.R. (1973/74) On extreme values in stationary sequences. *Z. Wahrscheinlichkeitstheorie und Verw. Gebiete*, **28**, 289–303.

**122** Leadbetter, M.R. (1983) Extremes and local dependence in stationary sequences. *Z. Wahrsch. Verw. Gebiete*, **65** (2), 291–306, doi:10.1007/BF00532484. URL `http://dx.doi.org/10.1007/BF00532484`.

**123** Albeverio, S. and Jentsch, V. and Kantz, H. (Eds.) (2005) *Extreme Events in Nature and Society*, John Wiley & Sons, New York-Chichester-Brisbane.

**124** Loynes, R.M. (1965) Extreme values in uniformly mixing stationary stochastic processes. *Ann. Math. Statist.*, **36**, 993–999.

**125** Leadbetter, M.R. and Nandagopalan, S. (1989) On exceedance point processes for stationary sequences under mild oscillation restrictions, in *Extreme value theory (Oberwolfach, 1987)*, *Lecture Notes in Statist.*, vol. 51, Springer, New York, pp. 69–80.

**126** Chernick, M.R., Hsing, T., and McCormick, W.P. (1991) Calculating the extremal index for a class of stationary sequences. *Adv. in Appl. Probab.*, **23** (4), 835–850, doi:10.2307/1427679. URL `http://dx.doi.org/10.2307/1427679`.

**127** Kallenberg, O. (1986) *Random measures*, Akademie-Verlag, Berlin, 4th edn..

**128** Hsing, T., Hüsler, J., and Leadbetter, M.R. (1988) On the exceedance point process for a stationary sequence. *Probab. Theory Related Fields*, **78** (1), 97–112, doi:10.1007/BF00718038. URL `http://dx.doi.org/10.1007/BF00718038`.

**129** Haydn, N. and Vaienti, S. (2009) The compound Poisson distribution and return times in dynamical systems. *Probab. Theory Related Fields*, **144** (3-4), 517–542, doi:10.1007/s00440-008-0153-y. URL `http://dx.doi.org/10.1007/`

`s00440-008-0153-y`.

**130** Smith, R.L. and Weissman, I. (1994) Estimating the extremal index. *J. Roy. Statist. Soc. Ser. B*, **56** (3), 515–528. URL `http://links.jstor.org/sici?sici=0035-9246(1994)56:3<515:ETEI>2.0.CO;2-5&origin=MSN`.

**131** Ferro, C.A.T. and Segers, J. (2003) Inference for clusters of extreme values. *J. R. Stat. Soc. Ser B Stat. Methodol.*, **65** (2), 545–556, doi:10.1111/1467-9868.00401. URL `http://dx.doi.org/10.1111/1467-9868.00401`.

**132** Robert, C.Y., Segers, J., and Ferro, C.A.T. (2009) A sliding blocks estimator for the extremal index. *Electron. J. Stat.*, **3**, 993–1020, doi:10.1214/08-EJS345. URL `http://dx.doi.org/10.1214/08-EJS345`.

**133** Robert, C.Y. (2013) Automatic declustering of rare events. *Biometrika*, **100** (3), 587–606, doi:10.1093/biomet/ast013. URL `http://dx.doi.org/10.1093/biomet/ast013`.

**134** Freitas, A.C.M. and Freitas, J.M. (2008) Extreme values for Benedicks-Carleson quadratic maps. *Ergodic Theory Dynam. Systems*, **28** (4), 1117–1133, doi:10.1017/S0143385707000624. URL `http://dx.doi.org/10.1017/S0143385707000624`.

**135** Vitolo, R., Holland, M.P., and Ferro, C.A.T. (2009) Robust extremes in chaotic deterministic systems. *Chaos*, **19** (4), 043127, doi:10.1063/1.3270389. URL `http://link.aip.org/link/?CHA/19/043127/1`.

**136** Freitas, A.C.M., Freitas, J.M., and Todd, M. (2011) Extreme value laws in dynamical systems for non-smooth observations. *J. Stat. Phys.*, **142** (1), 108–126, doi:10.1007/s10955-010-0096-4. URL `http://dx.doi.org/10.1007/s10955-010-0096-4`.

**137** Gupta, C. (2010) Extreme-value distributions for some classes of non-uniformly partially hyperbolic dynamical systems. *Ergodic Theory and Dynamical Systems*, **30** (03), 757–771.

**138** Gupta, C., Holland, M., and Nicol, M.

(2011) Extreme value theory and return time statistics for dispersing billiard maps and flows, Lozi maps and Lorenz-like maps. *Ergodic Theory Dynam. Systems*, **31** (5), 1363–1390, doi:10.1017/S014338571000057X. URL http://dx.doi.org/10.1017/S014338571000057X.

139 Holland, M., Nicol, M., and Török, A. (2012) Extreme value theory for non-uniformly expanding dynamical systems. *Trans. Amer. Math. Soc.*, **364**, 661–688. URL http://dx.doi.org/10.1090/S0002-9947-2011-05271-2.

140 Freitas, A.C.M., Freitas, J.M., and Todd, M. (2013) The compound Poisson limit ruling periodic extreme behaviour of non-uniformly hyperbolic dynamics. *Comm. Math. Phys.*, **321** (2), 483–527, doi:10.1007/s00220-013-1695-0. URL http://dx.doi.org/10.1007/s00220-013-1695-0.

141 Keller, G. (2012) Rare events, exponential hitting times and extremal indices via spectral perturbation. *Dynamical Systems*, **27** (1), 11–27, doi:10.1080/14689367.2011.653329. URL http://www.tandfonline.com/doi/abs/10.1080/14689367.2011.653329.

142 Aytaç, H., Freitas, J.M., and Vaienti, S. (2014) Laws of rare events for deterministic and random dynamical systems. *Trans. Amer. Math.*, pp. doi: 10.1090/S0002–9947–2014–06 300–9, doi:10.1090/S0002-9947-2014-06300-9. URL http://dx.doi.org/10.1090/S0002-9947-2014-06300-9.

143 Freitas, A.C.M., Freitas, J.M., and Todd, M. (2015) Speed of convergence for laws of rare events and escape rates. *Stochastic Process. Appl.*, **125** (4), 1653–1687, doi:10.1016/j.spa.2014.11.011. URL http://dx.doi.org/10.1016/j.spa.2014.11.011.

144 Holland, M. and Nicol, M. (2014), Speed of convergence to an extreme value distribution for non-uniformly hyperbolic dynamical systems, Preprint.

145 Bandt, C. (2007) Random fractals, in *Physics and theoretical computer science:*

*from numbers and languages to (quantum) cryptography security* (eds J. Gazeau, J. Nevšetřil, and B. Rovan), IOS Press, Amsterdam, pp. 91–112.

146 Ferguson, A. and Pollicott, M. (2012) Escape rates for Gibbs measures. *Ergodic Theory Dynam. Systems*, **32** (3), 961–988. URL http://dx.doi.org/10.1017/S0143385711000058.

147 Hirata, M. (1993) Poisson law for Axiom A diffeomorphisms. *Ergodic Theory Dynam. Systems*, **13** (3), 533–556, doi:10.1017/S0143385700007513. URL http://dx.doi.org/10.1017/S0143385700007513.

148 Keller, G. and Liverani, C. (2009) Rare events, escape rates and quasistationarity: some exact formulae. *J. Stat. Phys.*, **135** (3), 519–534, doi:10.1007/s10955-009-9747-8. URL http://dx.doi.org/10.1007/s10955-009-9747-8.

149 Kifer, Y. and Rapaport, A. (2012), Poisson and compound Poisson approximations in a nonconventional setup, Preprint arXiv:1211.5238.

150 Dolgopyat, D. (2004) Limit theorems for partially hyperbolic systems. *Trans. Amer. Math. Soc.*, **356** (4), 1637–1689 (electronic), doi:10.1090/S0002-9947-03-03335-X. URL http://dx.doi.org/10.1090/S0002-9947-03-03335-X.

151 Denker, M., Gordin, M., and Sharova, A. (2004) A Poisson limit theorem for toral automorphisms. *Illinois J. Math.*, **48** (1), 1–20. URL http://projecteuclid.org/getRecord?id=euclid.ijm/1258136170.

152 Abadi, M. and Vergne, N. (2009) Sharp error terms for return time statistics under mixing conditions. *J. Theoret. Probab.*, **22** (1), 18–37, doi:10.1007/s10959-008-0199-x. URL http://dx.doi.org/10.1007/s10959-008-0199-x.

153 Alves, J.F., Freitas, J.M., Luzzatto, S., and Vaienti, S. (2011) From rates of mixing to recurrence times via large deviations. *Adv. Math.*, **228** (2), 1203–1236, doi:10.1016/j.aim.2011.06.014. URL http://dx.doi.org/10.1016/j.aim.2011.06.014.

**154** Freitas, A.C.M., Freitas, J.M., and Rodrigues, F.B., The speed of convergence of rare events point processes in non-uniformly hyperbolic systems, in preparation.

**155** Bradley, R.C. (2005) Basic properties of strong mixing conditions. A survey and some open questions. *Probab. Surv.*, **2**, 107–144 (electronic), doi:10.1214/154957805100000104. URL `http://dx.doi.org/10.1214/15495780510000104`, update of, and a supplement to, the 1986 original.

**156** Young, L. (1999) Recurrence times and rates of mixing. *Israel Journal of Mathematics*, **110** (1), 153–188.

**157** Young, L. (1998) Statistical properties of dynamical systems with some hyperbolicity. *Annals of Mathematics*, **147** (3), 585–650.

**158** Saussol, B. (2000) Absolutely continuous invariant measures for multidimensional expanding maps. *Israel J. Math.*, **116**, 223–248, doi:10.1007/BF02773219. URL `http://dx.doi.org/10.1007/BF02773219`.

**159** Ledrappier, F. (1981) Some properties of absolutely continuous invariant measures on an interval. *Ergodic Theory Dynamical Systems*, **1** (1), 77–93.

**160** Rychlik, M. (1983) Bounded variation and invariant measures. *Studia Math.*, **76** (1), 69–80.

**161** Young, L. (2002) What are SRB measures, and which dynamical systems have them? *Journal of Statistical Physics*, **108**, 733–754.

**162** Sterk, A., Vitolo, R., Broer, H., Simó, C., and Dijkstra, H. (2010) New nonlinear mechanisms of midlatitude atmospheric low-frequency variability. *Physica D: Nonlinear Phenomena*, **239** (10), 702 – 718, doi:http://dx.doi.org/10.1016/j.physd.2010.02.003. URL `http://www.sciencedirect.com/science/article/pii/S0167278910000527`.

**163** Chazottes, J.R. and Collet, P. (2013) Poisson approximation for the number of visits to balls in non-uniformly hyperbolic dynamical systems. *Ergodic Theory Dynam. Systems*, **33** (1), 49–80.

**164** Pitskel', B. (1991) Poisson limit law for Markov chains. *Ergodic Theory Dynam. Systems*, **11** (3), 501–513, doi:10.1017/S0143385700006301. URL `http://dx.doi.org/10.1017/S0143385700006301`.

**165** Galves, A. and Schmitt, B. (1990) Occurrence times of rare events for mixing dynamical systems. *Ann. Inst. H. Poincaré Phys. Théor.*, **52** (3), 267–281. URL `http://www.numdam.org/item?id=AIHPA_1990__52_3_267_0`.

**166** Collet, P., Galves, A., and Schmitt, B. (1992) Unpredictability of the occurrence time of a long laminar period in a model of temporal intermittency. *Ann. Inst. H. Poincaré Phys. Théor.*, **57** (3), 319–331. URL `http://www.numdam.org/item?id=AIHPA_1992__57_3_319_0`.

**167** Einsiedler, M. and Ward, T. (2011) *Ergodic theory with a view towards number theory, Graduate Texts in Mathematics*, vol. 259, Springer-Verlag London Ltd., London, doi:10.1007/978-0-85729-021-2. URL `http://dx.doi.org/10.1007/978-0-85729-021-2`.

**168** Haydn, N., Lacroix, Y., and Vaienti, S. (2005) Hitting and return times in ergodic dynamical systems. *Ann. Probab.*, **33** (5), 2043–2050, doi:10.1214/009117905000000242. URL `http://dx.doi.org/10.1214/009117905000000242`.

**169** Lacroix, Y. (2002) Possible limit laws for entrance times of an ergodic aperiodic dynamical system. *Israel J. Math.*, **132**, 253–263, doi:10.1007/BF02784515. URL `http://dx.doi.org/10.1007/BF02784515`.

**170** Kupsa, M. and Lacroix, Y. (2005) Asymptotics for hitting times. *Ann. Probab.*, **33** (2), 610–619, doi:10.1214/009117904000000883. URL `http://dx.doi.org/10.1214/009117904000000883`.

**171** Collet, P. and Eckmann, J.P. (2006) *Concepts and results in chaotic dynamics: a short course*, Theoretical and Mathematical Physics, Springer-Verlag, Berlin.

**172** Coelho, Z. (2000) Asymptotic laws for symbolic dynamical systems, in *Topics in symbolic dynamics and applications (Temuco, 1997)*, London Math. Soc. Lecture Note Ser., vol. 279, Cambridge Univ. Press, Cambridge, pp. 123–165.

**173** Abadi, M. and Galves, A. (2001) Inequalities for the occurrence times of rare events in mixing processes. The state of the art. *Markov Process. Related Fields*, **7** (1), 97–112. Inhomogeneous random systems (Cergy-Pontoise, 2000).

**174** Sinai, Y.G. (1972) Gibbs measures in ergodic theory. *Uspehi Mat. Nauk*, **27** (4(166)), 21–64.

**175** Ruelle, D. (1973) Statistical mechanics on a compact set with $Z^v$ action satisfying expansiveness and specification. *Trans. Amer. Math. Soc.*, **187**, 237–251.

**176** Bowen, R. (2008) *Equilibrium states and the ergodic theory of Anosov diffeomorphisms*, *Lecture Notes in Mathematics*, vol. 470, Springer-Verlag, Berlin, revised edn.. With a preface by David Ruelle, Edited by Jean-René Chazottes.

**177** Galves, A. and Schmitt, B. (1997) Inequalities for hitting times in mixing dynamical systems. *Random Comput. Dynam.*, **5** (4), 337–347.

**178** Collet, P. and Galves, A. (1995) Asymptotic distribution of entrance times for expanding maps of the interval, in *Dynamical systems and applications*, *World Sci. Ser. Appl. Anal.*, vol. 4, World Sci. Publ., River Edge, NJ, pp. 139–152, doi:10.1142/9789812796417_0011. URL http://dx.doi.org/10.1142/9789812796417_0011.

**179** Abadi, M. and Vergne, N. (2008) Sharp errors for point-wise Poisson approximations in mixing processes. *Nonlinearity*, **21** (12), 2871–2885, doi:10.1088/0951-7715/21/12/008. URL http://dx.doi.org/10.1088/0951-7715/21/12/008.

**180** Abadi, M. and Saussol, B. (2011) Hitting and returning to rare events for all alpha-mixing processes. *Stochastic Process. Appl.*, **121** (2), 314–323, doi:10.1016/j.spa.2010.11.001. URL http://dx.doi.org/10.1016/j.spa.2010.11.001.

**181** Paccaut, F. (2000) Statistics of return times for weighted maps of the interval. *Ann. Inst. H. Poincaré Probab. Statist.*, **36** (3), 339–366.

**182** Robinson, C. (1999) *Dynamical systems*, Studies in Advanced Mathematics, CRC Press, Boca Raton, FL, 2nd edn.. Stability, symbolic dynamics, and chaos.

**183** Haydn, N.T., Winterberg, N., and Zweimüller, R. (2014) Return-time statistics, hitting-time statistics and inducing, in *Ergodic Theory, Open Dynamics, and Coherent Structures*, *Springer Proceedings in Mathematics & Statistics*, vol. 70 (eds W. Bahsoun, C. Bose and G. Froyland), Springer New York, pp. 217–227, doi:10.1007/978-1-4939-0419-8_10. URL http://dx.doi.org/10.1007/978-1-4939-0419-8_10.

**184** Bruin, H., Saussol, B., Troubetzkoy, S., and Vaienti, S. (2003) Return time statistics via inducing. *Ergodic Theory Dynam. Systems*, **23** (4), 991–1013, doi:10.1017/S0143385703000026. URL http://dx.doi.org/10.1017/S0143385703000026.

**185** Pomeau, Y. and Manneville, P. (1980) Intermittent transition to turbulence in dissipative dynamical systems. *Comm. Math. Phys.*, **74** (2), 189–197.

**186** Liverani, C., Saussol, B., and Vaienti, S. (1999) A probabilistic approach to intermittency. *Ergodic Theory Dynam. Systems*, **19** (3), 671–685, doi:10.1017/S0143385799133856. URL http://dx.doi.org/10.1017/S0143385799133856.

**187** Collet, P. and Galves, A. (1993) Statistics of close visits to the indifferent fixed point of an interval map. *J. Statist. Phys.*, **72** (3-4), 459–478.

**188** Hirata, M., Saussol, B., and Vaienti, S. (1999) Statistics of return times: a general framework and new applications. *Comm. Math. Phys.*, **206** (1), 33–55, doi:10.1007/s002200050697. URL http://dx.doi.org/10.1007/s002200050697.

**189** Lyubich, M. (1997) Dynamics of quadratic polynomials. I, II. *Acta Math.*, **178** (2), 185–247, 247–297, doi:10.1007/BF02392694. URL

http://dx.doi.org/10.1007/BF02392694.

190 Graczyk, J. and Świątek, G. (1997) Generic hyperbolicity in the logistic family. *Ann. of Math. (2)*, **146** (1), 1–52, doi:10.2307/2951831. URL http://dx.doi.org/10.2307/2951831.

191 Jakobson, M.V. (1981) Absolutely continuous invariant measures for one-parameter families of one-dimensional maps. *Communications in Mathematical Physics*, **81**, 39–88, doi:10.1007/BF01941800.

192 Misiurewicz, M. (1981) Absolutely continuous measures for certain maps of an interval. *Inst. Hautes Études Sci. Publ. Math.*, **53**, 17–51. URL http://www.numdam.org/item?id=PMIHES_1981__53__17_0.

193 Collet, P. and Eckmann, J.P. (1983) Positive Liapunov exponents and absolute continuity for maps of the interval. *Ergodic Theory Dynam. Systems*, **3** (1), 13–46, doi:10.1017/S0143385700001802. URL http://dx.doi.org/10.1017/S0143385700001802.

194 Bruin, H. and Vaienti, S. (2003) Return time statistics for unimodal maps. *Fund. Math.*, **176** (1), 77–94, doi:10.4064/fm176-1-6. URL http://dx.doi.org/10.4064/fm176-1-6.

195 Bruin, H. and Todd, M. (2009) Return time statistics of invariant measures for interval maps with positive Lyapunov exponent. *Stoch. Dyn.*, **9** (1), 81–100, doi:10.1142/S0219493709002567. URL http://dx.doi.org/10.1142/S0219493709002567.

196 Benedicks, M. and Carleson, L. (1991) The dynamics of the Hénon map. *Ann. of Math. (2)*, **133** (1), 73–169, doi:10.2307/2944326. URL http://dx.doi.org/10.2307/2944326.

197 Coelho, Z. and de Faria, E. (1996) Limit laws of entrance times for homeomorphisms of the circle. *Israel J. Math.*, **93**, 93–112, doi:10.1007/BF02761095. URL http://dx.doi.org/10.1007/BF02761095.

198 Coelho, Z. (2004) The loss of tightness of time distributions for homeomorphisms of the circle. *Trans. Amer. Math. Soc.*, **356** (11), 4427–4445 (electronic), doi:10.1090/S0002-9947-04-03386-0. URL http://dx.doi.org/10.1090/S0002-9947-04-03386-0.

199 Durand, F. and Maass, A. (2001) Limit laws of entrance times for low-complexity Cantor minimal systems. *Nonlinearity*, **14** (4), 683–700.

200 Chaumoître, V. and Kupsa, M. (2005) Asymptotics for return times of rank-one systems. *Stoch. Dyn.*, **5** (1), 65–73, doi:10.1142/S0219493705001298. URL http://dx.doi.org/10.1142/S0219493705001298.

201 Grzegorek, P. and Kupsa, M. (2012) Exponential return times in a zero-entropy process. *Commun. Pure Appl. Anal.*, **11** (3), 1361–1383.

202 Saussol, B. (2009) An introduction to quantitative Poincaré recurrence in dynamical systems. *Rev. Math. Phys.*, **21** (8), 949–979, doi:10.1142/S0129055X09003785. URL http://dx.doi.org/10.1142/S0129055X09003785.

203 Collet, P. (1996) Some ergodic properties of maps of the interval, in *Dynamical systems (Temuco, 1991/1992)*, *Travaux en Cours*, vol. 52, Hermann, Paris, pp. 55–91.

204 Abadi, M. (2004) Sharp error terms and necessary conditions for exponential hitting times in mixing processes. *Ann. Probab.*, **32** (1A), 243–264, doi:10.1214/aop/1078415835. URL http://dx.doi.org/10.1214/aop/1078415835.

205 Abadi, M. (2006) Hitting, returning and the short correlation function. *Bull. Braz. Math. Soc. (N.S.)*, **37** (4), 593–609, doi:10.1007/s00574-006-0030-1. URL http://dx.doi.org/10.1007/s00574-006-0030-1.

206 Abadi, M. and Vaienti, S. (2008) Large deviations for short recurrence. *Discrete Contin. Dyn. Syst.*, **21** (3), 729–747.

207 Freitas, J.M., Haydn, N., and Nicol, M. (2014) Convergence of rare event point processes to the Poisson process for planar billiards. *Nonlinearity*, **27** (7), 1669–1687, doi:10.1088/0951-7715/27/7/1669. URL http://dx.doi.org/10.1088/

0951-7715/27/7/1669.

208 Holland, M.P., Rabassa, P., and Sterk., A.E. (2013), On the convergence to an extreme value distribution for non-uniformly hyperbolic dynamical systems, Preprint.

209 Rudin, W. (1987) *Real and complex analysis*, McGraw-Hill Book Co., New York, 3rd edn..

210 Benedicks, M. and Carleson, L. (1985) On iterations of $1 - ax^2$ on $(-1, 1)$. *Ann. of Math. (2)*, **122** (1), 1–25, doi:10.2307/1971367. URL http://dx.doi.org/10.2307/1971367.

211 Guckenheimer, J. and Williams, R. (1979) Structural stability of lorenz attractors. *Publications Mathématiques de l'Institut des Hautes Études Scientifiques*, **50** (1), 59–72, doi:10.1007/BF02684769. URL http://dx.doi.org/10.1007/BF02684769.

212 Hadyn, N.T.A. and Wasilewska, K. (2015) Limiting distribution and error terms for the number of visits to balls in non-uniformly hyperbolic dynamical systems. Preprint.

213 Misiurewicz, M. (1980) Strange attractors for the lozi mappings. *Non Linear Dynamics, R.G. Helleman (ed), New York, The New York Academy of Sciences*.

214 Collet, P. and Levy, Y. (1984) Ergodic properties of the Lozi Mappings. *Commun. Math. Phys.*, **93**, 461–481.

215 Young, L.S. (1985) Bowen-ruelle measures for certain piecewise hyperbolic maps. *Trans. Amer. Math. Soc.*, **287**, 41–48.

216 Chernov, N. and Markarian, R. (2006) Chaotic billiards. *Mathematical surveys and monographs, American Mathematical Society*, **127**.

217 Benedicks, M. and Young, L.S. (2000) Markov extensions and decay of correlations for certain hénon maps. *Asterisque*, **xi** (261), 13–56.

218 Hall, W.J. and Wellner, J.A. (1979) The rate of convergence in law of the maximum of an exponential sample. *Statist. Neerlandica*, **33** (3), 151–154, doi:10.1111/j.1467-9574.1979.tb00671.x. URL http://dx.doi.org/10.1111/j.1467-9574.1979.tb00671.x.

219 Hall, P. (1979) On the rate of convergence of normal extremes. *Journal of Applied Probability*, pp. 433–439.

220 de Melo, W. and van Strien, S. (1993) *One-dimensional dynamics*, Ergebnisse der Mathematik und ihrer Grenzgebiete (3) [Results in Mathematics and Related Areas (3)], vol. 25, Springer-Verlag, Berlin.

221 Katok, A. and Hasselblatt, B. (1996) *Introduction to the Modern Theory of Dynamical Systems*, Encyclopedia of Mathematics and its Applications, Cambridge.

222 Melbourne, I. and Török, A. (2004) Statistical limit theorems for suspension flows. *Israel J. Math.*, **144**, 191–209, doi:10.1007/BF02916712. URL http://dx.doi.org/10.1007/BF02916712.

223 Eagleson, G.K. (1976) Some simple conditions for limit theorems to be mixing. *Teor. Verojatnost. i Primenen.*, **21** (3), 653–660.

224 Rényi, A. (1958) On mixing sequences of sets. *Acta Math. Acad. Sci. Hungar.*, **9**, 215–228.

225 Rényi, A. (1950) Contributions to the theory of independent random variables. (russian. english summary). *Acta Math. Acad. Sci. Hungar.*, **1**, 99–108.

226 Falconer, K. (2003) *Fractal geometry*, John Wiley & Sons, Inc., Hoboken, NJ, 2nd edn., doi:10.1002/0470013850. URL http://dx.doi.org/10.1002/0470013850, mathematical foundations and applications.

227 Chazottes, J.R. and Collet, P. (2010), Poisson approximation for the number of visits to balls in nonuniformly hyperbolic dynamical systems, Preprint arXiv:1007.0171v1. URL http://arxiv.org/abs/1007.0171v1.

228 Wang, Q. and Young, L.S. (2008) Toward a theory of rank one attractors. *Annals of Mathematics*, **167** (2), 349–480.

229 Saltzman, B. (1962) Finite amplitude free convection as an initial value problem - I. *J. Atmos. Sci.*, **19**, 329–341.

230 Tucker, W. (1999) The lorenz attractor exists. *C. R. Acad. Sci. Paris Sér. I Math.*,, **328** (12), 1197–1202.

231 Lucarini, V. and Fraedrich, K. (2009)

Symmetry breaking, mixing, instability, and low frequency variability in a minimal lorenz-like system. *Phys. Rev. E*, **80**, 026 313.

**232** Blender, R. and Lucarini, V. (2013) Nambu representation of an extended lorenz model with viscous heating. *Physica D*, **243**, 86–91.

**233** Afraimovich, V.S. and Pesin, Y.B. (1987) Dimension of lorenz type attractors. *Mathematical physics reviews*, **6**, 169–241.

**234** Zhang, L. (2015) Borel cantelli lemmas and extreme value theory for geometric lorenz models. Preprint.

**235** Baladi, V., Benedicks, M., and Maume-Deschamps, V. (2002) Almost sure rates of mixing for i.i.d. unimodal maps. *Ann. Sci. Ècole Norm. Sup.*, **35**, 126.

**236** Baladi, V., Benedicks, M., and Maume-Deschamps, V. (2003) Corrigendum: Almost sure rates of mixing for i.i.d. unimodal maps, *Ann. Sci. Ecole Norm. Sup.* **35**, (2002), 126. *Ann. Sci. Ècole Norm. Sup.*, **36**, 322.

**237** Alves, J.F. and Araùjo, V. (2003) Random perturbations of nonuniformly expanding maps. *Asterisque*, **286**, 25–62.

**238** Boyarsky, A. and Góra, P. (1997) *Laws of chaos*, Probability and its Applications, Birkhäuser Boston Inc., Boston, MA, doi:10.1007/978-1-4612-2024-4. URL `http://dx.doi.org/10.1007/978-1-4612-2024-4`, invariant measures and dynamical systems in one dimension.

**239** Freitas, J.M., Faranda, D., Guiraud, P., and Vaienti, S. (2014) Sampling local properties of attractors via extreme value theory. *Arxiv preprint arXiv:1407.0412*.

**240** Haydn, N., Nicol, M., Vaienti, S., and Zhang, L. (2013) Central limit theorems for the shrinking target problem. *Journal of Statistical Physics*, **153** (5), 864–887, doi:10.1007/s10955-013-0860-3. URL `http://dx.doi.org/10.1007/s10955-013-0860-3`.

**241** Haydn, N., Nicol, M., Tôrôk, A., and Vaienti, S. (2014) Almost sure invariance principle for sequential and non-stationary dynamical systems. *ArXiv e-prints*.

**242** Saussol, B. (2000) Absolutely continuous invariant measures for multidimensional expanding maps. *Israel Journal of Mathematics*, **116** (1), 223–248, doi:10.1007/BF02773219. URL `http://dx.doi.org/10.1007/BF02773219`.

**243** Knuth, D.E. (1973) *The art of computer programming , Volume 2: Seminumerical Algorithms*, Addison-Wesley.

**244** Faranda, D., Mestre, M., and Turchetti, G. (2012) Analysis of round off errors with reversibility test as a dynamical indicator. *Int. Jou. Bif. Chaos*, **22**.

**245** Aytac, H., Freitas, J., and Vaienti, S. (2014) Hitting time statistics for observations of dynamical systems. *Nonlinearity*, **27**, 2377–2392.

**246** Lasota, A. and Mackey, M. (1994) *Chaos, fractals, and noise: stochastic aspects of dynamics*, Springer.

**247** Freitas, J.M., Faranda, D., Guiraud, P., and Vaienti, S. (2015), Extreme value theory for piecewise contracting maps with randomly applied stochastic perturbations, Preprint arXiv:1501.02913.

**248** Catsigeras, E., Guiraud, P., Meyroneinc, A., and Ugalde, E. (2011) On the Asymptotic Properties of Piecewise Contracting Maps. *ArXiv e-prints*.

**249** Catsigeras, E. and Budelli, R. (2011) Topological dynamics of generic piecewise contracting maps in $n$ dimensions. *International Journal of Pure and Applied Mathematics*, **68**, 61–83.

**250** Coutinho, R., Fernandez, B., Lima, R., and Meyroneinc, A. (2006) Discrete time piecewise affine models of genetic regulatory networks. *Journal of Mathematical Biology*, **52** (4), 524–570.

**251** Strumik, M., Macek, W., and S., R. (2005) Discriminating addittive from dynamical noise for chaotic time series. *Phys. Rev. E*, **72**, 036 219.

**252** Husler, J. (1986) Extreme values of non-stationary random sequences. *Journal of Appl. Prob.*, **23**, 937–950.

**253** Freitas, A., Freitas, J., and Vaienti, S. (2014) Extreme value laws for non stationary processes generated by sequential dynamical systems. *in preparation*.

**254** Conze, J.P. and Raugi, A. (2007) Limit theorems for sequential expanding

dynamical systems on $[0, 1]$, in *Ergodic theory and related fields*, *Contemporary Math.*, vol. 430, Amer. Math. Soc., Providence, pp. 55–91.

255 Baladi, V. (2008) Linear response despite critical points. *Nonlinearity*, **21** (6), T81. URL http://stacks.iop.org/0951-7715/21/i=6/a=T01.

256 Baladi, V. and Smania, D. (2012) Linear response for smooth deformations of generic nonuniformly hyperbolic unimodal maps. *Ann. Sci. ENS*, **45**, 861–926.

257 Abramov, R.V. and Majda, A. (2008) New approximations and tests of linear fluctuation-response for chaotic nonlinear forced-dissipative dynamical systems. *Journal of Nonlinear Science*, **18**, 303–341. 10.1007/s00332-007-9011-9.

258 Lucarini, V. (2008) Response theory for equilibrium and non-equilibrium statistical mechanics: Causality and generalized Kramers-Kronig relations. *J. Stat. Phys.*, **131**.

259 Wang, Q. (2013) Forward and adjoint sensitivity computation of chaotic dynamical systems. *Journal of Computational Physics*, **235** (0), 1 – 13, doi:http://dx.doi.org/10.1016/j.jcp.2012.09.007. URL http://www.sciencedirect.com/science/article/pii/S0021999112005360.

260 Reick, C.H. (2002) Linear response of the Lorenz system. *Phys. Rev. E*, **66**, 036 103.

261 Cessac, B. and Sepulchre, J. (2007) Linear response, susceptibility and resonances in chaotic toy models. *Physica D: Nonlinear Phenomena*, **225** (1), 13–28.

262 Lucarini, V. (2009) Evidence of dispersion relations for the nonlinear response of the Lorenz 63 system. *J. Stat. Phys.*, **134**.

263 Lucarini, V. and Sarno, S. (2011) A statistical mechanical approach for the computation of the climatic response to general forcings. *Nonlin. Proc. Geophys.*, **18**, 7–28.

264 Wouters, J. and Lucarini, V. (2012) Disentangling multi-level systems: averaging, correlations and memory. *Journal of Statistical Mechanics: Theory and Experiment*, **2012** (03), P03 003.

265 Wouters, J. and Lucarini, V. (2013) Multi-level dynamical systems: Connecting the ruelle response theory and the mori-zwanzig approach. *Journal of Statistical Physics*, **151** (5), doi:10.1007/s10955-013-0726-8. URL http://link.springer.com/10.1007/s10955-013-0726-8.

266 Gallavotti, G. and Cohen, E.G.D. (1995) Dynamical ensembles in stationary states. *Journal of Statistical Physics*, **80** (5-6), 931–970.

267 Barreira, L. (1997) Pesin Theory, in *Encyclopaedia of Mathematics* (ed. M. Hazewinkel), Kluwer, pp. 406–411.

268 Barreira, L. and Pesin, Y. (2006) Smooth ergodic theory and nonuniformly hyperbolic dynamics, in *Handbook of dynamical systems. Vol. 1B*, Elsevier, pp. 57–263. With an appendix by O. Sarig.

269 Barreira, L., Pesin, Y., and Schmeling, J. (1999) Dimension and product structure of hyperbolic measures. *Ann. Math.*, **149**, 755–783.

270 Carletti, T. and Galatolo, S. (2006) Numerical estimates of local dimension by waiting time and quantitative recurrence. *Physica A*, **364**.

271 Grassberger, P., Badii, R., and Politi, A. (1988) Scaling laws for invariant measures on hyperbolic and nonhyperbolic atractors. *Journal of Statistical Physics*, **51**, 135–178, doi:10.1007/BF01015324.

272 Lucarini, V., Speranza, A., and Vitolo, R. (2007) Parametric smoothness and self-scaling of the statistical properties of a minimal climate model: What beyond the mean field theories? *Physica D*, **234**, 105 – 123.

273 Albers, D.J. and Sprott, J.C. (2006) Structural stability and hyperbolicity violation in high-dimensional dynamical systems. *Nonlinearity*, **19** (8), 1801–1847. URL http://stacks.iop.org/0951-7715/19/i=8/a=005.

274 Beirlant, J. (2004) *Statistics of extremes: theory and applications*, John Wiley & Sons Inc.

275 J. R. M. Hosking, J.R.M. and Wallis, J.R. (1987) Parameter and quantile estimation for the generalized pareto distribution. *Technometrics*, **29** (3), 339–349.

**276** Arnold, V. and Avez, A. (1968) *Ergodic problems of classical mechanics*, Benjamin New York.

**277** Hasselblatt, B. and Katok, A. (2003) *A first course in dynamics: with a panorama of recent developments*, Cambridge Univ Pr.

**278** Lilliefors, H. (1967) On the Kolmogorov-Smirnov test for normality with mean and variance unknown. *Journal of the American Statistical Association*, **62** (318), 399–402.

**279** Sprott, J. (2003) *Chaos and time-series analysis*, Oxford Univ Pr.

**280** Young, L. (1985) Bowen-Ruelle measures for certain piecewise hyperbolic maps. *Transactions of the American Mathematical Society*, **287** (1), 41–48.

**281** Badii, R. and Politi, A. (1987) Renyi dimensions from local expansion rates. *Physical Review A*, **35** (3), 1288.

**282** Grassberger, P. (1983) Generalized dimensions of strange attractors. *Physics Letters A*, **97** (6), 227–230.

**283** Wolf, A., Swift, J.B., Swinney, H.L., and Vastano, J.A. (1985) Determining Lyapunov exponents from a time series. *Physica D*, **16** (3), 285–317.

**284** Rosenstein, M.T., Collins, J.J., and De Luca, C.J. (1993) A practical method for calculating largest Lyapunov exponents from small data sets. *Physica D*, **65** (1-2), 117–134.

**285** Skokos, C. (2010) The lyapunov characteristic exponents and their computation, in *Lect. Notes Phys., Berlin Springer Verlag*, vol. 790, vol. 790, pp. 63–135.

**286** Kac, M. (1934) On the notion of recurrence in discrete stochastic processes. *Proc. Nat. Acad. Sci. USA*, **20**, 376–379.

**287** Gao, J.B. (1999) Recurrence time statistics for chaotic systems and their applications. *Phys. Rev. Lett.*, **83** (16), 3178–3181.

**288** Hu, H., Rampioni, A., Rossi, L., Turchetti, G., and Vaienti, S. (2004) Statistics of Poincaré recurrences for maps with integrable and ergodic components. *Chaos*, **14**, 160.

**289** Buric, N., Rampioni, A., and Turchetti, G. (2005) Statistics of Poincaré recurrences for a class of smooth circle maps. *Chaos. Soliton. Fract.*, **23** (5), 1829–1840.

**290** Skokos, C., Antonopoulos, C., Bountis, T., and Vrahatis, M. (2002) Smaller alignment index (SALI): Detecting order and chaos in conservative dynamical systems, in *Proc. of 4th GRACM Congress on Computational Mechanics*, vol. 4, Citeseer, vol. 4, pp. 1496–1502.

**291** Skokos, C., Antonopoulos, C., Bountis, T.C., and Vrahatis, M.N. (2004) Detecting order and chaos in Hamiltonian systems by the SALI method. *J. Phys. A-Math. Gen.F*, **37**, 6269.

**292** Skokos, C., Bountis, T.C., and Antonopoulos, C. (2007) Geometrical properties of local dynamics in Hamiltonian systems: The Generalized Alignment Index (GALI) method. *Physica D*, **231** (1), 30–54.

**293** Cincotta, P.M., Giordano, C.M., and Simó, C. (2003) Phase space structure of multi-dimensional systems by means of the mean exponential growth factor of nearby orbits. *Physica D*, **182** (3-4), 151–178.

**294** Goździewski, K., Bois, E., Maciejewski, A.J., and Kiseleva-Eggleton, L. (2001) Global dynamics of planetary systems with the MEGNO criterion. *Astron. Astrophys.*, **378** (2), 569–586.

**295** Chirikov, B. (1969) Research Concerning the Theory of Nonlinear Resonance and Stochasticity Preprint 267. *Institute of Nuclear Physics, Novosibirsk.*

**296** Ott, E. (2002) *Chaos in Dynamical Systems*, Cambridge University Press, New York.

**297** Castillo, E. and Hadi, A.S. (1997) Fitting the generalized pareto distribution to data. *Journal of the American Statistical Association*, **92**, 1609–1620.

**298** Martinez, W. and Martinez, A. (2002) *Computational statistics handbook with MATLAB*, CRC Press.

**299** Faranda, D., Leoncini, X., and Vaienti, S. (2014) Mixing properties in the advection of passive tracers via recurrences and extreme value theory. *Physical Review E*, **89** (5), 052901.

**300** Shin, Y. and Schmidt, P. (1992) The kpss stationarity test as a unit root test. *Economics Letters*, **38** (4), 387–392.

**301** Böhm, R., Auer, I., Brunetti, M., Maugeri, M., Nanni, T., and Schöner, W. (2001) Regional temperature variability in the european alps: 1760–1998 from homogenized instrumental time series. *International Journal of Climatology*, **21** (14), 1779–1801.

**302** Gonzalez, T., Sahni, S., and Franta, W.R. (1977) An efficient algorithm for the kolmogorov-smirnov and lilliefors tests. *ACM Transactions on Mathematical Software (TOMS)*, **3** (1), 60–64.

**303** Ott, E. (1993) *Chaos in Dynamical Systems*, Cambridge Universality Press.

**304** Stanley, E. (1988) *Introduction to phase transitions and critical phenomena*, Oxford University Press.

**305** Guttal, V. and Jayaprakash, C. (2008) Changing skewness: an early warning signal of regime shifts in ecosystems. *Ecology Letters*, **11** (5), 450–460.

**306** Prigent, A. and Dauchot, O. (2005) Transition to versus from turbulence in subcritical Couette flows, in *Laminar turbulent transition and finite amplitude solutions* (eds T. Mullin and R. Kerswell), Springer, pp. 195–219.

**307** Manneville, P. (2012) Coherent structures in transitional plane couette flow. *IMA Journal of Applied Mathematics*, doi:10.1093/imamat/hxs026. URL http://imamat.oxfordjournals.org/content/early/2012/06/08/imamat.hxs026.abstract.

**308** Waleffe, F. (1997) On a self-sustaining process in shear flows. *Phys. Fluids*, **9**, 883–900.

**309** Manneville, P. (2010) *Instabilities, Chaos and Turbulence*, Imperial College Press.

**310** Duguet, Y., Schlatter, P., and Henningson, D.S. (2010) Formation of turbulent patterns near the onset of transition in plane Couette flow. *J. Fluid Mech*, **650**, 119–129.

**311** Tuckerman, L.S. and Barkley, D. (2011) Patterns and dynamics in transitional plane Couette flow. *Phys. Fluids*, **23**.

**312** Manneville, P. and Rolland, J. (2011) On modelling transitional turbulent flows using under-resolved direct numerical simulations: the case of plane Couette flow. *Theoretical and Computational Fluid Dynamics*, **25** (6), 407–420.

**313** Gibson, J. (2007) Channelflow: a spectral Navier–Stokes simulator in C++. *Georgia Institute of Technology*.

**314** Rolland, J. and Manneville, P. (2011) Ginzburg-landau description of laminar-turbulent oblique band formation in transitional plane Couette flow. *Eur. Phys. J. B*, **80**, 529–544.

**315** Philip, J. and Manneville, P. (2011) From temporal to spatiotemporal dynamics in transitional plane Couette flow. *Phys. Rev. E*, **83**.

**316** Dauchot, O. and Manneville, P. (1997) Local versus global concepts in hydrodynamic stability theory. *Journal de Physique II*, **7** (2), 371–389.

**317** Barkley, D. (2011) Modeling the transition to turbulence in shear flows. *J. Phys.: Conf. Ser.*, **318**, 032 001.

**318** Manneville, P. (2005) Modeling the direct transition to turbulence, in *IUTAM Symposium on Laminar-Turbulent Transition and Finite Amplitude Solutions*, Springer, pp. 1–33.

**319** Shepherd, T.G. (2014) Atmospheric circulation as a source of uncertainty in climate change projections. *Nature Geosci*, **7** (10), 703–708. URL http://dx.doi.org/10.1038/ngeo2253.

**320** Smith, R.L. (1988) A counterexample concerning the extremal index. *Adv. in Appl. Probab.*, **20** (3), 681–683, doi:10.2307/1427042. URL http://dx.doi.org/10.2307/1427042.

**321** Lucarini, V. (2012) Stochastic perturbations to dynamical systems: A response theory approach. *Journal of Statistical Physics*, **146** (4), 774–786, doi:10.1007/s10955-012-0422-0. URL http://dx.doi.org/10.1007/s10955-012-0422-0.

**322** Varadhan, S. (1984) *Large Deviations and Applications*, SIAM, Philadelphia.

**323** Touchette, H. (2009) The large deviation approach to statistical mechanics. *Physics Reports*, **478**, 1–69.

**324** Gallavotti, G. (2014) *Nonequilibrium and irreversibility*, Springer, New York.

**325** Kifer, Y. (1992) Averaging in dynamical systems and large deviations. *Inventiones mathematicae*, **110** (1), 337–370, doi:10.1007/BF01231336. URL

http://dx.doi.org/10.1007/
BF01231336.

326 Imkeller, P. and von Storch, J.S. (2001)
*Stochastic Climate Models*, Birkhauser,
Basel.

327 Majda, A. and Timofeyev, I. and Vanden
Eijnden, E. (2001) A mathematical
framework for stochastic climate models.
*Communications on Pure and Applied
Mathematics*, **54**, 5891–5974.

328 Monahan, A.H. and Culina, J. (2011)
Stochastic averaging of idealized climate
models. *Journal of Climate*, **24** (12),
3068–3088,
doi:10.1175/2011JCLI3641.1. URL
http://dx.doi.org/10.1175/
2011JCLI3641.1.

329 Intergovernmental Panel on Climate
Change [Eds.: J. Houghton et al.] (2001)
*IPCC Third Assessment Report: Working
Group I Report "The Physical Science
Basis"*, Cambridge University Press.

330 Intergovernmental Panel on Climate
Change [Eds.: S. Solomon et al.] (2007)
*Climate Change 2007 - The Physical
Science Basis: Working Group I
Contribution to the Fourth Assessment
Report of the IPCC*, Cambridge
University Press, Cambridge, UK and
New York, NY, USA. URL
http://www.worldcat.org/
isbn/0521880092.

331 Buizza, R., Miller, M., and Palmer, T.N.
(1999) Stochastic representation of model
uncertainties in the ecmwf ensemble
prediction system. *Q. J. R. Meteorol. Soc.*,
**125**, 2887–2908.

332 Palmer, T.N. (2001) A nonlinear
dynamical perspective on model error: A
proposal for non-local stochastic-
dynamic parametrization in weather and
climate prediction models. *Q. J. R.
Meteorol. Soc.*, **127**, 279–304.

333 Palmer, T. and Williams, P. (2009)
*Stochastic Physics and Climate
Modelling*, Cambridge Univ. Press.

334 Franzke, C.L.E., O'Kane, T.J., Berner, J.,
Williams, P.D., and Lucarini, V. (2015)
Stochastic climate theory and modeling.
*Wiley Interdisciplinary Reviews: Climate
Change*, **6** (1), 63–78,
doi:10.1002/wcc.318. URL http://
dx.doi.org/10.1002/wcc.318.

335 Chekroun, M., Liu, H., and Wang, S.

(2015) *Approximation of Stochastic
Invariant Manifolds: Stochastic Manifolds
for Nonlinear SPDEs I*, SpringerBriefs in
Mathematics, Springer International
Publishing. URL
https://books.google.com/
books?id=1-7nBQAAQBAJ.

336 Chekroun, M., Liu, H., and Wang, S.
(2015) *Approximation of Stochastic
Invariant Manifolds: Stochastic Manifolds
for Nonlinear SPDEs II*, SpringerBriefs in
Mathematics, Springer International
Publishing.

337 Ruelle, D. (1978) *Thermodynamic
Formalism*, Addison-Wesley, Reading.

338 Daron, J.D. and Stainforth, D.A. (2013)
On predicting climate under climate
change. *Environmental Research Letters*,
**8** (3), 034021. URL
http://stacks.iop.org/
1748-9326/8/i=3/a=034021.

339 Lindenstrauss, J. and L., T. (1996)
*Classical Banach spaces*, Springer, New
York.

340 Ragone, F., Lucarini, V., and Lunkeit, F.
(2014) A new framework for climate
sensitivity and prediction: a modelling
perspective. *ArXiv e-prints*.

341 Lorenz, E.N. (1984) Irregularity: a
fundamental property of the atmosphere.
*Tellus A*, **36** (2), 98–110. URL http:
//www.tellusa.net/index.php/
tellusa/article/view/11473.

342 Lorenz, E.N. (1968) Climatic
determinism. *Meteor. Monog.*, **8**, 1–3.

343 Sornette, D. and Ouillon, G. (2012)
Dragon-kings: Mechanisms, statistical
methods and empirical evidence. *Eur.
Phys. J. Special Topics*, **205**, 1–26.

344 Kaneko, K. (1992) Overview of coupled
map lattices. *Chaos: An Interdisciplinary
Journal of Nonlinear Science*, **2** (3),
279–282.

345 Kaneko, K. (1993) *Theory and
applications of coupled map lattices*,
vol. 12, John Wiley & Son Ltd.

346 Cardy, J. and Ziff, R.M. (2003) Exact
results for the universal area distribution
of clusters in percolation, ising, and potts
models. *Journal of statistical physics*,
**110** (1-2), 1–33.

347 Süveges, M. (2007) Likelihood estimation
of the extremal index. *Extremes*, **10** (1-2),
41–55.